\documentclass[reqno,11pt]{amsbook}

\usepackage{amsmath}
\usepackage{amsthm}
\usepackage{amssymb}
 
\usepackage{amscd}
\usepackage[all]{xy}
\usepackage{verbatim}
\usepackage{stmaryrd}

\usepackage[plainpages=false,colorlinks,hyperindex,pdfpagemode=None,bookmarksopen,linkcolor=red,citecolor=blue,urlcolor=blue]{hyperref}
\usepackage{pdflscape}

\usepackage{minibox}

\swapnumbers
\theoremstyle{plain}	
\newtheorem{theorem}[subsubsection]{Theorem}

\newtheorem*{theorem*}{Theorem}
\newtheorem{proposition}[subsubsection]{Proposition}
\newtheorem*{proposition*}{Proposition}
\newtheorem{lemma}[subsubsection]{Lemma}
\newtheorem*{lemma*}{Lemma}
\newtheorem*{fact*}{Fact}
\newtheorem{corollary}[subsubsection]{Corollary}
\newtheorem*{corollary*}{Corollary}
\newtheorem{conjecture}[subsubsection]{Conjecture}
\newtheorem{discreteseriesconjecture}[subsubsection]{Discrete Series Conjecture}
\theoremstyle{definition}
\newtheorem*{definition}{Definition}

\theoremstyle{remark}
\newtheorem{remark}[subsubsection]{Remark}
\newtheorem*{remarks}{Remarks}

\newtheorem{example}[subsubsection]{Example}

\DeclareFontFamily{OT1}{rsfs}{}
\DeclareFontShape{OT1}{rsfs}{n}{it}{<-> rsfs10}{}
\DeclareMathAlphabet{\mathscr}{OT1}{rsfs}{n}{it}

\newcommand{\expmap}{\mathfrak{exp}}  

\newcommand{\eig}{\mathrm{eig}}

\newcommand{\bruhat}{C^{\infty}_c}
\newcommand{\SSp}{\operatorname{\mathbf{Sp}}}

\newcommand{\limS}{\lim_{S^+}}
\newcommand{\spectrum}{\mathfrak{S}}
\newcommand{\liechecka}{{\mathfrak{a}}^*}
\newcommand{\X}{\mathbf{X}}

\newcommand{\cinfc}{C^\infty_c}

\newcommand{\Vminus}{C_c^\infty(X_\Theta)_\sigma}

\newcommand{\G}{\mathbf{G}}

\newcommand{\Rep}{\mathrm{Rep}}
\newcommand{\aut}{\mathrm{aut}}

\newcommand{\slope}{\mathrm{slope}}

\newcommand{\codim}{\mathrm{codim}}

\newcommand{\su}{{\operatorname{su}}}

\newcommand{\Ad}{\mathrm{Ad}}

\renewcommand{\H}{\mathbf{H}}

\newcommand{\OO}{\mathrm{O}}

\newcommand{\Q}{\mathbf{Q}}
\newcommand{\Z}{\mathbb{Z}}
\newcommand{\C}{\mathbb{C}}

\newcommand{\adele}{{\mathbb{A}_K}}

\renewcommand{\AA}{\mathbf{A}}
\newcommand{\BB}{\mathbf{B}}
\newcommand{\R}{\mathbf{R}}
\newcommand{\CC}{\mathbb{C}}

\newcommand{\GG}{\mathbf{G}}
\newcommand{\Spin}{\mathrm{Spin}}
\newcommand{\HH}{\mathbf{H}}
\newcommand{\LL}{\mathbf{L}}

\newcommand{\MM}{\mathbf{M}}
\newcommand{\NN}{\mathbf{N}}
\newcommand{\PP}{\mathbf{P}}
\newcommand{\un}{\mathrm{un}}
\renewcommand{\P}{\mathcal{P}}
\newcommand{\RR}{\mathbb{R}}
\newcommand{\iR}{{i\mathbb{R}}}

\newcommand{\WW}{\mathbf{W}}
 
\newcommand{\F}{\mathbf{F}} 

\renewcommand{\SS}{\mathbf{S}}

\newcommand{\TT}{\mathbf{T}}
\newcommand{\UU}{\mathbf{U}}
\newcommand{\VV}{\mathbf{V}}
\newcommand{\XX}{\mathbf{X}}
\newcommand{\YY}{\mathbf{Y}}
\newcommand{\ZZ}{\mathbf{Z}}
\newcommand{\QQ}{\mathbb{Q}}

\newcommand{\ab}{{\operatorname{ab}}}

\newcommand{\aff}{{\operatorname{aff}}}

\newcommand{\Ind}{\operatorname{Ind}}
\newcommand{\Hom}{\operatorname{Hom}}
\newcommand{\End}{\operatorname{End}}
\newcommand{\Aut}{{\operatorname{Aut}}}
\newcommand{\Planch}{{\operatorname{Planch}}}
\newcommand{\varchi}{\mathcal{X}}
\newcommand{\GGm}{\mathbf{G}_{\rm m}}
\newcommand{\Gm}{G_{\rm m}}
\newcommand{\GGa}{\mathbf{G}_{\rm a}}
\newcommand{\Ga}{G_{\rm a}}
\newcommand{\supp}{\operatorname{supp}}
\newcommand{\GGL}{\operatorname{\mathbf{GL}}}
\newcommand{\GL}{\operatorname{GL}}
\newcommand{\PPGL}{\operatorname{\mathbf{\mathbf{PGL}}}}
\newcommand{\PGL}{\operatorname{PGL}}
\newcommand{\DD}{\mathbf{D}}
\newcommand{\SSL}{\operatorname{\mathbf{SL}}}
\newcommand{\SL}{\operatorname{SL}}
\newcommand{\ssl}{\mathfrak{sl}}
\newcommand{\Mat}{\operatorname{Mat}}
\newcommand{\Sp}{\operatorname{Sp}}

\newcommand{\GGSp}{\operatorname{\mathbf{GSp}}}

\newcommand{\SO}{{\operatorname{SO}}}
\newcommand{\SSO}{{\operatorname{\mathbf{SO}}}}

\newcommand{\Gal}{\operatorname{Gal}}

\newcommand{\spec}{\operatorname{spec}}
\newcommand{\Vol}{\operatorname{Vol}}
\newcommand{\diag}{{\operatorname{diag}}}
\newcommand{\ev}{\operatorname{ev}}
\newcommand{\im}{\operatorname{im}}
\newcommand{\disc}{{\operatorname{disc}}}
\newcommand{\cont}{{\operatorname{cont}}}
\newcommand{\proj}{{\operatorname{proj}}}

\newcommand{\rk}{{\operatorname{rk}}}
\newcommand{\Id}{\operatorname{Id}}

\newcommand{\pr}{{\operatorname{pr}}}

\newcommand{\val}{{\operatorname{val}}}

\newcommand{\N}{\mathbb{N}}
\newcommand{\cts}{\mathrm{cts}}

\newcommand{\scon}{{\operatorname{sc}}}
\newcommand{\Mac}{{\mathit{Mac}}}
\newcommand{\Unf}{{\operatorname{Unf}}}

 \newcommand{\bluetext}{}
 
\makeatletter
\def\iddots{\mathinner{\mkern1mu\raise\p@
\vbox{\kern7\p@\hbox{.}}\mkern2mu
\raise4\p@\hbox{.}\mkern2mu\raise7\p@\hbox{.}\mkern1mu}}
\makeatother

\begin{document}
\frontmatter
\numberwithin{equation}{section}
\setcounter{tocdepth}{1}
\title{Periods and harmonic analysis on spherical varieties.}
\author{Yiannis Sakellaridis and Akshay Venkatesh}

\begin{abstract}
Given a spherical variety $\XX$ for a group $\GG$ over a non-archimedean local field $k$, the Plancherel decomposition for $L^2(X)$ should be related to ``distinguished'' Arthur parameters into a dual group closely related to that defined by Gaitsgory and Nadler. Motivated by this, we develop, under some assumptions on the spherical variety, a Plancherel formula for $L^2(X)$ up to discrete (modulo center) spectra of its ``boundary degenerations'', certain $\GG$-varieties with more symmetries which model $\XX$ at infinity. Along the way, we discuss the asymptotic theory of subrepresentations of $C^\infty(X)$ and establish conjectures of Ichino--Ikeda and Lapid--Mao.   We finally discuss global analogues of our local conjectures, concerning the period integrals of automorphic forms over spherical subgroups. 
\end{abstract}
\maketitle

\tableofcontents

\mainmatter

\section{Introduction}\label{secintro}

Let $\HH \subset \GG$ be algebraic groups over a field $k$. 
If $k$ is {\em local}, an important problem of representation theory 
is to decompose the $\GG(k)$-action on the space of functions on $\HH\backslash\GG(k)$; if
$k$ is {\em global}, with ring of adeles $\mathbb A_k$,   the study of {\em automorphic period integrals}
\begin{equation} \label{api} \varphi \rightarrow \int_{\HH(k) \backslash \HH(\mathbb A_k)} \varphi,\end{equation}
where $\varphi$ is an automorphic form on $\GG(k) \backslash \GG(\mathbb A_k)$, 
is a central concern of the theory of automorphic forms.  

Our goal (continuing the program of \cite{SaSpc}, \cite{SaSph}) is to {\em relate these questions, and   to formulate a unified framework
in which they can be studied;} the problems discussed here have been previously
been studied largely on a case-by-case basis. 

We shall set up a general formalism in Part 1
 and then, in Parts 2, 3, we give evidence, in the local context,
 that our formalism is indeed the correct one. In Part 4 we formulate the conjectures and give evidence in the global setting.
 The resulting circle of ideas could be understood as part of  a {\em relative Langlands program.}

In most
cases where \eqref{api} is related to an $L$-function, $\HH$ acts with an open orbit on the full flag variety of $\GG$; 
equivalently, a Borel subgroup of $\GG$ acts with an open orbit on $\HH\backslash \GG$.  
This leads us to the starting point of the theory, spherical varieties.

\subsection{}

Let $\GG$ be a reductive group and $\XX$ a $\GG$-variety. In this paper the group will always be split over the base field $k$, and the base field will be of characteristic zero. 

The variety $\XX$ is called {\em spherical} if it is normal and a Borel subgroup $\BB \subset \GG$ acts with a Zariski dense orbit.  
It is a remarkable fact that spherical varieties have a uniform structure theory
and are classified by combinatorial data.  They include all {\em symmetric} varieties. 

We will consider the questions formulated above in the case when $\HH\backslash \GG$ is a spherical variety under $\GG$. 
Our goal will be to formulate
conjectural answers  in terms of the data attached to the spherical variety $\HH\backslash\GG$.

Throughout we will use the convention of denoting with boldface
letters $\GG, \XX$ algebraic groups or algebraic varieties, and by $G, X, \dots$ their points
over a local field $k$. 

\subsection{} \label{illustrations}

 As illustrations, we use the following classes of spherical varieties:
 \begin{itemize}
 \item[-] {\em Symmetric}:  stabilizers on $\XX$ are fixed
 points of an involution on $\GG$.
 \item[-] {\em Gross--Prasad}: $\GG = \GG_n \times \GG_{n+1}$ acting (by right and left multiplication)
 on $\XX = \GG_{n+1}$, where $\GG_n = \SSO_n$ or $\GGL_n$. 
 \item[-] {\em Whittaker.}  Here $\XX = \UU\backslash \GG$, where $\UU$ is the maximal unipotent subgroup; instead of functions on $X$ we consider sections of a line bundle defined by a nondegenerate additive character of $U$. (This does not fall strictly in the framework of spherical varieties, but nonetheless our results and methods apply unchanged to this case -- cf.\ \S \ref{ssWhittakertype}.) 
 \end{itemize}

In the Gross--Prasad and Whittaker cases, it is conjectured that the global automorphic period
is related to special values of $L$-functions; this is also believed in many, but not all, symmetric cases.

\subsection{}
We formulate our main local and global conjectures, and then discuss our results  in \S \ref{intro-results}. We point the reader to Part 4 for more precise formulations of the conjectures. 

 Our conjectures are phrased in terms of a dual group $\check G_X$ attached
to the spherical variety $\XX$. This is inspired and motivated by the work of Gaitsgory and Nadler \cite{GN}.  We define the root datum  of $\check G_X$ in \ref{dualgroup}; the dual group comes equipped with a canonical morphism of the distinguished Cartan subgroup of $\check G_X \times \SL_2$ to the distinguished Cartan subgroup of $\check{G}$.  A {\em distinguished morphism} is an  extension of this to a map \begin{equation} \label{Ext} \check G_X \times \SL_2 \rightarrow \check {G};\end{equation} 
that satisfies a certain constraint on root spaces formulated in \S \ref{distinguishedmorphisms}. 
We conjecture that {\em such an extension always exists},  
and prove it (for most spherical varieties, termed ``wavefront'') {\em assuming}
that the Gaitsgory-Nadler construction satisfies certain natural axioms e.g.\  compatibility   with boundary degeneration and parabolic induction. We should note here that our definition of the dual group leaves out some varieties -- for instance, the $\GGL_n$-variety of non-degenerate quadratic forms in $n$ variables. Our harmonic-analytic results still hold in this case, but formulating a Langlands-type conjecture about the spectrum is a very interesting problem whose answer we do not know.

What is important for applications is that one can rapidly compute \eqref{Ext} in any specific case;
for example, we give a table of rank one cases in Appendix \ref{primerankone}.

Now let $k$ be a local field. An Arthur parameter for $\GG$ is a homomorphism
 $\phi: \mathcal{L}_k \times \SL_2 \rightarrow \check{G}$, such that the image of the first factor is bounded and the restriction to the second factor is algebraic. Here $\mathcal L_k$ is the Weil-Deligne group of $k$ (Weil group in the archimedean case). We say that $\phi$ is $\XX$-distinguished if it factors through 
 a map $\tilde{\phi}: \mathcal{L}_k \longrightarrow \check{G}_X $, i.e.
 $\phi(w,g) = \rho (\tilde{\phi}(w), g)$, where $\rho$ is the map \eqref{Ext}.

 \begin{conjecture} \label{spectrum} The support of the Plancherel measure for $L^2(X)$, as a $G$-representation, is contained in the union of Arthur packets
 attached to  $\XX$-distinguished Arthur parameters.  \end{conjecture}
 
In fact, we may enunciate a more precise conjecture, predicting a direct integral decomposition:
\begin{equation}\label{Plancherel-intro}L^2(X) = \int_\phi \mathcal H_\phi \mu(\phi),\end{equation}
where $\psi$ ranges over $\check G_X$-conjugacy classes of $\XX$-distinguished Arthur parameters, the Hilbert space $\mathcal H_\phi$ is isotypic for a sum of representations belonging to the Arthur packet corresponding to $\phi$, and the measure $\mu(\phi)$ is absolutely continuous with respect to the natural ``Haar'' measure on Arthur parameters.
The sharpened conjecture implies that the unitary irreducible $G$-representations that occur as {\em subrepresentations} of $L^2(X)$ -- the so-called relative discrete series -- 
 are all contained in Arthur packets arising from {\em elliptic} parameters $\mathcal{L}_k \longrightarrow \check{G}_X$, i.e.\ maps that do not factor through a proper Levi subgroup. 
  
 \begin{example} Let $\VV$ be a $2n$-dimensional vector space over $k$, 
 $\GG = \GGL(V)$, and $\XX$ the space of alternating forms on $\VV$. Then $\XX$ is a spherical $\GG$-variety; the group $\check G_X$ is isomorphic to $\GL_n$, and the map
 $$\check G_X \times \SL_2 = \GL_n \times \SL_2 \longrightarrow \GL_{2n}$$
 is the tensor product of the standard representations. The content of the Conjecture is then that the unitary spectrum of $L^2(X)$ are precisely the {\em Speh representations}
 $J(2, \sigma)$, where $\sigma$ is a tempered representation of $\GL(n, k)$; moreover, such a represention embeds into $L^2(X)$ precisely if $\sigma$ is discrete series. We point to the work of Offen--Sayag \cite{OS} for work in this direction. \end{example}  
\begin{example}
 For many low-rank spherical varieties, Gan and Gomez have proven this conjecture recently using the theta correspondence, \cite{GG}.
\end{example}

It is desirable to refine Conjecture \ref{spectrum} to a precise Plancherel formula. 
 We will discuss a more precise version of this conjecture in section \ref{sec:localconj}.

 \subsection{}
 We now discuss its relationship with a global conjecture about the Euler factorization of periods of automorphic forms (cf.\ Section \ref{sec:globalconj}).  Let $K$ be a global field, with ring of adeles $\adele$. Let $\XX=\HH\backslash \GG$ be a spherical variety defined over $K$, and let $\pi=\otimes\pi_v\hookrightarrow C^\infty([\GG])$ (where $[\GG]=\GG(K)\backslash \GG(\adele)$) be an irreducible automorphic representation of $\GG$. Under some assumptions on $\HH$ (multiplicity-one is clearly sufficient, but not necessary as the work of Jacquet \cite{JEuler} shows),  it is expected that the \emph{period integral} against Tamagawa measure on $[\HH]$:
\begin{equation}\label{periodintegral}
 I_H:\phi\mapsto \int_{[\HH]} \phi(h) dh
\end{equation}
(whenever it makes sense) is an \emph{Eulerian} functional on the space of $\pi$, i.e.\ a pure tensor in the restricted tensor product:
$$\Hom_{\HH(\adele)} (\pi,\CC) = \otimes'_v \Hom_{\HH(K_v)}(\pi_v,\CC).$$ 
We will state a conjecture in the multiplicity-free case, for a treatment of Jacquet's example cf.\ \cite{FLO}.

The refined Gross--Prasad conjecture by Ichino and Ikeda \cite{II} gives an explicit Euler factorization of this functional -- or rather of the hermitian form $\mathcal P^\Aut := |I_H|^2$ -- in the case of $\GG=\SSO_n\times \SSO_{n+1}$, $\HH=\SSO_n$, at least when $\pi$ is tempered. 
They conjecture that (up to explicit, rational global constants) the form $\mathcal P^\Aut$
factorizes as products of the local Hermitian forms
$$\theta_{\pi_v}: (u,\tilde u) \in \pi_v \times \overline{\pi_v} \mapsto \int_{\HH(K_v)} \left< \pi_v(h) u, \tilde u\right> dh$$

An important motivating observation for us was that these Hermitian forms {\em play an important role in the Plancherel formula for the space $L^2(\XX(K_v))$:} the Hermitian form $\theta_{\pi_v}$ can be viewed, via Frobenius reciprocity, as a morphism $\pi_v\otimes\overline{\pi_v}\to C^\infty(\XX(K_v))\otimes C^\infty(\XX(K_v))$. Its dual, composed with the unitary pairing $\pi_v\otimes\bar\pi_v\to \CC$, defines an invariant hermitian form $H_{\pi_v}$ on $C_c^\infty(\XX(K_v))$. The $L^2$ inner product of functions on $X$ is the integral of these Hermitian forms against the {\em standard} Plancherel measure on the unitary dual $\widehat G$.

Let us note here an important subtlety of the Plancherel formula. In general,
the theory of unitary decomposition associates to $L^2(X)$ (here $X=\XX(K_v)$ etc.)
only a measure {\em class} on the unitary dual $\widehat{G}$. To choose a specific measure $\mu$
in this class is essentially equivalent to fixing an embedding
$\pi \otimes \overline{\pi} \hookrightarrow C^{\infty}(X \times X)$,
for almost all $\pi$ in the support of that measure. In the {\em group case} 
($\XX = \HH^\diag \backslash \HH \times \HH$) there is a canonical normalization
of such an embedding, coming from the theory of matrix coefficients. 

In general,
there is no corresponding normalization; however, our local conjecture gives a natural candidate
for  $\mu$:
If the Plancherel formula can be written in terms of parameters into $\check G_X$, then the measure that one would use for the Plancherel decomposition of $L^2(G_X)$, where $G_X$ is the split group with dual $\check G_X$, seems to be a natural choice. This measure was (conjecturally) described in \cite{HII} in terms of Langlands parameters, and it only depends on the parameter, up to a rational factor that may show up for ramified representations of exceptional groups. Since our global conjecture is only up to a rational factor, this ambiguity does not concern us here, and we can think of Plancherel measure on $G_X$ as a measure on the set of bounded Langlands parameters into $\check G_X$.

Fixing this measure {\em gives rise to normalized embeddings of $\pi \times \overline{\pi}$
into $C^{\infty}(X \times X)$}, for almost every $\pi$ in the support of Plancherel measure (and, by some continuity property, for all), and by evaluating at the identity we get a $H$-biinvariant Hermitian form $\mathcal P_v^\Planch$ on $\pi$. We conjecture that this is ``the correct normalization for global applications'', i.e.\ whenever $\mathcal P^\Aut$ is Eulerian these local forms are the correct generalization of the forms $\theta_{\pi_v}$ of Ichino and Ikeda: 

\begin{equation}
 \mathcal P^\Aut = q \prod'_v \mathcal P_v^\Planch,
 \end{equation}
 where $q$ is a nonzero rational factor that we don't specify. The Euler product is typically non-convergent, and the product of all but finitely many factors should be interpreted as a product/quotient of special values of $L$-functions. 
 
 There are many assumptions for this conjecture, besides the local Conjecture \ref{spectrum}, such as a multiplicity-one assumption and the assumption that the space of the automorphic representation in question is the one corresponding to an ``$\XX$-distinguished (global) Arthur parameter''. The existence of these parameters is of course highly conjectural, but in certain cases there are more down-to-earth versions of the conjecture that one can formulate. We point the reader to Section \ref{sec:globalconj}.

Clearly, our global conjecture lacks the precision of \cite{II}, and should be considered as a guiding principle for the time being. In any case, it provides access to the mysterious link between global periods and automorphic $L$-functions, via a computation of local Plancherel 
measures that was performed in \cite{SaSph}.

\subsection{} \label{intro-results}
This paper is divided into four main parts. 
All four bear on the main conjecture, but the details of individual parts are to a large extent independent
and can be read separately.

Some of the main results are  Proposition \ref{dualgrouprootsystem}/Theorem \ref{sl2} (identification of dual group),  Theorem \ref{mainasthm} (asymptotics of representations, implying finite multiplicity), Theorem \ref{propIIconj} (Ichino--Ikeda conjecture), 
Theorem \ref{finiteds} (finiteness of discrete series), Theorem \ref{Bernsteinmap} (existence of scattering morphisms), 
Theorem \ref{tilingtheorem} (abstract scattering theorem), Theorem \ref{advancedscattering} (in many cases, a complete description of scattering) Theorem \ref{explicitBernstein} (Plancherel decomposition in terms of ``normalized Eisenstein integrals'') and Theorem \ref{holdsbyunfolding} (compatibility of the global conjecture with ``unfolding'').

 Let $k$ be a local non-archimedean field. 
Practically all of our results are obtained under the assumption that $\XX$ is ``wavefront'' (see \S \ref{invariants} for the definition).  This includes the vast majority of spherical varieties (e.g., 
in Wasserman's tables \cite{Wa} of rank $2$ spherical varieties, only three fail to be wavefront),
and in particular covers the Whittaker, Gross--Prasad, and all symmetric cases.

\begin{enumerate}
\item {\em Part 1 ( \S \ref{sec:review} and \S \ref{sec:dualgroupproofs}). Dual groups of spherical varieties.}

   It is primarily concerned
 with defining the dual group $\check G_X$ and establishing -- as far as possible -- the existence
 of the morphism $\check G_X \times \SL_2 \rightarrow \check G$.  As mentioned, 
 we prove (Theorem \ref{sl2}) that this morphism exists {\em assuming} the compatibility
 of the Gaitsgory-Nadler construction with certain natural operations, such as boundary degeneration and parabolic induction; 
 a by-product  of this proof is an identification of the root system of the Gaitsgory-Nadler dual group.

An important feature of $\check G_X$ is its relation to the geometry of $\XX$ at $\infty$. 
To each conjugacy class $\Theta$ of parabolic subgroups of $\check G_X$, 
we associate a spherical variety $\XX_{\Theta}$ (which we call a ``boundary degeneration'') under $\GG$; it models
the structure of a certain part of $\XX$ at $\infty$. 
 The dual group $\check G_{X_{\Theta}}$ to $\XX_{\Theta}$
is isomorphic to a Levi subgroup of a parabolic subgroup in the class $\Theta$.

 The reader more interested in local or global theory could skip most of this section,
 reading only the parts on the boundary degenerations, and perhaps
 glancing at the table of examples in the Appendix.

\item {\em Part 2: Asymptotics and the Ichino--Ikeda conjecture  (\S \ref{sec:asymptotics} -- \ref{sec:sttempered}).}

 We verify (\S \ref{sec:asymptotics}) that the multiplicity of any irreducible $G$-representation in $C^\infty(X)$ is finite; we compute (also \S \ref{sec:asymptotics}) the asymptotic behavior\footnote{Such results on asymptotics, but expressed in the more traditional language of Jacquet modules, were proven for symmetric varieties by Lagier \cite{Lagier}, and independently by Kato and Takano \cite{KatoTakano}.} of ``eigenfunctions'' (i.e., functions on 
$X$ whose $G$-span is of finite length). 

The latter result is naturally expressed in terms of an ``asymptotics'' map
\begin{equation}  \label{intro-asymp}  e_\Theta:  C_c^{\infty}(X_\Theta) \longrightarrow C_c^{\infty}(X),\end{equation}  
 see Theorem \ref{mainasthm}.

We remark that these results are corollaries to an understanding of the geometry of $X$ at $\infty$; this
 understanding plays a fundamental role throughout the entire paper.\footnote{Recently Bezrukavnikov and Kazhdan \cite{BK} have given a geometric analysis of Bernstein's second adjunction, which is closely related to our analysis in the special case where $\X$ is the group
 variety $\Delta \G \backslash \G \times \G$. They use, in particular, the structure of the wonderful compactification of $\G$ itself; cf.\ also \S \ref{sssupport}.} 

 By elementary methods, we are able in \S \ref{sec:sttempered} to completely describe a Plancherel formula (Theorem \ref{ST:Plancherel}) 
for ``strongly tempered varieties''; this is a condition that implies $\check G_X = \check G$ and
includes the Gross--Prasad and Whittaker cases, although {\em not} most symmetric cases.  
This gives, in particular, a simple derivation of the Whittaker-Plancherel formula for $p$-adic groups\footnote{While our paper was being written, a complete description of the Whittaker-Plancherel formula was obtained by Delorme \cite{DeWhi}. Our proof is rather different.} (more precisely: a simple reduction to the usual Plancherel formula). 

Using these results, we verify conjectures of Ichino--Ikeda (Theorem \ref{propIIconj})
and Lapid--Mao (Corollary \ref{LMconj}). 
We mention only the former: if $(H, G)$ is as in the Gross--Prasad conjecture 
then for any tempered representation $\Pi$ of $G$, the form
$$v \otimes v' \mapsto \int_{h \in H} \langle h v, v'\rangle, \ \ v \in \Pi, v' \in \tilde{\Pi}$$
on $\Pi \otimes \tilde{\Pi}$ is nonvanishing {\em if and only if}
$\Pi$ is $H$-distinguished.     

\item {\em Part 3: Scattering theory. ( \S \ref{sec:discrete} -- \S \ref{sec:explicit})}

This is the core of the paper; the results of this section are summarized on page \pageref{pennylane}.

The set of conjugacy classes of $\XX$-distinguished Arthur parameters is partitioned 
into subsets indexed by conjugacy classes of Levi subgroups of $\check G_X$. We give evidence in 
\S \ref{sec:Bernstein} -- \S \ref{sec:explicit} that the unitary spectrum of $L^2(X)$
has a corresponding structure; in the most favorable cases (that is to say: satisfying an easy-to-check combinatorial criterion)
this amounts to a Plancherel formula modulo the knowledge of discrete (modulo center) series for $\XX$ and all its boundary degenerations $\XX_\Theta$. 
(The preceding sections cover preliminary ground: \S \ref{sec:discrete} discusses a somewhat subtle issue concerning the fact that one can have continuous families of relative discrete series, and \S \ref{sec:linearalgebra} contains some lemmas in linear algebra that are necessary to formulate the scattering arguments).

 The main tool is ``scattering theory,'' which relates the spectrum of a space and its boundary. 
We obtain in Theorem \ref{Bernsteinmap} a canonical $G$-equivariant map
\begin{equation} \label{l2-asymp-intro} L^2(X_{\Theta}) \stackrel{\iota_{\Theta}}{\longrightarrow} L^2(X). \end{equation}
This should be viewed as a unitary analog of the smooth asymptotics map (\ref{intro-asymp}). We call this map the ``Bernstein map'',\footnote{Unfortunately, the term ``Bernstein map'' is used in \cite{BK} for the smooth asymptotics map (\ref{intro-asymp}).} because its existence is essentially equivalent to an unpublished argument of Joseph Bernstein, which proves that the continuous part of the Plancherel formula for $X$ should resemble the Plancherel formula for the boundary degenerations of $X$ as one moves towards infinity.
Conjecturally, it corresponds to the evident map on Arthur parameters 
induced by $\check G_{X_\Theta} \hookrightarrow \check G_X$.

Let $L^2(X)_{\Theta}$ be the image of $L^2(X_\Theta)_\disc$ under $\iota_{\Theta}$.   We conjecture
that $L^2(X)_{\Theta} = L^2(X)_{\Theta'}$ when $\Theta, \Theta'$ are associate.
In favorable cases we are able to prove this in Theorem \ref{tilingtheorem}, and, in fact, precisely describe the kernel of the 
morphism $\bigoplus \iota_{\Theta}$.

Finally we discuss in \S \ref{sec:explicit} an ``explicit description'' of $\iota_{\Theta}$ in terms
of Mackey theory. The goal here, which we only partially achieve, is to describe the morphisms $\iota_{\Theta}$ in terms of explicit intertwining operators, commonly refered to as ``Eisenstein integrals''. In some combinatorially favorable cases we fully achieve this goal, including many symmetric varieties.

\item {\em Part 4.  Conjectures.}  
 
In this part we formulate the local and global conjectures discussed above, and give some evidence for the global ones. The formalism here relies on the local and global Arthur conjectures \cite{Arthur-unipotent}. The local Conjecture \ref{localconjecture-weak} states that the unitary representation $L^2(\XX(k))$ (where $k$ is a local field) admits a direct integral decomposition in terms of $\XX$-distinguished Arthur parameters (and the natural class of measures on them). A finer version (Conjecture \ref{localconjecture-refined}) introduces a notion of ``pure inner form'' for a spherical variety $\XX$, inspired from the relative trace formula and the local Gan--Gross--Prasad conjectures \cite{GGP}.

In Section \ref{sec:globalconj} we formulate the global conjecture on Euler factorization of period integrals (under several assumptions). Finally, in Section \ref{sec:examples} we prove the global conjecture in some cases: periods of principal Eisenstein series, Whittaker periods for $\GGL_n$, and all periods that ``unfold'' to Whittaker periods for $\GGL_n$. Of course, explicit Euler factorizations for these periods have been known in the past; what we do is verify that the local factors are the ones predicted by our conjectures, which are related to local Plancherel formulas.  Much of this is known to experts, and our goal is in part
to express this computation of local factors  in the language of this paper.

\end{enumerate}

  \subsection{Proofs}
  
We outline the ideas behind the results at the heart of this paper, the local Plancherel formula developed in Parts 2 and 3. As we have already mentioned, the basic ingredient in many of the proofs is a good understanding of the geometry of $\XX$ at $\infty$.   We will give briefly some examples of the type of ideas that enter. 

\subsubsection{Geometry at $\infty$}
    
A critical fact in the theory of spherical varieties is that there exists a parabolic $\PP(\XX)^-$ (the notation is such because it is in the opposite class of parabolics to one we will denote by $\PP(\XX)$) so that 
\emph{ the geometry of $\XX$
at $\infty$ is modelled by a torus bundle $\mathbf{Y}$ over $\PP(\XX)^-\backslash \GG$.}

  For example,
the hyperboloid $x^2-y^2-z^2 = 1$ is spherical under the group $\mathrm{SO}_3$;
at $\infty$ it becomes asymptotic to the cone $x^2 - y^2 - z^2 = 0$, which 
is a line bundle over the flag variety $\mathbb{P}^1$.  
In terms of the varieties $\XX_{\Theta}$ previously mentioned, 
the torus bundle $\YY$ is obtained by taking for $\Theta$ the class of Borel subgroups in the dual group $\check G_X$. 
 
 In fact, this is an overly simplified view of the geometry of $\XX$ at $\infty$; more accurately, the geometry of $\XX$ at $\infty$ is modelled 
 by the so-called {\em wonderful compactification} $\overline{\XX}$. 
 The $\GG$-orbits on $\overline{\XX}$ are canonically in correspondence with conjugacy classes of parabolic subgroups of $\check{G}_X$; for each such conjugacy class $\Theta$ we call ``$\Theta$-infinity'' the corresponding orbit at infinity, and define the variety $\XX_{\Theta}$ as (the open $\GG$-orbit in) the {\em normal bundle} to $\Theta$-infinity.
  In particular, there is -- in the sense of algebraic geometry -- a degeneration of $\XX$ to $\XX_{\Theta}$; in intuitive terms, $\XX_{\Theta}$
 models a part of the geometry of $\XX$ at $\infty$. The $\GG$-variety $\XX_\Theta$ is ``simpler'' than $\XX$ in a very important way: it carries the additional action of a torus $\AA_{X,\Theta}$, generated from the actions of the multiplicative group on the normal bundles to all $\GG$-stable divisors containing $\Theta$-infinity.

\subsubsection{Geometry at $\infty$ over a local field} \label{sss:geominftyintro}
The discussion above has the following consequence for points over a local field $k$:

Let $J$ be an open compact subgroup of $G = \GG(k)$. 
We construct a canonical identification of a certain subset of $X/J$ with a certain subset of $Y/J$, mirroring the fact that $X$ approximates $Y$ at $\infty$.   This fact remains true (with different subsets) if we replace $Y$ by any $X_{\Theta}$. We call this identification the ``exponential map'', because it is in fact induced by some kind of exponential map between the ($k$-points of the) normal bundle and the variety.
 
In the case of the hyperboloid we may describe this as follows:  Write $X = \{ (x,y,z) \in k^3: x^2-y^2-z^2=1\}$ and $ Y = \{ (x,y,z) \in k^3: x^2-y^2-z^2=0\}.$
We declare two $J$-orbits $x J \subset X,  \ \  y J \subset Y$ to be {\em $\varepsilon$-compatible}
if there exists $x' \in xJ, y' \in yJ$ that are at distance $< \varepsilon$ for the nonarchimedean metric on $k^3$. Then, for all sufficiently small $\varepsilon$ (this notion depending on $J$), there exist compact sets $\Omega_X \subset X, \Omega_Y \subset Y$
so that the relation of $\varepsilon$-compatibility gives a bijection  between $(X-\Omega_X)/J$
and $(Y-\Omega_Y)/J$.

 \subsubsection{Asymptotics of eigenfunctions}  \label{sss:asympintro} 
Call a (smooth) function on $X$ or $Y$ an {\em eigenfunction} if its translates under $G = \GG(k)$ span a $G$-representation of finite length.    Then the fundamental fact of interest to us is that, 
for every ($J$-invariant) eigenfunction $f$ on $X$, there exists an eigenfunction $f_Y$ on $Y$ so that
``$f$ is asymptotic to $f_Y$'': that is to say, $f$ and $f_Y$ are identified under
the ``partial bijection'' (exponential map) between $X/J$ and $Y/J$.  In fact, this is a fact that does not require admissibility or finite length (although we only need the finite length case for the Plancherel formula): there is a $G$-morphism: $e_\Theta^*: C^\infty(X)\to C^\infty(X_\Theta)$ (the dual of (\ref{intro-asymp}) with the property that functions coincide with their images in  neighborhoods of $\Theta$-infinity identified via the exponential map.

\subsubsection{The argument of Bernstein}

Now we turn our attention to the Plancherel decomposition: its existence and uniqueness is guaranteed by theorems involving $C^*$-algebras, and by \cite{BePl} it is known to be supported on ``$X$-tempered'' representations. This means that the norm of every Harish-Chandra--Schwartz function $f$ on $X$ admits a decomposition:
$$\Vert f\Vert^2_{L^2(X)} = \int_{\hat G} H_\pi(f) \mu(\pi),$$
where the positive semi-definite hermitian forms\footnote{We feel free to write $H(f):=H(f,f)$ for a hermitian form $H$.} $H_\pi$ factor through a quotient $G$-space of finite length, isomorphic to a number of copies of $\pi$. In other words, $H_\pi(f)$ can be written as a finite sum of terms of the form $|l_i^*(f)|^2$, where $l_i:\pi\to C^\infty(X)$ is an embedding with tempered image.

Given the theory of the exponential map, explained above, and the asymptotic theory for such embeddings, if the support of $f$ is concentrated close to``$\Theta$-infinity'' then $f$ can be identified with a function $f'$ on $X_\Theta$ and the expression $|l_i^*(f)|^2$ can be identified with the same expression for some $l_i':\pi\to C^\infty(X_\Theta)$. In other words, the   hermitian form $H_\pi$ give rise to a hermitian form $H_\pi'$ on functions on $X_\Theta$; in precise
terms this is simply the pullback of $H_{\pi}$ under the map $C^{\infty}_c(X_{\Theta}) 
\rightarrow C^{\infty}_c(X)$ of \eqref{intro-asymp}.

Could the forms $H_{\pi}'$ appear in the  Plancherel formula for $L^2(X_\Theta)$? Almost, but not quite. The reason is that the Plancherel formula for $X_\Theta$ is invariant under the additional torus symmetry group $A_{X,\Theta}$ of $X_\Theta$, and this is not, in general, the case with the asymptotic forms $H_\pi'$. There is, however, an elementary way to extract their ``$A_{X,\Theta}$-invariant part'' $H_\pi^\Theta$ (possibly zero), and then \emph{one obtains a Plancherel decomposition for $L^2(X_\Theta)$}:
$$\Vert f'\Vert^2_{L^2(X_\Theta)} = \int_{\hat G} H_\pi^\Theta(f') \mu(\pi).$$

We show that this fact is equivalent to the existence of ``Bernstein maps'': 
$$ \iota_\Theta: L^2(X_\Theta)\to L^2(X), $$
which are characterized by the fact that they are ``very close'' to the smooth asymptotics maps (\ref{intro-asymp}) for functions supported ``close to $\Theta$-infinity''.

For readers familiar with the Plancherel decomposition of the space of automorphic functions, we mention that the analog of this map in that case arises as follows: one decomposes a function $f\in L^2(N(\adele)A(K)\backslash \PGL_2(\adele))$ (this is the analog of $X_\Theta$ in that case) as an integral of (unitary) $A(\adele)$-eigenfunctions, and then $\iota_\Theta f$ will be the corresponding integral of Eisenstein series, i.e., replace each eigenfunction on $N(\adele) A(K) \backslash \PGL_2(\adele)$ by the corresponding Eisenstein series on $\PGL_2(K) \backslash \PGL_2(\adele)$. Notice that this is taken on the tempered line, without any discrete series contributions. 

In our case we will define $\iota_\Theta$ in a more abstract way, and only later \S \ref{sec:explicit} (and under additional conditions) will we identify it in terms of explicit morphisms (``normalized Eisenstein integrals'') analogous to the Eisenstein series. Thus, in the language often used in the literature of harmonic analysis on symmetric spaces, the Bernstein maps can be identified with ``normalized wave packets''.

\subsubsection{Scattering}
 
The construction of Bernstein maps also implies that their sum, restricted to discrete spectra:
$$\sum \iota_{\Theta,\disc}: \bigoplus_{\Theta\subset\Delta_X}  \iota_\Theta(L^2(X_\Theta)_\disc) \to L^2(X)$$
is surjective.

For a full description of $L^2(X)$ in terms of discrete spectra there remains to understand the kernel of this map. Based on the dual group conjecture, we expect the images of $L^2(X_\Theta)_\disc$ and $L^2(X_\Omega)_\disc$ are orthogonal if $\Theta$ and $\Omega$ are not $W_X$-conjugate, and coincide otherwise. We prove this under some combinatorial condition
(``generic injectivity'':  \S \ref{genericinjectivity}), which is easy to check and is known to hold, at least, for all symmetric varieties. 

More precisely, it is easy to show, first, that the images of $L^2(X_\Theta)_\disc$ and $L^2(X_\Omega)_\disc$ are orthogonal unless $\Theta$ and $\Omega$ are of the same size (i.e.\ the corresponding orbits of the wonderful compactification are of the same dimension). The combinatorial condition, on characters of the associated boundary degenerations, is used to rule out the possibility that the images are non-orthogonal when $\Theta$ and $\Omega$ are non-conjugate. Finally, a delicate analytic argument shows that when they are conjugate the images of $L^2(X_\Theta)_\disc$ and $L^2(X_\Omega)_\disc$ have to coincide.
This is encoded by certain maps $L^2(X_{\Theta}) \rightarrow L^2(X_{\Omega})$ whenever $\Theta, \Omega$
are $W_X$-conjugate, the so-called \emph{scattering morphisms}. For a more detailed introduction to the results of scattering theory, we point the reader to Section \ref{pennylane}.

\subsection{Notation and assumptions}\label{notation} 

We have made an effort to define or re-define most of our notation locally, in order to make the paper more readable. An exception are the notions and symbols introduced in section \ref{sec:review}. We note here a few of the conventions and symbols which are used throughout: We fix a locally compact $p$-adic field $k$ in characteristic zero, with ring of integers $\mathfrak o$; we denote varieties over $k$ by bold letters and the sets of their $k$-points by regular font.  For example, if $\YY$ is a $k$-variety, we denote $\YY(k)$ by $Y$ without special remark. On the other hand, for complex varieties (such as dual groups or character groups) we make no notational distinction between the abstract variety and its complex points, and use regular font for both.

We always use the words ``morphism'' or ``homomorphism'' in the appropriate category (which should be clear from the context), e.g.\ for topological groups a ``homomorphism'' is always continuous, even if not explicitly stated so.
We denote by $\varchi(\MM)=\Hom(\MM,\GGm)$ the character group of any algebraic group $\MM$ and for every finitely generated $\Z$-module $R$ we let $R^*=\Hom(R,\Z)$. Normalizers are denoted by $\mathcal N$, or $\mathcal N_\GG$ when we want to emphasize the ambient group; this is not to be confused with the notation $N_\ZZ\YY$, which denotes the normal bundle in a variety $\YY$ of a subvariety $\ZZ$. 

We denote throughout by $\GG$ a connected, reductive, \emph{split} group over $k$, and by $\XX$ a homogeneous, quasi-affine spherical variety for $\GG$. For most of the paper, we assume this variety to be wavefront (cf. \S \ref{invariants} -- see the list that follows for a full set of assumptions). We fix\footnote{\label{footnoteBorel}It is actually for convenience of language that we fix a Borel subgroup; one could adopt the language pertinent to ``universal Cartan groups'' and show that all constructions, such as the quotient torus $\AA_X$, are unique up to unique isomorphism because of the fact that Borel subgroups are self-normalizing and conjugate to each other.} a Borel subgroup $\BB$ and denote its Zariski open orbit on $\XX$ by $\mathring{\XX}$. Parabolics containing $\BB$ will be called ``standard'', the unipotent
radical of $\BB$ will be denoted by $\UU$ and the reductive quotient of $\BB$ will be denoted by $\AA$, although in some circumstances we identify it with a suitable maximal subtorus of $\BB$. We let Weyl groups act on the left on tori, root systems 
etc., the action denoted either by an exponent on the left or 
as $w\cdot$; 
for example, the action of $W$=the Weyl group of $\GG$ on the character group of $\AA$ is given by: ${^w\chi}(a)=\chi(w^{-1}\cdot a)$.

As noted above, we feel free to identify the Langlands dual group $\check{G}$ of $\GG$ with its $\CC$-points, and therefore on the dual side we avoid the boldface notation. It comes equipped with a canonical maximal torus $A^*$.

To any such variety is attached the following set of data (cf.\ Section \ref{sec:review}): 
\begin{itemize}
\item[-] $\mathcal Z(\XX):=\Aut_\GG(\XX)^0$, the neutral component of the $\GG$-automorphism group of $\XX$;
\item[-] a parabolic subgroup $\PP(\XX)\supset \BB$, namely the stabilizer of the open Borel orbit;
its reductive quotient (and, sometimes, a suitable Levi subgroup) is denoted by $L(\XX)$ and the simple roots of $\AA$ in $\LL(\XX)$ by $\Delta_{L(X)}$;
\item[-] a torus $\AA_X$, which is a quotient of $\AA$;  it is the analog of the universal Cartan for the group.
 \item[-] a finite group $W_X$ of automorphisms of $\AA_X$, the ``little Weyl group;''
\item[-] the set $\Delta_X \subset \varchi(\XX) = \Hom(\AA_X,\GGm)$ of  \emph{normalized (simple) spherical roots} -- see \S \ref{invariants} for our normalization of the spherical roots; they are the simple roots of a based root system with Weyl group $W_X$.
\item[-] the {\em valuation cone} $\mathcal{V}$ inside $\mathfrak{a}_X := \Lambda_X\otimes \mathbb Q$, where $\Lambda_X=\varchi(\XX)^*=\Hom(\GGm,\AA_X)$.
It is the \emph{anti-dominant} Weyl chamber for the based root system defined by $\Delta_X$, and moreover contains
the image of the negative Weyl chamber of the pair $(\GG,\BB)$ under the projection $\mathfrak{a} \twoheadrightarrow \mathfrak{a}_X$ (where $\mathfrak a=\varchi(\BB)^*\otimes \QQ$). These anti-dominant chambers will also be denoted by $\mathfrak a^+$, $\mathfrak a_X^+$, and the intersection of $\Lambda_X$ with $\mathcal V$ is denoted by $\Lambda_X^+$. 
\item[-] a submonoid and a subsemigroup $A_X^+, \mathring A_X^+ \subset A_X$ defined as:
$$A_X^+ := \{a \in A_X: |\gamma(a)| \geq 1 \mbox{ for all } \gamma \in \Delta_X\}$$
and 
$$\mathring A_X^+ =  \{a \in A_X: |\gamma(a)| > 1 \mbox{ for all } \gamma \in \Delta_X\}.$$

 \item[-]
The subscript $\Theta$ denotes a subset of the set $\Delta_X$ of simple spherical roots associated to the spherical variety, and thus corresponds to the face of $\mathcal V$ to which it is orthogonal (the word ``face'' means the intersection of $\mathcal V$ with the kernel of a linear functional which is non-negative on $\mathcal V$ -- hence, it can refer to $\mathcal V$ itself). 
For each such $\Theta$ there is a distinguished subtorus $\AA_{X,\Theta}$ of $\AA_X$, with cocharacter group the orthogonal complement of $\Theta$ in $\Lambda_X$ (hence, $\AA_X=\AA_{X,\emptyset}$). We define $A_{X,\Theta}^+$, $\mathring A_{X,\Theta}^+$ in a similar way as for $A_X$ (see 
\S \ref{positiveelements}), using only the elements $\gamma\in \Delta_X\smallsetminus\Theta$. The group $\AA_{X,\Theta}$ is canonically isomorphic to the connected component of the $\GG$-automorphism group of the ``boundary degeneration'' $\XX_\Theta$ of $\XX$ (cf.\ \S \ref{normalbundles}), and therefore \emph{we are invariably using the notations $\AA_{X,\Theta}$ and $\mathcal Z(\XX_\Theta)$}.

\item[-]  \label{sufficientlydeep}By the phrase ``for $a$ sufficiently deep in $A_{X, \Theta}^+$''
we mean, informally, that $a$ is sufficiently far from the ``walls'' of $A_{X, \Theta}^+$;
formally, ``there exists $T>1$ so that, whenever $|\gamma(a)| \geq T$ for all $\gamma\in \Delta_X\smallsetminus\Theta$,  \dots.''

We use similar phrasing (e.g.\ ``sufficiently large'', ``sufficiently positive'') in many similar contexts. 

\end{itemize}

We follow standard notation for denoting duals of representations, e.g.\ $\tilde\pi$ is the smooth dual of a representation $\pi$, $\tilde M$ is the adjoint of a morphism $M$ between representations, etc. (In general, we feel free to move between the unitary and smooth categories of representations of $G$, since it is clear from the context -- or unimportant -- which of the two we are refering to.) 

Induction, the right adjoint functor to restriction, is denoted by $\Ind$, but we usually use the symbol $I$ to denote unitary induction, e.g.\ $I_P^G (\sigma)= \Ind_P^G(\sigma\delta_P^{\frac{1}{2}})$. Here $\delta_P$ is the modular character of the subgroup $P$, by which we mean the quotient of the right by a left Haar measure. (In the literature, this is sometimes $\delta_P^{-1}$; e.g.\ in Bourbaki.) The corresponding algebraic modular character: $\PP\to \GGm$, which is equal to the quotient of a right- by a left-invariant volume form, will be denoted by $\mathfrak d_P$.

Similarly, if $U$ is the unipotent radical of $P$, we use the notation $\pi_U$ to denote the normalized Jacquet module, that is:
as a vector space it consists of the $U$-coinvariants on $G$, and
we twist the action of $P$ by $\delta_{P}^{-1/2}$. In particular, there is always a canonical morphism: $(I_P^G(\sigma))_U\twoheadrightarrow \sigma$. 

We will in \S \ref{levivarieties} define certain parabolics $\PP_\Theta$, $\PP_\Theta^-$ (where $\Theta\subset\Delta_X$), and we will denote their unipotent radical by $\UU_\Theta$, $\UU_\Theta^-$ and their Levi quotient (or subgroup) by $\LL_\Theta$. In that case, the corresponding normalized Jacquet module $\pi_{U_\Theta}$, $\pi_{U_\Theta^-}$ of a smooth representation $\pi$ will also be denoted by $\pi_\Theta$, resp.\ $\pi_{\Theta^-}$. We caution the reader that this notation is only applied when working in the smooth category of representations; thus, the notation $L^2(X)_\Theta$ is reserved for a different space, a closed subspace of $L^2(X)$ defined in Corollary \ref{corollarysurjection}.
 
When $Y$ is an $H$-space (where $H$ is a group) endowed with a positive $H$-eigenmeasure, with eigencharacter $\eta$, we define the \emph{normalized} action of $H$ on functions on $Y$ by:
$$ (h\cdot f)(y) = \sqrt{\eta(h)} f(yh).$$
This makes $L^2(Y)$ (with respect to the given measure) a unitary representation, and identifies (in the setting of homogeneous spaces for $p$-adic groups) the space $C^\infty(Y)$ (\emph{uniformly} locally constant functions on $Y$, i.e. invariant by an open compact subgroup for $G$) with the smooth dual of $C_c^\infty(Y)$.

For a locally compact group $H$, 
we denote by $\widehat{H}$ the unitary dual of $H$, endowed with the Fell topology. Notice the notational distinction to $\check H$, which is used for the Langlands dual group of $H$. When $H$ is abelian, $\widehat{H}$ is the Pontryagin dual group of $H$, and we denote by 
 $\widehat{H}_{\C}$ the group of continuous homomorphisms $\Hom(H, \C^\times)$, so that $\widehat{H} \subset \widehat{H}_{\C}$. 

We slightly abuse the term ``wonderful'' to apply it to some embeddings of our spherical variety which are smooth but not necessarily wonderful; see \S \ref{compactifications} for details. 
 
   Given a function $Q$ on a subset of $(\mathbb{Z}_+)^r$ 
 we say that $Q$ is ``decaying'' if it is bounded by a negative exponential: $|Q(x_1, \dots, x_r)| \leq a^{\sum x_i}$
 where $a < 1$. We apply this term to functions on $A_X^+, \mathring{A}_X^+$ etc.
 by means of the natural valuation maps, i.e. sending an element in $a\in A_X^+$ to the valuations of all $\gamma(a)$.

We also use the notation $A \ll B$ to mean
that there exists a constant $c$ such that $A \leq c B$.

In the final part  we use using exponential notation for characters of tori, i.e.
for $\mathbf{T}$ a torus,  any character $\chi: \mathbf{T} \rightarrow \mathbb{G}_m$
will also be denoted by the symbol $e^{\chi}$ when it is more suggestive.

\subsection{Important assumptions:}\label{ssassumptions}
Later in the paper we introduce further assumptions on the spherical variety $\XX$ (and its boundary degenerations $\XX_\Theta$, introduced in \S \ref{normalbundles}), which are used in all 
theorems without explicit mention: 
\begin{itemize}
\item[-] We assume that the action of $\mathcal Z(\XX)=\Aut_{\GG}(\XX)^0$ is induced by the action of the center of $\GG$ (as we may in every case by replacing $\GG$ by $\mathcal Z(\XX)\times \GG$).

\item[-]  From Section \ref{localfieldgeom} onwards we suppose that $X$ carries a non-zero $G$-eigenmeasure (which we fix), and endow its ``boundary degenerations'' $X_\Theta$ with compatible eigenmeasures (cf.\ \S 4.1). As discussed in that section, this is not a significant restriction.

\item[-]
From Section \ref{sec:asymptotics} onward we assume that $\XX$ is {\em wavefront}, 
i.e., that it satisfies the condition enunciated in \S \ref{invariants}.
This is a genuine restriction, but applies in the vast majority of cases. {\bluetext We also assume from that point on that the connected central torus of $\GG$ surjects onto $\mathcal Z(\XX)$, which causes no harm in generality.}

\item[-]
 From Section \ref{sec:Bernstein} we suppose that the Discrete Series Conjecture \ref{dsconjecture} holds for  $\XX$;  this holds in many cases, e.g.\ $\XX$ is a symmetric variety 
 (see \S \ref{factorizable}) and can be checked in many others (possibly in every case) by the methods
 of \S \ref{ssunfolding}. 
 
{ \item[-] In Section \ref{sec:explicit} we assume, in addition, that $\XX$ satisfies the \emph{generic injectivity condition} of \S \ref{genericinjectivity} and is \emph{strongly factorizable}. ``Strongly factorizable'' is a strong condition that holds, for example, for symmetric spaces but not for most other spherical varieties. The generic injectivity condition holds more generally, again for symmetric varieties (Proposition \ref{delorme-recent-label}) and in many other cases (discussion
in \S \ref{genericinjectivity}). 

Notice that the generic injectivity condition is also explicitly assumed in the main scattering theorem \ref{advancedscattering}, but in its proof (Section \ref{sec:scattering}) it is only used at the very end, so Section \ref{sec:scattering} is not based on that assumption. 

\item[-] The main theorems of Section \ref{sec:explicit}, \ref{explicitBernstein}, \ref{explicitPlancherel}, are conditional on another combinatorial condition (which, again, is known to hold for symmetric varieties), injectivity of the ``small Mackey restriction'', but this condition is explicitly stated in the theorems.
}

 \end{itemize}

\subsection{Some open problems}  

Our work does not resolve the global (respectively local) conjectures about periods (respectively: the support of Plancherel measure for $L^2(X)$) including the need to refine these conjectures. In many cases (whenever the main Scattering Theorem \ref{advancedscattering} applies), the local conjecture is essentially reduced, through our work, to discrete spectra.\footnote{This is not completely the case when there are parameters into $\check G$ which admit many lifts to $\XX$-distinguished parameters; in that case, one needs to know that scattering maps respect whichever parametrization of discrete spectra by $\XX$-distinguished parameters one has, in order to reduce the conjecture through our work to discrete spectra.} Let us now discuss a few more open questions:

First, although we phrase our results throughout this paper in such a way that they should apply as stated
to the non-wavefront cases, most of our proofs break down. There are not many classes of non-wavefront varieties that we know of (the example $\GGL_n\backslash \SSO_{2n+1}$ is the typical one), but it seems that the non-wavefront case requires, and will lead to, a better understanding of harmonic analysis for $p$-adic groups and their homogeneous spaces. Notice that the example mentioned is also used by Knop  \cite[\S 10]{KnHC}, to illustrate that his powerful theory of invariant differential operators on spherical varieties provides a genuine extension of the action of the center of the universal enveloping algebra.

Secondly, since our Plancherel decomposition is based on the understanding of how discrete series vary with the central character (cf.\ the Discrete Series Conjecture \ref{dsconjecture}), it would be desirable to show in general (and not case-by-case, which can be done ``by hand'') that this conjecture holds, for instance that the unfolding process of \S \ref{ssunfolding} proves it. { This is a problem which pertains to the cases which are not what we call ``strongly factorizable''; the explicit Plancherel decomposition of section \ref{sec:explicit} via ``Eisenstein integrals'' is also open in those cases.}

Third, it would be desirable to settle the combinatorial assumptions of some of the theorems (notably, Theorem \ref{advancedscattering}) in some generality.

And, finally, a large class of problems has to do with developing Paley-Wiener theorems for spaces of Harish-Chandra--Schwartz or compactly supported functions.

\subsection{Acknowledgements}

We are grateful to Laurent Clozel, Patrick Delorme, Wee Teck Gan, Atsushi Ichino, Erez Lapid and Andre Reznikov
for interesting suggestions and discussions. We owe a special debt to Joseph Bernstein, who explained to us the core idea of \S \ref{sec:Bernstein}. We are also very grateful to Patrick Delorme for detailed comments and corrections on an earlier preprint.

We are also very grateful to the referees of the manuscript. Their careful reading
has greatly improved the paper. 

The first author (Y.S.) was supported by NSF grants
DMS-1101471 and DMS-1502270.
The second author (A.V.) was supported by a Packard Foundation fellowship and by the 
NSF grant
DMS-0903110.  We thank these agencies for their support.

This article has been written and re-written many times over the past {\bluetext nine} years, repeatedly proving our expectations for finishing it unrealistic. Indeed, it has been written so many times that the word ``nine'' has been repeatedly changed in the prior sentence.  Some of the results (those of section \ref{sec:asymptotics}) were announced as early as in the summer of 2007 at the Hausdorff Center in Bonn; others were based on the assumption of the spectral decomposition of sections \ref{sec:Bernstein} and \ref{sec:scattering}, which we only completed {\bluetext in its original form in 2013. The present form of the paper, which includes several corrections and improvements over the previous ones, many suggested by the referees, was completed in 2016. We offer our apologies to anyone who was promised a faster completion of the paper at any stage.}

\part{The dual group of a spherical variety}

	\section{Review of spherical varieties}  
	\label{sec:review}

The purpose of this section is to collect the necessary facts about spherical varieties. Most of the proofs will be given in the next section.

A \emph{spherical variety} for a (split, connected) reductive group $\GG$ over a field $k$ is a normal variety $\XX$ together with a $\GG$-action, such that the Borel subgroup of $\GG$ has a dense orbit. \emph{We will assume throughout that $\XX$ is homogeneous and quasi-affine.} The assumption of quasi-affineness is not a serious one, since every homogeneous $\GG$-variety which is not quasi-affine is the quotient of a quasi-affine $\GGm\times \GG$-variety by the action of $\GGm$ (and the cover is $\GGm\times \GG$ spherical if the original one was $\GG$-spherical). We fix a (complex) dual group $\check{G}$ to $\GG$; it comes equipped with a canonical maximal torus $A^*$, and the group $\check G$ is canonical up to conjugation by elements of $A^*$. One of our primary goals is to attach to the spherical variety $\XX$
a reductive group $\check{G}_X$ together with a morphism:
\begin{equation} \label{Desid-1} \check{G}_X \times \SL_2 \rightarrow \check{G}.\end{equation}
We shall see that this morphism determines a great deal about the spherical variety, both its geometry and its representation theory. 

In \S \ref{invariants} we introduce basic combinatorial invariants, including the root system associated to a spherical variety by F.\ Knop.  We also give the definition of a wavefront variety; 
we assume at most points in this text that the varieties under consideration are wavefront. 

In \S \ref{dualgroup} we modify this root system and discuss the associated reductive group $\check{G}_X$, which we term the dual group of the spherical variety. This is expected to be related to the group constructed by Gaitsgory and Nadler in \cite{GN}; in the next section (\S \ref{sec:dualgroupproofs}) we shall discuss the morphism \eqref{Desid-1}; in particular, we will prove the existence of this morphism if one makes certain natural assumptions regarding the Gaitsgory-Nadler group. 

In \S \ref{compactifications} we review the theory of toroidal compactifications, and in \S \ref{normalbundles} and \S \ref{ssdegen} we  study the normal bundles of $\GG$-orbits in those. This study will be of importance later: we will interpret the asymptotics of special functions on $\XX(k)$ using these normal bundles. 

In \S \ref{ssWhittakertype} we present the modifications needed in order to treat cases such as the Whittaker model.

In \S \ref{levivarieties} and \ref{sshorocycles} we introduce some other homogeneous varieties associated to $\XX$: Levi varieties and horospherical varieties. The former are closer to the traditional harmonic-analytic approach of studying Levi subgroups of a group, and will not be important for our formulations; they will be useful, however, for some proofs. The latter will be essential in explicating our harmonic-analytic constructions through the language of Eisenstein integrals in Section \ref{sec:explicit}.

Finally, in \S \ref{examples} we discuss the example of $\XX=\PPGL_n$.

		\subsection{Invariants} \label{invariants} 
		
We fix throughout a Borel subgroup $\BB$ and denote the open $\BB$-orbit on $\XX$ by $\mathring \XX$. (See, however, footnote \ref{footnoteBorel}.)

Let $\HH$ be an algebraic group acting on a variety $\YY$. The multiplicative group of non-zero rational $\HH$-eigenfunctions (semiinvariants) will be denoted by $k(\YY)^{(\HH)}$. We will denote the group of $\HH$-eigencharacters on $k(\YY)^{(\HH)}$ by $\varchi_H(\YY)$, and if $\HH$ is our fixed Borel subgroup $\BB$ then we will denote $\varchi_B(\YY)$ simply by $\varchi(\YY)$. If $\YY$ has a dense $\BB$-orbit, then we have a short exact sequence: $1\to k^\times \to k(\YY)^{(\BB)} \to \varchi(\YY) \to 1$. 

As an intrinsic way of understanding the geometry at $\infty$ of the variety $\XX$,
one studies valuations of its function field. A point of basic interest is how the geometry of $\XX$ at $\infty$ interacts
with the simple operation of acting on $\XX$ by a one-parameter subgroup $\GGm$.
The discussion that follows formalizes the study of such matters:

For a finitely generated $\mathbb Z$-module $M$ we denote by $M^*$ the dual module $\Hom_{\mathbb Z}(M,\mathbb Z)$. For our spherical variety $\XX$, we let $\Lambda_X=\varchi(\XX)^*$ and  $\mathfrak a_X = \Lambda_X\otimes_{\mathbb Z} \QQ$. A $\BB$-invariant, $\QQ$-valued valuation on $k(\XX)$ which is trivial on $k^\times$ (triviality on $k^\times$) will be implicitly assumed from now on) induces by restriction to $k(\XX)^{(\BB)}$ an element of $\Lambda_X$.
 We let $\mathcal V\subset \mathfrak a_X$ be the cone generated by the images of \emph{$\GG$-invariant valuations}. (By \cite[Corollary 1.8]{KnLV}, the map from $\GG$-invariant valuations to $\mathfrak a_X$ is injective.) We denote by $\Lambda_X^+$ the intersection $\Lambda_X\cap \mathcal V$. This is precisely the monoid of $\mathbb Z$-valued valuations. Notice that $\mathcal V$ contains the image of the \emph{negative} Weyl chamber under the natural map $\mathfrak a\twoheadrightarrow \mathfrak a_X$, \cite[Corollary 5.3]{KnLV}. 
 (To get a sense for some of the geometry here, and in particular why there are ``distinguished directions'' in $A_X$ at all, the reader may wish to glance at the example in \S \ref{axpluspgl}.)

The notation $\mathcal V$ is compatible with the literature on spherical varieties, but in this paper we also denote it, invariably, by $\mathfrak a_X^+$. We say that $\XX$ is a \emph{wavefront} spherical variety if $\mathcal V$ is precisely equal to the image of the negative Weyl chamber.  The terminology  is due to the validity of the Wavefront Lemma \ref{wavefront}; this class of varieties
was not previously singled out, to our knowledge.
Symmetric varieties, in particular, are wavefront \cite{KnLV}, but not, for instance, the variety $\UU\backslash \GG$, where $\UU$ is a maximal unipotent subgroup. (However, the latter becomes wavefront if we consider the additional action of a maximal torus ``on the left''; i.e., it is wavefront as a homogeneous variety for $\GG \times \AA$. Perhaps  the simplest non-wavefront spherical variety which cannot be treated in this way is the variety $\GGL_2\backslash \SSO_5$.)

The \emph{associated parabolic} to $\XX$ is the standard parabolic $\PP(\XX):= \{ g\in \GG | \mathring \XX \cdot g = \mathring \XX\}$. Let us choose a point $x_0\in \mathring \XX(k)$ and let $\HH$ denote its stabilizer; hence  $\XX=\HH\backslash \GG$, and $\HH\BB$ is open in $\GG$. There is the following ``good'' way of choosing a Levi subgroup $\LL(\XX)$ of $\PP(\XX)$, depending on the choice of $x_0$: Pick $f\in k[\XX]$, considered by restriction as an element of $k[\GG]^{\HH}$, such that the set-theoretic zero locus of $f$ is $\XX \smallsetminus \mathring \XX$.  Then $f$ is a $\PP(\XX)$-eigenfunction, but not an eigenfunction for a larger parabolic. Thus, its differential $df$ at $1\in \GG$ defines an element in the coadjoint representation of $\GG$, whose centralizer $\LL(\XX)$ is a Levi subgroup of $\PP(\XX)$.  The intersection $\LL(\XX)\cap \HH$ is known to contain the derived group of $\LL(X)$ \cite[Proposition 2.4]{KnAs}.

We fix throughout a maximal torus $\AA$ in $\BB \cap \LL(\XX)$. We define $\AA_X$ to be the torus: $\LL(\XX)/(\LL(\XX)\cap \HH) = \AA/ (\AA\cap \HH)$; its cocharacter group is the lattice $\Lambda_X$ defined above. We identify $\AA_X$ with a subvariety of $\mathring \XX$ through the orbit map: $a\mapsto x_0\cdot a$.

We can also think of $\AA_X$ as a canonical, abstract torus associated to $\XX$, in a similar way that the \emph{universal Cartan group} is associated to the group $\GG$. More precisely, if for every Borel subgroup $\BB$ we think of $\AA_X$ as the quotient by which $\BB$ acts on $\mathring \XX/\UU$ (here in the sense of geometric quotient), then for every two choices of Borel subgroups, any element conjugating one to the other induces the same isomorphism between the corresponding $\AA_X$'s. 

\begin{remark} A \emph{symmetric subgroup} $\HH$ is usually defined as the fixed point group of an involution $\theta$ on $\GG$, and a symmetric variety as the space $\XX=\HH\backslash\GG$; therefore it comes with a chosen point $x_0=\HH\cdot 1$. On the other hand, in the treatment of symmetric varieties one usually doesn't choose a Borel subgroup a priori. In that case, for a Borel $\BB$ such that $x_0\cdot \BB$ is open, the group $\PP(\XX)$ is what is called a \emph{minimal $\theta$-split parabolic}. Moreover, the Levi subgroup $\LL(\XX)$ constructed above is the unique $\theta$-stable Levi subgroup of $\PP(\XX)$, $\LL(\XX)=\PP(\XX) \cap {^\theta\PP(\XX)}$.
\end{remark}

The cone $\mathcal V=\mathfrak a_X^+$ is the fundamental domain for a finite reflection group $W_X\subset \End(\mathfrak a_X)$, called the \emph{little Weyl group} of $X$. If we denote by $W$, resp.\ $W_{L(X)}$, the Weyl groups of $\GG$ and $\LL(\XX)$ with respect to $\AA$ then there is a canonical way to identify $W_X$ with a subgroup 
of $W$, normalizing $W_{L(X)}$ (which it intersects trivially) \cite[\S 6.5]{KnHC}. The set of simple roots of $\GG$ corresponding to $\BB$ and the maximal torus $\AA\subset \BB$ will be denoted by $\Delta$. 

Consider the (strictly convex) cone negative-dual to $\mathcal V$: $\mathcal V^\perp=\{\chi\in\varchi(\XX)\otimes \RR | \left<\chi,v\right>\le 0 \text{ for every } v \in \mathcal V\}$. The generators of the intersections of its extremal rays with $\varchi(\XX)$ are called the \emph{spherical roots} of $\XX$. They are known to form the set of simple roots of a based root system with Weyl group $W_X$ \cite{KnAs}. This root system will be called the \emph{spherical root system} of $\XX$. We will denote this set of simple roots by $\Sigma_X$, following the notation of \cite{Lu}, in order to distinguish it from the set of ``(simple) normalized spherical roots'', which will be defined later and denoted by $\Delta_X$.

Finally, we discuss the group of $\GG$-automorphisms of $\XX$, based on \cite{KnAu}: Of course, for any homogeneous variety $\XX=\HH\backslash \GG$ we have $\Aut_\GG(\XX)=\mathcal N(\HH)/\HH$.\footnote{This also holds for any quasi-affine spherical variety, if $\HH$ denotes the stabilizer of a point on the open $\GG$-orbit -- cf.\ the remark after Lemma 6.6 in \emph{loc.cit.}.} As it is known \cite[5.2, Corollaire]{BP}, the quotient $\mathcal N(\HH)/\HH$ is diagonalizable. We will be denoting the torus $\Aut_\GG(\XX)^0$ by\footnote{\label{ftntbundle} In cases like the Whittaker model, where $\XX=\HH\backslash \GG$ and we have a morphism $\Lambda:\HH\to \GGa$ whose composition with a complex character of $k$ gives rise to the line bundle under consideration, the definition of $\mathcal Z(\XX)$ should be as the connected component of the subgroup of $\mathcal N(\HH)/\HH$ which stabilizes $\Lambda$, cf.\ \S \ref{ssWhittakertype}.} $\mathcal Z(\XX)$. {\bluetext There is no harm in assuming that the connected center $\mathcal Z(\GG)^0$ of $\GG$ surjects onto $\mathcal Z(\XX)$ under its natural action on $\XX$ (by replacing $\GG$ by $\mathcal Z(\XX)\times\GG$, if necessary). This will be a standing assumption from section \S \ref{sec:asymptotics} onward.}

Knop defines yet another root system\footnote{We caution the reader against confusing the notation of \cite{KnAu} with ours. We have reserved the letters $\Sigma$, $\Delta$ to denote sets of simple roots, as is customary in the literature, and more precisely the notation $\Delta_X$ for our normalized spherical roots. On the contrary, Knop uses the letters $\Delta$ and $\Delta_X$ to denote the full sets of roots of the root systems generated by what we denote here by $\Sigma_X$, resp.\ $\Sigma_X'$.} with a distinguished set $\Sigma_X'$ of simple roots which is proportional to $\Sigma_X$ (but does not coincide with our $\Delta_X$). Notice that every $\GG$-automorphism acts by a non-zero constant on each element of $k(\XX)^{(\BB)}$, and this defines a map: $\Aut_{\GG}(\XX)\to \Hom(\varchi(\XX),\GGm)=\AA_X$. Knop proves that the corresponding map of character groups belongs to a short exact sequence:
\begin{equation}
 0\to\left< \Sigma_X'\right>_\Z \to \varchi(\XX)=\varchi(\AA_X)\to \varchi(\Aut_\GG(\XX))\to 0.
\end{equation}
In particular, we have: 
\begin{equation}\label{autom}\varchi(\mathcal Z(\XX))=\varchi(\XX)/(\varchi(\XX)\cap \left<\Sigma_X\right>_\QQ).
\end{equation}

Finally, we include the following useful lemma of Knop:
\begin{lemma}[Non-degeneracy] \label{nondeg}
 For every quasi-affine spherical variety $\XX$, the set of coroots $\check\alpha$ of the group such that $\check\alpha$ is perpendicular to $\varchi(\XX)$ is precisely the set of coroots in the span of $\Delta_{L(X)}$.
\end{lemma}

\begin{proof}
 This is \cite[Lemma 3.1]{KnAs}. 
\end{proof}

\subsection{The dual group of a spherical variety}  \label{dualgroup}

  Here we formulate, without proofs, our results on the dual group of a spherical variety and the associated ``Arthur $\SL_2$''. All results will be proven in section \ref{sec:dualgroupproofs}.

Let $A^*$ denote the dual torus of $\AA$. By definition, it is the complex torus whose cocharacter group is $\varchi(\AA)$. Similarly, let $A_X^*$ be the dual torus of $\AA_X$. The map $\AA \twoheadrightarrow \AA_X$ dualizes to a map $A_X^* \to A^*$ (with finite kernel). Ideally, we would like this to extend to a map of connected reductive groups $\check G_X \to \check{G}$, where $\check G_X$ has Weyl group $W_X$. Unfortunately this is not possible in general. For example, in the case of the variety $\operatorname{\mathbf{PO}_2}\backslash \PPGL_2$ the kernel of $A_X^* \to A^*$ is of order $2$, but the group $\check G=\SL_2$ does not have any non-trivial connected covers.

Nonetheless, Gaitsgory and Nadler \cite{GN} have canonically associated to any spherical affine embedding of $\XX$ a reductive subgroup $\check G_{X,GN}\subset\check G$ whose maximal torus is the image of $A_X^*$ in $A^*$ and whose Weyl group is, conjecturally, equal to $W_X$. (In the case of $\operatorname{\mathbf{PO}_2}\backslash \PPGL_2$, this is just equal to $\SL_2$.) In our notation we suppress the dependence on the affine embedding: it is expected that $\check G_{X,GN}$ is independent of the choice of embedding, and in fact this follows from our arguments below, based on some natural assumptions (GN2)--(GN5) on the Gaitsgory-Nadler dual group. 

In this section, we shall construct -- under assumptions on $\XX$ that rule
out an example such as $\operatorname{\mathbf{PO}_2}\backslash \PPGL_2$-- a slightly different dual group $\check G_X$ whose
maximal torus is literally $A_X^*$.

Let $\gamma$ be a spherical root, i.e.\ $\gamma\in \Sigma_X$. It is known (either from the work of Akhiezer \cite{Ak} classifying spherical varieties of rank one, or from the work of Brion \cite{BrO} for a classification-free argument) that either:
\begin{enumerate}
\item $\gamma$ is proportional to a positive root $\alpha$ of $\GG$, or ;
\item $\gamma$ is proportional to the sum $\alpha+\beta$ of two positive roots which are orthogonal to each other and part of some system of simple roots (not necessarily the one corresponding to $\BB$).
 \end{enumerate}
\label{gammaprimedef}
Notice that in the second case $\gamma$ is not proportional to a root of $\GG$, and therefore the two cases are mutually exclusive.
In the first case we set $\gamma'=\alpha$, and in the second we set $\gamma'=\alpha+\beta$. Equivalently, $\gamma'$ is the primitive element in the intersection of the $\RR_+$-span of $\gamma$ with the root lattice of $\GG$. The motivation for these choices will be explained in the next section, where we will revisit the work of \cite{BrO}. Notice that in the second case the roots $\alpha,\beta$ are not unique; however, we will see that there is a canonical choice (whose elements will be called the roots \emph{associated to} $\gamma$). The set $\{\gamma'|\gamma\in\Sigma_X\}$ will be denoted by $\Delta_X$. Let $\Phi_X$ denote the set of $W_X$-translates of the elements in $\Delta_X$.  We will see in the next section:

\begin{proposition}
 The pair $(\Phi_X,W_X)$ defines a root system, and $\Delta_X$ constitutes a set of simple roots for it.
\end{proposition}

This root system will henceforth be called the \emph{normalized spherical root system} of $\XX$, and the elements of $\Delta_X$ will be the \emph{(simple) normalized spherical roots}.  When it doesn't matter if we are working with $\Sigma_X$ or $\Delta_X$ (for instance, choosing subsets of either of them), we may be abusing language and talking about ``the set of spherical roots'' $\Delta_X$. We will use the notation $\check \Phi_X$, $\check\Delta_X$ -- both subsets of $\varchi(\XX)^*$ -- for the dual root system and the corresponding set of simple coroots.

\begin{proposition} \label{dualgrouprootsystem} 
 Assume that $\Sigma_X$ does not contain any elements of the form $2\alpha$, where $\alpha$ is a root of $\GG$. Then the set $(\varchi(\XX)^*,\check\Phi_X,\varchi(\XX),\Phi_X)$ is a root datum. 
\end{proposition}

We refer to the corresponding complex reductive group $\check G_X$ with maximal torus $A_X^*$ as {\em the dual group of the spherical variety $\XX$.}  

The restrictions on $\Sigma_X$ guarantee, roughly speaking, that  the case of $\operatorname{\mathbf{PO}_2}\backslash \PPGL_2$ does not ``appear'' in the ``rank-one degenerations'' of $\XX$.

{ Our main result regarding the dual group of the spherical variety relies on the following statements. To formulate them, we must choose an affine embedding $\XX^a$ of $\XX$, in order to apply the work of Gaitsgory and Nadler and attach to it a group $\check G_{X^a,GN}$; our description (Theorem \ref{sl2}) of the Gaitsgory-Nadler group based on the following statements will eventually prove that $\check G_{X^a,GN}$ is independent of the choice of embedding.  The first statement is included in \cite{GN}, although not explicitly stated. We will treat the rest -- (GN2),  (GN3), (GN4) and (GN5) -- as axioms. It seems likely that these could be checked, and we discuss them briefly in \S \ref{GNaxioms}, but, as we are not specialists in the technical details, we prefer to phrase them as hypotheses.   

\begin{enumerate} 
\item[(GN1)] The image of $\check G_{X^a,GN}$ commutes with ${2\rho_{L(X)}}(\CC^\times)\subset A^*$, where $2\cdot\rho_{L(X)}$ denotes the sum of positive roots of $\LL(\XX)$, considered as a morphism: $\Gm\to A^*$.
\item[(GN2)] The Weyl group of $\check G_{X^a,GN}$ equals $W_X$.
\item[(GN3)] For any $\Theta\subset \Sigma_X$ the dual group of $\XX^a_{\Theta}$ is canonically a subgroup of $\check G_{X^a,GN}$. Here, $\XX^a_{\Theta}$ is a certain ``affine degeneration'' of $\XX^a$, to be introduced in \S \ref{ssdegen}.
\item[(GN4)] If the open $\GG$-orbit $\XX\subset \XX^a$ is parabolically induced, $\XX = \XX_L \times^{\PP^-} \GG$, where $\XX_L$ is spherical for the reductive quotient $\LL$ of $\PP^-$, then the dual group $\check{G}_{X^a,GN}$
belongs to the standard Levi subgroup $\check L$ of $\check G$ corresponding to the class of parabolics\footnote{We clarify the correspondence between $\PP^-$ and $\check L$. Recall that $A^*\subset \check G$ is the dual torus of the ``universal Cartan'' of $\GG$, i.e.\ the reductive quotient of any Borel subgroup of $\GG$. One chooses a Borel $\BB$ opposite to $\PP^-$, and the positive coroots of the Levi $\check L$ are the roots in the Lie algebra of $\BB\cap \PP^-$.} opposite to $\PP^-$.  Moreover, if a connected normal subgroup $\LL_1$ of $\LL$ acts trivially on $\XX_L$, then $\check{G}_{X^a,GN}$ belongs to the dual group of $\LL/\LL_1$ (which is canonically a subgroup of $\check L$).
\item[(GN5)] If $\XX_1^+$ is a spherical homogeneous $\GG$-variety, $\TT$ a torus of $\GG$-automorphisms and $\XX_2^+=\XX_1^+/\TT$, and if $\XX_1, \XX_2$ are affine embeddings of $\XX_1^+,\XX_2^+$ with $\XX_2=\spec k[\XX_1]^\TT$,  then there is a canonical inclusion $\check{G}_{X_2, GN} \hookrightarrow \check{G}_{X_1, GN}$ which restricts to the natural inclusion of maximal tori: $A_{X_2, GN}^* \hookrightarrow A_{X_1, GN}^*$ (arising from $\varchi(\XX_2)\hookrightarrow\varchi(\XX_1)$).
\end{enumerate}
}

The formulation of the result is based on the notion of a ``distinguished morphism'' from the group $\check G_X$ to $\check G$. By a \emph{distinguished morphism} we mean a morphism which restricts to the canonical map of maximal tori: $A_X^*\to A^*$, and moreover satisfies a condition on the image of simple root spaces, which will be formulated in \S \ref{distinguishedmorphisms}. A \emph{distinguished morphism} from $\check G_X\times \SL_2$ to $\check G$ is one which restricts to a distinguished morphism on $\check G_X$ and, moreover, under the ``standard'' diagonal embedding of $\Gm$ to $\SL_2$ it restricts to the map ${2\rho_{L(X)}}: \Gm\to A^*$.

\begin{theorem} \label{sl2}

 Assuming that $\Sigma_X$ does not contain any elements of the form $2\alpha$ (where $\alpha$ is a root of $\GG$):
 
 \begin{enumerate}
 \item\label{unique} There is at most one $A^*$-conjugacy class of distinguished morphisms $$\check G_X \times \SL_2 \rightarrow \check G.$$
 
 \item Assume Axioms (GN2), (GN3), (GN4), (GN5). Then distinguished morphisms exist. 
 
 \item Assume (GN2), (GN3), (GN4), (GN5). Then the root system of the Gaitsgory-Nadler dual group is given by $(\check \Phi_X, W_X)$; in particular, it is independent of the choice of affine embedding of $\XX$.
 \end{enumerate}
 \end{theorem}

\begin{remark}\label{sl2after}
 We stated this theorem for the case that $\Sigma_X$ does not contain elements of the form $2\alpha$. However, in the general case the same statements are true if one replaces $\check G_X$ with the abstract group $\check G_X'$ defined by the same root system and with maximal torus equal to the image of $A_X^*\to A^*$. In particular, $\check G_{X,GN}\simeq \check G_X'$. However, this does not seem to be the correct dual group for all purposes of representation theory.
\end{remark}
  We give a rather clumsy proof of this theorem, which boils down to case-by-case considerations, in the next section. 
\begin{example}
 For the spherical variety $\XX=\TT\backslash\SSL_2$ (where $\TT$ is a non-trivial torus) the dual group is $\check G_X=\SL_2$ with its natural isogeny: $\SL_2\to \check G=\PGL_2$. Notice that here the action of $\SSL_2$ on $\XX$ factors through the quotient $\SSL_2\to\PPGL_2$. This ``explains'' why the representation theory of $\XX$ should be described in terms of $\SL_2$, which is the dual group of $\PPGL_2$. However,  this argument may not always work when $\check G_X\ne\check G_{X,GN}$.
\end{example}

		\subsection{Toroidal compactifications} \label{compactifications} 

We will freely use the word ``compactification'' of $\XX$ for what should more correctly be termed ``spherical embedding'' or ``partial spherical completion'', namely a normal $\GG$-variety $\bar\XX$ containing $\XX$ as a dense $\GG$-subvariety. A compactification $\bar\XX$ is said to be \emph{simple} if it contains a unique closed $\GG$-orbit, and \emph{toroidal} if no $\BB$-stable divisor in $\XX$ contains a $\GG$-orbit of $\bar\XX$ in its closure. To every simple toroidal embedding $\bar\XX$ we associate the cone $\mathcal C(\bar\XX)\subset\mathcal V$ spanned by the valuations defined by all $\GG$-stable divisors in $\bar\XX$. The main theorem in the classification of such embeddings states: 

\begin{theorem}[Luna and Vust, cf.\ {\cite[Theorem 3.3]{KnLV}}]\label{thmlunavust} The map $\bar\XX\mapsto \mathcal C(\bar\XX)$ induces a bijection between isomorphism classes of simple toroidal compactifications of $\XX$ and strictly convex, finitely generated cones in $\mathcal V$.
\end{theorem}

\begin{remark}
 Although this theorem, and the other theorems of Brion, Luna, Pauer and Vust which we are going to recall, have been stated for an algebraically closed field in characteristic zero, their proofs carry through over an arbitrary field in characteristic zero, as long as the group $\GG$ is split.
\end{remark}

Notice that for every $\GG$-orbit $\ZZ$ in a simple toroidal embedding $\bar\XX$, the union of all $\GG$-orbits whose closure contains $\ZZ$ is also a simple toroidal embedding. This way, we get a bijection between faces\footnote{The word ``face'' is used for the intersection of $\mathcal C(\bar\XX)$ with the kernel of a linear functional which is non-negative on $\mathcal C(\bar\XX)$; hence $\mathcal C(\bar\XX)$ is the face corresponding to the closed orbit, and $\{0\}$ is the face corresponding to $\XX$.} of $\mathcal C(\bar\XX)$ and orbits of $\GG$ on $\XX$. 

Now observe that when $\mathcal V$ itself is strictly convex (equivalently: $\Aut_\GG(\XX)$ is finite), this implies the existence of a canonical compactification with $\mathcal C(\bar\XX)=\mathcal V$. Its existence can also be characterized by the following conditions:

\begin{theorem}[{\cite[5.3, Corollaire]{BP}}] \label{BPtheorem} The following are equivalent:
\begin{enumerate}
 \item  There exists a simple complete toroidal compactification of $\XX$.
\item The $\GG$-automorphism group $\Aut_{\GG}(\XX)=\mathcal N(\HH)/\HH$ is finite.
\item The cone $\mathcal V$ is strictly convex.
\end{enumerate}
\end{theorem}

The corresponding complete variety $\bar\XX$ is sometimes (e.g.\ in \cite{KnAu}) called the \emph{wonderful compactification} of $\XX$, though the term ``wonderful'' (``magnifique'') is more often (e.g.\ in \cite{Lu}) reserved for the case when $\bar\XX$ is smooth. To understand the difference between the two, let us recall the Local Structure Theorem of Brion, Luna and Vust: Let $\bar\XX$ be a simple toroidal embedding of $\XX$. The complement of all $\BB$-stable divisors of $\bar\XX$ which are not $\GG$-stable is denoted by $\bar\XX_B$, and it is a $\BB$-stable open affine subvariety. Let $\YY$ be the closure of $\AA_X$ in $\bar\XX_B$; it is the toric compactification of $\AA_X$ characterized by the property that for $\lambda\in \Lambda_X=\Hom(\GGm,\AA_X)$ we have $\lim_{t\to 0} \lambda(t) \in \YY$ if and only if $\lambda\in \Lambda_X^+$. 

\begin{theorem}[{\cite[Th\'eor\`eme 3.5]{BLV}}] \label{localstructure} The action map $$\YY \times \UU_{P(X)}\to \bar\XX$$ is an open embedding and its image $\bar\XX_B$ meets each $\GG$-orbit in $\bar\XX$ along its open $\BB$-orbit.
\end{theorem}

Hence, the variety $\bar\XX$ will be smooth if and only if the toric variety $\YY$ is. Since $\mathcal V$ is a simplicial cone, it follows from this theorem that 
the  compactification  of Theorem \ref{BPtheorem} has, if any, only finite quotient singularities, and that it is smooth if and only if the monoid $\Lambda_X^+$ is generated by its intersections with the extremal rays of $\mathcal V$; equivalently: if and only if $\Sigma_X$ generates $\varchi(\XX)$. In that case, all $\GG$-orbits are smooth and the complement of the open $\GG$-orbit is a union of $\GG$-stable divisors intersecting transversely. 

\label{sctnotwavefront}
{ In this paper we will adopt the following convention: We will, by abuse of language, \textbf{refer to any smooth, complete, toroidal compactification of $\XX$ as the ``wonderful compactification''}; and this term will also be extended to certain non-complete embeddings when we consider Whittaker-type induction in \S \ref{ssWhittakertype}. We will attach certain ``boundary degenerations'' to $\XX$ indexed by subsets of the set of spherical roots, which although by construction seem to depend on the choice of such a compactification, an important result -- Proposition \ref{thetaindependent} -- states that they are actually completely canonical. Our representation-theoretic results, then, will be formulated in terms of these degenerations. Whenever proofs are identical for the  wonderful compactification of Theorem \ref{BPtheorem} (if it exists and is smooth) and a general smooth, complete, toroidal compactification, for simplicity and clarity we only formulate them for the former; whenever the general case needs extra arguments, 
we provide them.}
 
Let us now discuss the general case, to which one will necessarily resort when $\Aut_\GG(\XX)$ is not finite, or it is finite but the canonical embedding of Theorem \ref{BPtheorem} is not smooth.
 
\subsubsection{The non-wonderful case} \label{sssnonwonderful}  
 
A toroidal embedding which is not necessarily simple is described by a strictly convex fan, instead of a strictly convex cone, that is: a collection $\mathfrak F$ of (distinct) strictly convex subcones of $\mathcal V$ as above, such that each face of a cone of $\mathfrak F$ is also contained in $\mathfrak F$, and each point of $\mathcal V$ belongs to the relative interior of at most one cone in $\mathfrak F$. We point the reader to \cite[Theorem 3.3 and discussion after Corollary 5.3]{KnLV} for details. 

To get a \emph{complete, smooth} toroidal embedding we need to subdivide $\mathcal V$ into a finite union of (strictly convex) simplicial cones $\mathcal C_i$, such that $\Lambda_X\cap \mathcal C_i$ is a free monoid for every $i$; then these cones and their faces form the fan of the embedding. Each such embedding is the union of simple (non-complete) smooth toroidal embeddings; hence, the complement of the open $\GG$-orbit is a divisor with normal crossings. 
 
The Local Structure Theorem \ref{localstructure}, applies equally well to this case with the understanding that the toric variety $\YY$
is that associated to the fan defining the toroidal embedding. 

The set of $\GG$-orbits in such an embedding is in bijection with the set of cones in the fan. The relative interior of the cone corresponding to an orbit $\ZZ$ consists of those valuations $v$ whose \emph{center} is the closure of $\ZZ$, that is: $\mathfrak o_v \supset \mathfrak o_Z$ and $\mathfrak m_v \supset \mathfrak m_Z$, where by $\mathfrak o$ and $\mathfrak m$ we denote the subring of $k(\XX)$ and its ideal, respectively, defined by the valuation $v$ or by the closure of $\ZZ$. This face, in turn, corresponds (non-injectively, in general) to the subset $\Theta\subset \Sigma_X$ of spherical roots to which it is orthogonal. (Since the elements of $\Sigma_X$ are proportional to the normalized spherical roots, i.e.\ the elements of $\Delta_X$, we will interchangeably be identifying $\Theta$ with a subset of either of the two.) We say that ``$\ZZ$ corresponds to $\Theta$''.

\subsubsection{The notion of $\Theta$-infinity}

\label{Thetainftydisc}
For each smooth toroidal embedding $\bar \XX$ and every $\Theta\subset\Delta_X$ (or $\Theta\subset\Sigma_X$) we let $\infty_\Theta$, the ``\emph{$\Theta$-infinity}'', denote the closure of the union of all $\GG$-orbits which correspond to $\Theta$, i.e.\ whose corresponding face is orthogonal to the set $\Theta$ of spherical roots. We have $\infty_{\Delta_X} = \bar\XX$.   As we will see in a moment, all of those orbits correspond to isomorphic ``boundary degenerations''. 

When working at the level of $k$-points,  we use the notion of ``neighborhood of $\Theta$-infinity'', which we explicate for clarity: A ``neighborhood of $\Theta$-infinity'' in $X$ is, by definition, the intersection of $X$
with a  neighborhood of $\infty_\Theta$ in $\overline{X}$ (for the Hausdorff topology induced by that of $k$). 
Note that the fact that $\infty_\Theta$ is the {\em closure} of the union of all orbits corresponding to $\Theta$ affects the meaning of this notion { When the embedding $\overline{\XX}$ is not explicated, we will mean a ``wonderful'' embedding (in the above sense). It is clear that this abstract notion of $\Theta$-infinity does not depend on the choice of a wonderful (i.e.\ complete, toroidal) embedding.}

  For a given $\GG$-orbit on $\overline{\XX}$, we will say that it ``belongs to $\Theta$-infinity'' if it is contained in $\infty_\Theta$, but not in $\infty_\Omega$, for any $\Omega\subsetneq\Theta$.

\subsubsection{Spherical system of a $\GG$-orbit}

For the following proposition we remind that a ``color'' of a spherical variety is a prime $\BB$-stable divisor which is not $\GG$-stable (hence colors are in bijection with $\BB$-orbits of codimension one in the open $\GG$-orbit, and by abuse of language we will be calling these $\BB$-orbits ``colors'' as well), and that each color $\DD$ defines a ``valuation'' $\check v_D\in \varchi(\XX)^*$, defined exactly as in the case of $\GG$-stable divisors that we discussed above (i.e.\ by restriction to $k(\XX)^{(\BB)}$). 

\begin{proposition}\label{sphericalsystem-orbit}
 Let $\ZZ$ be a $\GG$-orbit  in a toroidal compactification, and let $\Theta$ be the set of (unnormalized) spherical roots to which the corresponding face $\mathcal F$ is orthogonal. Then:
\begin{enumerate}
 \item $\varchi(\ZZ) = \mathcal F^\perp\subset \varchi(\XX)$;
 \item $\PP(\ZZ)=\PP(\XX)$;
 \item $\Sigma_\ZZ = \Theta$;
 \item for each simple root $\alpha$ of $\GG$, which belongs to $\Theta$, there are precisely two colors $\DD_1',\DD_2'$ in $\mathring \ZZ\PP_\alpha$, obtained as the (multiplicity-free) intersection with $\ZZ$ of the closures of the two colors $\DD_1,\DD_2$ in $\mathring \XX\PP_\alpha$; the valuations $\check v_{D_i'}$ are the images of the valuations $\check v_{D_i}$ under the restriction map: $\varchi(\XX)^*\to \varchi(\ZZ)^*$. 
\end{enumerate}
\end{proposition}

\begin{proof}
 
The first and second statement follow from the Local Structure Theorem \ref{localstructure}. 

The third statement can be proven by induction on the codimension of $\ZZ$; hence, we may assume that $\mathcal F$ is a half-line. If this half-line belongs to $\mathcal V\cap (-\mathcal V)$ then the Local Structure Theorem \ref{localstructure} implies that $\ZZ$ is isomorphic to the quotient of $\XX$ by the subtorus of $\mathcal Z(\XX)$ corresponding to this line. Otherwise, the theory of toroidal embeddings implies that there is a wonderful compactification of $\XX/\Aut_\GG(\XX)$ which contains a quotient of $\ZZ$ by its automorphism group, and the statement follows by a characterization of simple spherical roots in terms of the $\GG$-action on the normal bundle to the closed orbit, cf.\ \cite[\S 1.3]{Lu}. 

If $\alpha$ is a simple root of $\GG$ which belongs to $\Theta=\Sigma_\ZZ$ then there are two colors in $\mathring \ZZ\PP_\alpha$, and they are precisely those colors on which the Borel eigenfunction $f_\alpha$ of weight $\alpha$ vanishes -- the valuation of $f_\alpha$ on each of them is one \cite[\S 1.4]{Lu}, \cite[Proposition 3.4]{LuGro}. Every nonzero Borel eigenfunction on $\ZZ$ extends (uniquely) to a Borel eigenfunction on $\XX$ \cite[Theorem 1.3]{KnLV}, and similarly the extension of $f_\alpha$ has simple zeroes on the two colors of $\mathring\XX\PP_\alpha$. This shows that their closures intersect $\ZZ$ along the colors of $\mathring \ZZ\PP_\alpha$ without multiplicity; the fact that all eigenfunctions extend means that the valuations induced by the latter are equal to the valuations induced by the former, restricted to the character group of $\ZZ$.
\end{proof}

\subsection{Normal bundles and boundary degenerations}	  \label{normalbundles}
		 
We are interested in understanding normal bundles of $\GG$-orbits in wonderful embeddings; as we will see in Section \ref{sec:asymptotics}, the asymptotics of eigenfunctions on $X$ can be understood in terms of eigenfunctions on the normal bundles. 

Let $\ZZ$ be a $\GG$-orbit in a smooth toroidal compactification $\bar\XX$ of $\XX$, corresponding to a subset $\Theta$ of the set of spherical roots (\S \ref{Thetainftydisc}; in the propositions that follow, $\Theta$ will be considered as a subset of either the unnormalized spherical roots $\Sigma_X$ or the normalized spherical roots $\Delta_X$, according to what is appropriate in each case). The normal bundle $N_\ZZ{\bar\XX}$ of $\ZZ$ in $\bar\XX$ carries an action of $\GG$, under which it is spherical; this follows immediately from the Local Structure Theorem \ref{localstructure}. The open $\GG$-orbit on $N_\ZZ{\bar \XX}$ will be denoted by $\XX_\Theta$ and called a \emph{boundary degeneration} of $\XX$.   

 We remark  that if $\XX$ is wavefront, then so is $\XX_\Theta$ {\bluetext under the action of $\mathcal Z(\XX_\Theta)\times \GG$}, for every $\Theta\subset\Delta_X$. That will follow from Proposition \ref{wavefrontlevi}. 

\begin{remark} 
We notice that in the case of general toroidal compactifications there may be many $\GG$-orbits corresponding to the same $\Theta$; as we will see, all varieties $\XX_\Theta$ will be canonically isomorphic to each other, which will allow us to use this notation indistinguishably.
\end{remark}

Let $\ZZ$ be a $\GG$-orbit in $\bar\XX$ corresponding to $\Theta\subset\Sigma_X$. 
The fact that $\XX_\Theta$ belongs to a normal bundle gives rise to a torus of $\GG$-equivariant automorphisms of it:

 \begin{lemma} \label{torusacting} 
Let $\Gamma$ denote the set of $\GG$-stable divisors in $\bar\XX$ which contain $\ZZ$, then there is a canonical action of the torus $\GGm^{\Gamma}$ on $\XX_\Theta$ by $\GG$-automorphisms, such that the quotient is isomorphic to $\ZZ$.
\end{lemma}

\begin{proof}
{ Since these $\GG$-stable divisors intersect transversely, the normal bundle $N_\ZZ \bar\XX$ splits canonically into a direct
sum of line bundles $L_\gamma$ indexed by the spherical roots $\gamma\in \Gamma$. Moreover, $\XX_\Theta$ is
isomorphic to the total space of the direct sum $\bigoplus_{\gamma\in \Gamma} L_\gamma$, minus the union
of the smaller direct sums; this follows, for example, from the Local Structure Theorem \ref{localstructure}. This is a principal bundle over $\ZZ$ with structure group the torus $\GG_m^\Gamma$.}
\end{proof}

{ Proposition \ref{thetaindependent} in the next subsection will say that for different orbits $\ZZ$ in a smooth toroidal compactification, corresponding to the same $\Theta$, the resulting boundary degenerations $\XX_\Theta$ are \emph{canonically} isomorphic. Removing the word ``canonically'',} this could also be inferred from the uniqueness theorem of Losev \cite[Theorem 1]{Lo} and the following proposition:

\begin{proposition}\label{sphericalsystem} \label{Cartanident}
 Let $\ZZ$ be a $\GG$-orbit in a toroidal compactification, and let $\Theta\subset \Sigma_X$ be the set of (unnormalized) spherical roots to which the corresponding face $\mathcal F$ is orthogonal, and let $\XX_\Theta$ be defined as above. Then:
\begin{enumerate}
 \item \label{samechar} $\varchi(\XX_\Theta) = \varchi(\XX)$; more precisely,  there is a canonical isomorphism of ``universal tori'': 
\begin{equation}\label{Cartanidenteq}
 \AA_X\simeq \AA_{X_\Theta}.
\end{equation}
 \item \label{sameparab} $\PP(\XX_\Theta)=\PP(\XX)$;
 \item \label{sameroots} $\Sigma_{\XX_\Theta} = \Theta$;
 \item \label{samecolors} for each simple root $\alpha$ of $\GG$, which belongs to $\Theta$, there are precisely two colors $\DD_1',\DD_2'$ in $\mathring \XX_\Theta\PP_\alpha$, and they induce the same valuations $\check v_{D_i}$ as the two colors in $\mathring\XX\PP_\alpha$.
\end{enumerate}
\end{proposition}

\begin{proof}
 Fix a Borel subgroup $\BB$, and consider the local $\PP(\XX)$-equivariant isomorphism of the Local Structure Theorem \ref{localstructure}. Notice that this isomorphism is not canonical, but its quotient by the action of $\UU_{P(X)}$ is. For the smooth toric variety $\YY$, if $\YY'\subset\YY$ denotes the closure of the orbit corresponding to $\ZZ$, there is a unique isomorphism: $\phi:N_{\YY'}\YY\xrightarrow{\sim}\YY$ with the following properties:
\begin{enumerate}
 \item $\phi$ is $\LL(\XX)$-equivariant;
 \item its ``partial'' differential induces the identity on $N_{\YY'}\YY$.
\end{enumerate}

Restricted to the open $\LL(\XX)$-orbits, this map gives a canonical isomorphism: $\AA_{X_\Theta}\xrightarrow{\sim}\AA_X$. The second statement also follows from Theorem \ref{localstructure}. The third follows from Proposition \ref{sphericalsystem-orbit} and Lemma \ref{torusacting}.  

The last assertion is proven as in Proposition \ref{sphericalsystem-orbit}, except that now one has to consider an affine   \emph{degeneration} of $\XX$ to $\XX_\Theta$, which will be discussed in \S \ref{ssdegen}: it is a $\GG$-variety $\mathscr X^a$ over a base $\mathscr B$, whose generic fiber is isomorphic to an affine completion of $\XX$ and its special fiber is isomorphic to an affine completion of $\XX_\Theta$. The inclusion of the open Borel orbit on each fiber is a direct product: $\mathring{\mathscr X}^a\simeq \mathscr B\times \mathring\XX$, and therefore $\BB$-eigenfunctions extend non-trivially to the special fiber. We omit the details, since this result will not be used in the sequel. 
\end{proof}

In terms of dual groups, the following is implied by Proposition \ref{Cartanident} under the assumptions of Theorem \ref{sl2}. (When $\Sigma_X$ contains elements of the form $2\alpha$, the analogous statement applies to $\check G_{X,GN}$, cf.\ Remark \ref{sl2after}.)

\begin{corollary}
 The dual group of $\XX_\Theta$ is the Levi of $\check G_X$ with simple roots $\check\Theta$ (where $\check\Theta$ denotes the set of coroots of the elements of $\Theta$ -- now considered as a subset of $\Delta_X$).
\end{corollary}

\subsubsection{Identification of Borel orbits}

Notice that the map (\ref{Cartanidenteq}), together with the Local Structure Theorem \ref{localstructure}, gives rise to a $\BB$-equivariant isomorphism:
\begin{equation}\label{Borbitident}
\mathring\XX_\Theta \xrightarrow{\sim} \mathring\XX
\end{equation}
inducing the identity on normal bundles. Such isomorphism is not completely canonical, as it depends on the choice involved in the Local Structure Theorem. However, it is canonical up to $\BB$-automorphisms of $\mathring\XX$ which induce the identity on $\mathring\XX/\UU$, which means up to a morphism of the form: 
$$(a,u)\mapsto (a, a^{-1}u_1a u)$$ 
(in the setting of Theorem \ref{localstructure}), where $u_1\in U_{P(X)}$ is fixed by the kernel of the map: $L(X)\to A_X$.

\subsubsection{Automorphisms}
It follows now from (\ref{autom}) that $\mathcal Z(\XX_\Theta)=\Aut_\GG(\XX_\Theta)^0$ can be canonically identified with the maximal subtorus\footnote{We will be using the notation $\mathcal Z(\XX_\Theta)$ and $\AA_{X,\Theta}$ interchangeably. On the other hand, $\AA_{X,\Theta}$ is not to be confused with $\AA_{X_\Theta}$; the latter, as we saw in Lemma \ref{Cartanident}, is isomorphic to $\AA_X$, which allows us never to use the notation $\AA_{X_\Theta}$ again.} $\AA_{X,\Theta}$ of $\AA_X$ whose cocharacter group is perpendicular to $\Theta$. We wish to explicate the embedding $\GGm^{\Delta_X\smallsetminus\Theta}\to \AA_{X,\Theta}$ obtained from Lemma \ref{torusacting} (we restrict ourselves to the wonderful case, where the $\Gamma$ of Lemma \ref{torusacting} is equal to $\Delta_X\smallsetminus\Theta$); in fact, it is enough to do so when $\Delta_X\smallsetminus\Theta$ only has one element, since these embeddings are obviously compatible with each other. 
Choose a Borel subgroup $\BB$ and a $\gamma\in \Delta_X$. Let $\gamma'$ be the corresponding unnormalized spherical root (i.e.\ $\gamma'\in \Sigma_X$). It follows from the Local Structure Theorem \ref{localstructure} that the valuation induced by the orbit corresponding to $\Delta_X\smallsetminus\{\gamma\}$ is equal to $-{\gamma'}^*$, where $-{\gamma'}^*$ is the element of $\mathcal V$ with $\left<\gamma',-{\gamma'}^*\right>=-1$ and $\left<\delta,-{\gamma'}^*\right>=0$ for all $\delta\in \Sigma_X\smallsetminus \{\gamma'\}$. (If $\mathcal Z(\XX)\ne 1$ then we take $-{\gamma'}^*\in \mathcal V'$.) Notice that under our assumption that $\bar\XX$ is smooth, the elements $-{\gamma'}^*, \gamma'\in \Sigma_X$ form a basis for the monoid $\Lambda_\XX^+$. Hence, the action of $m\in \GGm^{\{\gamma\}}$ multiplies a $\BB$-eigenfunction with eigencharacter $\chi$ by $\left<\chi,-{\gamma'}^*\right>$, and hence we obtain:

\begin{lemma}\label{torilemma}
 The composition of (\ref{autom}) with the identification of $\mathcal Z(\XX_\Theta)$ with a subtorus of $\AA_X$ is given by the maps:
\begin{equation}
 -{\gamma'}^*: \GGm^{\{\gamma\}} \to \AA_X.\end{equation}
\end{lemma}

\subsubsection{Positive elements} \label{positiveelements}

For every $\Theta\subset\Delta_X$ we denote: \begin{equation}
A_{X,\Theta}^+ := \{a \in A_{X,\Theta}: |\gamma(a)| \geq 1 \mbox{ for all } \gamma \in \Delta_X\smallsetminus\Theta\},
\end{equation}
and:
\begin{equation}
\mathring A_{X,\Theta}^+ =  \{a \in A_{X,\Theta}: |\gamma(a)| > 1 \mbox{ for all } \gamma \in \Delta_X\smallsetminus\Theta\}. 
\end{equation}

Let $\Lambda_{X,\Theta}^+$ (resp.\ $\mathring \Lambda_{X,\Theta}^+$) denote the intersection of $\Lambda_X (=\varchi(\XX)^*)$ with the face (resp.\ the relative interior of the face) of the cone $\mathcal V$ which is orthogonal to $\Theta$.
We can alternatively describe $A_{X,\Theta}^+$ (resp. $\mathring A_{X,\Theta}^+$) as the submonoid generated
by all  $m(x)$,  where $m\in \Lambda_{X,\Theta}^+$ is a cocharacter in $\mathcal V$ and $x\in k^\times\cap \mathfrak o$ 
 (resp.\ $m\in \mathring\Lambda_{X,\Theta}^+$ and $x \in k^{\times}  \cap \mathfrak p$).
 Hence, $A_{X,\Theta}^+$ (resp.\ $\mathring A_{X,\Theta}^+$ is the preimage of $\Lambda_{X,\Theta}^+$ (resp.\ of $\mathring \Lambda_{X,\Theta}^+$) under the ``valuation'' maps: 
$$A_{X,\Theta}\to A_{X,\Theta}/\AA_{X,\Theta}(\mathfrak o)\simeq \Lambda_{X,\Theta},$$
normalized so that for $\check\lambda\in \Lambda_{X,\Theta}$ the valuation of $\check\lambda(\varpi)$ is $\check\lambda$.

 When $\Theta=\emptyset$, $\AA_{X,\emptyset}=\AA_X$ and the notation is compatible with the notation $\mathfrak a_X^+$ that we have been invariably using for the cone $\mathcal V$ of invariant valuations. In this case we will denote $A_{X,\Theta}^+,\mathring A_{X,\Theta}^+$ by $A_X^+,\mathring A_X^+$ (since $A_{X,\emptyset}=A_X$). We remark that under the map: $A\to A_X$, the anti-dominant elements of $A$ are contained in $A_X^+$; indeed, we have already seen that $\mathcal V$ contains the image of the negative Weyl chamber of the group.

The following is an easy corollary of the definitions and of Lemma \ref{torilemma}:

\begin{lemma} \label{limit-lemma-goo}
Under the map $\GGm^{\{\gamma\}} \to \AA_X$ of Lemma \ref{torilemma}, $A_{X,\Theta}^+$ is precisely the image of $(\mathfrak o\smallsetminus\{0\})^{\{\gamma\}}$.

The elements of $\mathring A_{X,\Theta}^+$ are precisely those $a\in A_{X,\Theta}=\mathcal Z(X_\Theta)$ with the property that $\lim_{n\to \infty} a^n\cdot x$ is contained in $\Theta$-infinity, for some (hence all) $x\in X_\Theta$. The elements of $A_{X,\Theta}^+$ are precisely those $a\in A_{X,\Theta}$ with the property that $\lim_{n\to \infty} a^n\cdot x$ is contained in a $G$-orbit on $\bar X$ whose closure contains $\Theta$-infinity, for some (hence all) $x\in X_\Theta$.  
\end{lemma}

\subsection{Degeneration to the normal bundle; affine degeneration}\label{ssdegen} 

It is well-known \cite[\S 2.6]{Fu} that given a closed embedding $\YY\subset \ZZ$ of varieties, there is a canonical way to deform $\ZZ$ over $\GGa$ to the normal cone over $\YY$. A multi-dimensional version of this, in the case of a simple smooth toroidal embedding $\overline\XX$ of a spherical variety $\XX$ -- hence corresponding to a cone $\mathcal C(\overline{\XX})\subset \mathcal V$ such that $\Lambda_X\cap \mathcal C(\overline{\XX})$ is a free monoid with basis $(\check\lambda_i)_i$ indexed by a set $I$ -- is a morphism:
$$\overline{\mathscr X}^n\to \GGa^I$$
carrying a $\GGm^I$-action compatible with the action on the basis. Over $\GGm^I$ it is canonically $\GG\times \GGm^I$-isomorphic to $\overline{\XX}\times\GGm^I$, and more generally the fiber over a point on the base whose non-zero coordinates are those in the subset $J\subset I$ is canonically isomorphic to the normal bundle of the orbit in $\overline{\XX}$ corresponding to $J$, with the $\GGm^{I\smallsetminus J}$-factors stabilizing that point acting on the fiber via the inverse of their canonical action on the normal bundle.\footnote{We clarify the need for ``inverse'': as we approach a non-open $\GGm^I$-orbit on $\GGa^I$, we ``stretch'' the space $\XX$ away from the corresponding $\GG$-orbit closure. It seems actually better to think of the usual $\GGm^I$-action on the normal bundle, and to compactify $\GGm^I$ ``at $\infty$'', instead of at zero; however, this would be cumbersome notationally.}

In fact, we are only interested in the union of open orbits on the fibers. As an analysis starting from the Local Structure Theorem \ref{localstructure} easily shows, $\overline{\mathscr X}^n$ contains an open $\PP(\XX)$-stable subset which is $\PP(\XX)$-equivariantly isomorphic to $\mathring\XX \times \GGa^I$ over the base, which leads to two conclusions:
\begin{lemma}\label{innormal}
\begin{enumerate}
 \item The union of open $\GG$-orbits in the fibers is an open subset $\mathscr X^n\subset \overline{\mathscr X}^n$.
 \item The $\GG\times \GG_m^I$-variety $\mathscr X^n$ is simple toroidal, with associated cone equal to the diagonal of $\mathcal C(\overline{\XX})$, where we identify the $\check\lambda_i$'s both as elements of $\Lambda_X$ and as the generators of $\mathbb N^I \subset \varchi(\GGm^I)$.
\end{enumerate}
\end{lemma}

\begin{proof}
 Indeed, the first follows immediately by acting by $\GG$ on the open $\PP(\XX)$-stable subset described before, or just by observing that $\mathscr X^n$ is just the union of $\GG$-orbits of maximal dimension on $\overline{\mathscr X}^n$, and hence open.
 
 The second follows from an easy inversion of the local structure theorem: there is a map from the toroidal embedding claimed in the lemma to $\mathscr X^n$, and since it is an isomorphism on $\PP(\XX)$-stable open subsets generating both under the $\GG$-action, it has to be an isomorphism.
\end{proof}

On the other hand, if $\XX$ is spherical and quasi-affine, and $\XX^a$ is any affine embedding of $\XX$, then $k[\XX^a]$ carries a filtration by dominant weights in $\varchi(\XX)$, partially ordered with respect to the cone  dual to $\mathcal V$ \cite[\S 6]{KnAu}, more precisely: If we decompose $k[\XX^a]$ into a direct sum of highest weight spaces, then the $\lambda$-th piece $\mathcal F_\lambda$ of the filtration is the sum of spaces with highest weights $\mu$ such that $\left< \lambda-\mu,\mathcal V\right> \le 0$. (Recall that $\mathcal V$ contains the image of the \emph{negative} Weyl chamber.) The important element here is that the maximal possible cone $\mathcal V$ with this property is precisely the same cone as that governing embeddings of $\XX$, namely the cone of invariant valuations.

There is a well-known degeneration of filtered modules to their associated graded, and multidimensional versions of it, which in the literature of spherical varieties are described over various different bases, cf.\ \cite{Popov} or \cite[\S 5.1]{GN}.  We find it preferable to define our preferred degeneration of $\XX^a$ as a morphism:
$$\mathscr X^a\to \overline{\AA_{X,ss}},$$ 
where $\overline{\AA_{X,ss}}$ is the affine embedding of the quotient\footnote{We thank Jonathan Wang for pointing out a mistake in an earlier version.} $\AA_{X,ss}$ of $\AA_X$ determined by the cone $(-\mathcal V)$ (and ``ss'' stands for ``semisimple''). In other words, the open orbit in $\overline{\AA_{X,ss}}$ is the quotient of $\AA_X$ by the subtorus generated by cocharacters in $\mathcal V\cap (-\mathcal V)$; then $(-\mathcal V)$ maps to a strictly convex cone $(-\overline{\mathcal V})$ of cocharacters of this quotient, and $\overline{\AA_{X,ss}}$ is the corresponding affine embedding. (It can be non-smooth, but this does not matter for our purposes.) 

The variety $\mathscr X^a$ is by definition the spectrum of the ring:  $$ \oplus_\lambda t^{\lambda}\mathcal F_\lambda \subset k[\XX^a \times \AA_X],$$
where the $t^\lambda$'s are symbols for the canonical basis of the group ring of $\varchi(\XX)$. Notice that the ring contains
$$ \oplus_{\left<\lambda,\mathcal V\right>\le 0} t^{\lambda} \mathcal F_0,$$
which is the coordinate ring of $\overline{\AA_{X,ss}}$. The variety $\mathscr X^a$ carries an action of $\AA_X$ compatible with its action on the base, and the coordinate ring of the fiber over a point $a\in \overline{\AA_{X,ss}}$ fixed under a subtorus $\AA_a\subset \AA_X$ is graded with respect to the character group of $\AA_a$. 

 The fiber over $1\in \AA_{X,ss}$ is canonically isomorphic to $\XX^a$, and, more generally, the restriction of the defining map:
\begin{equation}\label{sea}\XX^a\times \AA_X \to\mathscr X^a\end{equation}
to $\AA_{X,ss}$ is isomorphic to the quotient map $\XX^a\times \AA_X \to \XX^a \times \AA_{X,ss}$, although this isomorphism depends on a choice of section $\AA_{X,ss}\to \AA_X$.  On the other hand, completely canonically, we have an isomorphism:
\begin{equation}\label{modUisom} \mathscr X^a\sslash \UU \simeq \XX^a\sslash\UU \times\overline{\AA_{X,ss}},\end{equation}
simply by embedding the $\lambda$-eigenspace of the Borel in $k[\XX]$ into the $t^\lambda$ summand of $k[\mathscr X^a]$. The map $\XX\times\AA_X\to \mathscr X$ descends to a map
$$ \XX\sslash \UU \times \AA_X \to \mathscr X^a\sslash \UU = \XX\sslash \UU \times \overline{\AA_{X,ss}}$$
which is given by 
$$ (\bar x ,a)\mapsto (\bar x\cdot a, \bar a)$$
(where the bar denotes the obvious images in the quotients).

The two degenerations are closely related. More precisely, consider as above a smooth toroidal embedding determined by cocharacters $\check\lambda_i$ $(i\in I)$. We use the inverses $(-\check\lambda_i)$ of these cocharacters to obtain an injective map
$$ \GGm^I \hookrightarrow \AA_X.$$
Since the $(-\check\lambda_i)$'s are in the interior of $(-\mathcal V)$, the composition with the quotient map $\AA_X\to \AA_{X,ss}$ extends to a map
\begin{equation}\label{jj}\GGa^I \to \overline{\AA_{X,ss}}.
\end{equation}

\begin{proposition}\label{degencompatibility}
The composition 
$$ \XX\times \GGm^I \to \XX\times \AA_X \to \mathscr X^a$$
extends to a morphism:
$$ \mathscr X^n \to \mathscr X^a$$
over \eqref{jj}, which identifies: 
$$\mathscr X^n  \simeq \mathscr X \times_{\overline{\AA_{X,ss}}} \GGa^I,$$ where $\mathscr X$ is the open subset of $\mathscr X^a$ consisting of the union of open $\GG$-orbits on the fibers.
\end{proposition}

\begin{proof}
As was the case for $\mathscr X^n$, $\mathscr X^a$ also has a $\PP(\XX)$-stable open subset $\mathring{\mathscr X}$ which meets each fiber over $\overline{\AA_{X,ss}}$ in precisely the open $\PP(\XX)$-orbit and is $\PP(\XX)$-equivariantly isomorphic to $\mathring\XX\times\overline{\AA_{X,ss}}$, cf.\ \cite[Proposition 5.2.2]{GN}. (The isomorphism depends again on the choice of a section $\AA_{X,ss}\to \AA_X$.) Thus, the union of its $\GG$-translates $\mathscr X$ (the open subset consisting of the union of all open $\GG$-orbits on the fibers) is a simple toroidal embedding of the open orbit of $\GG\times \AA_X$ on $\mathscr X^a$. This embedding is described by the cone that describes $\mathring{\mathscr X}\sslash \UU$ (equivalently: $\sslash \UU_{P(X)}$) as a toric $(\AA_X\times\AA_X)/\AA_1$-variety, where $\AA_1$ is the kernel of $\AA_X\to \AA_{X,ss}$, embedded anti-diagonally. Notice that the restriction of the filtration to $k[\XX^a]^\UU$ is actually a grading. (Thus, if we choose a section $\AA_{X,ss}\to \AA_X$, we have
$$\mathscr X^a\sslash \UU \simeq (\XX^a\sslash \UU) \times \overline{\AA_{X,ss}}$$
and
$$\mathring{\mathscr X}\sslash \UU \simeq (\mathring\XX\sslash \UU) \times \overline{\AA_{X,ss}}.)$$
It is easy to see, comparing with \eqref{sea}, that the cone of the toric embedding $\mathring {\mathscr X}\sslash U$ is equal to the antidiagonal copy of any section of the quotient $\mathcal V \to \overline{\mathcal V}$ (the latter denoting the image of $\mathcal V$ in the cocharacter space of $\AA_{X,ss}$); notice that modulo cocharacters of the diagonal of $\AA_1$, the choice of section does not matter.

It now follows easily from the theory of toroidal embeddings that the map $\XX\times \GGm^I\to \mathscr X\subset \mathscr X^a$ extends to a morphism $\mathscr X^n\to\mathscr X$; it extends because the cone of $\mathscr X^n$ (spanned by the antidiagonal of the $\check\lambda_i$'s) maps to the cone of $\mathscr X$. By the same theory (or, if you prefer, by the description of open $\PP(\XX)$-orbits), the induced map 
$$\mathscr X^n  \simeq \mathscr X \times_{\overline{\AA_{X,ss}}} \GGa^I$$
is an isomorphism.
\end{proof}

{ In the previous subsection we associated a boundary degeneration $\XX_\Theta$ to any $\GG$-orbit $\ZZ$, belonging to $\Theta$-infinity, of a smooth complete toroidal embedding. The proposition above provides an alternative way to define $\XX_\Theta$, in terms of the affine degeneration, which allows us to show that the variety $\XX_\Theta$ does not depend on the choice of $\ZZ$:

\begin{proposition} \label{thetaindependent}  
For any $\GG$-orbit $\ZZ$ belonging to $\Theta$-infinity, in a smooth toroidal embedding of $\XX$, the boundary degeneration $\XX_\Theta$ defined in \S \ref{normalbundles} is canonically isomorphic to the open $\GG$-orbit in a fiber of $\mathscr X^a$ over a point of $\overline{\AA_{X,ss}}$, namely: the point $\lim_{t\to 0} \check\lambda(t)$, where $\check\lambda$ is any cocharacter into $\AA_{X,ss}$ in the interior of the face of $\overline{\mathcal V}$ corresponding to $\Theta$.
\end{proposition}

\begin{proof}
Indeed, if we let $\overline{\XX}$ be the simple toroidal subembedding where $\ZZ$ is the unique closed orbit (i.e.\ $\overline{\XX}$ consists of all orbits in the given embedding that contain $\ZZ$ in their closure), and apply Proposition \ref{degencompatibility} to the corresponding normal bundle degeneration $\overline{\mathscr X}^n$, we obtain an identification of $\XX_\Theta$ with the fiber stated in the proposition. Notice that it does not matter which affine embedding of $\XX$ we use to construct the affine degeneration $\mathscr X^a$, as the statement only uses the open $\GG$-orbit. 
\end{proof}} 

Finally, we note that the identification of open Borel orbits \eqref{Borbitident} obtained from the toroidal embedding (canonical modulo $\UU$) is compatible with the identification \eqref{modUisom} obtained from the affine degeneration. We leave the verification to the reader.

\subsection{Whittaker-type induction} \label{ssWhittakertype}  

The harmonic-analytic results of this paper apply equally well to the space of Whittaker functions, and similar spaces induced from complex characters of additive subgroups, although the ``dual group'' formalism in this case is slightly lacking at the moment. While, for simplicity, we mostly ignore these cases in our notation (for example, we write $C_c^\infty(X)$ instead of $C_c^\infty(X,\mathcal L_\Psi)$, where $\mathcal L_\Psi$ could denote the complex line bundle defined by a character of the stabilizers), the arguments carry over verbatim. Therefore, we present here the conventions that need to be used in order to translate our results to that setup.

This is not be the most general setup possible, but we will restrict ourselves to it because we do not know how to describe the ``spherical roots'' in all cases. We give ourselves a parabolic subgroup $\PP^-$ and a Levi subgroup $\LL$ of $\PP^-$. Denote by $\VV$ the vector space of homomorphisms:
$$ \UU_{P^-}\to \GGa$$

{\bluetext  We give ourselves a \emph{wavefront} homogeneous spherical variety $\XX^L$ for $\LL$, together with an equivariant map: 
$$\Lambda: \XX^L \to \VV.$$
There is a group subscheme $\ker\Lambda$ of $\UU_{P^-}\times \XX^L$ over $\XX^L$, whose fiber over $x\in \XX^L$ is $\ker(\Lambda(x))$.

Finally, we set $\XX=\XX^L \times^{\PP^-} \GG$, the spherical variety ``parabolically induced'' from $\XX^L$ to $\GG$. The map $\Lambda$ defines a principal $\GGa$-bundle with a compatible $\GG$-action over $\XX$ (by induction from the principal $\GGa$-bundle $(\UU_{P^-}\times \XX^L)/\ker \Lambda$ over $\XX^L$, which is in fact equipped with an $\LL$-equivariant trivialization). Now we fix a nontrivial additive character $\psi: k\to \CC^\times$, which defines a reduction of the $\GGa$-bundle over $\XX$ (or, rather, the associated $k$-bundle over $X$) to a $\CC^\times$-bundle, and hence a complex line bundle $\mathcal L_\Psi$ over $X$. If $x\in X^L$ and $M\subset L$ denotes its stabilizer, then sections of $\mathcal L_\Psi$, restricted to the $G$-orbit of $x$, can be identified with functions $f$ on $G$ such that: 
$$f(umg) = \psi(\Lambda(x)(u)) f(g)$$
for $u\in U_{P^-}$, $m\in M$.

Everything that follows depends on the pair $(\XX, \Lambda)$, even though there is only $\XX$ appearing in the notation; as we will see, the same variety with more degenerate characters will appear as a ``boundary degeneration'' of the pair $(\XX,\Lambda)$, so it will be denoted by the letters $\XX_\Theta$, etc. In other words, the reader should consider $\XX$ as a symbol for the pair, not just the variety. In particular, the dual group $\check G_X$ that we are about to describe \emph{is not} the dual group of $\XX$ viewed as a spherical variety; rather, it is associated to the pair $(\XX,\Lambda)$.

We let $\XX_0$ denote the \emph{total space} of the $\GGa$-bundle; again, if we fix a point $x\in X^L$ and let $\MM$ be its stabilizer in $\LL$ and $\UU_1\subset\UU_{P_0}$ the kernel of $\Lambda(x)$, then $\XX_0\simeq \UU_1\MM\backslash \GG$. It is a non-spherical variety. Friedrich Knop has associated in \cite{KnAs} a ``little Weyl group'' to $\XX_0$, which is a finite crystallographic reflection group of automorphisms of $\varchi(\XX)$; we will denote it by $W_{X_0}$, or by $W_X$ since $X$ really stands for the pair $(\XX,\Lambda)$.

We will describe a root system associated to this Weyl group. To do that, we recall from \cite{ABS} a few facts about the $\LL$-representation $\VV = \Hom(\UU_{P^-},\GGa)$. First of all, it is \emph{prehomogeneous}, i.e.\ $\LL$ acts with an open orbit. (If the image of the map $\Lambda$ lies in the open orbit, then $\Lambda$ is called \emph{generic}.) Secondly, it decomposes as a direct sum of the irreducible modules of $\LL$ with \emph{lowest} weights $\Delta\smallsetminus\Delta_L$, i.e.\ the simple roots of $\GG$ in the unipotent radical of a parabolic opposite to $\PP^-$: 
\begin{equation}\label{Vdecomp}\VV = \bigoplus_{\alpha\in \Delta\smallsetminus\Delta_L} V_{-\alpha^\vee}\end{equation}
(where $V_\gamma$ denotes the irreducible $\LL$-module with highest weight $\gamma$, and $\alpha^\vee$ is the dual weight to $\alpha$).
Equivalently, the abelianization of $\UU_{P^-}$ has highest weights $-\Delta\smallsetminus\Delta_L$, cf.\ \cite{ABS}.

\begin{lemma}\label{equiv-Whittaker}
The following are equivalent for a morphism $\Lambda$ as above and an $\alpha\in\Delta\smallsetminus\Delta_L$:
\begin{enumerate}
\item \label{Vone} In the decomposition \eqref{Vdecomp}, the image of $\Lambda$ has zero component in the $\alpha$-summand.
\item \label{Vtwo} $\Lambda$ is trivial on the unipotent radical of the parabolic (containing $\PP^-$) whose Levi has simple roots $\Delta\smallsetminus\{\alpha\}$.
\item \label{Vthree} $\Lambda$ is trivial on the subgroup $\widetilde{\UU}_{-\alpha}$ of $\UU_{P^-}$ on which the center of $\LL$ acts by the restriction of the character $-\alpha$.
\item \label{Vfour} Let $\BB$ be a Borel subgroup opposite to $\PP^-$ (i.e.\ $\BB \PP^{-}$ is open in $\GG$), $\TT$ a Cartan subgroup of $\BB$, and let $x\in \XX^L$ belong to the open $\LL\cap \BB$-orbit. 
Then the character $\Lambda(x)$ is trivial on $\UU_{-\alpha}$, the one-dimensional root subgroup of $\UU_{P^-}$ (with respect to the chosen torus) corresponding to the simple root $-\alpha$.
\end{enumerate}
\end{lemma}

We will say that $\Lambda$ is $\alpha$-trivial if the equivalent conditions of this lemma are met, and $\alpha$-generic if not.

\begin{proof}
The implications \eqref{Vtwo} $\Rightarrow$ \eqref{Vthree} $\Rightarrow$ \eqref{Vfour} are immediately clear, since the unipotent subgroup of each statement is contained in the unipotent subgroup of the previous one.

By \cite{ABS}, the subgroups of $\UU_{P^-}$ on which the center of $\LL$ does not act by the restriction of an element of $-\Delta\smallsetminus\Delta_L$ belong to the derived subgroup of $\UU_{P^-}$. This shows that \eqref{Vthree} $\Rightarrow$ \eqref{Vtwo}. Moreover, under the quotient map $\UU_{P^-} \to \UU_{P^-}^\ab$, the subgroup $\widetilde{\UU}_{-\alpha}$ maps isomorphically onto the highest weight module $V_{-\alpha}$ (dual to the module $V_{-\alpha^\vee}$ of \eqref{Vdecomp}). This shows the equivalence \eqref{Vthree} $\Leftrightarrow$ \eqref{Vone}.

Finally, \eqref{Vfour} $\Rightarrow$ \eqref{Vthree} follows from the fact that the stabilizer $\MM$ of $x$ in $\LL$ stabilizes the additive character $\Lambda(x)$, the Borel $\BB_L = \BB\cap \LL$ normalizes $\UU_{-\alpha}$ (indeed, $[U_{-\alpha}, U_\beta] = 1$ for all distinct $\alpha,\beta\in\Delta$) and $\MM\BB_L$ is open in $\LL$. Thus, $\Lambda(x)(\ell u\ell^{-1}) = 0$ for all $\ell\in \LL$, $u\in \UU_{-\alpha}$, and the $\LL$-span of $\UU_{-\alpha}$ is $\widetilde{\UU}_{-\alpha}$.
(Compare with \cite[Lemma 6.1.1]{SaSph}.) Thus, $\Lambda(x)|_{\widetilde{\UU}_{-\alpha}} = 0$, and by homogeneity the same holds for $\Lambda(x')$, for all $x'\in \XX^L$.
\end{proof}

\begin{definition}\label{defWhittakerroots}
 Define the set $\Delta_X$ of (normalized) spherical roots of $\XX$ (really, of the pair $(\XX,\Lambda)$ or of the non-spherical variety $\XX_0$) as the union of the set of (normalized) spherical roots of $\XX^L$ and the set of $\alpha\in \Delta\smallsetminus\Delta_X$ for which $\Lambda$ is $\alpha$-generic. 
\end{definition}

\begin{proposition}\label{sphericalWhittaker}
 The set $\Delta_X$ belongs to the character group $\varchi(\XX)$, and forms a set of simple positive roots for a root system with little Weyl group $W_X$, i.e.\ Knop's Weyl group for the non-spherical variety $\XX_0$.
\end{proposition}

Notice that the character group $\varchi(\XX)$ used here is the same as the one we have defined for the spherical variety $\XX$. (The morphism $\Lambda$ plays no role.)

\begin{proof}
If $\Lambda$ is $\alpha$-trivial for some $\alpha\in \Delta\smallsetminus \Delta_L$, the variety $\XX_0$ is parabolically induced from the parabolic containing $\PP^-$ whose Levi has simple roots $\Delta\smallsetminus\{\alpha\}$, by property \eqref{Vtwo} of Lemma \ref{equiv-Whittaker}. 
Knop's little Weyl group $W_X=W_{X_0}$ generalizes the little Weyl groups of spherical varieties and, in particular, has the property that if a variety is parabolically induced, it is equal to the little Weyl group of the inducing variety. Thus, we are reduced to the case where $\Lambda$ is $\alpha$-generic for all $\alpha\in \Delta\smallsetminus \Delta_L$. In this case, $\Delta_X = \Delta_{X^L} \cup (\Delta\smallsetminus \Delta_L)$.

First we prove that  the elements of $\Delta\smallsetminus \Delta_L$ belong to $\varchi(\XX)$, i.e.\ are characters of $\AA_X$. Let $\BB$ be a Borel subgroup opposite to $\PP^-$, and $x\in \XX^L$ in the open $\LL\cap \BB$-orbit. Choose $\TT=\AA$ a Cartan subgroup chosen as in \S \ref{invariants}, i.e.\  such that it acts via the quotient $\AA_X$ on the orbit of $x$. By property \eqref{Vfour} in Lemma \ref{equiv-Whittaker}, the character $\Lambda(x)$ is nontrivial on the root subgroup $\UU_{-\alpha}$. This character is stabilized by the kernel of $\AA\to\AA_X$, hence the character $\alpha$ is trivial on this kernel, i.e.\ $\alpha\in \varchi(\XX)$.

To prove that $\Delta_X$ forms a set of simple positive roots for a root system with little Weyl group $W_X$, we first show that the linearly independent set of roots $\Delta_X$ determines a Weyl chamber for the action of $W_X$ on $\mathfrak a_X  = \varchi(\XX)\otimes\QQ$. 

Friedrich Knop has defined in \cite{KnOrbits} an action of the full Weyl group on a certain set of $\BB$-stable subsets of $\XX_0$, which includes $\XX_0$ itself. As was remarked in \cite[\S 5.4]{SaSpc}, in the current setting the simple reflection corresponding to the root $\alpha$ belongs to the stabilizer of $\XX_0$, and more precisely to the little Weyl group $W_X$.  From Knop's construction, it is immediate to see that the same is true for the simple reflections associated to elements of $\Delta_{X^L}$. Thus, a Weyl chamber for $W_X$ is contained in the cone $\mathcal C$ negative-dual to the set of characters $\Delta_X$.

For the converse, since $\XX^L$ was assumed to be wavefront, the map of anti-dominant chambers: $\mathfrak a^+_L \to \mathfrak a_{X^L}^+$ is surjective; here $\mathfrak a^+_L$ denotes the $L$-anti-dominant elements of $\mathfrak a$. The subset $\mathfrak a^+ \subset \mathfrak a^+_L$ of $\GG$-anti-dominant elements is defined by the additional conditions: $\left<\alpha,a\right>\le 0$ for all $\alpha\in\Delta\smallsetminus\Delta_L$, and similarly the subcone $\mathcal C \subset \mathfrak a_{X^L}^+$ is defined by the same additional conditions. Therefore, the map:
$$ \mathfrak a^+ \to \mathcal C$$
is surjective. On the other hand, the image of $\mathfrak a^+ $ in $\mathfrak a_X$ is contained in a Weyl chamber for the little Weyl group $W_X$; therefore, $\mathcal C$ coincides with that Weyl chamber.

We have shown that the linearly independent set of roots $\Delta_X$ determines a Weyl chamber for the action of $W_X$ on $\mathfrak a_X  = \varchi(\XX)\otimes\QQ$. The fact that the $W_X$-translates of $\Delta_X$ form a root system with this Weyl group follows if we show that whenever an element of $W_X$ carries one of these half-lines to another, it has to map the corresponding spherical roots to each other. This, in turn, relies on the following fact, which we will use more extensively in \S \ref{ssrootdatum}, and therefore we point the reader there for proofs and definitions: each element of $\Delta_X$ is either a simple root of $\GG$ or the sum of two strongly orthogonal roots. This fact uniquely determines it on the half-line that it spans.
\end{proof} 

In accordance with this definition of $\Delta_X$ we define all the corresponding invariants of $\XX$ (that is, of the pair $(\XX,\Lambda)$) as in \S \ref{notation}. For example, $\mathfrak a_X^+$ denotes the subset of $\varchi(\XX)^*\otimes \mathbb Q$ of elements which are anti-dominant with respect to $\Delta_X$; equivalently, the subset of $\mathfrak a_{X^L}^+$ of elements which are $\le 0$ on all $\alpha\in \Delta\smallsetminus\Delta_L$ for which $\Lambda$ is $\alpha$-generic. As we saw in the last proof, since $\mathfrak a_L^+\to \mathfrak a_{X^L}^+$ is assumed to be surjective (wavefront property), the same is true for the map: $\mathfrak a^+\to \mathfrak a_X^+$ \emph{when $\Lambda$ is generic} (that is, when $\Delta\smallsetminus\Delta_L\subset\Delta_X$), i.e.\ generic Whittaker-induction in this sense preserves the wavefront property.}

Now we discuss ``wonderful compactifications''. Again, the name will be applied more generally to smooth toroidal embeddings -- however, they will not be complete, since, as we will see smooth sections of $\mathcal L_\Psi$ vanish in certain directions. Hence, in our setting, a ``wonderful compactification'' of $\XX$ (taking into account the character $\Lambda$) will be a smooth toroidal embedding $\overline{\XX}$ of $\XX$ whose fan has support (=the union of its cones) equal to $\mathfrak a_X^+$. It is easy to see that such an embedding is of the form:
$$ \overline{\XX} = \overline{\XX^L}\times^{\PP^-}\GG,$$
where $\overline{\XX^L}$ is a toroidal embedding of $\XX^L$ defined by the same fan. The reason for this definition is the following lemma, which will be proven in \S \ref{ssCartan}.

\begin{lemma} \label{supportLpsi}
 The support of any element of $C^\infty(X,\mathcal L_\Psi)$ has compact closure in $\bar X$, where $\bar X$ denotes any ``wonderful'' embedding as described above. 
\end{lemma}
Hence, for our purposes the space $\bar X$ is as good as a compact space, since the support of smooth sections cannot escape in other directions. 

On the other hand, in the other directions the line bundle $\mathcal L_\Psi$ can be extended to the ``wonderful'' embedding:
\begin{lemma}
 The $\GGa$-bundle $\XX_0\to\XX$ extends to a $\GGa$-bundle (with a compatible action of $\GG$) over a wonderful embedding $\overline{\XX}$.
\end{lemma}

\begin{proof}
If $\overline{\XX^L}$ denotes the closure of $\XX^L$ in $\overline{\XX}$, then, as mentioned above, $\overline{\XX} = \overline{\XX^L} \times^{\PP^-}\GG$. Therefore, it is enough to show that the $\GGa$-bundle extends to $\overline{\XX^L}$ or, equivalently, that the morphism 
$$\Lambda: \XX^L\to \VV=\Hom(\UU_{P^-},\GGa)$$
extends.

Let $\mathfrak a_V$ be the dual to the vector space spanned by the $\BB$-eigencharacters on $k(\VV)$, and let $\mathcal C\subset \mathfrak a_V$ be the cone dual to the set of eigencharacters of \emph{regular} Borel eigenfunctions. This is the cone attached by Luna-Vust theory \cite{KnLV} to the spherical embedding $\VV$, since every $\BB$-stable divisor on $\VV$ is $\GGm$-stable and hence contains the closed orbit $\{0\}$. By \cite[Theorem 4.1]{KnLV}, to show that the map of spherical varieties $\XX^L\to\VV$ extends to $\overline{\XX^L}$, it is enough to show that the support of the fan of $\overline{\XX^L}$, i.e.\ the cone $\mathfrak a_X^+$, maps into $\mathcal C$ under the natural quotient map: $\mathfrak a_X\to\mathfrak a_V$.

By the decomposition \eqref{Vdecomp}, the coordinate ring $k[\VV]$ is the symmetric algebra $S^\bullet (\oplus_{\alpha \in -\Delta\smallsetminus\Delta_L} V_\alpha)$, where $V_\alpha$ is the irreducible $\LL$-module with highest weight $\alpha$. { Decomposing these symmetric powers in irreducible modules for $\GG$, we get highest weights that belong to the cone spanned by $-\Delta\smallsetminus\Delta_L$ and the relative root cone of $\LL$, i.e.\ all the highest weights in $k[\VV]$ are contained in the negative root cone of $\GG$. Hence, $\mathfrak a^+$ maps into $\mathcal C$, and since $\mathfrak a^+$ surjects onto $\mathfrak a_X^+$, the same holds for the latter.}
\end{proof} 

Based on this lemma, now, we can define the boundary degenerations $\XX_\Theta$, which are really symbols for a pair $(\XX_\Theta,\Lambda_\Theta)$, where $\XX_\Theta = \XX_\Theta^L\times^{\PP^-}\GG$ is a homogeneous spherical $\GG$-variety induced from $\PP^-$ (possibly: $\XX_\Theta\simeq \XX$ as $\GG$-varieties), and $\Lambda_\Theta: \XX_\Theta^L\to \VV$, where $\VV$ is as before.

{\bluetext Namely, as varieties $\XX_\Theta$ and $\XX_\Theta^L$ are defined as previously -- the open $\GG$-orbit, resp.\ $\LL$-orbit, in the normal bundle to a certain orbit on $\overline{\XX}$, resp.\ $\overline{\XX^L}$. If $\ZZ\subset \overline{\XX^L}$ is the orbit corresponding to $\XX_\Theta^L$ we have a canonical quotient map (quotient by $\AA_{X,\Theta}$): $\XX_\Theta^L\to\ZZ$, and we let $\Lambda_\Theta$ be its composition with the map $\ZZ\to \VV$ obtained, by the previous lemma, as the extension of $\Lambda$. }

If $\Omega = \Theta\sqcup\{\alpha\}$, where $\alpha$ is a simple root in the unipotent radical of a parabolic opposite to $\PP^-$, then $\XX_\Theta^L$ and $\XX_\Omega^L$ are isomorphic as varieties, but $\Lambda_\Theta$ is trivial on the summand $V_{-\alpha}$ under the decomposition
$$ (\UU_{P^-})^\ab \simeq \oplus_{\beta\in -\Delta\smallsetminus\Delta_L} V_\beta$$
dual to \eqref{Vdecomp} (equivalently, on the subgroup $\widetilde{\UU}_{-\alpha}$ of Lemma \ref{equiv-Whittaker}.)

\subsection{Levi varieties} \label{levivarieties} While the statements of our representation-theoretic results can be formulated without any reference to the internal structure of the boundary degenerations $\XX_\Theta$, the tools that we have at our disposal are unfortunately related to parabolic induction and restriction from Levi subgroups. This is a basic reason why we restrict ourselves to wavefront spherical varieties (s.\ Proposition \ref{wavefrontlevi} below). We start with a general statement that does not use the wavefront property:

\begin{lemma}\label{lemmalevi}
Let $\LL_\Theta,\PP_\Theta$ denote the standard Levi with simple roots $\tilde\Theta:=\Delta(\XX)\cup\supp(\Theta)$, and the corresponding standard parabolic. Let $\PP_\Theta^-$ denote the parabolic opposite to $\PP_\Theta$ which contains $\LL_\Theta$. There exists a spherical variety $\XX_\Theta^L$ of $\LL_\Theta$ such that $\XX_\Theta\simeq\XX_\Theta^L\times^{\PP_\Theta^-}\GG$.
\end{lemma}

Such a process of constructing a spherical variety of $\GG$ from one of a Levi subgroup is called ``parabolic induction''. 

\begin{proof}
 By \cite[Proposition 3.4]{Lu}, a homogeneous spherical variety $\YY$ for $\GG$ is induced from a parabolic $\PP^-$ (assumed opposite to a standard parabolic $\PP$) if and only if $\supp(\Delta_Y)\cup \Delta(\YY)$ is contained in the set of simple roots of the Levi subgroup of $\PP$.
\end{proof}

The variety $\XX_\Theta^L$ will be called a \emph{Levi variety} for $\XX$. It is a \emph{canonical} subvariety of $\XX_\Theta$ once the parabolic $\PP_\Theta^-$ has been chosen (in particular, once we choose a standard Borel and a Levi $\LL(\XX)$ as explained in \ref{invariants}), namely the subvariety of points stabilized by $\UU_\Theta$. Notice, moreover, that its inclusion in $\XX_\Theta$, indeed in $\mathring \XX_\Theta\cdot \PP_\Theta$, allows us to \emph{canonically} identify it with the quotient $\mathring \XX \PP_\Theta/\UU_{\PP_\Theta}$. In general, we could have $\LL_\Theta = \GG$; in the wavefront case, however, the boundary degenerations are parabolically induced from sufficiently small Levi subgroups:

\begin{proposition} \label{wavefrontlevi}
 Assume that $\XX$ is wavefront and that the map $\mathcal Z(\GG)^0 \rightarrow \mathcal{Z}(\XX)$
 is surjective.  Then the natural map: $$\mathcal Z(\LL_\Theta)^0\to \AA_{X,\Theta}=\Aut_\GG(\XX_\Theta)^0$$ (``action on the left'') is surjective.   In particular, for $\Theta\subsetneq \Delta_X$ the Levi $\LL_\Theta$ is proper.

This condition characterizes wavefront spherical varieties: if a spherical variety is not wavefront, then there exists a $\Theta\subset\Delta_X$ such that the above map is not surjective. In particular, if $\XX$ is wavefront then $\XX_\Theta$ is also wavefront {\bluetext under the action of $\mathcal Z(\LL_\Theta)^0\times \GG$}, for every $\Theta\subset\Delta_X$.
\end{proposition}

\begin{proof}

Since both groups are tori, it suffices to prove surjectivity for the corresponding map of their cocharacter groups, tensored by $\RR$.  We will prove somewhat more, namely, the map induces a surjection of positive chambers at this level. 

Let $\mathfrak{z}$ (resp. $\mathfrak{z}_X$) be the kernel of all roots on $\mathfrak{a}$ (resp. the kernel of all spherical roots on $\mathfrak{a}_X$).  The cone $\mathfrak{a}^+$ (resp. $\mathcal{V}$) is the preimage in $\mathfrak{a}/\mathfrak{z}$ (resp. $\mathfrak{a}_X/ \mathfrak{z}_X$)
of the convex hull of the half-lines spanned by the negative dual basis $\{-\alpha^*: \alpha \in \Delta\}$
(resp. $-\gamma^*: \gamma \in \Delta_X$) to $\Delta$ (resp. $\Delta_X$).  

Since the map $\mathfrak{a}^+ \rightarrow \mathfrak{a}_X^+$ is surjective, 
it follows that for every $\gamma \in \Delta$ there exists $\alpha \in \Delta$ so that 
$\mathbb{R}\alpha^* \twoheadrightarrow \mathbb{R}_{+} \gamma^*$ under the induced $\mathfrak{a}/\mathfrak{z} \rightarrow \mathfrak{a}_X/\mathfrak{z}_X$. In fact, this latter condition is {\em equivalent} to  the surjectivity of $\mathfrak{a}^+ \rightarrow \mathfrak{a}_X^+$ if we suppose (as we have been doing) that $\mathfrak{z} \rightarrow \mathfrak{z}_X$ is surjective; we can rephrase it thus:
\begin{equation}\label{wavefrontequiv}
\minibox{For every $\gamma \in \Delta_X$, there exists $\alpha_{\gamma} \in \Delta$ that is contained in the support of\\ $\gamma$,  but not in the support of any other simple root.}
\end{equation}

We claim that we may always choose $\alpha_{\gamma}$ so that it does not belong to $\Delta_{L(X)}$.
Suppose to the contrary: Let $S$ be the set of all roots in the support of $\gamma$, but not in the support of any other root, and suppose $S \subset \Delta_{L(X)}$. Every $\alpha \in S$
is orthogonal to the support of every spherical $\gamma' \neq \gamma$: 
Indeed, by non-degeneracy (Lemma \ref{nondeg}), $\langle \alpha^{\vee}, \gamma' \rangle = 0$,  
but $\langle \alpha^{\vee}, \beta \rangle \leq 0$ for all $\beta \in \mathrm{supp}(\gamma')$. 
Consequently, every $\alpha \in S$ is orthogonal to every element of $\supp(\gamma) - S$. 
If we write $\gamma = \sum_{\alpha} n_{\alpha} \alpha$, then 
$\gamma_S := \sum_{\alpha \in S} n_{\alpha} \alpha$ has the property that 
$\langle \alpha^{\vee}, \gamma_S \rangle = 0$ for all $\alpha \in S$. 
This contradicts the fact that the matrix $\langle \alpha^{\vee},\beta \rangle_{\alpha , \beta \in S}$ is nondegenerate.

The claims of the proposition follow immediately. Indeed, 
the convex hull of  $\{ \mathbb{R}_+ \alpha_{\gamma}^*: \gamma \notin \Theta\}$ 
is orthogonal to all roots in $\tilde{\Theta}$ and surjects onto $\mathfrak{a}_{X, \Theta}^+$. 

On the other hand, if $\XX$ is not wavefront,  equivalently: if \eqref{wavefrontequiv} fails,  then there is a spherical root $\gamma$ such that every simple root in the support of $\gamma$ is also contained in the support of some other spherical root. Hence, $\LL_{\Delta_X\smallsetminus\{\gamma\}} = \GG$, but $\mathcal Z(\GG)^0$ has image $\mathcal Z(\XX)\subsetneq \mathcal Z(\XX_\Theta)$.
\end{proof}

\subsection{Horocycle space} \label{sshorocycles} 
Let $\Phi$ denote any conjugacy class of parabolic subgroups of $\GG$, whose representatives contain a conjugate of $\PP(\XX)$ as a subgroup. We define the space of $\Phi$-horocycles on $\XX$ as the $\GG$-variety $\XX^h_\Phi$ classifying pairs: 
$$(\Q,\mathfrak O),$$
where $\Q$ is a parabolic in the class $\Phi$ and $\mathfrak O$ is an orbit of $\UU_Q$ contained in the open $\Q$-orbit on $\XX$ (it will also be denoted by $\XX_Q^h$). More explicitly, if we choose a Borel and hence a standard representative $\Q$ for the class $\Phi$ then, canonically: 
$$\XX_{Q}^h :=  \XX_\Phi^h \simeq \mathring \XX \cdot \Q / \UU_Q \times^\Q  \GG,$$
where, as usual, $\mathring\XX$ denotes the open orbit of the chosen Borel subgroup.

Now let $\Theta$ be a subset of spherical roots, take $\Q = \PP_{\Theta}$, and denote $\XX^h_{\Theta}:=\XX^h_Q$. The next lemma compares the space of $\Q$-horocycles with the analogous space for the boundary degeneration $\XX_\Theta$:

\begin{lemma}\label{importantobs} 
 If $\XX$ is a wavefront spherical variety, then there is a canonical identification:
\begin{equation} \label{importantobs-ID} \XX^h_{\Theta}:=\XX^h_Q \stackrel{\sim}{\rightarrow} (\XX_{\Theta})_Q^h, \end{equation} 
compatible with the identification of open Borel orbits (\ref{Borbitident}).
\end{lemma}

\begin{proof}
Let us start with $\Theta=\emptyset$, i.e.\ $\Q=\PP_\Theta=\PP(\XX)$. Then the Local Structure Theorem \ref{localstructure} and the canonical identification \eqref{Borbitident} immediately imply that:
$$ \mathring \XX/\UU_\Theta = \mathring \XX_\Theta/\UU_\Theta$$
canonically, and hence $\XX_Q^h\simeq (\XX_\Theta)_Q^h$.

In the general case, we first would like to exhibit $\YY_\Theta:=\mathring\XX_\Theta\PP_\Theta/\UU_\Theta$ as a boundary degeneration of $\YY:=\XX_\Theta^L=\mathring\XX \PP_\Theta/\UU_\Theta$. In a smooth toroidal embedding $\overline{\XX}$ of $\XX$, let $\XX_1$ be the $\PP_\Theta$-stable subset of those points whose $\PP_\Theta$-orbit closure contains $\Theta$-infinity. (For simplicity, we assume that $\Theta$-infinity is a unique orbit.) Then $\UU_\Theta$ acts freely on $\XX_1$, and $\XX_1/\UU_\Theta$ is a simple embedding of $\YY$ such that the open $\LL_\Theta$-orbit in the normal bundle to the closed orbit is canonically isomorphic to $\YY_\Theta$. In terms of the combinatorial description of spherical embeddings, this is the toroidal embedding of $\YY$ with valuation cone equal to the face of $\mathcal V$ that is orthogonal to $\Theta$.

Thus, according to \S \ref{ssdegen}, $\YY_\Theta$ is the open $\LL_\Theta$-orbit in the spectrum of the partial grading of $k[\YY]$ induced by those valuations. But as we saw in Proposition \ref{wavefrontlevi}, the grading induced by these valuations 
coincides (or can be refined by) the grading by characters of $\mathcal Z(\LL_\Theta)$. Therefore, $k[\YY]$ is already graded and there is nothing to degenerate: $k[\YY_\Theta]\simeq k[\YY]$, and hence $\YY_\Theta\simeq \YY$ (canonically). Hence, $\XX_Q^h=(\XX_\Theta)^h_Q$.
\end{proof}

\subsection{The example of $\PPGL(V)$ as a $\PPGL(V) \times \PPGL(V)$ variety } \label{examples}

We shall describe a basic example and compute much of the foregoing data as
an illustration.

Let $V$ be a vector space, and take $\XX = \PPGL(V)$ considered as 
$\GG = \PPGL(V) \times \PPGL(V)$ variety via $A  \cdot(g,h)= g^{-1} A h$.   Let $x_0 \in \XX$ be the homothety class of scalar multiplication.  In the discussion that follows, we understand $\GL(V)$ as acting on $V$ on the  right.

\subsubsection{The case $n=2$} \label{axpluspgl}

We use this as an example to explain the behavior of $A_X^+$ in geometric terms, 
i.e.\ in terms of the asymptotic behavior of horocycles. 

In this case, the compactification of $\PGL(2)$ is simply that induced from its embedding
into $\overline{\XX} := \mathbb{P}(M_2)$, where $M_2$ is the algebra of $2 \times 2$ matrices.

We take the Borel subgroup $\BB \subset \GG$ to be $\BB^+ \times \BB^-$, where $+$ and $-$ denote respectively upper and lower triangular matrices.  If we let $x_0 \in \XX$
be the identity automorphism, then $x_0$ is in the open $\BB$-orbit. This orbit
consists of the lines of elements $\left(\begin{array}{cc} x & y \\ z & w \end{array}\right) \in M_2(k)$ with $w \neq 0$.   This space is foliated by the $U$-orbits
  $$ \mathfrak{h}_c := \left\{  \left( \begin{array}{cc} * & * \\ * & d  \end{array} \right) : \frac{d^2}{\det} = c \right\}.$$ 
In terms of the geometry of $\overline{\XX}$, the horocycle $\mathfrak{h}_c$ has quite different behavior as $c \rightarrow 0$ and as $c \rightarrow \infty$:
\begin{itemize}
\item[-] As $c \rightarrow \infty$, 
the entire horocycle $\mathfrak{h}_c$ draws close {\bluetext to the divisor of singular matrices in $\overline{\XX}$, which is the closure of} the horocycle of (lines of) singular matrices of the form
  $\left( \begin{array}{cc} * & * \\ * & d \end{array}\right)$ with $d\ne 0$.

\item[-] On the other hand,
as $c \rightarrow 0$, the horocycle $\mathfrak{h}_c$  converges rather to the $\BB$-stable divisor (small Bruhat cell) on $\XX$ itself.
\end{itemize}
  In this way, the two ``directions'' in $A_X$ are distinguished from one another. 
  
  More precisely: 
  The action of $\BB$ on these horocycles factors through the quotient
  $\BB \rightarrow \AA \twoheadrightarrow \AA_X$ and defines a faithful $\AA_X$-action;
  in this case, it is explicitly given by $$ \left( \begin{array}{cc} x & * \\ 0 & 1 \end{array} \right), 
  \left(\begin{array}{cc} y^{-1} & 0\\ * & 1 \end{array}\right): \mathfrak{h}_c \mapsto \mathfrak{h}_{c \cdot (xy)}.$$
  Thus, the action of a positive one parameter subgroup $\lambda: \GGm \rightarrow \AA$
  (which projects to the {\em negative} of the valuation cone in $X_*(A_X)$) satisfies $\lambda(t) \mathfrak{h}_c \rightarrow \overline{\XX}-\XX$
  as $t \rightarrow \infty$, whereas a negative one parameter subgroup (which projects to the valuation cone of $X_*(A_X)$) has the opposite behavior.

\subsubsection{The compactification for general $n$} 
For $n > 2$ the wonderful compactification of $\XX$ is more subtle. It is classically known as the variety of {\em complete collineations}; a modern treatment is given by Thaddeus \cite{Th} 
or de Concini and Procesi \cite{DCP}.

  A set $\Theta\subset \Delta_X$ can be identified with a \emph{flag type} 
  in this case, namely an increasing sequence of dimensions: 
\begin{equation}\label{sequence}0= d_0 < d_1 < \dots < d_k < d_{k+1}=\dim(V).
\end{equation}

Then the ``boundary degeneration'' $\XX_\Theta$ can be described as classifying triples:
$$(K, I,  \phi),$$ 
where:
\begin{itemize}
\item[-]$K: V=K_0 \supset K_1 \supset K_2 \supset \dots \supset K_k \supset K_{k+1}=V$ is a decreasing flag with $\codim K_i=d_i$ -- the ``kernel flag'';
\item[-]$I: \{0\} = I_0\subset I_1 \subset I_2 \subset \dots \subset I_k \subset I_{k+1}=V$ is an increasing flag with $\dim I_i=d_i$ -- the ``image flag'';
\item[-]
  $\phi$ is the homothety class of a graded isomorphism:
$$\phi: \mathrm{gr}^K(V) \rightarrow \mathrm{gr}^{I}(V).$$
\end{itemize}

Explicitly, $\phi$ is a collection of linear isomorphisms:
$$\phi_i: K_{i-1}/K_i \rightarrow I_i/I_{i-1},$$
determined up to a common scalar multiple. The $\GG$-automorphism group of $\XX_\Theta$ is generated by scalar multiplications of the individual $\phi_i$'s.
 
On the other hand, the corresponding $\GG$-orbit $\ZZ_\Theta$ of the wonderful compactification can be identified with the variety classifying triples $(K,I,[\phi])$ as before, but with \emph{each} of the morphisms $\phi_i$ comprising $\phi$ defined up to homothety (and hence we denote them as $[\phi_i],\phi$).

We describe how these orbits are glued topologically in the wonderful compactification (either in the sense of Zariski topology, or in the sense of the usual topology of points over a $p$-adic field). It suffices to describe a basis of neighborhoods of a point $z\in \ZZ_\Theta$ inside of $\ZZ_{\Theta'}$, when $\ZZ_\Theta$ is contained in the closure of $\ZZ_{\Theta'}$ and of codimension one. In other words, if $\Theta$ corresponds to a sequence of integers as in (\ref{sequence}) then ${\Theta'}$ corresponds to the same sequence with a $d_i$ removed. But then, the description of neighborhoods is in complete analogy to the highest-dimensional case, namely when $\Theta$ corresponds to a sequence $0<d<\dim(V)$ and hence $\ZZ_{\Theta'}$ is the open $\GG$-orbit (i.e.\ $\PPGL(V)$).

In that case, for a point $z=(K,I,[\phi])\in \ZZ_\Theta$, $K$ is a sequence $V\supset K_1\supset 0$, and $I$ is a sequence $0\subset I_1\subset V$. Moreover the datum $[\phi]$ consists of homothety classes of isomorphisms: $[\phi_1: V/K_1 \xrightarrow{\sim} I_1]$ and $[\phi_2: K_1\xrightarrow{\sim} V/I_i]$.

A \emph{neighborhood of $z$} in $\PPGL(V)$ then consists of homothety classes $[g]$ of automorphisms $g: V\to V$ with the property that $[g]$ is in a neighborhood, in $\mathbf P\End(V)$, of the endomorphism $$V \to V/K_1 \xrightarrow{\phi_1} I_1 \hookrightarrow V,$$ (or rather its homothety class), and $[g^{-1}]$ is in a neighborhood, in $\mathbf P\End(V)$, of the endomorphism: $$V\to V/I_1 \xrightarrow{\phi_2} K_1\hookrightarrow V$$
(or rather its homothety class). 

To describe the algebraic structure of the whole wonderful variety is a little more complicated, but notice that these codimension-one orbits already show up in the closure of the embedding: $$ \PPGL(V) \ni [g]\to ([g],[g^{-1}]) \in \mathbf P \End(V) \times \mathbf P\End(V).$$
(Notice also that the algebraic structure of the wonderful compactification is easily obtained from a different construction, namely embedding $\PPGL(V)$ in $\mathrm{Gr}(\mathfrak g)$ -- the Grassmannian of the Lie algebra of $\GG$-- as the $\GG$-orbit of the diagonal: $\mathfrak{pgl}(V)\hookrightarrow \mathfrak g$.)

On the other hand, a neighborhood of the above point $z\in \ZZ_\Theta$ in $\XX_\Theta$ can be described as the set of triples $(I',K', \phi')$ with $I', K'$ in neighborhoods of $I$ and $K$ in the corresponding flag varieties, and $\phi'=(\phi_1',\phi_2')$ such that it(s homothety class) is in a neighborhood of (the homothety class of) $(\phi_1,0)$, and such that (the class of) ${\phi'}^{-1}$ is in a neighborhood of (the class of) $(0,\phi_2^{-1})$. Of course, when varying the point $(I',K')$ we take into account that $\phi'$ is naturally an element of some bundle over the corresponding product of flag varieties, and hence ``in a neighborhood'' makes sense for $\phi'$.

 \subsubsection{Identification of orbits}

Let $k$ be a $p$-adic field, and $J$ an open compact subgroup of $\GG(k)$. 
Later in this paper (\S \ref{subsec:expmap}) we describe how to identify $J$-orbits in a neighborhood of $\Theta$-infinity on $\XX_\Theta(k)$ with $J$-orbits on a neighborhood of $\Theta$-infinity on $\XX(k)$.  In preparation for that, let us describe now how to do this explicitly in the example of $\XX=\PPGL(V)$:

\begin{enumerate}
  \item  To go from $X(k)$ to $X_{\Theta}(k)$: 
 
Given an element $\bar{A} \in X(k)$, choose a representative $A \in \GL(n,k)$. According to the Cartan decomposition,
there exists bases $e_1, \dots, e_n$ and $f_1,\dots, f_n$ for $\mathfrak{o}^n$
with the property that $A e_i = \lambda_i f_i$. 
We order the $\lambda_i$ so that
\begin{equation} \label{CD} |\lambda_1| > |\lambda_2| > \dots |\lambda_n|.\end{equation}
The bases $e_i$, $f_i$ are unique
up to  ``corrections'' belonging to  $\GL(n, \mathfrak{o}) \cap D^{\pm 1} \GL(n, \mathfrak{o}) D^{\mp 1}$, where $D = \mathrm{diag}(\lambda_1, \dots, \lambda_n)$.

Fix $\Theta \subset \Delta_X$, which we identify as before
with a flag type $ \{0=d_0 < d_1 < d_2 < \dots <  d_k = n\}$. Say $A$ is $\Theta$-large
if $|\lambda_{d_s}/\lambda_{d_s+1}| \geq T$ for each $s$.

 Set 
$$K_1 = \langle e_n, \dots, e_{d_1+1} \rangle \supset K_2 = \langle e_n, \dots, e_{d_2} \rangle \supset \dots$$
$$ I_1 = \langle f_1, \dots, f_{d_1} \rangle \subset I_2 = \langle f_1, \dots, f_{d_2} \rangle \subset \dots$$
and we associate to $A$ the flags $K_*, I_*$ and the natural isomorphism $\mathrm{gr}^K \rightarrow \mathrm{gr}^I$
induced by $A$.  Although $e_i, f_i$ are not uniquely defined, the resulting map nonetheless gives a well-defined map
$$ \mbox{$J$-orbits on $\Theta$-large $A$} \longrightarrow \mbox{$J$-orbits on $X_{\Theta}$}.$$
so long as $T$ is sufficiently large relative to $J$; this follows easily from the uniqueness
of Cartan decomposition (in the sense noted above).

\item 
To go from $X_{\Theta}(k)$ to $X(k)$:

First the notion of ``$\Theta$-large'' on $X_{\Theta}(k)$: 
Let $z=(K,I,[\phi]) \in X_{\Theta}(k)$. $[\phi]$ is a class of morphisms
up to common homothety; we choose a representative $\phi$ from this class.

 We derive integral structures on $K_{i-1}/K_{i}$ and $I_i/I_{i-1}$
 as the images of $\mathfrak{o}^n \cap K_{i-1}$ and $\mathfrak{o}^n \cap I_i$; now
  let $\lambda_{\min}(\phi_i)$ and $\lambda_{\max}(\phi_i)$
be the smallest and largest  ``singular values'' 
of the Cartan decomposition of
 $\phi_i:  K_{i-1}/K_i \rightarrow I_i /I_{i-1}$.
  We then say that $z$ is $\Theta$-large if, for every $j$, 
$$|\lambda_{\min}(\phi_j) | > T \cdot | \lambda_{\max}(\phi_{j+1}) |.$$

  Choose once and for all a splitting 
of all flags of type $\Theta$, in such a way that
these splittings vary continuously in the (compact) space of such flags under the $k$-adic topology.  
 The chosen splittings give identifications
$\mathrm{gr}^K V \stackrel{\sim}{\rightarrow} V$ and similarly for $I$. In particular,
$\phi$ induces a linear map $V \rightarrow V$, that is to say,
an element $\tilde{z} \in X(k)$. 

Fix an open compact subgroup $J$. 
If we had chosen a different choice of splitting,
the resulting elements $\tilde{z}, \tilde{z}' \in X(k)$ nonetheless
still lie in the same $J$ orbit so  long as $T$ is chosen large enough. 

We have therefore obtained a map
$$ \mbox{$J$-orbits on $\Theta$-large $A \in X_{\Theta}(k)$} \longrightarrow \mbox{$J$-orbits on $X$}.$$

\end{enumerate}

	\section{Proofs of the results on the dual group} \label{sec:dualgroupproofs}  
 
 In this section we prove  the results of \S \ref{dualgroup}, including Theorem \ref{sl2}. We use heavily the work of F.~Knop and M.~Brion to define $\check G_X$, and the work of Gaitsgory-Nadler to construct a morphism
 $\check G_X \times \SL_2 \rightarrow \check{G}$. 

		\subsection{The root datum of a spherical variety} \label{ssrootdatum}

 First of all, we recall that F.\ Knop has defined \cite{KnOrbits} an action of the Weyl group on the set of Borel orbits on $\XX$; the stabilizer of the open orbit is equal to the group $W_{(X)}:=W_X\ltimes W_{L(X)}$. Here the group $W_X$, which is originally a subgroup of $\End(\varchi(\XX)\otimes\mathbb Q)$, is identified with its ``canonical lift'' to $W$, which consists of representatives of minimal length modulo $W_{L(X)}$ (the set of those representatives to be denoted by $[W/W_{L(X)}]$). In order to distinguish between the two, we will denote the subgroup of $\End(\varchi(\XX)\otimes\mathbb Q)$ by $\overline{W_X}$. Knop's action has the property \cite[Theorem 4.3]{KnOrbits} that $\varchi({^wY})=w\cdot \varchi(Y)$ (in characteristic zero).

Recall that two \emph{strongly orthogonal} roots in a root system are two roots whose sum and difference are not roots. We call two roots $\alpha$ and $\beta$ \emph{superstrongly orthogonal} if there is a choice of positive roots such that those two roots are simple and orthogonal. This is equivalent to the assertion that the only roots in the linear span of $\alpha$ and $\beta$ are $\pm \alpha, \pm \beta$.\footnote{Here is how to see the equivalence: If they are simple and orthogonal, the only root system they can generate is $A_1\times A_1$. Vice versa, if there are no more roots in their linear span, we can find real functionals $\ell_1$ and $\ell_2$, such that $\ell_1$ is positive on $\alpha, \beta$ and $\ell_2$ is zero on $\pm\alpha, \pm\beta$ and non-zero on all other roots. Then, for $s\gg 0$, the functional $\ell_1 + s \ell_2$ distinguishes a set of positive roots which must have $\alpha$ and $\beta$ as its simple elements, because it takes larger values on every other positive root. We thank Vladimir Drinfeld for pointing out this equivalence.} In \cite{BrO}, Brion proves the following:

\begin{theorem}[Brion]\label{Briontheorem}
 A set of generators for $W_{(X)}$ consists of elements $w$ which can be written as $w=w_1^{-1} w_2 w_1$ where:
\begin{itemize}
 \item $^{w_1}\mathring X =: Y$ with $\codim(Y)= l(w_1)$.
 \item $w_2$ is either of the following two:
\begin{enumerate}
 \item equal to the simple reflection $w_\alpha$  corresponding to a simple root $\alpha$ such that the $\PPGL_2$-spherical variety $\YY\PP_\alpha/\mathcal R(\PP_\alpha)$ is of the form $\TT\backslash \PPGL_2$ (where $\TT$ is a non-trivial torus) or $\mathcal N(\TT)\backslash \PPGL_2$;
 \item equal to $w_\alpha w_\beta$ where $\alpha,\beta$ are two orthogonal simple roots such that the $\PPGL_2$-spherical variety $\YY\PP_{\alpha\beta}/\mathcal R(\PP_{\alpha\beta})$ is of the form $\PPGL_2\backslash \PPGL_2\times \PPGL_2$.
\end{enumerate}
\end{itemize}
\end{theorem}

When we write, for instance, $\YY \PP_{\alpha} /\mathcal R(\PP_{\alpha})$, we mean
simply the homogeneous $\PP_{\alpha} /\mathcal R(\PP_{\alpha})$-variety
where a point stabilizer is given by the projection of a $\PP_{\alpha}$-stabilizer on $\YY$.

It is explained in \cite[\S 6.2]{SaSph}, based on an analysis of low-rank cases, that one can take the set of generators of the above theorem to be the reflections corresponding to (simple) spherical roots. However, we do not need and will not use this result.

\begin{corollary}
 Each spherical root $\gamma\in \Sigma_X$ is proportional to:
\begin{itemize}
 \item a sum of two superstrongly orthogonal roots of $G$, or
 \item a root of $G$.
\end{itemize}
\end{corollary}

\begin{proof}
By reduction to the varieties $\TT\backslash\PPGL_2$, $\mathcal N(\TT)\backslash\PPGL_2$, and $\PPGL_2\backslash\PPGL_2\times\PPGL_2$, together with the fact that $\varchi({^{w_1}\chi})=w_1\cdot \varchi(\XX)$, we deduce that the Weyl group elements described in both cases of the theorem induce hyperplane reflections on $\varchi(\XX)\otimes\mathbb Q$, and the vectors $w_1^{-1}\alpha$ (resp.\ $w_1^{-1}(\alpha+\beta)$) are inverted by those reflections, hence are proportional to the roots of the corresponding root system.

 Now consider the set $\mathcal{T}$ of all $\overline{W_X}$-conjugates of reflections thus obtained.
 The set of roots $\Sigma_{\mathcal{T}}$ associated to elements of $\mathcal{T}$ is a $\overline{W_X}$-stable subset
 of $\Sigma_X$ whose associated reflections generate $\overline{W_X}$.  In other words,
 $\Sigma_{\mathcal{T}}$ is a root subsystem with the same Weyl group.  This implies that $\Sigma_{\mathcal{T}} = \Sigma_X$, since root lines are characterized as the $-1$ eigenspaces of reflections. 
  \end{proof}

The two cases of the corollary are mutually exclusive, since the sum of two superstrongly orthogonal roots of $G$ cannot be proportional to a root. Hence, we can use them to define the \emph{type} of a root. Notice that, in the first case of the above theorem, if $\YY\PP_\alpha/\mathcal R(\PP_\alpha)$ is of the form  $\TT\backslash \PPGL_2$  then $w_1^{-1}\alpha\in \varchi(X)$, while if it is of the form $\mathcal N(\TT)\backslash \PPGL_2$ then $w_1^{-1}\alpha\notin \varchi(\XX)$ (while $2w_1^{-1}\alpha\in \varchi(\XX)$). 

\begin{definition}
 A spherical root $\gamma\in \Sigma_X$ is said to be:
\begin{itemize}
 \item of type $T$ if $\gamma$ is proportional to a root $\alpha$ of $\GG$ which belongs to $\varchi(\XX)$;
 \item of type $N$ if $\gamma$ is proportional to a root $\alpha$ of $\GG$ which does not belong to $\varchi(\XX)$;
 \item of type $G$ if $\gamma$ is proportional to the sum $\alpha+\beta$ of two strongly orthogonal roots of $\GG$.
\end{itemize}
 In this notation, the weight $\alpha$ (resp.\ the weight $\alpha+\beta$ for type $G$) will be called the \emph{normalized (simple) spherical root} corresponding to $\gamma$;
 it will sometimes be denoted by $\gamma'$.  The set of normalized (simple) spherical roots will be denoted by $\Delta_X$. 
\end{definition}

Notice that there is some issue here with the word ``simple'': while it should normally be used to distinguish elements of $\Delta_X$ from elements of the root system that they generate, it is customary in the theory of spherical varieties to call ``spherical roots'' only a set of simple roots. Therefore, we adopt this convention and feel free to drop the word ``simple'' when we talk about the set $\Delta_X$ or the set $\Sigma_X$. 

{\bluetext To say that $\gamma\in\Sigma_X$ is of type $N$ is equivalent to saying that $\gamma=2\alpha$, where $\alpha$ is a root of $\GG$. This follows, for example, from the classification of rank one wonderful varieties \cite{Ak}, but can also be deduced in a classification-free way from Theorem \ref{Briontheorem}.}

In the case of a normalized spherical root of type $G$, there is a canonical way to choose the roots $\alpha$ and $\beta$, which will be a useful fact later.

\begin{lemma}
 Consider a generator of $W_{(X)}$ as described in (and with the notation of) Theorem \ref{Briontheorem}. In cases $T$ and $N$ the root $w_1^{-1}\alpha$ is orthogonal to all the roots of $\LL(\XX)$, while in case $G$ the element $w_1^{-1} w_\alpha w_\beta w_1$ permutes the positive roots of $\LL(\XX)$. As a corollary, the element $w_1^{-1}w_\alpha w_1$, resp.\ $w_1^{-1} w_\alpha w_\beta w_1$, belongs to the canonical lift $W_X$ of $\overline{W_X}$.
\end{lemma}

\begin{proof}
 For types $T$ and $N$, we have $2 w_1^{-1}\alpha\in\varchi(X) \Rightarrow w_1^{-1}\alpha\perp \Delta_{\check L(X)}$. 
 
For type $G$, similarly, $w_1^{-1} (\alpha+\beta) \perp \Delta_{L(X)}$ implies that $\left< \alpha, w_1 \check\delta\right> + \left< \beta, w_1 \check\delta \right> =0 $ for every $\check\delta\in \Phi_{\check L(X)}$.

Let $\check\delta$ be such, then we claim that $w_1^{-1} w_\alpha w_\beta w_1 \check\delta \in \Phi_{\check L(X)}$ as well. By non-degeneracy (Lemma \ref{nondeg}), it suffices to show that it is perpendicular to $\varchi(X)$. Let $\chi\in\varchi(X)$. Then:
$$ \left< \chi, w_1^{-1} w_\alpha w_\beta w_1 \check\delta \right> = \left<\chi, \check\delta\right> - \left< \alpha, w_1\check\delta\right> \left< \chi , w_1^{-1}\check\alpha\right> - \left< \beta, w_1\check\delta\right> \left< \chi, w_1^{-1} \check\beta\right>.$$
From the description of ``type $G$'' in Theorem \ref{Briontheorem} we see that $\left<\chi, w_1^{-1}\check\alpha\right> = \left<\chi, w_1^{-1}\check\beta\right>$, which together with the above relation imply that the last expression is equal to zero.

Finally, we claim that if $\check\delta>0$ then $w_1^{-1} w_\alpha w_\beta w_1 \check\delta>0$. Indeed, since $\operatorname{dist}(\mathring X,Y)=l(w_1)$ and $W_{P(X)}$ stabilizes $\mathring X$ it follows that $w_1 \in [W/W_{P(X)}]$, and hence $w_1 \Phi_{\check L(X)}^+ \subset \Phi_{\check G}^+ \smallsetminus\{\check\alpha,\check\beta\}$. Hence, $w_\alpha w_1 \Phi_{\check L(X)}^+$ and $w_\beta w_1 \Phi_{\check L(X)}^+$ are contained in $\Phi_{\check G}^+$. We have also proved that $w_\alpha w_1 \Phi_{\check L(X)} = w_\beta w_1 \Phi_{\check L(X)}$, hence $w_\alpha w_1 \Phi_{\check L(X)}^+=w_\beta w_1 \Phi_{\check L(X)}^+$ and therefore $w_1^{-1} w_\alpha w_\beta w_1 \Phi_{\check L(X)}^+ = \Phi_{\check L(X)}^+$.

For the final conclusion, we remind that the canonical lift of the coset space $W/W_{L(X)}$ to $W$ consists of those elements which preserve the set of positive roots of $\LL(\XX)$.
\end{proof}

\begin{corollary} \label{onlytwo}
 For a spherical root $\gamma$ of type $G$, there are precisely two positive roots of $G$ in the $(-1)$-eigenspace of $w_\gamma$. They are both orthogonal to the weight $\rho_{L(X)}$.
\end{corollary}

We will call those the \emph{associated roots} of $\gamma$. For $\gamma$ of type $T$ or $N$ the \emph{associated root} will be the unique positive root of $\GG$ which is proportional to $\gamma$. The second statement of the lemma holds also for the associated root of a spherical root of type $T$ or $N$ for obvious reasons, namely that $\chi(\XX)$ is orthogonal to all roots of $\LL(\XX)$.

\begin{proof}
 The statement does not depend on whether $\gamma$ is simple or not, so it is enough to show it for generators $w_\gamma$ of the form $w_1^{-1} w_\alpha w_\beta w_1$ as in the previous lemma. Notice that if $\alpha'$ is an associated root then $w_\gamma \alpha' =-\alpha'$, while on the other hand, by the previous lemma, $w_\gamma \rho_{L(X)} = \rho_{L(X)}$. Hence, 
$$ \left<\rho_{L(X)},\alpha'\right> = \left<w_\gamma\rho_{L(X)},w_\gamma\alpha'\right> = -\left< \rho_{L(X)}, \alpha'\right> \Rightarrow \left< \rho_{L(X)}, \alpha'\right> =0.$$
\end{proof}

Now we return to defining the normalized root system:

\begin{proposition}
 Under the action of $W_X$ on $\varchi(\XX)$, the set $\Delta_X$ is the set of simple roots of a root system with Weyl group $W_X$.
\end{proposition}

\begin{proof}
 Since the same statement is true for $\Sigma_X$, it suffices to prove that if $\gamma_1,\gamma_2\in\Sigma_X$ and $w\in W_X$ are such that $w\gamma_1=\gamma_2$ then for the corresponding normalized spherical roots $\gamma_1', \gamma_2'$ we still have: $w\gamma_1'=\gamma_2'$. But this is obvious from the definitions.
\end{proof}

We denote by $\Phi_X$ the root system generated by $\Delta_X$.

Now we come to the root data which, conjecturally, correspond to the dual group of Gaitsgory-Nadler $\check G_{X,GN}\subset \check G$ (to be recalled in the next subsection) and its central isogeny $\check G_X\twoheadrightarrow \check G_{X,GN}$ (whenever it can be defined). Notice that in the first case, in order for a root datum to define a subgroup of $\check G$, its lattice should be a sublattice of $\varchi(\AA)$ without co-torsion.

\begin{proposition}
 The set $(\varchi(\AA)\cap \mathbb Q\cdot \varchi(\XX), \Phi_X, W_X)$ gives rise to\footnote{Usually a root datum is described in terms of a pair $L, \check{L}$ of finite free $\Z$-modules, together with subsets $\Phi \subset L, \check{\Phi} \subset \check{L}$ of roots and coroots.  However, this is determined up to isomorphism by the triple $(L, \Phi, W_L)$, where $W_L$
 is the Weyl group.} a root datum. If there are no spherical roots of type $N$ then the set $(\varchi(\XX),\Phi_X,W_X)$ also gives rise to a root datum.
\end{proposition}

\begin{proof}
By definition, the elements of $\Phi_X$ belong to $\varchi(\AA)\cap \mathbb Q\cdot \varchi(\XX)$, and if there are no spherical roots of type $N$ then they also belong to $\varchi(\XX)$ (as we deduce, again, by reduction to the varieties $\TT\backslash\PPGL_2$ and $\PPGL_2\backslash\PPGL_2\times\PPGL_2$ together with the fact that $\varchi({^{w_1}\chi})=w_1\cdot \varchi(\XX)$).   Therefore, it remains to check that the corresponding coroots are integral on the given lattices. 

Let $\gamma\in\Phi_X$ correspond to a hyperplane reflection of the form described in the  Theorem \ref{Briontheorem}. The associated coroot equals: 
\begin{itemize}
 \item the image of $w_1^{-1}\check\alpha$ in $\varchi(\XX)^*$, if $\gamma=w_1^{-1}\alpha$ is of type $T$ or $N$ (in the notation of the theorem);
 \item the image of $w_1^{-1}\check\alpha$ (which coincides with the image of $w_1^{-1}\check\beta$), if $\gamma=w_1^{-1}(\alpha+\beta)$ is of type $G$.
\end{itemize}
Those are integral on the given lattices, which completes the proof of the proposition.
\end{proof}

\subsection{Distinguished morphisms}  \label{distinguishedmorphisms}

We introduce the following notation:

\begin{itemize}
\item[-]
$\check G_X'$ the abstract complex reductive group defined by the root datum of $(\varchi(\AA)\cap \mathbb Q\cdot \varchi(\XX), \Phi_X, W_X)$.
\item[-] If there are no spherical roots of type $N$,  $\check G_X$ is the abstract complex reductive group defined by the root datum of  $(\varchi(\XX), \Phi_X, W_X)$.
\end{itemize}
These groups come with preferred maximal tori $A_X^*,A_X'^*$ and are unique up to the inner action of this torus. Moreover, since the root data used to define them are actually based (i.e.\ have a preferred choice of positive roots), the groups $\check G_X, \check G_X'$ also have a preferred choice of Borel subgroup containing the canonical maximal tori.
The obvious isogeny between their root data gives rise to a canonical central isogeny $\check G_X\twoheadrightarrow \check G_X'$. We conjecture that the group $\check G_X'$ is isomorphic to the one constructed by Gaitsgory and Nadler, which we denote by $\check G_{X,GN}$.  { A priori, the group $\check G_{X,GN}$ depends on the choice of an affine embedding of $X$; conditional on some assumptions on the Gaitsgory-Nadler dual group which will be discussed in \S \ref{GNaxioms}, we will show (Corollary
\ref{dualgroup:ID})  that $\check G_{X,GN}$ is indeed equal to $\check G_X'$ and hence independent of the affine embedding. More precisely, we will show that $\check G_{X, GN}$ is obtained from a \emph{distinguished} embedding of $\check G_X'$ into $\check G$.}

 Call a morphism $\check G_X \rightarrow \check G$ {\em distinguished}
 if: 
\begin{enumerate}\item it extends the canonical map $A_X^* \rightarrow A^*$;
  \item for every simple (normalized) spherical root $\gamma$, the corresponding root space of $\check{\mathfrak g}_X$ maps into the root spaces of its associated roots.
\end{enumerate}

By Lemma \ref{onlytwo} (and the comment which follows it), the second condition implies that the image of a distinguished morphism commutes with the image of ${2 {\rho_{L(X)}}}$.  
 
We will call a morphism $\check G_X \times \SL_2 \rightarrow \check G$ \emph{distinguished}
 if its restriction to $\check G_X$ is distinguished, and its restriction to $\SL_2$ is a { \emph{principal} morphism into $\check L(X)$: $$\SL_2\to \check L(X)\subset \check G$$ 
 with weight: 
 $$\Gm \stackrel{{2  \rho_{L(X)}}}{\longrightarrow} \check{G},$$
 where $\Gm$ is identified as a subgroup of $\SL_2$ in the standard way: $a\mapsto \left(\begin{array}{cc} a \\ & a^{-1} \end{array}\right)$.} 
 
We apply similar terminology in related situations: for instance, 
 there is a corresponding notion of distinguished map: $\check G_X' \rightarrow \check{G}$, or a distinguished map when we deal with a standard Levi subgroup of $\check G_X$ containing $A_X^*$.

		\subsection{The work of Gaitsgory and Nadler} \label{GNaxioms}

In this subsection { we fix an affine embedding $\XX^a$ of $\XX$}.

Let us denote by $A_{X,GN}^*$ the image of the dual torus $A_X^*$ inside $A^*$. 
Recall as before that we regard    the sum $2\rho_{L(X)}$ of positive roots of $\LL(\XX)$ as
defining a character ${2 \rho_{L(X)}}: \Gm \rightarrow A^*$. 

\begin{theorem}[Gaitsgory and Nadler] \label{gnconj}
 To every affine spherical variety $\XX^a$ one can associate a connected reductive subgroup $\check G_{X^a,GN}$ of $\check G$ with maximal torus $A_{X,GN}^*$. The group $\check G_{X^a,GN}$ is canonical up to $A^*$-conjugacy. 
\end{theorem}

This is not very informative as stated, but the assertions (GN1)--(GN5) of \S \ref{dualgroup}, which we recall here with a few extra comments for the convenience of the reader, give more information about the group $\check G_{X,GN}$. To formulate them, let $\XX_\Theta^a$ be, for every $\Theta\subset\Sigma_X$, the affine embedding of $\XX_\Theta$ obtained by partially grading the  coordinate ring of $\XX^a$. In terms of the affine degeneration $\mathscr X^a \to \overline{\AA_{X,ss}}$
discussed in \S\ref{ssdegen}, it is the fiber over $\lim_{t\to 0}\check\lambda(t)$, where $\check\lambda$ is any cocharacter in the cone of $\overline{\AA_{X,ss}}$ which lies in the interior of the face determined by $\Theta$. Since we are no experts in the technical details of \cite{GN}, we will only prove the first of the following assertions, and treat the remaining as hypotheses: \label{Gnaxiomspage} 

\begin{enumerate} 
\item[(GN1)] The image of $\check G_{X^a,GN}$ commutes with ${2\rho_{L(X)}}(\CC^\times)\subset A^*$.  
\item[(GN2)] The Weyl group of $\check G_{X^a,GN}$ equals $W_X$. (This is
a consequence of \cite[Conjecture 7.3.2]{GN}, as discussed there.) 
\item[(GN3)] For any $\Theta\subset \Sigma_X$ the dual group of $\XX^a_{\Theta}$ is canonically a subgroup of $\check G_{X^a,GN}$. (Our identification of $\check G_{X^a,GN}$, based on these axioms, shows that it is the Levi subgroup of $\check G_{X^a,GN}$ associated to $\Theta$.)
\item[(GN4)] If the open $\GG$-orbit $\XX\subset \XX^a$ is parabolically induced, $\XX = \XX_L \times^{\PP^-} \GG$, where $\XX_L$ is spherical for the reductive quotient $\LL$ of $\PP^-$, then the dual group $\check{G}_{X^a,GN}$
belongs to the standard Levi subgroup $\check L$ of $\check G$ corresponding to the class of parabolic subgroups opposite to $\PP^-$.  Moreover, if a connected normal subgroup $\LL_1$ of $\LL$ acts trivially on $\XX_L$, then $\check{G}_{X^a,GN}$ belongs to the dual group of $\LL/\LL_1$ (which is canonically a subgroup of $\check L$).
\item[(GN5)] If $\XX_1^+$ is a spherical homogeneous $\GG$-variety, $\TT$ a torus of $\GG$-automorphisms and $\XX_2^+=\XX_1^+/\TT$, and if $\XX_1, \XX_2$ are affine embeddings of $\XX_1^+,\XX_2^+$ with $\XX_2=\spec k[\XX_1]^\TT$,  then there is a canonical inclusion $\check{G}_{X_2, GN} \hookrightarrow \check{G}_{X_1, GN}$ which restricts to the natural inclusion of maximal tori: $A_{X_2, GN}^* \hookrightarrow A_{X_1, GN}^*$ (arising from $\varchi(\XX_2)\hookrightarrow\varchi(\XX_1)$).
\end{enumerate}

We shall give a proof of (GN1) or rather ``angle'' it out of the articles of Gaitsgory and Nadler.
In \cite{GN} a certain tensor category $\mathbf Q(X)$ (denoted $\mathbf Q(Z)$ in \emph{loc.cit.}) is constructed, together with adequate functors: 
\begin{equation} \label{GN-functor} \Rep(\check{G}) \stackrel{\mathrm{Conv}}{\rightarrow} \mathbf Q(X)  \stackrel{\mathrm{fib}}{\rightarrow} \mathrm{Vect}.\end{equation} 
The first category, as is usual, is constructed as $\GG(\C[[[t]])$ equivariant sheaves on the affine Grassmannian $\GG(\C((t)))/\GG(\C[[t]])$, and $\mathrm{Vect}$ denotes the category of vector spaces. The category $\mathbf Q(X)$ is constructed via a certain substitute for $\GG(\C[[t]])$-equivariant sheaves of $X(\mathbb{C}((t)))$.

\begin{proof}[Proof of (GN1)]

To show that the image of ${2\rho_{L(X)}}$ commutes with $\check G_X$, it suffices to show that there is a $\mathbb Z$-grading of the tensor category $\mathbf Q(X)$ such that under the ``convolution'' functor $\Rep(\check{G}) \stackrel{\mathrm{Conv}}{\rightarrow} \mathbf Q(X)$ and the equivalence of $\mathbf Q(X)$ with $\Rep(\check G_{X,GN})$ the grading corresponds to the decomposition of representations in $\CC^\times$-eigenspaces, where $\CC^\times$ acts via ${2\rho_{L(X)}}$. This grading is explicit, in the form of a cohomological shift, in \cite{GNhoro}, but implicit in \cite{GN}. More precisely, it is shown in \cite{GNhoro}[Theorem 1.2.1], where the special case of horospherical varieties is studied, that for a horospherical variety $\XX_0$ the irreducible objects of $\mathbf Q(X_0)$ can be identified with intersection cohomology sheaves of certain strata, shifted in cohomological degree, and the shift is precisely the grading that we want. In \cite{GN} the authors choose to forget about the 
cohomological shift, however, this shift has to be compatible with the fiber functor $\mathbf Q(X)\to\mathbf Q(X_0)$ (where $\XX_0$ is the boundary degeneration that we denoted before by $\XX_\Theta$ for $\Theta=\emptyset$) because the fiber functor is obtained via a nearby cycles functor; hence, the category $\mathbf Q(X)$ carries the grading corresponding to ${2\rho_{L(X)}}$, as well.
\end{proof}

\begin{remark} Concerning (GN2) (GN3), (GN4) and (GN5): The statement of (GN2) is easy to deduce from the results of Gaitsgory and Nadler in the most interesting cases. First of all, we claim that $\mathcal N_{\check G_{X, GN}}(A_X^*)$ is contained in $\mathcal N_{\check G}(A^*)$. Indeed, $\mathcal N_{\check G_{X, GN}}(A_X^*)$ centralizes the image of ${2\rho_{L(X)}}$ (because it belongs to $\check G_{X, GN}$) and normalizes the centralizer of $A_X^*$ inside of the centralizer of ${2\rho_{L(X)}}$. By non-degeneracy (Lemma \ref{nondeg}), the common centralizer of $A_X^*$ and ${2\rho_{L(X)}}$ is $A^*$, hence $\mathcal N_{\check G_{X, GN}}(A_X^*)\subset \mathcal N_{\check G}(A^*)$. Now, the combination of Theorem 4.2.1 and Proposition 5.4.1 of \cite{GN} imply that the restriction to $A_X^*$ of any irreducible representation of $\check G_{X, GN}$ contains a character in $\Lambda_X^+$. In the cases where $W_X$ coincides with the normalizer of $A_X^*$ in $\mathcal Z_{\check G}({2\rho_{L(X)}})$ (such as for symmetric spaces), it follows 
immediately that the Weyl group of the dual group of Gaitsgory and Nadler has to be the whole $W_X$, for otherwise any Weyl chamber of it would be larger than $\Lambda_X^+$. The requirement that $\mathcal N_{\mathcal Z_{\check G}({2\rho_{L(X)}})}(A_X^*)/A_X^* = W_X$ can be understood   representation theoretically as follows: it was proven in \cite{SaSpc} that the multiplicity of a generic unramified representation in the spectrum of $X$ is equal to the product of the ``geometric factor'' $(N_{\mathcal Z_{\check G}({2\rho_{L(X)}})}(A_X^*)/A_X^* : W_X)$ by the ``arithmetic factor'' of the number of open $\BB(k)$-orbits on $\XX$. Thus, $\mathcal N_{\mathcal Z_{\check G}({2\rho_{L(X)}})}(A_X^*)/A_X^* = W_X$ means that the geometric factor of unramified multiplicities is 1.

(GN3) should follow in the same way as the ``fiber functor'' construction in \cite{GN}, except that the fiber functor was constructed through a full degeneration of the spherical variety (i.e.\ a degeneration to $\XX_\emptyset$), while for (GN3) one would only perform a partial degeneration. (GN4) should also be feasible along the lines of \cite{GNhoro}, by interpreting geometrically the action of the center of $\LL$ ``on the left''. Finally, (GN5) should follow from the behavior of intersection cohomology under such quotients by toric actions. However, since we are not specialists in the subject we treat (GN2)--(GN5) as hypotheses.
\end{remark}

\subsection{Uniqueness of a distinguished morphism}

 In what follows we denote by $A_{X,GN}^*$ the (canonical) maximal torus of the Gaitsgory-Nadler dual group, i.e.\ the image of $A_X^*$ in $A^*$.  We fix throughout a standard basis $\{h,e,f\}$ for the Lie algebra $\mathfrak{sl}_2$. 
By the `weight'' of a morphism: $f: \SL_2\to \check G$ we understand either its restriction to $\Gm$, or the derivative of this: $\left<h\right>\simeq \mathfrak g_m \to \check{\mathfrak g}$. We will repeatedly use the following fact: if $f, f'$ have the same weight, they are conjugate
by an element of the centralizer of this weight.

 \begin{lemma} \label{rk1uniqueexists}
 Let $\gamma\in \Delta_X$, and let $\check G_\gamma$ be the corresponding subgroup of $\hat G_X'$, i.e.\ a connected reductive group with a canonical maximal torus $A_{X,GN}^*$ and unique simple coroot $\gamma: \Gm \rightarrow A_{X, GN}^*$. 
 
 \begin{enumerate}
 \item A distinguished morphism $\psi:\check G_\gamma\to \check G$, through which the root space of $\check \gamma$ maps into the sum of root spaces of the associated coroots, always exists;
 \item Any two such are conjugate by $A^*$; 
\item
If $\psi$ is such,  its centralizer in $A^*$ contains the common kernel of the associated roots to $\gamma$. 
\end{enumerate}
 \end{lemma}

\begin{proof}
Let $\AA_{\gamma}$ be the identity component of the kernel of $\gamma$. 
Then $\check G_{\gamma}$ is the almost-direct product
$A_{\gamma}^* \cdot f_{\gamma}(\SL_2)$, where $f_{\gamma} :\SL_2 \rightarrow \check{G}_X'$
has weight $\gamma$.  

Let $\check{M}$ be the centralizer of $A_{\gamma}^* \cdot \mathrm{image}({2  \rho_{L(X)}})$ inside $\check{G}$.   It is a reductive group.

If $\psi$ is distinguished, the map $ \psi \circ f_{\gamma} $ has image in $\check{M}$ and has weight $\gamma$; 
the association $\psi \mapsto \psi \circ f_{\gamma}$ gives a bijection between
distinguished morphisms  and the set of:  \begin{equation}\label{wet} m: \SL_2 \rightarrow \check{M},  \ \ m|\Gm = \gamma. \end{equation}

  If $\gamma$ is a root, the existence of $m$ as in \eqref{wet} is clear. 
 Otherwise,
it is the sum $\alpha + \beta$ of two superstrongly orthogonal positive roots;
choosing a positive system in which $\alpha, \beta$ are simple, we see there
are associated morphisms $f_{\alpha}, f_{\beta} : \SL_2 \rightarrow  \check{G}$ (corresponding to $\alpha, \beta$ thought of as coroots of $\check{G}$). These morphisms have commuting image,
since $\alpha, \beta$ are strongly orthogonal; therefore, we obtain a product morphism
$$f_{\alpha} \times f_{\beta}: \SL_2 \times \SL_2 \longrightarrow \check{G},$$
whose diagonal is a morphism $m$ as in (\ref{wet}). This proves the first claim.

To check that any two such morphisms are $A^*$-conjugate, it suffices
to check that any two $m$ as in \eqref{wet} are $A^*$-conjugate.  However,
any two such $m$ are conjugate by the centralizer of $m|\Gm$, 
i.e., by the centralizer of $\gamma: \Gm \rightarrow \check{M}$. 
That is the same as the centralizer of $A_{X, GN}^* \times \mathrm{image}({2 \check \rho_{L(X)}})$. 
Since the spherical variety $\XX$ is quasi-affine and, hence, non-degenerate, this centralizer is equal to $A^*$.   

The final assertion follows from the explicit construction of a distinguished morphism. \end{proof}

\begin{lemma} \label{Associatedlinindep}
 The associated roots to all $\gamma\in \Delta_X$ are linearly independent.
\end{lemma}

\begin{proof}[Proof (sketch)] 
This is a clumsy argument reducing to the low-rank cases: Consider a linear relation between associated roots with non-zero coefficients, let $R$ denote the support of all roots appearing and let $\alpha\in R$ be a simple root which is connected, in the Dynkin diagram, to at most one more element of $R$, and is not the shorter of the two.  (Let's call such a simple root ``extreme'' for the given collection of associated roots.) Necessarily, the root $\alpha$ has to be contained in the support of at least two associated roots in the linear relation, say $\gamma$ and $\delta$. These, in turn, should be associated to spherical roots $\varepsilon$ and $\zeta$ (not necessarily distinct). By inspection of spherical varieties of rank one or two \cite{Wa}, we see that an extreme simple root cannot be in the support of two associated roots.  \end{proof} 

\begin{proposition} \label{prop-uni}
Distinguished morphisms $\check G_X' \rightarrow \check{G}$
and $\check G_X' \times \SL_2 \rightarrow \check{G}$, if they exist,
are unique up to $A^*$-conjugacy. 
\end{proposition}

\proof
Let $\psi_1, \psi_2: \check G_{X}' \rightarrow \check{G}$ be distinguished. 
We have seen that, for every spherical root $\gamma$,
there exists $a_{\gamma} \in A^*$ so that
$$\Ad(a_{\gamma}) \psi_1|)_{\check G_{X_{\gamma}}'} = \psi_2 |_{\check G_{X_{\gamma}}'}.$$
The action of $A^*$ on the image of $\check G_{X_{\gamma}}'$ factors through
the morphism $A^* \rightarrow \Gm$ or $A^* \rightarrow \Gm^2$ induced by
the associated roots for $\gamma$, by part (3) of Lemma \ref{rk1uniqueexists}. Thus, by Lemma \ref{Associatedlinindep}, we may find
$a \in A^*$ so that $\Ad(a) $ and $\Ad(a_{\gamma})$ have the same action on 
$ \psi_1|_{\check G_{X_{\gamma}}'} $ for all $\gamma$.  In particular, there exists $a$ so that
$$\Ad(a) \psi_1|_{\check G_{X_{\gamma}}'} = \psi_2 |_{ \check G_{X_{\gamma}}'}$$
for all $\gamma$. Since the $\check G_{X_{\gamma}}'$ generate $\check G_X$,
it follows that $\psi_1, \psi_2$ are $A^*$-conjugate, as desired. 

This completes the proof of the assertion for $\check G_X'$.

The assertion for $\check G_X' \times \SL_2$ follows from the first: Fix a distinguished
embedding $\psi$ of $\check{G}_X'$. Now  ${2 \rho_{L(X)}}$
defines a morphism from $\Gm$ to the connected centralizer of $\psi(\check{G}_X')$;
any two $\SL_2$-morphisms with this restriction to the diagonal $\Gm$
must be conjugate under the connected centralizer of $\psi(\check{G}_X') \times {2 \rho_{L(X)}}$. 
The latter is a subgroup of $A^*$ (since $X$ is nondegenerate, Lemma \ref{nondeg}), commuting with $\psi(\check{G}_X')$.

\qed 

		\subsection{The identification of the dual group} 

In this subsection we will use the axioms (GN) in order to identify the (based) root datum of the Gaitsgory-Nadler dual group with the (based) root datum of the abstract group which we denoted by $\check G_X'$, and in fact to identify $\check G_{X,GN}$ as a subgroup of $\check G$ uniquely up to $A^*$-conjugacy. (Notice that in the classical setting, as opposed to the geometric one, the group $\check G$ itself is only canonical up to $A^*$-conjugacy.) Using Proposition \ref{prop-uni}, the only thing that we need to prove is that simple roots of $\check G_{X,GN}$ and $\check G_X'$ have the same length and, in fact, for every simple root $\gamma$ the embedding of the corresponding standard Levi $\check G_\gamma\subset\check G_{X,GN}\hookrightarrow \check G$ is distinguished. The argument will eventually boil down to the classification of rank-one wonderful varieties by Akhiezer \cite{Ak}. 

\begin{proposition} \label{prrkone}  
 Assume axioms (GN2)--(GN5). For spherical varieties of rank one there
 is an isomorphism of Gaitsgory-Nadler dual group with $\check G_X'$, inducing
 the identity on $A_{X,GN}^*$, and for any such isomorphism the embedding $\check G_X'\xrightarrow{\sim}\check G_{X,GN}\hookrightarrow \check G$ is distinguished. 
\end{proposition}

\begin{proof}

A \emph{homogeneous} spherical $\GG$-variety of rank one is of the form:
\begin{equation}
 \XX=\XX_1 \times^{\PP}\GG,
\end{equation}
where:
\begin{enumerate}
 \item $\XX_1$ is a $\GG_1$-torus bundle\footnote{By $\GG_1$-torus bundle we mean a principal torus bundle with an action of $\GG_1$ commuting with the action of the torus.}  over a variety $\HH\backslash\GG_1$ from Table 1 of \cite{Wa} (we denote by $\GG_1$ the group $G$ of \emph{loc.cit.});
 \item $\PP$ is a parabolic subgroup with a homomorphism: $\PP\to\Aut(\XX_1)$ whose image coincides, up to central subgroups, with the image of $\GG_1$.
\end{enumerate}

By inspection of this table, and using the axioms (GN1)--(GN5), one can show that the Gaitsgory-Nadler dual group is unambiguously equal to $\check G_{X_\gamma}'$. The table of Wasserman, together with more details of this argument, are given in Appendix \ref{primerankone}.  \end{proof}

\begin{corollary} \label{dualgroup:ID}
Assume axioms (GN2)--(GN5).
Then there exists a distinguished embedding $\check G_X' \hookrightarrow \check{G}$
with image $\check G_{X,GN}$. 
In particular,  the group $\check{G}_{X,GN}$ is canonically isomorphic to $\check G_X'$ up to $A_{X,GN}^*$-conjugacy.
\end{corollary}

\proof 
By (GN2), the coroots of $\check{G}_{X,GN}$ are proportional to elements of the set $W_X\cdot \Delta_X$.
(In fact, the lines through coroots are characterized as the $-1$ eigenspaces
of reflections in $W_X$.) Therefore, there exists a system of simple positive coroots for $\check{G}_{X,GN}$, each of which is proportional to one of the $\gamma\in \Delta_X$. 

On the other hand, by (GN3), the group 
$\check G_{X_{\gamma}, GN}$ is contained in $\check G_{X, GN}$; the former group,
as a subgroup of $\check G$, is identified through Lemma \ref{rk1uniqueexists} and Proposition \ref{prrkone}. It follows that we may suppose that $\check G_{X,GN}$
contains the image of a distinguished homomorphism
$$\check G_{\gamma}' \longrightarrow \check{G},$$
as in Lemma \ref{rk1uniqueexists}. Hence, the coroot of $\check G_{X,GN}$ proportional to $\gamma\in \Delta_X$ actually {\em equals} $\gamma$.

It follows that the coroots of $\check{G}_{X,GN}$ are precisely the elements of the set $W_X\cdot \Delta_X$.  Hence, the root data of $\check{G}_{X,GN}$ and $\check{G}_X'$ coincide \emph{canonically} (recall that from the work of Gaitsgory and Nadler the group $\check{G}_{X,GN}$ is canonical up to $A^*$-conjugacy) and so the two are canonically isomorphic up to $A_{X,GN}^*$-conjugacy.
 \qed 

		\subsection{Commuting $SL_2$}

In this subsection we will prove that there is a principal $\SL_2$ inside of $\check L(X)$ commuting with $\check G_{X,GN}$, assuming (GN3) and (GN4). Our proof will be quite clumsy, using combinatorial arguments to reduce the problem to the case of spherical varieties of small rank, where it is checked case-by-case.

The basic result, which will be established case-by-case in Appendix \ref{primerankone}, is the following:

\begin{proposition}\label{rankonesl2commute}
 Let $\XX$ be a spherical variety of rank one and assume (GN4). Then there is a principal map: $\SL_2\to \check L(X)$ which commutes with $\check G_{X,GN}$.
\end{proposition}

Using this, we can now show:

\begin{theorem}
 Assume that there is a distinguished embedding: $\check G_X' \to \check G$ (as we have proven under the (GN) assumptions in Corollary \ref{dualgroup:ID}). Then there is a principal $\SL_2\to \check L(X)$ whose image commutes with $\check G_X'$.
\end{theorem}

Recall also that, by Proposition \ref{prop-uni}, the resulting distinguished morphism: $\check G_X'\times \SL_2\to \check G$ is unique up to $A^*$-conjugacy.

\begin{proof}
 Fix a principal $\SL_2$ into $\check L(X)$ with weight ${2 \rho_{L(X)}}$, and denote its image by $S$; all such subgroups are $A^*$-conjugate. Fix a distinguished embedding of $\check G_X'$ into $\check G$. By Proposition \ref{rankonesl2commute}, for every $\gamma\in\check\Delta_X$ there is an $A^*$-conjugate of $S$ which commutes with $\check G_\gamma$. Equivalently, there is an $A^*$-conjugate of $\check G_\gamma$ which commutes with $S$. 
 Arguing as in Proposition \ref{prop-uni},  we may find $a \in A^*$ which conjugates all $\check G_\gamma$ simultaneously into the centralizer of $S$.  
\end{proof}

In the case that there are no spherical roots of type $N$ (equivalently, as we mentioned, no element of $\Sigma_X$ is of the form $2\alpha$, for $\alpha$ some root of $\GG$), composing this with the central isogeny: $\check G_X\to \check G_X'$ we get the desired distinguished morphism:
\begin{equation}
 \check G_X\times \SL_2 \to \check G.
\end{equation}
The proof of Theorem \ref{sl2} is now complete.

 \part{Local theory and the Ichino--Ikeda conjecture}

\section{Geometry over a local field} \label{localfieldgeom}
In this section we shall examine certain general features of the geometry of $X = \XX(k)$, where $k$ is a $p$-adic field. In particular, we shall establish the relationship between $G$-invariant
measures (or $G$-eigenmeasures) on $X$ and $X_{\Theta}$; this will lead us to fixing compatible measures on $X$ and $X_\Theta$ for the rest of the paper, as we indicated in \S \ref{ssassumptions}.

More importantly, we shall
establish the ``exponential map'' which relates the structure of $X$ and $X_{\Theta}$ near infinity.

\subsection{Measures} \label{ssmeasures} We may assume, without serious loss of generality, that $X$ carries a positive $G$-eigenmeasure $\mu$. Indeed, this is the case if the modular character of $H$ (the quotient of its right by its left Haar measures) extends to a character of $G$. For a given $\XX=\HH\backslash\GG$, the algebraic modular character $\mathfrak d_H$ of $\HH$ is defined over $k$, and it either has finite image or surjects onto $\GGm$. In the latter case, we may replace $\HH$ by the kernel $\HH_0$ of $\mathfrak d_H$; then there is an $\HH/\HH_0\times\GG$-eigen-volume form on $\HH_0\backslash\GG$, and its absolute value gives the invariant measure. (Notice also that $\HH_0\backslash \GG\to \HH\backslash\GG$ is surjective on $k$-points.) In the former case, there is an invariant volume form $\omega$ valued in the bundle over $\HH\backslash\GG$ defined by $\mathfrak d_H$. Since $\mathfrak d_H$ has image in the $n$-th roots of unity $\mathbb \mu_n$, for some $n$, the associated complex line 
bundle is trivial and therefore the absolute value of $\omega$ defines an invariant measure on $\HH\backslash \GG(k)$. (For economy of language, we will never again mention the possibility that our eigenforms are valued in a torsion line bundle, instead of the trivial one.)

We fix from now on such an eigenmeasure $\mu$, i.e.\ an eigenmeasure which is the absolute value of a volume form.  We define $L^2(X):=L^2(X,\mu)$, considered as a \emph{unitary} representation of $G$ by twisting the right regular representation by the square root of the eigencharacter of $\mu$, that is: 
\begin{equation} \label{normalizedaction} (g\cdot \Phi)(x)=\sqrt{\eta(g)}\Phi(xg),
\end{equation}
 where $\eta$ denotes the eigencharacter of $\mu$.

\begin{example}
 When $\XX=\UU\backslash \GG$, with $\UU$ a maximal unipotent subgroup, the action of $A\times G$ on $L^2(X)$ is defined as: $((a,g)\cdot f) (x)=\delta^{-\frac{1}{2}}(a) f(a\cdot x\cdot g)$, where $\delta$ is the modular character of the Borel subgroup (the quotient of a right by a left invariant Haar measure).
\end{example}

\subsection{The measure on $X_\Theta$} 
Our concern in this section is to relate the measures on $X$ and $X_{\Theta}$. 
 Specifically:

\begin{proposition}\label{propmeasure}
For every $\GG$-eigen-volume form $\omega$ on $\XX$ there is a canonical $\AA_{X,\Theta}\times \GG$-eigen-volume form $\omega_\Theta$ on $\XX_\Theta$, with the same $\GG$-eigencharacter, characterized by the property that for every Borel subgroup $\BB$ the restrictions of the forms $\omega$ and $\omega_\Theta$ to $\mathring\XX$, resp.\ $\mathring\XX_\Theta$ correspond to each other under the isomorphism (\ref{Borbitident}).
\end{proposition}

Again, of course, we will twist the action of $A_{X,\Theta}\times G$ on functions on $X_\Theta$ as in \eqref{normalizedaction}.

Indeed, if $\omega$ is a differential form of top degree on $\overline{\XX}$, then it induces
in a natural way a differential form of top degree on the normal bundle to any subvariety. 
Unfortunately, the $\GG$-eigen-volume forms are rarely regular at the boundary. 
Nonetheless (formalizing the intuition that normal bundles
model small neighborhoods of submanifolds) we may associate to any {\em rational} form 
on $\overline{\XX}$ a rational form on each degeneration $\XX_{\Theta}$.
We explain how to do this:

\subsubsection{Obtaining measures by degeneration/residues} \label{nbdegeneration}

 Let $\bar \XX$ be a smooth variety and let $\ZZ\subset\bar\XX$ be a closed subvariety, obtained as the intersection of a fixed set of reduced divisors $\DD_i: 1 \leq i \leq m$, with simple normal crossings and non-empty intersection. The choice of divisors $\DD_i$ makes the normal bundle of $\ZZ$ into a $\TT:=\GGm^m$-space:
 the associated grading is the decomposition of the normal bundle at a point $z \in \ZZ$ as the sum of  normal bundles at $z$
 to each $\DD_i$.

Now suppose $\omega $ is a nonzero  differential form of top degree on $\XX = \overline\XX \smallsetminus \cup_i \DD_i$. Let $-n_i - 1$ be the valuation of $\omega$ at $\DD_i$; let $f_i =0$ be a local equation for $\DD_i$. Then
we obtain a rational differential form on the normal bundle via
\begin{equation} \label{baromegadef} \bar{\omega} :=  \frac{ (\pi^* \mathrm{Res}  (\omega \cdot \prod f_i^{n_i}))  \wedge df_1' \wedge df_2' \wedge \dots \wedge df_m'} {\prod (f_j')^{n_j+1} }. \end{equation}
where $\mathrm{Res}$ denotes  the iterated residue,  $\pi$ is the projection from normal bundle to $\ZZ$, 
and $f_j'$ denotes the derivative of $f_j$, considered as a function on the normal bundle of
$\ZZ$.     It is possible that $\bar{\omega} \equiv 0$: consider the
case when $\dim \XX = 2$ and $\omega = (f_1^{-1} + f_2^{-1}) df_1 \wedge df_2$. 
 
 In fact, $\bar{\omega}$ is independent of choices.
Indeed, there  is a more intrinsic way of understanding  \eqref{baromegadef} via
  degeneration to the normal bundle:
$$\mathscr X\to \GGa^{m}$$
that we discussed in \S \ref{ssdegen} (where it was denoted by $\mathscr X^n$ to distinguish it from the affine degeneration).  
 
Let  $\chi_n :  (t_1,\dots,t_m)\mapsto \prod t_i^{n_i} $,  and define a differential form $\tilde \omega$ on $\GGm^n \times\XX \subset \mathscr{X}$ --
via \begin{equation}
 \tilde \omega := \chi_n(t) \cdot p^*\omega,
\end{equation}
where $p:\GGm^m \times\XX\to\XX$ is the natural projection. 
We regard $\tilde \omega$ as a rational differential form on $\mathscr{X}$. 
Then:

\begin{lemma}
 The restriction of the form $\tilde\omega$ to any fiber of the map: $\mathscr X\to \GGa^m$, is well defined as a rational differential form on that fiber. If $\omega$ is regular everywhere, then $\tilde\omega$ is also regular everywhere. Finally, $\tilde\omega$ is an eigenform for the action of $\GGm^m$ on $\mathscr X$, with eigencharacter $\chi_n^{-1}$.
\end{lemma}

\begin{proof}
 This is easy to see for codimension-one orbits of $\TT$ in $\mathcal B$, and then by Hartogs' principle it extends to $\mathcal X$. The fact that it is an eigenform for $\GGm^n$ is obvious from the definitions.
\end{proof}

The restriction to the fiber over zero coincides, up to sign, with the form $\bar{\omega}$ defined previously. 
 
\subsubsection{The case of spherical varieties}

We now specialize to the case $\XX=$ our spherical variety and $\bar\XX=$ its wonderful compactification. 
The set $\Delta$ parametrizing divisors at infinity can now be identified with $\Delta_X$.

Given a $\GG$-eigenform $\omega$ on $\XX$, the prior discussion
applied to $\overline{\XX}, \omega$ and the boundary divisors 
yields a form $\bar{\omega}$ -- henceforth denoted 
  $\omega_\Theta$ on $X_{\Theta}$; let us  verify it has the properties stated in Proposition \ref{propmeasure} and is in particular non-zero. In fact, the only property that needs to be verified is that it ``coincides'' with $\omega$ on open $\BB$-orbits under an isomorphism as in (\ref{Borbitident}); since $\BB$ is arbitrary, this implies that $\omega_\Theta$ is a $\GG$-eigenform, and the fact that it is an eigenform for $\AA_{X,\Theta}$ follows from the fact that $\tilde\omega$ is an eigenform for $\GGm^m$ and Lemma \ref{torilemma}. 

But for the restriction of $\omega$ to the open $\BB$-orbit, it is easy to see that $\omega_\Theta$ coincides with $\omega$ under (\ref{Borbitident}), using the Local Structure Theorem \ref{localstructure}.

This concludes the proof of Proposition \ref{propmeasure}.  \qed

	\subsection{Exponential map} \label{subsec:expmap}

This section plays a critical role in the paper: it shows that the asymptotic geometry of $X$ and $X_{\Theta}$ over a local field are ``the same'' in a suitable sense. Specifically, for any open compact subgroup $J$ of $G$ there is an identification between $J$-orbits in suitable neighborhoods of $\Theta$-infinity.  This is based on the usual idea that a normal bundle of a submanifold is diffeomorphic to a tubular neighborhood of that submanifold.

\begin{definition} Suppose $B,C$ are topological spaces\footnote{possibly with extra structure, 
e.g.\ locally ringed spaces, so that in particular the notion of ``morphism'' from an open subset of $B$ to $C$ is defined}. 

For any closed subspace $A \subset B$, and an open subset $U_A\subset A$ 
a {\em germ} at $U_A$ of a morphism to $C$ is
an equivalence class of pairs $(U_B,f)$: 
$$\{U_B: \mbox{ neighborhood of } U_A \text{ in } B, f: U_B \rightarrow C \mbox{ morphism} \}$$
under the equivalence relation $(U_B, f) \sim (U_B', f')$ if $f$ and $f'$ agree on a neighborhood of $U_A$.

The set of such germs, as $U_A$ is varying, forms a sheaf on $A$, which we will denote by $\underline{\mathrm{Mor}}_A(B,C)$. Its global sections over $A$ will be denoted by $\mathrm{Mor}_A(B,C)$. 
\end{definition}

In the setting that we are interested in, namely morphisms which are locally $p$-adic analytic on the $p$-adic points of smooth varieties, there is actually no difference between global sections of $\underline{\mathrm{Mor}}_A(B,C)$ and germs of morphisms: $U\to C$, where $U$ is a neighborhood of $A$ itself. 

 Let $\bar\XX$ be the wonderful compactification of $\XX$   or any smooth toroidal embedding (not necessarily complete),  and let $\ZZ\subset \bar\XX$ be the closure of a $\GG$-orbit belonging to $\Theta$-infinity (cf.\ \S \ref{Thetainftydisc}).  In this section, 
  we shall construct a canonical collection of elements
    $$  \expmap_{\Theta,J} \in \mathrm{Mor}_{Z/J}(N_Z{\bar{X}}/J , \bar{X}/J) $$ 
  where ``morphisms'' means (germs of) measure-preserving, continuous maps and
  $J$ ranges over all open compact subgroups. 
 In other words, fixing an open compact $J$, we have a way of transferring $J$ orbits
in a neighborhood of $Z$ in $N_Z{\bar{X}}$, to $J$-orbits in $\bar{X}$. 

{ By construction, the restriction of $\expmap_{\Theta, J}$ to $X_\Theta/J\subset N_Z{\bar{X}}/J $ will have image in $X/J\subset \bar X /J$ (in the sense that any representative of this germ has this property in a neighborhood of $Z$), and in Proposition \ref{expproperties} we will see that it is independent from the embedding and choice of orbit closure $\ZZ$, in the sense that all these germs of maps, obtained from different embeddings and orbits, glue together to give a well-defined germ of maps from a neighborhood of $\Theta$-infinity (cf.\ \S \ref{Thetainftydisc}) in $X_\Theta/J$ to a corresponding neighborhood in $X/J$. Thus, we have a well-defined element:
$$ \expmap_{\Theta, J} \in \mathrm{Mor}_{\infty_\Theta}(X_\Theta/J,X/J),$$
where now the notation stands for germs of maps defined in a neighborhood of $\Theta$-infinity. \emph{This is the sense in which the exponential map will be used in the largest part of the paper, i.e.\ without reference to a specific embedding}, but the geometric picture will be used in the construction and proofs. The collection of the elements $\expmap_{\Theta, J}$, as $J$ varies over a set of open-compact neighborhoods of the identity, will be denoted by $\expmap_\Theta$. We shall informally refer to $\expmap_\Theta$ and the various maps it induces on $J$-invariants as ``the exponential map'', because of its construction:  
}

Namely, the germ of $\expmap_{\Theta,J}$ will be induced by {\em 
any} $p$-adic analytic map $N_Z\bar X \rightarrow X$ inducing the identity
  on the normal bundle to $Z$ and respecting $G$-orbits.

\begin{proposition}\label{expdefinition}
Let $\ZZ$ be the closure of a $\GG$-orbit in $\bar\XX$, and let $\XX_1$ and $\XX_2$ be either of the varieties $N_{\bar\XX}(\ZZ)$ or $\overline{\XX}$.   There are locally $p$-adic analytic maps\footnote{Recall that this means that, in a neighborhood of every point, the map is given by a convergent power series
with respect to systems of local coordinates.}  $$\phi: U_1\to U_2,$$ where $U_i$ is a neighborhood of $Z$ in $X_i$ (henceforth called ``distinguished''), with the property that $\phi$ induces the identity between the normal bundles $N_\ZZ{\XX_1} \simeq N_\ZZ{\bar\XX} $ and $N_\ZZ{\XX_2} \simeq N_\ZZ{\bar\XX} $,
and that $\phi$ maps the intersection of every $\GG$-orbit with $U_1$ to the corresponding $\GG$-orbit on $\XX_2$.

Any such $\phi$ has the following property: 
Given an open compact subgroup $J\subset G$, there are $J$-invariant neighborhoods $U_1'\subset U_1, U_2'\subset U_2$ of $Z$ such that $\phi$ descends to a map: $U_1'/J\to U_2'/J$.

Finally, consider the open-compact topology on the space of such maps (with fixed domain). For every compact subset $\mathcal M$ of such maps, there are $J$-invariant neighborhoods $U_1'', U_2''$ of $Z$ in $X_1$, resp.\ $X_2$, such that all $\phi\in \mathcal M$ are defined and descend to the same map:
$$U_1''/J\to U_2''/J.$$
\end{proposition}

The germ of a map as in the proposition will be denoted by $\expmap_{X_1,X_2,J} \in \mathrm{Mor}_{Z/J}(X_1/J , X_2/J) $. 
When $\XX_1= N_{\bar\XX}(\ZZ)$ and $\XX_2=\overline{\XX}$, we also denote this by $\expmap_{\Theta,J}$ (where, as before, $\ZZ$ is the closure of a $\GG$-orbit 
on $\bar \XX$ belonging to $\Theta$-infinity). 

\begin{proof}
 
For a fixed neighborhood $U_1$ of $Z$ in $X_1$, consider the set $\mathcal K$ of locally $p$-adic analytic maps $\phi: U_1\to X_2$ with the properties of the proposition, that is: the differential of $\phi$ induces the identity on normal bundles, and $\phi$ maps points of a given $\GG$-orbit on $\XX_1$ to the corresponding orbit on $\XX_2$. We endow $\mathcal K$ with the open compact topology.  

\begin{lemma}\label{expstatement}
$\mathcal{K}$ is nonempty if $U_1$ is small enough. Moreover,  
 For every compact open subgroup $J$ of $G$ and any compact subset $\mathcal M$ of $\mathcal K$,  there is a neighborhood $U_1'\subset U_1$ of $Z$ such that the composites: $U_1' \overset{\phi}{\underset{\phi'}{\rightrightarrows}} X_2 \to X_2/J$ coincide, for any $\phi, \phi'\in\mathcal M$.
\end{lemma} 

Let us first see why this implies Proposition \ref{expdefinition}. Let $\mathcal M$ be a compact subset of such maps. The group $G$ acts on such maps by: $g\cdot \phi = g \circ \phi \circ g^{-1}$, and we may assume that $\mathcal M$ is $J$-invariant. But then, for $x$ in a subset $U_1'$ as in the lemma and all $j\in J$, we have: $j\circ \phi \circ j^{-1}(x) \in \phi(x) J = \phi'(x)J$ for all $\phi,\phi'\in \mathcal M$, and therefore the restrictions of all elements of $\mathcal M$ to $U_1'$ factors through $U_1'/J$ and give rise to the same map: $U_1'/J\to X_2/J$.

We now come to the proof of the lemma: Since we can glue {\em locally analytic} maps   
 we may replace $Z$ by an arbitrarily small open subset $Z'$ of it.  In other words,  
 if we have constructed maps $\phi_i: U_{1,i} \rightarrow X_2$ with the desired properties on an open covering of $U_1$ by open compact sets $U_{1,i}$,
 then we can refine the $U_{1,i}$ to a partition of $U_1$, and glue the restrictions of the $\phi_i$ to obtain $\phi \in \mathcal{K}$ as desired.
 
 We may assume that there is a Borel subgroup $\BB$ such that $Z'$ is contained in the open $\BB$-orbit, and use the Local Structure Theorem \ref{localstructure} to understand neighborhoods of $Z'$. Finally, we may replace $J$ by the its subgroup $J\cap B$, since it is stronger to prove that the projections to $X_2/J \cap B$ coincide.  Let $\YY$ be the toric variety of the Local Structure Theorem \ref{localstructure}, fix an isomorphism of the distinguished open $\BB$-subset $\bar\XX_B$ of $\bar\XX$ with $\YY\times \UU_{P(X)}$, and denote by $\ZZ':=\bar\XX_B\cap \ZZ$ (it is the closure of a $\BB$-orbit). We are left with proving:

\begin{quotation}
If $\XX_1', \XX_2'$ are either of the varieties $\YY\times\UU_{P(X)}$ or $N_{\YY\times\UU_{P(X)}}(\ZZ')$, then there are locally $p$-adic analytic maps: $U_1\to X_2'$, where $U_1$ is a neighborhood of $Z'$ in $X_1'$, inducing the identity on normal bundles and preserving the points of corresponding $\BB$-orbits; moreover, for any compact subgroup $J_B$ of $B$ and any compact subset $\mathcal M$ of such maps (defined on a fixed $U_1$) there is a smaller neighborhood $U_1'\subset U_1$ where the composites: $U_1' \overset{\phi}{\underset{\phi'}{\rightrightarrows}} X_2' \to X_2'/J_B$ coincide for any $\phi, \phi'\in\mathcal M$.
\end{quotation}
 
The statement is now easily reduced to the analogous statement about smooth toric varieties, indeed eventually to the case where
$\XX_1'', \XX_2''$ are either of the spaces  $\mathbf{V}$ or $N_{\ZZ}(\mathbf{V})$, 
when $\mathbf{V}$ is an affine space $(\mathbf{A}^1)^n$ and $\ZZ$  is the intersection of all coordinate hyperplanes, i.e. the origin of $\mathbf{V}$,
and we require the map to preserve all coordinate hyperplanes.    

The other assertion -- concerning the induced maps $U_1' \rightarrow X_2'/J_B$ --  reduces similarly to the following assertion: 
Given a locally analytic morphism $f: k^n \rightarrow k^n$ which preserves coordinate
axes, so that $f(0,0, \dots, 0) = (0,\dots, 0)$, and so that the derivative
of $f$ at zero is the identity map, then $f$ maps
each $J$-orbit near $\underline{0}$ to itself, if $J \subset (k^{\times})^n$ is open compact
acting on $k^n$ by coordinate multiplication.  Moreover, the notion of ``near'' can be taken
to be uniform if $f$ lies in a compact set of such maps. 
To see this, one notes that the Taylor expansion $f_j(x_1, \dots, x_n) = x_j (1 + \dots)$:
all higher order terms are divisible by $x_j$ because of preservation of coordinate axes.

(Notice that the requirement that the maps $\phi$ preserve orbits is made necessary by the fact that open compact subsets of the group do not provide a good uniform structure in the neighborhood of non-open orbits.)\end{proof}

We now come to the properties of $\expmap_{\Theta,J}$, preserving the notation of the previous proposition.

\begin{proposition}\label{expproperties}
 Any representative $\phi_J$ of $\expmap_{X_1,X_2,J} \in \mathrm{Mor}_{Z/J}(X_1/J,X_2/J)$ has the following properties: 
\begin{enumerate}
 \item it is eventually (that is: after restricting its domain to a smaller neighborhood of $Z/J$, if necessary) a measure-preserving bijection;

 \item it is eventually equivariant: for any $h\in \mathcal H(G,J)$ (the Hecke algebra of $G$ with respect to $J$), if $U_1\subset X_1/J$ denotes the domain of definition of $\phi_J$ then  there is a smaller neighborhood $U_1'$ of $Z/J$ such that for any $J$-invariant function $f$ on $U_1'$ we have:
$$ h*\phi_{J*} f = \phi_{J*} (h*f),$$
where $\phi_{J*}f$ denotes the push-forward of $f$ through the bijection $\phi_{J}$. 
 { \item Its restriction to $X_\Theta\subset N_Z\bar X$ (cf.\ Proposition \ref{thetaindependent}) does not depend on the embedding $\bar X$, in the following sense: for any two embeddings $\overline{\XX}$, $\overline{\XX}'$, orbit closures $\ZZ, \ZZ'$ as before and representatives $\phi_J, \phi_J'$ for the corresponding germs $\expmap_{\Theta,J}$, there are $J$-stable neighborhoods $N_\Theta, N_\Theta'$ of $Z$, resp.\ $Z'$ in $X_\Theta$, such that $\phi_J|_{(N_\Theta\cap N_\Theta')/J} = \phi_J'|_{(N_\Theta\cap N_\Theta')/J}$. Hence, by working with all orbit closures belonging to $\Theta$-infinity in a wonderful compactification, we get a well-defined germ $\expmap_{\Theta,J} \in \mathrm{Mor}_{\infty_\Theta} (X_\Theta/J,X/J)$ which does not depend on choices.  } 
\end{enumerate}

{Notice that, by the third statement, the ``eventual equivariance'' property of the second statement extends to a neighborhood of $\Theta$-infinity, when applied to smooth functions on $X_\Theta$.}

\end{proposition}
 
\begin{proof}
 
The preservation of volume follows immediately from the existence of distinguished equivariant morphisms between open $\BB$-orbits (\ref{Borbitident}) and the characterizing property of the measure on $X_\Theta$ (Proposition \ref{propmeasure}). 

We come to the proof of eventual equivariance: It is enough to prove it for elements of the Hecke algebra of the form $h=JgJ$. Write the double coset $JgJ$ as a (finite) disjoint union of right cosets: 
$$JgJ = \sqcup_i g_i J.$$
Let $J' = J\cap \bigcap_i g_i J g_i^{-1}$ (an open compact subgroup of $G$), and choose a distinguished map $\phi$ giving rise to representatives $\phi_J$ and $\phi_{J'}$ for $\expmap_{\Theta,J}$, $\expmap_{\Theta,J'}$.

We obviously have, eventually:
$$ \phi_{J*} (h*f) = \phi_{J*} (\sum_i g_i\cdot f) = \sum_i \phi_{J'*} (g_i\cdot f).$$
Therefore, it is enough to show that
$$ \phi_{J'*} g_i \cdot f = g_i \cdot \phi_{J'*} f.$$
But the maps: $\phi^{-1} g_i^{-1} \phi g_i$ are also distinguished, so the result follows from Proposition \ref{expdefinition}.

{The independence from the orbit closure for a \emph{given} embedding is seen as follows: If two orbit closures $\ZZ, \ZZ'$ belonging to $\Theta$-infinity do not intersect, then they have disjoint neighborhoods, so there is nothing to prove. Otherwise, their intersection also belongs to $\Theta$-infinity, and it is therefore enough to prove independence when one (say, $\ZZ$) is contained in the other. In that case we have natural identifications of normal bundles:
$$ N_\ZZ \overline{\XX} = N_\ZZ(N_{\ZZ'}\XX).$$
If we set $\YY= N_{\ZZ'}\overline{\XX}$, there is an ``exponential map'':
$$\phi: N_\ZZ\YY\to \YY,$$
i.e.\  $p$-adic analytic map
fixing $\ZZ$ and inducing the identity on its normal bundle, which is the identity on the open $\GG$-orbits, both identified with $\XX_\Theta$ via Proposition \ref{thetaindependent}; indeed, this can easily be seen by invoking the normal bundle degeneration and the Local Structure Theorem \ref{localstructure}. Composing with an exponential map from $\YY$ to $\XX$ we see that in a neighborhood of $\ZZ$ the two exponential maps coincide.

Thus we have a well-defined germ $\expmap_{\Theta,J}$ of maps in a neighborhood   of \emph{the $\Theta$-infinity of the given embedding} (a priori, depending on the embedding). Given two smooth toroidal embeddings $\overline{\XX}, \overline{\XX}'$, now, by the Luna-Vust theory we can find a third one $\overline{\XX}''$, open $\GG$-invariant subsets $\UU,\UU'\subset \overline{\XX}''$ and proper morphisms: $\UU\to \overline{\XX}, \UU'\to \overline{\XX'}$. Indeed, such an embedding can be obtained by constructing a fan, as described in \S \ref{sssnonwonderful}, which contains a partition of $\mathcal C$, for every cone $\mathcal C$ in the fan of $\overline{\XX}$ or $\overline{\XX}'$.
 
 It is easy to see that any representative for $\expmap_{\Theta,J}$ on $\overline\XX$, resp.\ $\overline{\XX'}$ pulls back to a representative for $\expmap_{\Theta,J}$ on $\UU$, resp.\ $\UU'$, and their germs glue together to the germ $\expmap_{\Theta,J}$ for $\overline\XX''$.}

\end{proof}

{
\subsubsection{Transitivity}\label{ssstransitivity}

Let $\overline{\XX}$ be a smooth toroidal embedding of $\XX$, $\ZZ$ an orbit closure belonging to $\Theta$-infinity, and $\ZZ'\subset\ZZ$ an orbit closure belonging to $\Omega$-infinity, for some $\Omega\subset\Theta$. Then $N_{\ZZ} \overline{\XX}$ is a smooth toroidal embedding of $\XX_\Theta$, and hence it is clear from the identity:
$$N_{\ZZ'} \left(N_{\ZZ} \overline{\XX}\right) = N_{\ZZ'} \overline{\XX} $$ 
that $\XX_\Omega$ is canonically identified with the corresponding boundary degeneration of $\XX_\Theta$. (One can also argue that using the affine degeneration of \S \ref{ssdegen}.)

If $\phi:N_{Z'}\bar X \rightarrow X$ is a $p$-adic analytic map inducing the identity on the normal bundle to $Z'$ and respecting $G$-orbits, then its partial differential along $Z$:
$$ N_{Z'} \bar X \to N_Z \bar X$$
also has the same properties. Hence, the exponential maps satisfy the transitivity property:

\begin{equation}\label{eqtransitivity}
 \expmap_\Omega = \expmap_\Theta \circ \expmap_\Omega^\Theta,  
\end{equation}
where by $\expmap_\Omega^\Theta$ we denote the corresponding exponential map for $X_\Theta$ (i.e.\ a compatible collection, over $J$, of elements of $\mathrm{Mor}_{Z/J}(X_\Theta/J,X/J)$).
}

\section{Asymptotics}\label{sec:asymptotics}

\textbf{From now on we assume that $\XX$ is \emph{wavefront}, {\bluetext with $\mathcal Z(\GG)^0$ surjecting onto $\mathcal Z(\XX)$.}}

\subsection{The main result} 
The main goal of this section is to relate the asymptotic behavior of eigenfunctions on $X$ to eigenfunctions on the boundary degenerations. Here we use the word ``eigenfunctions'' freely, meaning elements of $C^\infty(X)$ generating irreducible representations. As we saw in the previous section  \S\ref{ssmeasures}, the measure on $X$ canonically induces an $A_{X,\Theta}\times G$-eigenmeasure, with the same $G$-eigencharacter, on each boundary degeneration $X_\Theta$, which allows us to formulate the main theorem.

Recall that the notion of $\Theta$-infinity has been introduced in \S \ref{Thetainftydisc}.

\begin{theorem}\label{mainasthm}  
For every $\Theta\subset\Delta_X$ there exists a unique $G$-morphism
\begin{equation} e_\Theta: C_c^{\infty}(X_\Theta) \longrightarrow C_c^{\infty}(X)\end{equation}  with the property that for every open compact subgroup $J\subset G$ and any representative $\phi_J$  of $\expmap_{\Theta,J}$ there is a ($J$-stable) neighborhood $N_\Theta$ of $\Theta$-infinity such that for all $f\in C_c^\infty(N_\Theta)^J$ we have:
\begin{equation}\label{eqasymptotics}
 e_\Theta(f) = \phi_{J*} (f).
\end{equation}
{ In fact, the morphism is characterized by the validity of \eqref{eqasymptotics} in a neighborhood of \emph{some} orbit closure $\ZZ$ belonging to $\Theta$-infinity in \emph{some} smooth toroidal embedding $\overline{\XX}$}.
\end{theorem}

The notation $e_{\Theta}$ derives from  ``exponential map'', but also from the name ``Eisenstein series'', because the global analog of $e_\Theta$ is the construction of pseudo-Eisenstein series.
A neighborhood $N_\Theta$ of $\Theta$-infinity as in the theorem will be called \emph{``$J$-good''}, { and $N_\Theta/J$ will be identified with a subset of both $X/J$ and $X_\Theta/J$ via the map $\phi_J$ as above (which is now a canonical representative of the exponential map, restricted to $N_\Theta/J$). We will sometimes, by abuse of language, treat $N_\Theta$ itself as a subset of both $X$ and $X_\Theta$, when the statements that we are making are really about $N_\Theta/J$.}

We observe that, if $N_\Theta$ denotes a neighborhood of $\Theta$-infinity
for each $\Theta \subset \Delta_X$, then $\bigcup_{\Theta\ne \Delta_X} N_\Theta$  
necessarily has compact-modulo-center complement inside $X$, as follows from the compactness of a wonderful embedding
$\overline{X}$.  Therefore, the theorem
indeed controls the asymptotics in all directions simultaneously.

{ On the other hand, for given $\Theta$, the last statement of the theorem shows that the map $e_\Theta$ is characterized by its restriction along a \emph{unique} direction towards $\Theta$-infinity, in the following sense: Recall from the Luna-Vust theorem \ref{thmlunavust} that to any half-line in $\mathcal V=$ the cone of invariant valuations for $\XX$ we can attach a smooth toroidal embedding (where smoothness follows from the local structure theorem \ref{localstructure}). This embedding has a unique non-open $\GG$-orbit, and by choosing the half-line in the interior of the face corresponding to $\Theta$, this $\GG$-orbit will belong to $\Theta$-infinity.}

Dually, we have a morphism: 
\begin{equation} e_\Theta^*: C^\infty(X)\to C^\infty(X_\Theta).
\end{equation}

Theorem \ref{mainasthm} is equivalent to:
\begin{theorem} \label{mainasthm2} There is a unique $G$-morphism
$$e_\Theta^*: C^\infty(X) \to C^\infty(X_\Theta)$$
{ with the property that for every open compact subgroup $J\subset G$ and any representative $\phi_J$  of $\expmap_{\Theta,J}$ there is a ($J$-stable) neighborhood $N_\Theta$ of $\Theta$-infinity such that for all $f\in C^\infty(X)^J$ we have:
\begin{equation}\label{eqasymptotics2}
 e_\Theta^*(f)|_{N_\Theta} = \phi_J^* (f|_{N_\Theta}).
\end{equation}
 In fact, the morphism is characterized by the validity of \eqref{eqasymptotics2} in a neighborhood of \emph{some} orbit closure $\ZZ$ belonging to $\Theta$-infinity in \emph{some} smooth toroidal embedding $\overline{\XX}$}.
\end{theorem}

For every smooth representation $\pi$ of $G$ and any $G$-equivariant map: $M:\pi\to C^\infty(X)$, the composition with this morphism gives rise to a $G$-morphism: 
\begin{equation}\label{asymptoticsmap}\Hom_G(\pi,C^\infty(X))\ni M\mapsto M_\Theta\in \Hom_G(\pi,C^\infty(X_\Theta))
\end{equation}
which will be called the ``asymptotics'' map.
It has the property that for any $v\in \pi^J$ we have:
\begin{equation}\label{asympcond}
 M(v)|_{N_\Theta} = M_\Theta(v)|_{N_\Theta'}.
\end{equation} 
Moreover, \eqref{asympcond} uniquely characterizes $M_{\Theta}$.

We shall give two proofs of this result:
\begin{enumerate}
\item In \S \ref{SSasymptotics} we formulate a ``morally satisfactory'' proof based on the ``stabilization theorem''
of Bernstein and the $\exp$-map;
\item In \S \ref{ssCartan} we give a ``quick and dirty'' proof, using the known results about asymptotics of smooth matrix coefficients (of course, this sweeps under the carpet all arguments of the previous approach, which are used to establish the result in the group case). This method requires some additional piece of information, namely a (weak) generalized form of the Cartan decomposition which can be derived from the geometry of the wonderful embedding, in order to show that $G_x$-invariant functionals (where $x\in X$) can be computed using \emph{smooth} matrix coefficients.\footnote{This observation has already appeared in work of Lagier \cite{Lagier} and Kato-Takano \cite{KatoTakano}.} 
 \end{enumerate}
 
 \begin{remark}
 \begin{itemize}
 \item[(i)]  Recently, Bezrukavnikov and Kazhdan gave a  proof of  
 second adjunction \cite{BK} using the geometry of the wavefront compactification
 in the group case. In particular, they give in \S 4 a beautiful abstract approach to 
 essentially the same problem (although phrased in a special case, their method adapts without change to the current situation).  In the current context, it gives another proof  of asymptotics, without using Bernstein's stabilization theorem. 
 It uses as input certain finite generation statements such as Remark \ref{remark AAG}; in our context, we obtain these {\em a posteriori} from the asymptotics and the knoweldge that $\X$ is wavefront. Although 
 the argument can be reordered so that the proof of \cite{BK} goes through, it does not allow us to bypass the requirement that $\X$ be wavefront.

\item[(ii)] Also, we observe that our proofs do not require Bernstein's results if, for instance, one is interested
 only in the case of $\pi$ admissible (as is the case with the unitary theory in this paper). In that case, one can easily see that the usual facts about the Jacquet module suffice. Bernstein's results are used to generalize from the case of admissible representations to general smooth representations -- a generalization that, although conceptually very pleasing, we do not strictly require for the  the largest part of the paper (except for section \ref{sec:explicit}). 
\end{itemize}
\end{remark}

\begin{remark} \label{transitivity}
 If $\Theta\supset \Omega$ then we can apply the theorem to the variety $\XX_\Omega$ in order to get a morphism $e_\Theta^\Omega: C_c^\infty(X_\Theta)\to C_c^\infty(X_\Omega)$. Clearly, we have the transitivity property: $e_\Omega\circ e_\Theta^\Omega =  e_\Theta$, since the exponential maps have the same transitivity property \eqref{eqtransitivity}.
\end{remark}

The theorem implies the following:

\begin{theorem}[Finiteness of multiplicities.] \label{finiteness}
Let $\pi$ be an irreducible smooth representation of $G$, and let $\XX$ be a wavefront spherical variety {\bluetext with $\mathcal Z(\GG)^0\twoheadrightarrow\mathcal Z(\XX)$}. Then: $$\dim\Hom_G(\pi,C^\infty(X))< \infty.$$
\end{theorem}

\begin{proof}
First, we claim:

\begin{quote}
 For any $\Theta\subset\Delta_X$, we have $\dim\Hom_G(\pi,C^\infty(X_\Theta))< \infty$ if and only if $\dim\Hom_{A_{X,\Theta} \times G}(\chi\otimes \pi,C^\infty(X_\Theta))< \infty$ for every character $\chi$ of $A_{X,\Theta}$.
\end{quote}

Indeed, recall that $C^\infty(X_\Theta)$ is induced from the $P_\Theta^-$ representation $C^\infty(X_\Theta^L)$, and by Proposition \ref{wavefrontlevi} the action of $A_{X,\Theta}$ is induced by the center of the corresponding Levi $L_{X,\Theta}$ (up to possibly a finite index due to the fact that the map of $k$-points: $\mathcal Z(L_\Theta)^0\to A_{X,\Theta}$ may not be surjective). Since $\pi$ is irreducible, there is only a finite number of distinct characters $\chi$ of $\mathcal Z(L_{X,\Theta})^0$ such that $\pi$ could be embedded in a representation induced via $P_\Theta^-$ from a representation of $L_\Theta$ with $\mathcal Z(L_\Theta)^0$-character $\chi$. Therefore, finiteness for every character $\chi$ implies finiteness for $\pi$, forgetting the $A_{X,\Theta}$-action.

{\bluetext Now we may assume, by induction, that the theorem is true for $\XX_\Theta$ under the $\AA_{X,\Theta}\times\GG$-action, for every $\Theta\subsetneq \Delta_X$. Recall from Proposition \ref{wavefrontlevi} that $\XX_\Theta$ is also wavefront under this action.}

Now, the common kernel of all the asymptotics maps (\ref{asymptoticsmap}) (excluding $\Theta = \Delta_X$), consists of morphisms: $\pi\to C^\infty(X)$ such that the image of $\pi^J$ is supported on the complement of all neighborhoods $N_\Theta$ of Theorem \ref{mainasthm}. This complement has a finite number of $J$-orbits modulo the action of $\mathcal Z(X)$, and therefore, since $\mathcal Z(G)^0\twoheadrightarrow\mathcal Z(X)$), $\dim\Hom_G(\pi,C^\infty(X))< \infty.$ 
\end{proof}

\begin{remark} This proof only gives a bound for the dimension of $\Hom_G(\pi, C^{\infty}(X))$ which depends on the level of $\pi$, i.e.\ on which subgroup $J$ is such that $\pi^J\ne 0$. It is natural to ask whether there exists a bound {\em independent of $\pi$}. 
A plausible such upper bound would be the generic multiplicity of unramified principal series, computed in \cite[Theorem 5.3.2]{SaSpc}. 
\end{remark}

\begin{remark} \label{remark AAG}
 Aizenbud-Avni-Gourevitch \cite{AAG}
 have established the following result (using, as input, the finiteness of multiplicities):
If $J$ is an open compact subgroup of $G$, then
\begin{quote} $C^{\infty}_c(X)^J$ is finitely generated as a module
under the $J$-Hecke algebra $\mathcal{H}(G,J)$.
\end{quote}

Let us sketch how the theorem may be used to give another proof of this result: one argues just as in the previous result, but using the dual maps
$C^{\infty}_c(X) \longrightarrow C^{\infty}_c(X_{\Theta})$, {\em together with} the fact
(due to Bernstein) \cite[3.11]{BeCentre} that if $P$ is a parabolic with associated
Levi decomposition $P=MN$ and $J$ admits the Iwahori decomposition $J = J_{-} J_M J_{+}$, 
then parabolic induction maps finitely generated $\mathcal{H}(M, J_M)$-modules
to finitely generated $\mathcal{H}(G, J)$-modules. 
\end{remark}

Finally, we introduce some language in order to describe another corollary of the asymptotics. 
This will be only used much later (Proposition \ref{proofextend}) to show meromorphic continuation
of certain intertwiners related to ``Eisenstein integrals.'' 

A function on an abelian group is said to be {\em finite} if its translates (under the action of that abelian group on itself by multiplication) span a finite-dimensional vector space. 
For a normal $k$-variety $\VV$ with a distinguished divisor $\DD$, we will say that a complex-valued function $F$ on $V$ is $D$-finite if for every $x\in V$ there is a neighborhood of $V_x$ of $x$, 
and 
rational functions $f_1, \dots, f_m$, whose zero and polar divisors are contained in $\DD$ such that,  
on $V_x - V_x \cap D$, 
the function $F$ agrees with the pullback of a finite   function 
  on $\mathbb{G}_m^m$ by $(f_1, \dots, f_m)$.   
  
  The notion of finite function is stable under pullback, i.e.
given a morphism $\pi: \VV_1 \rightarrow \VV_2$ of algebraic varieties,
a divisor $\DD_2 \subset \VV_2$, and a function $f_2$ on
$V_2$ that is $D_2$-finite, the pullback $\pi^* f_2$ is $D_1$-finite,
where $\DD_1 = \pi^{-1} \DD_2$. 

\begin{corollary} \label{finiteasymp} 
Let $\bar\XX$ be a smooth toroidal embedding of $\XX$, and let $\nu:\pi\to C^\infty(X)$ be a morphism from an \emph{admissible} representation. Then the image of $\nu$ consists of $(\bar X\smallsetminus X)$-finite functions.
\end{corollary}

The statement makes sense since, as we have seen from the Local Structure Theorem \ref{localstructure}, the complement of the open $\GG$-orbit in a smooth toroidal embedding is a union of divisors intersecting transversely.

\begin{proof}
 For any point $x\in \bar X$, assumed to belong to a $\GG$-orbit $\ZZ$ corresponding to $\Theta\subset\Delta_X$, we can find by (\ref{Borbitident}) algebraic  identifications of its neighborhoods in $X$ and $X_\Theta$ which are compatible with the exponential map. Therefore, it is enough to show that the image of $e_\Theta^*\circ \nu$ consists of finite functions (with respect to the complement of $X_\Theta$ in $N_{\bar Z}\bar X$). But then this notion of finiteness can be seen to be equivalent  to finiteness under the action of the torus of Lemma \ref{torusacting}, which is a subgroup of $\mathcal Z(X_\Theta)$, and the image of admissible representations is clearly $\mathcal Z(X_\Theta)$-finite.
\end{proof}

\subsection{Proof of asymptotics}  \label{SSasymptotics}

In this subsection we prove Theorem \ref{mainasthm}, or rather its adjoint Theorem \ref{mainasthm2}. We fix a $\Theta\subset\Delta_X$, a smooth toroidal embedding $\overline{\XX}$ and a $\GG$-orbit closure $\ZZ$ belonging to $\Theta$-infinity. From now on we will be denoting by $\exp_\Theta$ a representative for $\expmap_{\Theta,J}$, where $J$ is an open compact subgroup under which the functions under consideration are invariant. The dependence on $J$ is suppressed from the notation, since the choice of $J$ does not matter for the statements.

We assume for notational simplicity that $G$ has a unique orbit on $X_\Theta$. Fix a parabolic subgroup in the class of $\PP_\Theta^-$ (denoted thus), and recall from \S \ref{levivarieties} that this identifies a subvariety of $\XX_\Theta$ as the Levi variety $\XX_\Theta^L$; its closure in $\overline{\XX}$ intersects $\ZZ$ along a $\PP_\Theta^-$-stable set $\ZZ_1$. (Both the $\XX_\Theta^L$ and $\ZZ_1$ are the sets of all elements of $\XX_\Theta$, resp.\ $\ZZ$ whose stabilizer $\GG_{z}$ belongs to $\PP_\Theta^-$, or equivalently: $\UU_\Theta^- \subset \GG_z \subset \PP_\Theta^- $.)
We have a map: $\XX_\Theta^L\to \ZZ_1$ arising from the structure of $\XX_{\Theta}$ as a subset of the normal bundle to $\ZZ$, and we denote by $\mathrm{fiber}(z)$ the preimage of a given point $z\in Z_1$. Note that $\mathrm{fiber}(z)$ is homogeneous under $A_{X, \Theta}=\mathcal{Z}(X_{\Theta})$.

Then Theorem \ref{mainasthm2} is equivalent to the following:
\begin{proposition} \label{asymptoticstool} (Assuming a unique $G$-orbit on $X_\Theta$.)

There is a system $\left(\mathcal{N}_{ J}\right)_J$ of $J\cap P_\Theta^-$-stable neighborhoods  of $Z_1$ in \footnote{\label{ftntnbhd} By this, we more properly mean that $\mathcal{N}_J$ is of the form
$\mathcal{N}'_J \cap X_\Theta^L$, where $\mathcal{N}'_J$ is a neighborhood of $Z$ in its normal bundle.} $X_\Theta^L$, as $J$ varies over a basis of open compact subgroup-neighborhoods of the identity in $G$, with the following property:

For any $z$ in the open $P_\Theta^-$-orbit (equivalently: $L_\Theta$-orbit) in $Z_1$, there is a $G$-morphism: $$C^\infty(X) \ni f \mapsto f_\Theta \in C^\infty(X_\Theta)$$
(a priori, depending on $z$) such that, for any $J$-invariant $f$, 
 \begin{equation} \label{charprop} f_{\Theta}  = \exp_{\Theta}^* f 
 \mbox{ on $\mathrm{fiber}(z) \cap   \mathcal{N}_{J}$.} \end{equation}
 
Moreover,  there is a unique $G$-morphism which has property \eqref{charprop} for \emph{some} $z$ and for \emph{some} such system of neighborhoods of $z$ in $\mathrm{fiber(z)}$.
    \end{proposition}
    
{ When there are multiple $G$-orbits on $X_\Theta$ or, rather, on the open orbit in $Z$, one just needs to choose more than one $z$'s representing all orbits.}

Of course, the proposition follows from Theorem \ref{mainasthm2} (applied to the given embedding $\overline{\XX}$ and orbit closure $\ZZ$), but we will see that, vice versa, the validity of the proposition for \emph{every} pair $(\overline{\XX},\ZZ)$ implies the full theorem: 

\begin{proof}[{Proof of Theorem \ref{mainasthm2} assuming Proposition \ref{asymptoticstool}}]

We first check that the morphism supplied by the proposition
is independent of choice of $z$, so let us temporarily denote the $f_\Theta$ provided by the proposition by $f_\Theta^z$. Notice that $P$ acts transitively on the $k$-points of the open $\PP$-orbit of $\ZZ_1$, since it is self-normalizing in $G$ (and we are assuming that $G$ acts transitively on $X_{\Theta}$, hence also on the points of the open $\GG$-orbit on $\ZZ$). For $z' := zg$, $g \in P$, and $f$ a $g^{-1} J g$-invariant function, the functions $(g f)_{\Theta}^z = g f_{\Theta}^z $ and $\exp_{\Theta}^* (gf)$ coincide on  
$\mathcal{N}_{J} \cap \mathrm{fiber}(z)$; so $f_{\Theta}^z$
and $g^{-1} \exp_{\Theta}^* (gf)$ coincide on $\mathcal{N}_{J}\cdot g \cap \mathrm{fiber}(z')$. 
Now, the system $\left(\mathcal{N}_{J} \cdot g\right)_J$ is still a system of neighborhoods of $Z_1$, and we can replace it by a system of smaller neighborhoods $\left(\mathcal{N}'_{J}\right)_J$ so that
 $g^{-1} \exp_{\Theta}^*(gf) = f$ there (for any such $f$).
Then 
$f_{\Theta}^z$ and $\exp_{\Theta}^*(f)$ coincide on $\mathcal{N}_{J}' \cap \mathrm{fiber}(z')$. 
The uniqueness statement of the Proposition implies that
$f^z_{\Theta}=f^{z'}_\Theta$.

Now we show that the map $f\mapsto f_\Theta$ has the properties of Theorem \ref{mainasthm2}, {  first for the given pair $(\overline{\XX},\ZZ)$. First of all, by continuity of the map $\exp_\Theta$ and the validity of the above proposition in a neighborhood (s. footnote ${\ref{ftntnbhd}}$) of $Z_1$, property \eqref{charprop} holds for any $z\in Z_1$, not only in the open orbit.} The space $P_\Theta^-\backslash G$ is compact; let $K_1\subset G$ be a compact preimage of it. Given an open compact subgroup $J$, let $J_0$ be the intersection of $kJk^{-1}$, $k\in K_1$, and denote by $\tilde{\mathcal N}_{ J}$ the union of the sets $\mathcal N_{ J_0} k$, $k\in K_1$; it is a neighborhood of $\Theta$-infinity (the closure of $Z$) in $X_\Theta$. If $f\in C^\infty(X)^J$ and $z\in \tilde{\mathcal N}_{ J}$, there is a $k\in K_1$ such that $zk^{-1}\in \mathcal N_{ J_0} $, and the function $k\cdot f$ is $kJk^{-1}$-invariant, in particular $J_0$-invariant. Hence we have:
$$ f_\Theta(z) = k f_\Theta (zk^{-1}) = \exp_{\Theta}^*(kf) (z k^{-1}) = k^{-1} \exp_{\Theta}^*(kf) (z).$$
Again, by the compactness of $K_1$ we may replace $\tilde{\mathcal N}_{ J}$ by a smaller neighborhood of $\Theta$-infinity on which $k^{-1} \exp_{\Theta}^*(kf) = f$ for every $k\in K_1$ and $f\in C^\infty(X)^J$. This shows that \eqref{eqasymptotics2} holds in a neighborhood $N_\Theta$ of the orbit closure $\ZZ$ in the given embedding $\XX$.

{ 
Now, as in the proof of Proposition \ref{expproperties}, for two different pairs $(\overline{\XX},\ZZ)$, $(\overline{\XX}',\ZZ')$ we will work with a third smooth toroidal embedding $\overline{\XX}''$ which contains open $\GG$-stable subsets $\UU,\UU'$ properly dominating $\overline{\XX}$, resp.\ $\overline{\XX}'$. If the morphisms $f\mapsto f_\Theta$ obtained from different orbit closures belonging to $\Theta$-infinity in $\overline{\XX}''$ are equal, then so are the morphisms obtained from the pairs $(\overline{\XX},\ZZ)$, $(\overline{\XX}',\ZZ')$. Thus, the existence statement of the theorem is reduced to the case when $\overline{\XX}=\overline{\XX}'$ but the orbit closures $\ZZ$, $\ZZ'$ are different. Moreover, again by passing to another embedding, we may assume that \emph{$\Theta$-infinity is connected}; this amounts to saying that the support of the fan of $\overline{\XX}$ with the relative interior of the face of $\mathcal V$ (the cone of invariant valuations) corresponding to $\Theta$ is connected. This includes, for example, the elementary case $\XX=\GGm\subset \XX = \PP^1$, $\ZZ=\{0\}, \ZZ'=\{\infty\}$ (where $\Theta= \Delta_X=\emptyset$ and the open orbit $\GG_m$ also belongs to $\Theta$-infinity).  

Under these assumptions, if we consider the equivalence relation between all orbit closures belonging to $\Theta$-infinity generated by: $\ZZ\sim \ZZ' \iff \ZZ\supset\ZZ'\mbox{ or }\ZZ'\supset\ZZ$, then all orbits are equivalent, and we are reduced to considering the case of two pairs $(\XX,\ZZ)$ and $(\XX,\ZZ')$ with $\ZZ\subset\ZZ'$. We need to prove that the corresponding morphisms: $C^\infty(X)\to C^\infty(X_\Theta)$ (temporarily to be denoted by $f\mapsto f_\Theta^Z$ and $f\mapsto f_\Theta^{Z'}$) are identical. 

In this case, the validity of the theorem for $N_\Theta=$ a neighborhood of $\ZZ$ is evidently weaker than its validity for a neighborhood of $\ZZ'$. However, identifying $\XX_\Theta$ with the open $\GG$-orbit in both $N_{\ZZ}\overline{\XX}$ and $N_{\ZZ'}\overline{\XX}$, a neighborhood of $Z$ in $X_\Theta$ \emph{does} include a neighborhood of $z$ in $\mathrm{fiber}(z)$, for some point $z$ in the open $G$-orbit on $Z'$. Thus, the map $f\mapsto f_\Theta^{Z}$ satisfies \eqref{charprop} in that neighborhood of $z$ in $\mathrm{fiber}(z)$, and by the uniqueness statement of Proposition \ref{asymptoticstool} it has to coincide with the map $f\mapsto f_\Theta^{Z'}$. This proves Theorem \ref{mainasthm2}.

The characterization (uniqueness) statement follows from the uniqueness statement of Proposition \ref{asymptoticstool}.

}

 \end{proof}

 \subsubsection{Setup for the proof of Proposition \ref{asymptoticstool}}

Let us fix $z\in Z$. The idea in the proof of the proposition is to replace the action of the monoid 
$A_{X, \Theta}^+$ along the orbit of a point $x\in \mathrm{fiber}(z)$ by the right action of a subtorus of $G$, or rather a subalgebra of the Hecke algebra of $G$. 

For that reason, choose a Levi subgroup $\LL_\Theta$ of $\PP_\Theta^-$.
We have seen in Proposition \ref{wavefrontlevi} that $\mathcal{Z}(\LL_{\Theta})^0$ surjects onto $\AA_{X, \Theta}$
and this induces a surjection of positive chambers at the level of Lie algebras (as noted  at the end of the proof of Proposition \ref{wavefrontlevi}). 
In what follows we use $a,b$ to denote elements of $\mathcal{Z}(\LL_{\Theta})^0$.   For $a \in \mathcal{Z}(L_{\Theta})^0$,
we write $|a| = \max_{\delta} |\delta(a)|$, the maximum being taken over all {\em negative} roots for
$\mathcal Z(\LL_{\Theta})$ with respect to $\PP_{\Theta}$, or what is the same, the positive roots with respect to $\PP_\Theta^-$. We let $\mathcal{Z}(L_{\Theta})^+\subset \mathcal{Z}(L_{\Theta})^0$ denote the set of elements where $|a| \leq 1$,
i.e.\ $|\chi(a)| \leq 1$ for all negative roots with respect to $\PP_{\Theta}$. 
 
 Note that $\mathcal{Z}(\LL_{\Theta})^0 \twoheadrightarrow \AA_{X, \Theta}$
induces $\mathcal{Z}(L_{\Theta})^+ \rightarrow A_{X, \Theta}^+ $ (this map may not be surjective
because of the operation of taking $k$-points, however).  

As elsewhere, we write, for $a \in \mathcal{Z}(L_{\Theta})^+$, that ``$a$ is sufficiently deep'' in place of `` there exists $\varepsilon > 0$ so that,
whenever $|a| \leq \varepsilon$, ... '' 
Finally, for $x\in \mathrm{fiber}(z)$ we put
$$x_a := x \cdot a.$$

\begin{lemma} 
 For every open compact subgroup $J$ with Iwahori factorization with respect to $P_\Theta$ and $P_\Theta^-$, the elements of the Hecke algebra:
$$h_a:=1_{J aJ}, a \in \mathcal{Z}(L_{\Theta}), |a| \leq 1.$$  
satisfy: $h_a\star h_b = h_{ab}$ (where $b$ is also an element of $\mathcal{Z}(L_{\Theta})$ satisfying $|b| \leq 1$). 
\end{lemma}

This is straightforward and well-known. As a consequence, the vector space spanned by these elements is a subalgebra $\mathcal H$ of the Hecke algebra of $G$. For now we fix such a subgroup $J$ and the notation of the lemma. 

We notice the following   fact: for $|a| \leq 1$  
and arbitrary $b$ we have $x_b\cdot JaJ = x_{ba} J\subset X_\Theta$
 because, recalling that $x$ is stabilized by the unipotent radical of $\PP_\Theta^-$ and taking a corresponding Iwahori factorization $J = \underbrace{ (J \cap U)}_{J^+} \underbrace{ (J \cap P^-) }_{J^-}$, we have $x_b \cdot Ja J = x \cdot b J^{+} J^{-} a J =  x \cdot b a (a^{-1} J^- a) J$, and  we have $a^{-1} J^- a \subset J^{-}$ because $|\delta(a^{-1}) | \leq 1$
for all negative roots for $\PP_\Theta^-$. 

 Equivalently, if we denote by $h_{a^{-1}}$ the adjoint of $h_a$ (which is equal to the characteristic measure of $J\check\lambda(a^{-1})J$) then:  \begin{equation}\label{Heckeaction}
 h_{a^{-1}} \star 1_{x_b J} = 1_{x_{ab} J}.
\end{equation}
Let us denote by $\mathcal H'$ the algebra spanned by the elements $h_{a^{-1}}$. Obviously, the elements $h_{a^{-1}}$ also satisfy the analogous statement of the previous lemma, i.e.\ $h_{a^{-1}} \star h_{b^{-1}} = h_{a^{-1}b^{-1}}$.

\begin{lemma}\label{expalongcochar} The exponential map is ``eventually equivariant'' in $\mathrm{fiber}(Z_1)$ with respect to the action of the algebra $\mathcal H'$; that is,
there exists a $J$-stable neighborhood $\mathcal N_{J}$ of ${Z_1}$ in $\mathrm{fiber}(Z_1)\cdot J$ such that, for any $x \in \mathcal{N}_J$ and any $|b| \leq 1$ we have: 
$$\exp_\Theta (h_{b^{-1}} \star x J) = h_{b^{-1}} \star \exp_\Theta(x J).$$

\end{lemma}

Of course, we have in our notation identified sets with their characteristic functions.

\begin{proof}

This follows from the facts:
\begin{enumerate}
\item[(i)] $\mathcal H'$ is finitely generated; 
\item[(ii)] the eventual equivariance of Proposition \ref{expproperties};
\item[(iii)] For each $x \in \mathrm{fiber}(Z_1)$ and every $\varepsilon>0$ the linear span of the characteristic functions of the sets $x_a J$ with $a\le\varepsilon$ is $\mathcal H'$-stable -- cf.\ (\ref{Heckeaction}).
\end{enumerate}
\end{proof}

As we have mentioned, the validity of the Proposition is independent of the choice of $\mathcal{N}_J$, and so we indeed take the neighborhood $\mathcal{N}_J$ so that the prior Lemma is valid.

\subsubsection{Inverting elements in the Hecke algebra} \label{subsub:invhecke}

Now recall Bernstein's ``stabilization theorem'' (see \cite[p 65]{BePadic}):  For any smooth representation $\pi$ and $a$ sufficiently close to zero (how close depends only on $J$), the action of $h_a$ on $\pi^{J}$ is stable, i.e.: 
\begin{equation}\label{stabilization}\pi^{J} = \ker(\pi(h_a)) \oplus \im (\pi(h_a)).
\end{equation}

This  stabilization theorem implies the generalization of ``Jacquet's lemma'' to the smooth case. Let $\NN^-$ denote the unipotent radical of $\PP_\Theta^-$, and denote by the subscript $_{N^-}$ the Jacquet module (coinvariants) of a representation with respect to $N^-$. The map   
$$\pi^J \to \pi_{N^-}^{J\cap P_\Theta^-}$$ 
which intertwines
the action of $h_b$  with the action\footnote{because we may write $h_b v = (J \cap N^{-}) b 
 \cdot b^{-1} (J \cap P_{\Theta}) b \cdot v$;
but $b^{-1} (J \cap P_{\Theta}) b$ fixes $v$ for sufficiently small $b$} of $\pi_{N^-}(b)$ --
has for kernel exactly $\ker(\pi(h_a))$, thereby inducing a bijection of $\im(\pi(h_a))$
onto $\pi_{N^-}^{J\cap P_\Theta^-}$.

The decomposition \eqref{stabilization} is independent of the choice of a sufficiently small $a$. It follows that the Hecke elements $h_b$ act invertibly on $\im(\pi(h_a))$ so long as $a,b$ are sufficiently small.  We extend this inverse to an operator $\tilde h_b$ on $\pi^J$ by defining $\tilde h_b$ to be zero on $\ker(\pi(h_a))$.

 Let us denote by $\mathfrak l_{J}$ the inverse of the induced bijection
 $\im(\pi(h_a)) \rightarrow \pi_{N^{-}}^{J \cap P_{\Theta}^-}$, the ``canonical lift'' to $\pi^{J}$. In these terms, we have:
\begin{equation} \label{canlift} \tilde h_b(f) = \mathfrak l_J(\pi_{N^-}(b^{-1})f_{N^-}), f\in \pi^J, \end{equation} 
 where $f_{N^{-}}$ is the image of $f$ in the Jacquet module. Let us note in particular that
 $\tilde h_b  \left(b f\right)$ is then $\mathfrak l_J f_{N^{-}}$.

\subsubsection{The proof of Proposition \ref{asymptoticstool}} 
{\bluetext We take the subgroups $J$ in our basis to admit Iwahori factorization with respect to $P_\Theta$ and $P_\Theta^-$,} and for a given $J$ choose $\mathcal{N}_J$ as in   Lemma 
\ref{expalongcochar}.  Take $x \in \mathcal{N}_J \cap \mathrm{fiber}(z)$; recall that $z$ is stabilized by $N^-$. 

Applying the remarks of \S \ref{subsub:invhecke} to the representation $\pi=C^\infty(X)$, we can define a functional:
\begin{equation}\label{limiteq}
\Lambda(f) = \lim_{|a|\to 0} (\tilde h_a f) (\exp_\Theta(x_a J)).
\end{equation}
 
This functional {\em a priori} depends on choices. Our goal is to show that it does not, and that it is in fact $G_x$-invariant, defining by Frobenius reciprocity the morphism: $e_\Theta^*: C^\infty(X)\to C^\infty(X_\Theta)$ that we are aiming at. Moreover, the resulting morphism $f \mapsto f_{\Theta}$
has the characterizing property \eqref{charprop} if we take
$\mathcal{N}_{ J} := a  \mathcal{N}_J$ where $a \in \mathcal{Z}(L_{\Theta})$
satisfies $|a| < \varepsilon'$, $\varepsilon'$ as in (3) below, and $\mathcal{N}_J$ was the neighborhood
from Lemma \ref{expalongcochar}. 

All this follows from the numbered statements following. In what follows, elements $a,b$
are always in $\mathcal{Z}(L_{\Theta})^+$. 

\begin{enumerate}
 \item \label{one} The limit (\ref{limiteq}) stabilizes for $|a|< \varepsilon$, where $\varepsilon>0$ depends only on $J$.

Indeed, for $a$ small as in Lemma \ref{expalongcochar} and any $a'$ with $|a'|\le 1$ we have:
$$ \tilde h_{aa'}  f(\exp_\Theta(x_{aa'} J)) = $$
$$\tilde h_{a'} \tilde h_{a} f (\exp_\Theta(h_{{a'}^{-1}}\star x_a J)) = $$
$$h_{a'} \star (\tilde h_{a'} \tilde h_a f)(\exp_\Theta(x_a J)) = $$
$$\tilde h_a f (\exp_\Theta(x_a J)),$$
the first step by \eqref{Heckeaction}, the second step due to Lemma \ref{expalongcochar}, the third because
$h_{a'}$ and $\tilde h_{a'}$ are inverse on the image of $h_a$.

\item The limit (\ref{limiteq}) is independent of $J$, and hence extends to a well-defined functional on $\pi=C^\infty(X)$.

Indeed, for $J'\subset J$ we clearly have $\im(\mathfrak l_J)\subset \im(\mathfrak l_{J'})$, and therefore $\mathfrak l_{J'}v = \mathfrak l_J v$ for every $v\in \pi_{N^-}^{J\cap P_\Theta^-}$.

 \item \label{three} 
For every $f\in \pi^J$ and $|a|<\varepsilon'$ (where $\varepsilon'>0$ depends only on $J$) we have: $$\Lambda(\pi(a) f) = f(\exp_\Theta(x_aJ)).$$ 

(This gives \eqref{charprop} for a suitable choice of neighborhood, as mentioned above.) 

Indeed, by the definition we have (for $|a|<\varepsilon$, where $\varepsilon$ is as in (\ref{one})): 
$$\Lambda(\pi(a) f) = (\proj_{\im (\mathfrak l_J)} f) (\exp_\Theta(x_a J)),$$
where the projection is with respect to the direct sum (\ref{stabilization}). Hence, if $f\in \im(\mathfrak l_J)$ then we are done. For general $f$, write $\bar{f} := \proj_{\im (\mathfrak l_J)} f$ for the projection of $f$ onto the image of the canonical lift; take $b$ small enough so that $h_b$ is stable, and take $a$ with $|a|<\varepsilon$. Then:
\begin{eqnarray*} \Lambda(\pi(ab) f) &=& \bar{f} (\exp_\Theta(x_{ab} J))  = 
  \bar{f}(h_{b^{-1}}\star\exp_\Theta(x_{a} J))  \\ &=& h_b\star \bar{f}(\exp_\Theta(x_{a} J)) = 
  (\proj_{\im (\mathfrak l_J)} h_b\star f) (\exp_\Theta(x_{a} J))  \\ &=& h_b\star f (\exp_\Theta(x_a J)) = f(\exp_\Theta(x_{ab} J)), \end{eqnarray*}
so the statement holds for $\varepsilon'= |b|\cdot \varepsilon$.

\item \label{four} Property (\ref{three}) characterizes $\Lambda$ among functionals on $\pi$ which factor through $\pi_{N^-}$. {\bluetext (The fact that $\Lambda$ factors through $\pi_{N^-}$ follows from \eqref{canlift}.)} 

Indeed, let us fix a $J$ in order to show that any $\Lambda'$ with the same property coincides with $\Lambda$ on $\pi^J$. For any $f$ in the image of $\pi^J$ and $ a \in Z(L_{\Theta})$ with $|a|\le \varepsilon'$ (where $\varepsilon'$ is such that property (\ref{three}) holds for both $\Lambda$ and $\Lambda'$) we have: 
$$\Lambda' (\pi(a)f)= f(\exp_\Theta(x_a J)) = \Lambda(\pi(a)f).$$
But, by assumption, $\Lambda'$ factors through $\pi_{N^-}$; since $\pi^J$ surjects onto $\pi_{N^-}^{P_\Theta^-\cap J}$ and $\pi_{N^-}(a)$ acts invertibly on $\pi_{N^-}^{P_\Theta^-\cap J}$, it follows that $\Lambda$ and $\Lambda'$ coincide there.

\item The functional $\Lambda$ is $G_x$-invariant. 

Indeed, for any $g\in G_x$ the functional $f\mapsto \Lambda(\pi(g) f)$ also satisfies (\ref{three}): We may assume that $g \in L_{\Theta}$ -- 
recall that we are supposing that the stabilizer of $x$ lies between $N^-$ and $P_{\Theta}^-$. 

Then:
$$\Lambda(\pi(g) \pi(a) f) = \Lambda(\pi(a) \pi(g) f) $$
(since $a\in \mathcal Z(L_\Theta)$)
$$= (\pi(g) f)(\exp_\Theta(x_a (gJg^{-1}))) \mbox{ for small } a$$
(notice that since $a$ belongs to the center of $L_\Theta$, $\pi(g) f$ is in the image of the canonical lift for $gJg^{-1}$ if $f$ is in the image of the canonical lift for $J$, hence the corresponding $\varepsilon'$ remains dependent only on the open compact subgroup for such elements)
$$ = f(\exp_\Theta( x_a gJ))$$
(by eventual equivariance; again, the implicit estimate is independent of the particular $f$)
$$ = f(\exp_\Theta( x_a J)).$$

\end{enumerate}

This proves Proposition \ref{asymptoticstool}, and hence Theorems \ref{mainasthm}, \ref{mainasthm2}.
 
{\bluetext We note that the uniqueness statement also follows directly for the following, which is proven like part (\ref{four}) in the above proof and is more general than the setting of spherical varieties. (It applies if we replace $X_\Theta$ by $N^-\backslash G$.)

\begin{lemma}
If $J$ is an open compact subgroup, and $N_\Theta$ a $J$-stable neighborhood of $\Theta$-infinity in $X_\Theta$, then the elements of $C_c^\infty(X_\Theta)^J$ which are supported in $N_\Theta$ generate $C_c^\infty(X_\Theta)^J$ under the Hecke algebra of $J$-biinvariant functions. 
\end{lemma}
 
\begin{proof} 
Let $V$ denote the quotient of $C_c^\infty(X_\Theta)$ by the $G$-subspace generated by the elements of $C_c^\infty(X_\Theta)^J$ which are supported on $N_\Theta$. If $V^J\ne 0$, then the smooth dual $\pi$ of $V$, a subspace of $C^\infty(X_\Theta)$, will also have a non-zero space of $J$-invariant elements, with the property that, as functions on $X_\Theta$, $f(a\cdot x) =0$ for any $x$ and for $a\in A_{X,\Theta}^+$ ``small'' enough.

If now $x$ is a point where $f(x)\ne 0$ for some $f\in \pi^J$, $N^-$ denotes the unipotent radical of the parabolic of type $P_\Theta^-$ contained in the stabilizer of $x$, and $\Lambda'$ denotes the functional ``evaluation at $x$'',  then on one hand $\Lambda'(f) \ne 0$ and on the other $\Lambda'(\pi(a) f') = f'(x_a) = 0$ for every  $f'\in \pi^J$ and $a \in \mathcal Z(L_\Theta)$ sufficiently ``small'' (depending only on $x$). On the other hand, $\Lambda'$ factors through the Jacquet module $\pi_{N^-}$, the map $\pi^J\to \pi_{N^-}^{J\cap P_\Theta^-}$ is surjective, as mentioned above (as a corollary of the stabilization theorem), and $a$ acts invertibly on $\pi_{N^-}^{J\cap P_\Theta^-}$, from which we get that $f_{N^-} = \pi_{N^-}(a) f'_{N^-}$ for some $f'\in \pi^J$, and hence
$$ \Lambda'(f) = \Lambda'(\pi(a) f') =0,$$
a contradiction.
\end{proof} 
 
} 

\subsection{Cartan decomposition and matrix coefficients} \label{ssCartan}
The above were just a reformulation of arguments due to Casselman and Bernstein, using the wavefront assumption for the variety $X$ which allowed us to ``push'' a point on $X$ to infinity using only ``anti-dominant'' cocharacters (with respect to a parabolic whose open orbit includes this point).

On the other hand, we can also present the above argument in a way where we  reduce everything to the (known) case of smooth matrix coefficients on the group. This argument first appeared in \cite{Lagier}, \cite{KatoTakano}. The reduction to smooth matrix coefficients is based on the following observation, for which we assume as fixed a Borel subgroup $\BB$ and we consider the ``universal Cartan'' $\AA_X$ of $\XX$ as a subvariety of $\XX$ as explained in \S \ref{invariants}. We also fix a maximal torus $\AA\subset\BB$ such that $\AA_X$ is an $\AA$-orbit, and denote by $A^+$ the anti-dominant elements of $A$. We denote by $A_X^+$ the elements of $A_X$ corresponding to the cone $\mathcal V$ of invariant valuations under the canonical map: $A_X^+\to \Lambda_X^+$; in the wavefront case, the map: $A^+/\AA(\mathfrak o) \to A_X^+/\AA_X(\mathfrak o)$ is surjective.

\begin{lemma}\label{coverlemma}
 For any spherical variety $X$, there is a compact subset $U\subset G$ such that $A_X^+  U= X$.
\end{lemma}

For symmetric spaces see the paper \cite{BO}. 

\begin{proof}  
{ For the purpose of this proof, we may replace $\XX$ by its quotient by $\mathcal{Z}(\XX)$.
After all, if we have found such a set $U$ which works for that quotient, 
then any point $x \in X$ differs from an element of $A_X^+ U$
by an element of $\mathcal{Z}(X)$, but $A_X^+$ is $\mathcal{Z}(X)$-invariant.

 Let $\bar X$ be a wonderful embedding of $X$, and let $\YY$ be the toric embedding of $\AA_X$ of the Local Structure Theorem \ref{localstructure}. 
  Being a toric variety under the split torus $\AA_X$, it admits a canonical structure over $\mathfrak{o}$  -- that is to say, the toric scheme
  defined by the same combinatorial data --  and $\YY(\mathfrak{o}) \cap A_X =  A_X^+$.
  
 We may assume that $U$ is open; then, by the Local Structure Theorem again, $\YY(\mathfrak o) U$ will be open in $\bar X$.     
 
 We now prove, for each $G$-orbit $Z \subset \bar{X}$, that there exists a compact open subset $U$ such that $(\YY(\mathfrak{o}) \cap Z) U = Z$. 
 In particular, taking $Z=X$, this implies the desired assertion. 
 We proceed by induction on the dimension of $Z$.  

The orbit $Z$  of minimal dimension is closed and thus compact. So the assertion is clear in this case. 
  
  Take a general orbit $Z$. By inductive assumption, we may assume that $U$ has been chosen so large
  that $\YY(\mathfrak{o}) U$ contains all orbits $Z'$ of lower dimension. Because $\YY(\mathfrak{o}) U$ is open,
  it contains an open neighborhood of all these $Z'$. But, because the closure of $Z$ is compact, 
  this means that $Z - \YY(\mathfrak{o}) U$ is compact. Enlarging $U$ appropriately, then, we may suppose that
  $\YY(\mathfrak{o}) U$ contains $Z$, as desired. }

    \end{proof}

In particular, if $X$ is wavefront, there are a finite subset $\{x_1,\dots, x_n\}\subset \mathring X$ and a compact subset $U\subset G$ such that:
$$ \cup_{i} x_i A^+ U = X.$$ 
Based on this, we can prove:

\begin{corollary}[The Wavefront Lemma] \label{wavefront}
Let $X$ be a wavefront spherical variety and $x_0,x_1,\dots$ representatives for the $G$-orbits on $X$. Let $o_i:G\to X$ be the corresponding orbit maps. There is a subset $G^+\subset G$ such that $\cup_i o_i(G^+)=X$ and with the property: For every open $K_1\subset G$ there is an open $K_2\subset G$ such that for every $i$ and every $g\in G^+$ we have $o_i(g)K_1\supset x_i K_2\cdot g$. 
\end{corollary}

To reformulate, assuming without loss of generality that $K_1$ and $K_2$ are subgroups:  If we take a double coset $K_2 g K_1$  with $g \in G^+$,
its image under $o_i$ consists of a single $K_1$-orbit. Informally, the orbit map does not ``smear out'' a double coset too much. Although in some ways inelegant, this is a very useful tool for reducing questions about $X$ to questions about $G$. 

 In particular, $o_i$ defines a map $K_2 \backslash K_2 G^+ K_1 / K_1 \rightarrow X/K_1$.  In
other words, if we enlarge $G^+$ to be  left $K_2$-invariant and right $K_1$-invariant,
the orbit map defines
$$o: K_2 \backslash G^+ / K_1 \longrightarrow X/K_1.$$

Such a result was used by Eskin and McMullen \cite{Eskin-McMullen} -- where $X$ is a symmetric variety under a semisimple real group $G$ -- in order to establish certain equidistribution and counting results.
These arguments partly inspired an early version of this section and we owe an intellectual debt to these authors. 

\begin{proof}[Proof of the corollary]
For notational simplicity, let us assume that there is only one $G$-orbit represented by $x_0\in \mathring X$, and that the map $A\to A_X$ (hence also the map $A^+\to A_X^+$, where $A^+$ denotes anti-dominant elements in $A$) is surjective. We can take the subset $G^+=A^+ U$ and then it is clearly enough to prove the lemma for $g\in A^+$; indeed, one is reduced to this case by replacing $K_1$ by the intersection of all of its $U$-conjugates (which is open since $U$ is compact). By shrinking $K_1$ further, if necessary, we may assume that it admits a decomposition: $K_1=N^- M^+$ where $N^-$ belongs to the unipotent radical of the parabolic opposite to $B$ and $M^+=K_1\cap B$ is completely decomposable: $M^+=(M^+\cap A)\cdot \prod_{\alpha>0} M_\alpha$, with the $M_\alpha$'s belonging to the corresponding root subspaces. 

For $g\in A^+$, we have $g^{-1} M^+ g\subset M^+$. Therefore, $x_0 g K_1\supset x_0 g M^+ \supset x_0 M^+ g$. Since $x_0\BB$ is Zariski open and hence $x_0 M^+$ is open in the Hausdorff topology, we can find a compact open subgroup $K_2$ of $G$ such that $x_0 K_2\subset x_0 M^+$.
\end{proof}

We also mention the following strengthening of Lemma \ref{coverlemma}, which has been proven under additional assumptions. It generalizes the Cartan and Iwasawa decompositions. In what follows, we assume smooth integral models over $\mathfrak o$ for the groups and varieties involved (and reductive, for $\GG$; in particular, $K=\GG(\mathfrak o)$ is a hyperspecial maximal compact subgroup), and we assume that the point $x_0$ used to define $\AA_X$ as a subset of $\XX$ in \S \ref{invariants} belongs to $\mathring\XX(\mathfrak o)$.

\begin{theorem}[Under additional assumptions]
The set $A_X^+\subset X$ contains a complete set of representatives for $K$-orbits on $X$; elements of $A_X^+$ which map to distinct elements of $\Lambda_X^+$ belong to different $K$-orbits.
\end{theorem}

This theorem was proven in \cite{SaRS} using an argument of \cite{GN}, under assumptions that are satisfied at almost every place, if $\XX$ is defined over a global field. The first part of the theorem    was proved in the symmetric case by Delorme and S\'echerre \cite{DS}; see also \cite{BO} for related results. The implication of Lemma \ref{coverlemma} is obvious, using the fact that wavefront varieties are precisely those for which $A_X^+$ can be covered by a finite number of orbits of the monoid $A^+$.

We are ready to apply the Wavefront Lemma in order to obtain the desired results on asymptotics:

\begin{proposition}\label{smoothfnls}
Let $X$ be a wavefront spherical variety and $x_0,x_1,\dots$ representatives for the $G$-orbits on $X$. Let $o_i:G\to X$ be the corresponding orbit maps, and assume that $K_1, K_2, G^+$ are as in Corollary \ref{wavefront}. Let $M:\pi\to C^\infty(X)$ be a morphism from a smooth representation of $G$, and denote by $L_i$ its composition with ``evaluation at $x_i$'', considered as a functional on $\pi$. Then for every $v\in \pi^{K_1}$, and for every point $x\in X$ with $x=x_i g$, $g\in G^+$ we have $M(v) (x) = \left<\pi(g) v, K_2* L_i\right>$, where we denote by $K_2*$ the operator ``convolution by the characteristic measure of $K_2$''.
\end{proposition}
In other words, the morphism $M$ is determined by the smooth functionals $K_2* L_i$.
Using known results about the asymptotics of matrix coefficients, then, 
we see that we can understand completely the asymptotics of $M(x)$. 
In principle, this could be used to give a second proof of Theorem \ref{mainasthm}; however,   
as pointed out to us by a referee, this is not a formality. Rather it 
requires a further analysis of how to choose $g$ when $x$ 
  is near a given ``wall'' of the compactification.  We will not carry this out here,
 and simply present the Proposition as an alternate way to understand asymptotics,
 albeit one that is less well adapted to the geometry of the spherical variety.

\begin{proof}
 Immediate corollary of the Wavefront Lemma \ref{wavefront}.
\end{proof}

Finally, we give the proof of Lemma \ref{supportLpsi}. We point the reader to \S \ref{ssWhittakertype} for the notation.  

\begin{proof}[Proof of Lemma \ref{supportLpsi}]
 Let us, for clarity, denote by $\YY$ the space $\XX$ without the character. {Since $\YY$ is parabolically induced from the variety $\XX^L$, its data $\Lambda_Y^+\subset \mathfrak a_Y^+$, $A_Y^+$, its spherical roots etc.\ are those of the variety $\XX^L$. 
 
We choose a point $x\in \XX^L$ with stabilizer $\MM\subset\LL$ and a torus $\AA$ as in \S \ref{notation} so that its quotient $\AA_Y$ can be identified with the $\AA$-orbit of the point of $x$. By Lemma \ref{coverlemma}, it is enough to show that for every $f\in C^\infty(X,\mathcal L_\Psi)$ the support of $f|_{A_Y^+}$ has compact closure in $\overline X$.

To understand what this support condition means, choose a smooth, complete toroidal embedding $\overline{\YY}$ of $\YY$ which contains $\overline{\XX}$, and apply the Local Structure Theorem \ref{localstructure}. The closure of $A_Y^+$ in $\overline Y$ is a compact subset of a smooth toric variety. To describe the subsets which are compact in $\overline X$, we use the ``valuation'' map: $A_Y^+\to \Lambda_Y^+$ -- notation as in \S \ref{invariants}. The subsets of $A_Y^+$ which have compact closure in $\overline X$ are precisely those whose valuations are contained in a finite number of translates of $\Lambda_X^+ = \Lambda_Y^+\cap \mathfrak a_X^+$.

Recall that the condition defining $\mathfrak a_X^+$ inside of $\mathfrak a_Y^+$ was determined by adding the simple roots of a parabolic opposite to $\PP^-$ to the spherical roots of $\YY$; thus, a sequence $\check\lambda_n$ of elements of $\Lambda_Y^+$ which does not belong to any finite number of translates of $\Lambda_X^+$ has the property that for some such root $\alpha$, 
\begin{equation}\label{Haight} \left<\alpha,\check\lambda_n\right>\to \infty.\end{equation}

Fix $f\in C^\infty(X,\mathcal L_\Psi)$ which violates the lemma, and let $\check \lambda_n$ be such a sequence of elements in the image of the support of $f$ under the valuation map. 
The $\GG$-stabilizer of a point on $A_Y\subset Y$ contains the unipotent radical of $\PP^-$.  we decompose its Lie algebra into root spaces for the torus $\AA$. The meaning of \eqref{Haight} is that for $u$ in the root space $U_{-\alpha}\simeq k$ in the unipotent radical of $P^-$:
$$ \lim_n \check\lambda_n(\varpi)^{-1} u \check\lambda_n(\varpi) = 0,$$
where $\varpi$ is any uniformizing element. (To be precise, the torus $\AA_Y$ does not act on $\UU_{P}$; but we saw in the proof of Proposition \ref{sphericalWhittaker} that it acts on $\UU_{-\alpha}$ through this quotient.)}

This implies that $f$ cannot be invariant by a compact open subgroup $U_0\subset U_{-\alpha}$, because then, for $u$ in that subgroup, denoting for simplicity $\check\lambda_n(\varpi)$ by $a_n$:
$$ f(a_n) = f(a_n u) = f((a_n u a_n^{-1}) a_n) = \Psi(a_n u a_n^{-1}) f(a_n),$$ and for large $n$ and some $u\in U_0$,  $\Psi(a_n u a_n^{-1})\ne 1$.
\end{proof}

\subsection{Mackey theory, the Radon transform and asymptotics} \label{ssRadon}

Thus far we have constructed a canonical ``asymptotics'' map
 $$ e_\Theta: C_c^{\infty}(X_\Theta) \longrightarrow C_c^\infty(X).$$ 
Here we will show how, in some instances, this map may be made more explicit.

The basic principle is as follows:  In order to write down any ``explication'', we need, first of all,
some common ``context'' to compare the varieties $X$ and $X_{\Theta}$. 
 Although the varieties $\mathbf{X}$ and $\mathbf{X}_{\Theta}$
look quite different, there is by Lemma \ref{importantobs} a {\em canonical identification} of the space of $P_{\Theta}$-horocycles for them. We shall prove that the adjoint asymptotics map commutes with integration along $\Theta$-horocycles.    

The operation of integration along horocycles is a classical concern of integral geometry, at least in certain other contexts, the so-called Radon transform. For example
in the case of the quotient of $\SL_2$ by its unipotent subgroup this reduces to the most classical case: integration of a function on an affine plane over lines.

\subsubsection{The Radon transform } \label{sssRadon}
Let $\Theta\subset\Delta_X$, and recall the horocycle space $\XX_{\Theta}^h$  defined in \S \ref{sshorocycles}.

Notice that $\UU_\Theta$ acts freely on $\mathring\XX \PP_\Theta$ and that, since $\XX$ is assumed quasiaffine, the orbits of any unipotent subgroup on $\XX$ are all closed. Therefore, we have a well-defined ``Radon transform'', defined by integration over generic $U_\Theta$-orbits:
\begin{equation} \label{radon} C^{\infty}_c(X) \stackrel{R_\Theta}{\longrightarrow} C^{\infty}(X_\Theta^h, \delta_\Theta).\end{equation}

 Here $C^{\infty}(X^h_{\Theta}, \delta_{\Theta})$ denotes
 the space of smooth sections of the complex line bundle $\delta_{\Theta}$
obtained thus: 
 Let $\mathcal{L}_{\Theta}$ be the algebraic line bundle over $\XX^h_{\Theta}$
 whose fiber at a point is the line of invariant volume forms on the corresponding
 unipotent group (recall there is a defining morphism $\XX^h_{\Theta} \rightarrow \PP_{\Theta} \backslash \GG$, which may be thought of as the variety of conjugates of $\UU_{\Theta}$;
 so to each point of $\XX^h_{\Theta}$ there is an associated unipotent group.) 
Now $\delta_{\Theta}$ is obtained from $\mathcal{L}_{\Theta}^{-1}$ via 
 reduction through $k^\times\xrightarrow{|\bullet|} \RR^\times_+ \hookrightarrow \CC^\times$.

Since we twist the action of $G$ on $C_c^\infty(X)$ by $\sqrt{\eta}$, where $\eta$ denotes the eigencharacter of the chosen measure on $X$, we will do the same for $C^\infty(X_\Theta^h,\delta_\Theta)$, in order for $R$ to be equivariant.

\subsubsection{Asymptotics via second adjunction}

Fix a ``standard'' parabolic $\PP_\Theta$ and consider the ``Levi quotient'':
\begin{equation}\label{Leviquotient}
 \XX_\Theta^L = \mathring\XX\cdot\PP_\Theta / \UU_\Theta.
\end{equation}
Recall that in the wavefront case this is canonically isomorphic to the Levi variety defined in \S \ref{levivarieties}, which is why we are using the same notation. We will assume that $\XX$ is of wavefront type from now on, so that our results on asymptotics hold, and will show how these results translate in terms of Mackey theory and the Radon transform. For simplicity, we denote the Jacquet module of any representation $\pi$ with respect to $U_\Theta$ by $\pi_\Theta$.

The Radon transform previously defined gives an identification:
\begin{equation} \label{blahblah} C^{\infty}_c(\mathring{X} P_{\Theta})_\Theta \stackrel{\sim}{\rightarrow} C^{\infty}_c(X_\Theta^L, \delta_{\Theta})\otimes \delta_\Theta^{-1}\end{equation}
Here (owing to the normalizing factor in the definition of Jacquet module) we twist the action of $L_\Theta$ on functions on $X_\Theta^L$ in such a way that $L^2(X_\Theta)$ is unitarily induced from $L^2(X_\Theta^L)$, as we are about to explain:

On the space $X_\Theta^L = \mathring X P_\Theta/ U_\Theta$ the measure on $X$ gives rise to an $L_\Theta$-eigenmeasure for which the following is true:
$$\int_{\mathring X P_\Theta} f(x) dx = \int_{X_\Theta^L} \int_{U_\Theta} f(ux) du dx.$$
This depends on the choice of Haar measure on $U_\Theta$. The character by which $L_\Theta$ acts on this measure is $\delta_\Theta \eta$ (recall that $\eta$ is the eigencharacter of the measure on $X$). Thus, we need to twist the unnormalized action of $L_\Theta$ on functions by $(\eta\delta_\Theta)^{\frac{1}{2}}(l)$ in order to obtain a unitary representation.

Another way to describe this twisting is the following: if we identify $X_\Theta^L$ as a subvariety of $X_\Theta$ as before, and $g\in P_\Theta^-$ with image $l\in L_\Theta^-$, then for a function $f$ a function on $X_\Theta$ we have:
\begin{equation} \label{Ltwist}  l\cdot (f|_{X_\Theta^L}) := \delta_\Theta^{\frac{1}{2}}(l) (g\cdot f)|_{X_\Theta^L}.\end{equation} 
(The twisting by $\sqrt\eta$ is already contained in the $G$-action on $X_\Theta$.)

We leave it to the reader to check that there is a choice of invariant measure, valued in the line bundle defined by $\delta_\Theta$, over $P_\Theta^-\backslash G$ such that:
$$L^2(X_\Theta) = I_{P_\Theta^-}^G \left(L^2(X_\Theta^L)\right),$$
where $I_{P_\Theta^-}^G$ denotes unitary induction with respect to that measure.

On the other hand, when we consider $\XX_\Theta^L$ as a subvariety of $\XX_\Theta^h$ and taking into account the twisting of the action we have by restriction of sections a map:
$$C_c^\infty(X_\Theta,\delta_\Theta)\twoheadrightarrow C^\infty(X_\Theta^L,\delta_\Theta)\otimes \delta_\Theta^{-1},$$ where $C^\infty(X_\Theta^L,\delta_\Theta)$ denotes smooth sections of the restriction of the above line bundle to $X_\Theta^L$. We have a (non-canonical) isomorphism of $L_\Theta$-representations: $$C^\infty(X_\Theta^L,\delta_\Theta)\simeq C^\infty(X_\Theta^L)\otimes \delta_\Theta.$$ We will denote the 
representation $C^{\infty}_c(X_\Theta^L, \delta_{\Theta})\otimes \delta_\Theta^{-1}$ which appears in (\ref{blahblah}) by $C^{\infty}_c(X_{\Theta}^L)'$; it is non-canonically isomorphic to $C^{\infty}_c(X_{\Theta}^L)$.

Hence we get an embedding:
\begin{equation}\label{Mackey}
 m_\Theta: C^{\infty}_c(X_{\Theta}^L)' \hookrightarrow C_c^\infty(X)_\Theta.
\end{equation}

The same considerations for $X_\Theta$ give an embedding, to be denoted by the same symbol:
\begin{equation}\label{MackeyTheta}
 m_\Theta: C_c^\infty(X_\Theta^L)'  \hookrightarrow C_c^\infty(X_\Theta)_\Theta.
\end{equation}
Recall that the quotients $\mathring X P_\Theta/U_\Theta$ and $\mathring X_\Theta P_\Theta/U_\Theta$ are both canonically isomorphic to $X_\Theta^L$ by \eqref{importantobs-ID}.

The analysis of the Jacquet module of $C_c^\infty(X)$ in terms of $P_\Theta$-orbits is usually called ``Mackey theory'' or ``the geometric lemma'', and therefore we will call the embeddings $m_\Theta$ ``Mackey embeddings''.

The following result {\em in principle} identifies the asymptotics (but is rather hard to use in practice):

 \begin{proposition}\label{propositionMackey}
 \begin{enumerate} 
 \item
Any $G$-morphism $M: C^{\infty}_c(X_{\Theta}) \rightarrow \pi$
is uniquely determined by its ``Mackey restriction'', i.e.\ by the induced map: 
\begin{equation}\label{Mackeyrestriction}
C^{\infty}_c(X^L_{\Theta})'  \rightarrow  \pi_\Theta
\end{equation}
obtained as the composition of the Mackey embedding \eqref{MackeyTheta} with the map of Jacquet modules induced by $M$.
In other words, the Mackey restriction map is injective (in fact, bijective):
$$\Hom_G(C_c^\infty(X_\Theta),\pi)\to \Hom_{L_\Theta}(C_c^\infty(X_\Theta^L)',\pi_U).$$
 \item
The diagram
\begin{equation}\label{diagrams1} 
\xymatrix{
& C_c^\infty(X_\Theta)_\Theta \ar[dd]_{e_\Theta}  \\
C^{\infty}_c(X^L_{\Theta})' \ar[ur]^{m_\Theta} \ar[dr]^{m_\Theta}\\
& C_c^\infty(X)_\Theta
}
\end{equation}
commutes, where the slanted maps are as above. 
\end{enumerate}
  \end{proposition}
 
\begin{proof}
Given that $C_c^\infty(X_\Theta) = I_{P_\Theta^-}^G C_c^\infty(X_\Theta^L)$, the first statement is precisely the statement of the second adjunction of Bernstein \cite[p61]{BePadic}: for any smooth representations $\sigma$ of $L_\Theta$ and $\pi$ of $G$, restriction with respect to the Mackey embedding $\sigma\hookrightarrow I_{P_\Theta^-}^G(\sigma)_\Theta$ gives rise to a bijection:
$$\Hom_G(I_{P_\Theta^-}^G(\sigma),\pi)\to \Hom_{L_\Theta}(\sigma,\pi_\Theta).$$
 
 To check commutativity, note that 
the space $C_c^\infty(X_\Theta^L)'$ is generated over $L_\Theta$ by the images of functions $f\in C_c^\infty(\mathring X P_\Theta)$ supported ``close enough to $\Theta$-infinity'', i.e.\ in a ``good'' neighborhood with respect to their stabilizer. (Indeed, any given compactly supported function on $X_\Theta^L$ can be translated by the center of $L_\Theta$ in the desired direction.) 
 It suffices to check commutativity on the images of those elements, which follows from the compatibility (in the sense of Lemma \ref{importantobs}) of the isomorphism: $\mathring XP_\Theta/U_\Theta\xrightarrow{\sim} \mathring XP_\Theta/U_\Theta$ with the exponential map.
\end{proof}

\subsubsection{Asymptotics and the Radon transform}

The dual asymptotics map
$$e_\Theta^*:C^\infty(X)\to C^\infty(X_\Theta),$$ when restricted to compactly supported functions, does not, in general, preserve compact support. However, elements in its image are compactly supported along unipotent orbits:

\begin{proposition}\label{affinesupport}
There is an affine equivariant embedding $\XX_\Theta\hookrightarrow \YY$ such that for every $\Phi\in C_c^\infty(X)$ the support of $e_\Theta^*\Phi$ has compact closure in $Y$.
\end{proposition}

This is a result of Bezrukavnikov and Kazhdan in the group case \cite[Prop 7.1]{BK}, and we prove it for the general case in the next subsection extending their argument. Given that orbits of unipotent groups on affine varieties are closed, it follows that the support of $e_\Theta^*\Phi$ intersects each unipotent orbit on a compact set, and hence Radon transform converges absolutely on $e_\Theta^*(C_c^\infty(X))$. Using this, we will explicate here the dual asymptotics map:

 \begin{proposition}\label{propositionRadon}
The square: 
\begin{equation}\label{Radoncommutes}
\begin{CD} 
C^{\infty}_c(X) @>{R_\Theta}>> C^{\infty}(X^h_\Theta,\delta_\Theta)\\
 @VV{e_\Theta^*}V  @VV{=}V \\
C^{\infty}(X_{\Theta}) \supset e_\Theta^*\left(C^{\infty}_c(X)\right) @>{R_\Theta}>> C^{\infty}(X_{\Theta}^h,\delta_\Theta) 
\end{CD}
\end{equation}
is commutative.  Here the right vertical arrow is induced from the canonical identification \eqref{importantobs-ID} of $P_{\Theta}$-horocycles on $X$ and $X_{\Theta}$. 
\end{proposition}

In principle, one can invert $R$ on the {\em parabolically induced} spherical variety $X_{\Theta}$, by using the theory of intertwining operators. We will do this in \S \ref{ssexplicitsmooth}.

\begin{proof} 
 
To avoid being overwhelmed by $\Theta$-subscripts, we write $\mathbf{Q} = \mathbf{P}_{\Theta}$, with Levi factorization $\mathbf{Q} = \LL \UU$;
let $\UU^-$ be opposite to $\UU$.  Note that $\LL = \LL_{\Theta}$ with prior notation.
 
Fix a compact subgroup $J$ that admits an Iwahori factorization $J = J_{U^-} J_L J_{U}$; we have written $J_S$ for $J \cap S$ whenever $S$ is any subgroup of $G$.   We verify the commutativity at the level of $J$-invariants. 
 
 Note that, for any $x \in \mathring{X}$, and any $a \in \mathcal{Z}(L_{\Theta})$ that is sufficiently large (with respect to the roots of $\PP_\Theta$),  we have $$ x (a J_{U^-} a^{-1}) \subset x J_Q \Rightarrow x_a J  = x a J_{U^-} J_Q = x a J_{U^-} a^{-1} a J_Q $$ $$\subset x J_Q a J_Q \subset x_a a^{-1} J_Q a J_Q\subset x_a J_Q,$$
 where we have written $x_a := xa$ for short. 
 From this we deduce that  $\int_{U} f(x_a k u) du$ is independent of $k \in J$, whenever $f$ is itself $J$-invariant.  In particular,   
 $$\int_{U} f(x_a u) = \int_{U} ([JuJ] \star f) (x_a) du,$$
 whenever  $a$ is sufficiently large (how large depends on $x$ and $J$); here
 $[J u J]$ is the $J$-bi-invariant measure of total mass $1$ supported on $JuJ$. The same conclusion holds if $x\in \mathring X_\Theta$ and $f$ is a $J$-invariant function on $X_\Theta$.
  
 Now let us compare the Radon transforms of $f \in C_c^\infty(X)$ and $f_\Theta := e_\Theta^*(f) \in C_R^\infty(X_{\Theta})$. 
 Fix a $J$-good neighborhood $N_\Theta$ of $\Theta$-infinity in $X$, let $x \in N_{\Theta} \cap \mathring{X}$ and let $x' J$ be the corresponding $J$-orbit
on $X_{\Theta}$, i.e.\ the image under the $\exp$-map.    Then,  if $N_\Theta$ is taken sufficiently small the orbits $xa J$ and $x' aJ$ are also matching under the $\exp$-map for $a \in \mathcal{Z}(L_{\Theta})$
sufficiently large. Indeed, it is enough to show this when $J$ is replaced by any smaller subgroup $J'$ with Iwahori factorization, e.g.\ a subgroup $J'$ such that $xJ' = xJ'_Q$; then $xa J'= x J' a J'$, and the claim follows from the eventual equivariance (statement (2)) of Proposition \ref{expproperties} and the fact that the characteristic measures of the sets $J'aJ'$, with $a$ anti-dominant, form a finitely generated subalgebra of the Hecke algebra (so by choosing a finite number of generators, we can find a neighborhood $N_\Theta$ where the exp-map is  equivariant with respect to that subalgebra).

According to what we noted above,  for $a \in \mathcal{Z}(L_{\Theta})$ sufficiently large we have:
$$ \int_{U} f_\Theta (x'_a u) du = \int_{U}  ([JuJ]\star f_\Theta)(x'_a) du = \int_{U}  ([JuJ]\star f)(x_a) du = \int_U f(x_a u) du,$$
the middle step because of the equivariance of the asymptotics map.
  
Write for $\mathscr{H}(G,J)$ the Hecke algebra of $G$ with respect to $J$. 
We have seen that the image of $C_c^\infty(X)^J$ under $e_{\Theta}^* R-R e_{\Theta}^*$ is a $\mathscr{H}(G,J)$-invariant subspace $W \subset C^\infty(X_{\Theta}^h,\delta_\Theta)^J$ with the following property: for any $y\in X_{\Theta,Q}$ there exists a sufficiently positive element $z$ of $A_{X,\Theta}$ (notice that the action of $\AA_{X,\Theta}\simeq \mathcal Z(\XX_\Theta)$ on $\XX_\Theta$ gives rise to a $\GG$-commuting action on $\XX_{\Theta}^h = (\XX_{\Theta})^h_{\Q}$, as well) such that $f(z\cdot y)=0$ for every $f\in W$. 
 
As is the case for $X_\Theta$, though, the characteristic functions of the sets $z \cdot y J$, where $y$ varies over $X_{\Theta}^h$ and $z\in A_{X,\Theta}$ is a sufficiently positive element (where the notion of ``sufficiently positive'' depends on $y$) generate $C_c^\infty(X_{\Theta}^h,\delta_\Theta)^J$ over $\mathscr H(G,J)$. Thus, any $f\in C^\infty (X_{\Theta}^h,\delta_\Theta)^J$ is given by the analogous formula of (\ref{limiteq}): $f(y) = \lim_{z\to 0} (\tilde h_z f)(z\cdot y)$ for every $f\in C^\infty (X_{\Theta}^h,\delta_\Theta)^J$.

Therefore, $W=0$.
\end{proof}

\subsection{Support of elements in $e_\Theta^*(C_c^\infty(X))$}\label{sssupport}

We now prove Proposition \ref{affinesupport}.

Let $\varchi(\XX)^+$ denote the weights of \emph{regular} Borel eigenfunctions on $\XX$; in other words, $\varchi(\XX)^+$ is the set of highest weights in the decomposition of $k[\XX]$ into a multiplicity-free direct sum of highest weight modules. This monoid induces a partial order $\succeq$ on the set $\mathfrak a_X = \mathfrak a_X = \varchi(\XX)^*\otimes \QQ$ by: $\mu\succeq\lambda \iff \left<\mu,\chi\right>\ge \left<\lambda,\chi\right>$ for all $\chi\in \varchi(\XX)^+$. We will extend this order to $\mathfrak a=\varchi(\BB)^*\otimes\QQ$ by pull-back ($\mu\succeq\lambda\iff \bar\mu\succeq\bar\lambda$, where $\bar\mu,\bar\lambda$ denote the images in $\mathfrak a_X$), and to the tori $A,A_X$ via the natural maps (which we denote by ``$\log$''):
\begin{eqnarray}\label{logmaps}
 \log: A & \to  &\varchi(\BB)^* \subset \mathfrak a, \\
  \log: A_X & \to &\varchi(\XX)^* \subset \mathfrak a_X. \nonumber
\end{eqnarray}
Since $\varchi(\XX)^+$ is contained in the set of dominant weights of $\varchi(\BB)$, the set of elements $\succeq 0 $ in $\mathfrak a$ contains the coroot cone. Finally, we will use, as we have done so far, the notation $\mathfrak a^+, A^+$ (resp.\ $\mathring{\mathfrak a}^+, \mathring A^+$ for anti-dominant (resp.\ strictly anti-dominant) elements, 
and similarly for $\mathfrak a_X$, $\mathfrak a_{X,\Theta}$, etc.\ (with respect to the root system of $\XX$). We are still in the wavefront case, so we have a surjection: $\mathfrak a^+\twoheadrightarrow\mathfrak a_X^+$.

Choose a splitting of the exact sequence of groups: 
\begin{equation}\label{charsequence}
  1 \to k^\times \to k(\XX)^{(\BB)}\to \varchi(\XX)\to 1,
\end{equation}
denoted by $\varchi(\XX)\ni \chi\mapsto f_\chi \in k(\XX)^{(\BB)}$. (Such a splitting exists because $\varchi(\XX)$ is free.)

Having fixed a maximal torus $\AA\subset \BB$, we choose a special maximal compact subgroup $K\subset G$ stabilizing a special point in the apartment of $A$ in the building of $G$. Hence, $K$ acts transitively on the points of the flag variety of $G$ and satisfies the Cartan decomposition: $G = K\Lambda^+ K$, where $\Lambda^+$, denotes the \emph{anti-dominant} elements in $\Lambda := \varchi(\BB)^*$, considered as elements of $G$ via any choice of uniformizer in $F$, which leads to a splitting of the $\log$ map \eqref{logmaps}.

Define, for every $\chi\in \varchi(\XX)^+$, a $K$-invariant ``norm'' on $X$:
$$ \Vert x \Vert_\chi := \max_{k\in K} |f_\chi(xk)|.$$

\begin{lemma}
 The norms $\Vert x \Vert _\chi$ have the following properties:
 \begin{enumerate}
  \item $\Vert x\Vert_{\chi_1+\chi_2} \le \Vert x\Vert_{\chi_1} \cdot \Vert x \Vert_{\chi_2}$;
  \item For any $a\in A^+$, $\Vert xa\Vert_\chi \le |\chi(a)| \cdot \Vert x\Vert_\chi$.
 \end{enumerate}
\end{lemma}

\begin{proof}
 The first property is obvious. The second follows by viewing the points of $X$, via the evaluation maps, as elements in the dual $V_\chi^*$ of the highest weight module $V_\chi$ spanned by $f_\chi$. Then $V_\chi^*$ can be realized as sections of the line bundle over the flag variety of $\GG$ induced from the character $\chi$ 
 of $\BB$, that is: regular functions on $\GG$ satisfying $F(bg) = \chi(b) F(g)$ for $b\in \BB$, and for such a section $F$ the norm we previously defined is equivalent to the norm: $\Vert F\Vert = \max_{k\in K} |F(k)|$. From this realization it is easy to see that $\Vert a\cdot F\Vert \le |\chi(a)| \Vert F\Vert$ for $a\in B$ dominant, where by $a\cdot F$ we denote the right regular action of $a$ on $F$: This follows from the Bruhat-Tits comparison (\cite[10.2.1]{HainesRostami}, cf.  \cite[p. 148]{CartierP}) of the Iwasawa and 
 Cartan decompositions: $K a_\lambda K \subset \cup_{\mu\ge \lambda} Ua_\mu K$, where $\mu\ge \lambda$ means that $\lambda-\mu$ belongs to the coroot cone (and recall that every $\chi\in \varchi(\XX)^+$ is $\ge 0$ on the coroot cone). An equivalent way to state that is that the norms defined extend to $K$-invariant norms on $V_\chi^*$, and the spectral norm of $a\in B$ is $|\chi(a)|$.
\end{proof}

For every $\mu\in \mathfrak a_X$, let $X_{\succeq \mu}$ denote the set of elements $x\in X$ with $\Vert x\Vert_\chi \le q^{-\left<\chi,\mu\right>}$ for all $\chi\in \varchi(\XX)^+$; obviously, by the first statement of the previous lemma, only a finite number of $\chi$'s is needed to define this set.
 These sets have the following properties:
\begin{enumerate}
\item $X_{\succeq \mu}$ is a set with compact closure in the points of the affine closure $\overline{\XX}^\aff:= \spec k[\XX]$; indeed, the $K$-translates of a generating set of $\BB$-eigenfunctions generate the coordinate ring $k[\XX]$ (since $BK=G$), and therefore any finite set of generators of $k[\XX]$ has bounded evaluations on elements of $X_{\succeq \mu}$.
\item The sets $X_{\succeq \mu}$ define a filtration of $X$ by $K$-invariant sets, decreasing with respect to the $\succeq$ ordering on $\mathfrak a_X$ (i.e.\ $X_{\succeq \mu}\subset X_{\succeq \lambda}$ if $\lambda \preceq \mu$). 
\end{enumerate}

\begin{lemma} \label{boundedaction}
 For any $\lambda\in \Lambda^+$, $\mu\in \mathfrak a_X$, we have: $$X_{\succeq \mu} \cdot K a_\lambda K \subset X_{\succeq \mu+\bar\lambda},$$ where $\bar\lambda$ is the image of $\lambda$ in $\mathfrak a_X$.
\end{lemma}

\begin{proof}
 In fact, since the sets are $K$-invariant, this reduces to the statement: $X_{\succeq \mu} \cdot  a_\lambda \subset X_{\succeq \mu+\bar\lambda}$. This follows immediately from the second statement of the previous lemma.
\end{proof}

Now we compare these sets for $X$ and $X_\Theta$. First of all, we clarify that \emph{we will use for the definition of the sets $X_{\Theta,\succeq \mu}$ only the set of weights $\varchi(\XX)^+$}, despite the fact that $X_\Theta$ might have more regular functions; this corresponds to the affine embedding of $\XX_\Theta$ obtained from the ``affine degeneration'' of the affine closure of $\XX$, cf.\ \S \ref{ssdegen}. This embedding will play the role of $\YY$ in the proof of Proposition \ref{affinesupport}. 
Secondly, we choose a splitting of the sequence analogous to \eqref{charsequence} for $\XX_\Theta$ as follows: 

Let $\XX^a$ denote the affine closure of $\XX$ (i.e.\ $\XX^a=\spec k[\XX]$), and consider the affine degeneration $\mathscr X^a$ of \S \ref{ssdegen} over the base $\overline{\AA_{X,ss}}$ (same notation as in \S \ref{ssdegen}). Notice that $\mathscr X^a \times_{\overline{\AA_{X,ss}}}\AA_X\simeq \XX^a\times\AA_X$ canonically, so we can define the regular function $F_\chi:(x,a)\mapsto f_\chi(x) \chi(a)$ on it. Choose an affine embedding $\overline{\AA_X}$ of $\AA_X$ where the kernel $\AA_1$ of $\AA_X\to\AA_{X,ss}$ acts freely, and such that $\overline{\AA_X}/\AA_1 = \overline{\AA_{X,ss}}$. This corresponds to a lifting of cocharacters from $\AA_{X,ss}$ to $\AA_X$, and gives rise to a base change:
$$ \widetilde{\mathscr X^a} = \mathscr X^a \times_{\overline{\AA_{X,ss}}}\overline{\AA_X}.$$

It follows immediately from the definition of $\mathscr X^a$ that the function $F_\chi$ extends to a regular function on $\mathscr X^a$, and its restriction to $\XX_\Theta$ will be denoted by $f_\chi^\Theta$. By the compatibility between the normal bundle and the affine degenerations (Proposition \ref{degencompatibility}), it can be seen that the correspondence between $f_\chi$ and $f_\chi^\Theta$ also arises from the almost canonical identification of Borel orbits \eqref{Borbitident}.

In the case of $X_\Theta$ we obviously have:
\begin{lemma}\label{white}
For every $\mu\in \mathfrak a_X$ and every $a\in A_{X,\Theta}$ we have: 
 $$a\cdot X_{\Theta,\succeq\mu} = X_{\Theta,\succeq \mu+\log a}.$$
\end{lemma}

\begin{proof}
 Indeed, for every $\chi\in \varchi(\XX)^+$, the function $f_\chi^\Theta$ is $\chi$-equivariant with respect to the action of $\AA_{X,\Theta}$: $f_\chi^\Theta (a\cdot x) = \chi(a) f_\chi^\Theta(x)$.
\end{proof}

Now we discuss compatibility of these sets with the exponential map:

\begin{lemma}\label{normscoincide} 
For every $\mu\in \mathfrak a_X$, and a $K$-good neighborhood $N_\Theta$ of $\Theta$-infinity, and for all $\kappa$ sufficiently deep in $\mathfrak a_{X,\Theta}^+$, the sets $X_{\succeq \mu+\kappa}\cap N_\Theta$ and $X_{\Theta,\succeq \mu+\kappa}\cap N_\Theta$ coincide.
\end{lemma}

Recall that, by abuse of language, we are treating here $N_\Theta$ as a subset of both $X$ and $X_\Theta$, when we really mean that $N_\Theta/K$ is identified as a subset of both $X/K$ and $X_\Theta/K$.

\begin{proof}
Let $\overline{\XX}$ be any simple smooth toroidal embedding of $\XX$, and let $\ZZ\subset\overline{\ZZ}$ be its closed $\GG$-orbit. The corresponding normal bundle degeneration $\overline{\mathscr X}^n\to \GGm^I$ (discussed in \S \ref{ssdegen}) has the following property, essentially by construction: For every distinguished (cf.\ \S \ref{expdefinition}) $p$-adic analytic map $\phi$ from the $k$-points of the normal bundle $N_\ZZ{\overline{\XX}}$ to $\overline{X}$, every point $x\in N_\ZZ{\overline{\XX}}(k)$, and every ``strictly positive'' (i.e.\ positive on every coordinate) cocharacter  $\check\lambda:\GGm\to \GGm^I$, we have: 
\begin{equation}\label{battleofbackyard}\lim_{m\to 0} (\phi(\check\lambda(m^{-1})x), \check\lambda(m)) = x.\end{equation} 
(Recall that the algebrogeometric meaning of $\lim_{m\to 0}$ is that the map extends from $\GGm$ to $\GGa$.) While a result of type
\eqref{battleofbackyard} can be proved in the setting of a normal crossing divisor on a general variety, in our current setting it can be proved directly using the Local Structure Theorem \ref{localstructure} to reduce to the case when $\overline{\XX}$ 
is an affine space and $\overline{\ZZ}$ is an intersection of coordinate hyperplanes and then computing explicitly. 

The property of ``distinguished'' maps to preserve $\GG$-orbits is actually irrelevant for this, but if $\phi$ also preserves $\GG$-orbits then both sides of the limit will be in the open set denoted by $\mathscr X^n$ in \S \ref{ssdegen} if the right-hand side is, and by Proposition \ref{degencompatibility} they can be regarded as points on the affine degeneration $\mathscr X^a$.

Hence, given such a $\phi$, for every $x\in X_\Theta$, considered as a subvariety of $\mathscr X^a$, we have: 
\begin{equation}\label{limit1} x =  \lim_{a \in A_{X,\Theta}^+} (\phi(ax),a^{-1}).
\end{equation}

Now, viewing $\mathscr X^a$ as an affine $\GG\times\AA_X$-spherical variety, we can define the sets $\mathscr X^a_{\succeq \mu}$ as in the case of $\XX$ and $\XX_\Theta$. Here $\mu$ can be in the sum $\mathfrak a_X\oplus \mathfrak a_X$, but we will be interested in the antidiagonal of $\mathfrak a_X$ only, so we assume that $\mu\in\mathfrak a_X$. The definition is through the functions $F_\chi$ as above (and their multiples by characters of $\AA_X$ extending to the base $\overline{\AA_X}$), which specialize to both $f_\chi$ and $f^\Theta_\chi$. 

For any fixed $\mu\in\mathfrak a_X$, the set $\mathscr X^a_{\succeq \mu}$ is open and compact (as in the discussion prior to Lemma \ref{boundedaction}), it intersects $X_\Theta$ on the set $X_{\Theta,\succeq \mu}$, and it intersects $X\times \{a\}$ on the set $X_{\succeq \mu-\log a} \times\{a\}$.  This, together with \eqref{limit1}, means that for $\log a$ sufficiently deep in $\mathfrak a_{X,\Theta}^+$, the point:
$$(\phi(a x),a^{-1})$$
belongs to $X_{\succeq \mu+\log a} \times\{a^{-1}\}$ if and only if $x\in X_{\Theta,\succeq \mu}$ . 

By Lemma \ref{white}, the point $y=ax$ belongs to $X_{\Theta,\succeq \mu+\log a}$ if and only if $x\in X_{\Theta,\succeq \mu}$. We deduce that: 
$$\phi( X_{\Theta,\succeq \mu+\log a}) = X_{\succeq \mu+\log a}$$
as long as $\log a$ is sufficiently deep in $\mathfrak a_{X,\Theta}^+$. In particular, the intersections of the sets $X_{\Theta,\succeq \mu+\log a}/K$ and  $X_{\succeq \mu+\log a}/K$ with $N_\Theta/K$ (identified as a subset of both $X/K$ and $X_\Theta/K$ via $\phi$) coincide.

\end{proof}

We are now ready to prove Proposition \ref{affinesupport} on the support of elements of the form $e_\Theta^*\Phi$. As mentioned above, $\YY$ will be the affine embedding of $\XX_\Theta$ such that the set of highest weights of $k[\YY]$ is $\varchi(\XX)^+$. For every $\lambda\in \Lambda^+$, let $\mathcal H_{\ge \lambda}$ denote the set of elements of the (full) Hecke algebra of $G$ supported in the union of cosets $Ka_\mu K$ of the Cartan decomposition, where $\mu-\lambda$ belongs to the coroot cone.  As in \cite[Lemma 8.8]{BK}, one proves:

\begin{lemma}\label{Heckeelement}
Given an open compact subgroup $J$, there is a finite subset $S$ of $A_{X,\Theta}$ such that for all $f\in C_c^\infty(X_\Theta)^J$ and every $a\in A^+$ with image  $\bar a \in A_{X,\Theta}^+$ there is an $F\in C_c^\infty(X_\Theta)^J$ whose support lies in $\bar a S \supp f$, and a Hecke element $h \in \mathcal H_{\ge \log a'}^J$, such that $f = h\star F$. Here $\log a'$ denotes the dual weight of $\log a$ ( $=-w\log a$, where $w$ is the longest Weyl group element).
\end{lemma}

The argument is identical to that of \emph{loc.cit}.\ and we omit it. We will now prove a stronger and more precise statement than that of Proposition \ref{affinesupport}: Let $\Phi\in C^\infty(X)^J$ whose support lies in $X_{\succeq \mu}$ for some $\mu$ in $\mathfrak a_X$. If $S$ is as in the previous Lemma and $\lambda \in \mathfrak a_X$ is such that $\lambda + \log S\preceq \mu $ then we will prove that $e_\Theta^* \Phi$ is supported in $X_{\Theta,\succeq \lambda}$. 

Indeed, let $f\in C_c^\infty(X_\Theta)^J$ be supported in the complement of $X_{\Theta,\succeq\lambda}$.  
Choose a $J$-good neighborhood $N_\Theta$ of $\Theta$-infinity such that the sets $X_{\succeq \mu+\kappa}$ and $X_{\Theta,\succeq \mu+\kappa}$ coincide for $\kappa$ deep enough in $\mathfrak a_{X,\Theta}^+$, according to Lemma \ref{normscoincide}. Choose an element $a\in A^+$ with image $\bar a \in \mathring A_{X,\Theta}^+$ such that $\bar aS\supp(f)\subset N_\Theta $, where $S$ is as in Lemma \ref{Heckeelement}. According to that lemma, $f = h\star F$, where $F$ is supported in $\bar aS \supp(f)$ and $h \in \mathcal H^J_{\ge \log a'}$ . Then the support of $F$ does not meet $X_{\Theta, \succeq \mu+\log \bar a}$, while by Lemma \ref{boundedaction} the support of $h^\vee \star \Phi$ (where $h^\vee$ denotes the dual Hecke element of $h$, which belongs to $\mathcal H^J_{\ge \log a}$) is contained in $X_{\succeq \mu+\log \bar a}$. If $a$ has been chosen so that $\log\bar a$ is sufficiently deep in $\mathfrak a_{X,\Theta}^+$ (which we may assume), then $\supp e_\Theta F$ does not meet $X_{\succeq \mu+\log\bar a}$ 
by Lemma \ref{normscoincide}. Therefore:

$$ \left<e_\Theta^*\Phi,f\right> = \left< e_\Theta^* \Phi,h\star F\right> = \left< e_\Theta^* h^\vee \star \Phi, F\right> = 0.$$
 \qed

\section{Strongly tempered varieties} \label{sec:sttempered}

It is clear that, if $G$ is a compact group and $X$ a compact homogeneous $X$-space, a Plancherel formula for $L^2(X)$ is a formal consequence of a Plancherel formula for $L^2(G)$, together with an understanding of which representations are $\XX$-distinguished.
Indeed, this is so even if we suppose only that point stabilizers are compact. 

What is perhaps surprising is that a corresponding phenomenon -- the Plancherel measure for $X$
is determined by a Plancherel formula for $L^2(G)$ -- persists even when point stabilizers are noncompact, so long as they are ``not too big.'' 
 (As a reference for the Plancherel formula for $L^2(G)$ itself, for $G$ a $p$-adic group, see Waldspurger's paper \cite{WaPl}.)
 
We term the spherical varieties for which this is so {\em strongly tempered,} 
and discuss their general theory. 
As a consequence of our general discussion, we will prove a conjecture of Ichino and Ikeda,
as well as a conjecture of Lapid and Mao. (For the latter, we give a short proof of a Whittaker-Plancherel formula.)

\subsection{Abstract Plancherel decomposition} 
\label{Plancherelgeneralities}

A Plancherel formula for $L^2(X)$ is, by definition, an isomorphism $L^2(X)\cong\int_{\hat G} \mathcal H_\pi \mu(\pi)$ of unitary $G$-representations; here $\hat G$ denotes the unitary dual of $G$ and the Hilbert space $\mathcal H_\pi$ is $\pi$-isotypic, i.e.\ isomorphic to a direct sum (in our case, finite) of copies of the unitary representation $\pi$. 
We describe $\mu$ as being a Plancherel measure\footnote{In some treatments,
a Plancherel formula is described by a collection of measures $\nu_n$, for $n \in \{1,2, \dots,\} \cup \{\infty\}$, together with an isomorphism $L^2(X) \cong \sum_{n} \int \pi^{\oplus n} \nu_n(\pi)$. In this language, $\mu$ is in the same measure class as $\sum_{n} \nu_n$. }
 for $L^2(X)$; any other Plancherel measure $\mu'$ belongs to the same measure class as $\mu$.

We recall from (cf.\ \cite{BePl}) how to describe such a decomposition, and more generally any morphism from $L^2(X)$ to a direct integral of unitary representations.  
The subspace $C^{\infty}_c(X)$ of smooth, compactly supported, functions is countable-dimensional, 
so there is a family of morphisms $$L_\pi: C^{\infty}_c(X)\to \mathcal H_\pi$$ (defined for $\mu$-almost every $\pi$) such that $\alpha\mapsto L_\pi(\Phi)$ represents $\Phi$ for every $\Phi\in C^{\infty}_c(X)$.

By pull-back,  we obtain seminorms $\Vert\bullet\Vert_\pi$ on $\bruhat(X)$;    the spaces $\mathcal H_\pi$ can be identified with the completions of $\bruhat(X)$ with respect to the seminorms $\Vert\bullet\Vert_\pi$. In particular, the spaces $\mathcal H_\pi$ are completions of the spaces of \emph{$\pi$-coinvariants}:
\begin{equation} \label{coinvariants} \bruhat(X)_\pi := \left(\Hom_G(\bruhat(X), \pi)\right)^* \otimes \pi.
\end{equation}
Notice that there is a canonical quotient map: $\bruhat(X)\twoheadrightarrow \bruhat(X)_\pi$ (surjectivity follows from the irreducibility of $\pi$).

Therefore, by a \emph{Plancherel decomposition} of $L^2(X)$ (or a quotient thereof), we will mean the following set of data: a positive  measure\footnote{More precisely: a positive Borel measure, with respect to the standard Borel structure on $\hat{G}$; see \cite[Prop. 4.6.1]{Dixmier}.} $\mu$ on $\hat G$; and a measurable set of invariant, non-zero seminorms $\Vert\bullet\Vert_\pi$ on the spaces $\bruhat(X)_\pi$, for $\mu$-almost every $\pi$, so that for every $\Phi \in \bruhat(X)$:
\begin{equation}  \label{bsn} \Vert\Phi\Vert^2= \int_{\hat G} \Vert\Phi\Vert_\pi^2 \mu(\pi). \end{equation}
(Here we denoted the image of $\Phi$ in $\bruhat(X)_\pi$ again by $\Phi$.)

The data $\mu,(\Vert\bullet\Vert_\pi)_\pi$ are uniquely determined up to the obvious operation of multiplying $\mu$ by a non-negative measurable function which is $\mu$-almost everywhere non-zero and dividing $(\Vert\bullet\Vert_\pi)_\pi$ by the square root of that function. This is the content of ``uniqueness of Plancherel decomposition.'' In fact, we shall need a slightly stronger uniqueness, even when we allow certain norms to be possibly negative:

\begin{proposition} \label{positivedense}  Suppose given a positive measure $\mu$ on $\hat{G}$ 
as well as a family $\pi \mapsto H_{\pi}$ of Hermitian forms  on $\bruhat(X)_{\pi}$ so that
$\Phi \mapsto H_{\pi}(\Phi)$ is $\mu$-measurable, for every $\Phi$ in $\bruhat(X)$, and moreover
$\|\Phi\|^2=\int H_{\pi}(\Phi) \mu(\pi) $ for all $\Phi$. Then $H_{\pi}$ are positive semidefinite for $\mu$-almost every $\pi$, so that $(\mu, H_{\pi})$ define a Plancherel formula. 
\end{proposition}
\begin{proof}
Let $(M, \sigma)$ be a Levi subgroup and a supercuspidal representation.  Let $Y'$ 
be the set of all unramified characters of $M$ modulo the finite subgroup of those for which $\sigma\otimes\chi\simeq \sigma$ (it is a complex algebraic variety). Let $W_M$ be the normalizer of $M$ in the Weyl group -- it acts on $Y'$, and we set $Y=Y'\sslash W_M = \spec \CC[Y']^{W_M}$.

The theory of the Bernstein center has the following consequence: For
a ``good'' basis of compact subgroups $J$ and any $\alpha \in \mathbb{C}[Y]$, 
there exists an element $f \in \mathcal H(G,J)$ (the Hecke algebra of $J$-biinvariant compactly supported measures on $G$) so that $f$ acts on 
the $J$-invariant space of any representation
$i_P^G(\sigma \cdot \chi)$ as scalar multiplication by $\alpha(\chi)$, and $f$ acts on every other Bernstein component as $0$.

Let $\hat{G}_0$ be the subset of $\hat{G}$ consisting of irreducible representations
that occur as a subquotient of some $i_P^G(\sigma \otimes \chi)$.  
Now $\hat{G}_0$ is closed and open in $\hat{G}$; 
  the measure $\mu$ induces a measure $\mu_0$ on $\hat{G}_0$, and it is sufficient 
to show that $H_\pi$ is positive semidefinite for $\mu_0$-almost every $\pi$, because the union of sets $\hat{G}_0$
as we vary $M$ is all of $\hat{G}$. 

If $\pi \in \hat{G}_0$ is a subquotient of $i_P^G(\sigma \otimes \chi)$
then $\chi$ is uniquely determined modulo $W_M$, i.e.\ the image of $\chi$ in $Y$ is uniquely determined. This gives a map $\hat{G}_0 \rightarrow Y$. Let $Y_0$
be the closure (in the usual topology on $Y$) of the image of $\hat{G}_0$. 
The   induced map $\pi:  \hat{G}_0 \rightarrow Y_0$ is a Borel map. 

Let $\bar{\mu}$ be the push-forward of the measure $\mu_0$
to $Y$.
The disintegration of measure implies that we may disintegrate the measure $  \mu(\pi)$ as an integral $\int_{y}  \mu_y \cdot  d\bar{\mu}(y) $, where $y \mapsto \mu_y$
is a measurable mapping from $Y_0$ to the space of measures on $\hat{G}_0$,
and each $\mu_y$ is entirely supported on $(\hat{G}_0)_y$, 
the fiber of the mapping above $y$. 
The space $(\hat{G}_0)_y$ is finite, and thus $\mu_y$ is nothing more than a function
on this finite set. 

 Fix $\Phi \in C^{\infty}_c$, and
let $$F^{\Phi}(y) = \int H_{\pi}(\Phi)d \mu_y  \left( = 
\sum_{\pi \in (\hat{G}_0)_y}  H_{\pi}(\Phi) \mu_y(\{\pi\}) \right) .$$ It is measurable on $Y_0$. 
 Then
for any $z \in \C[Y]$, the theory of the Bernstein center implies
$$\int_{y}  |z(y)|^2  F^{\Phi}(y) d\bar{\mu} \geq 0,$$
The Bernstein center induces a {\em dense} subalgebra of $C(Y_0)$; indeed
it separates points on $Y$ and it is closed under complex conjugation.
Therefore, $F^{\Phi}(y) \geq 0$ for $\bar{\mu}$-almost all $y$. 
The space $C^{\infty}_c$ being of countable dimension, this can be said simultaneously for all $\Phi$. 
That is to say, away from a set $S \subset Y_0$
with $\mu(S) = 0$, we have:
$$\sum_{\pi \in (\hat{G}_0)_y}  H_{\pi}(\Phi) \mu_y(\{\pi\}) \geq 0,$$
for all $\Phi \in C^{\infty}_c(X)$.

The left-hand side is a Hermitian form on the {\em finite length} $G$-representation
$\bigoplus_{\pi \in (\hat{G}_0)_y} C^{\infty}_c(X)_{\pi}$.  It follows that, whenever
$y \notin S$, we have  $H_{\pi}(\Phi) \geq 0$ 
for every $\pi$ in the fiber above $y$ with $\mu_y(\{\pi\}) >  0$.

Let $B$ be the set of $\pi \in \hat{G}_0$ for which $H_{\pi}$ fails to be positive semidefinite.
Then $B$ is measurable, since one can test the failure of positive semidefiniteness by a countable number of evaluations. But $\mu(B) = \int_{y} \mu_y(B) d\bar{\mu}(y)$. 
According to the discussion above, $\mu_y(B) = 0$ for $\mu$-almost all $y$,
so $\mu(B) = 0$, concluding the proof. 
   \end{proof}

To give invariant norms on $\bruhat(X)_{\pi}$ is equivalent to giving
an equivariant morphism 
$$M_\pi : \bruhat(X\times X)\to\pi\otimes\bar\pi$$
(where $\pi$ is assumed to have a unitary structure $\tilde\pi \xrightarrow{\sim}\bar\pi$).
 The hermitian forms associated to the above norms are the so-called \emph{spherical characters}:
$$\theta_\pi: \bruhat(X) \otimes \bruhat(X) = \bruhat(X\times X)\to\pi\otimes\bar\pi\to \CC$$
where the last arrow denotes the unitary pairing.  (This map is $G$-invariant.) Thus, to be explicit, 
these have the property that:
\begin{equation}\left<\Phi_1,\Phi_2\right>_{\mathcal H}= \int_{\hat G} \theta_\pi(\Phi_1\otimes\Phi_2) \mu(\pi). \end{equation}
Notice that, automatically, for $\mu$-almost every $\pi$ the spherical characters $\theta_\pi$ are positive semi-definite.

	\subsection{Definition; the canonical hermitian form}  \label{canhermform}

In this section we require, for simplicity, that stabilizers of points on $X$ are unimodular; in other words, $X$ admits an invariant measure. The modifications necessary to remove this assumption are straightforward.
We say that $\X$ is strongly tempered if, for any $x \in \X(k)$, the restriction of any ($G$-)tempered matrix coefficient to the stabilizer $H$ of $x$ in $G = \G(k)$ is in $L^1(H)$. 

In checking this, the following remark is useful: If $\pi_0$ is the normalized induction to $G$ of the trivial representation on a Borel subgroup $B$,  and $v_0 \in \pi_0$ the spherical vector (i.e., $K$-invariant for  $K$ a good maximal compact subgroup of $G$, satisfying the Iwasawa decomposition $G=BK$) 
then every tempered matrix coefficient $\varphi(g) $ is majorized by the spherical one (see \cite[Theorem 2]{CHH})
\begin{equation} \label{mcbound} |\varphi(g)| \leq c  \ \langle g v_0, v_0 \rangle .\end{equation}
(The right hand side is positive.) Moreover, if we fix an open compact subgroup $U$,
there exists a constant $c =c(U)$ which works whenever $\varphi(g) = \langle g u_1, u_2 \rangle$ arises
from $U$-invariant $u_1, u_2$, with $\|u_1\| = \|u_2\| = \|v_0\|$.

Let $(\pi, V)$ be a ($G$-)tempered representation of $G$, and assume that $X$ is a strongly tempered variety.
We define the morphism: $M_\pi: \pi\otimes\tilde\pi\to C^\infty(X\times X)$ characterized by the property\footnote{This property defines the morphism in the $G\times G$-orbit of the diagonal $\Delta X\subset X\times X$. We can extend it by zero on the whole space; in fact, the extension plays no role in what follows.} that 
$$M_\pi(v\otimes u)(x,x) = \int_{G_x} \left< \pi(h) v, u\right> dh.$$
We let $$\theta_\pi: C_c^\infty(X\times X)\to \tilde\pi\otimes\pi\to \CC$$ denote the adjoint composed with the canonical pairing.

Let $\mu_G$ denote the canonical Plancherel measure for $L^2(G)$, normalized as usual,  i.e.\  the spherical characters are simply the usual characters.  Since the spectrum of $G$ as a $G \times G$-representation is supported on representations of the form $\pi\otimes\tilde\pi$, $\mu_G$ will be thought of, as usual, as a measure on $\hat G$.

\begin{theorem} \label{ST:Plancherel} Suppose that $(\G, \X)$ is strongly tempered. 
Then $\theta_\pi$ and $\mu_G$ define a Plancherel formula for $L^2(X)$, 
in the sense that for $\Phi_1,\Phi_2\in C_c^\infty(X)$ we have:
\begin{equation}
 \left<\Phi_1,\Phi_2\right> = \int_{\hat G} \theta_\pi(\Phi_1\otimes\overline{\Phi_2}) \mu_G(\pi).
\end{equation}
In particular, $M_{\pi}(v,v)|_{\Delta X} \geq 0$ for every $\pi$, and $L^2(X)$ is tempered as a $G$-representation; its Plancherel measure is absolutely continuous
with respect to the group Plancherel measure.  \end{theorem}

\begin{remark}
 The positivity assertion states simply that $\int_{h \in G_x} \langle h v, v \rangle \geq 0$ -- 
an assertion that is obvious when $G_x$ is compact.  \end{remark}

\proof 
Let $\pi$ be a unitary representation (endowed with a invariant Hilbert norm, hence with a fixed isomorphism: $\tilde\pi\simeq\bar\pi$) and let $i_\pi:C_c^\infty(G)\to\bar\pi\otimes\pi$ denote the dual of matrix coefficient $m_\pi$.  (If we identify $\bar \pi \otimes \pi$ with a 
subspace of $\mathrm{End}(\pi)$, the morphism $i_{\pi}$ simply maps $f \in C_c^{\infty}(G)$
to $\pi(f)$.)

The Plancherel formula on $G$ can be written as: 
$$\left< f_1, f_2\right>_{L^2(G)} = \int \left< i_\pi(f_1), i_\pi(f_2)\right>_{HS} \mu_G(\pi)$$
where $\left< \,\, , \,\, \right>_{HS}$ denotes the Hilbert-Schmidt hermitian form on $\bar\pi\otimes\pi\subset \End(\pi)$.

We are going to assume, for simplicity, that $G$ has a single orbit on $X=H\backslash G$, but the general case follows in the identical fashion. Notice that the map:
$$C_c^\infty(G)\ni f\mapsto \Phi(x)=\int_H f(hx) dh\in C_c^\infty(X)$$ is surjective. Let $f$ and $\Phi$ be such, then:
\begin{equation} \label{ForRefAA1} \Vert \Phi \Vert_{L^2(X)}^2 = \int_G \int_H f(hg) \bar f(g) dh dg = \int_H \left< \mathcal L_{h^{-1}}(f), f \right>_{L^2(G)} dh\end{equation} 
where $\mathcal L_{\bullet}$ denotes the left regular representation of $G$.

We will use the following explication of $\theta_{\pi}$:  
\begin{equation} \label{ThetaExplication}  \theta_\pi (\Phi\otimes\bar \Phi)  = \int_H \left< \pi(h) i_\pi(f) , i_\pi(f) \right>_{HS} \end{equation} 
Indeed, suppose that $f$ is $J$-invariant,  and choose dual bases $v_1, \dots, v_n$ for $\pi^J$
and $v_1^*, \dots, v_n^*$ for $(\overline \pi)^J$. The definition of $\theta_{\pi}$
says that $\theta_{\pi}(\Phi \otimes \bar \Phi)$ is   the integral,
over $(g_1, g_2) \in X \times X$, of 
$\Phi(g_1) \overline{ \Phi}(g_2) , \int_{H} \langle \pi(h g_1) v_i,  \overline{ \pi}(g_2) v_i^*  \rangle$. 
Unfolding, this equals  
$$ \int_{G \times G} f(g_1) \overline{ f(g_2)} \int_{H} \langle \pi(h g_1) v_i, \overline{ \pi}(g_2) v_i^*  \rangle ,$$
which in turn equals 
  $ \int_{H}  \langle \pi(h) i_{\pi}(f)  v_i, i_{\overline{\pi}}(f) v_i^*  \rangle$
  as desired.

  Keeping in mind that $i_\pi(\mathcal L_{h^{-1} f})= \pi(h^{-1}) i_\pi(f)$ we get:
\begin{equation} \label{ForRefAA2} \Vert \Phi \Vert_{L^2(X)}^2 = \int_H \int_\pi \left< \pi(h) i_\pi(f) , i_\pi(f) \right>_{HS} \mu(\pi) dh = \end{equation}
$$ = \int_\pi \int_H \left< \pi(h) i_\pi(f) , i_\pi(f) \right>_{HS} \mu(\pi) dh = $$
$$ = \int_\pi \theta_\pi (\Phi\otimes\bar \Phi) \mu(\pi).$$

Notice that at all stages these integrals are {\em absolutely convergent}, justifying our
application of Fubini.  Indeed, because of \eqref{mcbound},
the integrand is bounded in absolute value by a constant multiple of $  \|i_{\pi}(f)\|^2_{HS} \cdot  \langle h v_0, v_0 \rangle$.  Then $\langle h v_0, v_0 \rangle$ is integrable over $H$ by assumption,
and $  \|i_{\pi}(f)\|^2_{HS}$ is $\mu$-integrable by the Plancherel theorem.

In particular, we have established the statement for $\Phi_1, \Phi_2 \in C^{\infty}_c(X)$. 
That the remaining statements follow is a consequence of Proposition \ref{positivedense}. 
\qed 

\subsection{The Whittaker case and the Lapid--Mao conjecture} \label{ssWhittaker}
A case of particular interest which does not literally fall under the strongly tempered is the ``Whittaker-Plancherel formula''. 
 
We shall now give a short reduction of this formula to the usual Plancherel formula.   The proof is largely the same as the previous, but the integrals are only conditionally convergent and we need to interpret them suitably; we therefore treat this case separately.

In fact, the treatment that follows covers many more cases of `Whittaker-induced'' models than the Whittaker model itself. We set up the notation very generally, but the key assumption is given in the paragraph below, and in practice is fulfilled only in ``nondegenerate'' cases. The general setup is as follows:   Let $\HH=\MM\ltimes\UU^-$ be a spherical subgroup, where $\MM$ is contained in a Levi subgroup $\LL$ of a parabolic $\PP^-$ and $\UU^-$ is the unipotent radical of this parabolic. Assume that $\Lambda_0:\HH\to \GGa$ is a homomorphism and let $\Psi=\Psi_{\Lambda_0}$ be the composition of $\Lambda_0$ with a nontrivial unitary character of $k$; let $\HH_0=\ker \Lambda_0$. 
Finally, assume that there is a cocharacter $\check\mu:\GGm\to\mathcal Z(\LL)$ which is nonpositive under the right adjoint action on $\mathfrak u^-$ (that is, its eigencharacters on $\mathfrak u^-$ are of the form $a\mapsto a^n$ with $n\le 0$), normalizes $\HH_0$, and acts nontrivially on the quotient $\HH/\HH_0$; identifying this quotient with $\Ga$, the action of $\check\mu$ is by a positive character: $x \cdot \check\mu(a) = a^{n_0}x$ with $n_0 > 0$. The crucial assumption is: 

\begin{quote}
 Let $v_s$ be a holomorphic section in the principal series $I_B^G(\delta_B^s)$ obtained by normalized induction from a power of the modular character of the Borel subgroup. The ``Jacquet integral'': 
 \begin{equation}\label{Jacquet integral}
 \int_H v_s(h) \Psi^{-1}(h) dh
\end{equation} 
 is convergent for $\Re s>0$, and extends continuously to $s=0$.
\end{quote}

This assumption, and in particular the extension of the integral to $s=0$, typically forces some nondegeneracy condition on $\Lambda_0$. One can calculate that the following examples satisfy it:
\begin{enumerate}
 \item $\HH$ is a maximal unipotent subgroup of $\GG$, $\Psi=$ a generic character. 
 
 \item $\GG= \GGSp_4$, $H=$ the Bessel/Novodvorsky subgroup $\TT\ltimes \UU$, where $\UU$ is the unipotent radical of the Siegel parabolic ($\simeq S^2 V$, where $V$ is a two-dimensional vector space), $\Lambda: \UU\to \GGa$ (i.e.\ $\Lambda \in S^2 V^\vee$) and $\TT \supset \SSO(\Lambda)$ is the stabilizer of $\Lambda$ in a Siegel Levi subgroup. 
 
 \item $\GG=\SSO_n\times \SSO_{n+2m+1}$, $\HH = \SSO_n \ltimes $the unipotent radical of the parabolic with Levi $\GGm^m\times \SSO_{n+1}$, $\Lambda$: a nondegenerate additive character of $\HH\cap \GGL_m$ (where $\GGL_m$ is in the Levi $\GGL_m\times \SSO_{n+1}$) This case has been considered by Waldspurger \cite[Chapter 5]{Waldspurger-GP2}. 
\end{enumerate}

For simplicity, we will only present the Whittaker case here, though with notation suggestive of the general case; the steps required for the other cases are completely analogous. Hence, we are in the setting of (1): $\HH$ is a maximal unipotent subgroup of $\GG$, $\Psi=$ a generic character, $\HH_0$ is the kernel of the corresponding algebraic morphism into $\GGa$.

\begin{proposition} \label{Whittakerconvergence}  The restriction of any tempered matrix coefficient of $G$ to $H_0$ belongs to $L^1 (H_0)$.
\end{proposition}

\begin{proof} 
We fix a Borel subgroup $\BB$ with $\check\mu(\GGm)\subset\BB$ so that $\HH\BB$ is open and $\HH \cap \BB$ is trivial. Let $\pi_s = I_B^G(\delta_B^s)$, the normalized principal series induced from the $s$-power of the modular character of $B$. The representation $\pi_s$ is tempered if $\Re s = 0$. By \eqref{mcbound}, it suffices to prove the proposition for matrix coefficients of $\pi_0$.

We fix the following invariant inner product on $\pi_s$ ($\Re s =0$):
\begin{equation}\label{innerproduct}\left< v_1, v_2\right> = \int_{H} v_1(u)\overline{v_2(u)} du,
\end{equation}
 where $ v_1,v_2\in \pi_s$, considered as functions on $G$. This integral converges; it is the restriction to the open $H$-orbit of the compact, $G$-invariant integral of $v_1(u)\overline{v_2(u)}$ over $B\backslash G$.
 
The underlying vector spaces of all representations $\pi_s$ can be identified with one another \emph{by restriction of functions to $K$}, a maximal compact subgroup satisfying the Iwasawa decomposition $G=BK$. In particular, $K$-invariant elements for all representations are identified in this common vector space -- let $v$ be a non-zero $K$-invariant element with $v|_K >0$, and $v_s$ its ``realization'' in $\pi_s$. (The inner product that we chose on $\pi_s$ is not compatible with this identification, but of course it varies continuously in $s$.) 

Let $f_s$ be the matrix coefficient $\left<\pi_s(g)v_s,v_s\right>$, for $\Re s=0$. To prove convergence of the integral $\int_{H_0} f(n) dn$ we extend $f_s$ to all $s$ with $\Re(s)\ge 0$ using expression (\ref{innerproduct}), that is:
$$f_s(g) = \int_{H} v_s(ug)\overline{v_s(u)} du.$$
 (Of course, for $\Re s\ne 0$ this does not represent a matrix coefficient.) {\bluetext  It will follow from the argument below that this expression converges for $\Re s \ge 0$, but for the moment we can restrict our attention to positive real $s$, and treat this as a possibly infinite expression.}
Notice that for $\Im s= 0$ we have $f_s(g)\ge 0$ for all $g$; in fact, for such $s$, we have $v_s(g)>0$ for all $g$. 

{\bluetext Thus, by expanding the definitions and using Fatou's lemma,
\begin{eqnarray}\label{limit2}\int_{H_0} f_0(n)dn =  \int_{H_0}\int_{H} v_0(un)\overline{v_0(u)} du dn \le  \nonumber \\  \le \limsup_{s\to 0^+, s\in \RR} \int_{H_0} \int_{H} v_s(un)\overline{v_s(u)} du = \limsup_{s\to 0^+, s\in \RR} \int_{H_0} f_s(n)dn,
\end{eqnarray}
be they finite or infinite.}

{ Now, for $\Re s>0$  the function $f_s$ is absolutely integrable over $H$ (for this statement, the character $\Lambda$ could be trivial). Indeed, we have
$ \int_{H} |f_s(g)|  \leq  \left( \int_H |v_s(u)| du \right)^2,$ and the integral $\int_H v_s(u)$ is known to be absolutely convergent from the study of standard intertwining operators.}
We therefore have
\begin{equation} \label{Jacquet} \int_{H} f_s(g)\Psi_\Lambda^{-1}(g) dg = \int_{H} \int_{H} v_s(ug)\overline{v_s(u)} du \Psi_\Lambda^{-1}(g) dg=\end{equation} 
\begin{equation*}=\left| \int_{H} v_s(u)\Psi_\Lambda^{-1}(u) du\right|^2,\end{equation*}
and $$W_s(\Lambda,g) = \int_{H} v_s(ug)\Psi_\Lambda^{-1}(u) du$$ is the Jacquet integral which converges absolutely for $\Re s>0$. (This is one of the assumptions that we made above, and is known to hold in the aforementioned cases, including the Whittaker case.)

We can let $\Lambda$ vary in the $k$-points of the one-dimensional vector space $\VV^*:=\Hom(\HH/\HH_0,\GGa)$, and then by \eqref{Jacquet} and inverse Fourier transform we get, 
still for $\Re(s) > 0$, 
\begin{equation} \label{gray2} \int_{H_0} f_s(n)dn = \int_{V^*} |W_s(\Lambda,1)|^2 d\Lambda,\end{equation}
for a suitable choice of Haar measure $d\Lambda$. 

Now -- as $\Lambda$ varies --  $W_s(\Lambda, 1)$ may be expressed in terms of the value of $W_s$ for a {\em fixed} character
at a varying point of the torus; 
by a routine computation with the known asymptotics of the spherical Whittaker function (s.\ the remark that follows),
  the integral (\ref{gray2}) is uniformly bounded for $s$ in a neighborhood of zero, and the right hand side of \eqref{limit2} is finite.
(We discuss this argument at more length below, phrased in a fashion where it generalizes more readily.) 
\end{proof}

\begin{remark}
Let us discuss in more detail how to phrase the final step of the proof -- bounding \eqref{gray2} -- in the language of this paper,
so that it may be readily generalized to other settings:

The asymptotics of the ``Whittaker function'' $W_s(\Lambda,g)$ can be derived from our earlier discussion of asymptotics in section \ref{sec:asymptotics}, by interpreting the value $W_s(\Lambda,g)$ as an element in a representation induced from a character of a spherical subgroup. Recall the cocharacter $\check\mu$ discussed before the statement of the proposition; in the Whittaker case, if $\HH$ is the subgroup corresponding to the \emph{negative} roots for some choice of Borel and Cartan subgroups, then $\check\mu$ is (in additive notation) a multiple of $-\check\rho$ (minus half the sum of positive coroots). In order to relate $W_s(\Lambda,g)$ to standard Whittaker functions, we only need to notice that changing $\Lambda$ corresponds to conjugating by an element of $\check\mu(k^\times)$; more precisely:
\begin{equation}\label{twoWhittaker}
 W_s(\Lambda,\check\mu(x)) = \delta_B^{s-\frac{1}{2}}(\check\mu(x))  \int_{H} v_s(u)\Psi_{\Lambda}^{-1}(\check\mu(x)u\check\mu^{-1}(x)) du.
\end{equation}

Indeed (the argument is quite simple, but we formulate it in some language that can be generalized):
\begin{itemize}
 \item The variety $\HH\backslash \GG$, equipped with the \emph{trivial} character of $\HH$, is a boundary degeneration of the same variety equipped with the line bundle $\mathcal L_\Psi$ corresponding to induction from $\Psi$; denote that boundary degeneration by $\XX_\Theta$. (In the Whittaker case, this will be the most degenerate case, so $\Theta=\emptyset$, but not in general.) The left action of the cocharacter $\check\mu$ has image in $\AA_{X,\Theta}$, with $\check\mu(\mathfrak o)$ mapping to $A_{X,\Theta}^+$.
 \item We may split the nonzero points of the one-dimensional vector space $V^*$ into a finite number of $\Gm$-orbits under the cocharacter $\check\mu$; denote them by $V_i$.
 \item Let $\delta_H$ be the modular character of the $k$-points of the algebraic group $\HH\check\mu(\GGm)$. We notice that the modular character $\delta_B$ on $\check\mu(k^\times)$ is inverse to the modular character $\delta_H$.
 \item If $\Lambda_i$ is a representative for $V^*_i$, then:
$$W_s(\Lambda_i,\check\mu(x)) = \int_{H} v_s(u\check\mu(x))\Psi_{\Lambda_i}^{-1}(u) du = $$
$$ = \delta_B^{s+\frac{1}{2}}(\check\mu(x)) \int_{H} v_s(\check\mu(x^{-1})u\check\mu(x))\Psi_{\Lambda_i}^{-1}(u) du =$$
$$ = \delta_B^{s+\frac{1}{2}}(\check\mu(x)) \delta_H(\check\mu(x)) \int_{H} v_s(u)\Psi_{\Lambda_i}^{-1}(\check\mu(x)u\check\mu^{-1}(x)) du = $$
$$\delta_B^{s-\frac{1}{2}}(\check\mu(x))  \int_{H} v_s(u)\Psi_{\Lambda_i}^{-1}(\check\mu(x)u\check\mu^{-1}(x)) du,$$
which shows \eqref{twoWhittaker}.
\end{itemize}

Therefore for a suitable choice of Haar measures we have: 
\begin{eqnarray}\int_{V_i^*} |W_s(\Lambda,1)|^2 d\Lambda = \int_{k^\times} \delta_B^{1-2s}(\check\mu(x))  |W_s(\Lambda_i, \check\mu(x))|^2 d(x^{-n_0}) = \nonumber\\
 = \int_{k^\times} \delta_B^{1-2s}(\check\mu(x)) |x|^{-n_0}  |W_s(\Lambda_i, \check\mu(x))|^2 d^\times x. \end{eqnarray}

Now we notice:\begin{itemize}
 \item For $|x|\gg 1$ we have $|W_s(\Lambda_i, \check\mu(x))|^2=0$.
 \item Since, by definition, $W_s(\Lambda,\bullet)$ is in the image of an operator: $I_B^G(\delta_B^s)\to C^\infty(H\backslash G,\Psi_\Lambda)$, and $\check\mu(\mathfrak o)\subset A_{X,\Theta}^+$, by the theory of asymptotics for $|x|\ll 1$ we will have that the function $x\mapsto W_s(\Lambda,\check\mu(x))$ is equal to a $k^\times$-finite function with generalized eigencharacters equal to the (unnormalized) exponents of $I_B^G(\delta_B^s)$.  For any given $\varepsilon>0$ the latter are bounded, for $s$ in a neighborhood of $0$, by $\delta_B^{-\frac{1}{2}+\varepsilon}(\check\mu(x))$. 
 \item Given this information on generalized eigencharacters to handle $|x|$ small, the vanishing for $|x| \gg 1$ to handle $|x|$ large,  and the fact that $W_s$ extends continuously to $s=0$ (and hence is pointwise bounded for $x$ in a compact set, for $s$ in a neighborhood of zero) to handle the remaining $x$, it follows that
$$ \delta_B^{1-\varepsilon}(\check\mu(x)) |W_s(\Lambda_i, \check\mu(x))|^2 $$
is (uniformly for $s$ close to $0$) integrable over $k^\times$. 
\end{itemize}

\end{remark}

\begin{corollary} \label{corollaryonnormalization}
 For every tempered representation $\pi$, there is a canonical normalization of the integral of matrix coefficients:
\begin{equation}\label{normalization}\int_{H}^* \left< \pi(u) v^1, v^2\right> \Psi_\Lambda(u)^{-1} du 
\end{equation}
as the evaluation at $\Lambda$ of the Fourier transform of the function:
$$u \in H/H_0 \mapsto \int_{H_0} \left<\pi(nu) v^1, v^2\right> dn.$$

\end{corollary}
Indeed, this Fourier transform is a {\em a priori} a distribution, but it is also invariant under an open compact subgroup of $A$, so it can be identified with a function in the complement 
of degenerate characters $\Lambda$.

Now, for any non-degenerate $\Lambda$ the normalized integral above defines a morphism: 
\begin{equation} \label{Mpidef} M_\pi^\Lambda: \pi \otimes\bar\pi\to C^\infty(H\backslash G,\Psi_\Lambda)\otimes C^\infty(H\backslash G,\Psi_\Lambda^{-1}), \end{equation} characterized by the property that $M_\pi^\Lambda (v_1 \otimes v_2)(1, 1) = \int_{H}^* \Psi_\Lambda^{-1}(u) \left<\pi(u)v_1, v_2\right> du$. (This can also be expressed in terms of Jacquet integrals, as we discuss in \S \ref{whit-explicit}.)  The hermitian dual of this, composed with the unitary pairing between $\pi$ and $\bar\pi$, will be denoted by $\theta_\pi^\Lambda$. Then Theorem \ref{ST:Plancherel} carries over:

\begin{theorem} \label{Whittakerplancherel}  The Plancherel formula for $L^2 (H \backslash G,\Psi_\Lambda)$ reads:  
$$\left<\Phi_1,\Phi_2\right> = \int_{\check G} \theta_\pi^\Lambda(\Phi_1\otimes\overline{\Phi_2})\mu_G(\pi).$$
In particular, for every tempered representation $\pi$ and $v\in \pi$ we have:
$$\int^*_{H} \left<\pi(u)v,v\right>  \Psi_\Lambda(u)  du \ge 0,$$
and $L^2 (H \backslash G,\Psi_\Lambda)$ is tempered\footnote{This is easy to deduce directly in the Whittaker case, since $U^{-}$ is amenable.} as a $G$-representation; its measure is absolutely continuous with respect to the Plancherel measure for $G$. 
\end{theorem}

\begin{proof} The proof that we saw in the strongly tempered case carries over almost verbatim,
with due care for the regularizations:

 Fix $f\in C_c^\infty(G)$, and let $\Phi_\Lambda(g)=\int_{H} f(ug) \Psi_\Lambda^{-1}(u) du$. Note that $\Phi_{\Lambda}$ is also compactly supported on $H \backslash G$,
 and so square integrable.  Let us consider:
           \begin{equation} \label{aleph-eqn} \aleph : \Lambda \mapsto  \Vert\Phi_\Lambda\Vert^2 = \int_{H} \Psi_\Lambda(u) \int_{\hat G}\left<\pi(u) i_\pi (f), i_\pi (f)\right>_{HS}    \mu_G(\pi)du.\end{equation} 
                     The equality here is proved in precisely the same way as \eqref{ForRefAA1} and \eqref{ForRefAA2}.

\eqref{aleph-eqn} is indeed a genuine, continuous function of $\Lambda$, locally constant on the non-degenerate locus. 
This statement follows easily from the facts that the support of $\Phi_{\Lambda}$ in $H \backslash G$ is compact, and the integrand
involved in the definition of $\Phi_{\Lambda}(g) $  is actually compactly supported (uniformly for $g$ in a compact set). 
On the other hand, the integral is no longer absolutely convergent in general as a double integral. To push through the previous computations, we shall study $\aleph$ as a distribution in $\Lambda$. 

Let $Q$ be a smooth function on $V^* = \Hom(\HH/\HH_0)(k)$, compactly supported away from the degenerate locus,  with Fourier transform $u\to \hat Q(u)$ (also compactly supported on $V$ and smooth). Then:
$$ \int_{V^*} Q(\Lambda) \cdot \aleph (\Lambda) d\Lambda = \int_{H} \hat Q(u) \int_{\hat G}\left<\pi(u) i_\pi (f), i_\pi (f)\right>_{HS}  \mu_G(\pi)du=$$
$$=  \int_{\hat G}\int_{H} \hat Q(u) \left<\pi(u) i_\pi (f), i_\pi (f)\right>_{HS}  du \,\, \mu_G(\pi)=$$ $$ = \int_{\hat G}  \int_{V^*} Q(\Lambda) \theta_\pi^\Lambda (\Phi \otimes \bar \Phi ) d\Lambda\,\, \mu_G(\pi)= \int_{V^*} Q(\Lambda) \int_{\hat G} \theta_\pi^\Lambda (\Phi\otimes\bar\Phi)\mu_G(\pi) \,\, d\Lambda,$$
where $\Phi(g)= \int_{H} f(ug) \Psi_\Lambda^{-1}(u) du$. The first equality on the final line, i.e. the introduction of $\theta_{\pi}^{\Lambda}$, is
just as in the proof of \eqref{ThetaExplication}.

To deduce the desired result from this, for any given non-degenerate $\Lambda_0$ there exists an open neighborhood $S$ of $\Lambda_0$ in $V^*$ so that $\aleph(\Lambda)$ and $\theta_\pi^\Lambda(\Phi \otimes\bar \Phi)^*$ are all constant on $S$ . We then choose $Q$ supported in $S$ to get the desired result.
\end{proof}
 
We note that a Plancherel formula for the Whittaker model has also been developed by Delorme \cite{DeWhi}, who used different methods. In the archimedean case, the Whittaker-Plancherel formula was developed by Wallach \cite{Wallach}.

A corollary to the last  sentence of the theorem was conjectured \cite[Conjecture 3.5]{LM}  by Lapid and Mao; this also follows from the results of Delorme.

 \begin{corollary} \label{LMconj}
 Suppose that $\pi$ is a generic irreducible representation of $G$ whose
 Whittaker functions are square integrable on $U^- \backslash G$; then $\pi$ is ($G$-)discrete series. (And similarly for all other cases of Theorem \ref{Whittakerplancherel}.)
 \end{corollary}
 
 { We note that the converse statement is also true: a generic, discrete series representation 
 has Whittaker functions which are square integrable on $U^{-} \backslash G$. This converse follows from the theory of asymptotics. 
 }
 \medskip

 \subsubsection{Explication} \label{whit-explicit}
The functionals of the theorem, when $\pi$ is an induced representation, can be described in terms of Jacquet integrals: 
\begin{equation} \label{fund}  \int_{U^-} \left<\pi(u)v^1,v^2\right> \Psi_\Lambda(u) du  = W^1_0(\Lambda) \overline{W^2_0(\Lambda)}, \end{equation} where $W^i_0(\Lambda)$ are constructed from $v^1, v^2$ by Jacquet integrals in an induced representation.  Precisely:

Any tempered representation is a direct summand of a representation of the form $I_P^G(\tau)$, where $\tau$ is a discrete series of the Levi quotient of $P$.   We equip it with the unitary structure
$\|v\|^2 = \int_{U_M^{-} \backslash U^{-}} \|v(u)\|_{\tau}^2 du$ (where $U_M^-=U^-\cap M$, $M$ a Levi subgroup of $P$).

If $\pi$ admits a Whittaker functional,
then $\tau$ does also.  In this case, we fix a Whittaker functional (unique up to a scalar of norm one)
$\eta_{\Lambda}:\tau \rightarrow \mathbb{C}$ so that
\begin{multline*} \eta_{\Lambda}( \tau(u) v) = \Psi_\Lambda(u) \eta_{\Lambda}(v), u \in U_M^{-}; 
\\  \eta_{\Lambda}(v_1) \overline{\eta_{\Lambda}(v_2)} 
= \int_{U_M^{-}}^* \left< \tau(u) v_1, v_2 \right> \Psi_{\Lambda}(u) du, \ \ v_1, v_2 \in \tau. \end{multline*}
(Indeed,  that the final integral can be thus factorized follows from multiplicity one for Whittaker functionals,
and 
the non-negativity of Theorem \ref{Whittakerplancherel}. It can also be verified that the resulting $\eta_{\Lambda}$ is nonzero
by means of the argument to be presented in the next section \S \ref{IIsubsec}, we briefly sketch the argument  -- namely, the known theory of asymptotics of Whittaker functions
mean that $\tau$ occurs inside the $L^2$-Whittaker space for $M$, and then Theorem \ref{Whittakerplancherel}
implies that the integral defining $\eta_{\Lambda}(v_1) \overline{\eta_{\Lambda}(v_2)} $ must have been nonzero.)

Now define the Jacquet integral  on  the family of representations $I_P^G(\tau\delta_P^s)$:
$$v \mapsto \int_{U_M^{-} \backslash U^{-}} \Psi_\Lambda(u) \eta_{\Lambda}(v(u))  du.$$
As before, this defines (after holomorphic continuation) a $(U^{-},\Psi_\Lambda)$-equivariant functional
$\Xi_{\Lambda} \in I_P^G(\tau \delta_P^s)^*$ without poles on the tempered axis $\Re s =0$. 
Let $v^1,v^2 \in \pi:=I_P^G(\tau)$, and let $W^1_0,W^2_0$ be the corresponding Whittaker functions $W^i_0 (g) = \Xi_{\Lambda}(g v^i)$ (we will also write: $W_0^i(\Lambda):=W_0^i(1)$).
The product $W^1_0(\Lambda) \overline{W^2_0(\Lambda)} $ is independent
of the choice involved in defining $\eta_{\Lambda}$.  Computing formally, 
the left-hand side of \eqref{fund} equals
\begin{multline}  \label{neededtorefertothis} \int_{(u,u') \in U_M^{-} \backslash (U^{-} \times U^{-})}  \Psi_\Lambda(u') 
\langle  v^1(uu'), v_2 (u) \rangle \ du \  du' 
\\ = 
 \int_{(u,u') \in U_M^{-} \backslash U^{-} \times U_M^{-} \backslash U^{-}}  \Psi_\Lambda(u') 
\eta_{\Lambda}(v^1(u u')) \overline{\eta_{\Lambda}(v^2(u))}  \ du \ du' 
\\  = W^1_0(\Lambda) \overline{W^2_0(\Lambda)}.
\end{multline} 

The conversion of this to a formal proof follows along the lines of the previous regularizations carried out in this section, and we omit it.

	\subsection{The Ichino--Ikeda conjecture} \label{IIsubsec}

The following establishes a conjecture of Ichino and Ikeda. 
\begin{theorem} \label{propIIconj}
Suppose that $\X$ is strongly tempered and wavefront, let $H$ be the stabilizer of a point on $X$.
Then:
\begin{equation}M_\pi: v\otimes \bar w\mapsto \label{intMC}\int_H \left< \pi(h) v, w\right> dh
\end{equation}
 defines a nondegenerate Hermitian form on $\pi_H$ (the $H$-coinvariants of $\pi$), for
every $G$-discrete series representation $\pi$. 

Moreover, if $\sigma$ is a discrete series representation of some Levi subgroup and $\pi=I_P^G(\sigma)$ (where $P$ is some corresponding parabolic), then the same expression defines a non-zero, non-negative hermitian form on $\pi_H$.
\end{theorem}

\begin{remark}
 In the setting of the ``Gross--Prasad variety'' $\SO_n\backslash\SO_n\times\SO_{n+1}$, it is known by recent work of Waldspurger \cite{Waldspurger-GP} that for any tempered $L$-packet of $G$ there is at most one element $\pi$ in the $L$-packet such that $\pi_H\ne 0$. Since any irreducible tempered representation is a subrepresentation of $I_P(\sigma)$ (where $\sigma$ is a discrete series of the pertinent Levi) and all subrepresentations of that belong to the same $L$-packet, it follows that in that case (\ref{intMC}) is non-zero for every irreducible tempered representation $\pi$ with $\pi_H\ne 0$. This is the conjecture originally formulated by Ichino and Ikeda. In fact, Beuzart-Plessis managed to refine the argument of our theorem \cite[Proposition 14.2.2]{B-P:SMF} to \emph{prove} multiplicity one in the induced discrete series, in the setting of the Gross--Prasad conjectures.
\end{remark}

 The idea of the proof:  By our discussion of asymptotics we can show that any embedding $\pi \hookrightarrow C^{\infty}(X/\mathcal Z(X),\omega_\pi)$ (where $\omega_\pi$ is the -- unitary -- central character of $\pi$) has image contained in $L^{2+\varepsilon}(X/\mathcal Z(X),\omega_\pi)$. This suggests that any $H$-distinguished $\pi$ must contribute to the Plancherel formula. 
 But we have already computed (Theorem \ref{ST:Plancherel}) the Plancherel formula for $X$
 and seen that $\pi$ occurs exactly when $M_{\pi}$ is nonzero. 
 
\proof

{ For $G$-discrete series $\pi$, by our discussion of asymptotics any morphism $\pi\to C^\infty(X/\mathcal Z(X),\omega_\pi)$ will have image in $L^2(X/\mathcal Z(X),\omega_\pi)$:

 Indeed, if $\Theta \subset \Delta_X$ then 
the $X_\Theta = X_\Theta^L \times^{P_\Theta^-} G$ and the asymptotic morphism $\pi \to C^{\infty}(X_{\Theta}) \to C^\infty(X_\Theta^L)$ factors through the corresponding Jacquet module $\pi_{{\Theta}^{-}}$.  Because $\pi$ is discrete series, all the characters of this Jacquet module
decay on the negative Weyl chamber of the center of $L_{\Theta}$. Since (a finite union of cosets of) this negative Weyl chamber surjects onto $A_{X, \Theta}^+$, by the proof of Proposition \ref{wavefrontlevi}, it follows that the image of any $v \in \pi$ ``decays in the $\Theta$-direction'';
since this is so for all $\Theta$ and the exponential map is measure-preserving, this image in fact belongs to $L^2$.  
\footnote{ To expand: Recall that, for this argument,   the action of $A_{X, \Theta}$ has been twisted already, as per
our general notational conventions from \S \ref{notation},
by the square root of the $A_{X, \Theta}$-eigencharacter for the measure on $X_{\Theta}$; thus the normalizations are precisely chosen so that decaying
exponents along $A_{X, \Theta}$ force square integrability. We also
used the fact that this twist is compatible with the twisting in the definition of the normalized Jacquet module.  This follows
because of the description of $X_{\Theta}$ as a parabolically induced variety, i.e. again from Proposition \ref{wavefrontlevi}.}

The claim on nondegeneracy of $M_{\pi}$, now, follows from Theorem \ref{ST:Plancherel}:  
indeed, the natural embedding: 
$$ \pi\otimes\Hom(\pi, C^\infty(X)) \hookrightarrow L^2(X/\mathcal Z(X),\omega_\pi),$$
$$ v\otimes M \mapsto M(v),$$
endows the space on the left with a non-degenerate hermitian form, of which the Plancherel form $\theta_\pi$ of Theorem \ref{ST:Plancherel} is just the hermitian dual (i.e.\ the dual via the pairing of $\pi\otimes\Hom(\pi, C^\infty(X))$ with $C_c^\infty(X)_{\tilde\pi} = \tilde\pi \otimes (\Hom(C_c^\infty(X),\tilde\pi))^*$). Thus $\theta_\pi$, and hence $M_\pi$, is non-degenerate.}

  We now turn to the claim on induced representations, which is more subtle; we need to pass from an ``almost-everywhere'' statement to a pointwise statement. We do this by establishing a uniform lower bound almost everywhere, and then we can specialize pointwise by a continuity argument. 

Let $\pi=I_P^G(\sigma)$ as in the statement with $\pi_H\ne 0$. We establish a series of claims:

\begin{itemize}
 \item Every $H$-distinguished direct summand of $\pi$ belongs to the support of Plancherel measure for $L^2(X)$. 
\end{itemize}

This follows by approximating matrix coefficients : Suppose that 
$\pi_0$ is an $H$-distinguished direct summand of $\pi$, so that we are given
a $G$-morphism $M: \pi_0 \rightarrow C^{\infty}(X)$.    The idea is now to approximate
the matrix coefficients of $\pi_0$ by truncating functions in the image of $M$. For a fixed open compact subgroup $J$, partition $X$ into $J$-fixed subsets:
$$ X=\sqcup_\Theta N_\Theta,$$
where $N_\Theta$ belongs to a $J$-good neighborhood of $\Theta$-infinity, is $A_{X,\Theta}^+$-stable and has compact image in $X_\Theta/A_{X,\Theta}$. { We let $\tilde N_\Theta= \bigcup_{\Omega\subset\Theta} N_\Omega$, a $J$-good, $A_{X,\Theta}^+$-stable neighborhood of $\Theta$-infinity.

For every $\gamma\in \Delta_X$, set $\hat\gamma:=\Delta_X\smallsetminus\{\gamma\}$, and choose an element $a_\gamma\in \mathring A_{X,\hat\gamma}^+$. Notice that $\mathring A_{X,\hat\gamma}/\mathcal Z(X)$ is a one-dimensional subtorus of $A_X/\mathcal Z(X)$. For $x\in X/\mathcal Z(X)$ define a \emph{radial function}:

$$ R(x):= \min\{ n\ge 1 | x\notin  a_\gamma^n \tilde N_{\hat\gamma} \mbox{ for any }\gamma\in \Delta_X\}.$$

It is then clear that the sets $X_n:=\{x\in X| R(x)\le n\}$ form an exhaustive, increasing filtration of $X$ by compact-mod-$\mathcal Z(X)$ sets. Moreover, they have the following property: 

\begin{quote}
 For any compact subset
$\Omega \subset G$, there exists an integer $n \geq 1$ so that $X_k\cdot \Omega \subset X_{k+n}$.
\end{quote}

Indeed, since the sets $X_n$ are by definition $J$-stable, one can replace $\Omega$ by a finite subset $\{g_i\}_i$. Then there is an $n$ such that for every $\gamma$ and $i$ we have: $a_\gamma^n \cdot \tilde N_{\hat\gamma} \cdot g_i^{-1} \subset \tilde N_{\hat\gamma}$ \emph{when $\tilde N_{\check\gamma}$ is considered as a subset of $X_\Theta$}, and hence also: 
$$ a_\gamma^{n+k} \cdot \tilde N_{\hat\gamma} \cdot g_i^{-1} \subset a_\gamma^k \cdot\tilde N_{\hat\gamma}$$ 
for all $k\ge 0$. Hence: $a_\gamma^{n+k} \cdot\tilde N_{\hat\gamma} \cdot \Omega^{-1} \subset a_\gamma^k \cdot\tilde N_{\hat\gamma}$, and by the equivariance property of $J$-good neighborhoods the same is true on $X$. Therefore, if a point $x$ is in $X_k$, i.e.\ does not lie in $a_\gamma^k \cdot\tilde N_{\hat\gamma}$ for any $\gamma$, then $x\Omega$ is in $X_{k+n}$, i.e.\ does not lie in $a_\gamma^{k+n} \cdot\tilde N_{\hat\gamma}$ for any $\gamma$.}

{ 
We also set $N_{\Theta,k}:= N_\Theta \cap X_k$ for any $k$.} {\bluetext By the theory of asymptotics, and the fact that $\pi_0$ is tempered, the quantity:
$$ \langle M(v_1), M(v_2) \rangle_{L^2(N_{\Theta,k}/\mathcal Z(X))} $$
either has a limit over $k$ -- this cannot happen for all $\Theta$ and all $v_1,v_2$ because $\pi_0$ is not discrete --  or is ``asymptotic to the integral of a generalized $A_{X,\Theta}$-eigenfunction with trivial generalized eigencharacter''. By the latter we mean that the function:
$$A_{X,\Theta}^+\ni a \mapsto M(v_1)(a\cdot x) \overline{M(v_2)(a\cdot x)} \Vol(a\cdot xJ), $$
for $x\in N_\Theta$, is (the restriction to $A_{X,\Theta}^+$ of)  a generalized eigenfunction with unitary and subunitary eigencharacters; and hence the integral of $ M(v_1)\overline{M(v_2)}$ over $N_{\Theta,k}$, as $k\to\infty$, is dominated by the integral of its summand corresponding to the trivial generalized $A_{X,\Theta}$-eigencharacter. By an easy calculation over finitely generated abelian groups, this means that 
there is a nonzero constant $c_\Theta(v_1,v_2)$ and a positive integer $r_\Theta (v_1,v_2)$ such that:
$$ \langle M(v_1), M(v_2) \rangle_{L^2(N_{\Theta,k})} \sim c_\Theta(v_1,v_2) k^{r_\Theta(v_1,v_2)}. $$}

{ Thus, there is a positive integer $r (= \max_{\Theta,v_1,v_2} r_\Theta(v_1,v_2))$ so that for every triple $(\Theta,v_1,v_2)$ the limit:
\begin{equation} \label{asymp-mc-ii} c(v_1,v_2):= \lim_k \frac{\langle M(v_1), M(v_2) \rangle_{L^2(X_k)}}{k^r}\end{equation}   exists, and moreover it is not zero for all such triples.}

{ 
Then the limit $c = c(v_1, v_2)$ defines a nonzero  Hermitian form on $\pi_0$.
We claim that it is $G$-invariant; indeed, for $g \in \Omega$, 
we have an inequality
$$\langle M(v_1), M(v_2) \rangle_{L^2(X_{k-n})} \leq
\langle M(gv_1), M(gv_2) \rangle_{L^2(X_k)} \leq \langle M(v_1), M(v_2) \rangle_{L^2(X_{k+n})}$$
 from which the invariance follows.   Therefore $c(v_1, v_2)$   is a multiple of the inner product.

Now, we can approximate (on compacta) matrix coefficients of $\pi_0$ by matrix coefficients of  $L^2(X)$:
 namely, given $v \in \pi_0$ and $g \in \Omega$ 
we may approximate $\langle g v, v \rangle$, which is a multiple of $c(gv, v)$,
by a multiple of $\langle M(g v), M(v) \rangle_{L^2(X_k)}$ for large $k$. Now this equals
$$\langle g M(v)|_{X_k}, M(v)|_{X_k} \rangle + \langle g M(v)|_{X_k} - M(gv)|_{X_k}, M(v)|_{X_k} \rangle$$
and the second term is bounded by Cauchy-Schwarz by the square root of
$$ \langle M(v), M(v) \rangle_{L^2(X_{k+n} - X_{k-n})} \cdot \langle M(v), M(v) \rangle_{L^2(X_k)}$$
which is of lower order than the main term as $k \rightarrow \infty$, because of  \eqref{asymp-mc-ii}.}
 Consequently, $\pi_0$ belongs to the support of Plancherel measure for $L^2(X)$.

\begin{itemize}
 \item There is a set of unitary unramified characters of $P$ of positive (Haar) measure such that $M_{\pi_\chi}\neq 0$, where $\pi_\chi = I_P(\sigma\otimes \chi)$.
\end{itemize}

Indeed, the only tempered representations in a neighborhood of $\pi$ under the Fell topology are of the form $\pi_\chi$. Since we know (from strong temperedness, i.e.\ Theorem \ref{ST:Plancherel}) that the Plancherel measure for $L^2(X)$ is supported in the set of tempered representations and is absolutely continuous with respect to the natural class of measures on them, it follows that there should be a set of unramified characters $\chi$ of positive Haar measure such that the Plancherel morphisms $M_{\pi_\chi}$ are non-zero.

As $\chi$ varies, we identify the underlying vector spaces of all $\pi_\chi$ in the natural way with a fixed vector space $V$. 

\begin{itemize}
 \item For $v\in V$, the expression $M_{\pi_\chi}(v,v) $ is a real analytic function of $\chi$.
\end{itemize}

Indeed, the integral of the matrix coefficient over $H$ is actually a countable sum of functions polynomial in $\chi$, which converges absolutely and uniformly in $\chi$.

It follows from the last two points that $M_{\pi_\chi}$ is non-zero for almost every $\chi$. There remains to show:
\begin{itemize} \item  $M_{\pi_{\chi}}$ is non-zero for every $\chi$.
\end{itemize} 
 For this, we will treat for simplicity the case of multiplicity one, i.e.\ $$M_{\pi_\chi}(v\otimes w)= L_{\pi_\chi}(v)\otimes \overline {L_{\pi_\chi}(v)},$$ where $L_{\pi_\chi} : \pi_\chi\to \CC$, an $H$-invariant functional.  It is easy to see that we may choose $L_{\pi_{\chi}}$ measurable in $\chi$. Moreover, we will assume that $G$ acts transitively on $X$. The proof is the same in the general case, but its formulation would obscure the argument.

Let $K_1$ be an open compact subgroup such that $L_{\pi_\chi}|_{\pi_\chi^{K_1}}\ne 0$ for a dense set of $\chi$'s. Let $v_\chi=\overline{K_1* L_{\pi_\chi}}\in \pi_\chi$, i.e.
$v_{\chi} \in \pi_{\chi}^{K_1}$ and $\langle u, v_{\chi} \rangle = L_{\pi_{\chi}}(u)$
for $u \in \pi_{\chi}^{K_1}$.  In particular, $v_{\chi}  \neq 0$ if and only if $M_{\pi_{\chi}}(v,w) \neq 0$
for some $v,w \in \pi_{\chi}^{K_1}$;  
so, the set of $\chi$ for which $\|v_{\chi}\| \neq 0$ is of full measure.

By Corollary \ref{smoothfnls}, there exists $K_2$ such that \begin{equation} \label{ctst} L_{\pi_\chi}(g v_\chi)= \left< g v_\chi,K_2* L_{\pi_\chi}\right> \ \ (g \in G^+).\end{equation} 
We are going to average both sides of \eqref{ctst} over $\chi$, and compute the $L^2$-norm. 
The left-hand side can be computed via the Plancherel formula for $X$; the right hand side, via the Plancherel formula for $G$. This will lead -- eventually -- to a ``almost-everywhere'' lower bound on the norm of $v_{\chi}$; by continuity, we will deduce
that $M_{\pi_{\chi}}$ is everywhere nonzero.

Let $Z$ be the set of distinct isomorphism classes of the unitary representations $\pi_\chi$; it has the canonical structure of an orbifold, and the restriction of the canonical Plancherel measure $\mu$ on $G$ is a well defined measure on $Z$. We will write $\pi_z, M_{\pi_z}, L_{\pi_z}, v_z$ etc.

For any $Z'\subset Z$ we consider the $K_1$-invariant $L^2$ function: 
$$\Phi: x\mapsto \int_{Z'} L_{\pi_z} (\pi_z(g) v_z)\mu(z),  \ \ x = Hg\in H \backslash G$$ on $X=H\backslash G$. 
For $x = Hg$, with $g \in G^+$, it coincides with the value of 
the  $K_2\times K_1$-invariant  function:
$$f: g\mapsto \int_{Z'} \left< \pi_z(g) v_z, K_2*L_{\pi_z}\right> \mu(z)$$
on $G$.
In what follows, it is harmless to replace $G^+$ by $K_2 \cdot G^+ \cdot K_1$. 

In order to compare norms we need the following fact (see discussion after
Corollary \ref{wavefront}):
\begin{quote}
 There is a constant $C>0$ such that the (surjective) orbit map
$$o: K_2\backslash G^+/K_1\to X/K_1$$
satisfies:
$$\Vol(o(S))\le C\cdot \Vol(S)$$
for any set $S$. (Here the volume of a subset $S \subset K_2\backslash G^+ / K_1$
is the volume of the corresponding $K_2, K_1$-invariant subset of $G$.) 
\end{quote}

Indeed, writing $x_0 \in X$ for the basepoint that corresponds to the identity coset in our identification $X \simeq H \backslash G$,  we have $x_0 g K_1  = x_0 K_2 g K_1$, which is covered by at most $[K_2 \backslash K_2 g K_1]$ translates of $x_0 K_2 $, all of which have equal measure. Therefore, the $X$-volume of the set $x_0 g K_1$ is bounded above by a constant multiple (not depending on $g\in G^+$) of 
the $G$-measure of $K_2 g K_1$ . 

Consequently we have:
$$\int_{Z'} \Vert v_z \Vert^2 \mu(z) \stackrel{(a)}{=} \Vert \Phi\Vert^2_{L^2(X)} = \int_{X/K_1} |\Phi(x)|^2 dx \le $$ $$\stackrel{(b)}{\le} C\cdot  \int_{K_2\backslash G^+/K_1} |f(g)|^2 dg \le C\cdot  \int_G |f(g)|^2 dg = $$
$$ \stackrel{(c)}{=}  C\cdot \int_{Z'} \Vert v_z\Vert^2 \Vert K_2* L_{\pi_z}\Vert^2 \mu(z).$$

Here equality (a) is the Plancherel formula for $X$; inequality (b) arises from the equality
$f(g) = \Phi(Hg)$ for $g \in G^+$ and our previous remark on measures; equality (c) is the Plancherel formula for $G$. 

Because $Z'$ is arbitrary, and the set of $z$ for which $\|v_z\| = 0$ has measure $0$, it follows that $\Vert K_2* L_{\pi_z}\Vert^2\ge C^{-1}$ for almost every $z$, and because it is continuous \footnote{For the continuity, note that $\Vert K_2* L_{\pi_z}\Vert|^2$ can be expressed  as the
sum of values of $M_{\pi_z}(w, w)$ over an orthonormal basis $w$ for $\pi_z^{K_2}$. In turn, this is expressed as a finite sum of
integrals $\int f_z(h) dh$ over $H$, where
each $f_z$ is a matrix coefficient, varying continously in $z$. 
Using the majorization  \eqref{mcbound} we see that, given $\varepsilon  > 0$, we can find a compact
subset $\Omega \subset H$ such that $\int_{h \notin \Omega} |f_z(h)| dh < \varepsilon$ for all $z$.
The continuity is now clear.} in $z$ it follows that the same holds for \emph{every} $z$.  
In particular, $L_{\pi_z}\ne 0$ for every $z$.
This completes the proof.
\qed 
 
 \newpage 
\part{Spectral decomposition and scattering theory}

\section{Results}

\label{pennylane}

This is the core part of the present paper, where we develop the Plancherel decomposition of $L^2(X)$.  By this we mean, more precisely, that we reduce the Plancherel decomposition
to the understanding of discrete series for $X$ and its boundary degenerations $X_{\Theta}$. 

Our approach to this is  different to previous works on similar topics -- in particular,
the series of papers establishing the Plancherel decomposition for real semisimple symmetric spaces \cite{BSP1,BSP2,De}.  Our viewpoint is close to that of Lax and Phillips on scattering theory, and further
from the viewpoint motivated by global Eisenstein series. The original idea, and a core argument in our approach, is due to Joseph Bernstein.

Our main theorem is Theorem \ref{advancedscattering}, and our approach to its proof
proceeds as follows:
 \begin{enumerate}
 \item[i.]  We first derive the existence of morphisms $L^2(X_{\Theta}) \rightarrow L^2(X)$, canonically characterized by their asymptotic properties,  by using general Hilbert space theory and the existence of asymptotics of eigenfunctions.  
 \item[ii.] 
 By using {\em a priori} information about  when a representation $\Pi$ can occur
simultaneously in $L^2(X_{\Theta}) $ and $L^2(X_\Omega)$, where $\Theta \neq \Omega$, we are able to analyze
  interactions between these maps from (i).  This allows us, in particular,  to decompose
  $L^2(X)$ in terms of the discrete parts of $L^2(X_{\Theta})$. 
  \item[iii.]  Finally,  we discuss
the question of writing an explicit formula for the morphism, using the Radon transform. (This corresponds to the study of ``Eisenstein integrals''.) 
\end{enumerate}

Our results are not complete in general.   Step (i) works very generally -- it uses only asymptotics of eigenfunctions as an input, and is performed in Sections \ref{sec:linearalgebra}--\ref{sec:Bernstein}.  Step (ii), performed in Sections \ref{sec:directintegrals}--\ref{sec:scattering}, requires a few ``non-formal'' inputs to work. 
Nonetheless, these non-formal inputs (``generic injectivity'', \S \ref{genericinjectivity} and  the Discrete Series Conjecture \ref{dsconjecture})
hold in a wide variety of cases (including symmetric spaces); and we expect the theorem to hold true in general, 
whether or not these inputs are valid. {\bluetext Step (iii) is performed in Section \ref{sec:explicit} under much more restrictive conditions, namely that the variety is strongly factorizable; the theory of Eisenstein integrals and their applications is open in the general case.}

 The reader who wishes to get an idea of our methods without diving into the details 
 may wish to read \S \ref{toymodel}. There we discuss the simple case of linear operators
acting on functions on the non-negative integers, and explain how our methods operate in this 
(very well-known) case. 

\subsection{Plancherel decomposition and direct integrals of Hilbert spaces} \label{ssdirectintegrals}

We recall first some generalities about direct integrals and Plancherel decomposition.

Any (separable) unitary representation $\mathcal H$ of $G$ admits a direct integral decomposition:
\begin{equation}
 \int_{\hat G} \mathcal{H}_{\pi} \mu(\pi)
\end{equation}
where $\hat G$ denotes the unitary dual of $G$, equipped with the Fell topology, and $\mu$ is a positive measure on the Borel $\sigma$-algebra of $\hat G$. (Note that these are standard Borel spaces, by \cite[Proposition 4.6.1]{Dixmier}, i.e.\ isomorphic as measurable spaces with the Borel space of a Polish space -- which can be taken to be the interval $[0,1]$ -- even though the topologies on $\widehat{G}$ will be usually non-Hausdorff.)

The Hilbert space $\mathcal H_{\pi}$ is $\pi$-isotypic (that is, a finite or countable direct sum of copies of $\pi$), and the direct integral makes sense only after we specify a family $F$ of ``measurable'' sections $\eta_:\hat G\ni \pi\mapsto \eta_\pi\in \mathcal H_\pi$ satisfying the following axioms:
\begin{enumerate}
 \item A section $\pi\mapsto\xi_\pi$ lies in $F$ if and only if for each $\eta\in F$ the function $\pi\mapsto\left<\eta_\pi,\xi_\pi\right>$ is measurable.
 \item There is a countable subfamily $\{\eta_i\}_i\subset  F$ such that for all $\pi\in \hat G$ the vectors $\{\eta_{i,\pi}\}_i$ span a dense subspace of $\mathcal H_\pi$.
\end{enumerate}
Here and later, we will be using the word ``measurable'' to mean measurable with respect to the completion of the Borel $\sigma$-algebra of $\hat G$ with respect to a given measure and a given family of measurable sections into Hilbert spaces, which should be clear from the context; when no measure is present, we mean Borel measurable. We will call this decomposition a \emph{Plancherel} decomposition (and the corresponding measure a Plancherel measure) if $\mu$-almost all spaces $\mathcal H_\pi$ are non-zero.

 \subsection{Discrete spectrum} \label{ss:dsdifficulties}
 
Before stating our results more precisely, we  describe the precise notion of ``discrete spectrum''
and discuss a difficulty that arises. It should be noted that this difficulty vanishes in the case of symmetric varieties, and thus (to our knowledge) does not arise in prior work. 

 The space $L^2(X)$ decomposes into a direct sum of a ``discrete'' and a ``continuous'' part:
\begin{equation*}
 L^2(X)=L^2(X)_\disc\oplus L^2(X)_\cont
\end{equation*}
More precisely, discrete means ``discrete modulo center'', where ``center'' is the connected component of $\Aut_G(X)$.
More formally, 
any $f \in L^2(X)$ can be disintegrated
as an integral $\int_{\omega} f_{\omega} d\omega$ indexed by characters $\omega$ of $\Aut(X)$; here, each $f_{\omega} \in L^2(X; \omega)$.   
\begin{definition}
$L^2(X)_{\disc}$
consists, by definition, of those $f$ for which almost every $f_{\omega}$ 
belongs to the direct sum of all irreducible subrepresentations of $L^2(X;\omega)$.
\end{definition}

It should be noted that it is by no means clear that this defines a closed subspace, although it follows from our later considerations.

Discrete series for a reductive group $G$ come in continuous families, with varying central character, which can be constructed by twisting matrix coefficients by characters of the group. A similar ``twisting'' is possible for $X$-discrete series when $\XX$ is a symmetric variety; however, for a general spherical variety $\XX$ with infinite automorphism group (hence, in the split case, $X$-discrete series appearing in continuous families), the variation of $X$-discrete series with ``central character'' presents serious challenges, which we analyze in section \ref{sec:discrete}. We expect that in every case the variation of $X$-discrete series can be described in terms of ``algebraically twisted'' families of representations, and we formulate this expectation as the ``Discrete Series Conjecture'' \ref{dsconjecture}. We show how this can proven in individual cases by the method of ``unfolding''; however, since we have not proven that this method always applies (although we know of no counterexample), the 
Discrete Series Conjecture remains a conjecture -- but an easy one to check in each individual case. Therefore, formulating theorems which are conditional on this conjecture seems to present no serious harm of generality or applicability.

At the first reading the reader might prefer to skip most of this section (\S \ref{sec:discrete}), using it only as a reference for terms and properties encountered in the rest of the paper. 
 In sections \ref{sec:Bernstein} and \ref{sec:scattering} we obtain our main results on the spectral decomposition, which are summarized in the theorem below.  

\subsection{Main result} \label{subsec:mainresult} 

  For any $\Theta, \Omega\subset\Delta_X$ (possibly the same),  
 
set $$W_X(\Omega,\Theta)= \{ w \in W_X: w \Theta=\Omega\}$$ 
Let $\mathfrak a_{X,\Theta} = \varchi(\AA_{X,\Theta})^* \otimes \mathbb Q$, $\mathfrak a_{X,\Theta}^*$ its $\mathbb Q$-linear dual. The vector space $\mathfrak a_{X,\Theta}$ has a ``dominant'' chamber $\mathfrak a_{X,\Theta}^+$ (namely, its intersection with the dominant chamber of $\mathfrak a_X$, which is a face of the latter), and it is known by properties of root systems that the union of the sets $w^{-1} \mathfrak a_{X,\Omega}^+$, $w\in W_X(\Omega,\Theta)$, where $\Omega$ varies over all possible subsets of $\Delta_X$ with $|\Omega|=|\Theta|$, is a perfect tiling for $\mathfrak a_{X,\Theta}$.  Here, by ``perfect tiling''
we mean that these sets cover $\mathfrak a_{X, \Theta}$ and their interiors are disjoint, cf. \cite[Lemma 1]{ArthurCorvallis}. \footnote{This ``perfect tiling'' claim is a simple property
of root systems; we give a proof. 
  In what follows, we discuss as if the system of spherical roots arises from a genuine Lie algebra; this would be the Lie algebra of the dual group to $\check{G}_X$. Now, any $\lambda \in \mathfrak{a}_{X, \Theta}$ defines  
a subset of roots, namely those which are non-negative on $\lambda$. This corresponds to the set of roots of a parabolic subalgebra $\mathfrak{p}$, with Levi factor $\mathfrak{m}$
given by the centralizer of $\lambda$.  Now there is a unique element $w$ of the Weyl group taking $\mathfrak{p}$
to a standard parabolic subalgebra $\mathfrak{q}$. This element $w$ carries the center $\mathfrak{a}_P$ of $\mathfrak{m}$ to the center 
$\mathfrak{a}_Q$ of the standard Levi of $\mathfrak{q}$. Moreover, $w$ is well-defined up to the Weyl group of $\mathfrak{q}$,
which acts trivially on $\mathfrak{a}_Q$, i.e. the map $\mathfrak{a}_P \rightarrow \mathfrak{a}_Q$ is uniquely determined.}
  Let:  
\begin{equation*}c(\Theta)=\mbox{ the number of those chambers in }\mathfrak a_{X,\Theta}=\sum_{\Omega} \# W_X(\Omega,\Theta)
\end{equation*}
where the second equality is a consequence of the tiling property. 

{ For the statement of the following theorem, recall that while the map: $\mathcal Z(\LL_\Theta)^0\to \AA_{X,\Theta}$ is surjective as a map of algebraic tori, it may not be surjective at the level of $k$-points.}

   We prove:

\begin{theorem}[Scattering theorem]\label{advancedscattering} Suppose that $\XX$ is a wavefront spherical variety and that all degenerations $X_\Theta$ (including $X_\Theta=X$) satisfy the Discrete Series Conjecture \ref{dsconjecture}. Also assume the condition of ``generic injectivity of the map: $\mathfrak a_X^*/W_X\to \mathfrak a^*/W$ on each face'' (cf.\ \S \ref{genericinjectivity})

 Then there exist canonical $G$-equivariant morphisms (\emph{``Bernstein morphisms''}) $$\iota_{\Theta}: L^2(X_\Theta) \rightarrow L^2(X)$$ and $A'_{X,\Theta}\times G$-equivariant isometries (\emph{``scattering morphisms''}):
$$S_w: L^2(X_{\Theta}) \longrightarrow L^2(X_\Omega), \ \ w \in W_X(\Omega,\Theta)$$
where $A_{X,\Theta}$ acts on $L^2(X_\Omega)$ via the isomorphism: $\AA_{X,\Theta}\xrightarrow{\sim} \AA_{X,\Omega}$ induced by $w$
{ and $A_{X,\Theta}'$ denotes the image of $\mathcal Z(L_\Theta)^0$ in $A_{X,\Theta}$,}
with the following properties:

 \begin{eqnarray} \label{compo1}
  \iota_{\Omega} \circ S_w & =&  \iota_{\Theta}, \\ 
 \label{composition}
S_{w'} \circ S_w &=& S_{w'w} \text{ for }  w\in W_X(\Omega,\Theta), w'\in W_X(Z, \Omega),
\\   \label{compo3} \iota_{\Omega}^* \circ \iota_{\Theta} &=& \sum_{w \in W_X(\Omega,\Theta)} S_w. 
\end{eqnarray}

Finally, the map: 
\begin{equation} \label{totalmapdef}
\sum_{\Theta} \frac{\iota_{\Theta,\disc}^*}{\sqrt{c(\Theta)}} : L^2(X)\to \bigoplus_{\Theta\subset\Delta_X} L^2(X_\Theta)_\disc\end{equation}
(where $\iota_{\Theta,\disc}^*$ denotes the composition of $\iota_{\Theta}^*$ with orthogonal projection to the discrete spectrum) is an isometric isomorphism onto the subspace of vectors $$(f_{\Theta})_\Theta \in \bigoplus_\Theta L^2(X_{\Theta})_\disc$$ satisfying:
$$S_w f_{\Theta} = f_{\Omega} \text{ for every }w \in W_X(\Omega,\Theta).$$
\end{theorem}

  We observe one way in which our approach differs from the usual one: 
for the proof of this Theorem, {\em we do not write down any explicit formula for $\iota_{\Theta}$;
instead, it is entirely characterized in terms of asymptotic properties. } 
The question of writing a formula for $\iota_{\Theta}$ is the concern of \S \ref{sec:explicit}; 
our main results are Theorem \ref{explicitBernstein} and Theorem \ref{explicitPlancherel}.
We do not summarize them here because they are more technical.  Our results  are complete in many cases (including the group case), but not in complete generality.

We do not know if every wavefront spherical variety satisfies either of the two algebraic multiplicity conditions of Theorem \ref{advancedscattering}, though we know no counterexample. However, they can easily be seen to fail in the case of the non-wavefront variety $\GGL_n\backslash \SSO_{2n+1}$. As was the case for the theory of asymptotics, we expect our scattering theory to extend to such cases as well:

\begin{conjecture}  
The conclusions of Theorem \ref{advancedscattering} are true for any spherical variety $\XX$.
\end{conjecture}

\section{Two toy models: the global picture and semi-infinite matrices} \label{toymodel}

We now consider two  ``toy models'' for our reasoning in this paper:
\begin{itemize}
\item[-] 

The first is not really a toy model: it is the global picture for automorphic forms. However,
it may be more familiar to readers of this paper than the local setting. 

The main point we wish to convey is that the relationship between ``smooth'' asymptotics $C^{\infty}_c(X_{\Theta}) \rightarrow C^{\infty}_c(X)$ and  the $L^2$-morphisms $L^2(X_{\Theta}) \rightarrow L^2(X)$
corresponds -- in this global setting --  to the difference to integrating Eisenstein series ``far from the unitary axis''
and ``on the unitary axis.''  This point of view is not central to us (although we establish, in the case of symmetric varieties and more generally, a corresponding result) but it  is closer to other developments and is useful to keep in mind. 

\item[-] 
The second is ``scattering theory on the non-negative integers'': 
Given a semi-infinite symmetric real matrix $A_{ij}$ with the property that
$A_{ij}$ depends only on $i-j$ when $i$ and $j$ are both large, what
can we say about its eigenvalues and the corresponding $L^2$-spectrum?

Here (unlike the global picture) we sketch proofs, so the reader may get
a simple idea of the techniques that we will use in the spherical variety case.
\end{itemize}

\subsection{Global picture}

Let $\mathbf{G}$ be a semisimple algebraic group over a number field $K$;
fix a maximal $F$-split torus with Weyl group $W$ and a minimal $K$-parabolic containing it. 
Let $\mathbf{P} = \mathbf{M} \mathbf{N}$ be a standard (i.e., containing the fixed minimal parabolic) parabolic $K$-subgroup
and denote by $\AA_P$ the center of $\mathbf{M}$.  
We write $w(\mathbf{P})$ for the number of chambers in $\Hom(\varchi(\mathbf{M}),\RR)$ (a {\em chamber} is a connected component of the space obtained by removing the kernel of roots.)
Write $X_{\mathbf{P}} = \mathbf{N}(\adele) \PP(K) \backslash \mathbf{G}(\adele)$
and set $X = X_{\mathbf{G}}$, so that 
$$X = \mathbf{G}(K) \backslash \mathbf{G}(\adele).$$

 By $C^{\infty}(X)$ we shall mean the set of functions that are invariant by an open compact
 of $\mathbf{G}(\adele_{\mathrm{finite}})$ and are smooth along each translate of $G_{\infty} := \G(K \otimes \RR)$. 
 By  $C^{\aut}(X)$ we mean the space of {\em automorphic forms} on $X$ (defined in the standard way, see e.g.\ \cite{BorelJacquet});  similarly for $X_{\PP}$. 
 
 \subsubsection{Smooth asymptotics}
 
 Recall the {\em constant term}:
 $$ c_{\PP}: C^{\infty}(X) \longrightarrow C^{\infty}(X_{\PP})$$
 obtained by integration along ``horocycles'':
 $c_{\PP} f(g) = \int_{\mathbf{N}(K) \backslash \mathbf{N}(\adele)} f(ng) dn,$
 the integral being taken with respect to the $\mathbf{N}(\adele)$-invariant probability measure. 
This is the  analog of our smooth asymptotics map $$e_{\Theta}^*: C^{\infty}(X) \rightarrow C^{\infty}(X_{\Theta});$$
 indeed under certain stronger smoothness assumptions on $f \in C^{\infty}(X)$ (for example, 
 $L^2$-boundedness of all elements $X f$, where $X$ belongs to the universal enveloping
 algbra of $\G(F \otimes \RR)$),  it is known
 that $f$ is asymptotically equal to $c(f)$, in a suitable sense and in a particular direction.\footnote{This has roughly the same content as the following fact, due to Harish-Chandra and Langlands: If $f$ is an automorphic form on $X$, then the truncation $\wedge f$ is of rapid decay.}

In the adjoint direction we have the pseudo-Eisenstein series:
$$  e_{\PP}: C^{\infty}_c(X_{\PP}) \longrightarrow C^{\infty}_c(X)$$
(with rapidly decaying image) defined by the rule $e_{\PP} f: g \mapsto  \sum_{P(K) \backslash G(K)} f(\gamma g).$
This is the analog of our smooth dual asymptotics map $$e_{\Theta}: C^{\infty}_c(X_{\Theta}) \rightarrow C^{\infty}_c(X).$$   

\subsubsection{Eisenstein series}
 
 We will be very brief, with the understanding that this section
 is targeted only at readers who have previous experience with the theory of Eisenstein series. 

In what follows, $\omega$ will denote a character of $A_P := \AA_P(\adele)/\AA_P(K)$, not necessarily unitary; 
the space  $ C^{\aut}(X_{P})_{\omega}$  denotes the subspace of $C^{\aut}(X_P)$ comprising functions that transform under the character $\omega$ under the {\em normalized} action of $A_P$.   Also if $\omega$ is unitary we
may form the space $L^2(X_{\PP})_{\omega}$
of functions $f$ which transform under $A_P$ and so that $|f|^2$ is integrable on $X_{\PP}/A_P$. 

We fix a Haar measure on $A_P$ and a dual measure on the Plancherel dual $\widehat{A_P}$. The set of all, not necessarily unitary, characters of $A_P$ will be denoted by $\widehat{A_P}_\CC$.

For all $\omega$ with ``sufficiently large real part'', and every $f \in C^{\aut}(X)_{\omega}$ the series defining
$e_{\PP}  f $ is absolutely convergent. Moreover, for a Zariski-open set of $\omega \in \widehat{A_P}_{\C}$,  this series can be regularized by a process of meromorphic continuation.   We refer to the result
as $e_{\PP, \omega}$.

The analog of our ``Bernstein'' $L^2$-morphism $\iota_\Theta: L^2(X_{\Theta}) \rightarrow L^2(X)$ is given by the unitary Eisenstein series:
\begin{eqnarray} \label{l2asymptotics} E_{\PP}: L^2(X_{\PP})  & \rightarrow  &L^2(X), \mbox{ characterized by:}   \\   \int_{\omega \in \widehat{A_P}} f_{\omega} d\omega & \mapsto &   \int_{\omega \in \widehat{A_P}}  e_{\PP, \omega} f_{\omega} d\omega \end{eqnarray}
for $f_{\omega} \in C^{\aut}(X_{\PP})_{\omega} \cap L^2(X_{\PP})_{\omega}$ that varies  measurably in $\omega$ (suitably interpreted).  
This map satisfies an asymptotic property analogous to that characterizing
$L^2(X_{\Theta}) \rightarrow L^2(X)$ (see Theorem \ref{Bernsteinmap}). 

How does this compare to the ``smooth'' asymptotics $C^{\infty}_c(X_{\Theta}) \rightarrow C^{\infty}(X)$? Given a ``sufficiently positive'' real character $\omega_0: A_P \rightarrow \mathbb{R}_{>0}$,  the morphism $e_{\PP}$ is characterized by: 
\begin{eqnarray}  e_{\PP}: C^{\infty}_c(X_{\Theta})  & \rightarrow & C^{\infty}(X) \\  
  \label{smoothasymptotics}  \int_{\mathrm{Re}(\omega)=\omega_0} f_{\omega}  & \mapsto &
\int_{\mathrm{Re}(\omega)=\omega_0} e_{\PP, \omega} f_{\omega}. \end{eqnarray} 
whenever the left-hand side belongs to $C^{\infty}_c(X_{\Theta})$.  (In fact, a function in $C^{\infty}_c(X_{\Theta})$ may be uniquely expressed as $\int_{\mathrm{Re}(\omega) = \omega_0} f_{\omega}$, so long as the real part of $\omega_0$ is sufficiently positive.)

The difference, then, between smooth and $L^2$-asymptotics is (from this point of view)
the ``line of integration.'' One can pass from smooth to $L^2$ by shifting contours:
on the other hand, this will introduce extra contributions coming from residues
of the map $e_{\PP, \omega}$. 
 
The maps $E_{\PP}$ satisfy an exact analog of Theorem \ref{advancedscattering}; see, for example, the ``main theorem'' on page 256 
of \cite{ArthurTraceFormula}, and for more details \cite{Langlands76, MoeglinWaldspurger}.

\subsection{Spectra of semi-infinite matrices: scattering theory on $\N$}

Let $C(\N)$ and $C_c(\N)$ denote the vector spaces of all functions, resp.\
all compactly supported functions, on the set $\N$ of non-negative integers.

Fix some real numbers $c_0, c_1, \dots, c_K$; we set $c_k = c_{|k|}$ if $|k| \leq K$ and $c_k = 0$ otherwise. 
We shall consider real self adjoint operators $T: L^2(\N) \rightarrow L^2(\N)$ that are given by the rule
\begin{equation} \label{Tasymp} T f(x) = \sum_{|k| \leq K} f(x+k) c_{|k|} \end{equation}
for all large enough $x$, i.e.\ there exists $M$ so that this equality holds whenever $x \geq M$. 
(Such operators are easily seen to exist: for example, \eqref{Tasymp} defines a self-adjoint
operator on $L^2(\Z)$, and then one composes with the orthogonal projection to $L^2(\N)$.)

We denote by $T_{\infty}$ the operator on $L^2(\Z)$ defined by rule \eqref{Tasymp} {\em for all $x$};

In what follows, 
let $P \in \mathbb{C}[z, z^{-1}]$ be defined by $\sum_{k=-K}^K  c_{|k|} z^k$;  we think of this as a meromorphic function of the complex variable $z$. 
We denote by $P^{-1}(\lambda)$ the set $ \{ z \in \C: P(z) = \lambda\}$.

\subsubsection{Smooth asymptotics}  \label{731}

There is a unique morphism
$$ e: C_c(\Z) \longrightarrow C_c(\N)$$
that intertwines the $T_{\infty}$ and $T$-action, and  carries the characteristic function $\delta_k$ of $k \in \Z$ (considered as an element of $C_c(\Z)$) to $\delta_k$ (now considered
as an element of $C_c(\N)$), if $k$ is sufficiently large.

Indeed,  $v_k := e(\delta_k)$ can be solved for inductively by means of the linear recurrence that they satisfy.  Indeed, 
$$  c_{-K} v_{m-K} + \dots + c_K v_{m+K} = e( T_{\infty} \delta_m) = T e(\delta_m) = T v_m$$
which is to say that we can determine $v_{s}$ given knowledge of $v_{s'}$  and $T v_{s'}$ when $s' > s$. 

The meaning of the phrase ``asymptotics'' is clearer for the
 dual morphism $e^*: C(\N) \rightarrow C(\Z)$.    Let $C^{\infty}(\N)^{\lambda}$ (resp. $C^{\infty}(\Z)^{\lambda}$) be  the space of $T$-eigenfunctions on $\N$ with eigenvalue $\lambda$
(resp. $T_{\infty}$-eigenfunctions on $\Z$ with eigenvalue $\lambda$). Then the dual asymptotics map  $e^*$ gives a natural (not always injective)
map 
\begin{equation} \label{asymp-lambda} C^{\infty}(\N)^{\lambda}  \rightarrow C^{\infty}(\Z)^{\lambda} \end{equation}
i.e.\ for any $v \in  C^{\infty}(\N)^{\lambda} $ there is a unique $w \in  C^{\infty}(\Z)^{\lambda} $
so that $v(n) = w(n)$ for all large $n$. 

Notice that any eigenfunction of $T$ (and of $T_\infty$) necessarily is of the form
\begin{equation}\label{Teigenfunction}
f(n) = \sum_{i=1}^{k} \alpha_i^n Q_i(n)  \ \ \text{ for }n > M+K,
\end{equation}
where  $\alpha_1, \dots, \alpha_k \in \mathbb{C}$ and $Q_i(n)$ are nonzero polynomials of total degree $\leq 2 K$; also, $M$ is as in \eqref{Tasymp}. 
If there indeed exists such an eigenfunction,  the eigenvalue is necessarily given by $P(\alpha_1) = P(\alpha_2) = \dots =P(\alpha_k)$.

\subsubsection{Finiteness of the discrete spectrum}    \label{732}

We show now that the eigenfunctions of $T$ in $L^2$
span a finite-dimensional space. 
(It is easier to see that the corresponding assertion for $C_c(\N)$:
indeed $C_c(\N)$ is a finitely generated module over $\C[T]$.)

Consider an eigenfunction of $T$ of the form (\ref{Teigenfunction}). If it is to belong to $L^2$ we must have $|\alpha_i| < 1$ for every $i$.

Fix for a moment the degrees $d_1, \dots, d_k$ of the polynomials $Q_i$. 
Then the condition that there exists an eigenfunction with the above ``asymptotic''
is a finite system of linear equations in the coefficients of the $Q_i$
as well as the $M+K$ values $f(0), \dots, f(M+K)$; 
the coefficients of this linear system  
depend algebraically on the $\alpha_i$. In particular the
set of $(\alpha_1, \dots, \alpha_n)$ for which there exist a solution with $\deg(Q_i) =d_i$
is a {\em constructible} subset of $\mathbb{C}^n$. (As usual, {\em constructible} means
that it is defined by a finite system of polynomial equalities and inequalities.)

Thus we have a constructible subset $Z \subset \mathbb{C}^k$ whose intersection with $\{ \underline{\alpha} : |\alpha_i| < 1\}$ 
is countable -- the corresponding finite-dimensional eigenspaces are orthogonal in the separable Hilbert space
$L^2(\mathbb{N})$. (Notice that in the above process we fixed the \emph{degree} of the polynomials, and not just an upper bound; so there are no redundant $\alpha_i$'s for each eigenfunction.) Therefore, this intersection is in fact finite:

 {  To verify the finiteness, we use the following properties of 
Zariski-closed subsets  $Y \subset \mathbf{C}^k$: There are only finitely many zero-dimensional (Zariski-)irreducible components of $Y$,
and if $p$ does not lie on such a component, then any neighborhood of $p$ in $Y$, for the standard topology on $\C^k$, is uncountable.  For the 
latter assertion we just note  that there exists
an analytic nonconstant map $f:\C \rightarrow Y$ with $f(0) = p$, which we may verify by slicing by a hyperplane to reduce to thc case of $\dim_p(Y)=1$, i.e. a curve, where
we can use the existence of a smooth neighborhood on a desingularization $\tilde{Y} \rightarrow Y$.
With this in mind, let $Z$ be a constructible set; it is a finite union of sets of the form $Y_i'$, where each $Y_i'$ is Zariski-open
in a Zariski-closed set $Y_i \subset \C^k$.  Let $B$ be any open subset of $\mathbb{C}^k$ for the usual topology.
If $ p \in Z \cap B$ is not one of the finitely many zero-dimensional components of some $Y_i$, 
our discussion above has shown that there exists
uncountably many points in $Z \cap B$. }

\subsubsection{$L^2$-decomposition}
 \label{734}

 There is a unique bounded map $$\iota: L^2(\Z)  \hookrightarrow L^2(\N)$$ 
 (in fact, with image in $L^2(\N)_{\cts}$) 
that intertwines  the actions of $T$ and $T_{\infty}$ and is ``asymptotically the natural identification'': 
\begin{equation} \label{iotatoydef} \| \iota \delta_k - \delta_k\|_{L^2(\N)} \rightarrow 0, \ k \rightarrow \infty. \end{equation} 
  We will not give a complete proof of this statement but we will outline the identification
  between the spectra of $T$ and $T_{\infty}$ that is the central ingredient.

   According to spectral theory, 
we may find a measure $\mu(\lambda)$ on the real line
together with an isomorphism 
\begin{equation}\label{Ndecomp} L^2(\N) \stackrel{\sim}{\rightarrow} \int_{\RR} \mathcal{H}_{\lambda} \mu(\lambda)\end{equation}
which carries the action of $T$ to multiplication by $\lambda$. Here $\mathcal{H}_{\lambda}$
is  a family of Hilbert spaces -- finite-dimensional in the present case -- over $\lambda \in \RR$.

 Now let $\nu = \mu - \sum_{\lambda} \mu(\{\lambda\}) \delta_{\lambda}$, the measure obtained from $\mu$ by removing all atoms. 
 This yields a corresponding decomposition of $L^2(\N)_{\cts}$:
 $$L^2(\N)_{\cts} \stackrel{\sim}{\rightarrow} \int_\RR \mathcal{H}_{\lambda} \nu(\lambda).$$

 In particular, when analyzing $L^2(\N)_{\cts}$, we can and do  neglect  the finite set of $\lambda$  for which there is a solution to $P(z) = \lambda$ {\em with multiplicity $>1$}, because
 that set has $\nu$-measure zero. 
  
 Because $T$ acts on $\mathcal{H}_{\lambda}$ as multiplication by $\lambda$, 
the morphism $C_c(\N) \rightarrow \mathcal{H}_{\lambda}$ 
must necessarily factor through the quotient $C_c(\N)_{\lambda}$ of $C_c(\N)  $
generated by all $ Tf - \lambda f$, for $f \in C_c(\N)$.   
In particular, $\mathcal{H}_{\lambda}$ is finite dimensional, since $C_c(\N)$ is a finitely generated $\C[T]$-module; 
also, any linear functional on  $\mathcal{H}_{\lambda}$ may be expressed
$f \mapsto \int f \phi_{\lambda}$, where $\phi_{\lambda} \in C(\N)^{\lambda}$ satisfies $T \phi_{\lambda} = \lambda \phi_{\lambda}$. 
(Here we write $\int$ for the functional $f \in C_c(\N) \rightarrow \sum_{x \in \N} f(x)$). 

Because of the asymptotics of eigenfunctions \eqref{Teigenfunction}, the 
image of $f$ in $C_c(\N)_{\lambda}$ 
 depends only on the values $f(0), \dots, f(M+K)$ together
with the values of the ``Fourier transform''  $\hat{f}(z)$, for 
$z \in P^{-1}(\lambda)$; here $\hat{f}(z) = \sum_{n} f(n) z^n$, the ``Fourier transform'' of $f$. 
It follows that $\mathcal{H}_{\lambda}$
can be identified with the completion of $C_c(\N)$ with respect to a Hermitian form of the type\footnote{Warning:  Neither $f \mapsto f(n)$, for $n \leq M_K$, nor the functionals $f \mapsto \hat{f}(z) \ \ (z \in P^{-1}(\lambda))$ have to descend to the quotient $C_c(\N)_{\lambda}$.}

\begin{equation} \label{HermitianFormToy} H_\lambda(f) = \sum_{1 \leq n,m \leq M+K} b_{n,m}(\lambda) f(n) \overline{f(m)} + \sum_{z, z' \in P^{-1}(\lambda)} a_{z, z'} \hat{f}(z) \overline{\hat{f}(z')},\end{equation}
It is not difficult to see that (away from a set of measure zero) 
$a_{z,z'} = 0$ unless $|z| \leq 1$ and $|z'| \leq 1$. To prove this, one works with a ``Schwartz'' space slightly larger than $C_c(\N)$,  allowing functions that have rapid (faster than any polynomial)
decay at $\infty$; then $H_{\lambda}$
must extend continuously to $C_c(\N)$, which translates into the desired vanishing. We do not give details here, but the idea
is due to Bernstein \cite{BePl} and the analogous step in our context is in  part (2) of Corollary \ref{HCform}.

The terms   $a_{z, z'} $ themselves include the ``diagonal'' case where $z'=z^{-1} = \bar z$; for those, set $a_z:= a_{z,z'}$. If we fix a function $f$ and consider an average of the hermitian forms of its (right) translates:
$$ \frac{1}{k+1} \sum_{i=l}^{l+k} H_\lambda(S^i f),$$
(where  $S$ is the translation operator $S f(x)=f(x-1)$ where we extend $f$ by zero off $\mathbb{N}$, and $l$ is arbitrary) then  in the limit as $k\to \infty$ only the ``diagonal'' terms survive. 
The reason for this is all other terms $a_{z,z'}$  occur with a coefficient involving $\frac{1}{k+1}\sum_{i=l}^{l+k} (z \overline{z}' )^i$, and this sum goes to zero for large $k$ in the non-diagonal case. 

On the other hand:

\begin{eqnarray} \label{limint} \Vert f\Vert^2_{L^2(\N)} &=& \frac{1}{k+1} \sum_{i=l}^{l+k} \Vert S^i f\Vert_{L^2(\N)} = \lim_{k\to\infty} \frac{1}{k+1} \sum_{i=l}^{l+k} \Vert S^i f\Vert_{L^2(\N)_\cts} =\nonumber\\
 &=& \lim_{k\to\infty} \frac{1}{k+1} \sum_{i=l}^{l+k} \int_\RR H_\lambda(S^i f) \nu(\lambda).
\end{eqnarray}

In the discussion that follows we will justify the fact that one can interchange the limit and the integral\footnote{We thank Joseph Bernstein for pointing out an omission in a previous version.} in order to arrive at the conclusion:

\begin{equation} \label{toy-conclusion} \|f\|^2 =  \int_{\RR} \nu(\lambda) \sum_{z \in P^{-1}(\lambda): |z| =1} |\hat{f}(z)|^2  a_z. 
\end{equation}

Notice that, if we think of $f$ as a function on $L^2(\Z)$, the above expression is invariant under left translation $S^{-1}$, in particular: it holds for \emph{every} element of $C_c^\infty(\Z)$. The uniqueness of the spectral decomposition for $T_\infty$ acting on $L^2(\Z)$ implies that it has to coincide with it. More precisely, writing $e(\theta) := e^{2 \pi i \theta}$ for $\theta \in \RR/\Z$,
$$\slope(\theta) := \frac{d}{d\theta} P(e(\theta)),$$
and $\slope(z)$ if $z=e(\theta)$, we can write the formula:
$$\|f\|^2 =  \int_{\RR/\Z} |\hat{f}(e(\theta))|^2 d \theta ,$$  (where the measure on $\theta \in \RR/\Z$ is the Haar probability measure) as: 
$$\|f\|^2 =  \int_\RR d\lambda \sum_{e(\theta) \in P^{-1}(\lambda)} |\hat{f}(e(\theta))|^2    |\slope(\theta)|^{-1};$$
this is the spectral decomposition for $L^2(\Z)$ under $T_\infty$. We conclude that:
$$ \int  d\lambda \sum_{z \in P^{-1}(\lambda) : |z|=1} |\hat{f}(z)|^2    |\slope(z)|^{-1} =  \int \nu(\lambda) \sum_{z \in P^{-1}(\lambda): |z| =1} |\hat{f}(z)|^2  a_z.$$
and from this we deduce
that $\nu(\lambda)$ is absolutely continuous with respect to the Lebesgue measure; without loss of generality, we may take it to  {\em equal} Lebesgue measure (restricted to the set where the expression on the left is supported), and then \begin{equation} \label{azz}  a_z  = |\mathrm{slope}(z)|^{-1}.\end{equation} 

This analysis shows, in effect,  that some part of the Plancherel formula for $T$ acting on $L^2(\N)$
is determined by the Plancherel formula for $T_{\infty}$ acting on $L^2(\Z)$.  This is the starting point
     for the construction of the morphism $\iota$ of \eqref{iotatoydef}, although we will not give details of that construction. We will use this type of idea in 
   \S \ref{sec:Bernstein}.

 \begin{remark} The Hermitian form $H_\lambda$ in fact satisfies more constraints, which force many values $a_{z, z'}$ to be zero; we don't discuss this here.  \end{remark}
  
We now proceed to justifying the interchange of limit and integral in \eqref{limint}. The argument will be more involved than the one used later in Section \ref{sec:Bernstein}, where we will have a priori knowledge of the fact that the values of $z$ in the above expression for $H_\lambda$ with $|z|<1$ are actually uniformly bounded, in absolute value,  by a constant $c<1$. Here we do not know that, and we will need to separate the set of $\lambda$'s in two.
 
First of all, fix the function $f$ and let us consider the effect of shifting the function $f$, i.e., replace it by the function 
$S f: n \mapsto f(n-1)$; we interpret $f(-1), f(-2)$ etc. as $0$.  

Associated to $f$ there is a (positive) spectral measure $\mu_f = H_\lambda(f) \nu(\lambda)$ on the set of $\lambda$'s, where $\nu(\lambda)$ is as in \eqref{Ndecomp}. 

We will need in particular the following key lemma:

\begin{quote} Let $L(\delta)$ be the set of eigenvalues $\lambda$
for which there is $z \in P^{-1}(\lambda)$
such that $1-\delta < |z|  < 1$. Then:
\begin{equation}\label{deltai}
 \lim_{(\delta,i)\to (0,\infty)}\mu_{S^i f}(L(\delta)) = 0.
\end{equation}

\end{quote}

Granted that, we can split the integral in \eqref{limint} into an integral on $\RR\smallsetminus L(\delta)$ and an integral on $L(\delta)$, with the contribution of the latter to the right-hand side of the being less than any given $\varepsilon>0$ as long as $\delta$ is small enough and $l$ is large enough. A result of linear algebra (Proposition \ref{invub}) allows us to apply the dominated convergence theorem to the former as $k\to \infty$ -- cf.\ also the proof of Theorem \ref{bernstein-abstract}. Since $\varepsilon$ is arbitrary, we arrive at \eqref{toy-conclusion}.

Finally, we are left with proving \eqref{deltai}. Since $L(\delta)$ is the intersection with $\RR$ of $\{P(z)|\,\, 1-\delta<|z|<1\}$, it is  a semialgebraic set and so a union of intervals (which may or may not contain their endpoints); 
moreover, the union of the sets $\{\delta\}\times L(\delta)$ is a semialgebraic subset of $(0,1)\times\RR$. Also   the number of components of $L(\delta)$ is bounded independently of $\delta$, as follows from the following fact, known as Hardt triviality \cite{Hardt} (although there is surely an elementary proof): 
\begin{quote}  Given a map $f: X \rightarrow Y$ of semialgebraic sets, there exists a finite partition 
of $Y$ into semi-algebraic subsets so that, on each part $Y_i$, the map $f$ is (semialgebraically) trivial, i.e.,
there exists a semialgebraic homeomorphism of $f^{-1}(Y_i) \rightarrow Y$
with $F_i \times Y_i \rightarrow Y_i$, where $F_i$ is any fiber. 
\end{quote}  
Since $L(\delta)$ is decreasing as $\delta \rightarrow 0$ with $\cap L(\delta) = \emptyset$,
 these intervals must all have length that goes to zero as $\delta$ does. 

 So it is enough to show that,
for any $\varepsilon > 0$ there are $k$ and $N$ so that:
$$\left( \mbox{ the $\mu_{S^i f}$-measure of   any interval of the form $\left( \frac{m}{2^k}, \frac{m+1}{2^k} \right)$, with $m \in \Z$ }\right)    < \varepsilon $$
for all $i  > N$. Notice that the support of all $\mu_{S^i f}$'s lies in a compact set, namely $[-\|T\|, \|T\|]$ where $\|T\|$ is the operator norm of $T$,  so it is enough to consider a finite number of those intervals.
 
As we will see, the measures $\mu_{S^i f}$ converge weakly to $\mu_f^\infty$ (the spectral measure of $f$ considered as an element of $L^2(\Z)$), which is absolutely continuous with respect to Lebesgue measure. Being, in addition, positive, we can eventually bound the $\mu_{S^if}$-measures of the aforementioned intervals in terms of the $\mu_f^\infty$-measures of slightly larger intervals, which in turn tend to zero with the length of those. 

To show the weak convergence, since these are positive measures whose total mass is uniformly bounded above,  it is enough to show it on polynomials. We have:
$$\int \lambda^m \mu_{S^i f}  =  \langle T^m S^i f , S^i f \rangle$$
which for large $i$ is equal to:
$$\langle T_{\infty}^m S^i f, S^i f \rangle = \langle T_{\infty}^m f, f \rangle = \int \lambda^m \mu_f^\infty.$$
This shows that $\mu_{S^i f}$ converge to $\mu_{f}^{\infty}$, as we just explained above, and concludes our argument 
 for the interchange 
 of limit and integral in \eqref{limint}.

\section{The discrete spectrum}  \label{sec:discrete}

This section addresses the discrete spectrum of $L^2(X)$, and in particular its variation with central character; see \S \ref{ss:dsdifficulties} for a general discussion of the difficulties involved, and \S \ref{newissue} for more details.

\subsection{Decomposition according to the center} \label{decomp-center}
We let $\widehat{\mathcal Z(X)}_\CC$ denote the group of complex characters of $\mathcal Z(X):=\Aut_G(X)^0$; it is a commutative complex group with infinitely many components (unless $\mathcal Z(X)$ is trivial). The identity component is the full torus of unramified characters of $\mathcal Z(X)$, and (of course) each connected component is a torsor thereof. The subgroup of unitary unramified characters acts transitively on the unitary characters of each component, inducing a canonical ``imaginary'' structure on each of them (i.e.\ the structure of a real algebraic variety which we will call ``imaginary''). The subgroup of unitary characters of $\widehat{\mathcal Z(X)}_\CC$ will simply be denoted by $\widehat{\mathcal Z(X)}$. We {\bluetext keep assuming, as we have done since Section \ref{sec:asymptotics}, that $\mathcal Z(\GG)^0$ surjects onto $\mathcal Z(\XX)$.} Notice that when we replace $\XX$ by a boundary degeneration $\XX_\Theta$, this means that we also 
have to ``enlarge'' the group $\GG$ in order to take into account the additional action of $\AA_{X,\Theta}=\mathcal Z(\XX_\Theta)$.

If we fix a Haar measure on $\mathcal Z(X)$, the maps
$$p_\omega(\Phi)(x)= \int_{\mathcal Z(X)} (z\cdot \Phi)(x) \omega^{-1}(z) dz$$
are surjective maps: $p_\omega: C_c^\infty(X)\to C_c^\infty(X, \omega)$, for every $\omega\in \widehat{\mathcal Z(X)}_\CC$. Here we denote by $C_c^\infty(X, \omega)$ the space of smooth functions on $X$, whose support has compact image under $X\to X/\mathcal Z(X)$ and which are eigenfunctions with eigencharacter $\omega$ under the (normalized) action of $\mathcal Z(X)$.  The fixed measures on $X$ and $\mathcal Z(X)$ induce a measure on $\mathcal Z(X)\backslash X$ and hence $L^2$-norms on the spaces: $C_c^\infty(X, \omega)$; the corresponding Hilbert space completions will be denoted by $L^2(X, \omega)$.

The Plancherel formula for $L^2(X)$ viewed as a unitary representation of $\mathcal Z(X)$ reads:
\begin{equation}
 L^2(X)=\int_{\widehat{\mathcal Z(X)}} L^2(X,\omega) d\omega;
\end{equation}
\begin{equation}\label{centerdecomp}
 \Vert \Phi \Vert_{L^2(X)}^2 = \int_{\widehat{\mathcal Z(X)}} \Vert p_\omega \Phi\Vert^2_{L^2(X,\omega)} d\omega,\,\,\text{ for all }\Phi\in C_c^\infty(X)
\end{equation}
where $d\omega$ is the Haar measure on $\widehat{\mathcal Z(X)}$ dual to the chosen measure on $\mathcal Z(X)$. From now on we will consider as fixed a measure on $\mathcal Z(X)$.

A \emph{relative discrete series (RDS) representation} for $X$, or an \emph{$X$-discrete series representation}, is a pair $(\pi, M)$ where $\pi$ is an irreducible smooth representation of $G$ with unitary central character $\omega$ and $M$ is a morphism: $\pi\to C^\infty(X)$ whose image lies in $L^2(X, \omega)$. The morphism $M$ induces a canonical unitary structure on $\pi$, pulled back from $L^2(X, \omega)$. The images of all such $M$, for given $\omega$, span the \emph{discrete spectrum} $L^2(X,\omega)_\disc$ of $L^2(X,\omega)$.

 Assuming that the orthogonal projections: $L^2(X,\omega)\to L^2(X,\omega)_\disc$ are measurable in $\omega$ (we will prove this in Proposition \ref{measurable}), they define a direct summand $L^2(X)_\disc$, the \emph{discrete spectrum} of $X$, which has a Plancherel decomposition according to the center: 
\begin{equation}
 L^2(X)_\disc=\int_{\widehat{\mathcal Z(X)}} L^2(X,\omega)_\disc d\omega;
\end{equation}

It is obvious that $X$-discrete series belong to $L^2(X)$ in the sense of Fell topology:

\begin{lemma}
 If $(L,\pi)$ is an RDS then $\pi$ belongs weakly (i.e.\ under the Fell topology) in $L^2(X)$.
\end{lemma}

\begin{proof}
 Choose a continuous section $Y$ of $X\to X/\mathcal Z(X)$, a ``radial'' function: $s:\mathcal Z(X)\to \RR_+$ (for example, project
 $\mathcal{Z}(X)$ to its quotient by its maximal compact, giving a free abelian group $\Lambda$, and let $s$ be the restriction
 of a Euclidean norm on $\Lambda \otimes \mathbb{R}$)
  and let $\wedge^T: C^\infty(X)\to L^2(X)$ denote ``normalized truncation'': 
$$ \wedge^T(\Phi) (x)= \left\{\begin{array}{l} \frac{\sqrt{\eta(z)}}{\sqrt{\Vol(s^{-1}[0,T])}} \Phi(x),\, \text{ if } x=y\cdot z\in Y\cdot s^{-1}[0,T] \\
0,\, \text{ otherwise }
\end{array}\right.$$
(Recall that $\eta$ denotes the eigencharacter of the fixed measure on $X$.) Then each diagonal matrix coefficient $\left<\pi(g) v, v\right>$ of $\pi$ is approximated, uniformly on compacta, by the diagonal matrix coefficients $\left< g\cdot \wedge^T(Lv),\wedge^T(Lv)\right>$ of $L^2(X)$, as $T\to\infty$.
\end{proof}

However, it is not as clear how the relative discrete series vary as one changes the central character $\omega$, and hence how to decompose the space $L^2(X)_\disc$. The reader can jump ahead to subsection \ref{newissue} for a discussion of the difficulties that one faces. We will first discuss the issue of relative discrete series abstractly, in order to obtain the uniform boundedness of their exponents which we use for the spectral decomposition.

\subsection{A finiteness result} \label{subsec:finiteness}

The purpose of this section is to show that the problem of constructing relative discrete series for $X$ is, essentially, described by polynomial equations, and to draw certain very general conclusions from this fact, in particular:

\begin{theorem}\label{finiteds}
 For a fixed open compact subgroup $J$ of $G$ and a unitary character $\chi$ of $\mathcal Z(X)$ the space $L^2(X,\chi)_\disc^J$ is finite-dimensional.
\end{theorem}

In the proof we shall use the following simple lemma:

\begin{lemma}  \label{linalg} Suppose $V,W$ are finite dimensional vector spaces,  $R$
a complex algebraic variety (identified with the set of its complex points), and
$\{T_i(r)\}_{i \in I} \subset \Hom(V,W) \ (r \in R)$ a (possibly infinite) collection of linear maps which vary polynomially in $r$ (that is, $T_i\in \Hom(V,W)\otimes \CC[R]$).
\begin{enumerate}
\item For each $r$, write $ W_r :=  \mathrm{span} \{ T_i(r): \ i \in I\}$. 
There exists a finite subset $J \subset I$ so that $W_r$ is spanned
by $\{T_j(r):j\in J\}$, for all $r$. 
\item Assume that $I$ is finite, and write $V_r$ to be the common kernel of all $T_i(r)$. 
There exists a constructible partition $R = \bigcup R_{\alpha}$
so that over each $R_\alpha$, $V_r$ is a trivializable sub-bundle of $V$ (i.e.\  there exists a vector space $V_1$
and isomorphisms $f_r: V_1 \stackrel{\sim}{\rightarrow} V_r$ varying algebraically with $r\in R_\alpha$.). 
\item Assume that $I$ is finite, suppose $V \subset W$, and let $\mathcal G$ denote the Grassmannian of all subspaces of $V$. The subset in $R\times\mathcal G$ consisting of a point $r\in R$ and a subspace of $V$ stable under each $T_i(r)$ is Zariski-closed. 
\end{enumerate}
\end{lemma}
\begin{proof} 
For the first and second assertions:
 Let $r_0$ be a point at which $\mathrm{rank} \  W_r $ (resp. $\mathrm{rank} \ V_{r}$)
 is maximized (resp.\ minimized).  It is easy to see that, in the first case, any finite collection of operators whose image spans $W_{r_0}$ at $r=r_0$ also spans $W_r$ for $r$ in a Zariski open neighborhood $U$ of $r_0$. In the second case, we can find a set of rational sections of vectors in $W$ which belong to the kernel of all $T_i$, are regular at $r_0$ and their specializations form a basis for the common kernel of all $T_i(r_0)$. Then there is a Zariski open neighborhood $U$ of $r_0$ where all these sections remain regular and linearly independent, and the dimension of the common kernel remains the same.  We then replace $R$ by $R-U$ and argue by induction on dimension.

As for the final assertion: 
the  set of $r$ in question is the inverse image in $R\times\mathcal G$ of a Zariski-closed subset
of $(\Hom(W,V))^J \times \mathcal G$, under the map $(r,g)\mapsto (\prod_j T_j(r),g)$.
 \end{proof}

\begin{proof}[Proof of the theorem]
The idea of the proof is to exploit the a priori knowledge of the fact that the set of $X$-discrete series is countable. On the other hand, we can describe the set of $J$-fixed vectors in $C^\infty(X)$ generating admissible subrepresentations as a constructible set. Constructibility and countability, together, will imply finiteness. See \S \ref{732} for this argument in a simpler setting,
and in particular see the final paragraph of \S \ref{732} for an explanation of why ``constructability+countability $\implies$ finiteness.'' 

For notational simplicity, assume at first that $\mathcal Z(X)$ is trivial.
For an assignment $$R: \Delta_X\supset\Theta\mapsto r_\Theta\in \mathbb N,$$ define the ``space of $R$-exponents'':
\begin{equation} \label{spectrumRdef} \spectrum_R = \prod_{\Theta} (\widehat{\mathcal{Z}(X_{\Theta})}_{\C}^J)^{r_\Theta}, \end{equation} where $J$ denotes  
``$J$-fixed.'''\footnote{The notion of $J$-fixed means the group of characters of $\mathcal Z(X_\Theta)$ which appear as eigencharacters on $C^\infty(X_\Theta)^J$.
{ Because our $\mathcal Z(X_{\Theta})$ is a quotient of $\mathcal Z(L_{\Theta})$,  by Proposition
\ref{wavefrontlevi},
this amounts to the character group of the quotient of $\mathcal Z(X_{\Theta})$ by an open compact subgroup; in particular, the set of $J$-fixed characters
is the group of homomorphisms from a finitely generated abelian group to $\mathbb{C}^{\times}$, and thus has a natural structure of algebraic variety.}}
   In words: an element of $\spectrum_R$
is a collection of characters of $\mathcal Z(X_{\Theta})$ for all $\Theta$, in fact, $r_{\Theta}$ of them for given $\Theta$.
 We denote an element of $\spectrum_R$ by $x=(x_j)$,   where the subscripts $j$ are of the form $(\Theta, i)$ with $i \leq r_{\Theta}$ an integer;  we denote by  
$\chi_{x_{\Theta, 1}},  \dots,  \chi_{x_{\Theta, r}}$ for $r \leq r_{\Theta}$ the various characters of $\mathcal Z(X_{\Theta})$ 
indexed by $x$.

The set $\spectrum_R$
has the natural structure of a complex variety. 
Let us fix a
$J$-good neighborhood $N'_{\Theta}$ for each $\Theta$, and let $N_\Theta=N_\Theta'\smallsetminus \cup_{\Omega\subset\Theta} N_\Omega'$. This includes $N_\Theta'=X$ for $\Theta=\Delta_X$, therefore the sets $N_\Theta$ cover $X$. Notice that, if we identify $N_\Theta/J$ with a subset of $X_\Theta/J$ via the exponential map, then this subset is compact (finite) modulo the action of $\mathcal Z(X_\Theta)$. 
For this reason, the set of functions on $\mathcal Z(X_\Theta) \cdot N_\Theta/J$ (as a subset of $X_\Theta/J$) annihilated by: 
\begin{equation}\label{annihilation}\prod_{i=1}^{r_\Theta}(z- \chi_{x_{\Theta, i}}(z)),\end{equation}
for all $z \in \mathcal{Z}(X_{\Theta})$, forms a finite dimensional vector space.  (Here we denote, for notational simplicity, simply by $z$ the normalized action of $z \in \mathcal{Z}(X_{\Theta})$ on functions on $X_\Theta$, which before was denoted by $\mathcal L_z$.)
 
Hence, we can form for each $x = (x_j) \in \spectrum_R$, the finite-dimensional subspace 
$V_x \subset C(X)^J$ consisting of functions $f$ with the following property: for every $\Theta$, there is a function $f_\Theta\in C(X_\Theta)^J$ that is annihilated,  for all $z \in \mathcal{Z}(X_{\Theta})$,  by (\ref{annihilation}) and such that $f|_{N_\Theta}=f_\Theta|_{N_\Theta}$; here we are using the identification of $N_{\Theta}/J$ with  subset of $X_{\Theta}/J$. 
 We shall say, for simplicity, that elements of $V_x$ are functions ``asymptotically annihilated'' by  
 the ideal $I_{x, \Theta}  \subset \C[\mathcal{Z}(X_{\Theta})]$
 generated by (\ref{annihilation}).

Notice that for any admissible subrepresentation $\pi\subset C^\infty(X)$ there is an $R$ and a $x\in \spectrum_R$ such that the image of $\pi^J$ belongs to $V_x$; more precisely, the image of $\pi^J$ is a finite-dimensional subspace of $V_x$ which is stable under the Hecke algebra for $J\backslash G/J$. { Vice versa, any $\mathcal H(G,J)$-stable subspace of $V_x$ is a finite module over the Bernstein center and hence generates a finite length (hence admissible) subrepresentation of $C^\infty(X)$.} 

There is a finite subset $S$ of $X/J$, depending only on $R$, such that for any given $x\in\spectrum_R$, each element of $V_x$ is determined by its restriction on $S$; indeed, the point $x$ determines ``characteristic polynomials'' for the action of all elements of $\mathcal Z(X_\Theta)$ on the asymptotics of $V_x$, which amount to recursive relations by which these asymptotics are determined by a finite number of evaluations. In particular, $V_x$ is finite-dimensional;
we may indeed consider each $V_x$ as a subspace of $\CC^S$ specified by a finite set of linear constraints which vary polynomially in $x \in \spectrum_R$. An application of Lemma \ref{linalg} shows that we may partition $\spectrum_R$ into constructible subsets so that the restriction of $\{V_x\}_x$ to each subset is a trivializable sub-bundle of the trivial bundle $\underline{\CC^S}$.

 Let us denote by $\mathcal G$ the Grassmannian of all linear subspaces of $\CC^S$.
 
Let $h_1, \dots, h_k$ generate the Hecke algebra for $J \backslash G  / J$. 
Then, by the ``eventual equivariance'' of Proposition \ref{expproperties}, $h_i V_x\subset \tilde{V}_x$, where $\tilde{V}_x$
is defined as $V_x$ but replacing $N_{\Theta}$
by certain smaller neighborhoods $\tilde{N}_{\Theta}$;
in particular, $V_x \subset \tilde{V}_x$, and it still makes sense for a subspace of $V_x$ to be ``Hecke stable.'' 

  It follows from Lemma \ref{linalg} again that the subset $Z_R\subset \spectrum_R\times\mathcal G$ defined as: 
\begin{equation}\label{zrdef} Z_R = \{ (x,M) \in \spectrum_R \times\mathcal G:  \mbox{$M \subset V_x$, and $M$ is Hecke-stable}\}
\end{equation}
 is a constructible subset of $\spectrum_R \times \mathcal G$.  

By definition, each vector in $V_x$ is ``asymptotically annihilated'' in the $\Theta$-direction 
by the ideal $I_{x, \Theta} \subset \C[\mathcal Z(X_\Theta)]$ which is determined by $x$. 
We consider the subset $Z'_R \subset Z_R$ consisting of pairs $(x, M)$ where 
where $M$ is an $\mathcal H(G,J)$-stable subspace of $V_x$ whose annihilator is {\em precisely}
$I_{x, \Theta}$. This is a constructible set, since it is obtained from $Z_R$ by removing the preimages of the $Z_{R'}$ for all $R'<R$. (Notice that we can use the same set $S$ for those, so we have a finite number of maps: $\spectrum_R \times\mathcal G \to \spectrum_{R'}\times\mathcal G$ by forgetting certain subsets of the exponents.)  

We now introduce the notion of \emph{subunitary} exponents, which are the exponents that $X$-discrete series have (see also \S \ref{sssubunitary}.) We say that a (complex) character $\chi$ of $A_{X,\Theta}$ is \emph{subunitary} { if $|\chi(a)|<1$ for all $a\in \mathring A_{X,\Theta}^+$. We say that it is \emph{strictly subunitary} if it is also subunitary when restricted to $A_{X,\Omega}$, for all $\Omega\supset\Theta$ -- i.e.\ when it does not become unitary after restriction to a ``wall'' of $A_{X,\Theta}^+$.} 
 
We say that an element of $\mathfrak S_R$ is subunitary if all its components are. We denote by $\mathfrak S_R^\su$ the semialgebraic, and open in the Hausdorff topology, subset of subunitary exponents. 

\label{newredstuff}
Let $x\in \spectrum_R^\su$ and let $f\in V_x$ be a function which generates under $\mathcal H(G,J)$ a Hecke-stable subspace of $V_x$. Then we claim that $f\in L^2(X)_\disc$. {\bluetext This is a generalization of Casselman's square integrability criterion, \cite[Theorem 4.4.6]{Casselman-notes}.} 

{ First of all, we claim that the exponents of $f$ are in fact \emph{strictly subunitary}, not just unitary, in every direction $\Theta\subsetneq\Delta_X$.
Indeed, as mentioned above, $f$ has to generate an admissible subrepresentation of $C^\infty(X)$. This implies that for every $\Theta\subset\Omega\subset\Delta_X$, its $\Omega$-exponents contain the restrictions of its $\Theta$-exponents to $A_{X,\Omega}$, and since the $\Omega$-exponents are assumed to be subunitary, for every $\Omega\supset\Theta$, it follows that the $\Theta$-exponents cannot be unitary on any ``wall'' of $A_{X,\Theta}^+$, i.e.\ they are strictly subunitary.

This implies that $f$ is in $L^2(X)$: recall from \S \ref{notation} that we are using {\em normalized} actions to define ``exponents'', wherein
we twist the action of $A_{X,\Theta}$ by the square root of the eigenmeasure, 
so that  the condition of ``strictly subunitary'' on the $A_{X,\Theta}$-exponents of every $f$ forces $f$ to be square integrable  on $N_{\Theta}$, i.e. the growth of measure on $X$ has been built into the normalization.}

Since $f\in L^2(X)$ and generates an admissible subrepresentation of $C^\infty(X)$, hence generates a subrepresentation of $L^2(X)$ which belongs to a finite sum of irreducibles, it follows that $f\in L^2(X)_\disc$.

Now we claim:
\begin{quote}
 The projection of $Z_R'$ to $\mathfrak S_R$ intersects $\mathfrak S_R^\su$ in a finite set.
\end{quote}

Indeed, let $(x,M)\in Z'_R$, with $x\in \mathfrak S_R^\su$. The space $M$ contains elements of $L^2(X)_\disc^J$ which are asymptotically annihilated by the
ideal corresponding to $x$, but not by any larger ideal.  Since $L^2(X)_{\disc}$ is a countable direct sum of irreducible subrepresentations, there are only countably many $x\in \mathfrak S_R^\su$ which admit a subspace $M$ with that property.

On the other hand, the projection of $Z'_R$ to $\mathfrak S_R $ is constructible, which implies that its intersection with $\spectrum_R^\su$ is either uncountable or finite. So this intersection is finite. Let us denote by $L^2(X)^J_R$ the subspace of $L^2(X)^J_\disc$ spanned by all such subspaces $M$.

Finally, observe:
\begin{quote}
 There is a positive integer $N$ such that every $f\in L^2(X)_\disc^J$ spanning an irreducible representation is contained in $L^2(X)_R^J$, for some $R$ with $|R|<N$,
\end{quote}
where we write $|R| = \sum_{\Theta} r_{\Theta}$ for short.

Indeed, for an embedding $\pi\to C^\infty(X)$ the asymptotics give $\pi\to C^\infty(X_\Theta)$ or, equivalently, an $L_\Theta$-equivariant map from the Jacquet module with respect to $P_\Theta^-$: $\pi_{\Theta^-}\to C^\infty(X_\Theta^L)$, where the action of $\mathcal Z(X_\Theta)$ coincides\footnote{Recall that in the wavefront case $\mathcal Z(\XX_\Theta)$ is a quotient of $\mathcal Z(\LL_\Theta)^0$, Proposition \ref{wavefrontlevi}.} with the action of $\mathcal Z(L_\Theta)^0$. The degree of an element of $\pi_{\Theta^-}$ as a finite $\mathcal Z(L_\Theta)^0$-vector is uniformly bounded as $\pi$ varies over all irreducible representations, i.e.

\label{BZbound}
{ \begin{quote} For any irreducible representation of $G$, the length (as $M$-representation) of the Jacquet functor  associated to
the parabolic subgroup $P=MN$  is bounded above by the order of the Weyl group of $G$.
\end{quote} 
 This follows from the exactness of the Jacquet functor, 
embedding $\pi$ inside a representation induced from a supercuspidal, and \label{BZref} then applying the ``geometrical lemma'' of Bernstein and Zelevinsky \cite[\S 2]{BZ}. }

Therefore, to ``detect'' exponents of discrete series we only need to work with a finite number of $R$'s. This implies that $L^2(X)_\disc^J$ is finite-dimensional.

We have until now discussed only the case where $\mathcal Z(X)$ is trivial. In the general case we will denote by $\mathfrak S_R$ the tuples of exponents which agree on $\mathcal Z(X)$, and repeat the same proof but considering, instead of $\mathfrak S_R$, its fiber $\mathfrak S_{R,\chi}$ over a given $\chi\in \widehat{\mathcal Z(X)}$. 

\end{proof}

 \subsection{Variation with the central character}  

\subsubsection{The problem}Now we will discuss the way that relative discrete series vary as their central character varies. Let us first explain through examples the difficulties that one faces:\label{newissue}

In the case of a reductive group ($\XX=\HH,\GG=\HH\times\HH$) we can twist any relative discrete series $\pi\simeq \tau\otimes\tilde\tau$ by characters of $G$ of the form $\eta\otimes\eta^{-1}$ (let us say: $\eta$ is unramified) in order to obtain a ``continuous family'' of relative discrete series. More precisely, if $M_\tau: \tau\otimes\tilde\tau\to C^\infty(X)$ denotes the ``matrix coefficient'' map, and we identify the underlying vector spaces of all the representations $\tau\otimes\eta$ (and those of their duals), then $M_\tau$ lives in the holomorphic (actually, polynomial) family of morphisms $M_{\tau\otimes\eta}$, parametrized by the complex torus $D^*$ of unramified characters of $H$. Moreover, for $\eta$ in the real subtorus $D_\iR^*$ of unitary characters these morphisms represent relative discrete series for $X=H$, and notice also that for every value of the parameter they are non-zero. There is a finite-to-one map $D^*\to\widehat{\mathcal Z(X)}_\CC$, and the Plancherel formula for the discrete 
spectrum of $X=H$ is a sum, over all such orbits $[\tau]$ of discrete series, of terms of the form:
$$ \Vert\Phi\Vert^2_{[\tau]} := \int_{D^*_\iR} C\circ \tilde M_{\tau\otimes\eta}(\Phi\otimes\bar\Phi) d\eta,$$
where $\tilde M$ denotes the adjoint of matrix coefficients, $C$ denotes the natural contraction map: $(\tau\otimes\eta)\otimes (\tilde\tau\otimes\eta^{-1}) \to \CC$,
and $d\eta$ is a suitable Haar measure on $D^*_\iR$. (Notice that it might happen here that a finite subgroup $F$ of $D^*_\iR$ stabilizes the isomorphism class of $\tau$; therefore this integral does not correspond to a direct integral decomposition of the space; for that we would have to write it as an integral over $D^*_\iR/F$.)

While in the group case (and more generally, as we shall see, in the case of a symmetric variety) one easily obtains ``continuous families'' of relative discrete series in this way, let us describe an interesting new issue which arises for general spherical varieties. Consider the action of $\GG = \GGm \times \PPGL_2$ on $\XX =\PPGL_2$,
where $\PPGL_2$ is acting on itself on the right, whereas $z \in \GGm$ acting via the left action of {\small $\left(\begin{array}{cc} z & 0 \\ 0 & 1 \end{array}\right)$.}
This is spherical; moreover, $L^2(X)$ decomposes as an integral indexed
by unitary characters of $k^{\times} = \GGm(k)$:
\begin{equation} \label{example} L^2(\XX(k)) = \int_{\widehat{k^{\times}}} L^2_{\omega} d\omega. \end{equation} 
Here $L^2_{\omega}$ is the unitary induction of $\omega$ from $\GGm(k)$
to $\PPGL_2(k)$.  It is known\footnote{Although this is only relevant as motivation, it follows from
Theorem \ref{propIIconj}, together with the fact -- immediate from the theory of the Kirillov model --  that {\em every} discrete series representation of $\PPGL_2$ is distinguished
for $(\GGm, \chi)$.  }
  that $L^2_{\omega}$ possesses discrete series  for all $\omega$;
however, {\em a priori}, it could vary wildly as $\omega$ changes. 
This leads to a difficulty in analyzing the most discrete part of the spectrum of $X$.

\subsubsection{Algebraicity and measurability}

Our goal here is to use similar arguments as in the proof of Theorem \ref{finiteds} in order to show that $L^2(X, \chi)_{\disc}^J$ varies measurably with $\chi\in \widehat{\mathcal Z(X)}$. 

Let us specify what this means. The family of spaces $H_\chi:=L^2(X, \chi)$ form a ``measurable family of Hilbert spaces parameterized by $\chi$''; by this we mean that we have a collection
of \emph{measurable sections} $\chi \mapsto f_{\chi}\in H_\chi$ satisfying certain natural axioms (cf.\ \cite{BePl}).
 This also gives rise to a notion of measurable sections: $\chi\mapsto T_\chi\in \End(H_\chi)$ as follows: $\chi\mapsto T_\chi$ is measurable if and only if for any two measurable sections $\chi\mapsto \eta_\chi,\chi\mapsto \eta'_\chi$ the inner product $\left<T_\chi\eta_\chi,\eta'_\chi\right>$ is a measurable function of $\chi$.

\begin{proposition}\label{measurable}
 The projections: $L^2(X,\chi)\to L^2(X,\chi)_\disc$ are measurable. 
\end{proposition}

\begin{proof}
The idea of the proof is to go over our previous proof of Theorem \ref{finiteds} and show that  
(suitably interpreted) the discrete spectrum in fact varies ``semi-algebraically'' with $\chi$. 

 Let us make a couple of reductions. First of all:
\begin{itemize}\item
 It suffices to show that there is a countable number of measurable sections $\chi\mapsto f_i(\chi)\in L^2(X,\chi)_\disc$ such that, for almost all $\chi$, $\{f_i(\chi)\}_i$ spans $L^2(X,\chi)_\disc$.
\end{itemize}

This is clear; one may construct the orthogonal projection in a ``measurable'' fashion from the $f_i(\chi)$, using the Gram-Schmidt process. 
 
  Now consider the spaces $V_x$, $x\in\spectrum_R$, discussed in the proof of Theorem \ref{finiteds}.  Recall (see the last paragraph of the proof of this Theorem) that we denote by $\spectrum_R$ the subset of $\prod_{\Theta} (\widehat{\mathcal{Z}(X_{\Theta})}_{\C}^J)^{r_\Theta}$ such that the restrictions of all factors to $\mathcal Z(X)$ coincide; the space $V_x$ then consists  of functions on $X$ whose ``asymptotics'' transform under the component characters of $x$, and in particular this notion depends on a choice of good neighborhoods of $\infty$. 
  
   As we remarked, there is a finite subset $S$ of $X/J$, depending only on $R$, such that the evaluation maps: $V_x\to \CC^S$ are injective, for every $x$. Hence, we consider the union of the spaces $V_x$ as a subspace of $\spectrum_R\times\CC^S$. It now suffices to show:
\begin{itemize}\item
 There is a finite number of measurable sections $$\widehat{\mathcal Z(X)}\ni \chi\mapsto f_i(\chi)\in \mathfrak S_{R,\chi}^\su\times \CC^S $$ such that:
\begin{itemize}
 \item The image lies in the union of the spaces $V_x$.
 \item For almost all $\chi$, the $f_i(\chi)$, considered as functions on $X$, actually lie inside $L^2(X,\chi)_\disc^J$ and moreover span this space.
\end{itemize}
\end{itemize}
 Indeed, if $f_i$ is measurable as a section into $\mathfrak S_{R,\chi}^\su\times \CC^S$ then it is also measurable as a section into $L^2(X,\chi)_\disc^J$; this follows by examining
 the proof (``by linear recurrence'') that the map $V_x \rightarrow \CC^S$ is injective for fixed $x$. 
 (Also, in what follows, the fact that the $f_i(\chi)$ actually lie inside $L^2(X, \chi)_{\disc}^J$ will follow as in the arguments
 on page \pageref{newredstuff}). 

Now we use the notation of the proof of Theorem \ref{finiteds}, adjusted to the present setting: Let $Z'$ be the subset of pairs $(x,M)$ such that $x \in \mathfrak S_R$ and $M$ an $\mathcal H(G,J)$-stable subspace of $V_x$, considered as a subspace of $\CC^S$,
with the property that $M$ is not asymptotically annihilated by any ideal larger than $I_{x, \Theta}$.  (``Hecke stable'' means that $V_x$ is Hecke stable when identified as a subspace
of $C^{\infty}(X)^J$, cf. discussion before \eqref{zrdef}). 
  Hence, again, $Z'$ is a constructible subset of $\mathfrak S_R\times \mathcal G$, where $\mathcal G$ denotes the Grassmannian all subspaces of $\CC^S$. 
  
  We define another constructible subset $Z''\subset Z'$ as follows: If $\mathcal G_N\subset \mathcal G$ denotes the Grassmannian of $N$-planes, a point $(x,M)\in Z'\cap (\mathfrak S_R\times \mathcal G_N)$ belongs to $Z''$ if and only if there is no $(x,M')\in Z'\cap (\mathfrak S_R\times \mathcal G_{N'})$, with $N'>N$.  In words, $Z''$ parameterizes pairs $(x, M)$
  where $M$ is of {\em maximal dimension} among Hecke-stable subspaces of $V_x$ annihilated exactly by $I_{x, \Theta}$. 
  
  Since the sum of two Hecke-stable subspaces is also Hecke stable, it is clear  that:
\begin{quote}
 The fibers of $Z''\to \mathfrak S_R$ are of size at most one.
\end{quote}

Finally, consider $Z''^\su$, the intersection of $Z''$ with $\mathfrak S_R^\su \times \mathcal G$; it is a semi-algebraic set. We have a canonical map: $Z''^\su\to\widehat{\mathcal Z(X)}$. By what we have already established in \S \ref{subsec:finiteness}, the fibers of this map are finite. 

We now utilize Hardt triviality, recalled in \S \ref{734}.  
It provides us with a finite number of semi-algebraic sections $$Y_i\ni \chi \mapsto g_i(\chi) \in \mathfrak S_{R,\chi}^\su\times\mathcal G,$$ where $Y_i\subset\widehat{\mathcal Z(X)}$ is semialgebraic.   Note that a semialgebraic set is a Polish space. It is known \cite[Appendix V]{DiVN} that a function $r: P_1 \rightarrow P_2$ from a locally compact second-countable space to a Polish space whose graph is a Polish space (more generally, a Souslin set) is Borel measurable.

  Choosing, locally on $\mathcal G$, a frame for the corresponding subspace of $\CC^S$, we get a finite number of measurable sections:
$$\widehat{\mathcal Z(X)}\ni \chi\mapsto f_i(\chi)\in \mathfrak S_{R,\chi}^\su\times \CC^S $$
(extended by zero away from the sets $Y_i$) with the property that, as $R$ varies over a finite set, their specializations span $L^2(X,\chi)_\disc^J$, for every $\chi$.
\end{proof}

\begin{corollary}
 There is a natural measurable structure on the family of Hilbert spaces $\left\{L^2(X,\chi)_\disc\right\}_{\chi\in \widehat{\mathcal Z(X)}}$, which identifies the direct integral:
\begin{equation}
 L^2(X)_\disc:=\int_{\widehat{\mathcal Z(X)}} L^2(X,\chi)_\disc d\chi
\end{equation}
with a closed subspace of $L^2(X)$.
\end{corollary}

\subsection{Toric families of relative discrete series}

\subsubsection{Factorizable spherical varieties} \label{factorizable}
 
To encode the difference between the examples that we saw in \S\ref{newissue}, namely that of the group and of $\PGL_2$ under the $k^\times\times\PGL_2$-action, we call a (homogeneous) spherical variety $\XX$ \emph{factorizable} if the Lie algebra $\mathfrak h$ of the stabilizer of a generic point on $\XX$ has a direct sum decomposition: $\mathfrak h=(\mathfrak h\cap \mathcal Z(\mathfrak g))\oplus (\mathfrak h\cap[\mathfrak g,\mathfrak g])$. We always apply this definition under the assumption that the connected center of the group $\GG$ surjects onto $\mathcal Z(\XX)$ -- for instance, the spherical variety $\SSL_2\backslash\SSL_3$ is not factorizable. We call a spherical variety $\XX$ \emph{strongly factorizable} if for every $\Theta\subset\Delta_X$ the Levi variety $\XX_\Theta^L$ is factorizable. Recall by Proposition \ref{wavefrontlevi} that in the wavefront case the action of $\mathcal Z(\XX_\Theta^L)$ is induced by the connected component of the center of $\LL_\Theta$ -- hence ``factorizable'' in this case 
means factorizable as an $\LL_\Theta$-variety. On the other hand, if $\XX$ is not wavefront then it cannot be strongly factorizable. Indeed, in that case there is a Levi variety $\XX_\Theta^L$ such that the connected center of $\LL_\Theta$ does not surject onto $\mathcal Z(\XX_\Theta^L)$; but $\varchi_{\LL_\Theta}(\XX_\Theta^L)$ has rank equal at most the rank of the quotient by which $\mathcal Z(\LL_\Theta)$ acts on $\XX_\Theta^L$, which in this case is strictly less than the rank of $\mathcal Z(\XX_\Theta^L)$.

While there are many interesting varieties which are not strongly factorizable, there are also many which are:

\begin{proposition}
If $\XX$ is symmetric, then $\XX$ is strongly factorizable. 
\end{proposition}

\begin{proof}
 First, we notice that the Levi varieties of a symmetric variety are also symmetric; more precisely, given an involution $\theta$ on $\GG$ and a $\theta$-split parabolic $\PP$ (i.e.\ a parabolic such that $\PP^\theta$ is opposite to $\PP$) the corresponding Levi variety is obtained from the restriction of $\theta$ to the Levi $\PP\cap \PP^\theta$. Therefore, it suffices to prove that every symmetric variety is factorizable. But any involution $\theta$ preserves the decomposition: $\mathfrak g=\mathcal Z(\mathfrak g)\oplus [\mathfrak g,\mathfrak g]$, and therefore the fixed subspace of $\theta$ on $\mathfrak g$ decomposes as the direct sum of $\mathcal Z(\mathfrak g)^\theta$ and $[\mathfrak g,\mathfrak g]^\theta$.
\end{proof}

{\bluetext A complete classification of strongly factorizable spherical varieties in terms of combinatorial data is given in \cite{DHS}.}

For factorizable varieties we can describe the discrete spectrum in terms of families of relative discrete series as in the group case. These families have much stronger properties than the measurability of Proposition \ref{measurable}. We will encode the basic properties of such families in what we will call \emph{toric families} of relative discrete series, and then try to extend them to the non-factorizable case.

\subsubsection{Definition of toric families}

In the discussion that follows, we shall often identify complex algebraic tori with their complex points, leaving the algebraic structure implicit. Let us introduce the following convention: By a \emph{torus of unramified characters} of the $k$-points of a connected algebraic group $\mathbf O$ we will mean the \emph{full} torus of unramified characters of a \emph{torus quotient} of $\mathbf O$. More precisely, to any subgroup
  $\Gamma\subset \varchi(\mathbf O)$ is associated a (torus) quotient $\mathbf O'$ of $\mathbf O$ with coordinate ring $k[\Gamma]$; to this quotient we associate the complex torus $D^*$ of characters of $O'$ that factor through the valuation map:
$$O'\to \Gamma^*:=\Hom(\Gamma,\mathbb Z).$$
  Hence, the complex points of $D^*$ are equal to $\Hom(\Gamma^*, \CC^\times)$ (here $\Gamma^*=\Hom(\Gamma,\Z)$). As noted above, we will be identifying such a torus with its group of $\CC$-points.

A torus of unramified characters for $O$ comes with a $\QQ$-structure, in particular an $\RR$-structure with Lie algebra canonically isomorphic to $\Hom(\Gamma,\RR)$; explicitly, the $\RR$-points of this structure are generated by characters of $O'$ the form $x\mapsto |\chi(x)|^r$ for $\chi\in \Gamma$ and $r\in \mathbb R$. We will be using the notation $D^*_\RR$ for the group of real points under this structure. 
 
There is a {\em second} canonical real structure on $D^*$, whose real points coincide with the maximal compact subgroup $D^*_\iR$ of unitary characters, and which we will call the \emph{imaginary} structure on $D^*$; explicitly, the real points of this structure are generated by characters of $O'$ of the form $x\mapsto |\chi(x)|^r$ for $\chi\in \Gamma$ and $r\in i\mathbb R$.

\begin{remark}
 By our definition, a torus of unramified characters of $\mathbf O$ is not necessarily a subgroup of the full group of unramified characters of $O$; in general, the map:
$$ D^*\to \{\textrm{unramified characters of }O\}$$
may have finite kernel because the $k$-points of $\mathbf O$ may not surject to the $k$-points of its torus quotient used to define $D^*$.
\end{remark}

\begin{definition} 
 A \emph{toric family} of relative discrete series for $X$ consists of the following data: a parabolic subgroup $P\subset G$, an irreducible unitary representation $\sigma$ of its Levi quotient $L$, a torus $D^*$ of unramified characters of $L$ and a family of morphisms, defined for almost every $\omega\in D_\iR^*$: $$M_\omega: \pi_\omega:=I_P^G(\sigma\otimes \omega)\to C^\infty(X),$$ where $I_P^G(\bullet)=\Ind_P^G(\delta_P^\frac{1}{2}\bullet)$ denotes the normalized induced representation, with the following properties:
\begin{enumerate}
\item For almost every $\omega\in D^*$ the representation $\pi_\omega$ is irreducible.
\item The morphism of complex varieties: $D^*\ni\omega \mapsto \chi_\omega\in \widehat{\mathcal Z(X)}_\CC$, obtained by taking the central character\footnote{For notational simplicity, for the rest of this section, we are assuming that the connected center of $G$ surjects onto $\mathcal Z(X)$; recall that we are assuming this to be true over the algebraic closure. If this is not true at the level of $k$-points, then $\mathcal Z(X)$ has to be replaced by the image of $\mathcal Z(G)^0$. For the purposes of this section, it would also not harm to replace $\GG$ by its quotient which acts faithfully on $\XX$; thus making $\mathcal Z(\XX) = \mathcal Z(\GG)^0$.}
of $\pi_\omega$, surjects onto one of the connected components of $\widehat{\mathcal Z(X)}_\CC$.
 
\item For $\omega\in D^*_\iR$ the  
image of $M_\omega$ lies in $L^2(X,\chi_\omega)$, and varies measurably with $\omega$.
\end{enumerate}
\end{definition}

\begin{remark}
To complete the analogy with the factorizable case, we expect the family of morphisms $M_\omega$ to extend to a rational family on all of $D^*$, in the usual sense: if we identify the restrictions of all $\pi_\omega$ to $K$ (and hence also the underlying vector spaces) as $V:=\Ind_{K\cap P}^K(\delta_P^\frac{1}{2}\sigma|_{K\cap P})$ then for all $x\in X$: $$(v\mapsto M_\omega(v)(x))\in \Hom_\CC(V,\CC(D^*)).$$ 

However, we will not need the rationality property here (except for a weak version of it which is covered by Theorem \ref{rationalds} which follows), and therefore we will not impose it or prove it in any case.
\end{remark}

We will sometimes condense notation and denote by $(\pi_\omega,M_\omega)_{\omega\in D^*}$ a toric family of relative discrete series as above. The notation conceals the fact that the morphisms are only defined for $\omega\in D_\iR^*$ (see, however, the previous remark). The condition ``for almost every'' means ``for almost every with respect to Haar measure on $D^*$'' but also, by means of the theory of reducibility of induced representations, ``in the complement of a finite number of divisors''.

Our basic goal in the rest of this section will be to prove instances of the following: 
\begin{discreteseriesconjecture}\label{dsconjecture}
 Given a spherical variety $\XX$, there is a parabolic subgroup $\PP=\LL\ltimes \UU\subset \GG$, a torus $D^*$ of unramified characters of $L$,
 and a countable collection $\{\mathcal D_i=(\pi^i_\omega,M^i_\omega)_{\omega\in D^*}\}$ of toric families of relative discrete series representations -- all associated to the \emph{same} $L$ and $D^*$, and finitely many of them containing non-zero $J$-fixed vectors, for any fixed open compact subgroup $J$ -- such that the norm of $L^2(X)_\disc$ admits a decomposition:
\begin{equation}
 \Vert\Phi\Vert_\disc^2 = \sum_i \int_{D_\iR^*} \Vert\tilde M_\omega^i(\Phi)\Vert^2 d\omega
\end{equation}
for a suitable Haar measure $d\omega$. 

Moreover, notice the canonical maps (with $A^*$ the canonical torus in $\check{G}$) 
$$D^*\to (\textrm{unramified characters of }P)^0 \rightarrow A^*$$
and 
$$ \mathcal Z(\check G_X)\hookrightarrow A_X^*\to A^*,$$
both of which\footnote{ The map  $ (\textrm{unramified characters of }P)^0 \hookrightarrow A^*$ is injective because, in fact,
we can identify  the group of unramified characters of $P$ with the center of the dual Levi subgroup $\check{L} \subset \check{G}$
canonically attached to the parabolic $P$.}  have finite kernel. The following is true:
\begin{quote}
 There is an isomorphism $D^*\simeq \mathcal Z(\check G_X)^0$ and an element $w\in \mathcal N_W(A_X^*)$ such that the diagram commutes: 
\begin{equation}\label{chargroups}\begin{CD}
 D^* @>>> A^* \\
@VVV  @VV{w}V \\
 \mathcal Z(\check G_X) @>>> A^*
\end{CD}
\end{equation}

\end{quote}

\end{discreteseriesconjecture}

Again, this is not necessarily a direct integral decomposition, as the images of the $M_\omega^i$'s could be non-orthogonal for different $i$'s and $\omega$'s corresponding to the same central character of $G$. However, it is easy to see that by an orthogonalization process we can make them orthogonal for different indices $i$, and as will follow from the proof (in the cases that we establish), the images will be orthogonal for $\omega$'s corresponding to different central characters of $L$. Of course, if there is multiplicity in the spectrum then this collection of families of RDSs is not unique, though our conjectures for the discrete spectrum in terms of Arthur parameters 
suggest that there might be a canonical way to pick mutually orthogonal toric families of relative discrete series spanning $L^2(X)_\disc$.

\subsubsection{Bounds on subunitary exponents} 

Before we discuss proofs of the Discrete Series Conjecture and related results, let us see a corollary which bounds the possible exponents by which an $X_\Omega$-discrete series can embed into $X$.

\begin{proposition}[Uniform boundedness of exponents] \label{uniformbound} 
 Let $J$ be a fixed open compact subgroup of $G$ and $\Theta, \Omega\subset \Delta_X$. Assume that $X_\Omega$ satisfies the Discrete Series Conjecture \ref{dsconjecture}.
 
 There exists a finite set $\mathcal{E}$ of homomorphisms $A_{X, \Theta} \rightarrow \mathbb{R}_+^{\times}$ so that for almost every (with respect to Plancherel measure on $X_\Omega$) relative discrete series representation $\pi$ for $X_\Omega$ { with $\pi^J \neq \{0\}$} and any morphism $M:\pi\to C^\infty(X)$ the exponents of $e_\Theta^* \circ M$ satisfy: $$|\chi| \in \mathcal{E}.$$

 In particular, there exists a constant $c < 1$ so that every exponent with $|\chi|<1$ on $\mathring{A}_{X, \Theta}^+$ satisfies:
$$ |\chi(a)| < c \textrm{ for }a \in \mathring{A}_{X, \Theta}^+.$$
\end{proposition}

\begin{proof} 
Since we are assuming that $X_\Omega$ satisfies the Discrete Series Conjecture \ref{dsconjecture}, almost every such $\pi$ is isomorphic to a subquotient of $I_P^G(\sigma \otimes \omega)$, where $\omega\in D_\iR^*$ and $(P,\sigma, D^*, (M_\omega)_\omega)$ varies in a \emph{finite} number of toric families of relative discrete series for $X_\Omega$.

Now consider the possible exponents of $I_P^G(\sigma \otimes \omega)$ along $X_{\Theta}$. Replacing $P$ by a smaller parabolic, we may assume that $\sigma$ is supercuspidal with central character $\chi$, and then the exponents of $I_P^G(\sigma \otimes \omega)$ along $\Theta$ are contained in the restrictions to $\mathcal{Z}(L_{\Theta})^0$ of the set of characters $\left\{^w(\chi\cdot \omega)\right\}$, where $w$ ranges over all elements of the Weyl group which map a Levi subgroup of $P$ into $P_\Theta$. The claim follows.
\end{proof} 

\subsubsection{A weaker result on $X$-discrete series}

By using general arguments similar to those in the proof of Theorem \ref{finiteds} and Proposition \ref{measurable}, we can easily prove a weaker result than the Discrete Series Conjecture. To formulate it, let us call ``algebraic family of representations'' any fixed vector space $V$, with a fixed, admissible action of a maximal compact subgroup $K$ of $G$, and a family of representations $\pi_\omega$ of $G$ on $V$ extending the action of $K$, such that $\omega$ varies on an algebraic variety $D^*$ and the action of $\mathcal H(G,J)$ on $V^J$ varies algebraically with $\omega$ for every open compact $J\subset K$. (That means, that for every $v\in V^J$ and $h\in \mathcal H(G,J)$ we have $\pi_\omega(h)v \in V^J\otimes \CC[D^*]$.) Such an algebraic family is obtained, for example, by starting from an admissible representation of a Levi subgroup, twisting it by characters of the Levi and parabolically inducing. 

Let such an algebraic family $(V, \pi_\omega)_{\omega\in D^*}$ be given. We assume that all members of the family are irreducible, and ask the question of which members of the family can appear as relative discrete series on $X$.

\begin{theorem}\label{rationalds}
Fix an open compact subgroup $J $. There is   a finite set of data $(Y_i, M_i)$, where  
$Y_i$ is a semi-algebraic subset of $D^*$ (considered as a real variety by restriction of scalars) and $M_i$ is a semi-algebraic family of relative discrete series $M_{i,\omega}:\pi_\omega\to C^\infty(X)^J$, $\omega\in Y_i$ with the property that for every $\omega$  
any relative discrete series $\pi_\omega\to C^\infty(X)$ is in the linear span of the $M_{i,\omega}$'s.
\end{theorem}

 We only sketch the main steps of the proof, since the arguments are similar to those encountered in the proof of Theorem \ref{finiteds}:

\begin{proof}[Sketch of proof]
The proof relies on showing that for every open compact $J\subset K$ there is a semi-algebraic subset $\mathcal A$ of the product:
$$\Hom(V^J, C(X)^J) \times D^*$$
with the property that $(M,\omega)\in \mathcal A$ if and only if $M$ is Hecke equivariant with respect to the $\pi_\omega$-action, and the image is in $L^2$-mod-center. 

Of course, since the space $\Hom(V^J, C(X)^J)$ is infinite (uncountable) dimensional, we need to fix neighborhoods of infinity and exponents in order to reduce it to finite dimensional spaces, as in the proof of Theorem \ref{finiteds}, in order to talk about a ``semialgebraic subset''.

Following that, we use again Hardt triviality in order to construct sections from semi-algebraic subsets of $D^*$ to $\mathcal A$. We can construct enough sections so that their specializations at all $\omega$ span the space of all relative discrete series from $\pi_\omega$.
\end{proof}

\subsubsection{Proof of Discrete Series Conjecture \ref{dsconjecture} in the factorizable case}
Now we focus on the easy case of the Discrete Series Conjecture \ref{dsconjecture}, and will discuss the general case in the next subsection.

Let $L=G$, let $\mathbf O$ be the quotient of $\GG^\ab$ by the image of $\HH$ (where $\HH$ denotes the connected component of the stabilizer of a point $x_0\in\XX$) and let $D^*$ be the torus of unramified characters of $\OO$. Hence, $\varchi(D^*)^*=\varchi_\GG(\XX)$. We claim, first of all, that $D^*\simeq \mathcal Z(\check G_X)^0$ canonically. Indeed, the cocharacter group of $\mathcal Z(\check G_X)^0$ is the set of all elements in $\varchi(\XX)$ which are perpendicular to the spherical coroots $\check \Delta_X$. The fact that $\XX$ is factorizable implies that $\varchi_\GG(\XX)$ and $\varchi_{\mathcal Z(\GG)}(\XX)$ have the same rank -- which is the rank of $\varchi(\XX)\cap \check \Delta_X^\perp$. Since $\varchi_\GG(\XX)\subset \varchi(\XX)\cap \check \Delta_X^\perp$, it follows that the two groups are commensurable and that the elements of $\varchi(\XX)\cap \check \Delta_X^\perp$ are characters of $\GG$. 
We claim that, in fact, they are trivial on an isotropy group $\HH$ or, equivalently, $\varchi_\GG(\XX)= \varchi(\XX)\cap \check \Delta_X^\perp$. Indeed, they have to be trivial on  the connected component $\HH^0$ of $\HH$ since they are commensurable with a group of characters which has this property. Now consider the morphism: $\XX^0:=\HH^0\backslash \GG\to \HH\backslash \GG=\XX$, an equivariant cover of degree $(\HH:\HH^0)$. We restrict this cover to the open $\BB$-orbits: $\mathring \XX^0\to \mathring \XX$, and see that $\HH/\HH^0$ acts faithfully on the fibers of this map. Therefore, every element of $\varchi(\XX)$ which extends to a character of $\GG$ has to be trivial on $\HH$. This proves the fact that $\varchi_\GG(\XX)= \varchi(\XX)\cap \check \Delta_X^\perp$ and, hence, $D^*=\mathcal Z(\check G_X)^0$. 

The choice of point $x_0\in X$ defines a map: $X\ni x\mapsto \bar x\in O$ (indeed, $x_0$ defines a morphism of algebraic varieties, and the map $X \rightarrow O$ is obtained by taking $k$-points),  and every relative discrete series $(\pi,M)$ can be twisted by elements of $D^*$ to obtain a toric family of RDS: $$(M_\omega)(v)(x):=M(v)\omega(\bar x).$$ Moreover the dimension of $D^*$ is equal to the dimension of $\mathcal Z(X)$, and therefore $D^*$ surjects onto the identity component of $\widehat{\mathcal Z(X)}_\CC$. 

Now choose a unitary central character $\chi$ and a direct sum decomposition of $L^2(X,\chi)$ in terms of relative discrete series:
$$L^2(X,\chi)=\oplus_i M^i(\pi^i).$$
This gives rise to a direct sum decomposition of $L^2(X,\chi\otimes\chi_\omega)$, for every $\omega\in D_\iR^*$:
$$L^2(X,\chi\otimes\chi_\omega)=\oplus_i (M^i_\omega)(\pi^i\otimes\omega),$$
and therefore from (\ref{centerdecomp}) we get:
$$ \Vert\Phi\Vert_\disc^2 = \sum_i \int_{D_\iR^*} \Vert\tilde M_\omega^i(\Phi)\Vert^2 d\omega$$
where $d\omega$ is a suitable Haar measure.
\qed

The next subsection will be devoted to reducing a much more general case to the factorizable one by a method called ``unfolding''.

\subsection{Unfolding}\label{ssunfolding}

\emph{In this section we take the right regular representation of any group on any space of functions to be \emph{unnormalized}; if we want to consider the normalized action, we will explicitly tensor by a character. As usual, our convention is that the group acts on the right on the given spaces and on the left on function- or measure spaces.}

\subsubsection{Introduction}
Let $\GG=\PPGL_2$, $\TT$= a nontrivial split torus in $\GG$. Let $\chi$ be a unitary character of $T$, and consider the representation $L^2(T\backslash G, \chi)$. By imitating a technique which has been used globally, in the theory of period integrals and the Rankin-Selberg method, we will show that there is a unitary, equivariant isomorphism: 
$$L^2(T\backslash G,\chi)\simeq L^2(U\backslash G, \psi),$$
where $U$ denotes a nontrivial unipotent subgroup and $\psi$ a nontrivial unitary complex character of it -- hence the right-hand-side is a Whittaker model. This is an amazing fact per se: For instance it shows that, as unitary representations, the spaces $L^2(T\backslash G,\chi)$ for all unitary $\chi$ are isomorphic! This is what will allow us to pass from a ``non-factorizable'' to a ``factorizable'' case: The spaces $L^2(T\backslash G,\chi)$ show up in the spectral decomposition of $X=G$ under the $T\times G$-action, which is not factorizable.  The unfolding technique shows that $L^2(X)\simeq L^2(k^\times\times U\backslash G,\psi)$, which is factorizable.\footnote{The meaning of factorizable here is exactly the same as in the case of the trivial line bundle: the Lie algebra of the stabilizer is a direct sum of its intersections with $[\mathfrak g,\mathfrak g]$ and with the Lie algebra of the center of $G$; however, notice that the center of $X$ is smaller when we have a non-trivial character, hence, for 
instance, in the case of the Whittaker model we do not enlarge the group to include the action of $\mathcal N(U)/U$ -- cf.\ footnote \ref{ftntbundle}.}  In the general case the passage is not to a factorizable space, but to one which is parabolically induced from a factorizable one.

We give a couple more examples: 

\begin{example}
 For $X=\SL_n\times\SL_n\backslash \GL_{2n}$, as a $k^\times\times k^\times\times \GL_{2n}$-space, we have:
$$ L^2(X) \simeq L^2(k^\times)\otimes L^2(\SL_n^{\diag}\ltimes \Mat_n \backslash \GL_{2n},\psi),$$
where $\psi$ is a non-degenerate character of $\Mat_n$ (the unipotent radical of the $n\times n$-parabolic) normalized by $\GL_n^\diag$. The space on the right is factorizable -- it is essentially the ``Shalika model''.
\end{example}

\begin{example}
 For $X= \SL_n\backslash \Sp_{2n}$ under the $k^\times\times\Sp_{2n}$-action, we have:
$$ L^2(X) \simeq L^2(k^\times)\otimes L^2(O_n\ltimes S^2(k^n) \backslash \GL_{2n},\psi),$$
where $\psi$ is a non-degenerate character of $S^2(k^n)$ (the unipotent radical of the Siegel parabolic) normalized by $O_n$. The space on the right is factorizable.
\end{example}

The examples above can be obtained by one step of ``unfolding'', i.e.\ one application of inverse Fourier transforms. The following one requires several steps (and we will encounter more of this kind later):

\begin{example}
 For $X=\SL_n\backslash \GL_n\times\GL_{n+1}$ under the $k^\times\times\GL_n\times\GL_{n+1}$ action, we have:
$$ L^2(X) \simeq L^2(k^\times) \otimes L^2(U_n\backslash \GL_n,\psi_n)\otimes L^2(U_{n+1}\backslash \GL_{n+1},\psi_{n+1}),$$
where the last two factors are Whittaker models for $\GL_n$ and $\GL_{n+1}$, respectively.
\end{example}

\subsubsection{Technicalities} To understand some technicalities that arise, replace in the example of $\TT\backslash \PPGL_2$ the group $\PPGL_2$ by $\SSL_2$, and consider the space $L^2(T\backslash \SL_2,\chi)$ for some unitary character $\chi$ of $T$. (The reader will notice that the space has not changed, only the group acting has.) Here we have to deal with the fact that the torus acts with multiple open orbits on the character group of a unipotent subgroup that it normalizes. In this case one needs to replace the space $L^2(U\backslash G, \psi)$ by $L^2((\{\pm I\}\UU\backslash\GG)(k),\chi|_{\{\pm I\}}\psi)$, which is not precisely the space of an induced representation, because $\SL_2$ does not act with a unique orbit on $(\{\pm I\}\backslash \SSL_2)(k)$. (Of course, in this case there is a larger group acting on the space, namely the $k$-points of $(\GGm\times\SSL_2)/\{\pm I\}^\diag$, but to deal with the general case we do not want to use this fact.)

Therefore one needs some careful definitions of the spaces that arise in the process of unfolding, not as representations induced from some character but as complex line bundles over the $k$-points of certain $\GG$-varieties. In the above examples, the space $\TT\backslash\GG$ will be understood as an affine bundle over $\BB\backslash\GG$, and the space $\UU\backslash\GG$ (resp.\ $\{\pm I\}\UU\backslash\GG$) as an open subvariety of the total space of a vector bundle over $\BB\backslash\GG$. The ``unfolding'' process is a Fourier transform from functions on the former to sections of a certain complex line bundle over the latter -- not the trivial line bundle, precisely because $\TT\backslash \GG$ is an affine bundle and not a vector bundle.

\subsubsection{Fourier transform on affine spaces}

For a vector space $\VV$ (considered as an algebraic variety) over $k$, Fourier transform is an isomorphism between $L^2(V)$ and $L^2(V^*)$, or $C_c^\infty(V)$ and $C_c^\infty(V^*)$, depending on a unitary complex character $\psi: k\to\CC^\times$ and assuming that the spaces are endowed with compatible Haar measures. Under the action of $\GL(V)$ the map is equivariant after we twist by the determinant character: 
$$C_c^\infty(V)\xrightarrow{\sim} C_c^\infty(V^*)\otimes |\det|^{-1}$$
or, if we want to work with unitary representations:
$$L^2(V)\otimes |\det|^\frac{1}{2} \xrightarrow{\sim} L^2(V^*)\otimes |\det|^{-\frac{1}{2}}.$$
Here by $\det$ we denote the determinant of $\GGL(\VV)$ on $\VV$, which is inverse to its determinant on $\VV^*$.

Now let $\F$ denote an affine space over $k$, that is, a variety equipped with a equivalence class of isomorphisms with $\mathbb A^n$ up to translations and linear automorphisms, for some $n$. We can still define Fourier transform of functions on $F$ as follows: First of all, let $\VV$ denote the unipotent radical of $\Aut(\F)$ (the group of affine automorphisms of $\F$); these are the ``translation'' automorphisms, i.e.\ those automorphisms which have no fixed point. The whole affine automorphism group is isomorphic to $\GGL(\VV)\ltimes\VV$, the isomorphism depending on the choice of a point on $F$. 

The dual vector space $\VV^*$ of $\VV$ can be thought of as the algebraic variety classifying reductions of $\F$ to a $\GGa$-torsor. Therefore, there is a ``universal'' principal $\GGa$-bundle over $\VV^*$, which we will denote by $\R$. It carries a canonical (right) $\Aut(\F)$-action.

Explicitly, we have:
\begin{equation}
 \R= (\F\times\VV^*\times\GGa )/\VV
\end{equation}
where $\VV$ acts on $\F\times\VV^*\times\GGa$ as:
\begin{equation}\label{Vaction}
 (f,v^*,x)\cdot u = (f\cdot u, v^*, x+\left<u,v^*\right>).
\end{equation}
The action of $\Aut(\F)$ descends to $\R$ from the following action on $\F\times\VV^*\times\GGa$:
\begin{equation}\label{Autaction}
 (f,v^*,x)\cdot g = (f\cdot g, v^*\cdot g, x).
\end{equation}
Notice that (\ref{Vaction}) and (\ref{Autaction}) do not coincide on $\VV$, i.e.\ $\VV\subset \Aut(\F)$ acts non-trivially on $\R$.

The choice of a point $f\in F$ defines a trivialization of this bundle: $\VV^*\times\GGa\ni (v^*,x)\mapsto [f,v^*,x]\in \R$. We will denote the corresponding trivialized bundle by $\R_f$. Notice that all these trivializations coincide over $0\in \VV^*$, but not over other points.

Now let $\psi: k\to \CC^\times$, a unitary character, be given. It defines a reduction of the corresponding $k$-bundle to a $\CC^\times$-bundle over $V^*$, and hence a complex line bundle over $V^*$, which we will denote by $R^\psi\to V^*$. Again, the choice of a point $f\in F$ defines a trivialization of this bundle, which we will denote by $R^\psi_f$. 

Fix throughout dual Haar measures on $F$ and $V^*$ with respect to $\psi$, where by ``Haar'' on $F$ we mean a $V$-invariant measure. Let $\mathfrak d_V$ denote the inverse of the character by which $\Aut(\VV)$ acts on an invariant volume form on $\VV$. Then $L^2(F)\otimes |\mathfrak d_V|^{-\frac{1}{2}}$ is a unitary representation of $\Aut(F)$, and notice also that $L^2(V^*, R^\psi)\otimes |\mathfrak d_V|^\frac{1}{2}$ is also well-defined and unitary, since the character $\psi$ is unitary and hence so are the isomorphisms between the different trivializations: $R^\psi_{f_1}\overset{Id}{\to}R^\psi\overset{Id}{\to}R^\psi_{f_2}$. (In other words, the line bundle $R^\psi\otimes\overline{R^\psi}$ is \emph{canonically} trivial.)

Choosing a point $f\in F$ identifies $F$ with $V$, and hence we can, using the character $\psi$, define Fourier transform:
\begin{equation}
 L^2(F)\otimes |\mathfrak d_V|^{-\frac{1}{2}} \to L^2(V^*)\otimes |\mathfrak d_V|^\frac{1}{2} = L^2(V^*, R_f^\psi)\otimes |\mathfrak d_V|^\frac{1}{2}.
\end{equation}
It is immediate to check that, independently of any choice of point:
\begin{lemma}\label{Fourierlemma}
 Fourier transform defines a canonical, $\Aut(F)$-equivariant isomorphism:
\begin{equation}
 L^2(F)\otimes |\mathfrak d_V|^{-\frac{1}{2}}\overset{\sim}{\to} L^2(V^*, R^\psi)\otimes |\mathfrak d_V|^{\frac{1}{2}}.
\end{equation}
\end{lemma}

Notice also that it preserves smooth, compactly supported sections (here it is preferable to tensor with the modular character only on one side):
\begin{equation}
 C_c^\infty(F)\overset{\sim}{\to} C_c^\infty(V^*, R^\psi)\otimes|\mathfrak d_V|.
\end{equation}

\subsubsection{Affine bundles and unfolding}

We generalize this to affine bundles: Let $\GG$ be an algebraic group over $k$, $\YY$ a homogeneous space for $\GG$ (with $\YY(k)\neq \emptyset$) and $\F$ an affine bundle over $\YY$. Let $\mathcal L$ be a complex $G$-line bundle over $Y$ (which, at first reading, the reader may take to be trivial). We want to define Fourier transform as an equivariant isomorphism:
\begin{equation}\label{Fouriertransform}
\mathcal F: C_c^\infty(F,\mathcal L)\otimes \eta^\frac{1}{2} \xrightarrow{\sim} C_c^\infty({V^*}, \mathcal L\otimes R^\psi \otimes |\mathfrak d_V|)\otimes \eta^\frac{1}{2}.
\end{equation}

We explain the notation: We denote by ${\VV^*}$ the underlying space of the vector bundle over $\YY$ whose fiber over $y\in \YY$ is the vector spaces of reductions of the corresponding fiber $\F_y$ to a $\GGa$-torsor.  The affine bundle $\F$ gives rise to a canonical principal $\GGa$-bundle $\R$ over ${\VV^*}$; the character $\psi$ defines a reduction of this to a complex line bundle $R^\psi$ over ${V^*}$. The complex line bundle $\mathcal L$ is also considered as a line bundle over $F$ or ${V^*}$ by pull-back (but denoted by the same letter). We let $\mathfrak d_V$ be the inverse of the character by which the stabilizer of a point $y\in \YY$ acts on a translation-invariant volume form on the fiber $\F_y$; its absolute value defines a complex line bundle over $Y$ (and, by pull-back, over $F$ and $V^*$), which explains the tensor product with $|\mathfrak d_V|$ on the right hand side. Finally, we \emph{assume} that there is a positive  $G$-eigenmeasure of full support on $Y$, with eigencharacter $\eta$, valued 
in $(\mathcal L\otimes \bar{\mathcal L} \otimes |\mathfrak d_V|)^{-1}$. This endows $C_c^\infty(F,\mathcal L)\otimes \eta^\frac{1}{2}$ with a $G$-invariant Hilbert norm; the corresponding Hilbert space completion will be denoted by $L^2(F,\mathcal L)\otimes \eta^\frac{1}{2}$. Recall that in this subsection only we use the \emph{unnormalized} action of the group on spaces of functions, so the twist by $\eta^\frac{1}{2}$ is the one that makes the representation unitary (and is absorbed in the normalized action for the rest of the paper).

We can now apply Fourier transform on compactly supported smooth sections, fiberwise with respect to the maps $\F\to \YY, {\VV^*}\to \YY$.  We can make a $G$-invariant choice of Haar measures\footnote{As usual, we demand that the invariant Haar measures arise from invariant differential forms; thus, the choice of a Haar measure on a fiber fully determines it on all fibers, even if $G$ does not act transitively on $Y$.} (i.e.\ translation-invariant measures) on the fibers of $F\to Y$, valued in the line bundle defined by $|\mathfrak d_V|$. Hence, if $f\in C_c^\infty(F,\mathcal L)$, then $|f|^2\in C_c^\infty(F,\mathcal L\otimes \bar{\mathcal L})$ and we can define an ``integration against Haar measure on the fiber'' morphism:
\begin{equation}
 C_c^\infty(F,\mathcal L\otimes\bar{\mathcal L})\to C_c^\infty (Y, \mathcal L\otimes \bar{\mathcal L} \otimes |\mathfrak d_V|)
\end{equation}

We claim that $C_c^\infty(V^*, \mathcal L\otimes R^\psi \otimes |\mathfrak d_V|)\otimes \eta^\frac{1}{2}$ also has a unitary structure. Indeed, for $f\in C_c^\infty(V^*, \mathcal L\otimes R^\psi \otimes |\mathfrak d_V|)$ we have $|f|^2\in C_c^\infty(V^*, \mathcal L\otimes\bar{\mathcal L} \otimes |\mathfrak d_V|^2)$, where we took into account the fact that $R^\psi\otimes\overline{R^\psi}$ is the trivial line bundle \emph{canonically}. Notice also that stabilizers of points act on the Haar measure on fibers of $\VV^*\to \YY$ by $\mathfrak d_V$. Therefore, integration along the fibers with respect to Haar measure gives:
\begin{equation}
C_c^\infty(V^*, \mathcal L\otimes\bar{\mathcal L} \otimes |\mathfrak d_V|^2) \to C_c^\infty(Y, \mathcal L\otimes\bar{\mathcal L}\otimes|\mathfrak d_V|), 
\end{equation}
i.e.\ \emph{the same line bundle whose dual we assumed admits a $G$-eigenmeasure with eigencharacter $\eta$}.

Having fixed these measures, it is now easy to see that fiberwise Fourier transform indeed represents an equivariant isomorphism as in (\ref{Fouriertransform}). The essence of what we call ``the unfolding step'' is the following:

\begin{theorem}\label{unfoldingisomorphism}
 With suitable choices of invariant Haar measures on the fibers of $F\to Y, {V^*}\to Y$, Fourier transform (\ref{Fouriertransform}) gives rise to a $G$-equivariant isometry of unitary representations:
\begin{equation}
\mathcal F: L^2(F,\mathcal L)\otimes \eta^\frac{1}{2} \xrightarrow{\sim} L^2({V^*}, \mathcal L\otimes R^\psi \otimes |\mathfrak d_V|)\otimes \eta^\frac{1}{2}.
\end{equation}
\end{theorem}

\begin{proof}
 Let $\mathcal F'$ denote the inverse of $\mathcal F$. Since the spaces of compactly supported smooth sections are dense in the corresponding $L^2$-spaces, and are mapped to each other under Fourier transform, the statement is equivalent to the following:
\begin{quote}
 The morphisms $\mathcal F, \mathcal F'$ are adjoint.
\end{quote}
 Since the morphisms are defined fiberwise and the $L^2$-hermitian forms are computed by ``integration on the fibers'' followed by (the same) ``integration over $Y$'', it suffices to prove this fact fiberwise, where ``adjoint'' will mean ``adjoint with respect to the $L^2$(Haar measure)-pairing on the fiber. But this is just Lemma \ref{Fourierlemma} (ignore tensoring by the characters, since the equivariance of the map has already been established). It is easy to see that there are compatible choices of Haar measures globally.
\end{proof}

\subsubsection{The goal of unfolding}

The method of unfolding provides a proof of Conjecture \ref{dsconjecture} in all cases that we have examined. However, we have not been able to establish the inductive step abstractly. The idea, in general, is to use this technique in order to prove the following conjecture. Recall that in this subsection we are not normalizing implicitly the right regular representations, so we denote by $\eta$ the eigencharacter of our measure on $X$ and tensor with its square root to make $L^2(X)$ unitary.
\begin{conjecture}\label{unfoldingconjecture}
 Given a homogeneous spherical variety $\XX$ there exists a Levi subgroup $\LL$ of a parabolic $\PP$ of $\GG$, a factorizable spherical variety $\WW$ of $\LL$ and a complex hermitian line bundle $\mathcal L_\Psi$ with an $L$-action over $W$, together with an isometric isomorphism of unitary representations (depending only on choices of measures on the various spaces):
\begin{equation}
L^2(X) \otimes\eta^\frac{1}{2}\xrightarrow{\sim} I_P^G \left( L^2(W,\mathcal L_\Psi)\right)\otimes\eta^\frac{1}{2}.
\end{equation}
\end{conjecture}

A proof of this conjecture proves the first part of the Discrete Series Conjecture \ref{dsconjecture}. The second part, regarding character groups, can easily be verified in each case, but we have not shown abstractly that the unfolding process establishes it.

\subsubsection{The inductive step} The proof of Conjecture \ref{unfoldingconjecture} in each particular case goes through a series of ``unfolding'' steps.

Let a homogeneous spherical variety $\XX=\HH\backslash \GG$ be given, where we keep assuming that $\mathcal Z(\GG)^0\twoheadrightarrow \mathcal Z(\XX)$, and let $\HH=\MM\ltimes\NN$ be a Levi decomposition. If $\XX$ is not factorizable or parabolically induced from a factorizable spherical variety then there  exists\footnote{The proof of existence of such a unipotent subgroup departs from the observation that there is a proper parabolic subgroup containing $\HH$; we omit the details since we do not know how to describe a canonical choice for $\UU$.} a unipotent subgroup $\UU$ of $\GG$, normalized by $\HH$ and such that $\NN\cap\UU$ is normal in $\UU$ and $\UU/\NN\cap\UU$ is non-trivial and abelian. 

Consider the variety $\YY=\MM\NN\UU\backslash\GG$. We will denote by $y$ the point of $\MM\NN\UU\backslash \GG$ represented by ``1''. Then $\XX\to\YY$ has a canonical structure of an affine bundle; indeed, the orbit map for the action of $\UU$ on a point over the fiber of $\MM\NN\UU$ identifies that orbit with the additive group $\UU/\NN\cap\UU$, the choice of point changes this identification by translations and the action of $\MM\NN$ is linear on $\UU/\NN\cap\UU$.
Call $\VV^*$ the corresponding dual vector bundle,  and $R^\psi$ the complex line bundle over $V^*$, as above. Applying Fourier transform, we get according to Theorem \ref{unfoldingisomorphism} an equivariant isomorphism:
\begin{equation}\label{BART}
 L^2(X)\otimes\eta^\frac{1}{2} \to L^2(V^*, R^\psi \otimes |\mathfrak d_V|)\otimes\eta^\frac{1}{2}.
\end{equation}

This is the first step of the unfolding process. We expect:
{
 \begin{conjecture}
$\VV^*$ is a spherical $\GG$-variety. 
\end{conjecture}

The $L^2$-isomorphism of \eqref{BART}, and the implied finiteness of multiplicities, should also imply the spherical property of the above conjecture. Clearly, however, a direct, geometric proof would be more desirable. In any case, since we do not have a complete recipe for proving Conjecture \ref{unfoldingconjecture}, one needs to check this in any specific case, which is easy to do.

In any case, the fact that $\XX$ is spherical implies that $\MM$ acts with an open orbit on $\XX_y=\UU/(\UU\cap\NN)$ (the fiber of $\XX$ over the point $y\in Y$).\footnote{ Indeed, if $\widetilde\PP$ is a parabolic with unipotent radical $\widetilde\UU$ containing $\NN\UU$, and with $\MM$ in its Levi $\LL$, then since $\MM\NN$ acts with an open orbit on the flag variety of $\GG$, it has to act with an open orbit on the open $\tilde\PP$-orbit of the flag variety, which is isomorphic to $\LL/\BB_L \times \widetilde\UU$. Equivalently, $\MM$ acts with an open orbit on $\LL/\BB_L \times \widetilde\UU/\NN$, and since $\widetilde\UU/\NN$ is fibered over $\widetilde\UU/\UU\NN$ with fiber $\UU/\NN\cap\UU$, it follows that $\MM$ acts with an open orbit on the latter.} This implies, in particular, that $\GG$ acts with an open orbit $\mathring \VV^*$ on $\VV^*$; this open orbit is isomorphic to $(\MM\NN)_{v^*}\UU\backslash \GG$, where $(\HH)_{v^*}$ is the stabilizer in $\HH$ of a generic element $v\in V^*_y$ (which is the dual of the ``translation'' automorphism group of $X_y$). 
}

Recall that $R^\psi$ is the complex line bundle over $V^*$ where $H_{v^*}U$ acts on the fiber over $v^*\in V^*_y$ by the character $\Psi:=\psi\circ v^*$ (considered as a character of the whole group $H_{v^*}U$), and the same description holds for the other points of $V^*$. In the case when $G$ acts transitively on $\mathring V^*$ we can write:
$$L^2(V^*, R^\psi \otimes |\mathfrak d_V|) \otimes\eta^\frac{1}{2}= L^2(H_{v^*}U\backslash G, |\mathfrak d_V|\Psi)\otimes\eta^\frac{1}{2}.$$

It is easy to see the following:
\begin{lemma} 
 Fix a point $v^*\in (V^*)_y$ (the fiber of $V^*$ over the point $y\in Y$), and the trivialization of the fiber of $R_\psi$ over $v^*$ corresponding to a chosen point $x\in X_y$. (Recall that the choice of a point $x\in X_y$ makes the fiber $X_y$ into a vector space, and hence makes Fourier transform scalar-valued.) Let $\Psi$ denote the character by which the stabilizer $H_{v^*}U$ acts on the fiber of $R_\psi$. Then, for suitable choices of measures, the map:
 $$C_c^\infty(X)\otimes\eta^\frac{1}{2} \to C^\infty(V^*, R^\psi \otimes |\mathfrak d_V|)\otimes\eta^\frac{1}{2} \xrightarrow{\ev_{v^*}} \CC$$
 inducing the above $L^2$-isomorphism is given by the integral:
$$\Phi\mapsto \int_{U/U\cap N} \Phi(u) \Psi^{-1}(u) du $$
 and its adjoint:
 $$C_c^\infty(V^*, R^\psi \otimes |\mathfrak d_V|)\otimes\eta^\frac{1}{2} \to C^\infty(X)\otimes\eta^\frac{1}{2} \xrightarrow{\ev_{x}} \CC$$
is given by the integral:
$$\Phi \mapsto \int_{(\HH_{v^*}\backslash \HH)(k)} \ev_{v^*}(R_h\Phi) dh, $$
where $R_h$ denotes the regular representation.
\end{lemma}

\begin{proof}
{ We follow the definitions and constructions that we have presented:
The first integral represents, simply, Fourier transform over the fiber ($\simeq (\UU/\UU\cap \NN)(k)=U/U\cap N$) of the affine bundle $X$($=F$ in the notation of Theorem \ref{unfoldingisomorphism}) over $y$. The fact that the second is its adjoint is completely formal; for simplicity, when $G$ acts with a unique orbit on $\mathring V^*$, and using the isomorphism $(H_{v^*}U\cap H) \backslash H_{v^*}U = U/U\cap N$, when $f\in C_c^\infty(X)$ and $\Phi\in C_c^\infty(V^*, R^\psi \otimes |\mathfrak d_V|) = C_c^\infty (H_{v^*}U\backslash G, |\mathfrak d_V|\Psi)$:

$$ \int_{H_{v^*}U\backslash G} \left(\int_{U/U\cap N} f(u g) \Psi^{-1}(u) du\right) \Phi(g) dg = \int_{(H_{v^*}U\cap H)\backslash G} f(g) \Phi(g) dg = $$
$$ =\int_{H\backslash G} \left(\int_{H_{v^*}\backslash H} \Phi(hg) dh \right) f(g) dg.$$
} 
\end{proof}

If the stabilizer $\HH_{v^*}$ of a generic point on the dual of $U/(N\cap U)$, modulo the center of $\GG$, has finite character group, then we are done with the proof of Conjecture \ref{unfoldingconjecture} in the given case: the variety to which we have unfolded is factorizable.

If not, it is convenient (for the purpose of this theoretical presentation) to ``fold back'' in order to eliminate the character $\psi$, i.e.\ to do the following:\footnote{Without choosing a base point $v^*$ on $V^*$, the ``folding back'' process admits the following description: We have obtained sections of a certain complex bundle $R^\psi$ over the $k$-points of a vector bundle $\VV^*$ over $\YY$. Pull back the $\GGa$-bundle $\mathbf R^\psi$ (and the corresponding complex vector bundle $R^\psi$) to the blow-up $\mathbb B\VV^*$ of $\VV^*$ along the zero section. Let $\YY':=\mathbb P\VV^*$ (the projectivization of $\VV^*$); then $\mathbb B\VV^*$ is a line bundle over $\mathbb P\VV^*$. There is a $\GGa$-bundle $\XX'$ over $\YY'$ such that the bundle $R^\psi$ over $\mathbb B\VV^*$ is the one obtained by Fourier transforms from the trivial complex bundle over $\XX'$. } if $\widetilde{(\MM\NN)_{v^*}}$ denotes the stabilizer of the line of $v^*$ (necessarily $\widetilde{(\MM\NN)_{v^*}}/{(\MM\NN)_{v^*}}\simeq \
\GGm$ since $\GG$ acts almost transitively on $\VV^*$), and $\UU_0$ denotes the kernel of $v^*$, let $\XX'= \widetilde{(\MM\NN)_{v^*}}\UU_0\backslash \GG$. Now we can ``fold back'' to $\XX'$, which brings us to the original setup of a spherical variety for $\GG$, without the line bundle defined by an additive character. What we cannot show in general is that the variety $\XX'$ is closer to being ``factorizable'' than the original variety $\XX$; for instance, that its stabilizers have a larger unipotent part. However, all of the examples that we have examined show that with correct choices this is indeed the case.

\begin{example}
It is instructive at this point to discuss the unfolding process for the variety $\XX=\GGa^2 \backslash \SSO_5$ under the $\GGm^2\times \SSO_5$-action. Here $\GGm^2$ is a Cartan subgroup and $\GGa^2$ is the subgroup containing the root subspaces of the two long roots. We want to show:
$$ L^2(X)\simeq L^2((k^\times)^2)\otimes L^2(U\backslash \SO_5,\psi),$$
the last factor being a Whittaker model.

One could choose for the first step the group $\UU$ to be the one corresponding to the sum of the long root spaces and one short root space, that is: the unipotent radical of a parabolic with Levi of type $\GGm\times \SSO_3$. The reader will see that ``unfolding'' this way will lead to a non-factorizable situation which we cannot unfold further (more precisely, ``folding back'' as was suggested right above we return to the original space).

On the other hand, one may take for $\UU$ the unipotent subgroup corresponding to the sum of the short root spaces and one long root space, that is: the unipotent radical of a parabolic with Levi of type $\GGL_2$. The second step now goes through with $\UU$= a maximal unipotent subgroup, and leads to the Whittaker model.
\end{example}

We finish by mentioning a few more examples.

\begin{example}
 Generalizing the example of $\PPGL_2$ as a $\GGm\times\PPGL_2$-variety that we mentioned earlier, let $\XX=\SSL_n\backslash \GGL_{n+1}$ as a $\GG=\GGm\times \GGL_{n+1}$-space, where $\GGm=\GGL_n/\SSL_n$, where $\GGL_n$ belongs to the mirabolic subgroup of $\GGL_{n+1}$ (the stabilizer of a non-zero point under the standard representation). Let $\NN$ be the unipotent radical of the mirabolic, and let $\Psi$ be a non-trivial complex character of $N$. Then we have a $G$-equivariant isometry (for $n\ge 2$):
$$ L^2(\SL_n\backslash\GL_{n+1})\simeq L^2(P_n\ltimes N\backslash\GL_{n+1},\Psi)$$
where $P_n\subset \SL_n$ denotes the stabilizer of $\Psi$. 

Here (for $n\ge 2$) the group $k^\times$ does not act trivially on $L^2(P_n\ltimes N\backslash\GL_{n+1},\delta_N\Psi)$. However, the variety $P_n\ltimes N\backslash\PGL_{n+1}$ is parabolically induced (and the character $\Psi$ is trivial on the unipotent radical of its parabolic), hence:
$$L^2(P_n\ltimes N\backslash\GL_{n+1},\Psi)=\Ind_P^G L^2((\SL_{n-1}\times N_2)\backslash (\GL_{n-1}\times\GL_2),\Psi) =$$ $$= \Ind_P^G \left( L^2(k^\times) \otimes L^2(N_2\backslash \GL_2,\Psi)\right).$$
Here $P$ denotes the parabolic of type $\Gm\times \GL_{n-1}\times \GL_2$, $N_2$ denotes a non-trivial unipotent subgroup of $\GL_2$, and $\Ind_P^G$ denotes unitary induction.
\end{example}

\begin{example}
 Let $G=\GL_{2n}$, let $\Sp_{2n}$ denote a symplectic subgroup of $G$ and let $H$ be the subgroup of $\Sp_{2n}$ stabilizing a point under the standard representation of $G$. Then by successive applications of ``unfolding'' one can show that:
$$L^2(H\backslash G)= \Ind_P^G L^2(N\times N\backslash \GL_n\times GL_n,\Psi)$$
where $P$ is the parabolic of type $\GL_n\times\GL_n$, $N\times N$ is a maximal unipotent subgroup in its Levi and $\Psi$ is a non-degenerate character of this subgroup.
\end{example}

\section{Preliminaries to the Bernstein morphisms: ``linear algebra''} \label{sec:linearalgebra}
This section collects some simple results in ``linear algebra'' (interpreted broadly) which will be used  in 
\S \ref{sec:Bernstein}. 

The reader may wish to refer to the contents only as needed. The main purpose of this section is to separate the ``abstract'' parts of arguments from the  parts that are specific to spherical varieties. 

\begin{itemize}
\item[-] The first sections (\S \ref{sssubunitary} -- \S \ref{sshermitian}) pertain to the following general question: given a Hermitian form on a vector space $V$,
and a group $S$ acting on $V$, how can one canonically replace the form by an $S$-invariant one? Assuming that the forms are $S$-finite, i.e.\ generate a finite-dimensional vector space under the action of $S$, this is possible. These constructions will be
used in \S \ref{sec:Bernstein}.

\item[-] In section \ref{meas-eigenpproj} we show that, given a family of $S$-finite linear functionals (or hermitian forms) on $V$ which vary in a measurable way over a parameter space, we may extract their eigenprojections to certain generalized eigencharacters (for instance, unitary ones) and still get a measurable family.

\end{itemize}

\subsection{Basic definitions} \label{sssubunitary}
Suppose $S$ is a finitely generated abelian group 
together with a finitely generated submonoid $S^+ \subset S$ that generates $S$. 
 
Thus there is a surjective homomorphism $\mathbb{Z}^k \twoheadrightarrow S$ so that $S^+$
is the image of $\mathbb{Z}^k_{\geq 0}:=\{\mathbf{x}: x_i \geq 0\}$.

 By a locally finite $S$-vector space $V$ we shall mean a vector space $V$ over $\C$ equipped with a locally finite action of $S$ (i.e.\ the $S$-span of any vector is finite dimensional).
 
For a vector $v\in V$, an \emph{exponent} of $v$ is any generalized eigencharacter of $S$ on the space of translates; the \emph{degree} of $v$ is the dimension of $\langle Sv \rangle$.

If $\chi$ is a character of $S$, we write $|\chi| \leq c$ (resp.\ $|\chi|<c$) if $|\chi(s)| \leq c$ (resp. $|\chi(s)| < c$) for all $s \in S^+ \smallsetminus \{0\}$;
similarly, we define $|\chi| \geq c$, $|\chi| > c$. Note that if $|\chi| = 1$ (i.e., $|\chi| \leq 1$ and $|\chi| \geq 1$) 
then $\chi$ is unitary. If $|\chi|<1$ we will say that $\chi$ is \emph{strictly subunitary}. {(We will say simply \emph{subunitary} when $|\chi|$ is $<1$ in the \emph{interior} of $S^+$, i.e.\ its elements which do not lie on the boundary of the cone spanned by $S^+$ in $S\otimes \RR$. In the rank-one case, where this notion will be used in later sections, subunitary and strictly subunitary coincide.)}

We warn that $|\chi| < 1$ is a stronger condition than ($|\chi|\le 1$ and $|\chi|\ne  1$). 
Indeed the statement $|\chi| < 1$ amounts to asking 
that $\chi$ ``decay in all directions'', and not merely in {\em some} directions. 
In practice, we will always be able to arrange this by shrinking $S^+$ if necessary. 
 
\label{simplefact} We often use the following observation: any $S^+$-stable subspace of a locally finite $S$-vector space $V$ is also $S$-stable. Indeed, the $S$-span of a vector being finite dimensional implies that the inverse of any invertible operator is a polynomial in the operator.

 \subsection{Finite and polynomial functions} 
\label{fpolsubsec}
Now we specialize to the case of functions on $S$; a function whose $S$-translates span a finite dimensional vector space will be called a \emph{finite function}. A finite function whose only exponent is the trivial character of $S$ is called a {\em polynomial.}
This coincides with the usual use of ``polynomial'' when $S=\mathbb Z^k$.
For any any finite function $f$, there exists characters $\chi_i$ and polynomials $P_i$ so that
$f(s) = \sum \chi_i(s) P_i(s)$. 

For every finite function $f$, we refer to the dimension of the space spanned by its $S$-translates 
as its {\em degree}. 

\begin{lemma} \label{polbdd}
A polynomial function that is bounded on $S^+$ is constant. 
\end{lemma}

\begin{proof} It is enough to consider the case $S = \mathbb{Z}^k, S^+ = \Z_{\geq 0}^k$. Our assertion 
reduces to the following: if a polynomial function on $\mathbb{R}^k$ is bounded
on $\Z_{\geq 0}^k$, then it is constant. \end{proof}

 \begin{lemma} \label{ffbound2}
Let $f$ be any finite function bounded on $S^+$. Then there
exists a unique $S$-invariant functional $\langle S f  \rangle \rightarrow \C$ that sends the constant function $1$ to $1$,
where $\langle S f \rangle$ is the span of all translates of $f$ by $S$. 
\end{lemma}

We refer to this functional as $\limS$. For instance, if $S = \mathbb{Z}$, $S^+ = \Z_{\geq 0}$, $t$ a nonzero complex number of absolute value $\le 1$ (not equal to $1$), and $f$ is the function $n \mapsto 3+ t^n$, then $\limS f = 3$. 

\begin{proof}
 If $f$ is any finite function, bounded on $S^+$, we may write it as a sum of generalized eigenfunctions, $f =\sum f_{\chi}$, where each $f_{\chi}$ belongs to the $S^+$-span of $f$ and is therefore itself bounded on $S^+$.  The putative functional must (by invariance) send $f_{\chi}$ to zero for $\chi \neq 1$. On the other hand $f_1$ is bounded polynomial and thus constant.  Therefore, the only possibility for the functional is
$$\sum f_{\chi} \mapsto [f_1],$$
where $[f_1]$ is the constant value of $f_1$. It is clear that this functional is $S$-invariant and has the desired normalization. 
\end{proof}

{ 
We refer to a sequence of positive probability measures $\nu_i$ on $S^+$, defined for all sufficiently large positive integers $i$, as an {\em averaging sequence}\label{averaging sequence} if
it is obtained in the following way:  Let $\ell_1, \ell_2$ be  monotone increasing affine functions $\mathbb{Z} \rightarrow \mathbb{Z}$  with $\ell_2 -\ell_1 \rightarrow \infty$ as $i \rightarrow \infty$;
for example, $\ell_1(i)=i, \ell_2(i) =2i$. 
Choose a surjection $\mathbb{Z}^k \rightarrow S$ mapping $\Z_{\geq 0}^k$ onto $S^+$,  and let $\nu_i$ be the image of the uniform probability measure on $[\ell_1(i),  \ell_2(i)]^k$. 

In particular, such a sequence has the following properties: 
 For arbitrary $s \in S$, the measure $s\nu_i$ is eventually (i.e.\ for any large enough $i$, depending on $s$) supported on $S^+$   and  the total mass of the difference $|s\cdot \nu_i-\nu_i|(S)$ approaches $0$.   }

 \begin{lemma}  \label{limandaveraging}
Let $f$ be a finite function and $\nu_i$ an averaging sequence. 
\begin{enumerate}
\item If $f$ bounded on $S^+$, we have $\int f \nu_i \longrightarrow \limS f$. 
 \item If $f$ is unbounded on $S^+$ and all exponents $\chi$ of $f$ satisfy 
 either $|\chi| = 1$ or $|\chi| < 1$,   then $\int |f|^2 \nu_i \longrightarrow \infty$. 
 \end{enumerate}
\end{lemma}

The restriction on exponents in the latter part is for simplicity, to avoid situations where $f$ grows in some directions and decays in others. This restriction will be satisfied in our applications. 

\begin{remark}
 The incongruence between the two statements (the first for $f$ and the second for $|f|^2$) will not appear in our applications, as we will apply both to the functions obtained by evaluating $S$-finite hermitian forms on translates of a vector in an $S$-vector space; evidently, however, we need a positivity assumption for the second statement to hold.
\end{remark}

\begin{proof}
Consider first the case of $f$ bounded. Note that $\int f \nu_i$ are bounded. We may choose a subsequence of the $i$ so that
all the integrals $\int g \nu_i$ converge for $g \in \langle S f\rangle$. Then 
$g \mapsto \lim_{i} \int g \nu_i$
defines an $S$-invariant functional on $\langle S f \rangle$, which is necessarily 
$\lim_{S^+}$.  Since the subsequence was arbitrary, the result follows.

Now suppose $f$ is unbounded; write $f= \sum f_{\chi}$, where $f_{\chi}$ has generalized character $\chi$.
At least one $f_{\chi}$ is unbounded. Therefore $|\chi| = 1$, since if $|\chi| < 1$
then certainly $f_{\chi}$ must be bounded.  
   Twisting by $\chi^{-1}$, we may suppose that $\chi=1$, i.e.\ $f_1$ is a polynomial.

     There exists an element $\Delta$ of the group algebra $\C[S]$ so that $\Delta \star f= f_1$. 
    Write $\Delta = \sum a_s s$. Let $\nu_i'$ be another averaging sequence, to be chosen momentarily.
    Let $\nu_i^* = \sum_{a_s \neq 0} s^{-1} \cdot  \nu_i'$, a sum of translates of $\nu_i'$. We  can and do choose $\nu_i'$
    in such a way that $\nu_i^* \leq  C \cdot \nu_i$  for some positive constant $C$ -- choose $\nu_i'$ by replacing the linear forms $\ell_1, \ell_2$  used to define $\nu_i$ by $\ell_1'=\ell_1+A, \ell_2'=\ell_2-A $ for a large enough integer $A$.   The integral  $ \int |f_1|^2 \nu_i'$  is a sum of terms of the form $\int (s f) \overline{(s' f)} d\nu_i'$, 
 and by Cauchy-Schwarz we may bound this by $\int |f|^2 \nu_i^*$. That is to say,  
$$\int |f_1|^2 \nu_i' \leq \mathrm{const} \cdot \int |f|^2 \nu_i^*.$$
Visibly, the left-hand side is unbounded (since $f_1$ is a polynomial). Since $\nu_i^* \leq  C \nu_i$  our assertion is proved. 
\end{proof}

\begin{lemma} \label{boundS}
Let $F(n,c)$ be the space of all finite functions on $S$ of degree $\leq n$ all of whose
exponents satisfy $|\chi| < c$, for some $c < 1$.  Then we can find a finite subset $\Lambda \subset S^+$, depending only on $(n,c)$, so that for all $f\in F(n,c)$:
$$\sup_{k \in S^+} |f(k)| \leq   \max_{\lambda \in \Lambda} |f(\lambda)|,$$
In fact, there exists a decaying   function $Q: S^+ \rightarrow \mathbb{R}$  depending only on $(n,c)$
so that, for any $f\in F$: \begin{equation} 
\label{finebound} |f(k)| \leq Q(k)  \max_{\lambda \in \Lambda} |f(\lambda)| . \end{equation}
\end{lemma}

\begin{remark}
Recall that the notion of ``decaying function'' was defined in the introduction \S \ref{notation}. By enlarging $\Lambda$ we may assume that $|Q|\le 1$. This way, the first statement becomes a special case of the second.
\end{remark}

\begin{proof}
We may again assume that $S = \mathbb{Z}^k, S^+ = \Z_{\geq 0}^k$. 
Let $s_1 = (1, 0, \dots, 0), \dots, s_k = (0,0,\dots, 1)$ be the standard generators for $S^+$. 
Fix an $f\in F(n,c)$, and let $\mathbf{P} = (P_1, \dots, P_k)$ be the characteristic polynomials of $s_1, \dots, s_k$ acting on $\langle Sf \rangle$, each of degree $n_i\le n$.  Note that all the coefficients of all $P_i$ are bounded in terms of $(n,c)$. 
Let $F(\mathbf{P})$ be the set of functions annihilated by $P_i(s_i)$ for each $i$.

 Put $\Lambda = \prod_{i=1}^k [0,n_i-1]$. 
 The evaluation map $\mathrm{ev}: F(\mathbf{P}) \rightarrow \C^{\Lambda}$ is a linear isomorphism. 
The action of translation by $s_i$ on $F(\mathbf P) \simeq \C^{\Lambda}$ is expressed by a certain
endomorphism $A_i \in \End(\C^{\Lambda})$ {\em whose matrix entries are bounded
in terms of $(n,c)$.} Therefore:
$$f(t_1, \dots, t_k) = \left( A_1^{t_1} \dots A_k^{t_k}  \mathrm{ev}(f) \right) (\mathbf{0}).$$

The second statement now follows from the following: Suppose that $\Omega$ is a compact subset of $\End(\CC^\Lambda)$ so that, for every $A \in \Omega$,
all of the eigenvalues of $A$ are $\leq c<1$; then there exists $N$ so that $\|A^N\| \leq \frac{1}{2}$ for all $A \in \Omega$. Here, $\|\bullet\|$ denotes any norm on $\End(\CC^\Lambda)$. To check that, take $1>c' > c$, use $A^k = \frac{1}{2\pi i } \int_{|z| = c'} \frac{z^k}{z-A} dz$ and
the fact that $$ \| (z-A)^{-1} \|, \ \ |z|=c', A \in \Omega,$$ being a continuous function on a compact set, is bounded. As mentioned in the remark, this implies the first assertion, as well (by enlarging $\Lambda$).
 \end{proof} 

\begin{lemma} \label{boundU}
Suppose $f$ is a  finite function on $S$ of degree $\leq n$, bounded on $S^+$, 
all of whose exponents are unitary. Then 
$$ \sup_{k \in S^+} |f|^2 \leq   n \limS |f|^2.$$
\end{lemma}
\begin{proof}
The boundedness forces $f$ to be a sum of {\em eigenfunctions} (with unitary character, by assumption): if we write $f = \sum_{\chi \in I} f_{\chi}$, then each $f_{\chi}$ is bounded.
We claim that $f_{\chi}$ is proportional to $\chi$. For $\chi=1$ that is 
  Lemma \ref{polbdd} and in the general case $\chi=\chi_0$ we apply the same reasoning to $f \chi_0^{-1}$.
 
The result easily follows. 
\end{proof}

\subsection{Hermitian forms}\label{sshermitian} We now turn to the properties of Hermitian forms on $S$-vector spaces. 
In what follows, where we speak of ``Hermitian forms'' we always mean  \textbf{positive semi-definite} Hermitian forms.
If $(u,v)\mapsto H(u,v)$ is a Hermitian form, we use the notation $H(v)$ for $H(v,v)$.

Notice that if $H$ is a hermitian form on a finite-dimensional $S$-vector space $V$, then it can be considered as an element of the tensor product representation $V^*\otimes \bar V^*$ of $S$, and hence the form itself is an $S$-finite vector in an $S$-vector space; in particular, it makes sense to talk about its exponents. For a hermitian form on a locally finite, possibly infinite-dimensional, $S$-vector space we call exponents of $H$ the union of its exponents on all $S$-stable, finite-dimensional subspaces. (The form itself might not be $S$-finite, in this case.)

Let $H$ be a Hermitian form on a locally finite $S$-vector space, and put $$ V_f = \{ v\in V:  H(S^+ v) \mbox{ is bounded.}\}$$
This is an $S$-invariant subspace of $V$.  That it is a subspace is a consequence of the inequality
 $$H(x_1 + \dots + x_m) \leq m \sum H(x_i),$$
 whereas the $S$-stability follows from the observation at the end of \S \ref{simplefact}.

The following is an obvious application of Lemma \ref{ffbound2}, with positivity following, for instance, from Lemma \ref{limandaveraging}:
\begin{lemma} \label{invariantform} Let $H$
be a Hermitian form on a locally finite $S$-vector space $V$. 
Then 
$$H^S(v) := \limS \left( s \mapsto H(sv) \right),$$
defines an $S$-invariant Hermitian form on $V_f$.  
\end{lemma}

We refer to $H^S$, extended by $\infty$
off $V_f$, as the {\em associated $S$-invariant form.}  Then the association $H \mapsto H^S$ is linear, i.e., given Hermitian forms $H_1, H_2$ and positive scalars $a_1, a_2$, we have
$(a_1 H_1 + a_2 H_2)^S = \sum a_i H_i^S$. The following is a corollary to what we have already proved in Lemma \ref{limandaveraging}, taking into account that we may express $H$ (on any finite dimensional subspace) as a sum of squares of linear forms:
  
\begin{lemma} \label{nui}
Let $H$ be a Hermitian form on a locally finite $S$-vector space $V$. 
Let $\nu_i$ be an averaging sequence; then, for any $v \in V_f$, 
$$\int H(a v) \nu_i = : \nu_i \star H(v) \longrightarrow H^S(v).$$
The same assertion holds for every $v\in V$ (with possible infinite values on the right hand side) if it is assumed that all exponents of $H$ satisfy $|\chi|=1$ or $|\chi^{-1}|<1$ on $S^+$.
\end{lemma}

We shall say that a form $H$ on a locally finite $S$-vector space
is {\em $c$-good} if  
\begin{equation} \label{goodexponents} 
\mbox{Every exponent $\chi$ of $H$ satisfies $|\chi| = 1$ or $|\chi^{-1}| \leq c^2$.}
\end{equation} 
The reason for inverting $\chi$ is that the exponents of $S$ on linear functionals are inverse to those on the vector space itself; we write $c^2$ for convenience when comparing to the case of the square of a linear functional. 

{ \begin{remark}
Notice that $|\chi^{-1}|\leq c^2$ rules out the possibility of subunitary exponents which are not strictly subunitary. We do that for simplicity, since 
in later sections we will only use the case where the rank of $S$ is $1$, so subunitary and strictly subunitary coincide.
\end{remark}
}

The form $H$ is simply {\em good} if it is $c$-good for some $c < 1$.  In particular, this excludes
the possibility of exponents that grow along some ``walls'' of $S^+$ but not along other. Similar terminology will be applied to a vector, if it applies to its $S$-span.

\begin{lemma}  \label{cgoodrad} 
Suppose $H$ is $c$-good.   Let $R$ be the sum of all generalized eigenspaces
corresponding to characters $\chi$ that satisfy neither
$|\chi| = 1$ nor $|\chi| \leq c$. Then $R$ lies in the radical of $H$ and
$H$ factors through $V/R$.  

In particular one may write (on any finite-dimensional, $S$-stable subspace): $H = \sum |\ell_i|^2$, where each $\ell_i$ is linear 
and each Hermitian form $|\ell_i|^2$ is itself $c$-good.  \end{lemma}

\begin{proof} 
 If $v \in V$ is a
eigenvector for $S$, with eigencharacter $\chi$, then
$$  \langle H^s \mbox{ considered in } V^* \otimes \overline{V^*},  v \otimes  v
\mbox{ considered in }  V \otimes \overline{V} \rangle  
 = |\chi(s)|^{-2} H(v) .$$
 
 That implies that $|\chi|^{-2}$ is an exponent of $H$ if $H(v) \neq 0$. 
 In particular, if $\chi$ doesn't satisfy $|\chi| = 1$ or $|\chi| \leq c$, 
 then $H(v) = 0$, i.e., $v$ lies in the radical of $H$ (because $H$ is semidefinite, 
 $H(v) = 0 \implies H(u+v) = H(u)$ for all $u$, i.e.\ $v$ is in the radical). 
 
   This conclusion remains valid if
$v$ were simply a {\em generalized} eigenvector. This proves the first assertion of the lemma.

For the second assertion (concerning linear forms):
Choose any expression of $H$ as a sum of squares on the vector space $V/R$. 
If $\ell$ is any linear functional on $V/R$ 
then all exponents of $\ell \otimes \overline{\ell} \in V \otimes V^*$ 
are of the form $\chi = (\psi_1 \overline{\psi_2})^{-1}$, where $\psi_{1}, \psi_2$
are exponents of $S$ on $V/R$. In particular, $ |\chi^{-1} | = |\psi_1| |\psi_2| $
satisfies $|\chi^{-1}| = 1$ or $|\chi^{-1}| \leq c$; that proves the second assertion.

\end{proof}

\begin{proposition} \label{invub} Let $H$ be a hermitian form on a locally finite $S$-vector space $V$, whose elements have degree bounded by $n$, and assume that $H$ is $c$-good. 
  
There exist a finite subset $\Lambda \subset S^+$ and a constant $C$ depending only on $(n,c)$ so that:
\begin{equation} \label{moo1} H^a \leq  C (H^S + \sum_{\lambda \in\Lambda} H^{\lambda}),\end{equation}
for any $a \in S^+$; here $H^a(v) := H(av)$.  
\end{proposition}

 \begin{proof}
  In fact, we may suppose $V$ to be finite-dimensional and $H$ to be the square
  of a linear form $\ell$, in view of the prior Lemma.

 In that case the function
 $a \mapsto \ell(av)$ is finite for any fixed $v \in V$. Write
 $$\ell(av) = f(a) + g(a)$$
where the functions $f,g$ possess only unitary (resp. sub-unitary) exponents, i.e.\ all exponents of $f$ satisfy $|\chi|=1$, and all exponents of $g$ satisfy $|\chi| \leq c <1$.  The degrees of $f,g$ are both bounded by $n$.
 
 If $f$ (equivalently $f+g$) is unbounded
 on $S^+$, then the result is obvious, as the right-hand side of the putative inequality is infinite.
 We suppose therefore that $f$ is bounded on $S^+$. 
   Apply Lemmas \ref{boundS} and \ref{boundU} (taking $\Lambda$ as in Lemma \ref{boundS}):
 \begin{eqnarray*}
  \begin{aligned} \sup_{S^+} |f+g|^2  
  & \leq &    2 \sup_{S^+} |f|^2 + 2  \sup_{S^+} |g|^2 
 \\  & \leq & 2 n \limS |f|^2 +2  \max_{\Lambda} |g|^2 
\\   & \leq &  2 n \limS |f|^2 +  4 \max_{\Lambda} \left( |g+f|^2 + |f|^2 \right)
\\  & \leq  & 6n \limS |f|^2+ 4 \max_{\Lambda} |f+g|^2 
\\ &=& 6n \limS |f+g|^2 +4 \max_{\Lambda} |f+g|^2, \end{aligned}
  \end{eqnarray*}
the last line since $\limS |f+g|^2 = \limS |f|^2$. 
\end{proof} 
 
\begin{corollary} \label{subunitarybound} Notation and assumptions as in the previous lemma, let $\Pi^{<1}$ be the $S$-equivariant projection
 of $V$ onto all generalized eigenspaces with subunitary   exponent, and $\Pi^{=1}$
  the $S$-equivariant projection onto generalized eigenspaces with unitary exponent. Let $H^{<1} = H \circ \Pi^{<1}$ and $H^{=1} = H \circ \Pi^{=1}$, and denote by $n_1$ the number of distinct unitary exponents of $V$. 
  
 Then:
\begin{eqnarray} H^{=1}(v)  & \leq &  n_1 H^S(v) \\ 
H^{<1}(av) & \leq  &  Q(a) (H^S(v) + \max_{\Lambda} H^{\lambda}(v)) \end{eqnarray}
where the decaying function $Q$ depends only on $(n,c)$.  
  \end{corollary}

 \begin{proof}

 Both assertions are linear in $H$, so we reduce, as in the previous proof, to the case where $H = |\ell|^2$.  Note that $\ell = \ell \circ \Pi^{=1} + \ell \circ \Pi^{<1}$;
 indeed (adopting the notation of Lemma \ref{cgoodrad})
the functional $\ell$ may be supposed to contain $R$ in its kernel,
but $\mathrm{id} - \Pi^{=1} - \Pi^{<1}$ is a projector to a subspace of $R$. 
 
 Note that the decomposition $\ell = \ell \circ \Pi^{=1} + \ell \circ \Pi^{<1}$
 induces the decomposition $\ell(av) = f +g$ of the prior proof, i.e.
 $H^{=1}(av) = |f|^2$
and $H^{<1}(av) = |g|^2$. For the first bound, now, apply  Lemma \ref{boundU}. (Notice that either the right hand side is infinite, or the eigenprojection of $v$ to the sum of subspaces with unitary exponents has degree bounded by the number of unitary exponents.)
For the second,  bound $|g|^2$ using \eqref{finebound}. 
 \end{proof}

\subsection{Measurability of eigenprojections} \label{meas-eigenpproj}
Here we discuss the following issue: Let $S^+$ be  a finitely generated abelian monoid, and let $W$ be an $\mathbb C[S^+]$-module of countable dimension. Let $(L_y)_{y\in Y}$ be a family of linear functionals on $W$, varying with a parameter $y$ in a measurable space $(Y,\mathcal B)$ (i.e.\ a set equipped with a $\sigma$-algebra). Assume that the family $L_y$ is measurable in the sense of evaluations (the function $y\mapsto \left<w,L_y\right>$ is measurable for every $w\in W$), and that \emph{each $L_y$ is $S^+$-finite}, i.e.\ generates a finite-dimensional $\CC[S^+]$-submodule of $W^*$.

Hence, each $L_y$ has a decomposition in generalized $S^+$-eigenspaces:
\begin{equation}
 L_y = \oplus_{\chi\in \widehat{S^+}_\CC} L_y^\chi,
\end{equation}
where $\widehat{S^+}_\CC$ denotes the space of all characters (not necessarily unitary) of $S^+$. Let $K\subset \widehat {S^+}_\CC$ be a measurable subset (with respect to the natural Borel $\sigma$-algebra on $\widehat{S^+}_\CC$), and denote:
\begin{equation}
 L_y^K = \oplus_{\chi\in K} L_y^\chi.
\end{equation}

We want to prove:
\begin{proposition}\label{measurable2}
 The functionals $L_y^K$ vary measurably with $y$.
\end{proposition}
Again, measurability here is meant with respect to evaluations at every vector of $W$, like above.

In practice, we will use this proposition to isolate the unitary part of an $S^+$-finite functional or a hermitian form.

First of all, we can reduce the proof of the proposition to the case $S^+=\mathbb N$. Indeed, first we replace $S^+$ by $\mathbb{N}^k$ by taking a surjective map: $\mathbb{N}^k \to S^+$, and then we observe that the Borel $\sigma$-algebra of $\widehat{\mathbb N^k}_\CC$ is the product $\sigma$-algebra of the Borel $\sigma$-algebras of $\widehat{\mathbb N}_\CC$, and hence we may assume without loss of generality that $K$ is a product subset. (Let us clarify what this means for the functionals $L_y^K$: writing $K$ as a countable union of subsets $K_n$ corresponds to writing $L_y^K$ as the weak limit of functionals $L_y^{K_n}$, and this limit stabilizes for a given $y$. Since for the complement $K'$ of $K$ we have $L_y^{K'}= L_y-L_y^K$, the same argument shows that the set of subsets $K$ which satisfy the Proposition is closed under countable unions and intersections, i.e.\ forms a $\sigma$-algebra, and hence it is enough to check for a set of subsets generating the $\sigma$-algebra.) Then we can obtain $L_y^K$ 
in a finite number of steps by restricting the generalized eigencharacters coordinate-by-coordinate. Hence, from now on we will assume that $S^+=\mathbb N$, and we will denote its generator by $x$.

The proposition will now be established via the following result: Identify $\CC[S^+]$ with the ring $\CC[x]$ of polynomials in one variable (where $x$ corresponds to the generator of $S^+\simeq \mathbb N$, and call \emph{minimal polynomial} of $L_y$ the monic generator of its annihilator in $\mathbb C[S^+]$. We will denote it by $\mathfrak m_y$. 

There is a natural measurable structure on $\CC[S^+]\simeq \CC[x]$. Namely, a set is measurable if for any $d$ its intersection with polynomials of degree $\le d$
is Borel-measurable with respect to the standard topological structure on that complex vector space. 

\begin{lemma}\label{measurablepoly}
 The minimal polynomials $\mathfrak m_y$ vary measurably in $y\in Y$.
\end{lemma}

Let us see why this implies Proposition \ref{measurable2}.

First of all, by partitioning $Y$ in a countable union we may assume that the degree of $\mathfrak m_y$ is constant in $Y$, say $\mathfrak m_y\in \CC[S^+]_N$ (degree $N$). The coefficients of each polynomial are elementary symmetric functions in its roots, and we may pick a measurable section: $\CC[S^+]_N\to \CC^N$ of the map $\CC^N\ni(\alpha_i)_i \mapsto \prod_i (x-\alpha_i)\in \CC[S^+]$. (We have continued with the prior notation, so that $x$ is an element of $\CC[S^+]$
corresponding to a generator $g$ for $S^+$).  Hence, we may index the roots $(\alpha_{i,y})_i$ of $\mathfrak m_y$ in a measurable way. Finally, for given measurable $K\subset \CC$, the set $A_y\subset\{1,\dots,N\}$ of indices such that $\alpha_{i,y}\notin K$ varies measurably with $y$. Hence, we may write the minimal polynomial $\mathfrak m_y$ in a measurable way as a product:
$$\mathfrak m_y = \mathfrak m_y^1 \mathfrak m_y^2,$$
where $\mathfrak m_y^1 = \prod_{i\in A_y} (x-a_{i,y})$ and $\mathfrak m_y^2= \prod_{i\notin A_y} (x-a_{i,y})$.
The polynomials $\mathfrak m_y^1$ and $\mathfrak m_y^2$ are relatively prime, hence we can find a polynomial $\mathfrak m_y^3$ which is inverse to $\mathfrak m_y^1$ in $\CC[x]/\mathfrak m_2$. Again, this can be done in a measurable way using the division algorithm.
 Then:
\begin{equation}  
L_y^K = \mathfrak m_y^3(x) \mathfrak m_y^1(x) L_y.
\end{equation}
Indeed, $\mathfrak m_y^1(x)$ annihilates the summands of $L_y$ with exponents outside of $K$, and since on the rest of the summands $\CC[x]$ acts via the quotient $\CC[x]/\mathfrak m_y^2$, the product $\mathfrak m_y^3(x) \mathfrak m_y^1(x)$ acts as the identity on them. This shows that $L_y^K$ is measurable, i.e., it concludes the proof that Lemma \ref{measurablepoly}
implies Proposition \ref{measurable2}. 

To prove Lemma \ref{measurablepoly}, we notice that $\mathfrak m_y (x) = x^N + a_{N-1} x^{N-1} +\dots+ a_0$ if and only if:
\begin{enumerate}
 \item for every $w\in W$ we have: $L_y(x^N\cdot w)+ a_{N-1} L_y(x^{N-1}\cdot w) + \dots + a_0 L_y(w)= 0$, and
 \item this is not the case for any polynomial of smaller degree.
\end{enumerate}
We may enumerate a vector space basis $(w_i)_{i\ge 1}$ of $W$, and for every $N, n$ we consider the following system of linear equations in the unknowns $a_0,\dots,a_{N-1}$:
\begin{equation}
 \mathbf S_{N,n}: \left( L_y(x^N\cdot w_i)+ a_{N-1} L_y(x^{N-1}\cdot w_i) + \dots + a_0 L_y(w_i)= 0\right)_{i=1}^n.
\end{equation}
The polynomial $x^N + a_{N-1} x^{N-1} +\dots+ a_0$ satisfies the above two conditions (i.e.\ is the minimal polynomial $\mathfrak m_y)$ if its coefficients are the \emph{unique} solution of the system: $\mathbf S_{N,\infty} = \cup_n \mathbf S_{N,n} $.

For given $n, N$ the set of $y\in Y$ such that the system has a solution is measurable; indeed we can attempt to solve the system by row operations, the order of which depends only on whether some coefficients vanish or not (which, of course, depends measurably on $y\in Y$). Among those, uniqueness of the solution is also a measurable condition, for the same reason. Finally, among the latter the unique solution $(a_0, a_1, \dots,  a_{N-1})$ varies measurably in $y\in Y$, again for the same reason. 

Hence, for given $N$ the set $Y_N$ of $y\in Y$ such that $\mathbf S_{N,\infty}$ has a unique solution is measurable (notice that for given $y$, if $\mathbf S_{N,\infty}$ has a unique solution then so does $\mathbf S_{N,n}$ for some $n$), and the solution varies measurably in $y\in Y$.  
This proves the lemma.

We will also use this result in the following form:

\begin{corollary}[Measurability of eigenspaces.] \label{cor:meas}   Let $W$ be a finite dimensional vector space 
 and $T(z) \in \mathrm{End}(W)$ a family of matrices
 varying measurably as $z$ varies in a measurable space $Z$.  
Then the $T(z)$-invariant  projection to the generalized $0$ eigenspace varies measurably with $z$.

More generally, suppose that $\alpha_z: \mathbb{Z}^N \rightarrow \mathrm{Aut}(W)$
 is a measurable family of actions and $\chi: \mathbb{Z}^N \rightarrow \mathbb{C}^{\times}$ a character. Then the canonical ($\alpha_z(\mathbb Z^n)$-invariant) projection of $W$ to the generalized $\chi$-eigenspace for $\alpha_z$
 varies measurably with $z$. 
\end{corollary}

\section{The Bernstein morphisms} \label{sec:Bernstein}

\textbf{From now on  we assume that the Discrete Series Conjecture \ref{dsconjecture} holds for  $\XX$ and all its degenerations $\XX_\Theta$ (for example, $\XX$ is strongly factorizable). } { Although the structure of discrete series will not be used explicitly in the present section, we will use its corollaries, such as the boundedness of subunitary exponents \ref{uniformbound}.}
Since this will be an ongoing assumption, it will not be explicitly included in the theorems.    

\subsection{} \label{bmp}

In the present section, we construct a canonical morphism $L^2(X_\Theta) \rightarrow L^2(X)$.
We may think of this morphism $L^2(X_{\Theta}) \rightarrow L^2(X)$ as 
\begin{itemize}
\item[-]  {\em the unique morphism asymptotic to the
``naive'' identification of functions on $X_{\Theta}$ and $X$.}
\end{itemize}

It may be worth beginning our section with the following easy Lemma, which
contains the germ of many of the ideas used in this section: 
\begin{lemma}  \label{berngerm} The support of Plancherel measure for $L^2(X_{\Theta})$ is contained in the support of  Plancherel measure for $L^2(X)$. 
\end{lemma} 

\proof  
It suffices to show that matrix coefficients of the form 
\begin{equation}\label{matrixcoeff}\left<g\cdot f, f\right>,
\end{equation}
 where $f\in C_c^\infty(X_\Theta)$ can be approximated, uniformly on compacta, by diagonal matrix coefficients of $L^2(X)$. Assume that $f$ is $J$-invariant, and given a compact, $J$-biinvariant subset $K$ of $G$ set $J'=\cap_{k\in K} kJk^{-1}$. We may translate $f$ by the action of $\mathcal Z(X_\Theta)$ into a $J'$-good neighborhood of infinity with the property that the identification with a neighborhood of $\Theta$-infinity on $X$ is equivariant\footnote{Equivariance implicitly assumes the identification of larger neighborhoods, of course.} under the action of the elements of $\mathcal H(G,J')$ whose support is in $K$. Then the matrix coefficients (\ref{matrixcoeff}) coincide on $K$, whether $f$ is considered as a function on $X_\Theta$ or on $X$.
\qed

We now formulate properties of the morphism more precisely.

Recall that 
\begin{equation}\label{etheta}
e_\Theta: C_c^\infty(X_\Theta) \to C_c^\infty(X) 
\end{equation}
denotes the equivariant ``asymptotics'' map which, for $J$-invariant functions supported in a $J$-good neighborhood of $\Theta$-infinity, coincides with the identification of $J$-orbits under the exponential map. The formulations that follow involve, actually, only functions supported close enough to $\Theta$-infinity, and therefore the equivariant extension of this identification is not being used in the statement of the theorem. 

We will need to ``push'' functions towards $\Theta$-infinity; recall that $\mathring A_{X,\Theta}^+$ denotes the subset of ``strictly anti-dominant'' elements of $A_{X,\Theta}$, i.e.\ those which push points on $X_\Theta$ towards $\Theta$-infinity.

Denote by $\mathcal L_a$ the (normalized) action of $a\in A_{X,\Theta}$ on functions on $X_\Theta$, 
which we understand as $a^{-1} \cdot f$ in the normalization of \S \ref{ssmeasures}; in other words, 
the $\mathcal{L}_a$ for $a \in \mathring A_{X, \Theta}^+$ push functions towards $\infty$.

\begin{theorem}\label{Bernsteinmap} For every $\Theta\subset\Delta_X$ there is a canonical $G$-equivariant morphism: $\iota_\Theta: L^2(X_\Theta)\to L^2(X)$, characterized by the property that for any $a\in \mathring A_{X,\Theta}^+$ and any $\Psi\in C_c^\infty(X_{\Theta})$ we have:
\begin{equation} \label{bchar}
 \lim_{n\to\infty} (\iota_\Theta \mathcal L_{a^n}\Psi- e_\Theta \mathcal L_{a^n}\Psi)=0.
\end{equation}
\end{theorem}

In words, \eqref{bchar} says that $\iota_{\Theta}$ is approximately equal to the identity furnished by the exponential map on functions supported near $\Theta$-infinity. 
Although results of this type are present in the scattering theory literature, 
the idea of the proof that we present here (essentially, the proof of Theorem \ref{bernstein-abstract}) is due to Joseph Bernstein and we will consequently refer to $\iota_{\Theta}$
as the ``Bernstein morphism''.

\begin{remark}  If we replace $e_\Theta$ by $\tau_\Theta:=$``truncation in a fixed $J$-good neighborhood $N_\Theta$ of $\Theta$-infinity'', we can generalize property (\ref{bchar}) to non-compactly supported smooth $L^2$-functions:  
\begin{equation}\label{bcharfornoncompact}
 \lim_{n\to\infty} (\iota_\Theta \mathcal L_{a^n}\Phi- \tau_\Theta \mathcal L_{a^n}\Phi)=0.
\end{equation}
As in the theorem, the limit is taken inside $L^2(X)$, { and we identify the function $\tau_{\Theta} \mathcal{L}_{a^n} \Phi$,
with a function on $X$ by means of the exponential map.}

The proof is very simple: For given $\varepsilon>0$ we can find $m$ and a function $\Psi\in C_c^\infty(N_\Theta)^J$ such that $\Vert \Psi - \mathcal L_{a^m}\Phi\Vert <\varepsilon$ (and hence $\Vert \mathcal L_{a^n}\Psi - \mathcal L_{a^{m+n}}\Phi\Vert <\varepsilon$ for every $n\ge 0$). Then, applying (\ref{bchar}) to $\Psi$, we get some $n$ such that:
$$ \Vert \iota_\Theta \mathcal L_{a^n}\Psi- e_\Theta \mathcal L_{a^n}\Psi\Vert<\varepsilon.$$
Therefore, for large enough $n$:
$$ \Vert \iota_\Theta \mathcal L_{a^{m+n}}\Phi- \tau_\Theta \mathcal L_{a^{m+n}}\Phi\Vert \le \Vert \iota_\Theta\left(\mathcal L_{a^{m+n}}\Phi - \mathcal L_{a^n}\Psi\right)\Vert +  $$
$$+ \Vert \iota_\Theta \mathcal L_{a^n}\Psi- e_\Theta \mathcal L_{a^n}\Psi\Vert + \Vert  e_\Theta \mathcal L_{a^n} \Psi- \tau_\Theta \mathcal L_{a^{m+n}}\Phi \Vert < (\Vert\iota_\Theta\Vert+2)\varepsilon.$$
{ For the last inequality, we note that for large enough $n$
that $e_{\Theta} \mathcal{L}_{a^n}\Phi$ is obtained from $\tau_{\Theta} \mathcal{L}_{a^n} \Phi$
via the identification of $J$-orbits arising from the exponential map; thus the last term has norm bounded by the norm of $\| \tau_{\Theta} \mathcal{L}_{a^n}
\Phi - \tau_{\Theta} \mathcal{L}_{a^{m+n}} \Phi \|$, which is at most $\varepsilon$.}

However, the proof of Theorem \ref{Bernsteinmap} will be via an estimate:
\begin{equation} \label{rateofconvergence}
 \left\Vert\iota_\Theta \mathcal L_{a^n}\Psi- e_\Theta \mathcal L_{a^n}\Psi\right\Vert^2 \le C_\Psi \cdot Q^J(a^n),
\end{equation}
for $\Psi\in C_c^\infty(X_\Theta)^J$, where $Q^J$ is a decaying function on $\mathring A_{X,\Theta}^+$ (cf.\ Lemma \ref{decaying}) which depends only on the open compact subgroup $J$, and $C_\Psi$ is a constant that depends on $\Psi$, see Lemma \ref{decaying}. Such an estimate is not valid for functions which are not compactly supported.
\end{remark}

\subsection{ Harish-Chandra--Schwartz space and temperedness of exponents}  \label{HCS}
In \cite{BePl} Bernstein explains how to prove that the Plancherel formula 
for a space \emph{of polynomial growth} like $X$ is supported on $X$-tempered\footnote{The reader of \cite{BePl} will notice that the notion of temperedness obtained there is slightly stronger than temperedness with respect to the Harish-Chandra--Schwartz space; indeed, we can replace that with $L^2(wdx)^\infty$ for any \emph{summable} weight $w$. However, since summability of weights depends on the rank of the variety, it is more convenient to work with the weaker condition provided by the Harish-Chandra--Schwartz space.} representations. 

We remind what this means. In order to do this we reprise some of the remarks of \S \ref{Plancherelgeneralities}, but replacing the role of $\bruhat(X)$ by the slightly larger space  
$\mathscr C(X)$ of \emph{Harish-Chandra--Schwartz} functions on $X$. 

A function $r: X\to \RR_+$ is called a \emph{radial function} if it is positive, locally bounded and proper, i.e.\ such that the balls $B(a):=\{x\in X| r(x)\le a\}$ are relatively compact, and has the property that for every compact $J\subset G$ there is a constant $C>0$ such that $|r(x\cdot g)-r(x)|<C$ for all $x\in X, g\in J$. Two radial functions $r$ and $r'$ are called equivalent if the quotient $ \frac{1+r}{1+r'}$ is bounded above and below by absoltue positive constants.

A space is called of \emph{polynomial growth} (for a given radial function) if for some compact $J\subset G$ there is a polynomial $a\mapsto P(a)$ such that for all $a>0$ the ball $B(a)$ can be covered by $\le P(a)$ sets of the form $x\cdot J$. This notion is clearly invariant under equivalence of radial functions.

By the generalized Cartan decomposition (see \S \ref{ssCartan}),\footnote{Using the Cartan decomposition is again not necessary: it is enough to know that the union of $J$-good neighborhoods of $\Theta$-infinity, for all $\Theta\subsetneq\Delta_X$, has a compact complement (modulo center).} the space $X$ of $k$-points of our spherical variety possesses a natural equivalence class of radial functions, with respect to which it is of polynomial growth. { In fact, if $\mathcal Z(X)=1$, such a radial function $R(x)$ was described in the proof of Theorem \ref{propIIconj}. We leave the details of the general case to the reader (the only difference being that one needs to quantify as well the ``distance'' from \emph{all} orbits in a smooth toroidal compactification of $X$, including ``orbits belonging to $\Delta_X$-infinity'').} We fix such a radial function $r$.

Then we define the Harish-Chandra--Schwartz space as the Fr\'echet space: 
\begin{equation}\mathscr C(X)=\lim_{\overset{\to}{J}}\bigcap_d L^2(X, (1+r)^d dx)^J
\end{equation}
the limit taken over a basis of neighborhoods of the identity.

 Bernstein proves that the embedding $\mathscr C(X)\to L^2(X)$ is \emph{fine} which implies that any Hilbert space morphism to a direct integral of Hilbert spaces: 
\begin{equation}\label{Hilbertintegral}
 m: L^2(X)\to \int \mathcal H_\alpha \mu(\alpha)
\end{equation}
 is \emph{pointwise defined} on $\mathscr C(X)$, i.e.\ there is a family of morphisms $$L_\alpha:\mathscr C(X)\to \mathcal H_\alpha$$ (defined for $\mu$-almost every $\alpha$) such that $\alpha\mapsto L_\alpha(\Phi)$ represents $m(\Phi)$ for every $\Phi\in \mathscr C(X)$. 

Notice that such a family of morphisms induces, by pull-back, seminorms $\Vert\bullet\Vert_\alpha$ on $\mathscr C(X)$. If the morphism $m$ is surjective (set-theoretically, hence open) then the spaces $\mathcal H_\alpha$ can be identified with the completions of $\mathscr C(X)$ with respect to the seminorms $\Vert\bullet\Vert_\alpha$. 

A Plancherel decomposition for $L^2(X)$ -- or, more generally, for some closed invariant subspace of $L^2(X)$ -- can be described by
the choice of a measure $\mu$ on $\hat{G}$ and a  measurably varying family $\Vert \bullet \Vert_{\pi}$ of norms on $\mathscr C(X)$,
with the properties that $\Vert \bullet \Vert_{\pi}$ factors through the natural morphism from $\mathscr C(X)$ to the space of $\pi$-coinvariants\footnote{In the case of the Harish-Chandra--Schwartz space, ``homomorphism'' will always mean ``continuous homomorphism''. The space of (smooth vectors on) $\pi$ is endowed with the discrete topology or, what amounts to the same for homomorphisms, the coarsest $\CC$-vector space topology. Similarly, ``linear functionals'' and ``hermitian forms'' will always be continuous.}
\begin{equation} \label{HCcoinvariants} \mathscr C(X)_\pi := \left(\Hom_G(\mathscr C(X), \pi)\right)^* \otimes \pi.
\end{equation} and also $\|\Phi\|_{L^2(X)} ^2= \int_{\pi} \Vert \Phi \Vert_{\pi} ^2 \mu(\pi)$ for every $\Phi \in \mathscr C(X)$. 
(Recall that in the case of wavefront spherical varieties, which we are discussing, the spaces $\Hom_G(\mathscr C(X), \pi)$ are finite-dimensional.)

We have a canonical map (recall that the space of $\pi$-coinvariants of $C_c^\infty(X)$ was defined in a completely analogous way in (\ref{coinvariants})):
\begin{equation}\label{SchwartztoHC}
C_c^\infty(X)_\pi \twoheadrightarrow \mathscr C(X)_\pi. 
\end{equation}

\begin{remark} \label{decomp-into-eigen} It is sometimes more convenient to think of the Plancherel formula as giving a ``decomposition into eigenfunctions'': 

The Hermitian norms $\Vert \bullet \Vert_{\pi}$ induce
$C^{\infty}_c(X)_{\pi} \rightarrow \overline{\left(C^{\infty}_c(X)_{\pi}\right)^\sim} = 
\overline{C^{\infty}(X)^{\pi}}$. For every $f \in C^{\infty}_c(X)$, 
the conjugate of the image of $f$ under the map\footnote{Here $C^{\infty}(X)^{\pi}$ denotes the $\pi$-isotypical subspace, which is to say, the image
 of $\pi \otimes \Hom(\pi, C^{\infty}(X)) \rightarrow C^{\infty}(X)$. }
$$C^{\infty}_c(X) \rightarrow C^{\infty}_c(X)_{\pi} \rightarrow \overline{C^{\infty}(X)^{\pi}}$$
will be denoted by $f^{\pi}$. 
 
Then we have a pointwise decomposition:
\begin{equation}\label{pointwisedecomposition}
 f(x) = \int_{\hat G} f^{\pi}(x)\mu(\pi),
\end{equation} 
which is another way of writing the Plancherel decomposition for the inner product $\left< f, \Vol(xJ)^{-1} 1_{xJ}\right>$ for a sufficiently small open compact subgroup $J$. Note that we also have for $f, g \in C^{\infty}_c(X)$ the equality
\begin{equation} \label{symmetryproperty} \langle f,  g^{\pi} \rangle_{X} = \langle f^{\pi}, g \rangle_X \end{equation}
since both are different ways to express the inner product $H_\pi(f,g)$ (where $H_\pi$ is the hermitian form corresponding to $\Vert\bullet\Vert_\pi$).

We remark, however, that although the image $f_\pi$ of a function $f$ in the space of $\pi$-coinvariants is canonically defined, this is not the case for $f^\pi$, which \emph{depends on the choice of Plancherel measure}.

 \end{remark}

By the asymptotics map $e_\Theta$ (cf.\ (\ref{etheta})) we get a canonical map:
$$ C_c^\infty(X_\Theta)_\pi\to C_c^\infty(X)_\pi.$$

The following follows directly from the definitions and the discussion of \ref{sshermitian}. Recall that we always consider the normalized action of $A_{X,\Theta}$ on functions on $X_\Theta$, and that the space $C_c^\infty(X_\Theta)_\pi$ is $A_{X,\Theta}$-finite for every irreducible representation $\pi$. Hence, we may decompose into sums of generalized eigenspaces:
\begin{equation}\label{leornleq} C_c^\infty(X_\Theta)_\pi = C_c^\infty(X_\Theta)_\pi^{< 1} \oplus C_c^\infty(X_\Theta)_\pi^{1} \oplus C_c^\infty(X_\Theta)_\pi^{\nleq 1}
\end{equation}
with exponents, respectively, satisfying\footnote{Our notation is explained as follows: We denote the sum of eigenspaces of $C_c^\infty(X_\Theta)_\pi$ by exponents satisfying $|\chi^{-1}|< 1$ by $C_c^\infty(X_\Theta)_\pi^{<1}$ because \emph{the hermitian forms will be subunitary there}. Another way to think of it is the following: if $l$ is a smooth linear functional on $C_c^\infty(X_\Theta)_\pi^{<1}$, then via the duality $C_c^\infty(X_\Theta)\otimes C^\infty(X_\Theta)\to \CC$ it can be considered as an element of $C^\infty(X_\Theta)$. \emph{That element will be decaying on $\mathring A_{X,\Theta}^+$.}}
 $|\chi^{-1}|< 1$, $|\chi|=1$ and $|\chi^{-1}|\nleq 1$ on $\mathring A_{X,\Theta}^+$.
 
 Here $|\chi^{-1}| \nleq 1$ means that there exists $a \in \mathring A_{X, \Theta}^+$ with $|\chi^{-1}(a)|  >  1$. 
 The three possibilities are mutually exclusive and one must occur for each $\chi$: 
 Recall that $\mathring A_{X, \Theta}^+$ denotes the ``strict interior'' of the cone $A_{X, \Theta}^+$. Now, 
if the final possibility does not occur, then  $|\chi^{-1}| \leq 1$ on $\mathring A_{X, \Theta}^+$;
then $|\chi^{-1}| $ is bounded above on $A_{X, \Theta}^+$, from which one sees that $|\chi^{-1}| \leq 1$
on $A_{X, \Theta}^+$; but that means that 
either  $|\chi|=1$ or $|\chi|< 1$ on the ``strict interior'' $\mathring A_{X, \Theta}^+$.

\begin{lemma} \label{HCform}
 Let $H_\pi$ be a $G$-invariant hermitian form on $C_c^\infty(X)_\pi$.   
\begin{enumerate}
 \item Through the map $e_\Theta$ it is pulled back to an $A_{X,\Theta}$-finite, $G$-invariant Hermitian form $e_\Theta^* H_\pi$ on $C_c^\infty(X_\Theta)_\pi$.
 \item Suppose  the form $H_{\pi}$ factors through the Harish-Chandra--Schwartz space, i.e.\ through \eqref{SchwartztoHC}. 
 Then, for any $\Theta$, if we decompose as in (\ref{leornleq}), the summand $C_c^\infty(X_\Theta)_\pi^{\nleq 1}$ lies inside the radical of $e_{\Theta}^* H_{\pi}$ (in particular, the form vanishes there).  
\end{enumerate}
\end{lemma}

\begin{proof}

For part (1), the only assertion to be proved is the $A_{X, \Theta}$-finiteness; however,
in our present case, the multiplicity of $\pi$ in $C^{\infty}(X_{\Theta})$ is finite, from which the result follows at once.

For part (2) it is enough to show
that $C_c^\infty(X_\Theta)_\pi^{\nleq 1}$
lies in the kernel of the composite
$$C_c(X_{\Theta})_{\pi} \longrightarrow C^{\infty}_c(X)_{\pi} \rightarrow \mathscr{C}(X)_{\pi}.$$

In other words, given a morphism $\lambda:\mathscr{C}(X) \rightarrow \pi$, 
we need to verify that $$C_c(X_{\Theta}) \stackrel{e_{\Theta}}{\rightarrow} C^\infty_c(X)
\rightarrow \mathscr{C}(X)  \stackrel{\lambda}{\rightarrow} \pi$$ vanishes on $C_c^\infty(X_\Theta)_\pi^{\nleq 1}.$

Suppose that the $\chi$-eigenspace on $C^{\infty}_c(X_{\Theta})_{\pi}$
is nonzero, and that there exists 
 $a \in \mathring A_{X, \Theta}^+$ so that $|\chi(a)| > 1$. 
  
 Take $\Psi \in C^{\infty}_c(X_{\Theta})^J$. 
 Then (by continuity of $\lambda$) the norm of the image of $\mathscr{L}_{a^n} \Psi$
 grows at most polynomially in $n$.   But the projection of $\mathscr{L}_{a^n} \Psi$
 to this $\chi$- generalized eigenspace -- if nonzero -- grows as $\chi(a^n)$, at least for a subsequence of $n$: 
  after replacing this projection by a linear combination of translates by various $a^k$,
 we may suppose that it  is a nonzero element of the genuine -- not just generalized --  $\chi$-eigenspace.
 \end{proof}

\subsection{Plancherel formula for $X_\Theta$ from Plancherel formula for $X$}

The theorem below is the heart of Theorem \ref{Bernsteinmap}.

\begin{theorem} \label{bernstein-abstract}
Consider a Plancherel decomposition for $L^2(X)$ (cf.\ \ref{Plancherelgeneralities}):
\begin{equation}\label{Pl}  \Vert\Phi\Vert^2= \int_{\hat G} H_\pi(\Phi) \mu(\pi). \end{equation}
{ Fix an open compact subgroup $J$}, and 
fix a  strictly positive cocharacter $s:  \mathbb{G}_m \rightarrow \AA_{X,\Theta}$ (i.e.\ a cocharacter in the strict interior, in $\Lambda_X^+$, of the face corresponding to $\Theta$),
and { let $S = s(\varpi^{\Z})$ for $\varpi$ a uniformizer of $k$. }
 
 Consider the pullback $e_\Theta^*H_{\pi}$ of $H_\pi$ to $C^{\infty}_c(X_{\Theta})_{\pi}$, and let $\left(e_\Theta^*H_{\pi}\right)^S$ be the associated  $S$-invariant form (by Lemma \ref{invariantform}).  
 Let $\Psi \in C^{\infty}_c(X_{\Theta})^J$; then
\begin{equation}\label{pullbackdecomp}
\|\Psi\|^2 = \int_{\hat G}  \left(e_\Theta^*H_{\pi}\right)^S(\Psi) \mu(\pi).
\end{equation}

Therefore, the hermitian forms $\left(e_\Theta^*H_{\pi}\right)^S$ define a Plancherel formula for $L^2(X_{\Theta})^J$.
\end{theorem}

\begin{remark}
 In particular, for almost every $\pi$ the invariant forms $\left(e_\Theta^*H_{\pi}\right)^S$ are $A_{X,\Theta}$-invariant and finite, and do not depend on the choice of $S$; we will be denoting them by $H_\pi^\Theta$. 
 
 Indeed, it is entirely  possible that these norms take infinite values for some $\pi$; but
 this must happen only on a set of measure zero: Clearly, for each individual
 $\Phi \in C_c(X_{\Theta})$, the set of $\pi$ for which $e_{\Theta}^* H_{\pi}(\Phi) = \infty$ has measure zero. Since $C_c(X_{\Theta})$ has countable dimension, 
 this implies the stronger statement, because, if a norm takes infinite values, it does so on at least
 one element of a basis. 
\end{remark}

\begin{remark} This theorem is roughly the analog of \eqref{azz} from the discussion of the toy model of scattering on $\mathbb{N}$.\end{remark}

\begin{proof} First of all, we notice that it suffices to prove the analogous statement to (\ref{pullbackdecomp}) for the Hilbert spaces $L^2(X,\chi), L^2(X_\Theta,\chi)$, for every $\chi\in\widehat{\mathcal Z(X)}$, and for a function $\Psi \in C^{\infty}_c(X, \chi)$. 

Let $S$ be as in the statement of the theorem, and set $S^+=S\cap \mathring A_{X,\Theta}^+$.
In order to simplify notation in what follows, we assume that $\mathcal Z(X)=1$, since the arguments are exactly the same in the general case.

Define  for $a \in S^+$  the function $\Phi_a =  e_{\Theta} \mathcal{L}_a   \Psi \in C^{\infty}_c(X)$.
Then, as $a$ approaches infinity inside $S^+$, we have
 $\|\Phi_a\|_{L^2(X)}  \rightarrow \|\Psi\|_{L^2(X_\Theta)}$;  indeed equality holds for sufficiently large $a$. 
Moreover, by definition, 	
 $H_\pi(\Phi_a)= e_\Theta^*H_{\pi}(\mathcal L_a \cdot \Psi)$.

The group  $S$ acts on $e_\Theta^* H_{\pi}$  through a finitely generated quotient, and we may therefore apply the results of \S \ref{sec:linearalgebra}.    
Let $\nu_i$ be an averaging sequence of probability measures on $S^+$ (\S \ref{averaging sequence}). Recall each of those is, by construction, of finite support. 
Then:
 
\begin{eqnarray} \label{eqnu} \lim_{i\to\infty} \int_{\hat G}  (\nu_i \star e_\Theta^*H_{\pi})  (\Psi)  \mu(\pi) =  \lim_{i\to\infty} \int_a \int_{\hat G}  e_{\Theta}^* H_{\pi} (\mathcal{L}_a \Psi)  \mu(\pi) d\nu_i(a)  \\
\nonumber =\lim_{i\to\infty} \int_a \|\Phi_a\|_{L^2(X)}^2 d\nu_i(a) = 
\|\Psi\|_{L^2(X_\Theta)}^2,\end{eqnarray}
the last step because $\|\Phi_a\|_{L^2(X)}$ and $\|\Psi\|_{L^2(X_\Theta)}$ are eventually equal.

Our task is to interchange the limit and the integral on the left hand side of \eqref{eqnu}. By Lemma \ref{nui}, $\lim_{i\to\infty}   (\nu_i \star e_\Theta^*H_{\pi})  (\Psi) = \left(e_\Theta^*H_{\pi}\right)^S(\Psi)$. Applying Fatou's lemma:
\begin{equation} \label{mf} \int_{\hat G} \left(e_\Theta^*H_{\pi}\right)^S(\Psi) \mu(\pi) \leq \|\Psi\|^2.\end{equation}

Before we continue with the proof of Theorem \ref{bernstein-abstract}, we draw a corollary from this inequality that will be used in the sequel:

\begin{corollary}\label{abscont}
 \begin{enumerate}
  \item The set of $\pi$ such that $\left(e_\Theta^*H_\pi\right)^S$ takes the value $\infty$ is of ($X$-Plancherel = $\mu$) measure zero. Therefore,  strengthening Lemma \ref{HCform} for the Plancherel forms, for almost all $\pi$ the restriction of $e_\Theta^*H_\pi$ to $C_c^\infty(X_\Theta)_\pi^1$ (the sum of \emph{generalized} eigenspaces with unitary exponents) factors through the maximal eigenquotient.
  \item The restriction of $\mu$ to the set of $\pi$ with $\left(e_\Theta^*H_\pi\right)^S\ne 0$ (equivalently: to the set of $\pi$ for which $e_\Theta^* H_\pi$ has unitary exponents) is absolutely continuous with respect to Plancherel measure on $X_\Theta$.

  \item The Plancherel measure for $X$ is absolutely continuous with respect to the sum, over all $\Theta\subset\Delta_X$, of $G$-Plancherel measures for the discrete spectra of $X_\Theta$.
 \end{enumerate}
\end{corollary}

\begin{proof}[Proof of the Corollary]
 \begin{enumerate}
  \item The first statement   is clear from (\ref{mf}), and the second follows from the fact that $S^+$ is arbitrary in $\mathring A_{X,\Theta}^+$, and that $\left(e_\Theta^*H_\pi\right)^S$ takes the value $\infty$ if its restriction to $S$-unitary generalized eigenspaces does not factor through the maximal eigenspace quotient.
  \item Both sides of (\ref{mf}) define $G$-invariant, positive semi-definite hermitian forms on $C_c^\infty(X_\Theta)$, and if $\mathcal H_l$ and $\mathcal H_r$ (for ``left'' and ``right'') denote the corresponding Hilbert spaces then we have a morphism of unitary representations: $\mathcal H_r\to \mathcal H_l$, necessarily surjective since the image of $C_c^\infty(X_\Theta)$ is dense. By \cite[\S 8]{Dixmier}, the Plancherel measure for $\mathcal H_l$ is absolutely continuous with respect to the Plancherel measure for $\mathcal H_r$.
  \item If $H_\pi\ne 0$ but has only subunitary exponents in all non-trivial directions, then $\pi$ is an $X$-discrete series. 
  
  Indeed, we show that  $H_{\pi}$ extends continuously to $L^2(X)$.For a fixed function
  $f \in C^{\infty}_c(X_{\Theta})$ and $a \in A_{X, \Theta}^+$ the quantity $H_{\pi}(e_{\Theta} \mathcal{L}_a f)$
  decays  rapidly with $a$, i.e. bounded above by $|\chi^{-1}(a)| $ where $|\chi^{-1}| < 1$ on $\mathring{A}_{X, \Theta}^+$. 
  In fact, more is true, namely  $|\chi^{-1}| < 1$ on $A_{X, \Theta}^+$.  If not, there is  a ``wall'' of $A_{X, \Theta}$
  along which $|\chi|=1$; this corresponds to an $\Omega \supset \Theta$ for which $e_{\Omega}^* H_{\pi}$
  has unitary exponents, contradicting our supposition. Since $|\chi^{-1}| < 1$ on $A_{X, \Theta}^+$, we deduce, 
  by taking $f$ to be the characteristic function of a single $J$-orbit, that $H_{\pi}$ is $L^2$-bounded on 
  the $A_{X, \Theta}^+$-span of $f$; by taking a finite set of such $\Theta$ and $f$, we  deduce that  $H_{\pi}$
  is bounded on $L^2(X)^J$, which implies that it is also bounded on $L^2(X)$ -- a $G$-invariant Hermitian form
  on a finite sum of copies of $\pi$ is uniquely determined by its restriction to $J$-invariants, for sufficiently small $J$.
  
    The restriction of Plancherel measure to the set of such $\pi$'s is, by definition, the Plancherel measure for $L^2(X)_\disc$.

 Otherwise, there is a $\Theta$ such that $\pi$ belongs to the set of representations for which $e_\Theta^* H_\pi$ has unitary exponents. Applying the second statement, for the set of those $\pi$ the statement is reduced to the analogous statement for the Plancherel measure of $X_\Theta$, and the claim follows by induction.
 \end{enumerate}
\end{proof}

  We continue with the proof of Theorem \ref{bernstein-abstract} -- we want to upgrade \eqref{mf} to an equality. 
  Let us discuss what might go wrong in order to better understand this.  Let us consider an increasing sequence
  $a_1, a_2, \dots , a_n , \dots $ in $S$ that ``go to $\infty$'' inside $S^+$;
  consider the functions $\Phi_{a_1}, \dots, \Phi_{a_n}, \dots \in C^{\infty}_c(X)$.   One could imagine that there existed a sequence of irreducible $G$-subrepresentations $\pi_1, \pi_2, \dots \subset L^2(X)$ 
  so that $\Phi_{a_j} \in \pi_j$ (or, more generally, such that a bounded below percentage of the norm of the $\Phi_{a_j}$'s is concentrated on those discrete series). In this case, the left-hand side of \eqref{mf} will be zero (or bounded away from $\Vert\Psi\Vert^2$). 
But we know that this cannot happen precisely because of finiteness of discrete series (Theorem \ref{finiteds}).   The input from Section \ref{sec:linearalgebra} generalizes this result and allows us to show
that, in all cases, \eqref{mf} may be replaced by equality:
\begin{equation} \label{mfequality} \int_{\hat G}\left(e_\Theta^*H_{\pi}\right)^S(\Psi)\mu(\pi) = \|\Psi\|^2.\end{equation}

Indeed, Corollary \ref{abscont} and Proposition \ref{uniformbound} imply that there is a uniform bound on the $S^+$-subunitary exponents for $\mu$-almost all $\pi$ with $\pi^J\ne 0$ (for some fixed open compact subgroup $J$). Moreover, the following easy lemma implies that the degree of elements of $C_c^\infty(X_\Theta)_\pi$ as $A_{X,\Theta}$-finite vectors is also uniformly bounded: 
\begin{lemma}\label{numberofexponents}
There is an integer $N$ (which, in fact, can be taken to be equal to the order of the Weyl group) such that for all irreducible representations $\pi$ the degree of all elements of $C_c^\infty(X_\Theta)_\pi$ as $A_{X,\Theta}$-finite vectors is $\le N$.
\end{lemma}

\begin{proof}[Proof of the lemma]
By the definition of $C_c^\infty(X_\Theta)_\pi$, this is the same as the degree of elements of $\Hom_G(\tilde \pi, C^\infty(X_\Theta))$. Since $X_\Theta$ is parabolically induced from $X_\Theta^L$, this space is isomorphic to $\Hom_{L_\Theta}(\tilde \pi_{P_\Theta^-}, C^\infty(X_\Theta^L))$, and since $X$ is wavefront an $L_\Theta$-morphism is also an $A_{X,\Theta}$-morphism (Proposition \ref{wavefrontlevi}). Therefore, the degree is bounded by the $\mathcal Z(L_\Theta)^0$-degree of elements of the Jacquet module $\tilde \pi_{P_\Theta^-}$.

But $\tilde\pi$ is a subquotient of a parabolically induced irreducible supercuspidal representation; therefore, the degree of any element of any Jacquet module of $\tilde\pi$ (with respect to the action of the center of the corresponding Levi) is bounded by the order of the Weyl group.
\end{proof}

Therefore, we may apply Proposition \ref{invub} to $\mu$-almost all $\pi$ with $\pi^J\ne 0$ to deduce:

\begin{quote}There are a finite set $\Lambda$ and a constant $C$ so that, {\em for all indices $i$}:
 $$\nu_i \star e_\Theta^*H_{\pi}( \Psi) \leq  C \left( \left(e_\Theta^*H_{\pi}\right)^S(\Psi) + \max_{a \in \Lambda} e_\Theta^*H_{\pi}(a \Psi)\right).$$
\end{quote}
 The right-hand side is, by \eqref{mf}, integrable. Therefore,  
we may apply the dominated convergence theorem to \eqref{eqnu}, arriving at
the desired conclusion.  This finishes the proof of Theorem \ref{bernstein-abstract}.
\end{proof}

\subsection{The Bernstein maps. Equivalence with Theorem \ref{Bernsteinmap}} \label{Bmapconstruc} 

We shall construct the desired maps of Theorem \ref{Bernsteinmap}, i.e.
$$\iota_{\Theta}: L^2(X_{\Theta} ) \longrightarrow L^2(X)$$  and prove that they have the required properties, using as input Theorem \ref{bernstein-abstract}.  Note that, 
although bounded, this map is usually not an isometry; however, see Proposition \ref{Bmapisometry}.

One should like to produce this by completing the asymptotics map $$e_{\Theta}: C^{\infty}_c(X_\Theta) \rightarrow C^{\infty}_c(X)$$ of \S \ref{sec:asymptotics} but it does not, in general, extend continuously to a map $L^2(X_{\Theta}) \rightarrow L^2(X)$ and must be modified.   
Roughly, this modification is to ``project out'' the part of $e_\Theta$ which is due to subunitary exponents (such as, but not restricted to, the projection of $e_\Theta$ to $L^2(X)_\disc$).\footnote{Here is some motivation, in a simple situation where there is only one exponent and it occurs with multiplicity onet: Suppose we are given a representation $\pi \overset{\nu}{\hookrightarrow} C^{\infty}(X)$ and 
  a {\em non-unitary} character
$\chi$ of $A_{X, \Theta}$ with the property that the embedding
$$ \pi \rightarrow C^{\infty}(X) \xrightarrow{e_{\Theta}^*} C^{\infty}(X_{\Theta})$$
transforms under  $A_{X, \Theta}$ (normalized action) by $\chi$.    Suppose, to the contrary, that $e_{\Theta}$ were $L^2$-bounded, with norm $\|e_{\Theta}\|_{\mathrm{op}}$. 
Then, for any $g \in C_c^\infty(X_{\Theta})$, $f\in \nu(\pi)$, the quantity $ \langle  f, e_{\Theta} \left(  \mathcal{L}_a g   \right) \rangle$
is bounded independent of $a \in A_{X, \Theta}$; indeed, it is bounded by
$$\|e_{\Theta}\|_{\mathrm{op}} \|f\|_{L^2(X)} \|g\|_{L^2(X_{\Theta})}.$$
But we may rewrite this expression as $\langle \mathcal{L}_a^{-1} e_{\Theta}^* f, g \rangle$, 
and this is proportional to $\chi(a^{-1})$, thus unbounded if nonzero. Contradiction. 
Thus, the failure of the asymptotics map to extend to $L^2$ is related
to the existence of discrete series, and, more generally, representations in $L^2(X)$ with subunitary exponents in the $\Theta$-direction. } 

Fix a Plancherel formula for $X$ as in (\ref{Pl}). Let $\mathcal H_\pi$ denote the completion of $C_c^\infty(X)_\pi$ under the Plancherel norms. Fix the corresponding Plancherel formula (\ref{pullbackdecomp}) for $X_\Theta$ according to Theorem \ref{bernstein-abstract}, and let $\mathcal H_\pi^\Theta$ denote the corresponding completion of $C_c^\infty(X_\Theta)_\pi$.

Recall that $\cinfc(X)_{\pi}$ splits as in (\ref{leornleq}), and obviously (by $A_{X,\Theta}$-invariance) the norm of $\mathcal H_{\pi}^{\Theta}$ factors through projection to $\cinfc(X_\Theta)_\pi^1$. 
Consider the map $\iota_{\Theta,\pi}$ obtained as the  composition:
\begin{equation}\label{iotatheta} \cinfc(X_\Theta)_\pi \stackrel{\Pi}{ \longrightarrow }\cinfc(X_\Theta)_\pi^1 \xrightarrow{e_{\Theta,\pi}} \cinfc(X)_{\pi}\end{equation}
where $\Pi$ is the $A_{X, \Theta}$-invariant projection onto $\cinfc(X_{\Theta})_{\pi}^1$, and 
we denote by $e_{\Theta,\pi}$ the map induced on $\pi$-coinvariants by the asymptotics $e_\Theta$.

The square of the norm of $\iota_{\Theta, \pi}$ with respect to $\mathcal H_{\pi}^{\Theta}, \mathcal H_{\pi}$ is  bounded by   the number of distinct exponents of $A_{X,\Theta}$ on $\cinfc(X_\Theta)_\pi^1$. 
 This is a consequence of  
Corollary \ref{subunitarybound}, applied to the hermitian forms $H=e_\Theta^* H_\pi$, $H^S= H_\pi^\Theta$.
In particular, this morphism extends to a bounded map on Hilbert spaces:
\begin{equation} \label{iotazero} \iota_{\Theta, \pi}: \mathcal{H}_{\pi}^{\Theta} \rightarrow \mathcal{H}_{\pi}.  \end{equation}

We now seek to integrate to obtain a map
 \begin{equation} \label{iotathetadef}
 \iota_\Theta= \int_{\hat G} \iota_{\Theta,\pi}.
\end{equation}

The relevant issues of measurability are handled through a straightforward application of Proposition \ref{measurable2}.

As for norms, recall that the number of generalized $A_{X,\Theta}$-eigencharacters by which an irreducible representation $\pi$ can be embedded into $C^\infty(X_\Theta)$ is bounded by the number of generalized eigencharacters in its Jacquet module with respect to a parabolic opposite to $P_\Theta$.  This is uniformly bounded (see the final quoted result 
at the end of \S  \ref{subsec:finiteness} on page \pageref{BZref}).   Thus the norms of the resulting maps $\mathcal{H}_{\pi}^{\Theta} \rightarrow \mathcal{H}_{\pi}$
are uniformly bounded.  Therefore, by Corollary \ref{subunitarybound},  \eqref{iotathetadef} gives a $G$-equivariant bounded map $$L^2(X_\Theta)\to L^2(X).$$

\begin{remark}  \label{DualBernsteinMap} Let us also discuss the description of the dual Bernstein maps.  Take $\Phi \in C^{\infty}_c(X)$. 

Recall (Remark \ref{decomp-into-eigen}) that fixing a Plancherel measure for $X$ induces  
a pointwise decomposition $\Phi = \int \Phi^{\pi} \mu(\pi)$, where each $\Phi^{\pi}\in C^\infty(X)$ is $\pi$-isotypical.  Similarly we may decompose, using the \emph{same} Plancherel measure: $$\iota_{\Theta}^* \Phi = \int_{\hat G} (\iota_{\Theta}^* \Phi)^{\pi} \mu(\pi).$$ Here we take $\Phi$ to be any smooth, $L^2$ function, and since $\iota_\Theta^*\Phi$ may not belong to a ``nice'' subspace (where the Plancherel decomposition is pointwise defined), for every $x\in X_\Theta$ the quantity $(\iota_{\Theta}^* \Phi)^{\pi}(x)$ should be thought of as an element of $L^1(\hat G, \mu)$. Then we have:

\begin{proposition} \label{cinfmap} (For $\mu$-almost all $\pi$), $(\iota_{\Theta}^* \Phi)^{\pi}$ is image of $e_{\Theta}^* (\Phi ^{\pi})$ under the $A_{X, \Theta}$-invariant projection $$C^{\infty}(X_{\Theta})^{\pi} \twoheadrightarrow C^{\infty}(X_{\Theta})^{\pi, 1}.$$ 
\end{proposition}

In words: Take the asymptotics of $\Phi^{\pi}$, and discard all ``decaying'' exponents. 

\begin{proof} Indeed, it is clear that $\iota_{\Theta}^*$ is obtained by integrating the adjoints $\iota_{\Theta, \pi}^*$ of the morphisms $\iota_{\Theta,\pi}$ of \eqref{iotazero} over $\hat{G}$.  

Now consider our construction of $\iota_{\Theta,\pi}$:

$$\xymatrix{
 C_c^\infty(X_\Theta)_\pi \ar[r]^{\Pi} \ar[d]&  C_c^\infty(X_\Theta)_\pi^1 \ar[r]^{e_{\Theta, \pi}} & C^{\infty}_c(X){\pi} \ar[d]_{e_{\Theta,\pi}}\\
 \mathcal{H}^{\Theta}_{\pi} \ar[rr]^{\iota_{\Theta, \pi}} & &  \mathcal H_\pi,
}
$$
where $\Pi$ is as before the projection to the some of generalized eigenspaces with unitary exponents.
The assertion in question follows by dualizing the entire diagram.

\end{proof}

\end{remark}

\subsection{Property characterizing the Bernstein maps}   
 Let us now verify that $\iota_{\Theta}$ has the property described in Theorem \ref{Bernsteinmap}.  More precisely, we will prove:
 
\begin{lemma}\label{decaying}
There is a decaying function $Q^J$ on $\mathbb{Z}^+$, depending only on $J$, such that:   \begin{equation} 
 \Vert \iota_\Theta \mathcal L_{a^n}\Psi- e_\Theta \mathcal L_{a^n} \Vert \le C_\Psi Q^J(n)
\end{equation}
for any $\Psi\in C_c^\infty(X_\Theta)^J$ and some constant $C_\Psi$ depending on $\Psi$.
\end{lemma}
 
\begin{proof}
 
Fix $\Psi\in C_c^\infty(X_\Theta)^J$ and for $a\in \mathring A_{X,\Theta}^+$, denote  $\Psi_a:=\mathcal L_a \cdot \Psi$.

Recall from Lemma \ref{HCform} that the pull-back $e_{\Theta}^* H_{\pi}$ factors through the sum: $C_c^\infty(X_\Theta)_\pi^{<1}\oplus C_c^\infty(X_\Theta)_\pi^1$.  Let $H_{\pi}^{<1}$ be the pull-back of $H_\pi$ to $C_c^\infty(X_\Theta)$ via the composition:
\begin{equation} \label{overcomplicated} \cinfc(X_{\Theta})_{\pi} \twoheadrightarrow \cinfc(X_{\Theta})_{\pi}^{<1} \hookrightarrow \cinfc(X_{\Theta})_{\pi} \xrightarrow{e_{\Theta,\pi}} \cinfc(X)_{\pi}\end{equation}
In other words, if we write $\Psi_a^1$ for the image of $\Psi_a$ via $C_c^\infty(X_\Theta)\twoheadrightarrow C_c^\infty(X_\Theta)^1$, we have $H_\pi^{<1}(\Psi_a) = e_{\Theta,\pi}^*H_\pi(\Psi_a - \Psi_a^1)$. Notice, moreover, that by definition $H_\pi(\iota_\Theta(\Psi_a)) = e_\Theta^* H_\pi (\Psi_a^1)$.   (But the latter is \emph{not} equal, in general, to $H_\pi^\Theta(\Psi_a)$ as different unitary eigenspaces of $C_c^\infty(X_\Theta)_\pi$ may \emph{not} be orthogonal under $e_{\Theta,\pi}^* H_\pi$.)

Then, by following the definitions:
 $$\Vert \iota_{\Theta} (\Psi_a) - e_\Theta(\Psi_a) \Vert_{L^2(X)}^2 = \int_{\hat G} H_\pi\left(\iota_\Theta(\Psi_a)-e_\Theta(\Psi_a)\right) \mu(\pi) = $$
$$ \int_{\hat G} e_\Theta^*H_\pi \left(\Psi_a^1-\Psi_a\right) \mu(\pi) = \int_{\hat G} H_{\pi}^{<1}(\Psi_a)\mu(\pi).$$
So, it suffices to check the latter integral is bounded by a multiple of a decaying function $Q^J$.
But that follows from (the second inequality of) Corollary \ref{subunitarybound}.
Notice that the corollary applies for the same reasons as in the proof of Theorem \ref{bernstein-abstract}, that is because of Proposition \ref{uniformbound} and Lemma \ref{numberofexponents}. In the present context  it asserts that $$H_{\pi}^{<1}(\Psi_a) \leq Q^J(a) \left( H^\Theta_{\pi}(\Psi) + \sum_{i} H_{\pi}(e_\Theta\Psi_{a_i})\right),$$ 
where the decaying function $Q^J$ depends only on $J$ (since this subgroup determines bounds for the number and growth of subunitary exponents) and for some finite collection $\{a_1, \dots, a_m\}$. By integrating over $\pi$, we get:
$$ \Vert \iota_{\Theta} (\Psi_a) - e_\Theta(\Psi_a) \Vert_{L^2(X)}^2  \le Q^J(a) \cdot \left(\Vert \Psi\Vert_{L^2(X_\Theta)}^2+ \sum_i \Vert e_\Theta\Psi_{a_i}\Vert_{L^2(X)}^2\right).$$
\end{proof}

To finish the proof of Theorem \ref{Bernsteinmap}, there remains to check the {\em uniqueness} of a map $\iota_{\Theta}$ with the (weaker than that of the previous lemma) property \eqref{bchar}. 

Suppose $\iota_{\Theta}, \iota_{\Theta}'$ were two morphisms with that property. Their difference $\delta_\Theta$ then has the property that:
$$\| \delta_\Theta \mathcal L_{a^n} \Psi\|_{L^2(X)}\xrightarrow{n} 0,$$
for every $\Psi \in C_c^\infty(X_{\Theta})$.   

The quantity $\|\delta_\Theta \Psi\|_{L^2(X)}$ defines a $G$-invariant Hilbert seminorm
on $L^2(X_{\Theta})$, bounded by $C \|\Psi\|_{L^2(X_\Theta)}$ for some positive $C$.  We may disintegrate it as $\int_{\pi} N_{\pi}(\Psi)\mu(\pi)$, where 
$N_{\pi}$ is a $G$-invariant square-seminorm on the Hilbert space $\mathcal H^\Theta_{\pi}$
satisfying $N_{\pi}(\Psi) \leq C \cdot H_\pi^\Theta(\Psi)$.

Now $\int_{\hat G} N_{\pi}(\mathcal L_{a^n} \Phi)\mu(\pi) \rightarrow 0$. Reasoning as for \eqref{mfequality}, 
the associated $S$-invariant norms satisfy $\int_{\hat G} N^S_{\pi}(\Psi)\mu(\pi) = 0$, i.e.
$N^S_{\pi} = 0$ for almost all $\pi$.   
Since the function $a \mapsto N_{\pi}\left( \mathcal L_a \Psi \right)$
is bounded (by $C \cdot H_\pi^{\Theta}(\mathcal{L}_a \Psi) = C \cdot H_{\pi}^{\Theta}(\Psi)$), we deduce by Lemma \ref{boundU} that $N_{\pi}(\Psi) = 0$ for almost all $\pi$, as desired.

\subsection{Compatibility with composition and inductive structure of $L^2(X)$}

\begin{proposition} \label{propcomposition}
 For each $\Omega\subset\Theta\subset\Delta_X$, let $\iota_\Omega^\Theta$ denote the analogous Bernstein morphism: $L^2(X_\Omega)\to L^2(X_\Theta)$. Then:
\begin{equation}
 \iota_\Theta\circ \iota_\Omega^\Theta = \iota_\Omega.
\end{equation}
\end{proposition}

\begin{proof}
 This follows from the analogous result on the ``naive'' asymptotics maps:
$$ C_c^\infty(X_\Omega) \xrightarrow{e_\Omega^\Theta} C_c^\infty(X_\Theta) \xrightarrow{e_\Theta} C_c^\infty(X).$$
The composition of these arrows is equal to $e_\Omega$, cf.\ Remark \ref{transitivity}. 

Specializing to $\pi$-coinvariants, and taking into account that $C_c^\infty(X_\Omega)_\pi^1$ maps into $C_c^\infty(X_\Theta)_\pi^1$ (the restriction of a unitary character of $A_{X,\Omega}$ to the subtorus $A_{X,\Theta}$ remains unitary), we get the result by the definition (\ref{iotatheta}) of $\iota_{\Theta,\pi}$.
\end{proof}

\begin{corollary}\label{corollarysurjection}
 Let $L^2(X)_\Theta$ be the image\footnote{For $L^2$-spaces the index $\Theta$ is being used to denote the image of the discrete spectrum of boundary degenerations, while otherwise it is used to denote Jacquet modules. We hope that this will not lead to any confusion, since we do not use Jacquet modules in the category of unitary representations.}\ of $L^2(X_\Theta)_\disc$ under $\iota_\Theta$. Then:
\begin{equation}
 \sum_{\Theta\subset\Delta_X} L^2(X)_\Theta = L^2(X).
\end{equation}
\end{corollary}

\begin{proof}
 Assume the statement to be true if we replace $X$ by any $X_\Theta$, $\Theta\subsetneq \Delta_X$. Then, by Proposition \ref{propcomposition}, the orthogonal complement $\mathcal H'$ of $ \sum_{\Theta\subsetneq\Delta_X} L^2(X)_\Theta $ is orthogonal to $\iota_\Theta\left(L^2(X_\Theta)\right)$ for all $\Theta\subsetneq \Delta_X$. By the definition of $\iota_\Theta$, this means that the space $\mathcal H'$ admits a Plancherel decomposition:
$$\mathcal H'= \int_{\hat G} \mathcal H_\pi'\mu(\pi)$$
where the norms for all $\mathcal H_\pi'$ are decaying in all directions at infinity (i.e.\ they have only subunitary, no unitary exponents, for every $\Theta\subsetneq \Delta_X$). But then $\mathcal H'\subset L^2(X)_\disc = L^2(X)_{\Delta_X}$ by the generalization of Casselman's square integrability criterion, cf.\ \S \ref{newredstuff}. 
\end{proof}

\subsection{Isometry} 

As we mentioned, the Bernstein map $i_\Theta:L^2(X_\Theta)\to L^2(X)$ is not, in general, an isometry; in section \ref{sec:scattering} we will examine its kernel. However, it is an isometry if we restrict to a small enough subspace of $L^2(X_\Theta)$:

\begin{proposition}\label{Bmapisometry}
 Let $\mathcal H'\subset L^2(X_\Theta)$ be an $A_{X,\Theta}\times G$-stable subspace. Fix a Plancherel measure $\mu$ for $L^2(X)$ and corresponding direct integral decompositions for $\mathcal H:= L^2(X)$ and $\mathcal H'$:
$$ \mathcal H = \int_{\hat G} \mathcal H_\pi \mu(\pi),$$ 
$$ \mathcal H' = \int_{\hat G} \mathcal H'_\pi \mu(\pi).$$
Assume that for almost all $\pi$ the following is true: if $\mathcal H'_\pi = \oplus_\chi \mathcal H'_{\pi,\chi}$ is a decomposition into $A_{X,\Theta}$-generalized eigenspaces (necessarily, for almost all $\pi$, honest eigenspaces with unitary characters), then for distinct characters $\chi_i$ the images of $\mathcal H'_{\pi,\chi_i}$ via the Bernstein maps $\left.\iota_{\Theta,\pi}\right|_{\mathcal H'_\pi}: \mathcal H'_\pi\to \mathcal H_\pi$ are mutually orthogonal.

Then the restriction of the Bernstein map $\iota_\Theta$ to $\mathcal H'$ is an isometry onto its image.
\end{proposition}

\begin{remark} \label{UniqueExponentSituation}
 The proposition applies, in particular, to the case that (almost) every $\mathcal H'_\pi$ has a unique exponent. {\em This is the only setting in which we will use it.}
\end{remark}

\begin{proof}
 Having established the existence (and boundedness) of the morphisms $\iota_\Theta$ (and their disintegrations into $\iota_{\Theta,\pi}: \mathcal H_{\Theta,\pi}\to \mathcal H_\pi$, where $\mathcal H_{\Theta,\pi}$ is the disintegration of $L^2(X_\Theta)$ with respect to $\mu$; it is equipped therefore with a ``Plancherel'' hermitian norm), we may, a posteriori, rephrase the conclusion of Theorem \ref{bernstein-abstract} in terms of them (here $S$ is as in that theorem):

\begin{quote} 
Let $H_\pi$ denote the Hermitian form on $\mathcal H_\pi$, and $\iota_{\Theta,\pi}^* H_\pi$ its pull-back to $\mathcal H_{\Theta,\pi}$. The associated $S$-invariant form $\left(\iota_{\Theta,\pi}^* H_\pi\right)^S$ is equal to the ``Plancherel'' hermitian form on $\mathcal H_{\Theta,\pi}$.
\end{quote}

Indeed, let us further pull back these norms to $C_c^\infty(X_\Theta)_\pi$ via the canonical map: $C_c^\infty(X_\Theta)_\pi \to \mathcal H_{\Theta,\pi}$. By definition of $\iota_\Theta$, we have a commutative diagram:
$$\begin{CD}
 C_c^\infty(X_\Theta)_\pi @>>> C_c^\infty(X_\Theta)_\pi^1 @>>> \mathcal H_{\Theta,\pi} \\
&& @V{e_{\Theta,\pi}}VV @VV{\iota_{\Theta,\pi}}V \\
&& C_c^\infty(X)_\pi @>>> \mathcal H_\pi
\end{CD}
$$ 

Hence, our current pull-backs are obtained from the pull-backs of Theorem \ref{bernstein-abstract} (induced by $e_{\Theta,\pi}: C_c^\infty(X_\Theta)_\pi\to C_c^\infty(X)_\pi$) by composing with the $A_{X,\Theta}$-equivariant projection to $C_c^\infty(X_\Theta)_\pi^1$. But the process of taking $S$-invariants also factors through this projection, hence $\left(\iota_{\Theta,\pi}^* H_\pi\right)^S$ coincides with $\left(e_{\Theta,\pi}^* H_\pi\right)^S$ (as a hermitian form on $C_c^\infty(X_\Theta)_\pi$) and hence, by Theorem \ref{bernstein-abstract}, with the Plancherel hermitian norm on $\mathcal H_{\Theta,\pi}$.

The assumptions on $\mathcal H'$ now imply that the restriction of $\iota_{\Theta_\pi}^* H_\pi$ to $\mathcal H'_\pi$ is already $A_{X,\Theta}$-invariant. Hence, it coincides with the Plancherel hermitian form on $\mathcal H'_\pi$ (we implicitly use here that the construction of the
 ``associated invariant norm'' from Lemma \ref{invariantform} is compatible with passage to $S$-invariant subspaces; this is clear from the definition), and therefore the Bernstein map is an isometry when restricted to $\mathcal H'$.
\end{proof}

\section{Preliminaries to scattering (I): direct integrals and norms} \label{sec:directintegrals}

In this section and the next we gather some useful results, presented in an abstract setting, that will be used in \S \ref{sec:scattering}. 

\begin{itemize}

\item[-] \S \ref{ssDixmierRecall} recalls the general formalism of direct integrals of Hilbert spaces, which is essential for the Plancherel decomposition.

\item[-]  \S \ref{sshsnorms} discusses certain norms on direct integrals of Hilbert spaces; these norms will be used extensively in
\S \ref{sec:scattering}, in particular, \S \ref{ssestimates}. 

Roughly speaking, in \S \ref{ssestimates}, we will have available pointwise bounds on eigenfunctions, and we obtain pointwise bounds on general functions by first decomposing
into eigenfunctions and then applying this pointwise bounds; the norms that we discuss  are abstractions of this process. 

\end{itemize} 
\subsection{General properties of the Plancherel decomposition} \label{ssDixmierRecall}

\subsubsection{}
\label{trivialityremark}
Let $\mathcal H$ be a unitary representation of $G$. We discussed in \S \ref{ssdirectintegrals} the meaning of a Plancherel decomposition:
$$\mathcal H = \int_{\hat G} \mathcal H_\pi \mu(\pi).$$
By \cite{Dixmier}[Theorem 8.6.6 and A72], under the assumption that almost all $\mathcal H_\pi$ are of finite multiplicity, there is a countable partition of the unitary dual $\hat G$ into measurable sets $Z_i$, and hence a corresponding direct sum decomposition $\mathcal H= \oplus_i \mathcal H_i$, such that each $\mathcal H_i \simeq \mathcal H_i' \otimes V_i$ as a $G$-representation, where:
\begin{enumerate}
 \item  $V_i$ is a finite dimensional vector space of dimension $i$, with trivial $G$-action;
 \item $\mathcal H_i'$ is multiplicity-free, that it it admits a direct integral decomposition: $\mathcal H_i'\simeq \int_{Z_i} \mathcal H_{\pi}' \mu_i(\pi)$ with $\mathcal H_\pi'$ \emph{irreducible};
 \item the measures $\mu_i$ are mutually singular;
 \item the measurable structure is trivializable,\footnote{We are assuming here that all $\mathcal H_\pi'$ are infinite-dimensional; in general, the measurable structure is trivializable over the (measurable) subsets where $\mathcal H_\pi'$ has fixed dimension.} i.e.\ there is a Hilbert space $H_0$ and isomorphisms of Hilbert spaces: $\mathcal H_\pi'\xrightarrow{\sim} H_0$ inducing a bijection between the collection of measurable sections $\pi\mapsto \eta'_\pi\in \mathcal H_\pi'$ and the collection of measurable sections $\pi\mapsto \eta_\pi\in H_0$.
\end{enumerate}

The results cited below concerning unitary decomposition may be found in \cite{Dixmier}, in particular, Theorem 8.5.2 and 8.6.6 (existence and uniqueness of unitary decomposition)
and Proposition 8.6.4 (characterization of $G$-endomorphisms).

\subsubsection{Uniqueness of unitary decomposition}

Suppose that two unitary representations with Plancherel decompositions: $\int_{\hat G} \mathcal{H}_{\pi} \mu(\pi)$ and $\int_{\hat G} \mathcal{H}_{\pi}' \nu(\pi)$ are isomorphic. Then the measure classes of $\mu$ and $\nu$ are equal, and moreover, there exists an isometric isomorphism $\mathcal{H}_{\pi} \rightarrow \mathcal{H}'_{\pi}$ for $\mu$-almost every (equivalently: $\nu$-almost every) $\pi \in \hat{G}$.

\subsubsection{Endomorphisms} \label{Endo}

Notation as before. A family of (bounded) endomorphisms $\pi\mapsto T_\pi:\mathcal H_\pi\to\mathcal H_\pi$ is called \emph{measurable} if for every measurable section $\pi\mapsto\eta_\pi\in\mathcal H_\pi$ the section $\pi\mapsto T_\pi \eta_\pi$ is measurable (see \cite[A78]{Dixmier}). Any $G$-endomorphism $f$ of $\int_{\hat G} \mathcal{H}_{\pi} \mu(\pi)$ is ``decomposable,'' that is to say, there is a measurable family of $G$-endomorphisms $f_\pi$ of $\mathcal H_\pi$ such that $f(v) = \int_{\hat G} f_\pi(v_\pi) \mu(\pi)$ for $v=\int v_\pi \mu(\pi)$. We will symbolically write:
$$f = \int_{\hat G} f_\pi.$$
This assertion follows\footnote{We outline the argument: In the setting of \S \ref{trivialityremark}, there are no non-trivial $G$-morphisms between the different summands $\mathcal H_i$ (\emph{loc.cit.} Proposition 8.4.7). Now, since $V_i$ is finite-dimensional, we have $\End(\mathcal H_i) = \End(\mathcal H_i')\hat\otimes \End(V_i)$, and hence $\End(\mathcal H_i)^G = \End(\mathcal H_i')^G\hat\otimes \End(V_i)$. It now follows from Proposition 8.6.4 that the first factor consists precisely of the ``diagonalizable'' endomorphisms, i.e.\ those which are direct integrals of scalars in the $\mathcal H_\pi'$'s.} from Proposition 8.6.4 and Theorem 8.6.6 of \cite{Dixmier}. We shall refer to this as ``disintegration of endomorphisms.''

\subsubsection{Disintegration of morphisms} \label{ss:dis} 
Let $\mathcal{J} := \mathcal{H} \oplus \mathcal{H}'$ be the direct sum of two unitary $G$-representations. Let $\mu$ be a Plancherel measure for $\mathcal{J}$, so we may disintegrate
$$\mathcal{J}= \int_{\hat G} \mathcal{J}_\pi \mu(\pi).$$

We claim that there are measurable subfields $\mathcal{H}_\pi, \mathcal{H}'_\pi \subset \mathcal{J}_\pi$ (that is, Hilbert subspaces so that the corresponding projections are measurable) so that:
\begin{itemize}
 \item  $\mathcal{J}_\pi = \mathcal{H}_\pi \oplus \mathcal{H}'_\pi$;
\item $\mathcal{H}=\int_{\hat G} \mathcal H_\pi \mu(\pi)$ and $\mathcal{H}'=\int_{\hat G} \mathcal H'_\pi\mu(\pi)$.
\end{itemize} 
Notice that $\mu$ is not necessarily a Plancherel measure for $\mathcal H$ or $\mathcal H'$, as the space $\mathcal H_\pi$, $\mathcal H'_\pi$ could be zero for $\pi$ in a non-zero set.

To see that this is true, we disintegrate the projections of $\mathcal{J}$ onto $\mathcal{H}$ and $\mathcal{H}'$ to obtain measurable families of projections of $\mathcal{J}_z$; we define $\mathcal{H}_z$ and $\mathcal{H}'_z$  as the images of these projections. 

Let us now discuss the analog of \S \ref{Endo} for morphisms with different source and target. 
Thus let $\mathcal{H}, \mathcal{H}'$ be unitary $G$-representations and $f: \mathcal{H} \rightarrow \mathcal{H}'$ a $G$-morphism. We wish to ``disintegrate'' $f$.

This is reduced to the previous setting by replacing $f$ by the endomorphism $f+0$ of $\mathcal{J} := \mathcal{H} \oplus \mathcal{H}'$. Then for a decomposition $\mathcal J = \int_{\hat G} (\mathcal H_\pi\oplus\mathcal H'_\pi) \mu(\pi)$ as above, we get morphisms $f_\pi: \mathcal J_\pi\to \mathcal J_\pi$ for almost every $\pi$, which disintegrate $f+0$. Since, however, $\mathcal H'$ is in the kernel of $f+0$ and the image is contained in $\mathcal H'$, it follows that for $\mu$-almost every $\pi$ the morphism $f_\pi$ factors as:
\begin{equation}
 f_\pi: \mathcal H_\pi \to \mathcal H'_\pi.
\end{equation}

This is what we will mean by disintegration of a morphism $f:\mathcal H\to\mathcal H'$; a disintegration with respect to a Plancherel measure for $\mathcal H \oplus \mathcal H'$.

\begin{remark}
 It can easily be shown that the class of Plancherel measure of $\mathcal H \oplus \mathcal H'$ is precisely the class of a sum of Plancherel measures for $\mathcal H$ and $\mathcal H'$.
\end{remark}

\subsection{Norms on direct integrals of Hilbert spaces} \label{sshsnorms} \footnote{We use the word ``norm'' freely in this section, to include seminorms that are not necessarily bounded -- i.e., can be zero or infinite on some nonzero vectors. In other words, a norm on a vector space
 $V$ is a pair of a subspace $V_f \subset V$ -- the ``space of vectors of finite norm''--
 and a seminorm $N: V_f \rightarrow \mathbb{R}_{\geq 0}$). By convention,
 in this setting, we write $\|v\| = \infty$ for $v \in V - V_f$. } 

The reader may wish to postpone this section until reading Lemma \ref{toy-lemma}, which 
gives the basic case where the norms discussed here arise.

\subsubsection{Basic example} \label{mixednormdefrankone}

Let $(Y,\mathfrak B, \mu)$ be a (positive) measure space, denote by $F(Y)$ the space of measurable functions modulo essentially zero functions, and consider the corresponding space $L^2(Y,\mu)$. Every $v\in L^2(Y,\mu)$ corresponds to a signed measure $v\mu$ and a  
positive measure $|v|\mu$, thus defining an $L^1$-seminorm $\Vert\bullet\Vert_{L^1(Y,|v|\mu)}$ which is continuous on $L^2(Y,\mu)$. Explicitly, this seminorm is defined by
$$ h \in L^2(Y, \mu) \mapsto \int |h| |v| \mu,$$
which is evidently bounded by $\|v\| \cdot \|h\|$.

We take this observation a step further, and consider instead a direct integral of (non-trivial) Hilbert spaces over a measure space:
$$ \mathcal H=\int_Y \mathcal H_\pi \mu(\pi).$$
Again, given $v=\int_Y v_\pi\mu(\pi)\in \mathcal H$ we can define a corresponding $L^1$-norm on $\mathcal H$, namely:
\begin{equation}\label{L1norm}
 \left\Vert \int_Y h_\pi \mu(\pi)\right\Vert_{L^1_v} := \int_Y |\left<h_\pi,v_\pi\right>| \mu(\pi).
\end{equation}

Again, this is bounded in operator norm by $\|v\|$.  
On the other hand we have $\langle h, v \rangle \leq \|h\|_{L^1_v}$.

Notice that this norm depends only on $Y,\mathfrak B, \mathcal H, v$ and not on $\mu$, in the following sense: If we multiply $\mu$ by an almost everywhere positive function and divide the hermitian forms on the spaces $\mathcal H_\pi$ by the same function, then we have a canonical isomorphism of the new direct integral with $\mathcal H$, and the element corresponding to $v$ defines the same norm on $\mathcal H$.

\subsubsection{The case of a $G$-representation}

Suppose  $\mathcal{H}$ is an arbitrary unitary $G$-representation, and $(Y, \mu)$ is defined by Plancherel decomposition for $\mathcal{H}$. In this setting the relative norm can be described thus:  

The ring of essentially bounded Borel measurable functions on $\hat{G}$  
acts on any unitary $G$-representation by bounded $G$-endomorphisms.   
Then, for $x \in \mathcal{H}$, 
\begin{equation}\|x \|_{L^1_v} = \sup_{E} \|\langle x, E v \rangle \|,
\end{equation}
where $E$ ranges through  Borel measurable functions on $\hat{G}$   satisfying $|E(\pi)| \leq 1$ for all $\pi\in\hat G$.

 This description makes manifest the following: 
given a bounded $G$-morphism $f: \mathcal{H}_1 \rightarrow \mathcal{H}_2$, then, 
for any $w \in \mathcal{H}_2$, the pull-back by $f$ of the norm $L^1_{w}$
	is equal to the norm $L^1_{f^* w}$, where $f^*:\mathcal H_2\to\mathcal H_1$ denotes the adjoint of $f$.

\subsubsection{Relativization} \label{mixednormdefrelative}

We also introduce ``relative'' versions of the above norms. Let $p$ be a morphism of measure spaces: 
\begin{equation}\label{mapmeasurespaces}(Y,\mathfrak B,\mu)\to (Y',\mathfrak B',\nu)
\end{equation}
Assume that the direct image of $\mu$ under $p$ is absolutely continuous with respect to $\nu$ (i.e., $\mu$ is  zero on inverse images of null sets). 

We also need to make certain assumptions
on the measures we are dealing with; for the applications that we have in mind, it suffices to assume that $\mu$ is compact (see below) and $\nu$ is $\sigma$-finite. 

\begin{remark} (Measure-theoretic details).  \label{compactmeasure}

Recall \label{compactdefinition} from \cite[\S 451]{Fremlin} that a compact measure on a given sigma-algebra is one which is \emph{inner regular} with respect to a \emph{compact class} of subsets. Inner regular means that the measure of a measurable set $A$ is the supremum of the measures of the subsets of $A$ in the given class, and a compact class is a collection $\mathcal K$ of (measurable) subsets such that $\bigcap_{K\in \mathcal K'} K\ne 0$ whenever $\mathcal K'\subset \mathcal K$ has the finite intersection property.

In our applications, these conditions will be almost automatic: 
  these measure spaces will arise from a Plancherel decomposition either of a subspace of $L^2(X_\Theta)$, or of a locally compact abelian group.  These are measures on
  standard Borel spaces (cf. \cite[4.6.1, 7.3.7]{Dixmier}). They may be assumed to be finite:
  by the definition of Dixmier, these measures are $\sigma$-finite, and therefore can be replaced by finite measures in the same measure class. Finally, a finite measure on a standard Borel space
  is automatically inner regular and so compact \cite[434J (g)]{Fremlin}.  
\end{remark}

Under the assumptions above,  there is a family of measures $\{\mu_\rho\}_{\rho\in Y'}$ on $(Y,\mathfrak B)$ such that $\mu_\rho(A)= \mu_\rho(A\cap p^{-1}\{\rho\})$ for every $A\in \mathfrak B$ and:
$$\int_Y f(\pi)\mu(\pi) = \int_{Y'} \left(\int_Y f(\pi)\mu_\rho(\pi)\right) \nu(\rho).$$
for every measurable function $f$ on $Y$. This is the \emph{disintegration of measures}, see, for instance, \cite{Fremlin}[Theorem 452I].

Let $\mathcal H=\int_Y \mathcal H_\pi \mu(\pi)$ as above, then from $Y'$ we get a coarser decomposition of $\mathcal H$:
$$\mathcal H=\int_{Y'} \mathcal H_\rho \nu(\rho)$$
where $\mathcal H_\rho = \int_Y \mathcal H_\pi \mu_\rho(\pi)$.

Given a vector $v\in \mathcal H$ we can now define a norm on $\mathcal H$ which is a mixture of the above norms along the fibers of $Y\to Y'$ and the $L^2$-norm along $Y'$, more precisely:
If $v=\int_{Y'} v_\rho \nu(\rho)$ then 
\begin{equation}
 \left\Vert \int_{Y'} h_\rho \nu(\rho) \right\Vert_{Y,Y',\nu, v} := \left(\int_{Y'} \left\Vert h_\rho\right\Vert_{L^1_{v_\rho}}^2 \nu(\rho)\right)^\frac{1}{2}.
\end{equation}

This norm is continuous on $\mathcal H$ if $\Vert v_\rho\Vert_{\mathcal H_\rho}$ is essentially bounded in $\rho$ -- indeed, it is bounded by $\mathrm{sup}_{\rho}  \Vert v_\rho\Vert_{\mathcal H_\rho}$ times the norm on $\mathcal H$, where $\sup$ means ``essential supremum'' -- 
but not in general. Notice also the following: the norms do not depend on $v$ itself, but rather the collection of $v_\rho$ up to multiplying each by a scalar of norm one. Finally, \emph{they do depend on the choice of $\nu$}, albeit not on the choice of $\mu$ (in the sense described above, i.e.\ modifying the measure and the hermitian norms accordingly).

\begin{lemma} \label{lem:bounds} Suppose that each $\mathcal{H}_\rho$, $\rho\in Y'$, is infinite-dimensional; then 
the norm on $\mathcal{H}$ is not bounded by any 
finite sum of norms of the form $\Vert\bullet\Vert_{Y,Y',\nu, v}$.
\end{lemma}

\begin{proof}
In fact, given any $v_1, \dots, v_N \in \mathcal{H}$, we may find a nonzero vector $w \in \mathcal{H}$
with the property that $w_\rho \perp (v_j)_\rho$ for all $\rho \in Y'$. 
Then $\|w\|_{Y, Y', \nu, v_j} =0$ for each $j$. 
\end{proof}

\subsubsection{The case of $G \times A$-representations}

Our use of relative norms will be in the following situation: $G$ our fixed reductive $p$-adic group, $A$ a discrete abelian group, 
and $\mathcal{H}$ a unitary $G \times A$-representation whose Plancherel measures under $G\times A$ and $A$ satisfy the assumptions for disintegration; we take $Y = \widehat{ G \times A}, 
Y' = \widehat{A}$.  For each choice of $A$-Plancherel measure $\nu$ and each $v\in\mathcal H$, the following is clear:
\begin{lemma}
The norm $\|\cdot\|_{\widehat{G \times A}, \widehat{A}, \nu, v}$ is $A$-invariant.  
\end{lemma}
We assume that the Plancherel measure for $\mathcal H$ as an $A$-representation is absolutely continuous with respect to Haar measure on $\hat A$, and equip $Y'$ with $\nu$=Haar probability measure.

{\em In this setting,  denote the relative norm  $\|\cdot\|_{\widehat{G \times A}, \widehat{A}, \nu, v}$  by $\| \cdot\|_{A, v}$ for short.}

The following Lemma will play a key role in our later proofs. 
\begin{lemma} \label{Sumofnorms}
  Given $x_1, \dots, x_r \in \mathcal{H}$, and any corresponding collection of \emph{proper} subgroups
$T_1, \dots, T_r \subset A$, 
the Hilbert norm on  $\mathcal{H}$ is not majorized by $\sum \| \bullet \|_{T_j,x_j}$. 
\end{lemma}
\begin{proof} 

Assume to the contrary, then by scaling the $x_j$'s we may assume that $\Vert \bullet\Vert_{\mathcal H}\le \sum \| \bullet \|_{T_j,x_j}$. 
We have already seen that the relative norm $\| \bullet \|_{T_j,x_j}$
is bounded by at most $\sup_{\rho \in \widehat{T_i}}  \|x_i\|_{\rho}$ times the norm of $\mathcal H$. 
Recall that 
$ \int_{\rho \in \widehat{T_i}} \|x_i\|_{\rho}^2 =  \|x\|^2$, but
this gives no control on the supremum needed to bound the relative norm. 

Now let $\mathcal{H}' \subset \mathcal{H}$ be a $G \times A$-invariant subspace.
Then the restriction of $\| \bullet \|_{T_j,x_j}$
to $\mathcal{H}'$ is simply given by $\| \bullet \|_{T_j,\overline{x_j} }$,
where $\overline{x_j}$ is the projection of $x_j$ to $\mathcal{H}$. (This follows from the subsequent Lemma \ref{lemma1067}, applied
to $f$ the inclusion $\mathcal{H}' \hookrightarrow \mathcal{H}$.) 

It follows that it suffices to construct a nonempty $G \times A$-invariant subspace $\mathcal{H}'$
with the property that
$$  \| \overline{x_j}\|_{\rho} < \frac{1}{r}  \ \ \mbox{ for all $j$ and all $\rho \in \widehat{T}_j$}.$$
Then we have: $\Vert \bullet \Vert_{\mathcal H} \le \sum \| \bullet \|_{T_j,x_j} < r \cdot \frac{1}{r} \Vert \bullet\Vert_{\mathcal H}$, a contradiction.

 The strategy is to take $\mathcal{H}'$ to be the image   of the orthogonal projection $1_S$ induced by a  measurable subset $S \subset \widehat{A}$. 
 Then the orthogonal projection $\overline{x_j}$ of $x_j$ to $\mathcal{H'}$
 is simply $1_S \cdot x_j$; on the other hand, $\mathcal{H}'$ is nonempty
 so long as $S$ has positive Haar measure, by virtue of our assumption on the Plancherel measure of $\mathcal{H}$ with respect to $A$.  

 It suffices, then, to construct a set $S$ with the property that
\begin{equation} \label{wanted:S}  \| 1_S x_j\|_{\rho}^2 < \frac{1}{r^2} \end{equation} 
for every  $1 \leq j \leq r$ and every $\rho \in \widehat{T}_j$. 

The function $x_j \mapsto \| x_j \|_{\rho}$
is (in the notation of \cite[Appendix A]{Dixmier}) ``$\mu$-measurable'' where $\mu$ is the Haar measure, that is to say, 
in the complection of the Borel $\sigma$-algebra with respect to $\mu$.

Notice that $\widehat{A}$ has the structure of a compact abelian Lie group. Fixing any Riemannian metric on it, we can speak of balls $S(\varepsilon)\subset \widehat A$ of radius $\varepsilon$ around a point. In what follows, let us fix the Haar measures on $\widehat{A}$ and $\widehat{T_j}$
to be probability measures, and then each fiber of $\widehat{A} \rightarrow \widehat{T_j}$
is also endowed with a natural fibral probability measure (indeed, this fiber may be identified with the dual of the discrete group $A/T_j$ in a natural way). 

By the Lebesgue differentiation theorem, for Haar-almost every point of $\widehat A$ (taken as the center of the balls $S(\varepsilon)$) there is a constant $C$ such that for $\varepsilon$ sufficiently small and every $j$: 
\begin{equation} \label{salem} \int_{S(\varepsilon)} \|x_j\|_\chi^2 d\chi \leq C \cdot  \Vol(S(\varepsilon)).\end{equation} 
The left hand side can be written as: 
$$ \int_{\widehat{T_j}}   \| 1_{S(\varepsilon)} x_j\|_{\rho}^2 d\rho = \int_{S(\varepsilon)}  \| 1_{S(\varepsilon)} x_j\|_{\rho(\chi)}^2 f_j(\varepsilon,\chi)^{-1} d\chi, $$
where $f_j(\varepsilon,\chi)$ is the fibral volume of $S(\varepsilon)$ over $\rho(\chi)\in \widehat{T_j}$.
Clearly, $f_j(\varepsilon,\chi)\le C' \varepsilon$ for some constant $C'$, hence:
\begin{quotation}
 \emph{there is a subset $S\subset S(\varepsilon)$ of positive measure with $\| 1_{S(\varepsilon)} x_j\|_{\rho(\chi)}<\frac{1}{r}$ for all $\chi\in S$, $j=1,\dots ,r$.} In particular, $\| 1_{S} x_j\|_{\rho(\chi)}<\frac{1}{r}$ for all $\chi\in S, j=1,\dots,r$.
\end{quotation}

Indeed, the estimates above show that: 
$$\sum_j \int_{S(\varepsilon)} \| 1_{S(\varepsilon)} x_j\|_{\rho(\chi)}^2 d\chi \le rCC'\varepsilon \Vol(S(\varepsilon))$$
and we can choose $\varepsilon$ small enough so that $rCC'\varepsilon<\frac{1}{r^2}$.
This provides the desired set and proves the lemma.
\end{proof}

Now $A_1, A_2$ be two discrete abelian groups and $T: A_2 \rightarrow A_1$ a morphism. 
Let $\mathcal{H}_1$ and $\mathcal{H}_2$ be, respectively, unitary $G \times A_i$ representations with Haar $A_i$-Plancherel measure, and let $f: \mathcal{H}_1 \rightarrow \mathcal{H}_2$ be
a morphism which is $(G,T,A_2)$-equivariant up to a character of $A_2$ (i.e.\ $f \circ T(a)$ and $ a \circ f$ differ by a -- necessarily unitary -- character $A_2$).  What can we say about pull-backs of those relative norms? First of all,

\begin{lemma} \label{lemma1067}
Let $w \in \mathcal{H}_2$; then:
\begin{equation} \label{normpull}  f^* \| \cdot\|_{A_2, w} = \| \cdot\|_{A_2,  f' w}, 
\end{equation} 
an equality of norms on $\mathcal{H}_1$. Here $f'$ denotes the adjoint of $f$, and in defining the latter norm, we consider $\mathcal{H}_1$ as an $A_2$-representation via $T$. 
\end{lemma}
 
\begin{proof}
Decompose $\mathcal{H}_j = \int_{\chi \in \widehat{A}_2} \mathcal{H}_{j,\chi} d\chi$
and  $f = \int f_{\chi}$, with $f_{\chi}: \mathcal{H}_{1,\chi} \rightarrow \mathcal{H}_{2,\chi}$ (see \S \ref{ss:dis} for discussion). 
Then for $v = \int_{\chi} v_{\chi} d\chi \in \mathcal{H}_1 $ (and using similar notation for further decompositions of $v$, $w$, and $f$), 
\begin{eqnarray} \nonumber  \| f(v)\|^2_{A_2 , w} & 
=& \int_{\chi \in \widehat{A_2}} 
\|f_{\chi}(v_{\chi})\|^2_{L^1( w_{\chi})} d\chi \\ \nonumber 
& =& \int_{\chi \in \widehat{A_2}} \left| \int_{\hat G} |\left<f_{\pi,\chi} v_{\pi,\chi}, w_{\chi,\pi}\right>| \mu_\chi(\pi)\right|^2 d\chi 
\\  \nonumber  & =&   \int_{\chi \in \widehat{A_2}} \left| \int_{\hat G} |\left< v_{\chi,\pi}, f_{\pi,\chi}' w_{\pi,\chi}\right>| \mu_\chi(\pi)\right|^2 d\chi 
 \\ \label{abcdef}
 &=& \int_{\chi \in \widehat{A_2}} \|v_{\chi}\|^2_{L^1(f'w_{\chi})} d\chi = \|v\|^2_{A_2, f'w}.
  \end{eqnarray}
 Here $\mu_\chi$ denotes the disintegration of $G\times A_2$-Plancherel measure on $\mathcal H_1\oplus \mathcal H_2$ with respect to the forgetful map $\widehat{G\times A_2}\to \widehat{A_2}$.
\end{proof}

The important result will be that if $A_1$ and $A_2$ have different rank, the Hilbert space norm on $\mathcal H_1$ cannot be majorized (not even at the level of $J$-invariants) by any finite sum of pullbacks of mixed norms from $\mathcal H_2$. We keep assuming, of course, that the $A_i$-Plancherel measure of $\mathcal H_i$ is in the class of Haar measure.

\begin{lemma} \label{norm contra}
Let notation be as above.
Assume that $\dim(A_2 \otimes \mathbb{Q}) < \dim(A_1  \otimes \mathbb{Q})$. Then, for any open compact subgroup $J \subset G$,  the Hilbert norm on $\mathcal{H}_1^J$ is not majorized by any finite sum of norms of the form $f^* \|\bullet\|_{A_2, w}$. 

More generally, suppose given a finite collection of spaces $\mathcal{H}_2^{(j)}$, 
for $j=1, 2, \dots, $, together with tori $A_2^{(j)}$ and morphisms $T^{(j)}: A_2^{(j)} \rightarrow A_1$. Assume that $\dim(A_2^{(j)} \otimes \mathbb{Q}) < \dim(A_1  \otimes \mathbb{Q})$ for all $j$.
 Let $f_{(j)}: \mathcal H_1\to \mathcal{H}_2^{(j)}$ be morphisms as above.   Then the Hilbert norm on $\mathcal{H}_1^J$ is not majorized
by any finite sum of norms $f_{(j)}^* \| \bullet \|_{A_2^{(j)}, w}$. 
\end{lemma}

\begin{proof}
(With a single torus $A_2$:) When we decompose $\mathcal H^J_1$
over $Y:=\widehat{G \times A_2}$, each fiber $\mathcal{H}_{1, \pi}$
has infinite multiplicity as a $G$-representation.  The conclusion follows from Lemma \ref{lem:bounds} and Lemma \ref{lemma1067}.

(With multiple  tori $A_2^{(j)}$:)  We expand on the argument of the previous case. We need to show that
the Hilbert norm on $\mathcal{H}_1^J$ is not majorized by a finite sum  of norms of the type
$\sum \| \bullet\|_{A_2^{(j)}, f_{(j)}' w}$, where we follow the notation of \eqref{normpull}. Now our claim follows from Lemma \ref{Sumofnorms}. 
\end{proof}

Finally, a lemma on the $A_1$-invariance of the relative norms. Here we use as imput a property of disintegration of a $G\times A_1$-Plancherel measure with respect to the map $\widehat{G\times A_1}\to \widehat{G\times A_2}$, i.e.\ as opposed to the previous situation we forget a torus action. The assumption is an injectivity assumption, i.e.\ that for a given $G\times A_2$-representation there is a unique $G\times A_1$-representation appearing (stated measure-theoretically). Notice that in order to be able to disintegrate, we need a $\sigma$-finiteness assumption with respect to $G\times A_2$-Plancherel measure, which in our examples will be provided by Remark \ref{compactmeasure}.

\begin{lemma} \label{injectivemeasure}
Assume that for a disintegration $\mu = \int \mu_\alpha$ of Plancherel measure on $\mathcal H_1$ with respect to the forgetful map: $\widehat{G\times A_1}\to \widehat{G\times A_2}$, {\bluetext almost each of the measures $\mu_\alpha$ is concentrated on one point}. Then the pulled-back norms: $$f^* \|\bullet\|_{A_2, w}$$ are $A_1$-invariant. 
\end{lemma}

\begin{proof}
We have seen in \eqref{abcdef} that $$ f^* \|v\|^2_{A_2, w} =
\int_{\chi \in \widehat{A_2}} \|v_{\chi}\|^2_{L^1(f'w_{\chi})} d\chi,$$
where the norms $\|v_{\chi}\|^2_{L^1(f'w_{\chi})}$ are densely defined on the spaces $\mathcal H_{1,\chi}$. Clearly, $A_2$ acts trivially on these norms, so we need to show that they are $\mathrm{coker}(A_2 \rightarrow A_1)$-invariant.

By twisting the space $\mathcal{H}_{1, \chi}$ by a character of $A_1$ that extends $\chi$, we may suppose that $A_2$ acts trivially on $\mathcal{H}_{1, \chi}$. 
In this way we are reduced to the case where $A_2$ is trivial, and need to prove that under the same assumption for the Plancherel measure of $\mathcal H_1$ with respect to the map:
$$\widehat{G\times A_1}\to \hat G,$$
given a vector $w\in \mathcal H_1$, the norm:
$$\|\bullet\|_{L^1_w}$$
(defined with respect to a $\hat G$-Plancherel decomposition) is $A_1$-invariant.

Let $a\in A_1$. By the definition (\ref{L1norm}):
\begin{equation}\label{multbya}\|a\cdot v\|_{L^1(w)} = \int_{\hat G} \left|\left< a\cdot v_\pi, w_\pi\right>\right| \nu(\pi),
\end{equation}
where we have used a $G$-Plancherel decomposition $\mathcal H_1 = \int_{\hat G} \mathcal H_\pi \nu(\pi)$.

We may disintegrate the inner product $\left< a\cdot v_\pi, w_\pi\right>$ with respect to the $A_1$-Plancherel decomposition of $\mathcal H_\pi$ in such a way that the corresponding Plancherel measure $\mu_{\pi}(\chi)$ will be a disintegration of a given Plancherel measure $\mu(\pi,\chi)$ with respect to $\nu(\pi)$:
$$\left< a\cdot v_\pi, w_\pi\right> = \int_{\hat A_1} \left< a\cdot v_{\pi,\chi}, w_{\pi,\chi}\right> \mu_\pi(\chi).$$

But, by assumption, the measures $\mu_\pi$ are atomic (for almost all $\pi$), i.e.\ the last integral is equal to a multiple (independent of $a$) of $\left< a\cdot v_{\pi,\chi}, w_{\pi,\chi}\right>$ for some $\chi$ (depending on $\pi$). In particular, $|\left< a\cdot v_\pi, w_\pi\right>| = |\chi(a)\left< v_\pi, w_\pi\right>| = |\left<  v_\pi, w_\pi\right>|$, and therefore we get that (\ref{multbya}) is independent of $a$.

\end{proof}

\section[Preliminaries to scattering (II)]{Preliminaries to scattering (II): consequences of the conjecture on discrete series } \label{sec:linalg2} 
As the previous section, this section works out certain results needed in \S \ref{sec:scattering}.  The results here all depend on the validity
of the Discrete Series Conjecture \ref{dsconjecture} for $X$ and its degenerations, as in the statement of
Theorem \ref{advancedscattering}.

{ Recall that the canonical map:
$$ \mathcal Z(\LL_\Theta)^0\to \AA_{X,\Theta}$$
is surjective \emph{as a morphism of algebraic tori}. However, it may not be surjective at the level of $k$-points; we will thus be denoting by $A_{X,\Theta}'$ the image of:
$$ \mathcal Z(L_\Theta)^0\to A_{X,\Theta}.$$

The space of smooth functions on $X_\Theta$ varying by a character $\chi$ of $A_{X,\Theta}'$ will be denoted by $C^\infty(X_\Theta,\chi)$, and similarly for $L^2$-spaces etc.} 
We  will analyze when a representation
$\pi$ can occur simultaneously in $C^{\infty}(X_{\Theta}, \chi)$
and $C^{\infty}(X_{\Omega}, \psi)$.
What restrictions does this put on $\chi$ and $\psi$?
 In favorable situations 
they are related by an ``affine'' map  between the character groups of $A_{X, \Theta}'$ and $A_{X, \Omega}'$. In more detail:

\begin{itemize}

\item[-] \S \ref{ss:isogaffine} discusses the notion of an ``affine map''
between character groups of $k$-points of algebraic tori { (or finite-index subgroups thereof)}; this notion
will be used in understanding when a representation $\pi$ can occur
simultaneously in $C^{\infty}(X_{\Theta}, \chi)$ and $C^{\infty}(X_{\Omega}, \psi)$. 

\item[-] \S \ref{relationbetweencentralchars} applies this notion of ``affine map'' to the problem discussed: If $\pi$ embeds into $C^{\infty}(X_{\Theta})$ and $C^{\infty}(X_{\Omega})$, what is the relationship between central characters? 

\item[-] \S \ref{decompositionl2maps} uses the result of \S \ref{relationbetweencentralchars} to give a canonical decomposition of a morphism
$$L^2(X_{\Theta}) \rightarrow L^2(X_{\Omega})$$
into ``equivariant'' summands (equivariant for suitable actions of $A_{X, \Theta}'$). This is a critical {\em a priori} input into our analysis of scattering. 

\item[-] \S \ref{genericinjectivity}
shows that the decomposition of \S \ref{decompositionl2maps} can be further refined under the assumption that the map $\mathfrak{a}_X^*/W_X \rightarrow \mathfrak{a}/W$ is ``generically injective.''
\end{itemize}

\subsection{Isogenies of tori and affine maps on their character groups} \label{ss:isogaffine}

{ In this section we introduce the notion of an affine map between character groups of tori, which will be used for the canonical decomposition of morphisms in Proposition \ref{decompf}. We will eventually need to apply this notion to finite-index subgroups of tori, which will not be the points of an algebraic subgroup. In order to not make the notation too heavy, we present the definitions for algebraic tori only; to obtain the general definitions when a torus $\AA(k)$ is replaced by a finite index subgroup $A'$, the reader only has to replace:
\begin{itemize}
\item any occurence of $A=\AA(k)$ (or its character group) by $A'$ (resp.\ its character group);
\item the maximal compact subgroup $A_0\subset A$ by the maximal compact subgroup $A'\cap A_0\subset A'$;
\item finally, the notion of ``morphisms modulo isogenies'' $T: \AA_1\dashrightarrow \AA_2$ that we are about to introduce does not change.
\end{itemize}

}

 Let us remember that the category of algebraic tori is equivalent to the (opposite) category of finitely generated, torsion-free $\mathbb Z$-modules. The functor from the latter to finite dimensional $\QQ$-vector spaces (tensoring by $\QQ$ over $\Z$) corresponds to the semisimple category obtained from tori by inverting isogenies. We will be denoting a morphism in the latter by $T:\AA_1 \dashrightarrow \AA_2$; explicitly, such a morphism corresponds to an equivalence class of pairs of homomorphisms of tori: 
\begin{equation}\label{torimorphism}(\AA_1\to\DD, \AA_2\to\DD),\end{equation}
 with the second one finite and surjective, where ``equivalence'' is by passing simultaneously to a further finite quotient of $\DD$.

Each such morphism $T$ defines a canonical subgroup $A_1^T$ of $A_1$, as follows: recall that there is a canonical ``valuation'' map: $A_1\to \varchi(\AA_1)^*$. The map $T$ induces: 
\begin{equation}\label{cocharmaps}\varchi(\AA_1)^*\to \varchi(\AA_2)^*\otimes \QQ,
\end{equation}
 and the subgroup $A_1^T$ is defined as the preimage of those elements which map into $\varchi(\AA_2)^*$. In terms of a presentation (\ref{torimorphism}), this is equivalent to saying that the elements of $A_1^T$ are those whose images in $D$ have the same ``valuation'' as elements of $A_2$. 
 
 {\em Example.} If $\AA_1 = \AA_2 = \mathbb{G}_m$ and $T$ is the isogeny ``$x \mapsto x^{2/3}$,''
 described more formally as the diagram $\AA_1 \stackrel{2}{\rightarrow} \DD=\mathbb{G}_m \stackrel{3}{\leftarrow} \mathbb{G}_m$, then
 $$A_1^T = \{\lambda \in A_1  = k^{\times} | \mbox{ valuation of $\lambda$ is divisible by $3$} \}.$$

 As one can see from this example, $T$ does not induce a map $A_1^T \rightarrow A_2$,
 but it at least does induces a map:\begin{equation}\label{atmap}
 A_1^T/A_{1,0} \to A_2/A_{2,0},
\end{equation}
 where the index $~_0$ denotes maximal compact subgroup.

In particular, there is a canonical way to pull back any unramified character $\chi$ of $A_2$ to an unramified character $T^*\chi$ of $A_1^T$; the fact that $\chi$ is unramified will be implicit whenever we write such a pull-back.

A {\em component} of $\widehat{A_2}$ will be a connected component in the natural topology, i.e.,
the set of all characters with the same restriction to $A_{2,0}$. Every component is a coset for the component of the identity, i.e.\ the subgroup of unramified characters.

An {\em affine map} $\widehat{A_2}  \dashrightarrow \widehat{A_1}$  
{\em compatible with the morphism} $T: \AA_1 \dashrightarrow \AA_2$ is a mapping
$$ f: \mbox{some component of }\widehat{A_2} \rightarrow \widehat{A_1^T}$$
which is equivariant with respect to the natural homomorphism of unramified character groups induced by $T$:
$$\widehat{A_2}^0 \rightarrow \widehat{A_1^T}^0.$$
(In this equation, the superscript $0$ denotes connected component of the identity; in this case,
it coincides with the group of unramified characters, e.g.\ $\widehat{A}^0$ is the dual of $A/A_0$). The term ``affine'' is due to the analogy with affine maps between vector spaces (i.e.\ translates of linear maps).

In other words, for every unramified character $\chi_2$ of $A_2$ we have:
$$f(\chi_1 \chi_2) = f(\chi_1) T^* \chi_2, $$
when $\chi_1$ is in the component of definition of $f$.

{\em Explication.} For every affine map $f$ compatible with $T$  we may find
a character $\eta$ of $A_2$ and a character $\eta'$ of $A_1^T$ with the following property: 
the domain consists of all characters $\chi \in \widehat{A_2}$ for which
$\chi \eta^{-1}$ is unramified, and the map $f$ satisfies
\begin{equation} \label{fexplicit} f(\chi) =  \eta'  T^*(\chi \eta^{-1}). \end{equation} 

In the same way, we may define affine maps compatible with $T$ on the space of all (not necessarily unitary) complex characters, without the requirement that they preserve unitarity -- they will be denoted as:
$$f:  \widehat{A_2}_{\C} \dashrightarrow \widehat{A_1}_{\C}.$$
As before, we require $f$ to be defined only on one connected component of $\widehat{A_2}_{\C}$, have image in $\widehat{A_1^T}_{\C}$, and be equivariant with respect to the natural homomorphism of the identity components. 

We will use this generalization only once.

\begin{remark} \label{remarktori1} 
Note the following: If $T$ is defined by the pair of maps  $ \AA_1 \rightarrow \DD \leftarrow \AA_2$ as in \eqref{torimorphism}, we have mappings
$$ \widehat{A_2} \leftarrow \widehat{D}  \rightarrow \widehat{A_1} \rightarrow \widehat{A_1^T}$$
this gives us particular a multivalued function
$\widehat{A_2} \longrightarrow \widehat{A_1^T}$, where each
``value'' is a (possibly empty)  finite subset of $\widehat{A_1^T}$:
the image in $\widehat{A_1^T}$ of all preimages in $\widehat{D}$. 
 
 Then this (set-valued) morphism
is given, on each component of $\widehat{A_2}$, by
$$ \chi \in \widehat{A_2} \mapsto \{f_1(\chi), \dots, f_r(\chi)\}$$
where the $f_i$ are a (possibly empty) collection of affine maps compatible with $T$.

Indeed, suppose we begin with $\alpha \in \widehat{A_2}$
with nonempty image $\{ \psi_1, \dots, \psi_r \} \in \widehat{A_1^T}$.  Then it is easy to verify that the image of $\alpha \chi$, for 
 $\chi$ an unramified character, is given by 
 $$ \{ \psi_1 T^* \chi, \dots, \psi_r T^* \chi \} \in \widehat{A_1^T}.$$
 The rule $f_i : \alpha \chi \mapsto \psi_i T^* \chi$
 define an affine map from the component of $\widehat{A_2}$ containing $\alpha$
 to $\widehat{A_1^T}$, and this collection $\{f_1, \dots, f_r\}$ has the desired property. 
  \end{remark}

 \subsubsection{Maps of Hilbert spaces} \label{sss:fequiv}
Notation as above; in particular we have a morphism $T: \AA_1 \dashrightarrow \AA_2$
in the isogeny category and an affine map $f: \widehat{A_2} \dashrightarrow \widehat{A_1}$
on character groups that covers $T$. Choose also $\eta$ and $\eta'$ as in \ref{fexplicit}.

 Suppose we are given Hilbert spaces
$\mathcal{H}_1, \mathcal{H}_2$ with actions of $A_1, A_2$ respectively.  We  say that a mapping
$$ S: \mathcal{H}_2 \rightarrow \mathcal{H}_1$$
is $f$-equivariant if:
\begin{quotation}
  $S$ factors through $(A_{2,0}, \eta)$-coinvariants\footnote{Because $A_{2,0}$ is compact, with discrete dual, the canonical map from invariants to coinvariants is an isomorphism.} and produces $(A_{1,0}, \eta')$-invariants (where by the index $~_0$ we denote the maximal compact subgroups), 
and if  we define twisted actions of $A_2/A_{2,0}$ and $A_1/A_{1,0}$ on these coinvariant and invariant spaces via the rules:
$$a * v = \eta^{-1}(a) a\cdot v, \  a' *  v' = \eta'(a')^{-1} a'\cdot v'$$
then \begin{equation}\label{Siequivariance}
  S( T(a') * v)=  a'  * S(v), \ \ a' \in A_1^T/A_{1,0}.
\end{equation}
\end{quotation}

Equivalently, $S$ is $f$-equivariant if we may disintegrate with respect to the $A_2$, resp.\ $A_1^T$-action:
 $$\mathcal{H}_2 = \int_{\chi} \mathcal{H}_{2,\chi},$$
$$\mathcal{H}_1 = \int_{\chi'} \mathcal{H}_{1,\chi'},$$
and $S = \int_{\chi} S_{\chi}$ with $S_\chi=0$ unless 
$\chi \in \widehat{A_2}$ belongs to the domain of $f$, 
in which case $S_\chi$ is a morphism: $\mathcal{H}_{2,\chi} \rightarrow \mathcal{H}_{1,f(\chi)}$.

\subsection{Relationship between central characters for $X_{\Theta}$ and $X_{\Omega}$} \label{relationbetweencentralchars}

Recall that we are assuming the validity of the Discrete Series Conjecture \ref{dsconjecture} for all boundary degenerations $\XX_\Theta$. { Recall also that the image of the map:
$$ \mathcal Z(L_\Theta)^0\to A_{X,\Theta}$$
is denoted by $A_{X,\Theta}'$.}

\begin{proposition}  \label{e2c} 
 Let $\Theta,\Omega$ be two (possibly the same) subsets of $\Delta_X$.

 \begin{enumerate} 
\item  Let $J$ be an open compact subgroup of $G$.  There is a finite collection of morphisms 
\begin{equation} \label{Ti} T_i: \AA_{X,\Omega}\dashrightarrow \AA_{X,\Theta},
\end{equation}
and\footnote{The affine map $f_i$ determines the morphism $T_i$; nonetheless, we prefer to keep 
both in our notation.} affine maps 
\begin{equation}\label{fi}  f_i: \widehat{A'_{X, \Theta}}_\CC \dashrightarrow \widehat{A'_{X, \Omega}}_\CC 
\end{equation}
compatible with $T_i$ so that, for almost every $\chi\in \widehat{A_{X, \Theta}'}$ and every representation $\pi \hookrightarrow L^2(X_{\Theta},\chi)_{\disc}$ with non-zero $J$-fixed vector the following is true:

\emph{If $\pi$ embeds in $C^\infty(X_{\Omega},\psi)$ then\footnote{(More precisely, given that $f_i$ is defined only on a component of $A_{X, \Theta}$, we should say
that $(\psi, \chi)$ belongs to the graph of $f_i$. We shall allow ourself this type of imprecision at several points.) }
$\psi = f_i(\chi)$ for some $i$. }

\item \label{e2cisog} The subcollection of those morphisms (\ref{Ti}) which are \emph{isogenies}, and those affine maps (\ref{fi}) which preserve unitarity:
$$f_i: \widehat{A_{X, \Theta}'} \dashrightarrow \widehat{A_{X, \Omega}'}, $$
 is enough in order for the statement to be true for \emph{almost all} $\chi$ and all $\pi$ (with nonzero $J$-invariant vectors) which embed both into $L^2(X_\Theta,\chi)_{\disc}$ and $L^2(X_\Omega,\psi)_{\disc}$.
 
\end{enumerate} 
\end{proposition}

\subsubsection{Proof of Proposition \ref{e2c}}

The Discrete Series Conjecture \ref{dsconjecture} applied to the Levi varieties $X_\Theta^L$, $X_\Omega^L$, together with the finiteness of relative discrete series with $J$-fixed vectors (Theorem \ref{finiteds}) imply that there is a finite number of triples $(P^-,\sigma, D^*_\iR)$ where $P^-$ is a parabolic subgroup of $P_\Theta^-$, $\sigma$ a supercuspidal representation of its Levi quotient $L$ and $D^*_\iR$ a torus of unitary unramified characters of $P^-$ (with $D^*_\iR\to\widehat{A_{X,\Theta}'}$ finite) such that, for almost every $\chi\in \widehat{A_{X, \Theta}'}$, any representation $\pi\in L^2(X_\Theta,\chi)_\disc$ which admits a $J$-invariant vector is a subquotient of $\pi'=I_{P^-}^G(\sigma\otimes\omega)$ for such a triple and some $\omega\in D^*_\iR$. For each such triple we have morphisms:
$$\mathcal Z(\LL_\Theta)^0 \hookrightarrow \mathcal Z(\LL)^0 \to \DD$$
(here $\DD$ denotes the torus quotient of $\PP^-$ defining this toric family of relative discrete series), whose composition is finite and surjective, and such that $\chi$ is the twist of $\eta:=$the central character of $\sigma$ (restricted to $\mathcal Z(\LL_\Theta)^0$) by the pull-back of an unramified character of $\DD$. Notice that, after possibly replacing $\DD$ by a finite quotient, the map $\mathcal Z(\LL_\Theta)^0\to \DD$ factors through a map:
\begin{equation}\label{map1}
\AA_{X,\Theta} \to \DD,
\end{equation} 
since unramified characters of $D$ have to be trivial on the kernel of $\mathcal Z(L_\Theta)^0\to \mathcal Z(X_\Theta)$.
($D$ is a torus quotient of $\mathbf{L}$, by definition, and the corresponding  representations parabolically induced from $\mathbf{P}$ 
embed in functions on $X_{\Theta}$; in particular, their central character   must factor through $\mathcal Z(L_\Theta)^0\to A_{X,\Theta}'$, proving the claim.)

Now, if $\pi$ embeds in $C^{\infty}(X_{\Omega},\psi)$ it must be a submodule of a representation parabolically induced from $P_{\Omega}^-$, because
of the description of the variety $X_{\Omega}$ itself as being parabolically induced (Lemma \ref{lemmalevi}). 
In order for $I_{P^-}^G(\sigma\otimes\omega)$ to have a common subquotient with a representation induced from $P_\Omega^-$, equivalently with a supercuspidal induced from a parabolic subgroup of $P_\Omega^-$, that supercuspidal should be a $w$-twist of $\sigma\otimes\omega$, for some element $w\in W$, the Weyl group of $G$.  Each $w\in W$ such that $wL\subset L_\Omega$ defines a morphism: 
\begin{equation}\label{map2}\mathcal Z(\LL_\Omega)^0\to \mathcal Z(\LL)^0 \to \DD.\end{equation}
 Let $T:\mathcal Z(\LL_\Omega)^0\dashrightarrow \mathcal Z(\XX_\Theta)= \AA_{X,\Theta}$ be the morphism defined by the equivalence class of the pair of maps (\ref{map1}),(\ref{map2}). The possible $\mathcal Z(L_\Omega)^T$-characters by which $\pi$ can be embedded into $C^\infty(X_\Omega)$ are thus ``images'' of $\chi \in 
  \widehat{ A_{X,\Theta}'}$ under
 the multivalued mapping arising from the diagram:
 $$  \widehat{ A_{X,\Theta}'}_{\C} \leftarrow \widehat{D}_{\C} \rightarrow \widehat{\mathcal{Z}(L_{\Omega})^0}_{\C}$$
 which, as we discussed in Remark \ref{remarktori1}, can be expressed
 as an affine mapping on each component of the character groups:
\begin{equation}\label{chgroupsmap} \widehat{A'_{X, \Theta}}_\CC \dashrightarrow \widehat{\mathcal{Z}(L_{\Omega})^0}_{\C}
\end{equation}

Now we verify that $T$ factors through the quotient: $\mathcal Z(\LL_\Omega)^0\twoheadrightarrow \AA_{X,\Omega}$: If we are given an affine mapping (\ref{chgroupsmap})
compatible with $T: \mathcal{Z}(\LL_{\Omega})^0 \dashrightarrow \AA_{X, \Theta}$, 
with the property that the image lies within
$\widehat{A_{X, \Omega}'}_\C$ (considered naturally as a subset of $ \widehat{\mathcal{Z}(L_{\Omega})^0}_{\C}$), 
then  in fact $T$ must factor\footnote{Here is the argument in an abstract context: 
Suppose given an isogeny $T : \AA_1 \dashrightarrow \AA_2$
and a corresponding affine map $f: \widehat{A_1} \rightarrow \widehat{A_2^T}$. (Our affine maps above are with respect to finite-index subgroups, but this doesn't make a difference for the argument.) 
We suppose that there is a quotient $\AA_1\twoheadrightarrow \BB$
such that the image of $f$ is contained in pullbacks of elements of $\BB$.
We claim, then, that $T$ factors through $\BB$ in the isogeny category.
Let $\mathbf{K}$ be the kernel of $\AA_1 \rightarrow \BB$. The assumption forces the pullback of any unramified character of $\AA$ to be trivial on $\mathbf{K}$;
equivalently, the image  of $\mathbf{K}(k) \cap A_1^T$ under 
$$ \mathbf{K}(k) \cap A_1^T \rightarrow A_1^T/A_{1,0} \rightarrow X_*(\AA_1)  \hookrightarrow
  X_*(\AA_2) \otimes \Q$$
is trivial.  In other words, the map  
$\varchi(A_1)^* \rightarrow \varchi(A_2)^* \otimes \Q$ induced by $T$ is trivial
on a finite index subgroup of the subgroup $\varchi(K)^*$ of $\varchi(A_1)^*$
corresponding to $\mathbf{K} \subset \AA_2$. This means that it is 
in fact trivial on all of $\varchi(K)^*$, and so
$T$ indeed factors in the isogeny category through the quotient $\AA_1/\mathbf{K}$.}
 through the mapping $\mathcal{Z}(L_{\Omega})^0 \rightarrow A_{X, \Omega}$, that is to say, $T$ determines a mapping 
$T: \AA_{X, \Omega} \dashrightarrow \AA_{X, \Theta}$
as desired.

For the second assertion,  if $\pi$ belongs to both $L^2(X_\Theta,\chi)_\disc$ and $L^2(X_\Omega,\psi)_\disc$ then applying the
 ``toric discrete series'' assumption to both we get triples $(P^-,\sigma, D^*_\iR)$  and $(Q^-,\sigma', D'^*_\iR)$ (with $P^-\subset P_\Theta^-$ and $Q^-\subset P_\Omega^-$) such that $\pi$ is a subquotient of the corresponding 
 induced representations.   It is known that, should $I_{P^-}(\sigma \cdot \chi)$
have a common subquotient as $I_{Q^-}(\sigma' \cdot \chi')$, this implies that there exists $w \in W$
carrying the Levi subgroup of $P^-$ to the Levi subgroup of $Q^-$
which carries $\sigma \cdot \chi$ to $\sigma' \cdot \chi'$.  Our assumption is that such $w$ exists for a set of positive measure; in particular, we may suppose that a {\em particular} $w$
works for a set of $\chi$ of positive measure. Twisting $\sigma, \sigma'$, we may suppose
that $w \sigma = \sigma'$ and that the following is true:
For a positive measure set $Z$ of unramified characters $\chi \in D_\iR^*$, the character $w \chi$ of $w P^{-}$ factors through the torus quotient corresponding to $D_\iR'^*$.

Now, given a set of positive measure (and thus Zariski-dense) of unramified unitary characters in $D^*_\iR$,
the intersection of their kernels is necessarily simply the maximal compact subgroup $D_0$ of the torus quotient $P^-\to D$ corresponding to $D_\iR^*$.  
Thus the intersection of the kernels of all $\chi \in Z$ is simply the preimage of $D_0$ in $P^-$.
Similarly, the intersection of the kernels of all $w \chi \ (\chi \in Z)$ is the preimage
of $D_0'$ in $Q^{-}$. 
So the map $w$ must then carry the preimage of $D_0$ in $P^-$ into the preimage of $D'_0$ in $Q^{-}$.  
The map $w$ must then carry this into the preimage of $D'_0$ in $Q^{-}$. In particular $w$ induces a mapping $D_\iR^* \rightarrow D'^*_\iR$ and this mapping
has the property that $\mathrm{ind}(\sigma \cdot \chi)$ and $\mathrm{ind}(\sigma' \cdot w(\chi))$
have a common subconstituent for all $\chi$. Since (as part of the assumption of relative discrete series) we suppose that the maps $D^*_{\iR} \rightarrow \widehat{A_{X,\Theta}'}$
and $D_{\iR}'^* \rightarrow \widehat{A_{X,\Omega}'}$ are finite and surjective, 
$w$ induces an isogeny $\mathcal{Z}(X_{\Theta}) \dashrightarrow
 \mathcal{Z}(X_{\Omega})$. The set of all such isogenies  $\mathcal{Z}(X_{\Theta})\dashrightarrow
 \mathcal{Z}(X_{\Omega})$ that arise in this fashion from some $w \in W$ then has the property stipulated by the proposition.

\qed
 
 \subsection{Canonical decomposition of maps $L^2(X_{\Theta}) \rightarrow L^2(X_{\Omega})$ } \label{decompositionl2maps}

\begin{proposition} \label{decompf} 
\begin{enumerate} 
 \item 
Suppose that $S: L^2(X_{\Theta})_{\disc}
\rightarrow L^2(X_{\Omega})$ is a $G$-equivariant morphism. Then there exists
a   unique (up to indexing) decomposition: 
\begin{equation}\label{seriesS} S=  \sum_{i=1}^\infty S_i
\end{equation}
such that each $S_i$ is a nonzero bounded morphism, and equivariant with respect to some (distinct) affine map between central character groups.

In other words, for each $i$ there is a pair $(T_i, f_i)$ (with $f_i\ne f_j$ when $i\ne j$), where
$$ T_i: \AA_{X,\Omega}\dashrightarrow \AA_{X,\Theta}$$
is a morphism in the isogeny category of tori and 
$$f_i: \widehat{A_{X, \Theta}'}\dashrightarrow \widehat{A_{X, \Omega}'}$$
is an affine map compatible with $T_i$, such that $S_i$ is 
$f_i$-equivariant (see \ref{sss:fequiv}). For each open compact subgroup $J$ only a finite number of summands in \eqref{seriesS} are non-zero on $L^2(X_\Theta)^J_\disc$.

\item  \label{orthogonality} If $|\Theta|\ne|\Omega|$ then $L^2(X)_\Theta\perp L^2(X)_\Omega$.

\end{enumerate}
\end{proposition}

\begin{proof}

 Given $S: \mathcal H_2:= L^2(X_\Theta)_\disc\to \mathcal H_1:=L^2(X_\Omega)$, let $\mu$ be a Plancherel measure for $\mathcal H_2\oplus\mathcal H_1$, and let $S=\int S_\pi$ be the corresponding decomposition of $S$, as in \S \ref{ss:dis}, where $\int_{\hat G} \mathcal H_{2,\pi} \mu(\pi)$ and $\int_{\hat G} \mathcal H_{1,\pi} \mu(\pi)$ are direct integral decompositions for $\mathcal H_2$ and $\mathcal H_1$, respectively.

Fix some $T:\AA_{X,\Omega}\dashrightarrow \AA_{X,\Theta}$, and unitary characters $\eta$ of $A_{X,\Theta}'$ and $\eta'$ of $A_{X,\Omega}'^T$. Let $\mathcal H_2^\eta, \mathcal H_1^{\eta'}$ be the eigenspaces where the maximal compact subgroups of $A_{X,\Theta}'$ and $A_{X,\Omega}'$ act via the characters $\eta$ and $\eta'$, respectively. The valuation maps give surjections:
\begin{equation}\label{mapp1} A_{X,\Omega}'^T \twoheadrightarrow \Gamma:=\varchi(\AA_{X,\Omega})^*\cap T^{-1}\varchi(\AA_{X,\Theta})^*
\end{equation}
and:
\begin{equation}\label{mapp2} A_{X,\Theta}' \twoheadrightarrow \varchi(\AA_{X,\Theta})^*;
\end{equation}
we may let the lattice $\Gamma$ act on $\mathcal J^{\eta,\eta'}:=\mathcal H_2^\eta \oplus \mathcal H_1^{\eta'}$ as:
\begin{equation}\label{gammaaction}\gamma\cdot (h + h') = \eta^{-1}(b) b\cdot h + \eta'^{-1}(a) a\cdot h',
\end{equation}
where $a$ is any lift of $\gamma$ via (\ref{mapp1}) and $b$ is any lift of $T(\gamma')$ via (\ref{mapp2}). Hence we get a homomorphism:
\begin{equation}\label{gammatoend}
 \Gamma\to \Aut_G( \mathcal J^{\eta,\eta'}).
\end{equation}

Let $\mathcal J^{\eta,\eta'}=\int_{\hat G} \mathcal J_\pi^{\eta,\eta'} \mu(\pi)$ be a Plancherel decomposition; we claim that (\ref{gammatoend}) decomposes into the direct integral of the analogous maps: 
\begin{equation}\label{gammatoendlocal}
 \Gamma\to \Aut_G( \mathcal J_\pi^{\eta,\eta'}),
\end{equation}
defined via the same formula (\ref{gammaaction}). Indeed, if we decompose
(by \S \ref{ss:dis})
 the action of $\gamma\in \Gamma$
as an integral of $\gamma_\pi: J_\pi^{\eta,\eta'}\to J_\pi^{\eta,\eta'}$ then it is clear that for almost\footnote{Recall that every $\pi$ can only appear with a finite number of $A_{X,\Theta}$-exponents (cf.\ the proof of Theorem \ref{finiteness}) in $C^\infty(X_\Theta)$. By the presumed validity of the Discrete Series Conjecture \ref{dsconjecture} for $X_\Theta$, statements that hold ``for almost all $\pi$'' in the spectrum of $L^2(X_\Theta)_\disc$ also hold ``for almost all $\chi$''.} all $\pi$ the endomorphism $\gamma_\pi$ of $\mathcal J_\pi^{\eta,\eta'} = \mathcal H_{2,\pi}^\eta\oplus\mathcal H^{\eta'}_{1,\pi}$ must coincide with the action of the element $\gamma$ defined as in (\ref{gammaaction}).  

The inner action of $\Aut_G( \mathcal J_\pi^{\eta,\eta'})$ on $\End_G( \mathcal J_\pi^{\eta,\eta'})$ defines by (\ref{gammatoendlocal}) an action of $\Gamma$ on the latter, and it is easy to see by its definition that this action preserves the subspace $\Hom_G(\mathcal H_{2,\pi}^{\eta},\mathcal H_{1,\pi}^{\eta'})$. In our setting, these spaces are finite-dimensional, 
and they carry a natural inner product that is preserved by $\Gamma$.  Let $S_\pi^{\Gamma}$ be the projection of $S_\pi$ to the eigenspace for the trivial character of $\Gamma$. 
We will show in a moment that  the function: $\pi\mapsto \Vert S_\pi^{\Gamma}\Vert$ is essentially bounded. Assuming that for a moment, by the Proposition
\ref{measurable2} on ``measurability of eigenprojections'' \footnote{By \S \ref{trivialityremark}, by decomposing $\hat G$ into a countable union of measurable sets, we may identify the measurable structure of the family of vector spaces $\End_G(\mathcal J_\pi^{\eta,\eta'})$ with that of a trivial family; hence, the Proposition applies. }  we may integrate  the morphisms $S_\pi^{\Gamma}$ in order to get a morphism:
 
\begin{equation}
 S^{\Gamma}:= \int_{\hat G} S_\pi^{\Gamma}: \mathcal H_2^\eta \dashrightarrow \mathcal H_1^{\eta'},
\end{equation}
where the dotted arrow means that it is well-defined on a dense subspace of $\mathcal H^\eta$. By construction, $S^{\Gamma}$ is equivariant with respect to any affine map: $f_i: \widehat{A_{X, \Theta}'}_\CC \dashrightarrow \widehat{A_{X, \Omega}'}_\CC$ which covers $T$ and maps $\eta$ to $\eta'$; indeed, it is clear by construction that \eqref{Siequivariance} is satisfied. 

Now, we verify that
 the function: $\pi\mapsto \Vert S_\pi^{\Gamma}\Vert$ is essentially bounded.To that end, decompose each $\mathcal J^{\eta,\eta'}$ into  $A_{X,\Omega}'^T \times A_{X, \Theta}'$-eigenspaces (these are genuine eigenspaces, rather than generalized ones, because
 this action preserves the natural inner product). From the Plancherel decomposition of $\mathcal J^{\eta,\eta'}$ \emph{as an $A_{X,\Omega}'^T\times G$-representation} it follows that the distinct generalized eigenspaces are, for almost all $\pi$, honest eigenspaces and orthogonal to each other. Since $S_\pi^{\Gamma}$ is the sum of some of the operators:
$$\pr_1\circ S_\pi \circ \pr_2$$
where $\pr_1$ and $\pr_2$ vary though all projections to $A_{X, \Omega}'^T \times A_{X,\Theta}'$- eigenspaces, the norm of $S_\pi^{\Gamma}$ can be bounded by the norm of $S_{\pi}$, multiplied by a number depending only on the number of distinct eigenspaces. Finally, we recall that this number is uniformly bounded by the order of the Weyl group (cf.\ Lemma \ref{numberofexponents}).

Now let $(T_i,f_i)$ vary over all those pairs of the first part of Proposition \ref{e2c}. Notice that they are finitely many if we restrict to representations with non-zero $J$-invariant vectors, so all together they will be at most countably many. Let $\Gamma_i$ denote the corresponding finitely generated abelian groups, defined as above. For each $i$ we get an operator $S_i= \int_{\hat G} S_{\pi,i}$ as above. We claim:

\begin{equation}\label{Ssum}
 S = \sum_i S_i.
\end{equation}

Indeed, for almost all $\pi$, by the second part of Proposition \ref{e2c} each eigenspace for the $A_{X,\Theta}' \times A_{X,\Omega}'$-action on $\Hom_G(\mathcal H_{2,\pi},\mathcal H_{1,\pi})$ is contained in the fixed (eigenvalue=$1$) subspace of $\Gamma_i$, for some $i$. Moreover, it is clear that for different $i$'s and almost all $\pi$ (in the setting of Proposition \ref{e2c}: almost all $\chi$) the fixed subspaces of $\Gamma_i$ and $\Gamma_j$ for $i\ne j$ are distinct.
If it is not so,  there is a positive measure set of $\chi \in \widehat{A_{X, \Theta}'}$
such that $f_i(\chi) = f_j(\chi) \in \widehat{A_{X, \Omega}'}$. This means
that $f_i, f_j$ {\em coincide}, a contradiction. This implies (\ref{Ssum}).

Regarding uniqueness: For any decomposition $S = \sum_i S_i$, 
where each $S_i$ is $(T_i, f_i)$-equivariant, 
we must have $S_{\pi} = \sum_{i} S_{\pi, i}$, where $S_{\pi, i}$ is obtained
by disintegrating the homomorphism $S_i$, in the same sense as we have seen above. 
But then (for almost all $\chi$) 
$S_{\pi, i}$  is {\em characterized} as the $\Gamma_i$-fixed part of $S$ (where $\Gamma_i$ and its action is defined as before), since
the measure of the set of $\chi$ where  $f_i(\chi) = f_j(\chi)$ for any $i \neq j$ is zero.

To prove the second statement, assume that $|\Theta|>|\Omega|$, hence $\dim \AA_{X,\Theta}<\dim \AA_{X,\Omega}$.

 From Proposition \ref{e2c} it follows that there is a subset $Z$ of $\widehat{A_{X,\Omega}'}$ of measure zero,
  and a subset $Z'$ of $\widehat{A_{X,\Theta}'}$ of measure zero,  such that if $\pi\in L^2(X_\Theta)_\disc$  does not have central character in $Z'$ and admits a non-zero morphism into $C^\infty(X_\Omega,\psi)$ for some $\psi\in \widehat{A_{X,\Omega}'}$ then $\psi$ belongs to $Z$.
(After all any morphism: $\AA_{X,\Omega}\dashrightarrow\AA_{X,\Theta}$ has positive-dimensional kernel.) Now, if $L^2(X)_\Omega$ is not orthogonal to $L^2(X)_\Theta$ then we get a non-zero morphism: $\iota_\Omega^*\iota_\Theta: L^2(X_\Theta)_\disc\to L^2(X_\Omega)_{\disc}$.   But this is impossible:  because we are assuming
the Discrete Series Conjecture for $X_{\Omega}$, the $L^2(X_{\Omega})_{\disc}$-Plancherel measure of representations with central character in $Z$ is also zero,
and similarly for $\Theta, Z'$. 
   It follows that $L^2(X)_\Omega\perp L^2(X)_\Theta$ if $|\Omega|\ne|\Theta|$.

\end{proof}

\section{Scattering theory} \label{sec:scattering}

\subsection{} In \S \ref{sec:Bernstein} we constructed canonical maps $\iota_\Theta: L^2(X_\Theta ) \rightarrow L^2(X)$, and we saw (Corollary \ref{corollarysurjection}) that $\sum_\Theta \iota_\Theta$ induces a surjection:
$$\bigoplus_{\Theta\subset\Delta_X} L^2(X_{\Theta})_{\disc} \twoheadrightarrow L^2(X).$$
We will be denoting the image of $L^2(X_\Theta)_\disc$ by $L^2(X)_\Theta$.

The question of \emph{scattering} is the description of the kernel of this map.
Our answer is described in Theorem \ref{advancedscattering} (which is conditional on 
Conjecture \ref{dsconjecture} and on generic injectivity, but we expect it to hold in full generality), which implies in particular that $L^2(X)_\Theta$ and $L^2(X)_\Omega$ coincide if $\Theta$ and $\Omega$ are $W_X$-associates (i.e.\ there is a $w\in W_X$ such that $w\Theta = \Omega$), and are orthogonal otherwise. 

The starting point for the proof of Theorem \ref{advancedscattering} is the relatively straightforward statement of Proposition \ref{decompf}:
\begin{quote}
Any morphism $L^2(X_{\Theta})_\disc \rightarrow L^2(X_{\Omega})_\disc$ decomposes uniquely as a sum:
\begin{equation}\label{SsumSi}\sum S_i, 
\end{equation}
where the morphism $S_i$, assumed non-zero, is equivariant with respect to an isogeny $T_i: \AA_{X,\Omega}\dashrightarrow \AA_{X, \Theta}$ and an affine map $f_i$ of character groups compatible with $T_i$, { and the pairs $(T_i,f_i)$ are assumed to be distinct}. \end{quote} 
Recall that this statement followed, essentially, from the assumption of validity of the Discrete Series Conjecture \ref{dsconjecture} and general facts about induced representations. 

For $\Theta,\Omega\subset\Delta_X$ and $w\in W_X(\Omega, \Theta)$ (the set of elements of $W_X$ taking $\Theta$ to $\Omega$), we say that a pair $(T_i,f_i)$ as above \emph{corresponds} to $w$ if $T_i$ is the isomorphism: 
\begin{equation}\label{gray}
\AA_{X,\Omega}\to \AA_{X,\Theta} 
\end{equation}
induced by $w^{-1}$ 
(and hence $f_i$ is the restriction to some component of $\widehat{A_{X,\Theta}'}$ of the map of character groups obtained by $T_i$). For a morphism $S:L^2(X_\Theta)\to L^2(X_\Omega)$ we will call \emph{the $w$-part of $S$} the sum of those summands in its decomposition (\ref{SsumSi}) for which $(T_i,f_i)$ is induced by (\ref{gray}). Notice that the $w$-part of such a morphism is $A_{X,\Theta}'$-equivariant when $A_{X,\Theta}'$ acts on $L^2(X_\Omega)$ via $w:A_{X,\Theta}'\to A_{X,\Omega}'$.

Applying this decomposition to the\label{discisok}\footnote{A minor remark:  The morphism $\iota_{\Omega}^* \iota_{\Theta}$ in fact
maps $L^2(X_{\Theta})_{\disc}$ into $L^2(X_{\Omega})_{\disc}$;
in particular, $\iota_{\Omega, \disc}^* \iota_{\Theta, \disc} $
and $\iota_{\Omega}^* \iota_{\Theta}$ coincide on $L^2(X_{\Theta})_{\disc}$. 
Indeed -- see the proof of Corollary 
\ref{corollarysurjection} -- we need only verify that, for $f \in L^2(X_{\Theta})_{\disc}$
that $\iota_{\Omega}^* \iota_{\Theta} f$ is perpendicular to all 
$\iota^{\Omega}_{\Omega'} f'$, where $f' \in L^2(X_{\Omega'})$, and $\Omega'$
contains $\Omega$. But:
$$\langle \iota_{\Omega}^* \iota_{\Theta} f, \iota_{\Omega'}^{\Omega} f' \rangle = 
\langle \iota_{\Theta} f, \iota_{\Omega} \iota_{\Omega'}^{\Omega} f'  \rangle 
= \langle \iota_{\Theta} f, \iota_{\Omega'} f' \rangle,$$
and we may now apply the fact that $L^2(X)_{\Theta}$
and $L^2(X)_{\Omega'}$ are perpendicular if $|\Theta| \neq |\Omega'|$, 
in view of Proposition \ref{decompf}. 
}\emph{scattering morphisms} $i_\Omega^* i_\Theta$, we additionally need to establish the following:
\begin{itemize}
 \item the only $(T_i,f_i)$'s that appear are those corresponding to elements of $W_X(\Omega,\Theta)$;
 \item for all $w\in W_X(\Omega,\Theta)$ the $w$-part of $i_\Omega^* i_\Theta$ is an isometry. 
\end{itemize}
As we will see in \S \ref{proof:advancedscattering}, these two facts imply Theorem \ref{advancedscattering}.

To establish these two properties we need some algebraic input, encoded the condition of ``generic multiplicity one'' of Theorem \ref{advancedscattering}, together with some hard analysis. The analysis leads to Theorem \ref{tilingtheorem}, which should be regarded as the main result of this section; let us first discuss the algebraic condition.

\subsection{Generic injectivity of the map: $\mathfrak a_X^*/W_X\to \mathfrak a^*/W$}

\label{genericinjectivity}  

Denote: $\mathfrak a_{X}^* := \varchi(\XX)\otimes\QQ\subset \mathfrak a^*:=\varchi(\BB)\otimes\QQ$. Recall that the subgroup $W_{L(X)}$ of $W$ is the pointwise stabilizer of $\mathfrak a_X^*$, the subgroup $W_X$ normalizes $\mathfrak a_X^*$, and its action on it is generated by simple reflections. Fix a Weyl chamber for $W_X$ on $\mathfrak a_X^*$; the image of a face $\mathcal F$ of that Weyl chamber\footnote{Recall that a ``face'' is the intersection of the Weyl chamber with the kernel of a linear functional which is non-negative on it; hence the whole chamber is also a face.} in $\mathfrak a_X^*/W_X$ will be called a ``face'' of $\mathfrak a_X^*/W_X$. 

The condition called ``generic injectivity of the map: $\mathfrak a_X^*/W_X\to \mathfrak a^*/W$ on each face'' in the statement of Theorem \ref{advancedscattering} is the following: 
\begin{quote} 
{For every integer $d$, the restriction of the map: $\mathfrak a_X^*/W_X\to \mathfrak a^*/W$  to the 
collection of $d$-dimensional faces of $\mathfrak a_X^*/W_X$ is generically injective.}
\end{quote}
By ``generically injective'' we mean injective on a subset of full measure, for the natural class of measures, but this is easily seen to be equivalent, in this case, to injectivity outside of the image of a finite number of hyperplanes in $\mathfrak a_X^*$. In other words, outside of a meager set two distinct elements of $\mathcal F_d$  (the union of all $d$-dimensional faces) cannot be $W$-conjugate.

An equivalent way to formulate this condition, in terms of dual groups, is the following: 
\begin{quote}
{For every pair of (standard) Levi subgroups $\check L_\Theta, \check L_{\Omega}$ of $\check G_X$, and any isomorphism of their centers $\mathcal Z(\check{L}_{\Theta}) \rightarrow \mathcal Z(\check{L}_{\Omega})$  induced by an element of the Weyl group $W$ of $\check G$, there is an element of the little Weyl group $W_X$ which induces the same  isomorphism.}
\end{quote}

It is evident that this condition is very easy to check in each particular case.
It is always true in both extreme cases: When the dual group $\check{G}_X$ is isomorphic to $\SL_2$, and when it is all of $\check{G}$.  Less trivially, { Delorme has shown that it holds for all symmetric varieties:

\begin{proposition} \label{delorme-recent-label}\cite[Lemma 15]{Delorme-Plancherel}  If $\XX$ is symmetric, then it satisfies the generic injectivity condition.
\end{proposition}}

\begin{example} 
 For $\XX=\SSp_{2n}\backslash \GGL_{2n}$ the dual group is $\check G_X = \GL_n \hookrightarrow \GL_{2n} = \check G$, with a diagonal element $\diag(\chi_1,\dots,\chi_n)$ embedded as:
$$ \diag(\chi_1,\chi_1,\chi_2,\chi_2,\dots, \chi_n,\chi_n).$$

Let $\Theta,\Omega$ be two subsets of the simple roots of $\check G_X$. Their kernels, are subtori of $A_X^*$ which, when considered as subtori of $\check G$ are the connected centers of standard Levi subgroups $\check L_{\tilde\Theta}, \check L_{\tilde\Omega}$ corresponding to subsets $\tilde\Theta, \tilde\Omega$ of the simple roots of $\check G$. (Explicitly: for the usual numbering $1,\dots,$ of the roots of $\GL_n$ and $\GL_{2n}$, $\tilde\Theta$ is the union of $2\cdot \Theta$ and all odd roots, and similarly with $\tilde\Omega$.) Any isomorphism between $\mathcal Z(\check L_{\tilde\Theta})^0, \mathcal Z(\check L_{\tilde\Omega})^0$ induced by an element of $W$ is actually induced by an element of $W(\tilde\Omega,\tilde\Theta)$. But $W(\tilde\Omega,\tilde\Theta) = W_X(\Omega,\Theta)$, therefore $\XX$ satisfies the injectivity assumption.
\end{example}

\begin{lemma} \label{injectivefaces}    Assume that the map: $$\mathfrak a_X^*/W_X\to \mathfrak a^*/W$$ is generically injective on each face. 
  
  Then  -- notation as in (\ref{SsumSi}) -- the only pairs $(T_i, f_i)$ that can appear in the decomposition of a morphism $S: L^2(X_{\Theta})_{\disc} \longrightarrow L^2(X_{\Omega})_{\disc}$ 
are  those corresponding to elements of $W_X(\Omega,\Theta)$.

Hence, any such morphism decomposes as a sum of its $w$-parts:
\begin{equation} \label{onlyw}
\sum_{w\in W_X(\Omega,\Theta)} S_w,
\end{equation}
where $S_w$ is the sum of those $S_i$'s in \eqref{SsumSi} for which the corresonding affine map $f_i$ of character groups is induced (by restriction to a connected component of $\widehat{A_{X,\Theta}'}$) by the element $w$.
\end{lemma}

\begin{proof} 
Revisiting the proof of the last assertion of Proposition \ref{e2c}, let us give ourselves two toric families of relative discrete series $(P,\sigma, D^*_\iR)$ and $(Q,\sigma',D'^*_\iR)$, for $X_\Theta$ and $X_\Omega$ respectively, such that $D^*_\iR$ is the group of unitary elements in a torus $D^*$ of unramified characters identified with $\mathcal Z(\check L_{X,\Theta})$ under (\ref{chargroups}) and $D'^*_\iR$ is the group of unitary elements of a torus $D'^*$ identified with $\mathcal Z(\check L_{X,\Omega})$ under (\ref{chargroups}). Recall that the Lie algebras of $\mathcal Z(\check L_{X,\Theta})$, $\mathcal Z(\check L_{X,\Omega})$ are the complexifications of the vector spaces $\mathfrak a_{X,\Theta}^*$, $\mathfrak a_{X,\Omega}^*$, respectively. As we saw in the proof of Proposition \ref{e2c}, any summand of the scattering map should arise from an element $w\in W$ which takes $\mathfrak a_{X,\Theta}^*$ to $\mathfrak a_{X,\Omega}^*$.

If the map: $\mathfrak a_X^*/W_X \to \mathfrak a^*/W$ is generically injective on every face, this means that any element $w\in W$ which carries one family into to the other induces the same map: $\mathfrak a_{X,\Theta}^* \to \mathfrak a_{X,\Omega}^*$ as an element of $W_X(\Omega,\Theta)$.  

This proves the lemma.
\end{proof}

\subsection{The scattering theorem}

Define $L^2(X)_{i}$ to be the image of $$\bigoplus_{|\Theta| =i} L^2(X_{\Theta})_{\disc}$$
in $L^2(X)$.   Part \ref{orthogonality} of Proposition \ref{decompf} implies that we have a direct sum decomposition:
\begin{equation} \label{directbydim}
 L^2(X)=\bigoplus_i L^2(X)_i.
\end{equation}

Denote by $\mathfrak a_{X,\Theta}$ the space $\varchi(\AA_{X,\Theta})^*\otimes \QQ\subset\mathfrak a_X$ and by $\mathfrak a_{X,\Theta}^+$ its ``anti-dominant chamber'', i.e.\ its intersection with the cone $\mathcal V$ of invariant valuations. We denote by $\mathring{\mathfrak a}_{X,\Theta}^+$ the interior of $\mathfrak a_{X,\Theta}^+$.  If a morphism $T_i: \AA_{X,\Omega}\dashrightarrow \AA_{X,\Theta}$ is an isogeny as in part \ref{e2cisog} of Proposition \ref{e2c}, it induces an isomorphism (again to be denoted by $T_i$): $\mathfrak a_{X,\Omega}\xrightarrow{\sim} \mathfrak a_{X,\Theta}$.

Let $\Theta\subset\Delta_X$ and let $\Omega$ range over the subsets of $\Delta_X$ of the same size as $\Theta$ (including $\Theta$). Let $\mathcal H$ be any $A_{X,\Theta}\times G$-invariant closed subspace of $L^2(X_\Theta)_\disc$ and consider the ``scattering'' morphisms $\iota_{\Omega}^*\iota_\Theta$ restricted to $\mathcal H$. If in their decomposition (\ref{SsumSi}) the summand $S_i$ is non-zero on $\mathcal H$, we will say that the pair $(T_i,f_i)$ (or to be absolutely complete the triple $(\Omega, T_i, f_i)$) {\em appears in the scattering of $\mathcal{H}$.} 

The main result of this section is the following:

\begin{theorem}[Tiling property of scattering morphisms.] \label{tilingtheorem} 
Let $\Theta\subset\Delta_X$ and let $\mathcal H$ be a  nonzero $A_{X,\Theta}\times G$-invariant closed subspace of $L^2(X_\Theta)_\disc$.
\begin{enumerate}
 \item If a triple $(\Omega,T,f)$ appears in the scattering of $\mathcal H$ and:
\begin{equation}\label{scatteringtwo}\mathring{\mathfrak a}_{X,\Theta}^+  \cap T \mathring{\mathfrak a}_{X,\Omega}^+ \ne \emptyset\end{equation}
then $\Omega=\Theta$, $T=\Id$ and $f$ is also the identity (on a connected component of $\widehat{A_{X,\Theta}'}$). 

 \item If $(\Omega_i,T_i, f_i)$ varies in all the triples  which appear in the scattering of $\mathcal H$, then:  
\begin{equation}\label{scatteringone}\bigcup_i T_i (\mathfrak a_{X,\Omega_i}^+) = \mathfrak a_{X,\Theta}.\end{equation}  

Indeed, 
 there exists a splitting $\mathcal{H} = \bigoplus \mathcal{H}_{\alpha} $ 
such that:
\begin{itemize}
\item[(i)] The $A_{X, \Theta}'$-Plancherel measure for different $\mathcal{H}_{\alpha}$
is mutually singular; in particular, $\Hom_{A_{X, \Theta}}(\mathcal{H}_{\alpha}, \mathcal{H}_{\beta}) = 0$ for $\alpha \neq \beta$;
\item[(ii)] For $a \in \mathfrak a_{X,\Theta}$,  let $J$ denote the set of indices $i$ such that $a\in T_i (\mathring{\mathfrak a}_{X,\Omega}^+)$, 
and assume that $a$ is generic in the sense that $a$ does not lie on any wall of $T_i \mathfrak{a}_{X, \Omega}^+$. 
 For any $\alpha$, any $v \in \mathcal{H}_{\alpha}$ and any generic  $a \in \mathfrak{a}_{X, \Theta}$ 
then   \begin{equation} \label{largescattering}\sum_{i\in J} \left\Vert  S_i(v)  \right\Vert^2  \ge \|v\|^2,
\end{equation}
  
\end{itemize}
 \end{enumerate}
\end{theorem}

This will be enough to prove the main Scattering Theorem \ref{advancedscattering}. Let us first discuss this proof. 
Actually, we only use the {\em second} statement of Theorem \ref{tilingtheorem} in this proof; 
the first statement won't be used,  because it is contained in the ``generic injectivity'' condition. 

\begin{proof}[Proof of Theorem \ref{advancedscattering}] \label{proof:advancedscattering}

The existence and characterization of Bernstein morphisms, together with the fact that the images of their restrictions to discrete spectra span the whole space $L^2(X)$, has already been established in Section \ref{sec:Bernstein} (see Corollary \ref{corollarysurjection}). The rest of the statements of the theorem will first be proved by restriction to discrete spectra, i.e.\ $S_w$, for $w\in W_X(\Omega,\Theta)$, will first be defined as a morphism with the stated properties from $L^2(X_\Theta)_\disc$ to $L^2(X_\Omega)_\disc$. Notice that (\ref{totalmapdef}) involves only discrete spectra. At the end we will extend $S_w$ to the whole $L^2(X_\Theta)$.

Take $\Omega, \Theta$ with $|\Omega| = |\Theta|$.  We write $\iota_{\Theta, \disc}$ for the restriction of $\iota_{\Theta}$
to $L^2(X_{\Theta})_{\disc}$ and $\iota_{\Theta, \disc}^*$ the adjoint of this restricted map, i.e. it is  $\iota_{\Theta}^*$ followed by
  the orthogonal projection to discrete spectrum. 
By Lemma \ref{injectivefaces}, the morphism 
$\iota_{\Omega,\disc}^* \iota_{\Theta}=\iota_{\Omega}^* \iota_{\Theta,\disc}=\iota_{\Omega,\disc}^* \iota_{\Theta,\disc} $ (the equalities follow from  claim \eqref{orthogonality} of Proposition \ref{decompf}\footnote{
For example, $\iota_{\Omega,\disc}^* (\iota_{\Theta}-\iota_{\Theta,\disc})$ is zero because $L^2(X)_{\Omega}$ is orthogonal
to $\iota_{\Theta} f$ for $f \perp L^2(X_{\Theta})_{\disc}$; and the latter follows because such $\iota_{\Theta} f$
can be expressed as a sum of $\iota_{\Theta'} f'$ for some $|\Theta'| \neq |\Omega|$, by Corollary \ref{corollarysurjection} and Proposition  \ref{propcomposition}.
 }) decomposes as a sum of morphisms 
$\sum_{w \in W_X(\Omega,\Theta)} S_w$, where
each $S_w$ is $A_{X, \Theta}'$-equivariant, with
$A_{X, \Theta}'$
acting on $L^2(X_{\Omega})$ via the isomorphism $A_{X, \Theta}' \stackrel{w}{\rightarrow} A_{X, \Omega}'$ induced by $w$. In particular, $\iota_{\Omega,\disc}^* \iota_{\Theta,\disc} $ is zero unless $\Theta$ and $\Omega$ are $W_X$-associate.

Now, root systems have the following tiling property:  
\begin{quote}   The collection of subsets of $\mathfrak{a}_{X, \Theta}$ given by
$w^{-1} a_{X, \Omega}^+$, where $\Omega$ varies through subsets of $\Delta_X$
with $|\Omega| = |\Theta|$, and $w$ ranges through $W_X(\Omega,\Theta)$, 
form a perfect tiling of $\mathfrak{a}_{X, \Theta}$, i.e.\ their union is $\mathfrak{a}_{X, \Theta}$
and their interiors are disjoint. 
\end{quote} 
A proof of this property  has been indicated in the footnote of \S \ref{subsec:mainresult}.

In particular, for each  $\Omega$ and $w \in W_X(\Omega,\Theta)$, if $a\in w^{-1} \mathring{\mathfrak a}_{X,\Omega}^+ $
then $a$ does not belong to the corresponding set for any other pair
$(\Omega', w' \in W_X(\Omega',\Theta))$.
 So the second statement  of Theorem \ref{tilingtheorem}  implies that
we may split $L^2(X_{\Theta})_{\disc} = \bigoplus \mathcal{H}_{\alpha}$ 
in such a way that 
\begin{equation} \label{Swnorm} \|S_w v\| \geq \|v\| \end{equation} whenever $v$ belongs to any $\mathcal{H}_{\alpha}$. {\bluetext (Recall that $S_w$ is the sum of $S_i$'s which are induced from the element $w$, and those have image on different orthogonal direct summands of $L^2(X_{\Omega})$, corresponding to distinct connected components of the character group $\widehat{A_{X,\Omega}'}$, hence $\| S_w v\|^2 = \sum_i \|S_i v\|^2$, the sum ranging over those $i$'s.)
}

{ We can also assume that Proposition \ref {Bmapisometry}} {\bluetext (which can easily be seen to hold when $A_{X,\Theta}$ is replaced by $A_{X,\Theta}'$)} { applies to any $\mathcal H_{\alpha}$ and its image under any $S_w$.
In fact, choose a partition $\widehat{A_{X, \Theta}'} = \coprod U_{\beta}$ up to sets of measure zero,  with the property that  if
 if $w \in U_{\beta}$ and   $1\ne w \in W_X(\Theta,\Theta)$ then $w \chi \notin U_{\beta}$.  We may assume, without loss of generality,
 that the decomposition $\mathcal{H} = \bigoplus \mathcal{H}_{\alpha}$ furnished by the Theorem is fine enough that
the Plancherel measure for $\mathcal{H}_{\alpha}$ as an $A_{X, \Theta}'$-representation is wholly supported on some $U_{\beta}$. This is enough
to ensure that Proposition \ref{Bmapisometry} applies to $\mathcal H_{\alpha}$ and its image under any $S_w$: 
our choice and the decomposition \eqref{SsumSi} mean that any $ G$-equivariant map $\mathcal{H}_{\alpha} \rightarrow \mathcal{H}_{\alpha}$
is also $A_{X, \Theta}'$-equivariant, which means that almost all $\pi \in \hat{G}$ 
for  $\mathcal{H}_{\alpha}$-Plancherel 
measure  have a unique exponent, in the language of Proposition \ref{Bmapisometry}.}

Fix $w_0 \in W_X(\Omega,\Theta)$. 
  Let $\mathcal{H}'_{\alpha}$
be the image of $\mathcal{H}_{\alpha}$ under $S_{w_0}$. 
 Let $p: \mathcal{H}'_{\alpha} \rightarrow L^2(X_{\Omega})_{\disc}$ be the natural inclusion.  Then $S_{w_0} = p^* S_{w_0}$.  We may write
$$ (\iota_{\Omega} \circ p)^* \circ \iota_{\Theta} =  S_{w_0} + \sum_{w \neq w_0} p^* S_w.$$
 Now $p^* S_w|_{\mathcal{H}_{\alpha}} = 0$ for $w \neq w_0$, as we see by a consideration
 of $A_{X, \Omega}'$-character.   { Indeed choose (as in the discussion after \eqref{Swnorm})
 a set $U_{\beta} \subset \widehat{A_{X, \Theta}'}$ containing the $A_{X, \Theta}'$-support of $\mathcal{H}_{\alpha}$;
 then $\mathcal{H}_{\alpha}'$ is supported on $A_{X, \Omega}'$-characters in $w_0 U_{\beta}$,
 also $S_w \mathcal{H}_{\alpha}$ is supported on $A_{X, \Omega}'$-characters in $w U_{\beta}$,
 and $w_0 U_{\beta} \cap w U_{\beta}=\emptyset$ because of the way the sets $U_{\beta}$ were chosen (discussion after \eqref{Swnorm}). }
 Therefore:
\begin{equation} \label{Sw0}  \left.(\iota_{\Omega} \circ p)^* \circ \iota_{\Theta}\right|_{\mathcal H_{\alpha}} =  S_{w_0}\end{equation} 
But $\iota_{\Theta}|_{\mathcal H_{\alpha}}$ and $\iota_{\Omega} \circ p$ are both isometries according to Proposition \ref{Bmapisometry}.
In order for \eqref{Swnorm} to hold, then, the images $\iota_{\Theta}(\mathcal{H}_{\alpha})$
must be contained in $\iota_{\Omega}(\mathcal{H}'_{\alpha})$, and then $S_{w_0}$ is an isometry
from $\mathcal{H}_{\alpha}$ to $\mathcal{H}'_{\alpha}$.

 In particular (taking the sum over $\alpha$)   the image $L^2(X)_{\Theta}$ of $L^2(X_{\Theta})_{\disc}$
 under $\iota_{\Theta}$ is contained in $L^2(X)_{\Omega}$; by symmetry, the two coincide:
 $$L^2(X)_{\Theta} = L^2(X)_{\Omega}.$$ 
 Also,  $\mathcal{H}_{\alpha}'$ and $\mathcal{H}_{\beta}'$ are 
mutually orthogonal if $\alpha \neq \beta$: that follows by a consideration
of the $A_{X, \Omega}'$-action, in particular using property 2(i) 
 from Theorem \ref{tilingtheorem} and the equivariance property of $S_{w_0}$.
  Finally, $\sum_{\alpha} \mathcal{H}'_{\alpha}=
L^2(X_{\Omega})_{\disc}$:
the orthogonal complement $\mathcal{K} \subset L^2(X_{\Omega})_{\disc}$
of all  $\mathcal{H}'_{\alpha}$ is an $A_{X, \Omega}' \times G$-stable space
which is exactly the kernel of
$S_{w_0}^* $. But it is not hard to see that the adjoint
of $S_{w_0}: L^2(X_{\Theta})_{\disc} \rightarrow L^2(X_{\Omega})_{\disc}$
is $S_{w_0^{-1}}: L^2(X_{\Omega})_{\disc} \rightarrow L^2(X_{\Theta})_{\disc}$;
so $\mathcal{K}$ would belong to the kernel of $S_{w_0^{-1}}$; that contradicts  the second part of Theorem \ref{tilingtheorem} applied  { to $\mathcal{H} =\mathcal{K}$, where we choose
 $a$  similarly to the discussion before \eqref{Swnorm} to produce a vector $v$ with $\|S_{w_0^{-1}}v\| \geq \|v\|$.}

Therefore, $S_{w_0}$ gives an isometry
$$L^2(X_{\Theta})_{\disc} = \bigoplus \mathcal{H}_{\alpha} \rightarrow L^2(X_{\Omega})_{\disc} =
\bigoplus \mathcal{H}'_{\alpha}.$$

We note for later use that
\begin{equation}\label{genA} \mbox{The $A_{X, \Omega}'$-span of $\iota_{\Omega,\disc}^* L^2(X)$ 
is dense in $L^2(X_{\Omega})_{\disc}$.}
\end{equation} 
Indeed  the $A_{X, \Omega}'$-span of
$\iota_{\Omega, \disc}^* \iota_{\Theta} \mathcal{H}_{\alpha}$
must contain $\mathcal{H}'_{\alpha}$.

  We also saw above that the images $L^2(X)_{\Theta}, L^2(X)_{\Omega}$ are orthogonal if they are not associate, hence we can refine the decomposition (\ref{directbydim}) as:
\begin{equation}
 L^2(X)= \oplus_{\Theta/\sim} L^2(X)_\Theta,
\end{equation}
where $\sim$ denotes the equivalence relation of being $W_X$-associate.

Now, with the same notation as what we have just proved: 

\begin{enumerate}
\item $\iota_{\Omega} p^* \iota_{\Omega}^*$ is the identity
on $\iota_{\Theta}(\mathcal{H}_{\alpha})$: 

As we just saw, $\iota_{\Omega} \circ p$ maps $\mathcal{H}'_{\alpha}$ 
isometrically onto a subspace containing $\iota_{\Theta}(\mathcal{H}_{\alpha})$. 
Therefore $\iota_{\Omega} \circ p  \circ ( \iota_{\Omega} \circ p)^*$ is the identity on $\iota_{\Theta}(\mathcal{H}_{\alpha})$,
which implies the claim.

 \item $\iota_\Omega\circ S_{w_0} = \iota_\Theta$. Indeed,  \label{pageref_for_invariance} 
 it is enough to check this on each $\mathcal{H}_{\alpha}$,
 and there, by \eqref{Sw0}, 
$$\iota_\Omega\circ S_{w_0}|_{\mathcal H_{\alpha}} = \iota_\Omega\circ p^*\circ \iota_\Omega^*\circ \iota_\Theta|_{\mathcal H_{\alpha}} = \iota_\Theta|_{\mathcal H_{\alpha}}$$
by what we just showed.
 \item $S_{w'}\circ S_w = S_{w'w}$ when $w\in W_X(\Omega,\Theta)$ and $w'\in W_X(Z,\Omega)$.

Indeed, write $w_0$ instead of $w$, for compatibility with our prior notation. Then by  what we just showed $$ \iota_Z^*\iota_\Omega S_{w_0}|_{\mathcal H_{\alpha}} = \iota_Z^*\iota_\Theta|_{\mathcal H_{\alpha}},$$
and thus $\iota_Z^* \iota_{\Omega} S_{w_0} = \iota_Z^* \iota_{\Theta}$
as morphisms $L^2(X_{\Theta})_{\disc} \rightarrow L^2(X_Z)_{\disc}$. 

Hence the $w'w$-equivariant part of $\iota_Z^*\iota_\Theta$  coincides with the $w'$-equivariant part of $\iota_Z^*\iota_\Omega $ composed with $S_w$.
\end{enumerate}

Define the endomorphism of $\sum_{\Omega\sim\Theta} L^2(X_\Omega)_\disc$.
$$S := \sum_{\Omega_1,\Omega_2\sim \Theta, w\in W_X(\Omega_2,\Omega_1)} S_w.$$  
 Note that $S^2 =  c(\Theta) \cdot S$; here $c(\Theta)$ is as in the statement of Theorem \ref{advancedscattering}: it is the number of chambers in $\mathfrak{a}_{X, \Theta}$, or, what is the same, the sum $\sum_{\Omega \sim \Theta} \# W_X(\Theta,\Omega)$ (the equality follows from the tiling result, just as in the discussion of \S \ref{subsec:mainresult}).

     In particular, 
$$\bar{S} := \frac{S}{c(\Theta)}$$ is a projection.  
Its image is the  ``$S_w$-invariants'' (i.e.\ invariants over all possible $S_w$ between different $\Omega$'s in this associate class).

Set $T=\oplus_{\Omega\sim\Theta} \iota_{\Omega,\disc}$ (as an operator from $\sum_{\Omega\sim\Theta} L^2(X_\Omega)_\disc$ to $L^2(X)$), hence $T^*=\oplus_{\Omega\sim\Theta} \iota_{\Omega,\disc}^*$.    Then $T^* T = S$;  it follows then that the image of $T^*$
contains the image of $\bar{S}$, i.e.\ the ``$S_w$-invariants.''  On the other hand,   
the reverse containment follows from what we have already shown.   That is to say,
$\bar{S} T^* = T^*$, which follows from (2) above together with the formula $c(\Theta) = \sum_{\Omega \sim \Theta} \# W_X(\Omega,\Theta)$ already noted.  
So the image of $T^*$ is \emph{precisely} as stated.

We have proven all statements of the theorem when the $S_w$'s are defined only on the discrete spectrum. We now turn to ``upgrading'' them so they apply to the entire spectrum.

We notice the following: $S_w$ can be characterized as the \emph{unique} $(A_{X,\Theta}', w)$-equivariant isometry from $L^2(X_\Theta)_\disc$ to $L^2(X_\Omega)_\disc$ which satisfies: $\iota_{\Theta,\disc}=\iota_{\Omega,\disc}\circ S_w$. Indeed, this condition identifies $S_w^*$ on the image of $\iota_{\Omega,\disc}^*$ and we apply \eqref{genA}.

Let now $\Theta\subset\Delta_X$, $Z_1,Z_2\subset \Theta$ and $w\in W_{X_\Theta}(Z_1,Z_2) \subset W_X(Z_2,Z_1)$. In particular, we have scattering morphisms: 
$$S_w, S_w^\Theta: L^2(X_{Z_1})_\disc\to L^2(X_{Z_2})_\disc,$$
from applying the part of the theorem which is already proven to the varieties $X$ and $X_\Theta$, respectively. We claim that they coincide, i.e.
\begin{equation} \label{SwSwT} S_w = S_w^{\Theta} .\end{equation}Indeed, we have $\iota_{Z_i} = \iota_\Theta\circ \iota_{Z_i}^\Theta$ (Proposition \ref{propcomposition}), and the claim follows from the above characterization of $S_w$:
$$ \iota^{\Theta}_{Z_1} = \iota^{\Theta}_{Z_2} \circ S_w^{\Theta} \implies \iota_{Z_1} = \iota_{Z_2} \circ S_w^{\Theta}.$$

Take now $w\in W_X(\Omega,\Theta)$.
We may now define $S_w$ in general. We have a decomposition:

$$L^2(X_{\Theta}) = \bigoplus_{\{Z|Z \subset \Theta\} / \sim} L^2(X_{\Theta})_Z$$
where $L^2(X_{\Theta})_Z$ is the image of $L^2(X_Z)_{\disc}$ in $L^2(X_{\Theta})$, 
by the map $\iota^{\Theta}_Z$. 

We need, then, to define $S_w$ on each space $L^2(X_\Theta)_Z$. 
Put $Y = w(Z) \subset \Omega$. We define $S_w$ 
by requiring the following diagram to commute
\begin{equation} 
 \begin{CD}
L^2(X_{\Theta})_Z @>{S_w}>> L^2(X_{\Omega})_Y \\ 
@V{\iota^{\Theta,*}_Z}VV @V{\iota^{\Omega,*}_Y}VV \\
L^2(X_Z)_{\disc} @>S_w>>  L^2(X_Y)_{\disc} \\ 
 \end{CD}
\end{equation}
In fact, this diagram can be made to commute: {\bluetext by what we have just proven for discrete spectra (applied to $X_\Theta$ instead of $X$)}, the left-hand vertical arrow identifies $L^2(X_{\Theta})_Z$
with that subspace of $L^2(X_Z)_\disc$ that is invariant by all $S_w^{\Theta}$,
$w \in W_{X_\Theta}(Z,Z)$.

Similarly, the right-hand vertical arrow
identifies $L^2(X_{\Omega})_Y$ with that subspace of $L^2(X_Y)_{\disc}$
that is invariant by all $S_w^{\Omega}$, for $w \in W_{X_\Omega}(Y,Y)$. 
Because of the composition property for the maps $S_w$ and \eqref{SwSwT},
$S_w$ maps the first space isomorphically to the second. 

 It remains only to verify that this does not depend on the choice of $Z$ within its associate class;
 this is routine and we omit it.

\end{proof}

\subsection{Proof of the first part of Theorem \ref{tilingtheorem}} \label{prooffirstpart} 

{\bluetext For notational simplicity, in this proof we denote $A_{X,\Theta}'$, $A_{X,\Omega}'$ simply by $A_{X,\Theta}$, $A_{X,\Omega}$.}

The basic idea is to consider the inner product $\langle \iota_{\Omega} \Phi, \iota_{\Theta} \Psi \rangle$ for suitable $\Phi \in L^2(X_{\Omega}), \Psi \in L^2(X_{\Theta})$.  We
then ``push'' $\Phi$ and $\Psi$ towards infinity using 
suitable elements of $A_{X, \Omega}$ and $A_{X, \Theta}$. The given assumption
(after some analysis of scattering) forces this inner product
to converge to zero. We then derive a contradiction by comparing with the decomposition
of $\iota_{\Omega}^* \iota_{\Theta}$.

Suppose that the first property is false;   that is, there exists a triple $(\Omega,T,f)$ appearing in the scattering of $\mathcal H$ such that: $\mathring{\mathfrak a}_{X,\Theta}^+  \cap T \mathring{\mathfrak a}_{X,\Omega}^+ \neq\emptyset$. 

This means that there is a finitely generated subsemigroup $M^+$ of $A_{X,\Omega}^T$ with the following properties:
\begin{enumerate}
 \item $M^+ \subset \mathring A_{X,\Omega}^+$ and $T(M^+)\subset \mathring A_{X,\Theta}^+$;
 \item $M^+$ and $A_{X,\Omega}^0$ generate $A_{X,\Omega}^T$ as a group.
\end{enumerate}

We are slightly abusing notation here, since $Ta$ does not always make sense as an element of $A_{X,\Theta}$ when $a\in A_{X,\Omega}^T$; however, it does make sense as an element of $A_{X,\Theta}/A_{X_\Theta}^0$ via the valuation map and (\ref{atmap}), therefore the statements above make sense. Let $M$ be the group generated by $M^+$. In order to avoid similar clarifications in the rest of the proof, let us choose a homomorphism: $M\to A_{X,\Theta}$ which lifts $T: M \rightarrow  A_{X,\Theta}/A_{X,\Theta}^0$. By abuse of notation, we will be denoting this homomorphism by $T$ again. If $\Theta=\Omega$, $T=\Id$, we take this homomorphism to be the identity.

Let $S_j$ denote the summand of $\iota_\Omega^*\iota_\Theta$ corresponding to the quadruple $(\Omega,T,f)$ in the decomposition (\ref{SsumSi}). If $\Phi\in L^2(X_\Theta)_\disc^J$, then the equivariance property of $S_j$ reads:
\begin{equation}\label{equiv}S_j(T(m) *  \Phi) = m * S_j \Phi,
\end{equation}
for $m\in M$, where we have twisted the actions of the tori as in (\ref{Siequivariance}). 

When $m\in M^+$, the fact that $M^+\subset \mathring A_{X,\Omega}^+$ implies, in particular, that the mass of the function $S_j  (T(m)^n * \Phi)\in L^2(X_\Omega)_\disc^J$ will be ``moving towards a $J$-good neighborhood of infinity in $X_\Omega$ as $n\to \infty$.'' { By the quoted phrase,  we mean the following: writing $f_j = S_j  (T(m)^n * \Phi)$
and $N_{\Omega}$ for the $J$-good neighourhood of infinity, the norm $\|f_j\|_{L^2(X-N_{\Omega})} \rightarrow 0$ as $j \rightarrow \infty$. 
This follows from the fact (in turn from Lemma \ref{limit-lemma-goo}) that for any compact set $O \subset X_{\Omega}$ and for $n \gg 1$,  we have $m^n O \subset N_{\Omega}$,
and then choosing $O$ so large that the $L^2$-norm of $f$ on $X_{\Omega}-O$ is arbitrarily small. }

Now let $\Phi\in L^2(X_\Theta)_\disc^J, \Psi\in  L^2(X_\Omega)_\disc^J$.  Let us choose $\Phi$ to be $(A_{X,\Theta}^0,\eta)$-equivariant and $\Psi$ to be $(A_{X,\Omega}^0,\eta')$-equivariant, where the restrictions of $\eta$ and $\eta'$ to the maximal compact subgroups are those that correspond to the domain and image of $f$ (in particular, $S_j\Phi$ is $(A_{X,\Omega}^0,\eta')$-equivariant).

Choose an averaging sequence of measures $\nu_n$ on $M^+$   (\S \ref{averaging sequence}) and consider the inner products:
\begin{equation}\label{iptoinfty}
P_n:= \sum_{M} \nu_n(m) \left< \iota_\Theta (T(m)* \Phi),\iota_\Omega (m* \Psi)\right>.
\end{equation}
where we regard $\nu_n(m)=0$ off $M^+$. The sum is convergent:
the inner products that appear are bounded independently of $m \in M^+$, by Cauchy-Schwarz
and the boundedness of $\iota_{\Theta}, \iota_{\Omega}$; and $\sum \nu_n(m)=1$. 
We now evaluate $P_n$ in two different ways -- (i) and (ii) below. 
We suppose now also  that we are  {\em not} in the case $\Omega=\Theta, T = \mathrm{Id}$
and $f$ the identity. 

\begin{itemize}
\item[(1).] We show that $P_n\to 0$ as $n\to\infty$ (for suitable choice of the groups $M$, $T(M)$).

As $n\to \infty$, by the property (\ref{bchar}) characterizing the Bernstein morphisms we know that for $a\in M^+$ and $b=T(a)$, if $\Phi$ was compactly supported then $\Vert \iota_\Theta {b^n}* \Phi-e_\Theta {b^n}*\Phi\Vert\to 0$, and similarly for ${a^n}*\Psi$, where $e_\Theta: C_c^\infty(X_\Theta)\to C_c^\infty(X)$ is the asymptotics map.  

We also noted in \eqref{bcharfornoncompact} a certain analog of this statement 
for functions that are not necessarily compactly supported. 
That implies that, if we fix $J$-good neighborhoods $N_\Theta$, $N_\Omega$ of $\Theta$- and $\Omega$-infinity, then for \emph{any} elements in $L^2(X_\Theta)^J$, resp.\ $L^2(X_\Omega)^J$, we get that (\ref{iptoinfty}) approaches the inner product:
\begin{equation}\label{iptoinfty2}   P_n' := \sum_{M} \nu_n(m) \left<  \tau_\Theta (T(m)* \Phi),  \tau_\Omega (m* \Psi) \right>,
\end{equation}
{ Here the meaning of ``approaches'' is that    $\lim_n(P_n -P_n' )= 0$ and 
 $\tau_\Theta$ and $\tau_\Omega$ denote truncation of the given $J$-invariant functions, even when they are not compactly supported, to $J$-good neighborhoods of infinity, and identification via the exponential map with functions on $X$.  Indeed,  because each $\nu_n$ is a probability measure, it is enough to verify that
 the pointwise difference $\left< \left. \tau_\Theta (T(m)* \Phi)\right|_{N_\Theta}, \left. \tau_\Omega (m* \Psi) \right|_{N_\Omega}\right> -  \left< \iota_\Theta (T(m)* \Phi),\iota_\Omega (m* \Psi)\right>$,
 restricted to the support of $\nu_n$, approaches $0$ as $n \rightarrow \infty$.   In turn, that follows from
a slight generalization of  \eqref{bcharfornoncompact} (from one-parameter groups to several-parameter groups),  
 the fact that the support of $\nu_n$ is contained deeper and deeper in the interior of $M^+$ (more precisely,  the first  noted property of an averaging sequence,
 see before   Lemma \ref{limandaveraging}), 
 and the assumed facts that $M^+ \subset \mathring A_{X,\Omega}^+$ and $T(M^+)\subset \mathring A_{X,\Theta}^+$. } 

There are  now two cases. For what follows, denote by $a$ an element of $M^+$, and by $b$ its image under $T$. 
{ If $\Omega=\Theta$ but $T$ is not the identity we can and do choose $a$ so that $b \neq a$.}
 \begin{itemize}
 \item[(i)] \label{random-i-ref} $\Omega\ne \Theta$, or $\Omega = \Theta$ and $T\ne \Id$, in which case an arbitrarily large percentage of the masses of $\left. \tau_\Theta {b^n}* \Phi\right|_{N_\Theta}$ and $ \left. \tau_\Omega {a^n}* \Psi\right|_{N_\Omega}$ is eventually concentrated on disjoint sets. {  In other words, for every $\varepsilon  > 0$, there are subsets $N_{\Theta}' \subset N_{\Theta}$
 and $N_{\Omega}' \subset N_{\Omega}$, disjoint when identified with subsets of $X$,  with the property that  \label{daggerpageref}
\begin{equation}  \label{dagger}  \ \ \  \| \tau_\Theta {b^n}* \Phi\|_{L^2(N_\Theta')} > (1-\varepsilon) \| \tau_\Theta {b^n}* \Phi\|_{L^2(N_\Theta)}\end{equation} 
 and similarly for $\Omega$:
\begin{equation}\label{dagger2} \ \ \  \| \tau_\Omega {a^n}* \Psi\|_{L^2(N_\Omega')} > (1-\varepsilon) \| \tau_\Omega {a^n}* \Psi\|_{L^2(N_\Omega)}.\end{equation}
 This, in turn, comes from the following: Suppose first that $\Theta \neq \Omega$.  Take $J$-stable compact subsets
 $O \subset X_{\Theta}, O' \subset X_{\Omega}$ which support a sufficiently large percentage of the $L^2$-norms of $\Phi$, resp.\ $\Psi$. Then $a^n O$ and $b^n O'$ become disjoint for $n \gg 1$ when
 identified with subsets of $X$.  Indeed, after covering $O, O'$ by finitely many orbits of a compact subgroup of $G$, 
  it is enough to note that, taking $x \in O$ and $y \in O'$,
 the limit points of $a^n x$ and $b^n y$, with this identification, belong to distinct $G$-orbits on $\overline{X}$.
 That follows from Lemma \ref{limit-lemma-goo}.  
 Next, suppose that $\Theta = \Omega$ but $T$ is not the identity, so that $b \neq a$.   We now use similarly the fact that, for any compact set $O \subset X_{\Theta}$, the sets $a^n O, b^n O$ are eventually disjoint, which is clear,
 working $A_{X, \Theta}$-orbit by $A_{X, \Theta}$-orbit on $X_\Theta/J$. 
 } 
 
{\bluetext Taking now $N_\Theta'$, resp.\ $N_\Omega'$, to be the union of all $a^n$- (resp.\ $b^n$-) translates of $O$ (resp.\ $O'$, when $\Omega\ne\Theta$) for all large enough $n$, this shows \eqref{dagger}, \eqref{dagger2}, where we remind that the point is that $N_\Theta'$ and $N_\Omega'$ are disjoint. By the Cauchy-Schwartz inequality, this shows that the inner products averaged in \eqref{iptoinfty2} become arbitrarily small, for large $n$.}
 
\item[(ii)] $\Omega = \Theta$, $T=\Id$, but $f\ne \Id$.

Let us use notation defined before \eqref{Siequivariance}, in particular the twisted  * action that was defined there, which
is simply the twist of the usual action by a unitary character. 
Let us also take $a=b$. The definition of the * action gives:
$$\left< \tau_\Theta (b^n * \Phi),  \tau_\Theta (a^n* \Psi) \right> = 
(\eta'\eta^{-1})(a^n) \cdot \left<  \tau_\Theta \mathcal L_{a^n} \Phi,  \tau_\Theta \mathcal L_{a^n} \Psi\right>,$$
where $\eta$ and $\eta'$ are the (unitary) characters defining the twisted  * action of $a$, as before \eqref{Siequivariance}. Since $f\ne \Id$, we can choose $M$ so that $(\eta'\eta^{-1})|_M\ne 1$.
{ Then $P_n' \rightarrow 0$:
the definition of $P_n'$ now takes the shape \label{daggerdaggerref} 
\begin{equation}  \label{daggerdagger} P_n' = \sum_{m \in M}  \nu_n(m) (\eta'\eta^{-1})(m)   \left<  \tau_\Theta \mathcal L_{m} \Phi,  \tau_\Theta \mathcal L_{m} \Psi\right> \end{equation} 
and, for any $\varepsilon  > 0$, we may choose $n$ so large so the inner products appearing here 
are almost constant, i.e. there is a constant $A$ such that $| \left<  \tau_\Theta \mathcal L_{m} \Phi, \tau_\Theta \mathcal L_{m} \Psi\right> - A| < \varepsilon$ whenever $m$ belongs to the support of $\nu_n$.  (Were it not for the truncation operators, we could even take $\varepsilon  =0$, i.e. the inner products
would be exactly constant.) 
Thus, $|P_n'| \leq \varepsilon + \sum_{M} \nu_n(m) (\eta' \eta^{-1})(m)$.  
The property of averaging sequences that $S_n := \sum_{M} \nu_n(m) \chi(m)$ approaches zero as $n \rightarrow \infty$ for any nontrivial character $\chi$: 
to see this, choose $s \in M$ with $\chi(s)\neq 1$ and observe that $(\chi(s) S_n - S_n) = \sum_{M} (\nu_n(ms^{-1} ) - \nu_n(m)) \chi(m) $, and the latter sum approaches
zero as $n \rightarrow \infty$ by the second  property of an averaging sequence (see \S \ref{fpolsubsec}). }

\end{itemize}

\item[(2)] $P_n \rightarrow  \langle S_j \Phi, \Psi \rangle$.

Each term of the sum (\ref{iptoinfty}) is equal to:
$$ \left< \iota_\Omega^*\iota_\Theta (T(m)* \Phi), m* \Psi\right> = \sum_i \left< S_i (T(m)* \Phi), m* \Psi\right>.$$
We now show that, upon applying the average (\ref{iptoinfty}) to this expression, only the $(T,f)$-equivariant summand survives, that is, $S_j$. 

First of all, the $A_{X,\Theta}^0$ and $A_{X,\Omega}^0$-equivariance properties of $\Phi$ and $\Psi$ kill all summands $S_i$ such that the domain or image of $f_i$ is different from that of $f$. 

Next, for those $S_i$ with $T_i \neq T$ one sees by a similar argument to those already given that 
 the inner product $\left< S_i(T(m)*\Phi), m* \Psi\right>$ eventually approaches zero. 
{  Let $U = T_i^{-1} T$, which is defined as a nontrivial isogeny 
 $\mathbf{A}_{X,\Omega} \dashrightarrow \mathbf{A}_{X, \Omega}$. 
 If $m$ lies in the subgroup of $M$ such that $U(m)$ is defined,  
then $\left< S_i (T(m)*\Phi), m * \Psi \right> $ has the same absolute value as $\langle \mathcal{L}_{U(m)} S_i \Phi, \mathcal{L}_m \Psi \rangle$, where we no longer
twist the actions (twisting only affects the result by a scalar of absolute value $1$).  Note that $U \neq 1$, and  we now reason in exactly the same as
as the $\Theta=\Omega$ case after \eqref{dagger} on page \pageref{daggerpageref}.}

Finally, if $T_i=T$, the associated affine map $f_i$ of character groups differs
from the affine map $f$ for $T$ simply by multiplication by a character $\eta$ of $A_{X, \Omega}$, i.e.\ $f_i(\chi) = f(\chi) \eta$. 
Then we have, by \eqref{Siequivariance}, 
$$S_i (T(m) * \Phi) = \eta(m) m * S_i( \Phi).$$  
By the argument presented after\footnote{Our situation is even simpler now, because there are no truncations involved. All that is needed is the final statement of that argument involving $S_n$.} \eqref{daggerdagger} on page \pageref{daggerdaggerref} the weighted average of terms $\langle S_i T(m) *\Phi, m*\Psi\rangle$ will therefore tend to zero, unless $\eta=1$ which implies that $(T_i,f_i)=(T,f)$, i.e.\ $S_i$ is the summand $S_j$.

But we now see that
$$P_n \longrightarrow  \left<S_j \Phi,\Psi\right> \ \ (n \rightarrow \infty)$$
because $  \left< S_j (T(m)* \Phi), m* \Psi\right> = \left<S_j \Phi,\Psi\right>$ by \eqref{equiv}. 
\end{itemize}

Taking (i) and (ii) together:   $\left<S_j\Phi,\Psi\right>  = 0$ for all $\Phi, \Psi$, which is impossible unless $S_j \equiv 0$. 
  \qed

 The second part will require a further ingredient: So far, we have not used the fact that
 the complement of all the neighborhoods of $\infty$ is actually compact modulo the center.
Roughly speaking, we will show that if the second part of the theorem were not the true then we would be able to ``push'' any function $\Phi\in L^2(X)_{\Theta}$ away from infinity, without changing its $L^2$ norm and keeping its $L^\infty$ norm under control, which will lead to a contradiction. The argument is delicate, and we complete it in the remaining part of this section.

\subsection{Estimates}  \label{ssestimates}  

In this section we fix an open compact subgroup $J$ of $G$ and develop an estimate   for the norm on $L^2(X)^J$ in terms of certain norms for the ``constant terms'' $\iota_\Omega^*(\Phi)$ (for all $\Omega\subset\Delta_X$). The only input that we use is the uniform bound of subunitary exponents for $L^2(X)^J$, Proposition \ref{uniformbound}, which allows us to apply some results of the ``linear algebra'' section \ref{sec:linearalgebra}.   

We fix, for each $\Omega\subset \Delta_X$, a $J$-good neighborhood $\tilde N_{\Omega}$ of $\Omega$-infinity which is stable under $A_{X,\Omega}^+$ (when $\tilde N_\Omega/J$ is identified with a subset of $X_\Omega/J$). Let $N_\Omega=\tilde N_\Omega\smallsetminus \cup_{\Theta \subsetneq\Omega} \tilde N_\Theta$; then the sets $N_\Omega$ partition $X$, and $N_\Omega$ is represented by a finite number of $J$-orbits modulo the action of $A_{X,\Omega}^+$. { (Indeed, $\tilde N_\Omega$ is the intersection with $X$ of an
actual neighborhood $\tilde N_\Omega'$ of $\infty_\Omega$ in a wonderful compactification $\tilde X$, and if we remove the corresponding neighborhoods for $\Theta\subsetneq\Omega$ then the remaining compact subset of $\tilde N_\Omega'$ only intersects $\infty_\Omega$ along the orbit (or orbits) corresponding to $\Omega$-infinity.)}

We denote by $\tau_\Omega$ the operator ``truncation to $N_\Omega$'' (which can be considered as an operator on both $C(X)^J$ and $C(X_\Omega)^J$). We feel are free to identify functions on $N_\Omega/J$ as functions on both $X/J$ and $X_\Omega/J$, hence expressions of the form $\Vert \Phi - \tau_\Omega\iota_\Omega^*\Phi\Vert_{L^2(X)}$ will make sense
for $\Phi \in L^2(X)$.

\subsubsection{Bounding the $L^2(X)$-norm in terms of the asymptotics: the rank one case}

Let us first discuss a toy case, namely assume that $\XX$ is a spherical variety of rank one with $\mathcal Z(\XX)$ trivial. Fix an open compact $J\subset G$ and a $J$-good neighborhood $N_\Theta$ of infinity (there is only one nontrivial direction to infinity, which we will denote by $\Theta$), and denote by $\tau_\Theta$ the ``truncation'' to this neighborhood. Then we claim:
\begin{lemma} \label{toy-lemma}
 There is a finite set of elements $v_i\in L^2(X)^J$ and a constant $C$ such that, for $\Phi \in L^2(X)^J$ we have:
\begin{equation}
 \Vert \Phi - \tau_\Theta\iota_\Theta^* \Phi\Vert_{L^2(X)^J} \leq C \sum_i \Vert \Phi\Vert_{L^1_{v_i}}.
\end{equation}
\end{lemma}

Recall that the norms that appear in the final term have been defined in \S \ref{mixednormdefrankone}.

\begin{proof}
Let $\Phi_\pi$ denote the image of $\Phi\in C^{\infty}_c(X)$ in $\mathcal H_\pi$. Fix a Plancherel measure and recall that by $\Phi^\pi(x)$ we denote the pairing of $\Phi_\pi$ with the characteristic measure of $xJ$ with respect to the corresponding Plancherel form (see Remark \ref{decomp-into-eigen}). 

On $N_\Theta$ the difference $\Phi-\iota_\Theta^*\Phi$ can be expressed pointwise in terms of their spectral decomposition :   
$$ (\Phi-\iota_\Theta^*\Phi)(x) = \int_{\hat G} (\Phi^\pi(x)-(\iota_\Theta^*\Phi)^\pi(x)) \mu(\pi),$$ 
Note that, on the right hand side, we use the identification of $J$-orbits on $X$ and $X_{\Theta}$ to make sense of $\Phi^{\pi}(x)$
and $(\iota_{\Theta}^* \Phi)^{\pi}(x)$ simultaneously.  We allow ourselves to abbreviate $(\iota_{\Theta}^* \Phi)^{\pi}$
to $\iota_{\Theta}^* \Phi^{\pi}$ in what follows, to avoid a plethora of bracketing. 

Hence:
$$\Vert\Phi-\iota_\Theta^*\Phi\Vert_{L^2(N_\Theta)}^2 = \int_{N_\Theta} \left| \int_{\hat G} (\Phi^\pi(x)-\iota_\Theta^*\Phi^\pi(x)) \mu(\pi)\right|^2 dx.$$

By the asymptotics, $\Phi^\pi|_{N_\Theta}$ is equal to $e_\Theta^*\Phi^\pi|_{N_\Theta}$, and the difference $e_\Theta^*\Phi_\pi(x)-\iota_\Theta^*\Phi_\pi(x)$ can be expressed as a sum of (uniformly, as in Proposition \ref{uniformbound}) subunitary exponents -- see 
Remark \ref{DualBernsteinMap}. Therefore, by Lemma \ref{boundS}, there is a $L^2$ function (see below) $\Omega$ on $N_\Theta$ and a finite set of points $x_i$ on $N_\Theta$ such that for every $\pi, x$ we have: 
$$ |\Phi^\pi(x)-\iota_\Theta^*\Phi^\pi(x)|\le \Omega(x) \sum_i |\Phi^\pi(x_i)|.$$

Note that the important feature of $N_{\Theta}/J$ that is used here is 
that (considered as a subset of $X_{\Theta}/J$): it is $A_{X, \Theta}^+$-stable and consists of a finite number of $A_{X, \Theta}^+$-orbits. 
One  applies Lemma \ref{boundS} by pulling back to each copy of $A_{X, \Theta}^+$, i.e.
apply it to the function $a \cdot e_{\Theta}^* \Phi^{\pi} (x_0)$ for fixed $x_0 \in N_{\Theta}$. Note again that
the action of $a$ includes a twist 
  by the square root of the $A_{X, \Theta}$-eigenmeasure. In particular, the function $\Omega(x)$ actually lies in $L^2$, because 
its exponents in the $A_{X, \Theta}^+$-direction are subunitary for this twisted action of $A_{X, \Theta}$.

Therefore:  $$\int_{N_\Theta} \left| \int_{\hat G} (\Phi^\pi(x)-\iota_\Theta^*\Phi^\pi(x)) \mu(\pi)\right|^2 dx\le \Vert \Omega\Vert_{L^2(N_\Theta)}^2 \left(\int_{\hat G} \sum_i |\Phi^\pi(x_i)| \mu(\pi)\right)^2$$
and hence:
\begin{equation} \label{tbe} \Vert\Phi-\iota_\Theta^*\Phi\Vert_{L^2(N_\Theta)} \ll \sum_i \int_{\hat G}  |\Phi^\pi(x_i)|\mu(\pi). \end{equation} 

If we set $v_i=$the characteristic function of $x_i J$ then the integrals appearing on the right hand side are precisely the norms $\Vert \Phi\Vert_{L^1_{v_i}}$. We complement this set of $v_i$'s with the characteristic functions of the $J$-orbits on the complement of $N_\Theta$ (there are only finitely many such since we are in the rank one case with $\mathcal Z(X)=1$), and then the statement of the lemma is true.
\end{proof}

Now we allow $\mathcal Z(\XX)$ to be non-trivial, but keeping the rank of $X$ equal to one. We need to modify the statement of the lemma according to the morphism: $\hat G \to \widehat{\mathcal Z(G)^0}$. To simplify notation, since Plancherel measure is supported in the preimage of $\widehat{\mathcal Z(X)}\subset \widehat{\mathcal Z(G)^0}$, we feel free to write maps: $\hat G\to \widehat{\mathcal Z(X)}$, while we should be replacing $\hat G$ by the preimage of $\widehat{\mathcal Z(X)}$. As we have been doing until now, we fix a Haar measure on $\widehat{\mathcal Z(X)}$, which will be used as the Plancherel measure in the definition of the relative norms that appear in the following lemmas and Proposition \ref{propscatteringestimate};
this is allowable because clearly the Plancherel measure for $X$ as a $\mathcal{Z}(X)$-representation lies in the measure class of the Haar measure on $\widehat{\mathcal{Z}(X)}$.  (These norms were defined in \S \ref{mixednormdefrelative} and they now depend on the choice of Plancherel measure on $\widehat{\mathcal Z(X)}$.)

\begin{lemma} \label{estimaterankone}
 There is a finite set of elements $v_i\in L^2(X)^J$ such that on $L^2(X)^J$ we have: 
\begin{equation}
 \Vert \Phi - \tau_\Theta\iota_\Theta^* \Phi\Vert_{L^2(X)^J} \ll \sum_i \Vert \Phi\Vert_{  \mathcal Z(X),v_i}.
\end{equation}
\end{lemma}

\begin{proof}
 The proof is like before, except that now we will estimate the norm of $\Phi-\iota_\Theta^*\Phi$ on $N_\Theta$ by first integrating over the action of $\mathcal Z(X)$ and then applying the above arguments:
$$\Vert\Phi-\iota_\Theta^*\Phi\Vert_{L^2(N_\Theta)}^2= \int_{\widehat{\mathcal Z(X)}} \Vert\Phi_\chi - \iota_\Theta^*\Phi\Vert_{L^2(\mathcal Z(X)\backslash N_\Theta,\chi)}^2 d\chi.$$

Just as in \eqref{tbe}, there is a finite set of $J$-orbits on $N_\Theta$ (independent of $\chi$) such that $\Vert\Phi_\chi - \iota_\Theta^*\Phi\Vert_{L^2(\mathcal Z(X)\backslash N_\Theta,\chi)}$ is bounded by a constant   
 
 times:
$$ \sum_i \int_{\widehat{G}_\chi} |\Phi^\pi(x_i)| \mu_\chi(\pi),$$
where $\widehat{G}_{\chi}$ denotes the fiber of $\widehat{G}$ over $\chi$, and $\mu_\chi$ is the corresponding Plancherel measure on this fiber. Note that 
Lemma \ref{boundS} can be applied uniformly in $\chi$, because it depends only on an upper bound for subunitary exponents (the constant $c$ in that Lemma) and 
such a bound indeed follows from Proposition \ref{uniformbound}.

Therefore:  
\begin{eqnarray*}\Vert\Phi-\iota_\Theta^*\Phi\Vert_{L^2(N_\Theta)}^2 \ll \int_{\widehat{\mathcal Z(X)}} \left(\sum_i \int_{\widehat{G}_{\chi}} |\Phi_\pi(x_i)| \mu_\chi(\pi)\right)^2 d\chi \\ \ll \sum_i \int_{\widehat{\mathcal Z(X)}} \left( \int_{\hat{G}_\chi} |\Phi_\pi(x_i)| \mu_\chi(\pi)\right)^2 d\chi .\end{eqnarray*}

If we set $v_i$=the characteristic function of $x_i J$ then the sum on the right hand side is: $$\sum_i \Vert \Phi\Vert_{\widehat{\mathcal Z(X)},v_i}^2.$$
and we have thereby shown:
\begin{equation} \label{ontheway}\Vert \tau_\Theta \Phi - \tau_\Theta \iota_{\Theta}^* \Phi\|_{L^2(N_{\Theta})} \ll \sum_i \Vert \Phi\Vert_{\widehat{\mathcal Z(X)},v_i} \end{equation}
Complementing  this set of $v_i$'s by the characteristic functions of a finite set of $J$-orbits representing all $J$-orbits in $(X\smallsetminus N_\Theta)/\mathcal Z(X)$, we are done.
\end{proof}

\subsubsection{Bounding the $L^2(X)$-norm in terms of the asymptotics: the general case}
Let $\XX$ be a spherical variety of arbitrary rank now, $J\subset G$ an open compact subgroup.

 The analog of   \eqref{ontheway}   requires slightly more involved combinatorics because
of the possibility of exponents being unitary along a wall.

{ 
For (almost all) $\pi\in \hat G$ we have a decomposition of the asymptotics $e_{\Theta}^* \Phi^\pi$ in the $\Theta$-direction,  into (generalized) $A_{X,\Theta}$-eigencharacters, each of which is either unitary or $A_{X,\Theta}^+$-subunitary. To any such character $\chi$ we attach the subset $\Omega_\chi\supset \Theta$ corresponding to the largest ``face'' of $A_{X,\Theta}^+$ where the restriction of $\chi$ is unitary; for example, if $\chi$ is unitary then $\Omega_\chi=\Theta$, while if $\chi$ is strictly subunitary then $\Omega_\chi=\Delta_X$. (Notice that $\chi|_{\mathcal Z(X)}$ is necessarily unitary.) Accordingly, we have a decomposition:
$$ \Phi^\pi|_{N_\Theta} = \sum_{\Omega\supset\Theta} \Phi^{\pi,\Omega},$$
where $\Phi^{\pi,\Omega}$ contains the summands with generalized eigencharacter $\chi$ such that $\Omega_\chi=\Omega$, and integrating over $\hat G$: 
\begin{equation} \label{newdecomp} \Phi|_{N_\Theta} = \sum_{\Omega\supset\Theta} \Phi^\Omega.\end{equation} 
Moreover, we have 
\begin{equation} \label{newdecomp2}  (\Phi - \iota_{\Theta}^* \Phi)_{N_{\Theta}} = \sum_{\Omega \supsetneq \Theta} \Phi^{\Omega}. \end{equation}

Now let us extend these notions to all of $X$ at once:

 Note  (recall the definitions at the start of \S   \ref{ssestimates})  that $\tilde{N_{\Theta}} =  \bigcup_{\Omega \subset \Theta} N_{\Omega}$
and the quantity $\Phi^{\Theta}$ is defined on each $N_{\Omega}$. Thus we may regard $\Phi^{\Theta}$
as being defined on all of $\tilde N_{\Theta}$ and will extend it by zero off  $\tilde{N} _{\Theta}$.
With this convention, $\Phi = \sum_{\Omega} \Phi^{\Omega}$ on $X$,
and when we restrict to any $N_{\Theta}$ only those terms with $\Theta \subset \Omega$ are nonvanishing.

\begin{lemma}  \label{estimategenera}
 There are a finite number of elements $v_i\in L^2(N_\Theta)^J$ and a  decaying function $Q$ on $N_\Theta$ such that, for each $\Phi \in C^{\infty}_c(X)^J$, we have $$ |\Phi^{\Delta_X} (x)|^2 \le Q(x) \sum_i \Vert \Phi\Vert^2_{\widehat{\mathcal Z(X)}, v_i}$$
 for $x\in N_\Theta$. 
 
 In particular, for a suitable constant $C$ \begin{equation} \label{norm-estimate-new}  \|\Phi^{\Delta_X}\|_{L^2(X)} \leq  C \sum_i \Vert \Phi\Vert^2_{\widehat{\mathcal Z(X)}, v_i}\end{equation} 
 for each such $\Phi$. 
\end{lemma}  

Here the meaning of ``decaying'' for a function $Q$ on $N_{\Theta}$ is that
the associated function $a \cdot Q(x_0)$ is decaying
for $a \in A_{X, \Theta}^+$  and fixed $x_0 \in N_{\Theta}$; note that the action of $a$ is twisted,
as always, by the square root of the eigenmeasure. 

\begin{proof}
The second estimate \eqref{norm-estimate-new} follows immediately from the first:  Compute first the norm on $N_{\Theta}$, and note that
a decaying function is square-integrable;   then sum over $\Theta$.

The first  statement  will be proved as in Lemma \ref{toy-lemma} and Lemma \ref{estimaterankone}, but
we will replace Lemma \ref{boundS} by a slight modification. 

Let $J' \subset A_{X, \Theta}$ be an open compact subgroup acting trivially on $J$-invariant functions on $X_{\Theta}$. 
Let   $f$ be the function on $A_{X, \Theta}^+/J'$ whose value at $a$ is given by evaluating $a \cdot e_{\Theta}^*\Phi^{\pi}$ 
at any fixed $x_0 \in N_{\Theta}$.   Let $S = A_{X, \Theta}/J'$, $S^+$ the positive cone corresponding to $A_{X, \Theta}^+$.
Thus, $f$ is a finite function on $S$ whose degree is absolutely bounded above (see quoted statement from page \pageref{BZbound}). 
Let $f^{<}$ be the projection of $f$ onto all exponents $\chi$ for $S$ which are {\em strictly subunitary} on $S^+$ (see \S \ref{sssubunitary} for definition). 

What is needed to adapt the previous statements is precisely:

\begin{quote} Each point evaluation of $f^{<}$
is bounded by a finite sum of point evaluations of $f$,
and the constants appearing in this bound can be taken
to depend only on the degree $n$ of $f$ as an $S$-finite function.  
\end{quote}

However, the projection can be effected by a countour integral,
exactly as in in the proof of Lemma \ref{boundS}.  Let us go through that proof with $c=1$; 
as in that proof, we obtain a vector space $F(\mathbf{P})$ and an evaluation map $\mathrm{ev}: F(\mathbf{P}) \rightarrow \mathbb{C}^{\Lambda}$,
together with endomorphisms $A_1, \dots, A_t$ of $F(\mathbf{P})$, such that 
$$f(t_1, \dots, t_k) = \left( A_1^{t_1} \dots A_k^{t_k}  \mathrm{ev}(f) \right) (\mathbf{0}).$$

Now choose $\delta > 0$ so that any eigenvalue $\lambda_i$ of any $A_i$ that 
satisfies $|\lambda_i|  <1$ actually satisfies $|\lambda_i< 1- 2\delta$. In our context,
such a bound exists because of  Proposition \ref{uniformbound}, which shows there
are only a finite possible number of exponents which can appear for $J$-invariant functions. 
Now set $ P_i =  \int_{|z| = 1-\delta} \frac{dz}{z-A_i}  $;
it gives a projection onto all eigenvalues of $A_i$ which are less than $1$ in absolute value.
Then $P_1 P_2 \dots P_t$ furnishes the desired projection to subunitary exponents, 
and just as in the proof of Lemma \ref{boundS} its norm 
is absolutely bounded in terms of $n$.
\end{proof}

This implies a similar estimate for arbitrary  $\Omega$:
  \begin{equation} \label{FixedLemmaSecondBound} \|\Phi^{\Omega}\|_{L^2(X)} \ll \sum_i \Vert \iota_\Omega^* \Phi\Vert_{\widehat{\mathcal Z(X_\Omega)}, v_i}.\end{equation} 
 To see this, we just identify $\tilde{N}_{\Omega}$
to a subset of $X_{\Omega}$ and therefore identify $\Phi^{\Omega}$ to a function on $X_{\Omega}$. 
This function is precisely $(\iota_{\Omega}^* \Phi)^{\Omega}$, restricted to $\tilde{N}_{\Omega}$.

Recall our notational convention about relative norms, defined   before Lemma \ref{Sumofnorms}.
For every $v\in L^2(X_\Omega)$ we can consider the ``relative'' norm:
$$\Vert\bullet\Vert_{   \mathcal Z(X_\Omega), v}$$
on $L^2(X_\Omega)$.  
In what follows, we will use these norms for $v$=the characteristic function of some $J$-orbit in $N_{\Omega'}$, $\Omega'\subset\Omega$, hence identified with a function on $X_\Omega$.

 \begin{proposition} \label{propscatteringestimate}
Let $\Phi\in L^2(X)_r^J$ (the image of all $L^2(X_\Theta)_\disc^J$ with $|\Theta|=r$).  Then 
there is a constant $C$ and a  finite list of vectors $v_{\Omega, i} \in L^2(X_{\Omega})$ such that 
\begin{equation}
 \Vert \Phi \Vert \le   \Vert  \sum_{|\Omega|=r}  \tilde{\tau_{\Omega}}\iota_\Omega^*\Phi\Vert_{L^2(\tilde N_\Omega)^J}   + C  \sum_{|\Omega|>r}  \sum_i \Vert \iota_\Omega^*\Phi\Vert_{ \mathcal Z(X_\Omega),v_{\Omega,i}}
\end{equation}
where the meaning of $\tilde{\tau}_{\Omega}$ is ``restriction to $\tilde N_{\Omega}$ and consider $\tilde{N}_{\Omega}$ as a subset of $X$.'' (Thus, $\tilde{\tau}_{\Omega} \iota_{\Omega}^* \Phi$
is a function on $X$, and the $\Omega$-sum is taken on $X$.)
\end{proposition}

\proof

In particular, note that $\Phi^{\Omega} = 0$ if $|\Omega| <  r$, because $\iota_{\Omega}^* \Phi = 0$. 
By the triangle inequality 
$$\| \Phi \|_{L^2(X)} \leq \|  \sum_{\Omega=r} \Phi^{\Omega}\|_{L^2(X)}  + \sum_{\Omega > r} \|\Phi^{\Omega}\|_{L^2(X)}.$$

and apply \eqref{FixedLemmaSecondBound} to all terms with $|\Omega| > r$.  
And (tracking through the notation) $\sum_{|\Omega|=r} \Phi^{\Omega}$ is the same
as $\sum \tilde{\tau}_{\Omega} \iota_{\Omega}^* \Phi$. 

\qed 
  }

\subsection{Proof of the second part of Theorem \ref{tilingtheorem}} \label{tilingtheoremIIproof}

{\bluetext Again, for notational simplicity, we denote $A_{X,\Theta}'$, $A_{X,\Omega}'$ simply by $A_{X,\Theta}$, $A_{X,\Omega}$.}

Fix a $\Theta\subset \Delta_X$.  
We may partition $\widehat {\mathcal Z(X_{\Theta})}$   into an almost disjoint (i.e.\ disjoint up to a set of measure zero) union of measurable subsets $Y_{\beta}$ with the following properties:
\begin{itemize}
\item[(A)]
 If $(T_i,f_i)$ are the isogenies  and affine maps of the decomposition \eqref{SsumSi} for the map $\iota_\Theta^*\circ \iota_\Theta$, then $Y_{\beta} \cap  f_i( Y_{\beta})$ is of measure zero
 { unless $T_i$ is the identity and $f_i$ is also the identity.} 
 \item[(B)] Each $Y_{\beta}$ is a subset of a single connected component of $\widehat{A_{X,\Theta}}$. 
 \end{itemize} 
 {To carry this out, we proceed one connected component at a time.  First  note that  the set of $\chi \in \widehat{A_{X,\Theta}}$ in this component
 fixed by some $f_i$ is a closed set of measure $0$ -- it is, in suitable coordinates, a union of translates of sub-tori. Now choose a countable dense set $\{P_i\}$ in the complement of the fixed locus,
 and for each  $P_i$ let $B_i$ be an open  ball around $P_i$ that is disjoint from each $f_i(B_i)$. Now take for our partition $B_1, B_2-B_1, B_3-(B_1 \cup B_2)$ and so on. In our application,
 where the isogenies arise from a Weyl group, one can easily in fact give an explicit choice of $Y_{\beta}$ using positive chambers. }

Applying the corresponding idempotents $1_{Y_\beta}$ gives rise to a direct sum decomposition of $\mathcal{H}$ into $A_{X,\Theta}\times G$-stable subspaces for which the assumptions of Proposition \ref{Bmapisometry}  -- and, indeed,  the stronger assumptions
of Remark \ref{UniqueExponentSituation} -- are satisfied.
{ Namely,   denote by $\mathcal{H}'$ any one of the resulting $A_{X, \Theta} \times G$-stable summands of $\mathcal{H}$, corresponding to the image of $1_{Y_{\beta}}$. 
Decompose $\mathcal{H}' = \int_{\pi} \mathcal{H}'_{\pi} d\pi$ as $G$-representation. We need to verify that for almost all $\pi$ appearing in this decomposition, 
there is a unique $A_{X, \Theta}$-exponent on $\mathcal{H}'_{\pi}$. After all,
choose one such exponent $\chi$, belonging to $Y_{\beta}$ say;
then Proposition \ref{e2c} shows that any other exponent must be of the form $f_i(\chi)$ for some $(T_i, f_i)$,
but we know by virtue of assumption (A) above that $f_{i}(\chi) \notin Y_{\beta}$ unless $f_{i}(\chi) = \chi$. So $\chi$ is actually the unique $A_{X, \Theta}$-exponent of $\mathcal{H}_{\pi}$. }

Hence $\iota_\Theta$ is an isometry onto the image when restricted to each of these summands.

Replacing $\mathcal{H}$ by one of these summands, we may suppose that $\iota_\Theta$ is an isometry on our given subspace $\mathcal H\subset L^2(X_\Theta)$. For any $a\in A_{X,\Theta}$,  and $v \in L^2(X_{\Theta})^J_{ \disc}$ -- for some fixed open compact subgroup $J$ --  Proposition \ref{propscatteringestimate} applied to the function $\Phi=\iota_\Theta \mathcal L_{a^n} v$ yields: 
 $$\|v\| =  \Vert \mathcal L_{a^n} v\Vert =   \Vert \iota_\Theta \mathcal L_{a^n} v\Vert   $$
$$\le \left( \Vert \sum_{|\Omega|=r} \tilde{ \tau}_\Omega \iota_\Omega^*\iota_\Theta \mathcal L_{a^n} v\Vert + C\cdot \sum_{|\Omega|>r} \sum_j \Vert \iota_\Omega^*\iota_{\Theta} \mathcal{L}_{a^n} v\Vert_{ \mathcal Z(X_\Omega),v_{\Omega,j}}\right),$$
where $C$ denotes the implicit constant of Proposition \ref{propscatteringestimate}.   The relative norms that appear on the right hand side will depend, of course, on $J$.  

Decomposing $\iota_\Omega^*\iota_\Theta$, for $|\Omega|=r$, according to (\ref{SsumSi}), we see that: $$\lim_{n \rightarrow \infty}\Vert \tilde{\tau}_\Omega S_i \mathcal L_{a^n} v\Vert =0$$ 
unless  $\val(a) \in T_i \mathfrak a_{X,\Omega}^+$.\footnote{Suppose -- for simplicity, the modifications in general are not difficult -- that the image of $a$ in 
$A_{X, \Theta}/A_{X, \Theta}^0$ coincides with $T b$, where $b \in A_{X, \Omega}^T$
but $b \notin A_{X, \Omega}^+$. 
By \eqref{Siequivariance} it is enough to show that
$$ \Vert \tilde{\tau}_{\Omega}^* \mathcal{L}_{b_i^n} S_i f \Vert \rightarrow 0$$
when $b \notin A_{X, \Omega}^+$;   in coordinates this amounts to the following: { Suppose that $n=n_1+n_2$. Then, for $g \in L^2(\Z^n)$ and $b = (b_1, \dots, b_n)$  
where some $b_i < 0$ for $1 \leq i \leq n_1$, we have 
then $$ \| \mbox{translate of $g$ by $m \cdot b$} \|_{L^2(\N^{n_1} \times \Z^{n_2})} \stackrel{m\rightarrow \infty}{\longrightarrow} 0.$$ }}
Here we denoted by $\val$ the natural ``valuation'' map: $$A_{X,\Theta} \to \varchi(\AA_{X,\Theta})^*\subset \mathfrak a_{X,\Theta}.$$
Hence we get:
\begin{equation}\label{inequality}\|v\| \le \liminf_{n \rightarrow \infty}\left(  \Vert \sum_{|\Omega|=r,  i :  \,   T_i {\mathfrak a}_{X,\Omega}^+ \ni \val(a)} \tilde{\tau}_\Omega S_i \mathcal L_{a^n} v\Vert + C\cdot \sum_{|\Omega|>r} \sum_j \Vert \iota_\Omega^*\iota_{\Theta} \mathcal{L}_{a^n} v\Vert_{ \mathcal Z(X_\Omega),v_{\Omega,j}}\right). 
\end{equation} 
{ We emphasize that the summations over $i$ and $j$ are {\em finite}, because the sums in both Proposition \ref{propscatteringestimate} and 
 \eqref{SsumSi} are finite.}

Now assume $a$ to be ``generic'' in the sense that $\val(a)\in T_i \mathfrak a_{X,\Omega}^+$ implies $\val(a)\in T_i\mathring{\mathfrak a}_{X,\Omega}^+$.
As in the proof of the first part of Theorem \ref{tilingtheorem} (\S \ref{prooffirstpart}),  
it is easy to see that the different summands $\tilde\tau_\Omega S_i \mathcal L_{a^n} v$ become orthogonal in the limit $n\to \infty$, i.e.\ for any given $\epsilon>0$ there is an $N>0$ such that for all $n\ge N$ and any distinct indices $i, j$ we have:
$$|\left<\tilde\tau_{\Omega_i} S_i \mathcal L_{a^n} v,\tilde\tau_{\Omega_j} S_j \mathcal L_{a^n} v\right>| <\epsilon.$$
{   Indeed, let $(\Omega_i,T_i,f_i)$ (with $|\Omega_i|=r$) be the triples corresponding to the morphisms $S_i$ appearing in the first sum of \eqref{inequality}. In what follows, we will say that ``almost all the mass of a function $f$ is concentrated in a set $M$'' if the square of the $L^2$-norm of $f$ restricted to the complement of $M$ is less than a certain multiple of $\epsilon$. (We omit the straightforward task of specifying which multiple is needed.)

For two distinct indices $i$, $j$ we have the following possibilities: \begin{enumerate}
 \item $\Omega_i\ne \Omega_j$. In this case, let $K_i\subset \infty_{\Omega_i}$, $K_j\subset \infty_{\Omega_j}$ be compact subsets such that  almost all the mass of $S_i v$ is concentrated in the preimage of  $K_i$ under the quotient map $X_{\Omega_i}\to \infty_{\Omega_i}$ (and similarly for $K_j$). Recall that $\infty_\Omega$ denotes ``$\Omega$-infinity'', i.e.\ the union of orbits in a toroidal compactification which correspond to $\Omega$. 
 There are neighborhoods of $K_i$, $K_j$ in $X$ which are disjoint, and for $n$ large enough almost all the mass of $S_i \mathcal L_{a^n} v$, $S_j \mathcal L_{a^n} v$ is concentrated in the respective neighborhoods.
 
 \item $\Omega_i=\Omega_j=\Omega$ but $T_i\ne T_j$. By choosing a section of the quotient map $X_\Omega\to \infty_\Omega$, we can find a compact subset $M$ of the image of that section and a compact subset $N$ of $A_{X,\Omega}$ such that almost all the mass of $S_i v$ and $S_j v$ is concentrated on $N\cdot M$. Then, for $n$ large enough\footnote{Here we are slightly abusing notation and treating the isogenies $T_i^{-1}$, $T_j^{-1}$ as actual morphisms. In reality, $T_i^{-1}a^n$ is well-defined only for $n$ in a finite-index subgroup of $\mathbb Z$. To be rigorous, enlarge the set $N$ so that $N\cdot M$ also contains most of the mass of $S_i \mathcal L_{a^n} v$, $S_j \mathcal L_{a^n} v$ for $n$ in a set of representatives of the cosets of this subgroup.} the sets $T_i^{-1}a^n N$ and $T_j^{-1} a^n N$ are disjoint; therefore, almost all of the mass of $S_i \mathcal L_{a^n} v$, $S_j \mathcal L_{a^n} v$ is concentrated on the disjoint sets $T_i^{-1}a^n N\cdot M$, $T_j^{-1}a^n N\cdot M$, respectively.
 
 \item $\Omega_i=\Omega_j=\Omega$ and $T_i=T_j$ but $f_i\ne f_j$. Then, by property (A), the $A_{X,\Omega}$-Plancherel supports of $S_i \mathcal L_{a^n} v$ and $S_j \mathcal L_{a^n} v$ intersect at a set of Plancherel measure zero, and hence these vectors are orthogonal for every $n$. 
\end{enumerate}

We have thus established the ``orthogonality in the limit''. 
Notice also that, in the limit, all the mass of $S_i   \mathcal L_{a^n} v$ is concentrated in $\tilde N_{\Omega_i}$, so we can get rid of the restriction operators $\tilde \tau_\Omega$.
Moreover, using the equivariance property \eqref{Siequivariance} (by replacing, if necessary, the element $a$ by a suitable power $a^k$ so that it is of the form $T_i(a')$ in the notation of \eqref{Siequivariance}, for all $i$), 
and the fact that the $L^2(X_\Omega)$-norm is $A_{X,\Omega}$-invariant we can get rid of $\mathcal{L}_{a^n}$. }
Therefore, we get:

\begin{equation}\label{inequality2}\Vert v\Vert  \leq \left(  \sum_{|\Omega|=r , \, i: \ T_i \mathring{\mathfrak a}_{X,\Omega}^+ \ni \val(a)} \Vert S_i   v\Vert^2 \right)^\frac{1}{2} + O \left( \liminf_{n \rightarrow \infty}  \sum_{|\Omega|>r} \sum_i \Vert \iota_\Omega^*\iota_{\Theta} \mathcal{L}_{a^n} v\Vert_{ \mathcal Z(X_\Omega),v_{\Omega,i}}\right).
\end{equation}

To complete the proof, consider the $G \times A_{X, \Theta}$-invariant Hermitian norm:  
 $$ H(v) := \sum_{|\Omega|=r , \,  i : \ T_i \mathring{\mathfrak a}_{X,\Omega}^+ \ni \val(a)}  \Vert  S_i v \Vert^2.$$
  It is $A_{X, \Theta} \times G$-invariant, and also absolutely continuous
 with respect to $\| \cdot \|_{L^2}$, since the $S_i$ are bounded morphisms.  By spectral theory, if $H(v) < \|v\|$
 for some nonzero $v$, we may find an $A_{X, \Theta} \times G$-invariant
 space $\mathcal{H}' \subset \mathcal{H}$ and $\delta>0$ with the property that $H(v) < (1-\delta) \|v\|^2$
 on $\mathcal{H}'$.   

Thus on $\mathcal{H}'^J$ the Hilbert norm $\Vert v\Vert$ is bounded by a finite sum of relative norms:
$$\|v\| \leq C' \liminf_{n \rightarrow \infty} \sum_{|\Omega|>r} \sum_i \Vert \iota_\Omega^*\iota_{\Theta} \mathcal{L}_{a^n} v\Vert_{ \mathcal Z(X_\Omega),v_{\Omega,i}}, \ \ v \in {\mathcal{H}'}^J.$$
Decomposing $\iota_\Omega^*\iota_{\Theta}$ with respect to Proposition \ref{decompf}, and denoting the summands by $S_{\Omega,j}$, we get: 
$$\|v\| \leq C' \liminf_{n \rightarrow \infty} \sum_{|\Omega|>r} \sum_{i,j} \Vert S_{\Omega,j} \mathcal{L}_{a^n} v\Vert_{ \mathcal Z(X_\Omega),v_{\Omega,i}}, \ \ v \in {\mathcal{H}'}^J.$$
By the equivariance property\footnote{The $S_{\Omega,j}$ are not quite equivariant with respect to a morphism $\AA_{X,\Omega}\to \AA_{X,\Theta}$; therefore, we are really applying variants of Lemmas \ref{injectivemeasure} and \ref{norm contra} which apply to the affine maps of character groups introduced in Section \ref{sec:linalg2}; we leave the details to the reader.}  of $S_{\Omega,j}$, and Lemma \ref{injectivemeasure} (which we can apply because, by choice, we are in the ``unique exponent'' situation of Remark \ref{UniqueExponentSituation})  we get a contradiction to Lemma \ref{norm contra}.    \qed

\section{Explicit Plancherel formula} \label{sec:explicit}

\emph{In this section we assume that $\XX$ is strongly factorizable, cf.\ \S\ref{factorizable}, for example: a symmetric variety.} 

Since the formalism of this section may appear quite involved, we begin by a rough description of its thrust:

We wish to write a formula for the ``smooth'' and ``unitary'' asymptotics maps: $e_\Theta$ and 
$\iota_{\Theta}$. As we have explained in \S \ref{sshorocycles}, {\em the varieties $\XX$ and $\XX_{\Theta}$ look quite different,
but there is a canonical identification of their varieties of horocycles $\XX_\Theta^h$. }
Our assertion is that the maps $e_\Theta$ and $\iota_{\Theta}$ can be obtained by a suitable ``interpretation'' of the diagram
$$\mbox{functions on $X$} \rightarrow \mbox{functions on $X_\Theta^h$} \leftarrow \mbox{functions on $X_{\Theta}$.}$$
where the arrows are obtained by Radon transform.   (By ``interpretation'', we mean, roughly speaking, disintegrating the arrows spectrally and making sense of convergence issues.)

If one carries out the analogs of our constructions 
in the case where $\XX$ is a real symmetric space,
we arrive at the theory of ``Eisenstein integrals''  developed by van den Ban, Schlichtkrull, and Delorme \cite{BSP1, BSP2, De}.  The relationship between $\XX, \XX_{\Theta}, \XX_\Theta^h$ helps to give a geometric
interpretation of this theory and, in particular, the correct normalization of Eisenstein integrals.

{\em A notational convention:} We have tried (possibly foolishly) to avoid choosing a measure on unipotent radicals in this section. To this end we introduce the following notation: In various contexts we shall denote by $V'$ a space
that is non-canonically isomorphic to a space $V$; but the isomorphism depends on the choice
of a measure on a certain unipotent group.  To a first approximation when reading, then, 
the reader can simply ignore the primes and replace each $V'$ by $V$.

\subsection{Goals}

Recall that Radon transform was defined in \ref{sssRadon} as a canonical morphism:
\begin{equation} C^{\infty}_c(X) \stackrel{R_\Theta}{\longrightarrow} C^{\infty}(X_\Theta^h, \delta_\Theta),\end{equation}
 where by $C^{\infty}(X_\Theta^h, \delta_\Theta)$ we denote sections of the dual line bundle to the line bundle of Haar measures on the unipotent radical of the parabolic corresponding to each point. Recall that $C^{\infty}(X_\Theta^h, \delta_\Theta)$ denotes smooth sections of a line bundle over $X_\Theta^h$ where the stabilizer of each point on $X_\Theta^h$ (contained in a parabolic of type $P_\Theta$) acts on the fiber via the modular character $\delta_\Theta$, hence the notation. The action of $G$ is twisted by the same character on both $C^\infty(X)$ and $C^\infty(X_\Theta^h,\delta_\Theta)$, in order to make this map equivariant.
 We will denote by $C^{\infty}(X^h_{\Theta}, \delta_{\Theta})_X$
 the Radon transform of the space $C_c^\infty(X)$. 

Our starting point for the explicit identification of smooth asymptotics is Proposition \ref{propositionRadon}, which states that the adjoint asymptotics maps $e_\Theta^*$ commutes with Radon transform, i.e.\ the following diagram commutes:

\begin{equation}\label{Radoncommutes2}
\xymatrix{
 C_c^\infty(X) \ar[dd]_{e_\Theta^*} \ar[dr]^{R_\Theta}  \\
& C^{\infty}(X^h_{\Theta}, \delta_\Theta)_X \\
 C_R^\infty(X_\Theta) \ar[ur]_{R_\Theta}
}
\end{equation}
where we have used $C^{\infty}_R(X_{\Theta})$ as an {\em ad hoc} notation for those elements
of $C^{\infty}(X_{\Theta})$ for which the transform defining $R_{\Theta}$ is absolutely convergent. 

This diagram suggests that by  ``inverting'' the lower occurence of Radon transform we can obtain an explicit formula for the smooth asymptotics; we do this in \S \ref{ssexplicitsmooth}, where we express $e_\Theta f$, for every $f\in C_c^\infty(X_\Theta)$, as a ``shifted wave packet of normalized Eisenstein integrals''.

We then proceed to do the same for unitary asymptotics (the Bernstein maps). In that case, we will ``filter out'' non-unitary exponents from the inversion of Radon transform, and the expansion of $\iota_\Theta f$ will be as a ``wave packet'' over a set of unitary representations of $L_\Theta$. Similar as they may be, the two goals are independent, and the main theorems in the smooth and unitary case do not rely on each other. It should be possible, and would be very interesting, to obtain the expression for $\iota_\Theta$ by shifting the contour in the expression of $e_\Theta$, but this is not the direction that we pursue here. For this reason, however, the appearance of Eisenstein integrals in the expression for $\iota_\Theta$ is conditional on some weak ``multiplicity one'' assumption which we are able to verify in many cases (including symmetric varieties).

Recall from \S \ref{decomp-into-eigen} that, fixing a Plancherel decomposition: $$L^2(X) = \int_{\hat G} \mathcal H_\pi\mu(\pi)$$ we get a pointwise decomposition of (smooth) functions: $$\Phi(x) = \int_{\hat G} \Phi^\pi(x) \mu(\pi).$$ Of course, it is  more intrinsic to think of $\Phi^\pi \mu(\pi)$, which is a measure on $\hat G$ valued in functions on $X$, instead of $\Phi^\pi$, because that doesn't depend on the choice of Plancherel measure.

Our main goal here is to obtain such a decomposition for $\iota_\Theta f$ (where $f\in L^2(X_\Theta)^\infty$), but starting from a Plancherel decomposition for $L^2(X_\Theta)$:

Fix a Plancherel measure $\nu$ for $L^2(X^L_{\Theta})$; as usual, $X_\Theta^L$ denotes the Levi variety. We will recall the necessary facts about Levi varieties in a moment, including a normalization of their measures such that: $$L^2(X_\Theta)_{\disc} = I_{\Theta^-} (L^2(X_\Theta^L)_{\disc})$$ (unitary induction). We have a Plancherel decomposition:
$$L^2(X^L_{\Theta})= \int \mathcal{I}_{\sigma} \nu(\sigma).$$
By induction:
\begin{equation}\label{Planchereltheta}
 L^2(X_\Theta) = \int_{\widehat{L_\Theta}} \mathcal H_\sigma \nu(\sigma).
\end{equation}
where $\mathcal{H}_{\sigma} = I_{\Theta^-} \mathcal{I}_{\sigma}$ is the parabolic induction of $\mathcal{I}_{\sigma}$
to $G$.

Our goal is to obtain a formula:
\begin{equation} \label{iotathetaexplicit}
 \iota_\Theta f(x) = \int_{\widehat{L_\Theta}} (\iota_\Theta^\sigma f)(x) \nu(\sigma),
\end{equation}
for some explicit $\iota_\Theta^\sigma: \mathcal H_\sigma^\infty \to C^\infty(X)$.The abstract existence of such a formula is a tautology: setting $\tilde{1}_{xJ } :=  \Vol(xJ)^{-1} 1_{xJ}$,
one sees formally that
$\iota_{\Theta} f(x) = \left<f,\iota_\Theta^* \tilde{1}_{xJ}\right> = \int_{\widehat{L_\Theta}} \left< f,\iota_\Theta^*\tilde{1}_{xJ}\right>_\sigma \nu(\sigma),$
and hence:
\begin{equation}\label{defNsigma}
 \iota_\Theta^\sigma (f)(x):= \left< f,\iota_\Theta^* \tilde{1}_{xJ} \right>_\sigma,
\end{equation}
which clearly factors through $\mathcal H_\sigma$. As the above formula shows, $ \iota_\Theta^\sigma$ is adjoint to the map:
\begin{equation} \label{iota*def} \iota^*_{\Theta,\sigma}: C_c^\infty(X) \to \mathcal H_\sigma\end{equation}
decomposing $\iota_\Theta^*$.

We are able to identify some invariant of the morphisms $\iota_{\Theta, \sigma}^*$, the so-called \emph{Mackey restriction} (see Theorem \ref{explicitMackey}), and under an additional assumption which is known to be true for symmetric varieties
through the work of Blanc-Delorme \cite{BlancDelorme} this is enough to identify them, again, with normalized Eisenstein integrals.

Finally, the expression (\ref{iotathetaexplicit}) can be reinterpreted as a precise and explicit Plancherel formula for $\iota_\Theta L^2(X_\Theta)$, to the extent that the full Scattering Theorem \ref{advancedscattering} is known; see Theorem \ref{explicitPlancherel}.

\subsection{Various spaces of coinvariants} \label{ss:explicitplanchereldefns}   
In this section we shall, roughly speaking, decompose
the Radon transform $C^{\infty}_c(X_{\Theta}) \rightarrow C^{\infty}(X^h_{\Theta},\delta_\Theta)$
over representations of $L_{\Theta}$. More precisely:

\subsubsection{Goal}

Given irreducible unitarizable representation $\sigma$ of $L_{\Theta}$ 
we shall define  two spaces and a morphism between them: 
\begin{equation} \label{twospaces} RT_{\Theta}: C^{\infty}_c(X_{\Theta})_{\sigma} \longrightarrow  C^{\infty}(X^h_{\Theta}, \delta_{\Theta})_{X, \sigma} \end{equation}

The spaces $C^{\infty}_c(X_{\Theta})_{\sigma} $ and $C^{\infty}(X^h_{\Theta}, \delta_{\Theta})_{X, \sigma}$
are certain $I_{\Theta^-}(\sigma)$-isotypical quotients of $C^{\infty}_c(X_{\Theta})$
and $C^{\infty}(X^h_{\Theta}, \delta_{\Theta})_X$. (Recall that $C^{\infty}(X^h_{\Theta}, \delta_{\Theta})_X$  is the image of $C_c^\infty(X)$ under Radon transform.)

The definition of these spaces makes sense only for $\nu$-almost every $\sigma$. At
several points we will abuse notation by omitting the phrase ``for $\nu$-almost every $\sigma$.'' 

  As for the morphism
between these spaces, it is a version of the Radon transform
but is also related to the standard intertwining operator: 
\begin{equation}\label{intops}T_\Theta: I_\Theta(\sigma^-) \to \int_{U_\Theta} f(u \bullet) du \in I_\Theta(\sigma)'.
\end{equation}
In fact, there will be
a commutative diagram     \begin{equation} \label{green}
 \begin{CD}
C^{\infty}_c(X_{\Theta}) @>{R_{\Theta}}>> C^\infty(X_\Theta^h, \delta_\Theta)_X \\
@VVV @VVV \\
C^{\infty}_c(X_{\Theta})_{\sigma} @>{RT_{\Theta}}>> C^{\infty}(X^h_{\Theta}, \delta_{\Theta})_{X, \sigma}\\
@V{\sim}VV @V{\sim}VV \\
\left({\Hom_{L_\Theta}(C^{\infty}_c(X_\Theta^L),\sigma )}\right)^*\otimes I_{\Theta^-}(\sigma)
@>{\mathrm{id} \otimes T_{\Theta}}>> \left({\Hom_{L_\Theta}(C^{\infty}_c(X_\Theta^L),\sigma )}\right)^*\otimes I_{\Theta}(\sigma)'
 \end{CD}
\end{equation}
   
The space $I_{\Theta}(\sigma)'$ is isomorphic to $I_{\Theta}(\sigma)$ once one fixes a measure
   on $U_{\Theta}$, and is defined in  \ref{sss:twistingclass}.
In any case,
   we will denote the morphism as $RT_{\Theta}$ because it can be thought
   of either as Radon transform or standard intertwining operator.

\subsubsection{Parabolics and Levi subgroups}

In this section we will not fix parabolics in the classes of $\PP_\Theta$, $\PP_\Theta^-$, except when needed. Hence, we will be thinking of $\LL_\Theta$ as a ``universal Levi group of type $\Theta$'', that is, not a subgroup of $\GG$, but rather the reductive quotient of any parabolic in the class of $\PP_\Theta$, which is a canonical abstract group up to inner automorphisms. We can also define it as the reductive quotient of any parabolic in the class of $\PP_\Theta^-$, and the two definitions give canonically isomorphic groups, up to inner conjugacy, by identifying the reductive quotients with the intersection $\PP_\Theta\cap \PP_\Theta^-$. A ``representation of $L_\Theta$'' will actually be only an isomorphism class of representations, but when we fix parabolics $P_\Theta$ and $P_\Theta^-$ then we will implicitly be fixing (compatible) identifications of their reductive quotients with $L_\Theta$, and a realization for the representations of $L_\Theta$ under consideration.

Similarly, we do not fix a ``Levi variety'' $\XX_\Theta^L$ as a subvariety of $\XX_\Theta$, except when explicitly saying so. A choice of parabolic in the class of $\PP_\Theta^-$ uniquely identifies such a Levi subvariety of $\XX_\Theta$ (consisting of all points whose stabilizers contain the unipotent radical of that parabolic), but in general we are only interested in the isomorphism class of $\XX_\Theta^L$ as a homogeneous $\LL_\Theta$-space.

We let $\LL_{\Theta,X}^\ab$ be the torus quotient of $\LL_\Theta$ whose character group consists of all characters of $\LL_\Theta$ which are trivial on the stabilizers of points on $\XX_\Theta^L$. By our assumption that $\XX$ is strongly factorizable, the rank of $\LL_{\Theta,X}^\ab$ is equal to the rank of $\AA_{X,\Theta}$.

\subsubsection{Twisting class} \label{sss:twistingclass} 
The class of representations of $L_\Theta$ obtained by twisting an irreducible  unitary representation  $\sigma$ by all unramified characters of $L_{\Theta,X}^\ab$ will be called, for short, a {\em twisting class.} We will say ``for almost every $\sigma$'' in a twisting class for statements that hold in an Zariski open and dense set of elements of a twisting class; notice the canonical algebraic structure, coming from the torus structure of the set of unramified characters of $L_{\Theta,X}^\ab$.

\subsubsection{Twists and half-twists; normalized and unnormalized induction; intertwiners} \label{AnalPedantry}
We use the symbol $I^{\un}$ to denote unnormalized induction.  It will be useful to talk about unnormalized induction at first,  because some of the algebraic structures are made clearer.

Recall the definition of the space $C^\infty(X_\Theta^h,\delta_\Theta)$ from \S \ref{sssRadon}. First $\delta_\Theta$ denotes the modular character of $P_\Theta$, which is inverse to the modular character of $P_\Theta^-$; for example, the usual induction $I_{\Theta}(\sigma)$  is $I_{\Theta}^{\un}(\sigma \delta_{\Theta}^\frac{1}{2})$. 
  
Secondly, $C^\infty(P_\Theta\backslash G,\delta_{\Theta})$ denotes smooth sections of the complex line bundle over $P_\Theta\backslash G$ whose fiber at a point is dual to the space of Haar measures
on the unipotent radical of the corresponding parabolic. The space of sections of this line bundle is non-canonically isomorphic to the representation parabolically induced(unnormalized) from $\delta_{\Theta}$.  Similarly, $C^\infty(X_\Theta^h,\delta_\Theta)$ denotes sections of the pull-back line bundle on $X_\Theta^h$, and $C^\infty(X_\Theta^L,\delta_\Theta)$ sections of the restriction of this line bundle to $X_\Theta^L\subset X_\Theta^h$ (see, however, the discussion below on normalization of the action of $L_\Theta$).

With these notations, for instance, the intertwining operator is {\em canonically}:
$$I_{\Theta^-}^{\un}(\sigma) \longrightarrow I_{\Theta}^{\un}(\sigma, \delta_{\Theta})$$
where the right hand side is defined as follows:  
Tensor the $G$-linear vector bundle over $P_\Theta\backslash G$ corresponding to the induced representation $I_\Theta^{\un}(\sigma)$ by 
the line bundle $\delta_{\Theta}$.   Smooth sections of this line bundle can be denoted by $I_\Theta(\sigma, \delta_\Theta)$.  For the normalized induction the corresponding morphism is:
\begin{equation} \label{somenotn} I_{\Theta^-}(\sigma) \longrightarrow I_{\Theta}^\un(\sigma \delta_{\Theta}^{-\frac{1}{2}}, \delta_{\Theta}) \end{equation} 
and we abridge the right-hand side to $I_{\Theta}(\sigma)'$; it is isomorphic to
$I_{\Theta}(\sigma)$, but it is sometimes clearer 
not to make this explicit.   Similarly the Radon transform is canonically:
 $$C^{\infty}_c(X) \longrightarrow C^{\infty}(X_\Theta^h, \delta_{\Theta}).$$

 We twist the action of $L_\Theta$ on functions on $X_\Theta^L$ (or sections in $C^\infty(X_\Theta^L,\delta_\Theta)$) in such a way that $L^2(X_\Theta)$ is the normalized induction of $L^2(X_\Theta^L)$
 (or, for that matter, $C^{\infty}_c(X_{\Theta})$
is the normalized induction of $C^{\infty}_c(X^L_{\Theta})$). The explicit formula was given in
\eqref{Ltwist}.

 \textbf{We caution the reader that this may not be the most natural-looking action}; for instance, if $X$ has a $G$-invariant measure and we consider the Levi variety $\XX_\emptyset^L\simeq \AA_X$, the usual action of $A$ on $C^\infty(A_X)$ is twisted by the square root of the modular character of $P(X)$.

Recall that the {\em normalized} Jacquet module (which we have been using, by convention, throughout this paper) 
 twists the action of $L_{\Theta}$ by $\delta_{\Theta}^{-1/2}$.
Therefore, with the definitions above we get a natural inclusion:
\begin{equation} \label{Jcareful} C^{\infty}_c(X^L_{\Theta}, \delta_{\Theta}) \otimes \delta_{\Theta}^{-1} \hookrightarrow
  C^{\infty}_c(X)_{\Theta}.\end{equation} 
We explicate the map from $ C^{\infty}_c(X^L_{\Theta}, \delta_{\Theta})$
to $C^{\infty}_c(X)_{\Theta}$, and then the factor $\delta_{\Theta}^{-1}$ comes by checking how $L_{\Theta}$ acts on both sides: 
An element $\nu \in C^{\infty}_c(X^L_{\Theta}, \delta_{\Theta})$ assigns
  to each point of $X^L_{\Theta}$  an element $\nu_x$ in the dual of the space of Haar measures on $U_{\Theta}$.
  We send $\nu$  to any function in $C^{\infty}_c(\mathring{X} P_{\Theta})$ whose 
  integral over $x U_{\Theta}$ against the Haar measure $du$ coincides with $\langle \nu_x, du \rangle$.

  {\em We  abridge the left-hand side to $C^{\infty}_c(X^L_{\Theta})'$:}
  \begin{equation} \label{Jprecise} 
  C^{\infty}_c(X^L_{\Theta})' := C^{\infty}_c(X^L_{\Theta}, \delta_{\Theta}) \otimes \delta_{\Theta}^{-1}.
  \end{equation}
  
A choice of a measure on    $U_{\Theta}$ identifies $C^{\infty}_c(X^L_{\Theta})'$ with  $C^{\infty}_c(X^L_{\Theta})$.

\subsubsection{The definition of $C_c^\infty(X_\Theta)_\sigma$}

Let $\sigma$ be an irreducible representation of $L_{\Theta}$.
  Any choice of parabolic in the class of $P_\Theta^-$, and hence of a Levi subvariety $X_\Theta^L$ of $X_\Theta$, gives rise to a map:
$$ \Hom_{L_\Theta}(C_c^\infty(X_\Theta^L),\sigma)\hookrightarrow  \Hom_G (C_c^\infty(X_\Theta), I_{\Theta^-}(\sigma)),$$
obtained by viewing $I_{\Theta^-}$ as a functor, and $C_c^\infty(X_\Theta)$ as $I_{\Theta^-}(C_c^\infty(X_\Theta^L))$.  Note that it is important for the validity of this statement
that we twisted the action on $C^{\infty}(X^L_{\Theta} )$.

 This embedding gives a corresponding quotient of the $I_{\Theta^-}(\sigma)$-coinvariants:\footnote{Note that we will prove that the representations $I_{\Theta^-}(\sigma)$ are irreducible for almost every $\sigma$ in the family, cf Corollary \ref{cor:genericbehavior}. }
\begin{eqnarray*}  C_c^\infty(X_\Theta)_{I_{\Theta^-}(\sigma)} &=& \left( \Hom_G (C_c^\infty(X_\Theta), I_{\Theta^-}(\sigma)) \right)^* \otimes I_{\Theta^-}(\sigma)  \\ & \twoheadrightarrow& \left(\Hom_{L_\Theta}(C^{\infty}_c(X_\Theta^L),\sigma )\right)^*\otimes I_{\Theta^-}(\sigma)\end{eqnarray*}
and the last quotient is what we call \emph{the
 $\sigma$-coinvariants} $C^{\infty}_c(X_{\Theta})_{\sigma} $.

There should be no confusion with ``$\pi$''-coinvariants, when $\pi$ is a representation of $G$, since the fact that $\sigma$ is a representation of the Levi suggests that we are using the structure of $C_c^\infty(X_\Theta)$ as an induced representation. As a quotient of $C_c^\infty(X_\Theta)$, the space $C_c^\infty(X_\Theta)_\sigma$ does not depend on any of the choices made, though the isomorphism:
\begin{equation}\label{Vminus} 
\Vminus \simeq \left(\Hom_{L_\Theta}(C^{\infty}_c(X_\Theta^L),\sigma)\right)^*\otimes I_{\Theta^-}(\sigma)
\end{equation}
does.

We denote the dual of $C^{\infty}_c(X_{\Theta})_{\sigma}$ in $C^\infty(X_\Theta)$ by $C^{\infty}(X_{\Theta})^{\tilde{\sigma}}$, and note that it is isomorphic to 
(writing $\tau = \tilde{\sigma}$):
\begin{equation} C^{\infty}(X_{\Theta})^{\tau} \simeq  \Hom_{L_{\Theta}}(\tau, C^{\infty}(X^L_{\Theta})) \otimes I_{\Theta^-}(\tau).\end{equation}

\subsubsection{The definition of $C^\infty(X_\Theta^h,\delta_\Theta)_{X,\sigma}$}

Let us recall that $C^{\infty}(X_{\Theta}^h, \delta_{\Theta})_X$ is the image of Radon transform of $C_c^\infty(X)$.  The space $C^\infty(X_\Theta^h,\delta_\Theta)_{X,\sigma}$
will be a certain $I_{\Theta}(\sigma)$-isotypical quotient of that space, defined for almost all $\sigma$ in each twisting class. As before, the term \emph{$\sigma$-coinvariants} will be used for that quotient.

Choosing a parabolic $P_\Theta$ gives rise in a similar way as above to a subvariety of $X_\Theta^h$ which is canonically isomorphic to $X_\Theta^L$. 
Our normalization of the action of $L_\Theta$ (\ref{Ltwist}) implies that the restriction maps give a $P_\Theta$-equivariant surjection: 
$$C_c^\infty(X_\Theta^h,\delta_\Theta) \to C_c^\infty(X_\Theta^L)'.$$ 
(Recall that $C_c^\infty(X_\Theta^L)'$ was defined in \eqref{Jprecise}.) Composing with maps into $\sigma$ we get a canonical embedding:
\begin{equation}\label{blue}\Hom_{L_\Theta}(C_c^\infty(X_\Theta^L),\sigma)\hookrightarrow \Hom_G (C_c^\infty(X_\Theta^h,\delta_\Theta),I_\Theta(\sigma)').\end{equation}
We leave it to the reader to check the canonicity of this embedding, just recalling here that both $C_c^\infty(X_\Theta^h,\delta_\Theta)$ and $I_\Theta(\sigma)'$ were defined by pulling back a certain line bundle over $P_\Theta\backslash G$.

The morphisms $C^{\infty}_c(X_\Theta^h,\delta_\Theta) \rightarrow I_{\Theta}(\sigma)'$ that arise
in the image of \eqref{blue} may not extend to $C^{\infty}(X_\Theta^h, \delta_{\Theta})_X$.  However, it will follow from Proposition \ref{proofextend} that they extend for \emph{almost every} $\sigma$; we state it in a vague form, which will be clarified by Definition \ref{defRadoncoinvts} and Proposition \ref{proofextend}:
\begin{proposition}[Proved as Proposition \ref{proofextend}.] \label{morphismsextend}
The map (\ref{blue}) extends naturally, for almost every $\sigma$ in every twisting class (cf.\ \S \ref{sss:twistingclass}), to a map:
\begin{equation}\label{extends}\Hom_{L_\Theta}(C_c^\infty(X_\Theta^L),\sigma)\hookrightarrow \Hom_G (C^\infty(X_\Theta^h,\delta_\Theta)_X,I_\Theta(\sigma)')
\end{equation}
\end{proposition}

We can finally define the desired quotient $ C^\infty(X_\Theta^h,\delta_\Theta)_{X,\sigma}$ as 
 the image of $C^{\infty}(X_{\Theta}^h, \delta_{\Theta})_X$ under the mapping 
\begin{equation} \label{Vplus} C^{\infty}(X_{\Theta}^h, \delta_{\Theta})_X \rightarrow  \left( 
\Hom_{L_\Theta}(C_c^\infty(X_\Theta^L),\sigma)\right)^* \otimes I_{\Theta}(\sigma)'.\end{equation} 
obtained by dualizing \eqref{extends}. 
Note that this quotient is $I_{\Theta}(\sigma)$-isotypical. 
 
Finally, the combination of (\ref{Vminus}), (\ref{Vplus}) shows that \emph{we have a canonical morphism}:
\begin{equation}\label{Tw}
RT_\Theta:  C^{\infty}_c(X_{\Theta}) \to C^{\infty}_c(X_\Theta^h, \delta_{\Theta})_{X, \sigma}
\end{equation}
 induced by the standard intertwining operators \eqref{intops}.  Indeed, it is immediate to check that this morphism does not depend on the choices of parabolic.     These morphisms are defined for almost every $\sigma$, and 
 will be seen to be invertible for almost every $\sigma$ by Corollary \ref{cor:genericbehavior}.

Thus, by definition, the bottom square of \eqref{green} commutes.   { Let us explain also the commutativity of the top square. Because the vertical maps for the bottom square are isomorphisms, it is enough to check that the ``big square'' commutes, that is to say, 
\begin{equation}  \begin{CD}
C^{\infty}_c(X_{\Theta}) @>{R_{\Theta}}>> C^\infty(X_\Theta^h, \delta_\Theta)_X \\
@VVV @VVV \\
\left({\Hom_{L_\Theta}(C^{\infty}_c(X_\Theta^L),\sigma )}\right)^*\otimes I_{\Theta^-}(\sigma)
@>{\mathrm{id} \otimes T_{\Theta}}>> \left({\Hom_{L_\Theta}(C^{\infty}_c(X_\Theta^L),\sigma )}\right)^*\otimes I_{\Theta}(\sigma)'
 \end{CD}
\end{equation}
 Ignoring issues of convergence, directions around this square are given
by taking a function on $X_{\Theta}$,   integrating along $U_{\Theta}$-orbits on the open orbit, and projecting
to $\sigma$-coinvariants.  In the next section we will see that these integrals are absolutely convergent in a Zariski-open subset of each twisting
class (\S \ref{sss:twistingclass} for definition of a twisting class, Proposition \ref{proofextend} and Corollary \ref{cor:genericbehavior} for the convergence). Hence, in this Zariski open set the diagram commutes, and this extends to all $\sigma$ for which both composites are defined. 
 }

\subsection{Convergence issues and affine embeddings} 
 
 Fix a $\Theta\subset\Delta_X$. We consider a $G$-stable class $\mathcal F$ of smooth functions $f$ on $X$ or $X_\Theta$ with the following properties (stated here with respect to $\XX$):\label{Fproperties}
\begin{enumerate}
 \item \label{condaffine} There is an affine embedding $\overline{\XX}$ of $\XX$ such that the support of all $f\in\mathcal F$ has compact closure in $\bar X$. In a slight variation of the language of \cite{BK}, having implicitly fixed the affine embedding $\overline{\XX}$, we will say that such functions have \emph{bounded support}.
 \item \label{condmoderate} For a given compact open $J\subset G$, the elements of $\mathcal F^J$ are uniformly of \emph{moderate growth}; i.e.:  {  
There is a completion $\widetilde{\XX}$ of $\XX$, a finite open cover (in the Hausdorff topology) $\tilde X = \bigcup_i U_i$, and, for each $i$, a rational function $F_i$ which is regular on $U_i \cap X$, such that each $f\in \mathcal F^J$ satisfies:
$$ |f|\le C_f |F_i|$$
on $U_i$, where $C_f$ is a constant that depends on $f$.
}
\end{enumerate}

Note, in particular, that condition \eqref{condaffine} guarantees that the Radon transform
$R f$ is defined for $f \in \mathcal{F}$, since  the orbits of a unipotent group on an affine variety are Zariski closed.

For a complex character $\omega$ of the $k$-points of an algebraic group $\MM$ we can write its absolute value in terms of absolute values of algebraic characters of $\MM$:
$$ |\omega| = \prod_i |\chi_i|^{s_i}$$
where each $\chi_i: \MM \rightarrow \mathbb{G}_m$ is algebraic; 
  then we define the \emph{real part} of $\omega$: 
\begin{equation}
 \Re \omega:= \sum_i \Re s_i \chi_i \in \mathfrak m^*:=\varchi(\MM)\otimes \RR
\end{equation}
which is independent of choices.

Recall that $\LL_{\Theta,X}^\ab$ be the torus quotient of $\LL_\Theta$ whose character group consists of all characters of $\LL_\Theta$ which are trivial on the stabilizers of points on $\XX_\Theta^L$. Let us choose a parabolic $\PP_\Theta$, giving rise to the quotient map:
\begin{equation}\label{quotient} \mathring \XX \PP_\Theta/\UU_\Theta \simeq \XX_\Theta^L.
\end{equation}
Also, choose a base point on $X_\Theta^L$ in order to identify characters of $\LL_{\Theta,X}^\ab$ with functions on $\XX_\Theta^L$; our statements will not depend on any of these choices.

\begin{lemma}\label{lemmaconvergence}
 Consider an algebraic character $\chi\in\varchi(\LL_{\Theta,X}^\ab)$ as a function on $\mathring\XX \cdot \PP_\Theta$ via the quotient map (\ref{quotient}). Let $\overline{\XX}$ be an affine embedding of $\XX$. Then for $\chi$ in an open subcone\footnote{i.e., in a generating, saturated submonoid of $\varchi(\LL_{\Theta,X}^\ab)$ -- the intersection of $\varchi(\LL_{\Theta,X}^\ab)$ with an open subcone of $\varchi(\LL_{\Theta,X}^\ab)\otimes \mathbb R$} of $\varchi(\LL_{\Theta,X}^\ab)$ this function extends to a regular function on $\overline{\XX}$ which vanishes on  $\overline{\XX}\smallsetminus \mathring\XX \cdot \PP_\Theta$.
\end{lemma}

\begin{proof}
Consider the quotient map of $\LL_\Theta$-spaces:
$$\overline{\XX}\to \overline{\XX}\sslash\UU_\Theta =\spec k[\overline{\XX}]^{\UU_\Theta}.$$

We claim that $\XX_\Theta^L$ embeds as the open $\LL_\Theta$-orbit in $\overline{\XX}\sslash\UU_\Theta$, and its preimage is precisely $\mathring\XX\cdot \PP_\Theta$.

If $k[\overline{\XX}]=\oplus_{\lambda\in \varchi(\XX)^+} V_\lambda$ denotes the decomposition of $k[\overline{\XX}]$ into a (multiplicity-free) sum of irreducible subrepresentations, where 
$\varchi(\XX)^+$ is a saturated (by normality), generating (by quasi-affineness) submonoid of $\varchi(\XX)$, depending on $\overline{\XX}$,
then highest weight theory implies that $k[\overline{\XX}]^{\UU_\Theta}$ has the following multiplicity-free decomposition into irreducible $\LL_\Theta$-representations:
\begin{equation} \label{pink} k[\overline{\XX}]^{\UU_\Theta} = \oplus_{\lambda \in \varchi(\XX)^+} V_\lambda^{\UU_\Theta}.
\end{equation}
In particular, it is finitely generated: indeed, it is generated by the sum of $V_\lambda^{\UU_{\Theta}}$'s for a set of $\lambda$'s generating $\varchi(\XX)^+$.\footnote{To see that, note that if  $\Lambda$ is such a set
of $\lambda$s, and $\chi = \sum_{\lambda \in \Lambda} n_{\lambda} \lambda$, 
then, for $v_{\lambda}$ a highest weight vector in $V_{\lambda}$,  the product $\prod_{\lambda} v_{\lambda}^{n_{\lambda}}$ is a $\BB$-invariant vector 
of weight $\chi$; that shows that $V_{\chi}^{\UU_{\Theta}}$ is contained
is contained in the image of a suitable tensor product of $V_{\lambda}^{\UU_{\Theta}}$s 
by a multiplication map.}

On the other hand, we have a decomposition: 
$$k[\XX^L_{\Theta}]=\oplus_{\lambda\in \varchi(\XX_\Theta^L)^+} V'_\lambda,$$
where $V'_\lambda$ now denotes the highest weight module of weight $\lambda$ for $\LL_\Theta$, and $\varchi(\XX)^+\subset \varchi(\XX_\Theta^L)^+ \subset \varchi(\XX)$. By choosing a finite set of generators of $\varchi(\XX_\Theta^L)^+$ and a suitable -- finite -- set of $G$-translates of the corresponding highest weight vectors, we obtain a finite set of elements in the fraction field of $k[\overline{\XX}]^{\UU_\Theta}$ generating $k[\XX_\Theta^L]$; hence, the morphism: $\XX_\Theta^L\to \overline{\XX}/\UU_\Theta$ is birational and dominant, and since $\XX_\Theta^L$ is homogeneous it is an open embedding. 

On the other hand, let us verify that the preimage of the open $\LL_{\Theta}$-orbit (call it $\OO$) in 
 $\overline{\XX}\sslash\UU_\Theta$ is precisely $\mathring\XX\cdot \PP_\Theta$. 
 Were this not so, there is another $\PP_{\Theta}$-orbit on $\overline{\XX}$ whose image
 is equal to $\OO$;  in particular, there is some $\BB$-orbit  $\ZZ \subset \overline{\XX}$, disjoint from $\mathring{\XX}$, whose image contains the open $\BB \cap \LL_{\Theta}$ orbit on $\overline{\XX} \sslash \UU_{\Theta}$. 
Then the map $\ZZ\to \overline{\XX}\sslash\UU_\Theta$ is dominant; that would imply that no non-zero element of $k[\overline{\XX}\sslash\UU_\Theta]^{(\BB)} = k[\overline{\XX}]^{(\BB)}$ vanishes on $\ZZ$. This cannot be the case: the complement of the open $\BB$-orbit in $\overline{\XX}$ is a closed, $\BB$-stable subvariety; consider the $\BB$-stable ideal of regular functions vanishing on it -- it must contain non-trivial $\BB$-semiinvariants.

Now, it suffices to prove that for any affine embedding $\YY$ of a factorizable spherical $\LL_\Theta$-variety $\XX_\Theta^L$ the cone of characters of $\LL_{\Theta,X}^\ab$ which vanish in the complement of the open orbit is non-trivial and, actually, of full rank. This is the case, of course, for affine toric varieties, and we will reduce to this case using the quotient map:
$$\YY\to \YY\sslash[\LL_\Theta,\LL_\Theta].$$

All we need to prove is that the preimage of the open $\LL_\Theta$-orbit on $\YY\sslash [\LL_\Theta,\LL_\Theta]$ is not larger than $\XX_\Theta^L$.
Recall that, when a reductive group acts on an affine variety, any two closed
sets are separated by an invariant function (see \cite{Mumford}). 
It suffices, then, to show that  all $[\LL_\Theta,\LL_\Theta]$-orbits on $\XX_\Theta^L$ are closed in $\YY$. We claim that these are spherical $[\LL_\Theta,\LL_\Theta]$-varieties with a finite number of automorphisms -- then by \cite[Corollary 7.9]{KnAs} they have no non-trivial affine embeddings, hence have to be closed in $\YY$. Finally, to show that the $[\LL_\Theta,\LL_\Theta]$-orbits on $\XX_\Theta^L$ are spherical without continuous group of automorphisms, we use the hypothesis that $\XX$ is strongly factorizable -- hence $\XX_\Theta^L$ is factorizable under the $\LL_\Theta$-action. Recall that the connected $\LL_\Theta$-automorphism group of $\XX_\Theta^L$ is induced by the action of $\mathcal Z(\LL_\Theta)$; hence, factorizability means that $\HH\cap [\LL_\Theta,\LL_\Theta]$ is spherical in $[\LL_\Theta,\LL_\Theta]$ and coincides with the connected component of its normalizer there.
\end{proof}

{ \begin{remark}
A knowledge of the combinatorial data describing $\overline{\XX}$ allows to read off the precise cone of characters which vanish on the complement of $\mathring\XX\PP_\Theta$. Indeed, the above proof shows that these are precisely the characters which vanish in the complement of the open orbit of the toric $\LL_{\Theta,X}^\ab$-variety:
$$ \overline{\XX}\sslash [\LL_\Theta,\LL_\Theta]\UU_\Theta = \spec k[\overline\XX]^{[\LL_\Theta,\LL_\Theta]\UU_\Theta} = \spec k[\varchi(\XX)^+\cap \varchi(\LL_\Theta)].
$$
Hence, the monoid of characters which extend to the complement is the set of those elements of $\varchi(\XX)^+$ which are characters of $\LL_\Theta$, and the characters that vanish on the complement are those in the ``interior'' of the monoid.
\end{remark}
}

Let us see what this lemma implies. 
Again, if we fix opposite parabolics $\PP_\Theta,\PP_\Theta^-$, we can regard $\XX_\Theta^L$ as a subvariety of both $\XX_\Theta$ and $\XX_\Theta^h$.  
Fix a point $x_0\in X_\Theta^L$ in order to define a quotient map:
\begin{equation}\label{Levitochars} \XX_\Theta^L\to \LL_{\Theta,X}^\ab,\end{equation}
and use it to consider characters of $L_{\Theta,X}^\ab$ as functions on $X_\Theta^L$. If $M:C_c^\infty(X_\Theta^L)\to \sigma$ is a morphism of $L_\Theta$-representations and $\omega$ is a character of $L_{\Theta,X}^\ab$, we get a morphism:  
$$ \omega^{-1} M: \omega^{-1} \otimes C_c^\infty(X_\Theta^L)\to \omega^{-1}\sigma$$
simply by twisting by $\omega$. Explicitly, if $f\in \omega^{-1} \otimes C_c^\infty(X_\Theta^L)$ we have:
$$ \omega^{-1} M (f) = M(\omega f) \in \omega^{-1}\sigma,$$
where the underlying vector spaces of $C_c^\infty(X_\Theta^L)$ and $\omega^{-1} \otimes C_c^\infty(X_\Theta^L)$, as well as those of $\sigma$ and $\omega^{-1}\sigma := \omega^{-1}\otimes \sigma$ have been identified.

{ 
The statement of the following corollary will make use of the concept of ``extension of a morphism by a convergent series'', by which we mean the following: Let $\YY\subset\overline\YY$ be smooth varieties with an action of a group $\LL$, with $\YY$ open, and let $M: C_c^\infty(Y)\to \sigma$ be a morphism to a smooth representation $\sigma$. We say that it ``extends by a convergent series'' to $C_c^\infty(\overline Y)$ if for every $\Phi\in C_c^\infty(\overline Y)$ and $v\in \tilde\sigma$ (the smooth dual of $\sigma$) 
the inner product $\left<\Phi, \tilde M(v)\right>$ (where $\tilde M$ is the adjoint of $M$) converges absolutely, defining a morphism: 
$C_c^\infty(\overline Y)\to \widetilde{\tilde\sigma}$. 
(In the corollary, $\sigma$ is admissible so $\widetilde{\tilde\sigma}=\sigma$.) 
Here we think of $\tilde M(v)$ as a smooth measure on $Y$, hence one could equivalently write $\Phi$ as a convergent sum $\sum_i \Phi_i$ with $\Phi_i\in C_c^\infty(Y)$ 
and require that $\sum_i \left<\Phi_i,\tilde M(v)\right>$ converges, hence the language.
}

\begin{corollary} \label{corollaryconvergence}
 Let $\mathcal F$ be a class of functions as on p.\ \pageref{Fproperties}. Let $M:C_c^\infty(X_\Theta^L)\to \sigma$ be a morphism to an admissible $L_\Theta$-representation. Then:
\begin{enumerate}
 \item If $\mathcal F$ consists of functions on $X_\Theta$, for $\Re(\omega)$ in a translate of an open cone as in the previous lemma, the morphism $\omega^{-1} M$ can be extended by a convergent series to functions of the form: $f|_{X_\Theta^L}$, $f\in\mathcal F$. 
 \item If $\mathcal F$ consists of either functions on $X$ or on $X_\Theta$, the analogous statement holds for the Radon transform $R_\Theta$ of $\mathcal F$:  the morphism $\omega^{-1} M \otimes \delta_\Theta: \omega^{-1} \otimes C_c^\infty(X_\Theta^L,\delta_\Theta)\to \sigma'$ can be extended by a convergent series to functions of the form: $\left.R_\Theta f\right|_{X_\Theta^L}$.
 \item In the case of $X_\Theta$, if $M:C_c^\infty(X_\Theta^L)\to \sigma$ is a morphism to an irreducible unitary representation of $L_\Theta$ and $T_\Theta: I_{\Theta^-}(\omega^{-1}\sigma)\to I_{\Theta}(\omega^{-1}\sigma)'$ is the standard intertwining operator, then $T_\Theta$ is defined by a convergent integral and the following diagram commutes: 
\begin{equation} 
 \begin{CD}
  \mathcal F @>{R_\Theta}>> C^\infty(X_\Theta^h, \delta_\Theta) \\
 @V{I_{\Theta^-} (\omega^{-1}M)}VV  @VV{I_\Theta (\omega^{-1}M)}V \\
 I_{\Theta^-}(\omega^{-1}\sigma) @>{T_\Theta}>> I_\Theta(\omega^{-1}\sigma)'
 \end{CD}
\end{equation}
Here $\omega^{-1}M$ denotes the twist of $M$ defined above, and $I_\Theta (\omega^{-1}M)$, $I_{\Theta^-} (\omega^{-1}M)$ are the maps obtained by functoriality of induction. (The right vertical arrow is defined on the image of $\mathcal F$ by the previous statement, not on the whole space $C^\infty(X_\Theta^h, \delta_\Theta)$.)
\end{enumerate}
\end{corollary}
 
{
In the arguments that follow, we will implicitly use the following easy fact: If $\YY$ is a homogeneous variety for a group $\HH$, and the set $Y$ of its $k$-points is equipped with an $H$-eigenmeasure $dy$, if $\overline{\YY}$ is an embedding of $\YY$ and if $P$ is a regular function on $\overline{\YY}$ which vanishes on $\overline{\YY}\smallsetminus \YY$, then for any moderate-growth function $f$ on $Y$ whose support has compact closure in $\bar Y$, the function $|P^n|f$ is in $L^1(Y,dy)$ for sufficiently large $n$. The reason is that the eigenmeasure itself is of moderate growth, i.e.\ in a neighborhood of a point of $\bar Y$, choosing local coordinates for the etale topology, the measure can be written as $h(x)dx$, where $dx$ is the usual Lebesgue measure in these coordinates and $h$ is a function of moderate growth. Replacing $\overline\YY$ by a blowup $\widetilde\YY$, so that $f$ and $h$ are bounded, locally in a (compact, without loss of generality) neighborhood $U_i$ around any point $y \in \widetilde Y\smallsetminus Y$ by a rational function $F_i$ defined on $U_i\cap Y$ as in the definition of ``moderate growth'' (\S \ref{Fproperties}), we get that the integral of $|P^n|f$ on $U_i$ is bounded by:
\begin{equation} \label{Cronut} \int |P^n F_i| (x)  dx,\end{equation} 
and for large enough $n$ the function $P^n F_i$ has no poles (and hence is bounded) on $U_i$.
}

\begin{proof}[Proof of the corollary]
{ \label{Radon_convergence_pageref} The first two statements follow from the fact that the support of the functions $f\in \mathcal F$ is contained in a compact subset of an affine embedding of $\XX$ or $\XX_\Theta$, and that they are uniformly of moderate growth. (Indeed, their asymptotics in every direction are governed by exponents, by the theory of asymptotics that we have developed and the fact that the Jacquet functor preserves admissibility.)
 
More precisely, let $\chi$ be an algebraic character of $\LL_{\Theta,X}^\ab$ as in Lemma \ref{lemmaconvergence}, considered as a function on $\mathring\XX \PP_\Theta$ or $\mathring \XX_\Theta\PP_\Theta$ by fixing a base point and extending by zero to $\overline\XX$. 
The given morphism $M:C_c^\infty(X_\Theta^L)\to \sigma$ has an adjoint $\tilde M:\tilde\sigma\to C^\infty(X_\Theta^L)$ and, because $\sigma$ is admissible, the functions in its image are also of uniformly moderate growth on $X_\Theta^L$. Notice that $X_\Theta^L$ is closed in $\mathring X \cdot P_\Theta$. Let $f_1$ be in the image of $\tilde M$, and let $f_2\in \mathcal F$. 

For the first claim, applying the remark before the proof to $\YY=\XX_\Theta^L$ and $\overline{\YY}=$its closure in the given (from the properties of $\mathcal F$) affine embedding $\overline\XX_\Theta$ of $\XX_\Theta$, we get that the product:
$$ |\chi^n| f_1 f_2 $$
is an integrable function on $X_\Theta^L$ for $n\gg 0$. This proves the first claim. 

Similarly, for the second claim, if $f_2\in \mathcal F$, $v_1\in \tilde\sigma$ and $\omega=|\chi|^n$, the pairing:
\begin{equation} \label{Radon_convergence_equation}  \left< \omega^{-1}M\otimes\delta_\Theta \circ R_\Theta (f_2), v_1\right>\end{equation}  can be written as the integral over $X_\Theta^L\subset X_\Theta^h$ of $R_\Theta |\chi|^n f_1$ against an element $f_1\in C^\infty(X_\Theta^L, \delta_\Theta^{-1})$ of moderate growth or, equivalently (using the definition of Radon transform), as the integral of $f_2$ over $\mathring X P_\Theta$ (or $\mathring X_\Theta P_\Theta$) of $f_2$ against a function of moderate growth on $\mathring XP_\Theta$ (or $\mathring X_\Theta P_\Theta$). The same argument now shows that this is convergent for $n\gg 0$.}

For the third statement, let us first verify that $T_{\Theta}$ is absolutely convergent. This is the case for $\Re(\omega)$ in a translate of the dominant cone inside of $\varchi(\LL_{\Theta}) \otimes \RR$. 
 We need to show that $\varchi(\LL_{\Theta,X}^\ab)$ intersects the interior of that cone nontrivially. 
 
 The interior of the dominant cone of $\chi \in \varchi(\LL_{\Theta})$ consists of $\RR_+$-multiples of 
 precisely those algebraic characters $\chi$ for which the function
 $$f_{\chi}:  u^{-} l u \in \UU_{\Theta}^- \times \LL_{\Theta} \times \UU_{\Theta} \mapsto \chi(l)$$
 extends to a regular function on $\GG$
 that vanishes on the complement of  the ``open cell''
 $\UU_{\Theta}^-   \LL_{\Theta}  \UU_{\Theta} $.  (Without the vanishing condition,
 we do not get {\em strictly} dominant.) 
 
 But elements of $\varchi(\LL_{\Theta, X})$  which belong to $\mathcal T$ have that property: \label{Tdominant}
 We have seen that (considering  $\omega$ can be regarded as a function on the open $\PP_{\Theta}$
 orbit on $\XX_{\Theta}$) that $\omega$ extends
 to a regular function on $\XX_{\Theta}$.
Now if we consider the orbit map:
 $$ \GG\ni g\mapsto x\cdot g\in \XX_{\Theta} \rightarrow  \simeq \XX^{L}_{\Theta} \times_{\PP_{\Theta}^-} \GG,$$
where $x\in X_\Theta$ is a point mapping to the chosen point $x_0$ of $X_\Theta^L$,
 the function $\omega$ on $\XX_{\Theta}$ pulls back to the function  $f_{\omega}$
 on $\GG$.
 So $f_{\omega}$
 extends to a function on $\GG$, and so $\omega$ is dominant. 
 But even better:  because {\em the preimage of the open orbit in $\XX_{\Theta}$ is the open cell  $\UU_{\Theta}^-   \LL_{\Theta}  \UU_{\Theta} $}, and so
$f_{\omega}$ vanishes on the complement of that open cell, and so $\omega$ is strictly dominant. 
 This concludes the proof that $T_{\Theta}$ is absolutely convergent. 
 
The rest of the third statement will be formal after we unwind its meaning;
 it is simply the fact that the standard intertwining operator is given (in the appropriate context)
 by a Radon transform. 
 
  Indeed, 
  we have:
$$ I_\Theta (\omega^{-1}M) \circ R_\Theta (f) (1) = \omega^{-1}M \left(\left.R_\Theta f\right|_{X_\Theta^L}\right) = M \left(\omega \left.R_\Theta f\right|_{X_\Theta^L}\right) =$$
$$ =  M \left(\omega \lim_{n} \left.R_{U_n}f\right|_{X_\Theta^L}\right).$$
Here $U_n$ denotes a sequence of compact subgroups exhausting $U_\Theta$, and $R_{U_n}$ denotes the partial Radon transform: $R_{U_n} f(x) = \int_{U_n} f(x\cdot u) du$. For $x$ in any compact subset of $X_\Theta^L$, this limit eventually stabilizes because unipotent orbits in affine varieties are closed, and $f$ has compact support in some affine embedding. 

Next, we can interchange $M$ and the limit here. To check that, 
we must  give a bound on $M( \omega \left( R_{U_n} - R_{U_{\infty}} f \right) )$ 
in absolute value that goes to zero with $n$.  But this follows by examining the reasoning
by which we verified the convergence in the first place:   As in \eqref{Radon_convergence_equation} 
we must bound the integral of $ \left( R_{U_n} - R_{U_{\infty}} f \right)|_{X^L_{\Theta}} \cdot \omega$
against a function of moderate growth on $X^L_{\Theta}$, and this bound should go to zero with $n$. As we noted above,  the quantity
$ \left( R_{U_n} - R_{U_{\infty}} f \right)|_{X^L_{\Theta}} \cdot \omega$ vanishes on any fixed compact subset of $X^L_{\Theta}$
for large enough $n$.   Now, $R_{U_n} f$ is bounded in absolute value by  $R_{U_{\infty}}|f|$; 
and the rest follows just as in \eqref{Radon_convergence_equation}, using 
again the discussion around \eqref{Cronut}. 

Thus, interchanging $M$ and the limit, we get
$$I_\Theta (\omega^{-1}M) \circ R_\Theta (f) (1) = \lim_n M \left(\omega \left.R_{U_n}f\right|_{X_\Theta^L}\right)=$$
$$ =\lim_{n} \int_{U_n} M \left(\left.\omega \cdot (u\cdot f)\right|_{X_\Theta^L}\right) du = \lim_{n} T_{U_n} \circ I_{\Theta^-} (\omega^{-1}M)(f) \ (1),$$
where $T_{U_n}$ is the analogous partial version of $T_\Theta$. Since we 
verified that $T_{\Theta}$ is absolutely convergent, the last integral equals $$T_\Theta \circ I_{\Theta^-} (\omega^{-1}M) (f)(1).$$

\end{proof}

\begin{remark}
 The reader will notice from the proofs of Lemma \ref{lemmaconvergence} and Corollary \ref{corollaryconvergence} that the subcone of $\varchi(\LL_{\Theta,X}^\ab)\otimes \RR$ of the corollary depends only on the toric variety $\YY\sslash [\LL_\Theta,\LL_\Theta]\UU_\Theta$, where $\YY$ is the affine $\GG$-variety containing the support of elements of the class $\mathcal F$.

In particular, for the cases $\mathcal F=C_c^\infty(X)$ and $\mathcal F=C_c^\infty(X_\Theta)$, where $\YY$ can be taken to be the affine closure of $\XX$, resp.\ $\XX_\Theta$, the same cone $\mathcal T$ will work. Indeed, there is an affine embedding of $\XX_\Theta$ whose coordinate ring is, as a $\GG$-module, isomorphic to $k[\XX]$ (\S \ref{ssdegen}), and then the corresponding categorical quotients by $[\LL_\Theta,\LL_\Theta]\UU_\Theta$ will coincide.

{\em We will write $\omega\gg  0$  for a character $\omega$ as in the Corollary, when $\mathcal F=C_c^\infty(X)$ or $C_c^\infty(X_\Theta)$.} 
\end{remark} 
\begin{definition}\label{defRadoncoinvts}
 For $\omega\gg 0$, we define the extension of \eqref{blue} to a map:
\begin{equation}\Hom_{L_\Theta}(C_c^\infty(X_\Theta^L),\omega^{-1}\sigma)\hookrightarrow \Hom_G (C^\infty(X_\Theta^h,\delta_\Theta)_X,I_\Theta(\omega^{-1}\sigma)'),\end{equation}
by inducing from the extension of a morphism $\omega^{-1}M$ to $R_\Theta f|_{X_\Theta^L}$ as guaranteed by Corollary \ref{corollaryconvergence}.
\end{definition}

It is this extension which enabled us to define the coinvariant space $C_c^\infty(X_\Theta^h,\delta_\Theta)_{X,\omega^{-1}\sigma}$ for $\omega\gg 0$. The extension to almost every $\omega$ will follow from Proposition \ref{proofextend}.

\begin{corollary} \label{cor:genericbehavior}
 For $\sigma$ in a dense, Zariski open subset of a twisting class:
 \begin{enumerate}
 \item[i.] the induced representation $I_{\Theta^-}(\sigma)$ is irreducible;
 \item[ii.]  $T_\Theta$ is an isomorphism;
  \end{enumerate}
 
\end{corollary}

\begin{proof}
It is enough to show that conditions (i) and (ii) are verified at a single point, 
because the twisting class is an irreducible variety and if conditions
(i) or (ii) hold at a point of such a variety, they hold at a Zariski-open set.\footnote{Indeed, irreducibility amounts to asking that certain elements of a Hecke algebra, 
with respect to an open compact $J$,
generate the endomorphisms of $J$-fixed vectors; and if an algebraic family
of matrices has full rank at a point, it has full rank at a Zariski-open set. Similarly for (ii): once the representation
is irreducible, the map $T_{\Theta}$ is an isomorphism if and only if it is nonzero, which can be checked
on a single $J$-fixed vector.}

   We have already seen that the natural map
 $$ \varchi(\LL_{\Theta, X}) \otimes \RR \rightarrow \varchi(\LL_{\Theta}) \otimes \RR$$
 carries the cone $\mathcal{T}$ of Corollary \ref{corollaryconvergence} into the strongly dominant cone (see p.\ \pageref{Tdominant}). 
But it is well-known \cite[Theorem 6.6.1]{Casselman-notes} that, if one twists $\sigma$ by a sufficiently dominant
character, the induced representation is irreducible and the intertwining operator is an isomorphism.

\end{proof}

\begin{proposition} \label{proofextend}
Let $\mathcal F= C_c^\infty(X)$ or $C_c^\infty(X_\Theta)$, and 
let $M$ be an element of $ \Hom_{L_\Theta}(C_c^\infty(X_\Theta^L),\sigma) $. The composition of:
\begin{equation}\label{Eisenstein-unnormalized}
 \mathcal F \xrightarrow{R_\Theta} C^\infty(X_\Theta^h,\delta_\Theta) \xrightarrow{I_\Theta(\omega^{-1}M)} I_\Theta(\omega^{-1}\sigma)',
\end{equation}
(the second arrow being defined only on the image of the first), which converges for $\omega\gg 0$ by Corollary \ref{corollaryconvergence}, extends to a rational family of morphisms for all $\omega$.
\end{proposition}

 For $\mathcal F = C_c^\infty(X)$ this composition is just the unnormalized Eisenstein integrals of the literature (or rather, their adjoints).
We  also note that, in the case of  $X$ symmetric, this Proposition proven by Blanc and Delorme \cite{BlancDelorme} (see also discussion of their paper in \S \ref{BDtranslated}). 

\begin{proof}
We will refer to $X$ in our notation, but the argument for $X_\Theta$ is verbatim the same.

Let us choose a measure on $U_\Theta$ in order to identify $I_\Theta(\omega^{-1}\sigma)'$ with $I_\Theta(\omega^{-1}\sigma)$. For any $v^*\in \tilde\sigma$ the composition of these arrows with ``restriction to the coset $P_\Theta 1$'' and ``pairing with $v^*$'':\footnote{Again, we are implicitly identifying the underlying vector spaces of all representations in the family $\omega^{-1}\sigma$ (as $\omega$ varies), and hence $v^*$ lives in the dual of all of them.}
$$\mathcal F \to I_{\Theta} (\omega^{-1}\sigma) \to \omega^{-1}\sigma \xrightarrow{\left<\bullet,v^*\right>} \CC$$
is given formally by an integral
$$f \mapsto \int_X f\cdot F_\omega dx$$
(resp.\ such an integral on $X_\Theta$), and that the function $F_\omega$ is the product of a fixed, locally constant, $U_\Theta$-invariant function $F$ on $\mathring X P_\Theta$ (namely, $M^* v^*$, via the map $\mathring X P_\Theta\to X_\Theta^L$), and the character $\omega$.  

Consider the quotient: $\XX\to \XX\sslash \UU_\Theta$, which induces an affine embedding of $\XX_\Theta^L$. Let $\YY\to \XX\sslash\UU_\Theta$ be a proper surjective morphism, where $\YY$ is a smooth toroidal embedding of $\XX_\Theta^L$. Such an embedding and morphism always exist; indeed, if $\mathcal C(\XX\sslash \UU_\Theta)$ is the cone of the affine embedding $\XX\sslash \UU_\Theta$ of $\XX_\Theta^L$ (we point to \cite[page 8]{KnLV} for the definition of this cone), then $\YY$ will be described by a fan whose support is  $\mathcal C(\XX\sslash \UU_\Theta) \cap \mathcal V$ (where $\mathcal V$ denotes the cone of invariant valuations of $\XX_\Theta^L$, cf.\ Section \ref{sec:review}). By \cite{KnLV} we have a morphism $\YY\to \XX/\UU_\Theta$, which is proper by \emph{loc.cit.\ }Theorem 4.2 and surjective by \emph{loc.cit.\ }Lemma 3.2.  

In particular, we know from the Local Structure Theorem \ref{localstructure} that the complement of $\XX_\Theta^L$ in $\YY$ is a union of divisors intersecting transversely. Let $\overline{\ZZ}$ be the closure of the image of the map:
$$ \mathring\XX\PP_\Theta\to \XX\times \YY$$
(natural inclusion times projection to $\XX_\Theta^L$). It comes with morphisms: $\overline{\ZZ}\to \XX$ and $\overline{\ZZ}\to \YY$, the former surjective and proper.
(In fact, $\overline{\ZZ}$ is contained in $\XX \times_{\XX /\UU_{\PP_{\Theta}}} \YY$,
and the latter has a proper map to $\XX$). 
Applying resolution of singularities, we may replace $\overline{\ZZ}$
by a nonsingular variety $\ZZ$ equipped with a proper birational morphism
$\ZZ \rightarrow \overline{\ZZ}$.   In particular the induced map $\ZZ \rightarrow \XX$
is proper also. 

We notice that a proper morphism of algebraic varieties over $k$ induces a proper map of the corresponding topological spaces of $k$-points, see \cite[Proposition 4.4]{Brian}. That means that the  
 pull-back of $f$ is a locally constant, compactly supported function on $Z$.
 
Hence it is enough to show that the integral of the pull-back of $F_\omega dx $ over a compact open neighborhood in $Z$ is rational in $\omega$. 
Here, we make sense of the  ``pullback of $dx$''  because $dx$ is the measure obtained as the absolute value of a volume form, which can be pulled back to $\ZZ$.
As for $F_{\omega}$, we extend it first of all by zero off $\mathring{\XX} \PP_{\Theta}$
and then pull back.   Note that, fixing $\omega_0$, we have, by definition,
$$F_{\omega} = F_{\omega_0} \cdot \prod |f_i|^{s_i(\omega)},$$
where the $f_i$ are rational functions on $\ZZ$ 
and the exponents $s_i(\omega)$ vary linearly with $\omega$.    This is a matter of the definitions: 
the $f_i$ are obtained  by pulling back  the coordinate functions on $\mathbb{G}_m^t$ under
the maps  $$\mathring{\XX} \PP_{\Theta} /\UU_{\Theta}  \rightarrow  \mathbf{X}^{L}_{\Theta} \rightarrow \mathbf{X}^{L}_{\Theta} / [\mathbf{L}_{\Theta}, \mathbf{L}_{\Theta}]
= \mathbf{L}_{\Theta, \X}^{\ab} \simeq \mathbb{G}_m^t.$$

For a normal $k$-variety $\VV$ with a distinguished divisor $\DD$, 
we have defined before
Corollary \ref{finiteasymp} the notion of a function on $V $ being ``$D$-finite.'' 
 By that Corollary \ref{finiteasymp}, $F$ is a $D$-finite function on $Y$, where $\DD$ is the complement of $\XX_\Theta^L$, and therefore (see discussion prior to quoted Corollary) its pull-back to $\ZZ$ is also  so (with respect to the complement of $\mathring{X} P_\Theta$).

We may now apply the following consequence of Igusa theory (\cite{Igusa}, see in particular Theorem 8.2.1 and the proof that follows;
that reference deals with a special case of a single $f_i$, but for a discussion of the modification necessary for many $f_i$
one can proceed in the fashion of \cite[p. 5]{Denef})
\begin{quote} Consider the integral 
$$ I(\omega) := \int_{V_0}  F \cdot |\Omega| \cdot \prod_{i} |f_i|^{s_i(\omega)},$$ 
where:
\begin{itemize}
\item[-] $V_0$ is an open compact subset of $V$: 
\item[-] $F$ is a $D$-finite function;
\item[-]  $\Omega$ an algebraic volume form, with polar divisor in $\DD$;
\item[-]  $f_i$ rational functions with polar divisor  in $\DD$;
\item[-] The exponents $s_i(\omega)$ vary linearly in $\omega$. 
\end{itemize} 
If $I(\omega)$ converges for a open set of the parameters in the Hausdorff topology, then it is rational in $\omega$.
\end{quote}

  This establishes the proposition.
\end{proof}

\subsection{Normalized Eisenstein integrals and smooth asymptotics}
 
\subsubsection{Definition of Eisenstein integrals}\label{ssexplicitsmooth}
We now define our normalized version of Eisenstein integrals. Recall
that  for $\sigma$ an irreducible unitarizable representation of $L_{\Theta}$, 
we have a diagram
    \begin{equation} \label{Bugle}
 \begin{CD}
C^{\infty}_c(X_{\Theta}) @.  C^\infty(X_\Theta^h, \delta_\Theta)_X  @<{R}<<  C^{\infty}_c(X) \\
@VVV @VVV  \\
C^{\infty}_c(X_{\Theta})_{\sigma} @<{RT_{\Theta}^{-1}}<< C^{\infty}(X^h_{\Theta}, \delta_{\Theta})_{X, \sigma}\\
 \end{CD}
\end{equation} 
Since $RT_{\Theta}^{-1}$ and $C^{\infty}(X^h_{\Theta}, \delta_{\Theta})_{X, \sigma}$ make sense only for generic $\sigma$ in a twisting class, we implicitly assume that we are referring to such $\sigma$.

 Note that the Radon transform maps $C^{\infty}_c(X_{\Theta})$
 into $C^{\infty}(X^h_{\Theta}, \delta_{\Theta})$, but we do not know, a priori, that
 the image is contained in $C^{\infty}(X^h_{\Theta}, \delta_{\Theta})_X$. 
 This is why we omit an arrow at the upper left. (It can be shown that the image is indeed in this space.)

 We denote by  
$$ R_{\Theta, \sigma}: C^{\infty}_c(X) \rightarrow
C^{\infty}(X^h_\Theta, \delta_{\Theta})_{X, \sigma} , \ \ E_{\Theta,\sigma}^*:C_c^\infty(X) \to C_c^\infty(X_\Theta)_\sigma,$$
the morphisms obtained by following the arrows in the above diagram. 
We refer to $E_{\Theta, \sigma}^*$ as the 
 \emph{adjoint normalized Eisenstein integral}. (The adjoints of \emph{unnormalized} Eisenstein integrals, already encountered in \eqref{Eisenstein-unnormalized}, are essentially the operators above without the last arrow representing $RT_\Theta^{-1}$.) The established language in the harmonic analysis of real symmetric spaces would suggest calling it (normalized) \emph{Fourier transform}, but we prefer to reserve the term ``Fourier transform'' for additive groups.  
 Now $E_{\Theta, \sigma}^*$ is defined a priori in a Zariski-dense subset of each twisting class; let us note that this automatically means
 that it is defined in a {\em full measure} subset of the unitary axis of the twisting class. 
 
The \emph{normalized Eisenstein integral} is the adjoint map: 
\begin{equation}
E_{\Theta,\sigma}: C^\infty(X_\Theta)^{\tilde\sigma} \to C^\infty(X),
\end{equation}
where $C^\infty(X_\Theta)^{\tilde\sigma}$ denotes the dual of $C_c^\infty(X_\Theta)_\sigma$, considered as a subspace of $C^\infty(X_\Theta)$. 
 
Note that our notion of Eisenstein integral involves in a crucial
way the ``intermediary'' of $X_\Theta^h$ between $X$ and $X_{\Theta}$.  From our point of view,
this is closely related to the appearance of factors related to the action of intertwining operators in the Plancherel formula.

The main theorem of this subsection is the following: The Plancherel formula for $X_\Theta$ gives rise, for every $f \in L^2(X_\Theta)^\infty$, to a $C^\infty(X_\Theta)$-valued measure $f^{\sigma} \nu(\sigma)$ on $\widehat {L_\Theta}$, by the Plancherel formula:
\begin{equation}\label{onaxis}
\left<f, \Phi\right>_{L^2(X_\Theta)} = \int_{\widehat{L_\Theta}} \int_{X_\Theta}  f^{\sigma}(x) \overline{\Phi(x)} dx  \nu(\sigma)
\end{equation}
(for all $\Phi\in C_c^\infty(X_\Theta)$).

For almost every $\sigma$, the function $f^{\sigma}$ belongs to $C^\infty(X_\Theta)^{\sigma}$. If $f\in C_c^\infty(X_\Theta)$, this measure has a natural ``translation'' to any translate of $\widehat{L_\Theta}$ (see the discussion following Lemma \ref{Plancherel-moderate}), and we have:
\begin{equation}\label{offaxis}
\left<f,\Phi\right>_{L^2(X_\Theta)} = \int_{\omega^{-1}\widehat{L_\Theta}} \int_{X_\Theta} f^{\sigma}(x) \overline{\Phi(x)} dx  \nu(\omega\sigma)
\end{equation}
for every character $\omega$ of $L_{\Theta,X}^\ab$.

\begin{theorem}\label{explicitsmooth}
 For any $\omega\gg 0$, if $f\in C_c^\infty(X_\Theta)$ admits the decomposition (\ref{offaxis}) then:
\begin{equation}
  e_\Theta f (x) = \int_{\omega^{-1}\widehat{L_\Theta}} E_{\Theta,\sigma} f^{\tilde\sigma}(x) \nu(\omega\sigma).
\end{equation}
\end{theorem}

Recall that in our notation $E_{\Theta,\sigma}: C^\infty(X_\Theta)^{\tilde\sigma}\to C^\infty(X)$; if $\sigma\in \omega^{-1}\widehat{L_\Theta}$ then $\tilde\sigma\in \omega\widehat{L_\Theta}$. We proceed with the proof of the theorem in several steps, including the explanation of \eqref{offaxis}.

\subsubsection{Moderate growth}

\begin{proposition*}\label{moderategrowth} 
 The image of $e_\Theta^*: C_c^\infty(X) \to C^\infty(X_\Theta)$ is a space $\mathcal F$ of functions satisfying the assumptions of \S \ref{Fproperties}; namely, for any open compact subgroup $J$ the $J$-invariants are of uniformly moderate growth, and their support has compact closure in an affine embedding of $X_\Theta$.
\end{proposition*}

\begin{proof}
The statement on the support is Proposition \ref{affinesupport}. For the moderate growth, { we may partition $X_\Theta$ into the union of subsets $N_\Omega$ belonging to $J$-good neighborhoods of $\Omega$-infinity, for all $\Omega\subset\Theta$, so that $N_\Omega$ is compact modulo $A_{X,\Omega}$. It is then the case that the functions $e_\Theta^* \Phi$, $\Phi\in C_c^\infty(X)^J$, are of uniformly moderate growth if and only if the same is true for the functions $e_\Omega^* \Phi|_{N_\Omega}$, $\Omega\subset\Theta$. Indeed, it is easy to see that ``moderate growth'' is compatible with our exponential map, i.e.\ the exponential map of \S \ref{subsec:expmap} between neighborhoods of $\Omega$-infinity of $X_\Theta/J$ and $X_\Omega/J$ carries functions of uniformly moderate growth in a neighborhood of a point of $\Omega$-infinity to functions of uniformly moderate growth.

Thus, assuming that uniformly moderate growth has been proven for all $\Omega\subsetneq\Theta$, it now suffices to prove that the restriction of $e_\Omega^*\Phi$, $\Phi\in C_c^\infty(X)^J$ on $A_{X,\Theta}$-orbits is of moderate growth. }

Thus, it suffices to show the following: given $x\in X_\Theta$ and an open compact $J\subset G$ there is a finite number of regular functions $\omega_i$ of $\AA_{X,\Theta}$ such that for all $\Phi\in C_c^\infty(X)^J$ we have: 
$$\left|e_\Theta^*\Phi(a\cdot x)\right| \ll \sum_i |\omega_i(a)|.$$ 
{
Indeed, the pairs $(U_i, F_i=\omega_i)$, where $U_i$ are the open-closed subsets $U_i=\{a\in A_{X,\Theta}| |\omega_i(a)|\ge |\omega_j(a)| \mbox{ for all }j\}$ provide a cover of $A_{X,\Theta}$ as in the definition of ``moderate growth'', showing that the pull-back of $e_\Theta^*\Phi$ to $A_{X,\Theta}$ under the action map is of moderate growth.

}

Using a Plancherel decomposition:
$$ \Vert\bullet\Vert^2 = \int_{\hat G} H_\pi \mu(\pi)$$
for $L^2(X)$, we get:
$$e_\Theta^*\Phi(a\cdot x) =\Vol(axJ)^{-1}\left<\Phi, e_\Theta 1_{a\cdot xJ}\right>_{L^2(X)} = \Vol(axJ)^{-1}\int_{\hat G} H_\pi(\Phi,e_\Theta 1_{a\cdot xJ}) \mu(\pi).$$
The sesquilinear forms:
\begin{equation}\label{tree} C_c^\infty(X)\otimes\overline{C_c^\infty(X_\Theta)} \ni \Phi\otimes \bar f\mapsto H_\pi(\Phi,e_\Theta f) \in \CC\end{equation}
are $A_{X,\Theta}$-finite with respect to the action of $A_{X,\Theta}$ on the second variable; this is because
they factor through the $\pi$-coinvariants, which are of finite length (by  finiteness of multiplicity, i.e. Theorem  \ref{finiteness}). 

Now recall (Corollary \ref{corollarysurjection})
 that for almost every $\pi$ there exists an $\Omega\subset \Delta_X$ such that $\pi$ is a relative discrete series for $X_\Omega$, and the conclusion of Proposition \ref{uniformbound} is satisfied: the absolute value of the exponents of any $\pi\to C^\infty(X)$ belong to a finite set of homomorphisms: $A_{X,\Theta}\to \mathbb R_+^\times$. By an analog of Lemma \ref{boundS}, this implies that there is a finite subset $\Lambda\subset A_{X,\Theta}$ and a finite set of characters $\omega_i$ of $\AA_{X,\Theta}$ such that: 
$$|H_\pi(\Phi,e_\Theta 1_{a\cdot xJ})| \le  \left(\sum_{\lambda\in \Lambda}|H_\pi(\Phi,e_\Theta 1_{\lambda\cdot xJ})|\right)\cdot \sum_i |\omega_i(a)|$$ 
for all $a\in A_{X,\Theta}$; thus:
$$\left| e_\Theta^*\Phi(a\cdot x)\right| \le \left(\sum_{\lambda\in \Lambda} |e_\Theta^*\Phi(\lambda\cdot x)|\right) \cdot \sum_i |\omega_i(a)|.$$
\end{proof}

\subsubsection{Plancherel decomposition of moderate growth functions with bounded support}

Fix a parabolic in the class $P_\Theta^-$, hence a subspace $X_\Theta^L$ of $X_\Theta$, and compatible measures so that:
$$\int_{X_\Theta} f(x) dx = \int_{P_\Theta^-\backslash G} \int_{X_\Theta^L} (g\cdot f)|_{X_\Theta^L} f(x) dx dg.$$
Also fix a point on $X_\Theta^L$ in order to consider characters of $L_{\Theta,X}^\ab$ as functions on $X_\Theta^L$, as before. 

Let us fix a Plancherel decomposition for $X_\Theta^L$:
 $$\int_{X_\Theta^L} f_1(x)f_2(x) dx = \int_{\widehat{L_\Theta}} H'_\sigma(f_1,f_2) \nu(\sigma)$$
for $X_\Theta^L$; we have written the inner product as a bilinear pairing, so the forms $H_\sigma$ are bilinear pairings of the spaces of coinvariants:
$$H'_\sigma: C_c^\infty(X_\Theta^L)_\sigma\otimes C_c^\infty(X_\Theta^L)_{\tilde\sigma}\to \CC.$$

We can ``twist'' the forms $H_\sigma$ to forms:
$$H'_{\omega\sigma}: C_c^\infty(X_\Theta^L)_{\omega\sigma}\otimes C_c^\infty(X_\Theta^L)_{\omega^{-1}\tilde\sigma}\to \CC,$$
for characters $\omega$ of $L_{\Theta,X}^\ab$ which are not necessarily unitary, simply by setting:
$$H'_{\omega\sigma} (f_1,f_2) = H_\sigma(\omega^{-1}f_1,\omega f_2).$$
This definition is consistent (for unitary $\omega$) if and only if the Plancherel measure $\nu$ chosen is $\widehat{L_{\Theta,X}^\ab}$-invariant (which we can assume).

The Plancherel formula for $L^2(X_\Theta)$ will involve the forms ``induced'' from the $H'_\sigma$:
 $$\int_{X_\Theta} f_1(x)f_2(x) dx = \int_{\widehat{L_\Theta}} H_\sigma(f_1,f_2) \nu(\sigma),$$
 $$H_\sigma (f_1,f_2) = \int_{P_\Theta^-\backslash G} H'_\sigma(g\cdot f_1|_{X_\Theta^L}, g\cdot f_2|_{X_\Theta^L}) dg.$$

\begin{lemma}\label{Plancherel-moderate}
Let $\mathcal F$ be a class of functions on $X_\Theta$ as before, then for $\omega\gg 0$, $f_1\in \mathcal F$ and $f_2\in C_c^\infty(X_\Theta)$ we have a ``Plancherel'' decomposition of the inner product:
\begin{equation}\label{offaxis-Hsigma}
 \int_{X_\Theta} f_1\cdot f_2 = \int_{\omega^{-1} \widehat{L_\Theta}} H_\sigma(f_1,f_2) \nu(\omega\sigma),
\end{equation}
for $\omega\gg0$. Notice that the image of $f_1$ in $C_c^\infty(X_\Theta)_\sigma$ makes sense for $\omega\gg 0$ by Corollary \ref{corollaryconvergence}, thus the right hand side is well-defined.
\end{lemma}

\begin{proof}
We compute:
$$ \int_{X_\Theta} f_1\cdot f_2 = \int_{P_\Theta^-\backslash G} \int_{X_\Theta^L} (g\cdot f_1)(x) (g\cdot f_2)(x) dx dg = $$ $$=\int_{P_\Theta^-\backslash G} \int_{X_\Theta^L} \omega(x)(g\cdot f_1)(x) \omega^{-1}(x)(g\cdot f_2)(x) dx dg.$$
For $\omega\gg 0$ both factors of the integrand are in $L^2(X_\Theta^L)$, so applying the Plancherel decomposition we get:
$$\int_{P_\Theta^-\backslash G} \int_{\widehat{L_\Theta}} H'_\sigma(\omega\cdot (g\cdot f_1), \omega^{-1}\cdot (g\cdot f_2))\nu(\sigma) dg=$$
$$=\int_{\widehat{L_\Theta}} \int_{P_\Theta^-\backslash G}  H'_\sigma(\omega\cdot (g\cdot f_1), \omega^{-1}\cdot (g\cdot f_2)) dg \nu(\sigma)=$$
$$ = \int_{\omega^{-1} \widehat{L_\Theta}} H_\sigma(f_1,f_2) \nu(\omega\sigma).$$
\end{proof}

An alternative way to state that is: 
Recall that $C^\infty(X_\Theta)^{\sigma}$ is the dual of $C_c^\infty(X_\Theta)_{\tilde\sigma}$, considered as a subspace of $C^\infty(X_\Theta)$. The Plancherel formula for $X_\Theta$ gives rise, for every $f \in C_c^\infty(X_\Theta)$, to a $C^\infty(X_\Theta)$-valued measure $f^{\sigma} \nu(\sigma)$ on $\widehat {L_\Theta}$, defined by \eqref{onaxis}.
For almost every $\sigma$, the function $f^{\sigma}$ belongs to $C^\infty(X_\Theta)^{\sigma}$. The above discussion shows that  the definition of this measure can be extended to translates of $\widehat {L_\Theta}$ by characters of $L_{\Theta,X}^\ab$: simply replace \eqref{onaxis} by the analogous expression coming from \eqref{offaxis-Hsigma} (valid for \emph{any} $\omega$, if we take $f_1$ and $f_2$ to be compactly supported).

Then the above result amounts to saying that the expression \eqref{offaxis} is valid for $f\in \mathcal F$, as long as $\omega\gg 0$. For such $\omega$, and $\sigma\in\omega^{-1}\widehat{L_\Theta}$, the map: $\mathcal F\to C_c^\infty(X_\Theta)_\sigma$ (or, equivalently, to $C^\infty(X_\Theta)^\sigma$) is defined by the ``convergent series'' extension of morphisms of Corollary \ref{corollaryconvergence}.

\subsubsection{Proof of Theorem \ref{explicitsmooth}}

Let $\Phi\in C_c^\infty(X)$. By Propositions \ref{affinesupport} and \ref{moderategrowth} we may apply Lemma \ref{Plancherel-moderate} to the function $e_\Theta^*\Phi$; if $f\in C_c^\infty(X_\Theta)$, we then get:
$$  \left<\Phi,\overline{e_\Theta f}\right>_{L^2(X)}  = \left<e_\Theta^*\Phi, \bar f\right>_{L^2(X_\Theta)} = \int_{\omega^{-1}\widehat{L_\Theta}} H_\sigma(e_\Theta^*\Phi,f) \nu(\omega\sigma).$$

By the third statement of Corollary \ref{corollaryconvergence} and the definition of the normalized Eisenstein integrals, the image of $e_\Theta^*\Phi$ in $C_c^\infty(X_\Theta)_\sigma$ is equal to the adjoint normalized Eisenstein integral $E_{\Theta,\sigma}^*\Phi$. Therefore:
$$  \left<\Phi,\overline{e_\Theta f}\right>_{L^2(X)}  = \int_{\omega^{-1}\widehat{L_\Theta}} H_\sigma(E_{\Theta,\sigma}^*\Phi,f) \nu(\omega\sigma) = $$
$$= \int_{\omega^{-1}\widehat{L_\Theta}} \int_{X_\Theta} (E_{\Theta,\sigma}^*\Phi)(x) f^{\tilde\sigma} (x)  dx \nu(\omega\sigma) 
$$
(by the definition of $f^{\tilde\sigma}$ in \eqref{onaxis})
$$
= \int_{\omega^{-1}\widehat{L_\Theta}} \int_X \Phi(x) (E_{\Theta,\sigma} f^{\tilde\sigma}) (x)  dx \nu(\omega\sigma) = $$
$$ =  \int_X \Phi(x) \int_{\omega^{-1}\widehat{L_\Theta}}(E_{\Theta,\sigma} f^{\tilde\sigma}) (x)\nu(\omega\sigma) dx.
$$
This proves Theorem \ref{explicitsmooth}. 
\qed

\subsection{The canonical quotient and the small Mackey restriction} 

We will now discuss certain technical prerequisites for the explicit decomposition of unitary asymptotics (i.e.\ the Bernstein maps). The basic issue here is the absence of a diagram analogous to \eqref{Radoncommutes2}; for this reason, Eisenstein integrals do not appear explicitly a priori, and their relevance has to be established via their properties -- more precisely, their asymptotics, and the notion of ``small Mackey restriction'' that we are about to define.

\subsubsection{The canonical quotient of an induced representation} \label{sss:canquot}
If $\tau$ is a smooth admissible representation of $L_{\Theta}$, and the intertwining operator $T_{\Theta}: I_{{\Theta}^-}^G \tau  \rightarrow  I_{{\Theta}}^G( \tau)' $ is regular at $\tau$, we can obtain a certain canonical quotient of  the (normalized) Jacquet module
of the (normalized) induced representation $I_{{\Theta}^-}^G \tau$ as the composition:
$$ I_{{\Theta}^-}^G (\tau)_{\Theta}  \xrightarrow{T_\Theta} I_{{\Theta}}^G( \tau)' _{\Theta} \rightarrow \tau'$$
where $\tau'$ denotes a representation that is isomorphic to $\tau$ once a measure on $U_{\Theta}$ is fixed.
We call this ``the canonical quotient,'' even though it is defined
only when the intertwining operator is regular.

Note that there is a canonical inclusion $\tau' \hookrightarrow I_{\Theta^-}^G(\tau)_{\Theta}$ (by considering those elements of $I_{\Theta^-}^G(\tau)$ which are supported on the open $P_\Theta$-orbit),
and when composed with the canonical quotient this gives the identity, i.e. the composite
\begin{equation} \label{composid}  \tau'  \hookrightarrow I_{\Theta^-}^G(\tau)_{\Theta} \rightarrow \tau' \end{equation}
is the identity.

In the case of $C^{\infty}_c(X_{\Theta})_{\sigma}$ (abstractly isomorphic to
an induced admissible representation; see \eqref{Vminus}) 
we denote 
the corresponding quotient by $C_c^\infty(X_\Theta)_\sigma[\sigma]$:
\begin{equation} \label{squaresigmadef} \left(C_c^\infty(X_\Theta)_\sigma\right)_\Theta \to C_c^\infty(X_\Theta)_\sigma[\sigma].\end{equation} 
Again, it is defined only for the set of $\sigma$ for which the intertwining operator
$I_{\Theta^-}(\sigma) \rightarrow I_{\Theta}(\sigma)'$ is an isomorphism; but we
have seen in Corollary \ref{cor:genericbehavior} that this includes $\nu_{\disc}$-almost every $\sigma$. We have a canonical isomorphism:
\begin{equation}\label{canonquotisom}
 C_c^\infty(X_\Theta)_\sigma[\sigma] = \left(\Hom_{L_\Theta}(C^{\infty}_c(X_\Theta^L),\sigma)\right)^* \otimes \sigma'.
\end{equation}

Later we shall use the following property of the canonical quotient. To state it, note that
the action of $\mathcal{Z}(L_{\Theta})$ on the flag variety $P_{\Theta}^- \backslash G$
induces an action of  $\mathcal{Z}(L_{\Theta})$ by equivalences on the functor $I_{\Theta^-}$. 
 In this way, we obtain an action of $\mathcal{Z}(L_{\Theta})$ on
 $I_{\Theta^-}(\tau)$. Therefore, the Jacquet module $I_{\Theta^-}(\tau)_\Theta$ has a $\mathcal Z(L_\Theta)\times L_\Theta$-action.

Consider the morphism
\begin{equation} \label{purple}
 I_{\Theta^-}(\tau)_{\Theta} \rightarrow \tau'. \end{equation}

\begin{lemma}\label{lemmainvariant} 
 For generic $\tau$, the antidiagonal copy $\mathcal Z(L_\Theta)\hookrightarrow \mathcal Z(L_\Theta)\times L_\Theta$ acts trivially on the quotient $\tau'$ of (\ref{purple}) in other words: the $\mathcal Z(L_\Theta)$-action on $I_{\Theta^-}(\tau)_{\Theta}$ commutes with the $\mathcal Z(L_\Theta)\hookrightarrow L_\Theta$-action on $\tau'$.
\end{lemma}

\begin{proof}
 This is just a consequence of the fact that (\ref{composid}) is the identity, with its first arrow being $\mathcal Z(L_\Theta)$-equivariant and the second being $L_\Theta$-equivariant. 
\end{proof}

\subsubsection{The small Mackey restriction} \label{sss:smM}

Fix a parabolic in the class $P_\Theta$, and recall \eqref{Mackeyrestriction} that for a representation $\pi$ of $G$, ``Mackey restriction'' is the natural morphism:
$$\Hom_G(C_c^\infty(X),\pi)\to \Hom_{L_\Theta}(C_c^\infty(X_\Theta^L)',\pi_\Theta)$$
obtained by applying the Jacquet functor and identifying $C_c^\infty(X_\Theta^L)'$ with a subspace of the Jacquet module of $C_c^\infty(X)$  as in \eqref{Jprecise}. 

We will use the term ``small Mackey restriction''
in the following situation: Suppose that $\pi$ is endowed with an isomorphism
to an induced representation $I_{\Theta^{-}}(\tau)$ and the intertwining operator 
$T_\Theta:I_{\Theta^-}(\tau) \rightarrow I_{\Theta}(\tau)'$ is an isomorphism. In that case, 
composition with the canonical quotient gives:
\begin{equation}\label{smallMackey}
\Mac: \Hom_G(C_c^\infty(X),\pi)\to \Hom_{L_\Theta}(C_c^\infty(X_\Theta^L)', \tau')=\Hom_{L_\Theta}(C_c^\infty(X_\Theta^L), \tau).
\end{equation}
and this morphism will be, by definition, the ``small Mackey restriction.''

We have proved that for almost all $\tau$ in a given twisting class  (\S \ref{sss:twistingclass})
the intertwining map $I_{\Theta^-}(\tau) \rightarrow I_{\Theta}(\tau)'$  is an isomorphism
(see  Corollary \ref{cor:genericbehavior}) and, therefore, the notion of ``small Mackey restriction''
makes sense for $\pi = I_{\Theta^-}(\tau)$ at least for $\tau$ generic in a twisting class. 

In particular, given an irreducible representation $\sigma$ of $L_\Theta$, and any morphism (not necessarily the canonical one):
$$M: C_c^\infty(X_\Theta)\to \pi := C_c^\infty(X_\Theta)_\sigma,$$
and taking into account that the right-hand side has the structure of an induced representation
with canonical quotient $\tau'=C_c^\infty(X_\Theta)_\sigma[\sigma]$ (see \eqref{squaresigmadef}) 
the small Mackey restriction of $M$ is a morphism:
\begin{equation} \label{li2eq} \Mac (M): C_c^\infty(X_\Theta^L)'  \to C_c^\infty(X_\Theta)_\sigma[\sigma]. \end{equation} This enjoys the following property:

\begin{lemma}
\label{lemmainvariant2} For any $G$-morphism: 
$$ M: C_c^\infty(X_\Theta)\to C_c^\infty(X_\Theta)_\sigma, $$ the small Mackey restriction (\ref{li2eq}) is $A_{X, \Theta}$-invariant. 
\end{lemma}
\begin{proof} 
In fact, we claim that the $A_{X, \Theta}$-action on source and target
of \eqref{li2eq} 
coincides with the restriction of the $L_{\Theta}$-action
to the center $\mathcal{Z}(L_{\Theta})$, by means of the surjection
$\mathcal{Z}(L_{\Theta}) \twoheadrightarrow A_{X, \Theta}$. 
That will prove the lemma, because \eqref{li2eq} is certainly $L_{\Theta}$-equivariant.

For the source, this is easy to see from the definitions.

For the target, this is Lemma \ref{lemmainvariant}.
\end{proof}
 
 \subsubsection{The work of Blanc and Delorme} \label{BDtranslated}  
 
We now translate a very useful result of Blanc and Delorme \cite{BlancDelorme} in the symmetric case,  into our language. 
 
 They prove that, if $\XX$ is a {\em symmetric} variety, the small Mackey restriction:
 \begin{equation}\label{bdvar} \Hom_G(C_c^\infty(X),I_{\Theta^-}(\sigma))\to \Hom_{L_\Theta}(C_c^\infty(X_\Theta^L),\sigma)
\end{equation}
is {\em injective}, generically for $\sigma$ within a twisting class.

The main theorem of Delorme and Blanc actually concerns
the composite 
 $$  \Hom(C^{\infty}_c(X), I_{\Theta}(\sigma))  \rightarrow
    \Hom_{L_{\Theta}}(C^{\infty}_c(X)_{\Theta}, 
   I_{\Theta}(\sigma)_{\Theta})    \rightarrow \Hom(C^{\infty}_c(X^L_{\Theta})', \sigma).$$ 
It asserts
that the composite map  is an isomorphism. In fact, their paper analyzes this in terms 
of distributions on the flag variety invariant by a point stabilizer on $X$, but this is easily translated 
to the stated claim using the isomorphism: $\Hom_G(C_c^\infty(X),\pi)\simeq \Hom_G(\tilde \pi,C^\infty(X))$ for admissible representations $\pi$. 
One passes to \eqref{bdvar} by applying the  intertwining operator.

\subsection{Unitary asymptotics (Bernstein maps)}

Our main theorem for unitary asymptotics is the following; its formulation uses the $C^\infty(X_\Theta)$-valued measure $f^\sigma\nu(\sigma)$ on $\widehat{L_\Theta}$, attached to $f\in L^2(X_\Theta)^\infty$ as in \eqref{onaxis}. This was also used in the explicit description of smooth asymptotics, Theorem \ref{explicitsmooth}, except that there we were restricting ourselves to $f\in C_c^\infty(X_\Theta)$, and here we will not need to translate off the unitary spectrum. We confine ourselves to describing the restriction of $\iota_\Theta$ to $L^2(X_\Theta)_\disc$, since the space $L^2(X)$ is built out of the images of those spaces; so, $\nu_\disc$ will denote the ``discrete'' part of a Plancherel measure for $L^2(X_\Theta^L)$. Restricting to discrete spectra will simplify somehow the proof of Theorem \ref{explicitMackey} in \S \ref{proofexplicitMackey}.

\begin{theorem}\label{explicitBernstein} 
 Assume  that for $\nu_\disc$-almost all $\sigma$ the small Mackey restriction map \eqref{smallMackey}:
$$\Hom_G(C_c^\infty(X),I_{\Theta^-}(\sigma))\to \Hom_{L_\Theta}(C_c^\infty(X_\Theta^L),\sigma)$$
is injective (cf.\ \S \ref{BDtranslated} for its validity for symmetric varieties).

Then for a function $f\in L^2(X_\Theta)_\disc^\infty$ with pointwise Plancherel decomposition:
\begin{equation}\label{sigmadecomp}
 f(x) = \int_{\widehat{L_\Theta}} f^\sigma(x)\nu_\disc(\sigma)
\end{equation}
(where $f^\sigma \in C^\infty(X_\Theta)^\sigma$), its image under the Bernstein morphism is given by:
\begin{equation}
 \iota_\Theta f(x) = \int_{\widehat{L_\Theta}} E_{\Theta,\sigma} f^\sigma(x) \nu(\sigma).
\end{equation}
\end{theorem}

The decomposition \eqref{sigmadecomp} is analogous to \eqref{pointwisedecomposition}, where the space $C^\infty(X_\Theta)^\sigma$ now denotes the dual of the space of $\sigma$-coinvariants $C_c^\infty(X_\Theta)_\sigma$. Thus, the normalized Eisenstein integral $E_{\Theta,\sigma}$ (the dual of $E_{\Theta,\sigma}^*$) maps $C^\infty(X_\Theta)^\sigma$ into $C^\infty(X)$.

In combination with the scattering theorem \ref{advancedscattering}, this implies the following explicit Plancherel decomposition; recall that under Theorem \ref{advancedscattering}, $L^2(X)$ is a direct sum of the spaces $L^2(X)_\Theta$, where $\Theta$ ranges over associate classes of subsets of $\Delta_X$ (that is, conjugacy classes of Levi subgroups of $\check G_X$).

\begin{theorem}\label{explicitPlancherel}
 Under the assumptions of Theorem \ref{explicitBernstein},   the norm on $L^2(X)_\Theta$ admits a decomposition:
\begin{equation}
 \Vert \Phi\Vert_\Theta^2 = \frac{1}{|W_X(\Theta,\Theta)|}\int_{\widehat{L_\Theta}} \Vert E_{\Theta,\sigma}^* \Phi\Vert^2_{\sigma} \nu_\disc(\sigma). 
\end{equation}
The measure and norms here are the discrete part of the Plancherel decomposition (\ref{Planchereltheta}) for $L^2(X_\Theta)$.
\end{theorem}

The real content of Theorem \ref{explicitBernstein} is an unconditional statement about the ``Mackey restrictions'' of the morphisms that decompose the Bernstein map $\iota_\Theta$ (the morphisms $\iota_{\Theta, \sigma}^*: C^{\infty}_c(X) \rightarrow \mathcal{H}_{\sigma}$
from \eqref{iota*def}). Notice that $\mathcal H_\sigma$ is canonically an induced representation from a unitary representation of $L_\Theta$ (since $L^2(X_\Theta) = I_{\Theta^-} L^2(X_\Theta^L)$). Therefore, its ``canonical quotient'': 
$$\mathcal H_\sigma^\infty\twoheadrightarrow H_\sigma^\infty[\sigma]$$
is well-defined, cf.\ \S \ref{sss:canquot}.

Moreover, since $\mathcal H_\sigma$ is a completion of the $\sigma$-coinvariants $C_c^\infty(X_\Theta)_\sigma$, we have a commuting diagram of canonical quotients:

\begin{equation}\label{Budapest}
 \xymatrix{
C_c^\infty(X_\Theta)_\sigma \ar[d]\ar[r] & C_c^\infty(X_\Theta)[\sigma] \ar[d]^{\rm compl}\\ 
\mathcal H_\sigma^\infty \ar[r] & \mathcal H_\sigma^\infty[\sigma].
}
\end{equation}
The label ``compl'' denotes ``completion''.

For the maps $\iota_{\Theta, \sigma}^*$ and $E_{\Theta, \sigma}^*$ we have the notion of ``small Mackey restriction'', with images, respectively, in $C_c^\infty(X_\Theta)[\sigma]$ and $\mathcal H_\sigma^\infty[\sigma]$. 
Our (unconditional) assertion is that these two maps coincide after completion:
 
\begin{theorem}\label{explicitMackey} 
 The small Mackey restrictions of both $\iota_{\Theta,\sigma}^*$ and $E_{\Theta,\sigma}^*$ coincide (for $\nu_\disc$-almost every $\sigma$) after composing the latter with the map ``compl'' of \eqref{Budapest}.
\end{theorem} 

As a prelude to the proof, let us actually compute this small Mackey restriction 
for $E_{\Theta, \sigma}^*$: Recalling that the canonical
quotient of $C^{\infty}_c(X_{\Theta})_{\sigma}$ is identified with:
$$\left(\Hom_{L_\Theta}(C^{\infty}_c(X_\Theta^L),\sigma)\right)^* \otimes \sigma'$$
(cf.\ \eqref{canonquotisom}), the small Mackey restriction of $E_{\Theta, \sigma}^*$ is the  {\em natural projection}:
\begin{equation}
C_c^\infty(X_\Theta^L)
 \rightarrow  \left({\Hom_{L_\Theta}(C^{\infty}_c(X_\Theta^L),\sigma )}\right)^*\otimes \sigma.
 \end{equation}

 This is just obtained by tracing the definitions. An equivalent formulation
is the following: The composite of the maps defining the small Mackey restriction of $E_{\Theta, \sigma}^*$: \begin{equation}\label{tennis}
  \begin{CD} C_c^\infty(X_\Theta^L)'
@>>> C_c^\infty(X_\Theta)_\Theta @>>>  \left(C_c^\infty(X_\Theta)_\sigma\right)_\Theta @>>> C_c^\infty(X_\Theta)_\sigma[\sigma],  \end{CD}
\end{equation} 
(we recall that the first arrow is the canonical embedding defining Mackey restriction, the second is induced by $E_{\Theta, \sigma}^*$, and the third is the canonical quotient) coincides with the natural projection; this comes down to \eqref{composid}.

\subsubsection{Proofs}\label{proofexplicitMackey}

Our goal here is to prove Theorem \ref{explicitMackey}, from which Theorem \ref{explicitBernstein} follows by a formal argument (which we omit).

  Thus, the statement of Theorem \ref{explicitMackey} is that the compositions of the following arrows coincide:
\begin{equation}\label{wheel1} 
 C_c^\infty(X_\Theta^L)' \xrightarrow{\sim} C_c^\infty(\mathring X P_\Theta)_\Theta \hookrightarrow C_c^\infty(X)_\Theta \overset{\iota^*_{\Theta,\sigma}}{\underset{{\rm compl}\circ E_{\Theta,\sigma}^*}{\rightrightarrows}} (\mathcal H^\infty_\sigma)_\Theta \to \mathcal H^\infty_\sigma[\sigma].
\end{equation}
Recall that the index ``$\Theta$'' denotes Jacquet module with respect to $P_\Theta$.
Here, the maps we have denoted $\iota_{\Theta, \sigma}^*$ and $E_{\Theta, \sigma}^*$
are more precisely the Jacquet functors applied to these maps.

These morphisms are, a priori, just $L_\Theta$-equivariant. However, Lemma 
\ref{lemmainvariant2}  implies that their composition is also $A_{X,\Theta}$-equivariant. In what follows, we denote by $\left< \,\, , \,\, \right>_\sigma$ the hermitian inner product on $\mathcal H_\sigma$.

\begin{lemma} \label{Bilinear Lemma}
  Consider the bilinear form on $C_c^\infty(X_\Theta)$ given by the formula:
\begin{equation} \label{pairingscattering}
 (f_1,\overline{f_2})\mapsto \left<f_1,\iota_\Theta^* e_\Theta f_2\right>_\sigma.
\end{equation}
 It carries a finite\footnote{Indeed, in {\em both} arguments, it factors through the quotient $C^{\infty}_c(X_{\Theta})_{I_{\Theta^-}\sigma}$. } diagonal action of $A_{X,\Theta}$;  the $A_{X,\Theta}$-invariant part\footnote{i.e., the $A_{X, \Theta}$-equivariant projection to the generalized eigenspace with eigenvalue 1 which, a posteriori, is $A_{X, \Theta}$-invariant for almost all $\sigma$} of this pairing is equal (for $\nu_\disc$-almost all $\sigma$) to the pairing:
\begin{equation} \label{pairingstandard}
 (f_1,f_2)\mapsto \left<f_1,f_2\right>_\sigma.
\end{equation}
\end{lemma}

Here is an equivalent phrasing:  the quotient map:  $$ C^{\infty}_c(X_{\Theta})   \to C^{\infty}_c(X_{\Theta})_\sigma \to \mathcal H^\infty_\sigma$$
coincides with the $A_{X, \Theta}$-invariant part of 
$$ C^{\infty}_c(X_{\Theta})   \ni \Phi \mapsto   \mbox{ the image of $\iota_{\Theta}^* e_{\Theta} \Phi$ in } \mathcal H_\sigma  .$$

\begin{proof}
Write $\tau  = I_{\Theta^-}(\sigma)$.  Recall that $\tau$ is irreducible for
almost every $\sigma$. Let $\Pi$ be  the natural projection $C^{\infty}_c(X_{\Theta})_{\tau}
\rightarrow C^{\infty}_c(X_{\Theta})_{\sigma}$ -- indeed
the latter space is  $\tau$-isotypical by construction, and therefore a quotient of the former. 
 
 Let $f_1,f_2\in C_c^\infty(X_\Theta)$.  
\begin{eqnarray} \nonumber   \left<f_1, \iota_\Theta^* e_\Theta f_2\right>_\sigma  =  \left< f_1, \iota_\Theta^* \iota_\Theta f_2\right>_\sigma + \left<  f_1,  \iota_{\Theta}^* (e_{\Theta} - \iota_{\Theta}) f_2\right>_\sigma \\    \label{rvw} = \left<f_1, \left(\iota_\Theta^* \iota_\Theta f_2\right)_\disc\right>_\sigma + \left< f_1,  \iota_{\Theta}^* (e_{\Theta} - \iota_{\Theta}) f_2\right>_\sigma.
 \end{eqnarray}
since $\left< \,\, , \,\, \right>_\sigma$ is a Plancherel Hermitian form for the discrete spectrum of $X_{\Theta}$. 

Let us examine the second term. 
Fix a Plancherel formula for $L^2(X)$:
$$L^2(X) = \int_{\hat G} \mathcal H_\pi\mu(\pi)$$
and the corresponding Plancherel formula for $X_{\Theta}$ according to Theorem \ref{bernstein-abstract}:
$$L^2(X_\Theta) = \int_{\hat G} \mathcal H_\pi^\Theta \mu(\pi).$$
Then, by the uniqueness of Plancherel decomposition, $\left< \,\, , \,\, \right>_\sigma$
factors (for $\nu_\disc$-almost every $\sigma$) through a norm on $\mathcal{H}_{\tau}^{\Theta}$.  Consequently, we
may re-express this second term as a linear function of
\begin{equation} \label{TenZor}  \left(  \mbox{image of  $f_1$ in  $\mathcal{H}^{\Theta}_{\tau}$ } \right)  \otimes 
 \left( \mbox{ image of $ \iota_{\Theta}^* (e_{\Theta} - \iota_{\Theta}) f_2$ in $\mathcal{H}^{\Theta}_{\tau}$}\right).\end{equation}

 Now, in \eqref{TenZor}, the action of $A_{X, \Theta}$ on the first argument here is through unitary exponents
(being a unitary action on a Hilbert space).
 
   On the other hand, the map
$$f_2 \mapsto \mbox{ image of } \iota_{\Theta}^* (e_{\Theta} - \iota_{\Theta}) f_2 \in \mathcal{H}^{\Theta}_{\tau}$$
is given (for almost all $\tau$) by the composite
$ \iota_{\Theta, \tau}^*  \circ (e_{\Theta, \tau} - \iota_{\Theta, \tau} ) $
where we regard $e_{\Theta, \tau}$ and $\iota_{\Theta, \tau}$
as mappings $C^{\infty}_c(X_{\Theta})_{\tau} \rightarrow  \mathcal{H}_{\tau}$ 
and we regard $\iota_{\Theta, \tau}^*$ as a mapping $\mathcal{H}_{\tau} \rightarrow \mathcal{H}_{\tau}^{\Theta}$.  By \eqref{overcomplicated}, for every $a \in A_{X, \Theta}^+$
we have
$$ \| (e_{\Theta, \tau} - \iota_{\Theta, \tau}) (a^n \cdot  f_2)\|_{\mathcal{H}_{\tau}} \rightarrow 0,$$
for $f_2 \in C^{\infty}_c(X_{\Theta})^J$.  Since $\iota_{\Theta, \tau}^*$ is bounded,
we have a similar property for  $\iota_{\Theta, \tau}^* (e_{\Theta, \tau} - \iota_{\Theta, \tau}) (a^n \cdot  f_2)$: it converges to zero inside $\mathcal{H}_{\tau}^{\Theta}$.  That shows the action of $A_{X, \Theta}$ on the second argument
must be through non-unitary exponents. 
 
We have now established that the $A_{X, \Theta}$-invariant part of the second term of \eqref{rvw}
is zero. 

Consider the first term of \eqref{rvw}.   
We may decompose according to Proposition \ref{decompf}
$$ \left(\iota_\Theta^* \iota_\Theta f_2\right)_\disc = \sum_{i}  S_i (f_2)_\disc,$$
where the morphisms $S_i$ are equivariant with respect to isogenies $T_i: \AA_{X, \Theta} \rightarrow \AA_{X, \Theta}$,  cf.\ Proposition \ref{decompf}. Note that, 
by the footnote on page \pageref{discisok}, one could replace
$(\iota_{\Theta}^* \iota_{\Theta} f_2)_{\disc}$ by
$\iota_{\Theta, \disc}^* \iota_{\Theta, \disc} f_{2, \disc}$
and therefore we only need the statement for the discrete spectrum. 

Therefore:
$$\left< f_1, \iota_\Theta^* \iota_\Theta f_2\right>_\sigma = \sum_i \left< f_1 , S_i (f_2)_\disc\right>_\sigma.$$
Now $S_{\mathrm{id}}$ is equal to the identity. 
Although that certainly follows from Theorem \ref{advancedscattering} when generic injectivity is known, this does not require generic injectivity;  it is a direct consequence
of Proposition \ref{Bmapisometry}, using a decomposition as in   the start of \S \ref{tilingtheoremIIproof}. 

 Therefore, the $A_{X,\Theta}$-invariant summand of this pairing is equal to: 
 $$\left< f_1 ,(f_2)_\disc\right>_\sigma = \left< f_1 ,f_2\right>_\sigma.$$
\end{proof}

\begin{corollary}
Given a choice of pararabolic in the class $P_\Theta$ the two $P_\Theta$-equivariant maps: $C_c^\infty(X_\Theta^L)' \to C_c^\infty(X_\Theta)_\sigma[\sigma]$ obtained from following diagram are identical, for almost all
  irreducible representations $\sigma$ of $L_{\Theta}$ within a fixed twisting class:
\begin{equation} \label{mther}
\xymatrix{
& C_c^\infty(X)_\Theta  \ar[dr]^{\iota^*_{\Theta,\sigma}}  \\
C_c^\infty(X_\Theta^L)' \ar[r] & C_c^\infty(X_\Theta)_\Theta \ar[u]_{e_\Theta}\ar[r] & \left(\mathcal H^\infty_\sigma\right)_\Theta \ar[r]& \mathcal H^\infty_\sigma[\sigma].
}
\end{equation}
Here the arrows on the horizontal row are as follows: the first arrow is the canonical embedding defining Mackey restriction, the second arrow is induced by passing to $\sigma$-coinvariants and completing, and the third is the canonical quotient.
\end{corollary}

Again, in the above diagram, the arrows denoted $e_{\Theta}$ and $\iota_{\Theta, \sigma}^*$
are more precisely the Jacquet functors applied to those morphisms; we do not
denote this explicitly for typographical reasons.

\begin{proof}

 We will use Lemma \ref{lemmainvariant2}.  

It follows from Lemma \ref{Bilinear Lemma} that the quotient map: $C_c^\infty(X_\Theta)\to \mathcal H^\infty_\sigma$ is simply the $A_{X,\Theta}$-equivariant part of the map $\iota_{\Theta,\sigma}^*\circ e_\Theta$.  In other words, the ``upper'' and ``lower'' composite
$$ C_c^\infty(X_\Theta)_\Theta  \rightarrow \left(\mathcal H^\infty_\sigma\right)_\Theta$$ 
both have the same $A_{X, \Theta}$-invariant part.

But Lemma \ref{lemmainvariant2} implies that the composite map from $C^{\infty}_c(X^L_{\Theta})'$  (either ``lower'' or ``upper'') 
to $\mathcal H^\infty_{\sigma}[\sigma]$ is $A_{X, \Theta}$-invariant. 
 
 Therefore the difference of the two maps in \eqref{mther} is, on the one hand, 
   $A_{X, \Theta}$-invariant; on the other hand, its $A_{X, \Theta}$-invariant part is zero. 
 Therefore this difference is zero, i.e., the two maps of \eqref{mther} coincide. 
\end{proof}

Finally, we recall by Proposition \ref{propositionMackey} that the Mackey embeddings: $C_c^\infty(X_\Theta^L)' \hookrightarrow C_c^\infty(X_\Theta)_\Theta$, $C_c^\infty(X_\Theta^L)'\hookrightarrow C_c^\infty(X)_\Theta$ commute with the map $e_\Theta$. Hence, altogether we get a diagram of Jacquet modules:
\begin{equation}
\xymatrix{
C_c^\infty(X_\Theta^L)' \ar[r] \ar[dr] & C_c^\infty(X)_\Theta  \ar[dr]^{\iota^*_{\Theta,\sigma}}  \\
& C_c^\infty(X_\Theta)_\Theta \ar[u]_{e_\Theta}\ar[r] & \left(\mathcal H^\infty_\sigma\right)_\Theta \ar[r]& \mathcal H^\infty_\sigma[\sigma],
}
\end{equation}
where the composed morphisms: $C_c^\infty(X_\Theta^L)' \to \mathcal H^\infty_\sigma[\sigma]$ agree. The ``lower'' map is the Mackey restriction of $E_\Theta^*$ composed with completion, by the discussion of \eqref{tennis};  the upper map is the Mackey restriction of $\iota_{\Theta,\sigma}^*$; their agreement is the assertion of Theorem \ref{explicitMackey}.

\subsection{The group case}

The case when $X=H$, a (split) connected reductive group under the action of $G=H\times H$, is a multiplicity-free  (and symmetric) example, therefore Theorem \ref{explicitBernstein} holds. The more familiar form of the Plancherel formula in this case is obtained by relating normalized Eisenstein integrals to the duals of matrix coefficients:

Let $L$ be a Levi subgroup of $H$; the corresponding boundary degeneration $X_\Theta$ is isomorphic to $L\backslash (U\backslash H\times U^-\backslash L)$, where $P=LU$ and $P^-=LU^-$ are two opposite parabolics with Levi $L$. 

Let  $\tau$ be an irreducible representation of $L$, and $\sigma := \tau\otimes \tilde\tau$, a representation of $L_\Theta=L\times L$. The matrix coefficient map is $M_\tau:\sigma\otimes \delta_{P_\Theta}^{-\frac{1}{2}}\to C^\infty(L)$, (where $P_\Theta=P^-\times P$), and by applying the induction functor we get a $G$-morphism:
\begin{equation}
 I_{P\times P^-}(\sigma) \xrightarrow{I_{\Theta^-}M_\tau} C^\infty(X_\Theta).
\end{equation}
The image lies in what was previously denoted by $C^\infty(X_\Theta)^\sigma$. Finally, we may apply the normalized Eisenstein integral $E_{\Theta,\tilde\sigma}$ to this to get a map into $C^\infty(X)$:
\begin{equation}\label{verdi}
 I_{P\times P^-}(\sigma) \xrightarrow{E_{\Theta,\tilde\sigma}\circ M_\tau} C^\infty(H).
\end{equation}

On the other hand, we may choose an invariant measure on the suitable line bundle over $P\backslash G$ in order to identify the representation $\rho=I_P(\tilde\tau)$ with the dual of $I_P(\tau)$, and then we have a matrix coefficient map
\begin{equation}\label{puccini}
 M_\rho: I_P(\tau)\otimes I_P(\tilde \tau)\to C^\infty(H).
\end{equation}

The question is what is the relationship between (\ref{verdi}) and (\ref{puccini}). In order to formulate the answer, let $T_2$ denote the standard intertwining operator in the second factor (similarly, $T_1$ will denote the corresponding operator in the first factor):
$$I_{P\times P}(\sigma) \to I_{P\times P^-}(\sigma).$$
To define it, \label{keepthislabel} we use the measure on $U^-$ which corresponds to the chosen measure on $P\backslash G$, i.e.\ such that if $f\in \Ind_P^G(\delta_P)$ then:
$$\int_{P\backslash G} f(g) dg = \int_{U^-} f(u) du.$$

\begin{lemma}
 The map (\ref{verdi}) is the composition of (\ref{puccini}) with $T_2^{-1}$.
\end{lemma}

\begin{proof}
We add to the picture a third map,
\begin{equation}
 I_{P^-\times P}(\sigma) \xrightarrow{I_{\Theta}M_\tau} C^\infty(X_\Theta^h),
\end{equation}
arising as well from induction of matrix coefficients. Note that the space $C^\infty(X_\Theta^h)$ is dual to $C_c^\infty(X_\Theta^h,\delta_\Theta)$, non-canonically. Everything is very explicit here, and isomorphisms can be chosen compatibly, so we will not worry about distinguishing $I_\Theta(\sigma)$ from $I_\Theta(\sigma)'$, $C^\infty(X_\Theta^h)$ from the dual of $C_c^\infty(X_\Theta^h,\delta_\Theta)$, etc, leaving the details to the reader.

Following the definitions, and using the fact that the adjoint of a standard intertwining operator $T:I_P(\tau)\to I_{P^-}(\tau)'$ is again $T$, this time as a map $I_{P^-}(\tilde\tau)\to I_P(\tilde\tau)'$, one can see that the following diagram commutes:
\begin{equation}\label{cooldiagram}
\xymatrix{
C^\infty(H)  											& I_P(\tau)\otimes I_{P}(\tilde\tau) \ar[l]_{M_\rho}\\
C^\infty(X_\Theta^h)^\sigma \ar[u]^{R_{\Theta,\sigma}^*} 				& I_{P^-\times P}(\tau\otimes \tilde\tau)  \ar[l]_{I_\Theta M_\tau} \ar[u]_{T_1} \\
C^\infty(X_\Theta)^\sigma \ar@/^3pc/[uu]^{E_{\Theta,\tilde\sigma}} \ar[u]^{T_\Theta^{-1}} 	& I_{P\times P^-} (\tau\otimes\tilde\tau) \ar[l]_{I_{\Theta^-} M_\tau} \ar[u]_{T_\Theta^{-1}}.
}
\end{equation}

Composing $T_1$ with $T_\Theta^{-1}$ we get $T_2^{-1}$, which implies the claim of the lemma.
\end{proof}

The compatibility of measures on $X=H$ and $X_\Theta$ can be expressed as follows: one chooses measures on $U^-,L$ and $U$ such that $d(u^-) dl du $, as a measure on the open Bruhat cell $U^-LU$, is equal to the measure on $H$. Then one defines a measure on $X_\Theta=L\times{^{P\times P^-}} (H\times H)$ by considering the measures on $U^-,U$ as measures on $P\backslash G,P^-\backslash G$, respectively, and integrating the measure on $L$ over $P\backslash G\times P^-\backslash G$. Given the Plancherel formula for $L^2(L)$ (with respect to this measure):
\begin{equation}
 \Vert f\Vert^2 = \int_{\hat L} \Vert M_\tau^* f\Vert^2_\tau \nu(\tau),
\end{equation}
we get a Plancherel formula for $X_\Theta$:
\begin{equation}
 \Vert f \Vert^2 = \int_{\hat L} \Vert f\Vert'^2_\tau \nu(\tau),
\end{equation}
where the norms $\Vert f\Vert'^2$ are obtained from the dual of the induced matrix coefficient $I_{\Theta^-}M_\tau$ (the last horizontal arrow of diagram (\ref{cooldiagram})) and the unitary structure on $I_{P\times P^-}(\tilde \tau\otimes \tau)$ obtained from the unitary structure on $\tau$ and the given measures on $U^-,U$.

For those fixed measures, let $c(\tau)$ denote the constant which makes the following diagram commute:
\begin{equation}\label{ctau}\begin{CD} I_{P^-}^H(\tau)\otimes I_{P^-}^H(\tilde\tau) @>{T\otimes T}>> I_P^H(\tau)\otimes I_P^H(\tilde\tau) \\
   @V{C}VV @VV{C}V \\
\CC @>{\cdot c(\tau)}>> \CC.
  \end{CD}
  \end{equation}

Then from Theorem \ref{explicitPlancherel} we deduce the Plancherel formula for the group up to discrete Plancherel measures:

\begin{theorem} \label{groupPlancherel}
 There is a direct sum decomposition: $L^2(H)\simeq \oplus_{L/\sim} L^2(H)_L$, where the sum is taken over conjugacy classes of Levi subgroups, and a Plancherel decomposition for $L^2(H)_L$:
$$\Vert\Phi\Vert^2= \int_{\hat L_\disc/W(L,L)} \Vert \Phi\Vert_{I_P(\tau)}^2 c(\tau)^{-1} \nu_\disc(\tau).$$
Here the measure $\nu_\disc$ is the above Plancherel measure for $L^2(L)_\disc$, and the norm $\Vert\Phi\Vert_{I_P(\tau)}$ is the Hilbert-Schmidt norm of $\Phi dg$ acting by convolution on $I_P(\tau)$, for any parabolic $P$ with Levi subgroup $L$.
\end{theorem}

Notice that this Hilbert-Schmidt norm is precisely the norm obtained by the adjoint of $M_\rho$ (first horizontal arrow of the diagram (\ref{cooldiagram}) and the unitary structure for $I_P(\tilde\tau)\otimes I_P(\tau)$.

\part{Conjectures} 
 
\section{The local $X$-distinguished spectrum} \label{sec:localconj}

Let $K$ be a number field. In principle, the discussion of this section should be valid for global function fields, but since we will appeal to results of the rest of this paper, which used theorems on the structure of spherical varieties that have been proven only in characteristic zero, we restrict ourselves to number fields. We adopt the following notation: for an algebraic group $\GG$ over $K$, we denote by
$[\GG]$ the adelic quotient $\GG(K) \backslash \GG(\adele)$. We keep assuming that $\GG$ is split over the global or the local field, unless otherwise stated.

\subsection{Recollection of the Arthur conjectures\cite{Arthur-unipotent}} \label{Artur} 

To each local field $k$ (resp.\ global field $K$) one associates a locally compact group $\mathcal L_k$ (resp. $\mathcal L_K$) (the ``Langlands group''), together with morphisms that fit into the following diagram:

\begin{equation} \begin{CD}
\mathcal L_{K_v} @>>> \mathcal W_{K_v} @>>> \RR_{>0} \\
@VVV @VVV @VVV \\ 
\mathcal{L}_K @>>> \mathcal W_K @>>> \RR_{>0}
\end{CD} \end{equation}
where $\mathcal{W}$ denotes the Weil group. 
If $k$ is nonarchimedean, the image of the map $\mathcal{L}_k \rightarrow \RR_{>0}$ takes values in $q_k^{\Z}$, where $q_k$
is the cardinality of the residue field of $k$. 
There is as yet no fully satisfactory definition
of $\mathcal{L}_K$ in the case of a number field; we use it primarily for motivational purposes. 

  For a summary of these these conjectural $L$-groups  we refer to \cite{ArthurAutomorphic}.
In particular, in the case of a local field, $\mathcal{L}_k$ can be taken to be the Weil group
of $k$ in the archimedean case, and its product with $\mathrm{SU}_2$ in the nonarchimedean case.

Functoriality implies the following property of $\mathcal{L}_K$:
\begin{equation} \label{fzd} \mbox{Frobenius elements are dense.} \end{equation}
By this we mean the following: Given any morphism $\varphi: \mathcal{L}_K
\rightarrow \mathrm{GL}_m(\C)$,  the associated morphism
$\mathcal{L}_{K_v} \rightarrow \mathcal{L}_K \rightarrow \GL_m(\C)$ factors through 
$q_v^{\Z}$ for almost all places $v$; the image of a generator
is a {\em Frobenius at $v$}, and these are necessarily Zariski dense
in $\mathrm{image}(\varphi)$.

\begin{enumerate} 
\item 
A {\em local},  resp.\ {\em global Arthur parameter} is a homomorphism
$\psi: \mathcal L_k \times \SL_2(\CC) \rightarrow \check{G}$, resp. $\psi: \mathcal L_K  \times \SL_2(\CC) \rightarrow \check{G}$,
so that the restriction to $\mathcal L_k$ ($\mathcal L_K$) {\em has bounded image} and the restriction to $\SL_2$ is \emph{algebraic}. Let us call the restriction of the Arthur parameter to $\SL_2(\CC)$ the \emph{type} or \emph{$\SL_2$-type} of the Arthur parameter.

Given an Arthur parameter, its composition with the morphism $\mathcal L_k \to \mathcal L_k\times \SL_2(\CC)$: $w \mapsto \left(w \times \left(\begin{array}{cc} |w|^\frac{1}{2} & 0 \\ 0 & |w|^{-\frac{1}{2}} \end{array} \right)\right)$ defines a Langlands parameter, which we shall refer to as the {\em associated} Langlands parameter.

\item Conjecturally, to each $\check G$-conjugacy class of local Arthur parameters $[\psi]$ we may ``naturally'' associate a finite set of unitary representations  
of $G(k)$, the {\em Arthur packet} of $\psi$. These representations should all have the same {\em infinitesimal character} (using the notation of \cite[\S 8]{Vo}: in the nonarchimedean case two representations $\pi, \pi'$ are said to have the same infinitesimal character
if and only if the restrictions of their Langlands parameters to the Weil group coincide)  and contain the $L$-packet of the associated Langlands parameter. They behave as expected with respect to parabolic induction: if the parameter $\psi$ has image in a Levi $\check M$ of the dual group, the associated packet consists of the irreducible summands of (unitarily) parabolically induced representations from the corresponding packet of $M$ (where $M$ is a Levi subgroup of $G$ whose conjugacy class corresponds to the conjugacy class of $\check M$) \cite[p.44]{Arthur-unipotent}.
 Notice that local Arthur packets are not, in general, mutually disjoint.

\item  \label{pagerefitem3} Over a local non-archimedean field $k$ fix a hyperspecial maximal compact subgroup $G_0\subset G$, if such exists. By the theory of principal series and the Satake transform, different isomorphism classes of $G_0$-unramified (for short: ``unramified'') representations have different infinitesimal characters. Therefore, every Arthur packet contains at most one unramified representation. 

In the reverse direction, suppose that $\psi_1, \psi_2$ are Arthur parameters
whose associated packets contain the {\em same} unramified representation. Then 
$\psi_1 | \SL_2$ and $\psi_2 | \SL_2$ are conjugate:

In fact, set $\alpha_i $ to be the restriction of $\psi_i$ to ${\rm G}_m \subset \SL_2$.  We may suppose that $\alpha_i({\rm G}_m) \subset A^*$. It suffices to check
that the derivative $d\alpha_1$ is conjugate to $d\alpha_2$.  
The assertion about  ``infinitesimal
character'' in $A$-packets  shows that there exists {\em bounded} elements $g_i \in \check G$ (i.e.\ elements spanning relatively compact subgroups) 
so that $g_1  \alpha_1(q^{1/2}) \mbox{ is conjugate to } g_2 \alpha_2(q^{1/2}).$

Let $W$ be the Weyl group of $A^* \subset \check{G}$ , and $\liechecka := \varchi(A^*)^* \otimes \mathbb{R}$. 
There is a natural projection $\eig: \check{G} \rightarrow \liechecka/W$:
 it is the unique conjugacy-invariant continuous assignment that coincides
 with the natural projection $A^* \stackrel{H}{\rightarrow} \liechecka \rightarrow \liechecka/W$, where $H$ is the ``logarithm map'' characterized by
$$|\alpha(t)| = e^{\langle \alpha,  H(t) \rangle}, 
\ \   \alpha \in \varchi (A^*), \ \ t \in A^*.$$
 Note also that if $g_1,s$ are semisimple commuting elements and
 $g_1^{\mathbb{Z}}$ is relatively compact, then $\eig(g_1 s) = \eig(s)$;
 this follows from the equality
$H(g) = \frac{H(g^n)}{n}$ for elements of $A^*$.   We conclude that $\eig (\alpha_1(q^{1/2})) = \eig(\alpha_2(q^{1/2}))$, whence the conclusion.

\item To every $\check G$-conjugacy class of global Arthur parameters  $[\psi]$ one should be able to associate a subspace $\mathcal A_{[\psi]}$ of the space of automorphic forms, such that 
$$L^2(\G(K) \backslash \G(\adele)) = \int \mathcal A_{[\psi]} \mu(\psi),$$
where the measure class of the above direct integral is   the natural measure class\footnote{An Arthur parameter $\psi$ can be twisted by parameters into the centralizer of its image. Those form a locally compact abelian group, and the natural measure class on the orbit of $\psi$ is the class of Haar measure.} on the set of conjugacy classes of Arthur parameters.

We note that Arthur did not phrase the conjectures in terms of the subspaces $A_{[\psi]}$; however, other formulations e.g. \cite{GanGurevich}
phrase the global conjecture in such a fashion, and the work of V. Lafforgue (see \cite[\S 2.2]{VL}, which summarizes results proved in \cite{VL2}) gives further evidence
for it over a global function field. 

 By means of Langlands' work on the spectral decomposition, this conjecture is equivalent to a description of the discrete spectrum (modulo center):  Fix a unitary central idele class character $\Omega$ and consider the set of conjugacy classes of Arthur parameters which correspond to this idele class character and, moreover, their image does not lie in any proper Levi subgroup of $\check G$. Then we should have:
\begin{equation} \label{globaldecomposition}  L^2(\G(K) \ZZ(\adele) \backslash \G(\adele), \Omega)_{\disc} = \oplus L^2_{[\psi]}\end{equation} 
where $L^2_{[\psi]}=\mathcal A_{[\psi]}$, with $[\psi]$ ranging over these classes. (For simplicity we will drop the brackets from $[\psi]$ from now on.)

The spaces $\mathcal A_{[\psi]}$ have the following properties: each irreducible subquotient $\pi \subset \mathcal A_{[\psi]}$ factors as a restricted tensor product $\otimes'_v \pi_v$, and $\pi_v$ belongs to the local Arthur packet
associated to the pullback $\psi_v$ of $\psi$ to $\mathcal L_{K_v}$.   For almost all $v$, $\pi_v$ is the unramified representation corresponding to the associated Langlands parameter of $\psi_v$. 
\end{enumerate}

Of course, the above is by no means a description of all the properties that the Arthur packets are expected to have; just those which we will use in this paper.

\subsection{The conjecture on the local spectrum (weak form)}

Let us recall that, to every quasi-affine spherical variety without roots of type N, we have associated
a distinguished class of morphisms $f: \check G_X \times \SL(2) \rightarrow \check{G}$;
the restriction $f|\SL(2)$ is a principal Levi for $\check L(X)$, and 
$f(\check G_X)$ is in the conjugacy class of the Gaitsgory-Nadler dual group $\check G_{X, GN}$.

\begin{definition}
 An $\XX$-distinguished Arthur parameter is a commutative diagram of the form: 
\begin{equation}\label{Xparameter}
\xymatrix{
& \check G_X \times \SL_2 \ar[dr]  \\
\mathcal L_k \times \SL_2 \ar[ur]^{\phi\otimes\Id}\ar[rr] & & \check G,
}
\end{equation}
where $\phi$ is a tempered (i.e.\ bounded on $\mathcal L_k$) Langlands parameter into $\check G_X$ and the right slanted arrow is the ``canonical''\footnote{Recall from Theorem \ref{sl2} that it is canonical up to $\check A^*$-conjugacy.} one.
\end{definition}

\begin{remark}
Notice that an ``$\XX$-distinguished Arthur parameter'' is not really an Arthur parameter into $\check G$, but rather such a parameter together with a lift of it to $\check G_X\times \SL_2$. However, given an Arthur parameter into $\check G$ we will say that it \emph{is} $\XX$-distinguished if it admits such a lift, and we will also call its Arthur packet $\XX$-distinguished.

If $X$ has roots of type $N$, we may, for the purposes of the weak conjecture, replace $\check G_X$ in the above definition by the Gaitsgory-Nadler group $\check G_{X,GN}$; however, this is not appropriate for the more refined Conjecture \ref{localconjecture-refined}.
\end{remark}

The group $\check G_X$ acts by conjugacy on the set of $\XX$-distinguished Arthur parameters, and, as happens with Langlands parameters, there is a natural class of measures on the set of $\check G_X$-conjugacy classes of $\XX$-distinguished Arthur parameters. Thus, we may for instance talk about ``$\XX$-discrete parameters'', meaning those diagrams (\ref{Xparameter}) for which they image of $\phi$ does not lie in a proper Levi subgroup of $\check G_X$.

\begin{conjecture}[Local Conjecture -- weak form] \label{localconjecture-weak} 
Let $k$ be a local field. There is a direct integral decomposition:
\begin{equation}
 L^2(\XX(k)) = \int_{[\psi]} \mathcal H_\psi \mu(\psi),
\end{equation}
where:
\begin{itemize}
 \item $[\psi]$ varies over $\check G_X$-conjugacy classes of $\XX$-distinguished Arthur parameters;
 \item $\mu$ is in the natural class of measures for $\XX$-distinguished Arthur parameters modulo conjugacy;
 \item $\mathcal H_\psi$ is isomorphic to a (possibly empty) direct sum of irreducible representations belonging to the Arthur packet associated to the image of $\psi$ in $\check G$.
\end{itemize}
\end{conjecture}

This conjecture states only \emph{necessary} conditions for a representation to belong weakly to $L^2(\X(k))$ (namely, it has to belong to the Fell closure of $\XX$-distinguished Arthur packets), and it also postulates that the $X$-discrete spectrum belongs to Arthur packets with $\XX$-discrete parameter. It may be, though, that the Hilbert space corresponding to such a parameter is zero; the refined Conjecture \ref{localconjecture-refined} will address that issue.

As far as we know, this conjecture was not anticipated elsewhere, e.g.\ in the (fairly extensive) study of the spectrum of symmetric varieties. 
The conjecture above gives (another) sense in which the Arthur packets are natural in representation theory.

\begin{remarks}

\begin{enumerate}

\item A corollary of the conjecture is this: the support of the discrete 
spectrum $L^2_{disc}(\XX(k))$ is contained in Arthur packets associated
to {\em discrete series parameters $\mathcal L_k \rightarrow \check{G}_X$.} More colloquially, if $\GG_X$ is a split $k$-group with dual group $\check G_X$, there should be a ``lifting'' from the discrete spectrum $L^2_{disc}(\XX(k))$ to the discrete series of $\GG_X(k)$.  

We have, however, no understanding of which portion of the relatively discrete spectrum
is in fact relatively supercuspidal, i.e.\ has compactly supported image in $C^\infty(X)$.

 \item  Another subtle issue is {\em which} elements of the Arthur packet for $\psi$ actually show up in $\mathcal{H}_{\psi}$ above;
 the conjecture as written gives no information on this. 
 For a start of a discussion of this point, see
 \S \ref{subsection:PIF}.

\item It is essential that the conjecture was formulated with Arthur packets, 
rather than the associated $L$-packets.  It is possible, for instance, for the Arthur type to be nontrivial
but yet $L^2(\X(k))$ contains weakly a tempered representation; an example is given by $\GG $ of type $G_2$ acting on the level set of an invariant quadratic form in its standard ($7$-dimensional) representation.

\item The conjecture predicts that $L^2(X)$ is tempered when
$\check G_X = \check G$. Although we do not have a general proof, 
Theorem \ref{ST:Plancherel} proves this in many cases. 

\item Suppose that $(G,X)$ is symmetric over $k = \RR$. Let $\mathfrak{a}$ 
be a ``maximally $\sigma$-split torus'' and let $\mathfrak{l}$
be its centralizer, the Lie algebra of the Levi subgroup associated to $X$.  It is known by the work of Flested-Jensen and Matsuki, Oshima
that any discrete series for $L^2(X)$ is the cohomological induction from $\mathfrak{l}$
of a one-dimensional representation (see \cite{VoganSymmetric} for a summary of these results.) 

The results of Adams and Johnson \cite{AdamsJohnson} likely
could be used to give evidence   that these indeed belong to Arthur packets as predicted by the conjecture, but we have not verified the details. 

\item Gan and Gomez have shown that the Howe duality correspondence  
can be used to exhibit establish the conjecture   for certain pairs $(G,X)$ (\cite{GG}).

 \item Conjecture 
\ref{localconjecture-weak} addresses the {\em unitary} spectrum; what of other $\XX$-distinguished
 representations?  In other words, what if we want to decompose $C^{\infty}(X)$ rather than $L^2(X)$?
 This is a more subtle issue; the trivial representation is always $\XX$-distinguished
 locally and globally, but does not belong to the Arthur parameter for an $\XX$-distinguished Arthur parameter (except, of course, when $\XX=$a point).  
 \end{enumerate}

\end{remarks}

\subsection{A  global to local argument} 

In this section, we admit the Arthur conjectures as formulated in \S \ref{Artur}, and establish some evidence -- 
Theorem \ref{ltg} below -- for our Conjecture \ref{localconjecture-weak} , using a local-global argument. 

We continue to denote by $K$ a global field, with ring of adeles $\adele$. We will only discuss here the case where $\check G_X$ is a subgroup of $\check G$, i.e.\ $\check G_X= \check G_{X,GN}$, so an $\XX$-distinguished Arthur parameter is indeed an Arthur parameter into $\check G$:
$$ \psi: \mathcal L_{(k\mbox{ or }K)} \times \SL_2 \to \check G$$
factoring through $\check G_X\times \SL_2$ in the specified way.

Call such a local Arthur parameter $\psi$ {\em weakly} $\XX$-distinguished if:
\begin{enumerate}
\item \label{wksl2} $\psi|\SL_2$ is conjugate to the $\SL_2$-type of $\XX$;
\item \label{wkgx} The semisimple part of every $\psi(g), \ g \in \mathcal{L}_k \times \SL_2$, 
is conjugate to an element in the image of $\check G_X \times \SL_2$. 
\end{enumerate}

In general this is strictly weaker than $\XX$-distinguished, but in some instances is equivalent to it.   
A typical example is: $\GG = \GGL_{2n}$ acting on $\XX = \{ \mbox{alternating $2$-forms.} \}$ In this case, 
$\check G_{X} = \check G_{X,GN}$ is the  centralizer of the $\SL_2$-type of $\XX$, and so 
requirement \eqref{wksl2} above already guarantees that any weakly distinguished parameter
is $\XX$-distinguished.

Let $\widehat{\G(K_v)}_{X-wkdist}$ be the set of unitary representations of $\GG(K_v)$  that belong to a  weakly $\XX$-distinguished Arthur packet. 

{
The following theorem will require some local input from \cite{SaSph}: It was proven in \cite[Theorem 9.0.1]{SaSph} that, under certain combinatorial assumptions on the spherical variety $\XX$, its unramified $L^2$-spectrum is supported on, what turns out to be, the set of $\XX$-distinguished unramified Arthur parameters. This is a result of explicit computation, and we do not know a conceptual proof or reasoning for it that does not invoke the Langlands dual. Similar results have been obtained in other (including non-split) cases by Hironaka, Offen and others, e.g.\ \cite{Hironaka, Offen}.

If we invoke Langlands duality, though, it is easy to see on the $L$-group side why this is so. The case $\check G_X\subset\check G$ that we are considering is particularly straightforward, because for the image of an $X$-distinguished Arthur parameter into $\check G$ to be unramified, the parameter itself had to be unramified in the first place. But an unramified parameter $\mathcal L_k \to \check G_X$ has image in a Cartan subgroup, hence unramified representations can only appear in the most continuous part of the spectrum.

Theorem 9.0.1 in \emph{loc.cit}.\ shows that, under assumptions on $\XX$, the Plancherel measure for $L^2(\XX(k))^{K_G}$ where $k$ is a non-archimedean place and $K_G$ is a hyperspecial maximal compact subgroup of $\GG(k)$, is supported on the set of unramified representations which are subquotients of $I_{P(X)}(\chi)$, the normalized-induced representation obtained from a (unitary) character $\chi\in \check A_X$ of $P(X)$. The conditions include that $\XX$ is homogeneous affine (that is, $\XX=\HH\backslash \GG$ with $\HH$ reductive) or Whittaker-induced from a homogeneous affine variety (in the sense of our \S \ref{Whitsubsection}). There are also certain combinatorial conditions to be checked, which are probably satisfied automatically under the previous conditions; a table of some varieties satisfying those appears at the end of \emph{loc.\ cit.}

In terms of Arthur parameters, this theorem states that there is a Plancherel decomposition the unramified spectrum of $\XX(k)$ in terms of represenations with $\XX$-distinguished unramified Arthur parameter. Notice that by property (3), \S \ref{Artur} of Arthur packets, this is the \emph{only} Arthur parameter that these representations admit.

}

\begin{theorem} \label{ltg}
Assume the Arthur conjectures (as formulated in \S \ref{Artur}). 
Suppose that:
\begin{itemize}
\item[(i)] the set of points in $\XX(K)$ with anisotropic stabilizers
is dense in $\XX(K_w)$;
\item[(ii)] $\XX_{K_v}$ satisfies the conditions of \cite[Theorem 9.0.1]{SaSph}  (see discussion before Theorem) for { almost all} places $v$. 
\end{itemize}
Then the support of $L^2(\XX(K_w))$ is contained in the closure
of $\widehat{\G(K_w)}_{X-wkdist}$ in the Fell topology.
\end{theorem} 

Note that assumption (i) is trivially
satisfied if $\GG$ is anisotropic (for this theorem we don't need to assume that $\GG$ itself is split);
as we will see in the course of proof, it 
   could also be replaced by either of the following two assumptions:
\begin{itemize}
\item[-] there exists a place $w$ for which there is a $\XX$-distinguished {\em supercuspidal} representation;
\item[-] the convolution of the invariant measure on  $[\G_{x_0}] \subset [\G]$ (where $x_0\in \XX(K)$) with a compactly supported measure on $\GG(\adele)$ is given by an $L^2$-density on $[\GG]$.
 \end{itemize}

The proof is inspired by the Burger-Sarnak principle as well as work of Clozel. It is based on a globalization result that is perhaps of independent interest (although it is very closely related to other results, 
  it has a slightly different range of applicability, since it is based on the use of the Fell topology). 

Let us denote by $\ZZ$ the center of $\GG$.
An \emph{automorphic discrete series} is an irreducible unitary representation $\pi$ of $\G(\adele)$ together with a non-zero morphism $\nu: \pi \to L^2(\G(K)\ZZ(\adele)\backslash \G(\adele), \Omega)$ (where $\Omega$ is the central character of $\pi$, an idele class character). 
Let $x_0 \in \mathbf{X}(K)$. An automorphic discrete series $(\pi, \nu)$ is \emph{$(\mathbf{X}, x_0)$-distinguished} if the functional $$f\mapsto \int_{[\mathbf{G}_{x_0}]} \nu(f)$$ is non-zero on the subspace of smooth vectors of $\pi$.\footnote{For this definition we need, of course, the integrals to converge. Ignoring analytic difficulties, we could also describe the $\XX$-distinguished spectrum as follows: There is a canonical intertwiner: $C_c^\infty(\X(\adele))\to C^\infty(\GG(K)\backslash \G(\adele))$ (summation over $\X(K)$). It seems natural to define {\em globally distinguished} as ``the image of the adjoint morphism.''   However, this leads to difficulties when not all point stabilizers on $\X$ are anisotropic.}  The restriction to smooth vectors is an important technical point: The elements of $\pi$
are $L^2$-functions on the automorphic space; but smooth vectors are given by genuine functions, and are defined on measure zero subsets such as $[\mathbf{G}_{x_0}]$.

If $S$ is any finite set of places of $K$, we write $K_S$ for $\prod_{v \in S} K_v$ and
  $K^S$ for the restricted direct product (with respect to integers) $\prod'_{v \notin S} K_v$.

\begin{theorem} \label{RefBurSar} Suppose that, for some $x_0\in \XX(K)$, $\GG_{x_0}$ is anisotropic; write  $X_{S}^{(0)}$ for $x_0\cdot \G(K_S) \subset \X(K_S)$. 
Suppose that $\sigma$ is an irreducible representation of $\G(K_S)$ which belongs to the support of Plancherel measure for
 $L^2(X_S^{(0)})$. Then there exist a sequence of $(\X, x_0)$-distinguished automorphic representations $\pi_i$ 
whose restrictions $\pi_{i,S}$ to $\G(K_S)$ converge, in the Fell topology, to $\sigma$. 
 \end{theorem}

 We now prove that Theorem \ref{RefBurSar} $\implies$ Theorem \ref{ltg}, and return
 to the proof of Theorem \ref{RefBurSar} in \S \ref{ssrbs}.

The key point of the argument is due to L. Clozel (\cite{Clozel}) and lies in the beautiful idea of considering two places simultaneously. 
Indeed, let  $S =\{v, w\}$ comprise two places, so that $v$ is a ``good place'' for $(\mathbf{G}, \mathbf{X})$, i.e.\ a place where \cite[Theorem 9.0.1]{SaSph} holds. 
That theorem computes precisely the decomposition of the {\em unramified part} 
of $L^2(\mathbf{X}(K_v))$.

Now let $\sigma_w$ be an arbitrary unitary representation that occurs weakly in $L^2(\mathbf{X}(K_w))$. 
Let $\sigma_v$ be an unramified representation occuring weakly in $L^2(\mathbf{X}(K_v))$.  Then $\sigma_S := \sigma_v \otimes \sigma_w$ occurs weakly in $L^2(\mathbf{X}(K_S))$. By Theorem \ref{RefBurSar}, 
  there exists a sequence of $\mathbf{X}$-distinguished automorphic forms 
$\pi_j$ whose local constituent at $v,w$ approach $\sigma_v, \sigma_w$ respectively, in the Fell topology.  Choose an Arthur parameter $\psi_j$ with the property that  the global packet $\mathcal{A}_{\psi_j}$
contains a representation isomorphic to $\pi_j$.

Fix now $\psi = \psi_j$. Let $Q$
be the Zariski-closure of the image of $\psi | \mathcal{L}_K$. 
 The map
 
\begin{equation} \label{mapaboveabove}  \psi: \mathcal{L}_K \rightarrow \mathcal{L}_K \times \SL_2 \longrightarrow \check{G}\end{equation} 
 visibly factors through  
\begin{equation} \label{mapabove} \left (\psi, \left( \begin{array}{cc} |\cdot|^{1/2} & 0 \\ 0 & |\cdot|^{-1/2} \end{array}\right) \right) : \mathcal{L}_K \longrightarrow 
\check G \end{equation}
where we regard the target $\Gm$ as the torus in $\SL_2$.   

By property (3), \S \ref{Artur} of Arthur packets,  together with the theorem  of \cite{SaSph}
that computes the unramified Plancherel measure, the $\SL_2$-type of $\psi_j$ (for large enough $j$) is the same as the $\SL_2$-type $\iota: \SL_2 \rightarrow \check G$ associated to the spherical variety (viz. the principal $\SL_2$ for $\check \LL_{\XX}$).  
 
 { Here, we have used a strengthening of (3) on page \pageref{pagerefitem3}, given a sequence of unramified representations  $\sigma_j$ of 
 $G$ belonging to Arthur packets with a given $\SL_2$-type $\iota_1$,
 and another unramified representation $\sigma_{\infty}$ which is a limit point for the Fell topology
 of the $\sigma_j$, then any Arthur packet containing $\sigma_{\infty}$ also has $\SL_2$-type $\iota$. 
  Indeed, it follows  (see e.g.  \cite[Theorem 2.2]{tadic})   the Satake parameters of $\sigma_j$
  approach that of $\sigma_{\infty}$, and then we argue as on page \pageref{pagerefitem3},
  where now equality is replaced by ``almost equality,'' which is enough for our purposes
  because the set of possibilities for ``$\eig (\alpha_i(q^{1/2}))$''  (as in the last line of the argument on page \pageref{pagerefitem3} is discrete. 
  }

Let $\mathrm{Fr}_s$ denote the image, under \eqref{mapabove},  of a Frobenius element at the unramified place $s$. 
The $\mathrm{Fr}_s$, when we vary the place $s$, are Zariski-dense in $Q \times\Gm$. Indeed, by assumption \eqref{fzd}, 
their Zariski closure gives a subgroup $Q' \subset Q \times \Gm$ projecting onto $Q$;
if $Q' \neq Q \times \Gm$, then there is an integer $m \geq 1$ and a character
$\chi: Q \rightarrow \Gm$ so that $Q' = \{(q, x): x^m = \chi(q)\}$.  But this cannot be:
if we write $(q_s, x_s)$ for the image of $\mathrm{Fr}_s$ in $Q(\C) \times \Gm(\C)$, 
then $q_s^{\mathbb{Z}}$ is precompact whereas $x_s^{\mathbb{Z}}$ is not.

On the other hand, the aforementioned theorem of \cite{SaSph} shows more precisely that the image of every $\mathrm{Fr}_s$
inside $\check G$, under \eqref{mapaboveabove}, is conjugate to an element of $f(\check G_X \times\SL_2)$. (Here, $f: \check{G}_X \times \SL_2 \rightarrow \check{G}$
is the distinguished morphism associated to $X$.)

Therefore, a Zariski-dense set of elements in the image $Q^*$ of $Q \times \Gm$ inside $\check {G}$ 
are conjugate to elements of $f(\check G_X \times \SL_2)$. 

 We have just seen that any conjugacy-invariant function on $\check G$,
zero on $f(\check G_X \times \SL_2)$, must be identically zero on $Q^*$. 
Since conjugacy-invariant algebraic functions separate semisimple conjugacy classes,
we conclude that every semisimple element in $Q^*$ is conjugate to $f(\check G_X \times \SL_2)$. 
The claimed assertion follows.  \qed

\subsection{Proof of Theorem \ref{RefBurSar}} \label{ssrbs}

The theorem follows immediately from:
\begin{quote} 
Let $\Omega \subset 
  G_S = \G(K_S)$ be compact and $f \in C_c(X_S^{(0)})$. 
 Then $y \mapsto \langle y f, f \rangle $, for $y \in \Omega$,  is a convex combination of diagonal matrix coefficients 
 of $(\X, x_0)$-distinguished automorphic representations. 
\end{quote}

Indeed, the quoted 
statement means that $L^2(X_S^{(0)})$, considered as a $G_S$-representation, is ``weakly contained''  in 
a direct sum $\bigoplus \pi_i$ of the various $(\X, x_0)$-distinguished automorphic representations.
According to \cite[Proposition 8.6.8]{Dixmier}  -- see also Theorem 3.4.4, {\em loc. cit.} for the definition of weak containment -- this means that the support of Plancherel measure
for $L^2(X_S^{(0)})$ is contained in the support of Plancherel measure for $\bigoplus \pi_i$, in particular, in the closure of the set of restrictions $\pi_i|G_S$.

  \proof
 We write the proof in the $S$-arithmetic language, rather than adelically. Fix a function $f \in C_c(X_S^{(0)})$.

 Let $H_S$ be the point stabilizer of $x_0 \in X_S^{(0)}$ in $G_S$, so that $H_S = \H(K_S)$. 
  Choose a congruence subgroup $\Gamma$  of $\G(K_S) = G_S$
and  $U$ a compact subset of $G_S$ so that,
with  $\Gamma_H :=
\H(K) \cap \Gamma$, we have: 
\begin{itemize}
\item[(i)] $U=U^{-1}$ and $\Omega \subset U$; 
\item[(ii)] $x_0 U$ contains the support of $gf$ for any $g \in \Omega$;
\item[(iii)]  $\Gamma_H \cdot U \supset H_S$
\item[(iv)] $U_4\cdot H_S \cap \Gamma =\Gamma_H$, where we write $U_2 = U \cdot U, U_3= U \cdot U \cdot U$ and so on. 
\end{itemize}

This can be done: First of all choose $\Gamma$ and choose $U$ satisfying (i) -- (iii). Now
we may shrink $\Gamma$, leaving $\Gamma_H$ unchanged, so that (iv) is satisfied,
by passing from $\Gamma$ to a subgroup of the form $\pi_N^{-1} \pi_N \Gamma_H$, 
where $\pi_N$ is the ``reduction modulo $N$'' map, for a suitable large ideal $N$. 

Now fix $g \in \Omega$ and
set $$ F (g):=   \sum_{\gamma \in \Gamma_H \backslash \Gamma} f(x_0 \gamma g),$$
a compactly supported function on $\Gamma \backslash G_S$. 
(Indeed, its support is contained in $\Gamma \backslash \Gamma U_2$). 

We are going to show that
\begin{equation} \label{ipequality}  \langle  yf, f \rangle_{X_0^{(S)}} =c \langle  yF, F \rangle_{\Gamma \backslash G_S}, \ \ (y \in \Omega). 
\end{equation}
where the positive constant $c$ depends only on normalization of measure. 
This will conclude the proof: 

If $\pi $ is an irreducible $G_S$-subrepresentation
of functions on $\Gamma \backslash G_S$ (we do not require square integrability!) 
and there exists $v \in \pi$ such that
$$ \langle F , v \rangle_{\Gamma \backslash G_S} \neq 0,$$
then $\pi$ is distinguished, because
$$ \int_{\Gamma \backslash G_S}  F (g) \overline{v(g)} dg 
= \int_{\Gamma_H \backslash G_S} f (x_0 g) \overline{v(g)} dg
\\ = \int_{h \in \Gamma_H \backslash H_S} \int_{g \in H_S \backslash G_S} f (x_0   g) \overline{v(hg)} dg,$$
so some translate $v^g$ has nonvanishing period over $\int_{\Gamma_H \backslash H_S}$. 
In particular, the expression $ \langle g F , F  \rangle$
is in fact  -- after spectrally expanding $F$ -- a convex combination of diagonal matrix coefficients
of $(\XX, x_0)$-distinguished representations, as required.    

Thus indeed \eqref{ipequality} will conclude the proof.

We compute, by unfolding, that
\begin{multline} \langle  yF, F \rangle_{L^2(\Gamma \backslash G_S)} =
\int_{g \in \Gamma \backslash G_S}yF(g ) \overline{ \sum_{\gamma \in \Gamma_H \backslash \Gamma}
f (x_0 \gamma g)}   dg \  \\ =  
\int_{g \in \Gamma_H \backslash G_S}    yF(g )  
\overline{f(x_0  g)}dg
\\ =  
\int_{\Gamma_H \backslash G_S}  dg \ f(x_0 gy) \overline{ f (x_0 g) } + \int_{\Gamma_H \backslash G_S} dg \
\sum_{\gamma \in \Gamma_H \backslash \Gamma - \{1\}}  f(x_0  \gamma gy) \overline{ f (x_0  g) }
 \end{multline}

 We claim that the final term is zero:   If not, there exists $ \gamma \in \Gamma - \Gamma_H$ 
 and $g \in G_S$ 
 such that $f( x_0 \gamma g y)$ and $f_2 (x_0   g)$ are both nonzero. 
 In particular, by (ii)
 $$ \gamma gy \in H_S U \mbox { and }   g \in H_S U. $$
Adjusting $\gamma$ on the left by $\Gamma_H$ and using (iii), 
we may suppose that $\gamma g y \in U_2$.  
 Therefore, 
 $$ \gamma = (\gamma g y) \cdot y^{-1} \cdot  g^{-1} \in U_4 H_S, $$  
a contradiction to (iv).

\qed
\begin{remark} 
This has the following corollary: 

\begin{quote} If $\sigma$ is {\em automorphically isolated} as well
as in the support of $L^2(X_S^{(0)})$,  it is the local constituent of an $(\X, x_0)$-distinguished
global representation. 
\end{quote}

Here, we say that 
a unitary irreducible representation $\sigma$ of $\G(K_S)$ is {\em automorphically isolated}
if there do not exist a sequence $\sigma_i$ of unitary $\G(K_S)$-representations, 
each of which occurs as the local constitutent of an automorphic representation,
which converge to $\sigma$ in the Fell topology.

  Results of a similar nature are well-known;
see e.g. \cite{SPP} when the $\sigma$ are supercuspidal.  The condition noted above (automorphically isolated) is 
very slightly  weaker.  For instance, every discrete series for $\mathbf{GL}_n$ is automorphically isolated, and it seems likely that discrete series representations are {\em always} automorphically isolated (although we do not know how to prove it -- it would be interesting to verify that this is a consequence of Arthur's conjectures,  or to verify it for the other classical groups). 
\end{remark}

\subsection{Pure inner forms} \label{subsection:PIF} 

We now formulate a refined version of the prior Local Conjecture
\ref{localconjecture-weak}.

For the purposes of the present subsection we assume:
\begin{enumerate}
 \item The spherical variety $\XX$ has no roots of type $N$; in particular, its dual group $\check G_X$ is defined.
 \item The center of $\GG$ acts faithfully on $\XX$. (If this is not the case, one should replace $\GG$ by its quotient by the kernel of the action of $\mathcal Z(\GG)$.)
\end{enumerate}

By a ``pure inner form'' of $\G$ we understand
  an isomorphism class $\alpha$ of left $\GG$-torsors. 
We call an isomorphism class $\alpha$ of left $\GG$-torsors, with the property that $\XX \times^{\G} \TT(k)\ne \emptyset$
for $\TT$ in the class $\alpha$, a ``pure inner form'' of $\XX$.  We denote by $\GG^{\alpha}$
the automorphism group of a torsor in the class $\alpha$, and we denote by $G^{\alpha}$ its $k$ points. 
We give further discussion and examples (see Examples \ref{pifex1}--\ref{pifex3}) after we formulate the conjecture.

Let us recall that Vogan \cite{Vogan} has proposed a version of Arthur's conjectures 
whereby 
  the Arthur packet should be considerd, in fact, as a collection of representations
of varying pure inner forms of $\GG$.  In particular each element
of the Arthur packet defines a representation of the group $\prod_{\beta} G^{\beta}$,
where $\beta$ ranges over pure inner forms of $G$, and the representation
is understood to be nontrivial on only one direct factor.

\begin{conjecture}[Local Conjecture -- strong form] \label{localconjecture-refined}  
There is a direct integral decomposition:
\begin{equation}\label{localstrongeq}
 \oplus_\alpha L^2(X^\alpha) = \int_{[\psi]} \mathcal H_\psi \mu(\psi),
\end{equation}
where
\begin{itemize}
 \item $\alpha$ parametrizes pure inner forms of $X$; 
 \item  we regard both sides as representations of the product: 
 \begin{equation}\label{pifproduct}\prod_{\beta}  G^{\beta}\end{equation} of all pure inner forms of $G$ -- the right-hand side  as discussed above, and the left-hand side
 by means of the evident map from pure inner forms of $X$  to pure inner forms of $G$; 
 \item $[\psi]$ varies over $\check G_X$-conjugacy classes of $\XX$-distinguished Arthur parameters;
 \item $\mu$ is in the natural class of measures for $\XX$-distinguished Arthur parameters modulo conjugacy;
 \item $\mathcal H_\psi$ is isomorphic to a \emph{multiplicity-free} direct sum of irreducible representations belonging to the (Vogan) Arthur packet associated to the class $[\psi]$;
 \item for $\mu$-almost all $\psi$, the spaces $\mathcal H_\psi$ are \emph{non-zero}.
\end{itemize}
\end{conjecture}

{ Note that, when multiple inner forms of $X$ correspond to the same inner form of $G$, the \emph{same} irreducible representations of this inner form may appear multiple times on both sides of \eqref{localstrongeq}; this doesn't contradict the multiplicity-freeness requirement, since we consider elements of the Vogan-Arthur packet as representations of \eqref{pifproduct} (where the isomorphic inner forms can appear as distinct factors).}

 In comparison to the weak version \ref{localconjecture-weak}, the present form states that the condition on Arthur parameters to be $\XX$-distinguished is also \emph{sufficient} for the A-packet to be distinguished, as long as we take pure inner forms into consideration. 

We have also postulated that the spaces $\mathcal H_\psi$ should be multiplicity-free. In the case of $\GL_n$, where A-packets are singletons, this means that the multiplicity of a given representation on $\oplus_\alpha L^2(X^\alpha)$, at least generically in the sense of Plancherel measure, is given by the number of lifts of its Arthur parameter to an $\XX$-distinguished Arthur parameter. For example, the multiplicity statement is true for the most continuous spectrum of $X$ under the assumptions of the Scattering Theorem \ref{advancedscattering}: indeed, the spectrum of the most degenerate boundary degeneration $X_\emptyset$ is a multiplicity-free direct integral over Arthur parameters with ``Langlands part'' into the maximal torus $A_X^*$ of $\check G_X$, and the corresponding ``most continuous spectrum'' $L^2(X)_\emptyset$ is a multiplicity-free direct integral over $W_X$-conjugacy classes of such parameters. 
However, we should point out that there is not enough evidence about whether the multiplicity statement is correctly formulated for ramified representations in the case of nontrivial Arthur-$\SL_2$.

In the remainder of this section, we examine more carefully the notion of pure inner form of $\XX$:

Consider the quotient stack: $[\XX\times \XX /\GG]$ (we understand the diagonal action of $\GG$ without putting brackets). We denote by $[\XX\times \XX /\GG](k)$ the \emph{set of isomorphism classes of $k$-objects} of the stack. By abuse of language, we will be calling them ``$k$-points''. They consist of isomorphism classes of diagrams: 
$$ \TT \to \XX\times \XX,$$
where $\TT$ is a (right) $\GG$-torsor and the map is $\GG$-equivariant. Two such diagrams are isomorphic if there is an isomorphism of torsors which makes the composite commute.

In what follows we will denote isomorphism classes of $\GG$-torsors by small Greek exponents (they correspond bijectively to elements of $H^1(k,\GG)$), and the exponent will appear on the opposite side from which $\GG$ acts. For instance, ${^\alpha \TT}$ denotes a right $\GG$-torsor ``in the class $\alpha$''; by composing with the inverse map we get a left $\GG$-torsor which will is denoted $\TT^\alpha$. The $\GG$-automorphism group of a torsor is an inner form $\GG$; for the torsors ${^\alpha\TT}$ and $\TT^\alpha$, this form will be denoted by $\GG^\alpha$ and will act on the left, resp.\ on the right. Notice the canonical $\GG^\alpha\times\GG^\alpha$-equivariant isomorphism: ${^\alpha\TT}\times^\GG \TT^\alpha \simeq \GG^\alpha$ (here, unlike the rest of the paper, the left multiplication of $\GG^\alpha$ on itself is defined as a left action).

Given a $\GG$-variety $\VV$ and a left $\GG$-torsor $\TT^\alpha$, we denote by $\VV^\alpha$ the $\GG^\alpha$-variety: $\VV\times^\GG \TT^\alpha$. The following is easy to see by applying the $\times^\GG \TT^\alpha$ operation:
\begin{lemma}\label{pointsonstacks}
 The set of isomorphism classes of $\GG$-morphisms: ${^\alpha\TT}\to \VV$ is in natural bijection with the set of $\GG^\alpha(k)$-orbits on $\VV^\alpha(k)$. In particular, the existence of such a morphism is equivalent to the statement: $\VV^\alpha(k)\ne \emptyset$.
\end{lemma}

To apply this to $\VV = \XX\times\XX$, where $\XX$ is our spherical $\GG$-variety and $\GG$ acts diagonally on $\VV$, we notice that $\VV^\alpha = \XX^\alpha\times\XX^\alpha$. Indeed, if $\VV$ carries an action of a larger group $\tilde\GG\supset \GG$ and $\F^\alpha:= \tilde\GG \times^\GG \TT^\alpha$, a $\tilde\GG$-torsor, then obviously $\VV\times^\GG \TT^\alpha = \VV\times^{\tilde\GG} \F^\alpha$; in our case, $\tilde\GG:=\GG\times\GG$ and $\F^\alpha = \TT^\alpha\times \TT^\alpha$, so we get:

\begin{lemma}
 The set of isomorphism classes of $\GG$-morphisms: ${^\alpha\TT}\to \XX\times\XX$ is in natural bijection with the set of $\GG^\alpha(k)$-orbits on $\XX^\alpha(k)\times\XX^\alpha(k)$. In particular, the existence of such a map is equivalent to the statement: $\XX^\alpha(k)\ne \emptyset$.
\end{lemma}

\begin{corollary}
 We have:
\begin{equation}
 [\XX\times \XX /\GG](k) = \sqcup_\alpha \XX^\alpha(k)\times \XX^\alpha(k) /\GG^\alpha(k),
\end{equation}
where $\alpha$ runs over all pure inner forms of $\XX$.
\end{corollary}

\begin{example} \label{pifex1} 

If $\XX = \HH\backslash \GG$, the pure inner forms of $\XX$ correspond to $\GG$-torsors obtained by reduction of $\HH$-torsors, i.e.\ to the image of the map: $H^1(k,\HH)\to H^1(k,\GG)$.

In fact, we have an isomorphism of stacks: $\XX\times\XX/\GG\simeq \HH\backslash\GG/\HH$ and $\HH$ has a fixed point on $\HH\backslash \GG$, one gets from Lemma \ref{pointsonstacks}:
 
\begin{equation} 
 [\XX\times \XX /\GG](k) = \sqcup_{\beta\in H^1(k,\HH)} (\HH^\beta\backslash \GG^\beta)(k)/\HH^\beta(k).
\end{equation}
 Here, similarly, we denote by $\HH^\beta$ the automorphism group of an $\HH$-torsor ${^\beta\SS}$ in the class of $\beta$ and by $\GG^\beta$ the isomorphism class of its reduction to a $\GG$-torsor ${^\beta\TT}$ (i.e.\ ${^\beta\TT}= {^\beta\SS}\times^\HH \GG$); the action of $\HH^\beta$ on ${^\beta\SS}$ gives rise to a natural injection: $\HH^\beta\hookrightarrow\GG^\beta$.  
\end{example} 

\begin{example} \label{pifex2} 
 The pure inner forms of $\XX=$ a point coincide with the isomorphism classes of $\GG$-torsors, despite the fact that all varieties $\XX^\alpha$ are isomorphic.
\end{example}
 \begin{example} \label{pifex3} 
 Let $V\subset W$ be two non-degenerate quadratic spaces of codimension one in each other and let $\XX=\HH\backslash \GG = \SSO(V)\backslash (\SSO(V)\times\SSO(W))$. 
 
 Isomorphism classes of $\SSO(V)$-torsors correspond canonically to isomorphism classes of quadratic spaces of the same dimension and discriminant as $V$, and similarly for $\SSO(W)$-torsors. The reduction of a $\SSO(V)$-torsor to an $\SSO(W)$-torsor corresponds to the operation $V^\alpha\mapsto V^\alpha\oplus k$ (orthogonal direct sum), where $V^\alpha$ is a quadratic space corresponding to the given torsor. Pure inner forms of $\XX$ correspond to isomorphism classes of quadratic spaces $V^\alpha\subset W^\alpha$, where $V^\alpha$ has the same dimension and discriminant as $V$ and $W^\alpha\simeq V^\alpha\oplus k$. This is the setting of the Gross--Prasad conjectures \cite{GP};
 these conjectures   have been now largely established  through the work of Waldspurger and (for the unitary analogue) Beuzart--Plessis (\cite{Waldspurger-GP, Waldspurger-GP2, Waldspurger:eps, B-P:SMF, B-P:Arch}).  
\end{example}

\begin{remark}
Suppose that $\XX = \HH\backslash \GG$, and that $\HH= \MM\ltimes \UU$ is the Levi decomposition of $\HH$. Let $\Lambda:\UU\to \GGa$ be a homomorphism which is fixed by $\MM$. 
Extend it to a homomorphism $\Lambda: \HH \to \GGa$ by making it trivial on $\MM$. 
We are interested in defining ``pure inner forms'' for the space $L^2(X,\mathcal L_\psi)$, where $\psi: k\to \CC^\times$ is a character and $\mathcal L_\psi$ is the complex line bundle defined by $\psi\circ\Lambda$. 

Since unipotent groups in characteristic zero have trivial Galois cohomology, all pure inner forms of $\XX$ correspond to $\MM$-torsors. Let $c:\Gal(\bar k/k)\to \MM(\bar k)$ be a cocycle; by inner automorphisms, it defines an inner twist $\HH^\alpha$ of $\HH$. Since the action of $\MM$ preserves the morphism $\Lambda:\UU\to \GGa$ (where $\MM$ acts trivially on $\GGa$), the ``inner twist'' of $\Lambda$ by $c$ is defined over $k$, that is, we have a morphism $\Lambda^\alpha:\HH^\alpha\to\GGa$. 

With this convention, then, the Conjecture also applies to ``Whittaker-type'' induction.
\end{remark}

 \section{Speculation on a global period formula} 
 \label{sec:globalconj}
 
 Throughout this section we adopt the following notation: For $K$ a global field -- which we understand as fixed -- and $\GG$ any algebraic group over $K$, we denote by
$[\GG]$ the adelic quotient $\GG(K) \backslash \GG(\adele)$. 

In this section we wish to discuss a potential generalization of the work of Ichino--Ikeda \cite{II} to all spherical varieties.   Our central conjecture, speaking somewhat imprecisely, gives a link between the {\em local Plancherel formula} and {\em global periods}. 

We assume that $\XX$ is \emph{homogeneous affine}, i.e.\ the stabilizers of points on $\XX$ are reductive; the discussion is also valid for ``parabolically-induced'' or ``Whittaker-induced'' varieties of this form. For a discussion of the general non-affine case, cf.\ \cite{SaRS}; our conjectures naturally extend to this more general case, but because of the speculative nature of this case we will ignore it in our formulations and only appeal to a special case in Theorem \ref{holdsbyunfolding}.

We will consider automorphic representations that embed weakly in $L^2([\GG])$, hence: abstract irreducible \emph{unitary} representations $\pi\simeq \otimes'_v \pi_v$ of $\GG(\adele)$ which admit a \emph{tempered} embedding: $$\mbox{ smooth subspace of $\pi$} \hookrightarrow C^\infty([\GG]),$$ i.e.\ the image is in $L^{2+\varepsilon}$ for every $\varepsilon>0$. The embedding will not be part of the data, but it will be assumed to be unitary whenever the image is discrete modulo center. The normalization of embeddings corresponding to the continuous spectrum will be discussed in \S \ref{ssexplanation}.   

For the automorphic representations $\pi$ that we will encounter, \textbf{we make the following multiplicity-one assumption}:
\begin{center} For all places $v$ of $K$, we have $\dim\Hom_{\GG(K_v)} (\pi_v,C^\infty(\XX(K_v))) \le 1.$
\end{center}
Strictly we should write $\pi_v^{\infty}$ above, instead of $\pi_v$ -- which is by definition unitary. We will allow ourselves
this imprecision at some points below.

As the work of Jacquet \cite{JEuler} shows (see also \cite{FLO}), the multiplicity-one assumption is too restrictive -- one could have Euler products even without it; however, we will contend ourselves to provide some conjectures in this more restrictive setting.

\subsection{Tamagawa measure} \label{Tamagawameasure}

We use throughout Tamagawa measures for $[\GG]$, $\GG(\adele)$ and, more generally, the adelic and local points of smooth homogeneous varieties. 

To define Tamagawa measure we proceed as follows:  
Let $\mu$ be the measure on $\adele$ that assigns
mass $1$ to the quotient $\adele/K$, and fix a factorization
$\mu=\mu_v$, where $\mu_v$ is a measure on $K_v$. 
Fix a $K$-rational top differential form $\omega$.  Then 
(using the choice of $\mu_v$) 
we obtain  a volume form
$|\omega|_v$ on $\GG(K_v)$. 
For $v$ a finite place, set  $c_v $ to be the $|\omega_v|$-mass of $\GG(\mathfrak o_v)$ and take $c_v=1$ otherwise. 
For all but finitely many places, $c_v$ is the value of a certain local $L$-factor,
and we can interpret the (non-convergent, in general) Euler product $C = \prod_v c_v$ accordingly. 
If $\GG$ has trivial $k$-character group, then $C \neq 0$ and we define the Tamagawa measure:
$$  C^{-1} \cdot \left( \prod_v c_v |\omega_v| \right). $$

 (If $\GG$ has a nontrivial character group, one usually regularizes the situation  
 by multiplying or dividing by the appropriate power of the (correspondingly partial) Dedekind zeta function. 
We will discuss in \S \ref{ssexplanation} how our conjecture is independent of such a choice.)

 We fix factorizations of the Tamagawa 
measures, e.g.\ if $dx$ denotes the invariant Tamagawa measure on $\XX(\adele)$, we fix an Euler product: $dx = \prod_v dx_v$, where $dx_v$ is a measure on $\XX(K_v)$ and $dx_v(\XX(\mathfrak o_v))=1$ for almost all $v$.

\subsection{Factorization and the Ichino--Ikeda conjecture}  \label{iisubsec}
Pick $x_0 \in \mathbf{X}(K)$; let $\mathbf{H}$ be the stabilizer of $x_0$. We assume throughout that the connected component of the center of $\GG$ acts faithfully on $\XX$.

Let $\nu:\pi = \otimes \pi_v \hookrightarrow C^\infty([\GG]) $ be an automorphic representation, together with an embedding. By multiplicity one, 
the global ``period'' Hermitian form 
\begin{equation}\label{autperiod} \P^\Aut :  \varphi \longrightarrow \left|  \int_{[\HH]}  \nu(\varphi ) \right|^2 \end{equation}
factorizes as a product of local $\HH(K_v)$-biinvariant Hermitian forms $\P_v$ on the representations $\pi_v$. Of course, this period integral is not always convergent, and has to be suitably regularized. \textbf{In the discussion that follows we assume such a regularization.}

\medskip

\underline{ Basic question.} {\em Is it possible to give a purely local expression for $\P_v$? } 
\medskip

Let us make this more precise by formulating the answer given by Ichino--Ikeda \cite{II}, based on the results of Waldspurger \cite{Waldspurger-torus} and others. First, assume that $\XX$ is ``strongly tempered'', and let
$\P^\Planch_v$ be the ``canonical hermitian form'' discussed in \S \ref{canhermform}:
\begin{equation} \label{Chf} \P^\Planch_v(u_1,u_2) = \int_{\mathbf{H}(K_v)} \langle \pi_v(h) u_1, u_2 \rangle dh .\end{equation} 

The case under consideration in \cite{II} is $\GG = \SSO_V \times \SSO_{V\oplus \Ga}$, where $V$ is a nondegenerate quadratic space, and $\HH=$the diagonal copy of $\SSO_V$ in $\GG$. If $\pi$ is a tempered and cuspidal automorphic representation of $\GG$ (in this case there is a unique, up to scalar, embedding of $\pi$ in $C^\infty([\GG])$), and we fix a \emph{unitary} such embedding, then by the Arthur conjectures it should be attached to a ``global Arthur parameter'' $\phi$ (cf.\ \S \ref{Artur}) whose centralizer in the connected dual group $\check G$ is a finite $2$-group $S_\phi$, and which is trivial on the ``Arthur $\SL_2$'' -- i.e., it is a global Langlands parameter. In that case, Ichino and Ikeda conjecture:
\begin{equation} \label{IIconj}  \P^\Aut = \frac{1}{|S_{\phi}|} \prod_{v}' \P^\Planch_v, \end{equation}

This Euler product is not absolutely convergent, and the symbol $\prod'$ denotes that it should be understood ``in the sense of $L$-functions'': more precisely, it is computed in \cite{II} that for local unramified data the local factors are equal to a certain quotient of special values of $L$-functions, and the meaning of $\prod'$ is that one should replace almost all Euler factors by the corresponding quotient of special values of the (analytically continued) partial $L$-functions.

The conjecture has been verified for $n \leq 3$, and special cases in higher rank. There is also evidence that {\em exactly the same conjecture} applies to the Whittaker period, with the regularized Plancherel hermitian forms that we defined in \S \ref{ssWhittaker}; in the case of $G=\GL_n$ the validity of the conjecture for the Whittaker case is known to experts, and we will recall the argument in Section \ref{sec:examples}; Lapid and Mao have recently proven it for automorphic representations of the double metaplectic cover of $\Sp_{2n}$ \cite{Lapid-Mao-Whittaker}. As we will see, by our interpretation of ``unfolding'' (\S \ref{ssunfolding}) the conjecture also holds whenever we can ``unfold'' the period integral to a known case, such as in the case of the Rankin-Selberg integral on $\GL_n \times \GL_{n+1}$ (i.e.\ the space $\GL_n\backslash\GL_n\times\GL_{n+1}$).

In the general case, the form  \eqref{Chf} does not converge, nor can it be regularized by the methods of \S \ref{ssWhittaker}.  A typical example where it diverges is the 
$\mathrm{Sp}_{2n}$-period inside $\GL_{2n}$, which has been studied by Jacquet--Rallis \cite{Jacquet-Rallis} and Offen
\cite{Offen-symplectic}. 
However, by Proposition \ref{ST:Plancherel}, the Hermitian form $\P^\Planch_v$ is intrinsically characterized -- at least off a set of Plancherel measure zero -- by its role in a Plancherel formula for $L^2(\mathbf{X}(K_v))$.  This suggests a reformulation of the definition of $\P^\Planch_v$
which has the possibility of working even when \eqref{Chf} is divergent.

\subsection{Local prerequisites for the conjecture}\label{sslocalpre}  We keep assuming multiplicity one at all places, a corollary of which is that $\check G_X\subset \check G$. We feel free to assume the validity of all conjectures in this paper for $X$, in particular Conjecture \ref{localconjecture-weak} on the local $L^2$-spectrum of $X$.  Hence, 
\begin{equation}\label{L2withApackets}
 L^2(\XX(K_v)) = \int_{\{\phi\}/\sim} \mathcal H_\phi \mu(\phi)
\end{equation}
is a decomposition of $L^2(\XX(K_v))$ in terms of $\check G_X$-conjugacy classes of $X$-distinguished Arthur parameters, where the measure $\mu$ belongs to the natural class of measures on Arthur parameters. The spaces $\mathcal H_\phi$ here may be zero.

We would like to fix a Plancherel measure in this class, in order to fix the associated norms on $\mathcal H_\phi$. The idea is to ``fix the natural Plancherel measure for $\GG_X(K_v)$'', where $\GG_X$ is the split group with dual group $\check G_X$. This is slightly problematic, in the sense that the Plancherel measure for $\GG_X(K_v)$ is, according to the conjectures of Hiraga--Ichino--Ikeda \cite{HII}, not quite a measure on the set of (tempered) Arthur parameters into $\check G_X$, but also depends on the representation in the corresponding packet. However, for the purposes of the present discussion, where we formulate our conjecture up to $\mathbb Q^\times$, this will not matter.

More precisely, it is expected that there is a measure $\mu_v$ on the set of local, tempered Arthur parameters (i.e.\ bounded Langlands parameters) into $\check G_X$ modulo conjugacy, such that $L^2(\GG_X(K_v))$ admits a direct integral decomposition analogous to (\ref{L2withApackets}), where for every unitary, tempered representation $\tau$ the Plancherel norm on the space spanned by its matrix coefficient is a multiple of the canonical (Hilbert-Schmidt) norm by an integer (which can be bounded independently of the representation). In fact, there is a \emph{minimal} such choice in the $\mathbb Q^\times$ class of $\mu_v$, in the sense that with that choice for some representations in the packet (those, conjecturally, corresponding to characters of the component group of the centralizer of the parameter) we will not need to multiply by an integer. Of course, this measure depends on the choice of a measure on $\GG_X(K_v)$, but 
again we may choose Tamagawa measures globally to eliminate the dependence on local choices in the conjecture that follows. A good choice of local Tamagawa measures is described in \cite{Gross-motive}, and for this choice there is a very precise conjecture on Plancherel measures in \cite{HII}. For discrete series:
\begin{equation}\label{formaldegrees}
 \mu^\Planch_v(\tau)=\frac{\left<1,\pi\right>}{|S_{\phi}^\natural|} \cdot |\gamma(0,\tau,\Ad,\psi)| 
\end{equation}
hence in that case we would take the measure for the corresponding Langlands parameter $\phi$ to be: $$\mu_v(\phi)= \frac{1}{|S_{\phi}^\natural|} \cdot |\gamma(0,\tau,\Ad,\psi)| .$$ We refer the reader to \cite[p.\ 287]{HII} for the notation and normalization.

By the obvious bijection between (local) $X$-distinguished Arthur parameters into $\check G$ and tempered Langlands parameters into $\check G_X$ we can consider this measure as a measure on the set of $\check G_X$-conjugacy classes of $X$-distinguished Arthur parameters. To this choice of measure corresponds, for almost every $\tau \in L^2(\XX(K_v))$, a generalized character:
\begin{equation}
 \theta_v^\Planch : C_c^\infty(\XX(K_v)\times \XX(K_v))\to \tau\otimes\bar\tau \to \CC
\end{equation}
or dually (and composing with evaluation at the chosen point $x_0$) a hermitian form:
\begin{equation}
 \mathcal P_v^\Planch: \tau\otimes\bar\tau \to C^\infty(\XX(K_v)\times\XX(K_v))\xrightarrow{\ev_{(x_0,x_0)}} \CC.
\end{equation}
Although this is, a priori, not well defined at any specific $\tau$ which is not in the discrete spectrum of $\XX(K_v)$, it is expected to be rational in $\tau$; this follows from Theorem \ref{explicitBernstein} under the assumptions of that theorem. Therefore, $\mathcal P_v^\Planch$ is uniquely defined wherever it is regular.

\underline{Question:} Is it true that $\mathcal P_v^\Planch$ is \emph{regular} on the set of $X$-tempered representations?

From now on we will assume this to be so, or we will assume the local components of our global representations to be on the regular set. The last local piece of input that we need to discuss is the value of $\mathcal P_v^\Planch$ on normalized unramified data.

Let $v$ be a non-archimedean place of $K$, unramified over $\mathbb Q$. Assume that $\GG$ and $\XX$ carry an integral model at $v$ with $x_0\in\XX(\mathfrak o_v)$ and such that $\XX(K_v)$ satisfies the ``generalized Cartan decomposition'' \cite[Axiom 2.4.1]{SaSph} with respect to the hyperspecial maximal compact subgroup $\GG(\mathfrak o_v)$ -- this is the case at almost every place under the multiplicity-one assumption. Let $\pi$ be an irreducible, unramified (with respect to $\GG(\mathfrak o_v)$) representation in $L^2(\XX(K_v))$, which is isomorphic to the unramified subquotient of the representation $I_{P(X)}^G(\chi)$. Then, up to a combinatorial condition which is easy to check and which is expected to hold for all affine homogeneous spherical varieties ({s.\ the statement of \cite[Theorem 7.2.1]{SaSph} --} from now on \textbf{we assume this condition to hold for our spherical variety}), the value of $\mathcal P_v^\Planch(u\otimes\bar u)$ where $u\in \pi^{\GG(\mathfrak o_v)}$ with $\|u\|=1$ 
follows from the Plancherel formula of \cite[Theorem 9.0.1]{SaSph}.
 More 
precisely, it is the quotient $L_X$ 
of $L$-values attached to the spherical variety $\XX$ in \emph{loc.cit.}, divided by the ``Plancherel measure for $G_X$''. Hence,\footnote{See the table at the end of \cite{SaSph} for some examples.} in the notation of \cite[Definition 7.2.3]{SaSph}  (in particular: using exponential notation $e^{\check\gamma}$ instead of $\check\gamma$ for characters of tori), but adding the index $v$:
\begin{equation}\label{defLX}
 L^\sharp_{X,v}(\pi) := \mathcal P_v^\Planch(u\otimes\bar u) = 
\end{equation}
\begin{equation*}
 =\frac{c_v^2}{Q^{P(X)}_v} \cdot \frac{\prod_{\check\gamma\in\check\Phi_X}(1-q_v^{-1}e^{\check\gamma})}{\prod_{\check\theta\in\Theta} (1-\sigma_{\check\theta}q_v^{-r_{\check\theta}}e^{\check\theta})}(\chi). \end{equation*}

We recall that $Q^{P(X)}_v=[\GG(\mathfrak o_v): \PP(\XX)^-(\mathfrak o_v) \PP(\XX)(\mathfrak o_v) ]$, where $\PP(\XX)^-$ is a parabolic opposite to $\PP(\XX)$. It is also equal to:
\begin{equation}\label{defQ}
 Q^{P(X)}_v=\prod_{\check\alpha \in \mathfrak u_{P(X)}} \frac{1-q_v^{-1}e^{\check\alpha}}{1-e^{\check\alpha}}(\delta^\frac{1}{2})
\end{equation}
(the product over all roots in the unipotent radical of the parabolic dual to $\PP(\XX)$). For affine varieties the constant $c_v$ is:
\begin{equation}\label{defc}
 \frac{\prod_{\check\theta>0} (1-\sigma_{\check\theta}q_v^{-r_{\check\theta}}e^{\check\theta})}{\prod_{\check\gamma>0}(1-e^{\check\gamma})}(\delta_{P(X)}^\frac{1}{2})
\end{equation}
but for a variety which is ``Whittaker-induced'' from a homogeneous affine spherical variety $\XX'$ of a Levi subgroup, the constant $c_v$ is the same as for $\XX'$. For example, for the Whittaker model itself we have $c_v=1$. 

We define $Q^{P(X)}_v(s)$ by replacing $\delta^\frac{1}{2}$ by $\delta^{\frac{1}{2}+s}$ in (\ref{defQ}), $c_v(s)$ by replacing $\delta_{P(X)}^\frac{1}{2}$ by $\delta_{P(X)}^{\frac{1}{2}+s}$ in (\ref{defc}) and:
\begin{equation}\label{defLXs}
L^\sharp_{X,v}(\pi,s):=\frac{c_v(s)^2}{Q^{P(X)}_v(s)} \cdot \frac{\prod_{\check\gamma\in\check\Phi_X}(1-q_v^{-1-s}e^{\check\gamma})}{\prod_{\check\theta\in\Theta} (1-\sigma_{\check\theta}q_v^{-r_{\check\theta}-s}e^{\check\theta})}(\chi). \end{equation}

\begin{remark}  The purpose of introducing the parameter $s$ is to make sense of the Euler product of the $L^\sharp_{X,v}$'s as the analytic continuation of a quotient of $L$-functions. 
 (We understand the existence of such an analytic continuation as {\em part of the conjecture.}) If the pertinent (global) $L$-functions turn out to have zeroes or poles when $s=0$ the way we have chosen the $s$-parameter plays an important role; for instance, changing some occurences of $s$ by $2s$ could introduce a power of $2$ as an extra factor. (Of course, this wouldn't matter at present, since we are only formulating conjectures up to $\mathbb Q^\times$.) It appears by \cite{II} that the definitions we have given here are the correct ones. 
\end{remark}

We explain how to deduce (\ref{defLX}) from \cite{SaSph}: By definition, $\mathcal P_v^\Planch(u\otimes\bar u)$ is the quotient of ``Plancherel measure for $\XX(K_v)$'' by ``Plancherel measure for $\GG_X(K_v)$'', where by ``Plancherel measure'' we mean, as in \emph{loc.cit.\ }the Plancherel measure corresponding to Hecke eigenfunctions normalized to have value $1$ at $x_0$. Notice that since we are using Tamagawa measures, the formulas that follow will differ from those of \emph{loc.cit.\ }by a factor of $(1-q^{-1})^{-\operatorname{rk} X}$, though this factor actually doesn't play a role since it is cancelled upon division. By Theorem 9.0.1 of \emph{loc.cit.\ }the unramified Plancherel measure for $\XX(K_v)$, considered as a measure on $A_X^*/W_X$, is:
$$\frac{1}{Q^{P(X)}_{v,G} (1-q^{-1})^{\operatorname{rk} X}} L_{X,v}(\chi) d\chi$$
where we write $Q^{P(X)}_{v,G}$ to emphasize that this factor is defined with the group $\GG$ in mind. On the contrary, for $\GG_X(K_v)$ the corresponding factor will be defined with respect to the group $\GG_X\times\GG_X$, and hence we will denote it by $\left(Q^{P(X)}_{v,G_X}\right)^2$. Hence, the unramified Plancherel measure for $\GG_X(K_v)$ is:
$$\frac{1}{\left(Q^{P(X)}_{v,G_X}\right)^2 (1-q^{-1})^{\operatorname{rk} X}} L_{G_X,v}(\chi) d\chi$$
which is equal to:
$$\prod_{\check\gamma\in\check\Phi_X} \frac{1-e^{\check\gamma}}{1-q_v^{-1}e^{\check\gamma}}(\chi).$$
By the definition of $L_{X,v}$ in 7.2.3 of \emph{loc.cit.}, the claim of (\ref{defLX}) follows.

\subsection{Global conjecture} 

Let $\mathcal A_{[\psi]}$ be the subspace of the space of automorphic forms corresponding to an \emph{$X$-distinguished Arthur parameter}. The notion of $X$-distinguished, here, is the same as locally: the restriction of the parameter to the hypothetical Langlands group is bounded and factors through $\check G_X$, while the $\SL_2$-type is the $X$-distinguished $\SL_2$-type. Hence, this is the analog of the ``tempered'' hypothesis in the Ichino--Ikeda conjecture (although it does not 
imply that the automorphic representations are tempered; rather they are ``tempered relative to $X$.'') 

One can speculate about extending all that follows to the general case along the lines of \cite{II}, but we have no reason to get into that here. The space $\mathcal A_{[\psi]}$ is a unitary space; for the discrete-modulo-center space we fix norms by integrating over $[\GG/\ZZ]$ (where $\ZZ$ denotes the center) against Tamagawa measure, and we will explain how to fix norms on the continuous spectrum in \S \ref{ssexplanation}, after we give a rough statement of the conjecture.

The following should be viewed more as a working hypothesis, rather than a solid conjecture. We feel more confident about it in the case where there is at most one $X$-distinguished representation in each local $A$-packet, which should be the case if and only if the relative trace formula for $H\backslash G/H$ is stable. A good conjecture, including an understanding of the unspecified rational constants, should be the result of a theory of endoscopy for the relative trace formula.

To formulate it, let $\mathcal A'_{[\psi]}$ denote a subspace of $\mathcal A_{[\psi]}$ with the properties:
\begin{itemize}
 \item $\mathcal A_{[\psi]}'$ contains \emph{with multiplicity one} all irreducible automorphic representations which occur in $\mathcal A_{[\psi]}$;
 \item the restriction of the hermitian form $\mathcal P^\Aut$ to the orthogonal complement of $\mathcal A_{[\psi]}'$ is zero.
\end{itemize}
 Such a subspace exists by the assumption that $X$ is multiplicity-free. It is not unique, as we may arbitrarily choose its component inside isotypic subspaces of $\mathcal A_{[\psi]}$ where $\mathcal P^\Aut$ is identically zero. However, such a choice allows us to formulate the conjecture uniformly.

\begin{conjecture}[Period conjecture]\label{periodconjecture}  
Let $\psi$ denote an $X$-distinguished global Arthur parameter. For each irreducible $\nu:\pi=\otimes_v \pi_v\hookrightarrow  \mathcal A_{[\psi]}'$ there is a rational number $q$ such that:
\begin{equation}\label{periodEuler}
 \left.\mathcal P^\Aut\right|_{\nu(\pi)} = q\cdot \prod'_v \mathcal P_v^\Planch.
\end{equation}
Here $\mathcal P_v^\Planch$ are the $\HH(K_v)$-biinvariant forms on $\pi_v$ ``normalized according to $\GG_X(K_v)$-Plancherel measure'', as explained in \S \ref{sslocalpre}, and conventions for the interpretation of the Euler product and Tamagawa measures will be explained in the next subsection.
\end{conjecture}

\subsection{How to understand the Euler product}\label{ssexplanation}

Again, the Euler product (\ref{periodEuler}) should be understood in the sense of $L$-functions; namely, for every $u=\otimes u_v\in \pi=\otimes' \pi_v$ there will be a large enough set $T$ of places (including the archimedean ones and the places of ramification of $K$ over $\mathbb Q$) such that:
\begin{itemize}
 \item there is (and we fix) a smooth integral model for $\GG$ and $\XX$ outside of $T$, with $x_0\in \XX(\mathfrak o_T)$ (where $\mathfrak o_T$ denotes the ring of $T$-integers);
 \item the formula of \cite[Theorem 7.2.1]{SaSph} for eigenvectors of the spherical Hecke algebra $\mathcal H(\GG(K_v),\GG(\mathfrak o_v))$ on $\XX(K_v)$ holds  
 for $v\notin T$.
 \item $u_v\in \pi_v^{\GG(\mathfrak o_v)}$ for $v\notin T$.
\end{itemize}

Then for $v\notin T$ we have:
$$P_v^\Planch(u_v)= L^\sharp_{X,v}(\pi).$$

The equality (\ref{periodEuler}) should be thought of as a \emph{formal equality}, whose real meaning is:  
\begin{equation}
\left|  \int_{[\HH]}  \nu(u)(h) |\omega|(h) \right|^2 = q  L_X^{\sharp\,\,(T)}(\pi) \cdot \prod_{v\in T} \mathcal P_v^\Planch(u_v)
\end{equation}
where $L_X^{(T)}(\pi)$ is the value at $s=0$ of the analytic continuation of the quotient of partial  $L$-functions (outside of $T$) whose Euler factors are (\ref{defLXs}).
Again, we emphasize that the existence of such an analytic continuation should be considered as part of the conjecture. 

Here on the \emph{left hand side} we have explicitly integrated against \emph{a rational differential form} in order to make the point that one may have to divide the whole expression by some zeta factors to make the two sides finite. Tamagawa measures are, by definition, defined by taking the absolute value of invariant, rational, volume forms, times local ``convergence factors'' which are cancelled, globally, by multiplying by a special value of a partial $L$-function. However, if a group $\HH$ has nontrivial $k$-character group, then the corresponding partial $L$-function has a pole at the desired point of evaluation. This is usually resolved by multiplying, instead, by the leading coefficient of its Laurent expansion, which leads to a non-canonical but quite standard choice.
In that case, we expect that $L_X^{\sharp\,\,(T)}(\pi)$ will also have a pole of at most the same order, and should be replaced by its leading term. Equivalently, one should treat both sides of the above equation as formal Euler products and ``cancel'' the same power of the (partial) Dedekind zeta function of $k$ from both.

\begin{example}
Let $\GG=\PPGL_2$, $\HH=\AA\subset\GG$ a split torus. 
This example is the original ``Hecke integral,'' which was reinterpreted adelically in the work of \cite{JacquetLanglands} of Jacquet and Langlands. 

Then $\AA(\mathfrak o_v)=1-q_v^{-1}$ at almost every place, and therefore the ``formal'' measure of a set $S=\prod_v S_v$ with $S_v=\AA(\mathfrak o_v)$ outside of a finite set of places $T$ is:
$$|\omega|(S) = \frac{1}{\zeta_K^{(T)}(1)} \prod_{v\in T} |\omega|_v(S_v),$$
which is zero. Therefore, the period integral should be computed with respect to the measure:
$$\zeta_K^{(T)}(1)|\omega|(S):=\prod_{v\in T} |\omega|_v(S_v).$$
(Notice that it is not standard to multiply by a \emph{partial} $\zeta$-function, but in fact the conjecture is independent of how exactly one chooses to normalize the Tamagawa measure!)

On the right hand side, correspondingly, we have, outside of a finite set of places:
$$L^\sharp_{X,v}(\pi) = \frac{(1-q_v^{-1})^2}{1-q_v^{-2}} \prod_{\check\alpha\in\check \Phi_G} \frac{1-q_v^{-1}e^{\check\alpha}}{(1-q_v^{-\frac{1}{2}}e^{\frac{\check\alpha}{2}})^2}(\chi_v) = 
\frac{1-q_v^{-1}}{1-q_v^{-2}} \cdot \frac{(L_v(\pi_v,\frac{1}{2}))^2}{L_v(\pi_v,\Ad,1)}.$$
Notice that the \emph{same} globally problematic factor of $\frac{1}{\zeta_{K,v}(1)}$ appears on the right hand side, as well!

Therefore, for a cuspidal representation $\pi=\otimes' \pi_v$ the conjecture says that the period integral with respect to the measure $\zeta_K^{(T)}(1)|\omega|(S)$ is equal, up to a rational factor, to:
$$\zeta_K^{(T)}(2)\frac{(L^{(T)}(\pi,\frac{1}{2}))^2}{L^{(T)}(\pi,\Ad,1)}\cdot \prod_{v\in T} \mathcal P_v^\Planch.$$
Of course, this is known to hold, with the implicity rational factor equal to 1. We will explain the meaning of the period integral on other parts of the spectrum below.
\end{example}

\begin{remark} It is not necessary that there exists a meaningful, finite regularization of the period integral for \emph{every} representation. For example, in the case of $\AA\subset\PPGL_2$ and $\pi=1$, the trivial representation, it is reasonable to think of the ``correct'' value of the period integral (with respect to a non-zero finite measure, such as $\zeta_K^{(T)}(1)|\omega|(S)$) as being ``infinity''. This is reflected on the right hand side, as well: indeed, for the trivial representation the partial $L$-factors on the right hand side are:
$$\frac{(\zeta_K^{(T)}(1))^2}{\zeta_K^{(T)}(2)},$$
which is infinite. 
\end{remark}

Finally, we explain how to make sense of the conjecture for the continuous spectrum. Recall that we are assuming an appropriate regularization of period integrals, so we will only explain how to choose a norm on the spaces of unitary Eisenstein series. The idea here is that the Eisenstein series morphism:
$$ \mathcal E_P: \Ind_{\PP(\adele)}^{\GG(\adele)}(\delta_P^\frac{1}{2}\sigma) \to C^\infty([\GG]),$$
where $\sigma$ is a discrete automorphic representation for the pertinent Levi subgroup $L$, should be an isometry. However, on $\PP\backslash\GG$ we have again the issue of making sense, globally, of Tamagawa measures. More precisely, let $\omega$ be a $K$-rational invariant volume form on $\PP\backslash\GG$ valued in the line bundle defined by $\mathfrak d_P^{-1}$. Then, locally (having fixed good integral models outside of a finite set of places $T$, and taking $v\notin T$ such that $\sigma_v$ is unramified), we consider the induced square-norm:
\begin{equation}\label{normstandard}\int_{\PP\backslash\GG(K_v)} \|\phi_v^0(g)\|^2 |\omega|_v(g)
\end{equation}
 where the unramified vector $\phi_v^0$ is defined by choosing $u\in \sigma_v^{\LL(\mathfrak o_v)}$ with $\|u\|=1$ and setting $\phi_v^0(pk)=\delta_P^{\frac{1}{2}}\sigma(p) u$ ($p\in \PP(K_v),k\in \GG(\mathfrak o_v)$). Then one computes that this integral, using measures coming from integral, residually non-vanishing volume forms, is equal to:
$$Q_v^{P} = [\GG(\mathfrak o_v): \PP^-(\mathfrak o_v) \PP(\mathfrak o_v) ] = \prod_{\check\alpha \in \mathfrak u_{P}} \frac{1-q_v^{-1}e^{\check\alpha}}{1-e^{\check\alpha}}(\delta^\frac{1}{2}).$$

The Euler product of the $Q_v^P$'s, understood as the quotient of special values of zeta functions, is ``infinite''. Therefore, for the conjecture to make sense we need to redefine the norm (\ref{normstandard}) of the ``standard vector'' $\phi_v^0$ to be equal to $1$ outside of a finite set $T$ of places and, correspondingly, divide the Euler factors on the right hand side, for $v\notin T$, by $Q_v^{P}$; that is:

$$\mathcal P^\Aut (\mathcal E_P(\prod_{v\notin T} \phi_v^0 \cdot \prod_{v\in T} \phi_v)) = \prod'_{v\notin T} \frac{L^\sharp_{X,v}(\pi_v)}{Q_v^{P}} \cdot \prod_{v\in T} \mathcal P_v^\Planch(u_v),$$
where the partial Euler product on the right is now expected to make sense as a quotient of $L$-values.

\begin{remark}
 If $\PP$ is not a self-associate parabolic, then the variety $\YY=\UU_P\backslash \GG$ is (spherical and) multiplicity-free for the group $\LL\times\GG$, and the requirement that the Eisenstein morphism $\mathcal E_P$ be an isometry is equivalent to the validity of our conjecture for the variety $\YY$ (i.e.\ for the constant term of the Eisenstein series. If $\PP$ is self-associate then $\YY$ is not multiplicity-free, and our conjecture holds, tautologically, for the ``first summand'' of the constant term of the Eisenstein series.
\end{remark}

\subsection{Everywhere discrete or unramified}

Because of the meager state of knowledge about the Arthur conjectures in general, it is useful to discuss a specific case which can be formulated without reference to them. 
 For any representation of $\GG(\adele)$ and a large enough collection of places $T$, write $\pi_T$ for the vectors that are unramified outside $T$. Here by ``large enough'' we mean, as before:
\begin{itemize}
 \item $T$ includes the archimedean places;
 \item $K$ is unramified over $\mathbb Q$ outside of $T$;
 \item there is (and we fix) a smooth integral model for $\GG$ and $\XX$ outside of $T$, with $x_0\in \XX(\mathfrak o_T)$;
 \item the formula of \cite[Theorem 7.2.1]{SaSph} for eigenvectors of the spherical Hecke algebra $\mathcal H(\GG(K_v),\GG(\mathfrak o_v))$ on $\XX(K_v)$ holds for $v\notin T$.
\end{itemize}
``Unramified'', of course, means ``fixed by $\GG(\mathfrak o_v)$''. Following standard notation, we denote $K_T=\prod_{v\in T} K_v$, and $L^{(T)}$ a partial $L$-function outside of $T$.

The conjecture that follows is stated with the help of a modification $L^\flat_X$ of $L^\sharp_X$, which will be defined afterwards:

\begin{conjecture}[$X$-variant] \label{Xvariant}
Endow $\XX(K_T)$ with the invariant measure $\mu_T$ such that $\mu^T\cdot \mu_T = $Tamagawa measure, where $\mu^T$ is the invariant measure on $\prod_{v\notin T} \XX(K_v)$ such that $\mu^T\left(\prod_{v\notin T} \XX(\mathfrak o_v)\right)=1$. 

If $\pi \in L^2([\G])$ is irreducible and $l_T:\pi_T\hookrightarrow L^2(\X(K_T))$ is an isometric embedding (i.e.\ $\pi_T$ is an $X$-discrete series), then
for $\phi\in \pi_T$:
$$ \left| \int_{[\HH]} \phi \right|^2 \in \mathbb Q^\times  \cdot 
L_X^{\flat\,\,(T)}(\pi) \cdot |l_T(\phi)(x_0)|^2. $$
\end{conjecture}

This does not follow in an entirely routine way from Conjecture \ref{periodconjecture}, because $|\ev_{x_0}\circ l_T|^2$  differs from $\prod_{v\in T}\mathcal P^\Planch_v$ by a factor which takes into account the normalization of Plancherel measures. According to the conjectural formula (\ref{formaldegrees}), the Plancherel measure used to define $\mathcal P^\Planch_v$ for $v\in T$ is, up to a rational number, equal to an adjoint $\gamma$-factor for the group $G_X$. Since we do not want to use any functoriality assumptions in order to formulate the conjecture at this point, we will substitute the product of these $\gamma$-factors at places $v\in T$ by the inverse of the corresponding partial gamma factor away from $T$, since we expect the product of these gamma factors over all places to be equal to $1$. We also have to take into account that the measure for $v\notin T$ is normalized, here, to give mass one to $\XX(\mathfrak o_v)$, while ``local Tamagawa measure'' gives mass:
\begin{equation}\frac{Q^{P(X)}_v \cdot (1-q^{-1})^{\rk X}}{c_v}
\end{equation}
to $\XX(\mathfrak o_v)$ \cite[Theorem 9.0.3]{SaSph}. 

Therefore, formally the local unramified factor for $L_{X,v}^\flat$ is:
\begin{equation}
L_{X,v}^\flat (\pi):= \frac{Q^{P(X)}_v \cdot (1-q^{-1})^{\rk X}}{c_v} \cdot L^\sharp_{X,v}(\pi) \cdot \gamma_{\check G_X}(\pi,\Ad, 0), 
\end{equation}
where $\gamma_{\check G_X}$ denotes the adjoint gamma factor for the unramified representation $\pi$ regarded as a semisimple conjugacy class in $\check G_X$. However, the latter is zero (since it has a numerator has the factor $(1-q^s)^{\rk X}$ evaluated at zero). But again, we are only interested in making sense of $L_{X,v}^\flat$ globally (i.e.\ making sense of the partial quotient of $L$-functions $L_X^{\flat\,\,(T)}$), and the corresponding partial $\gamma$-factor $\gamma_{\check G_X}^{(T)}(\pi,\Ad, 0)$ should be finite and non-zero, since the factors for $v\in T$ are all assumed to correspond to discrete parameters. Hence, we define:
\begin{equation}
 L_{X,v}^\flat (\pi,s):= c_v(s) \cdot \frac{(1-q_v^{-s})^{\rk A_X} \prod_{\check\gamma\in\check\Phi_X}(1-q_v^{-s}e^{\check\gamma})}{\prod_{\check\theta\in\Theta} (1-\sigma_{\check\theta}q_v^{-r_{\check\theta}-s}e^{\check\theta})}(\chi),
\end{equation}
and let $L_X^{\flat\,\,(T)}(\pi)$ denote the value, at $s=0$, of the (conjectural) meromorphic continuation of:
\begin{equation}
 L_X^{\flat\,\,(T)}(\pi,s) = \prod_{v\notin T} L_{X,v}^\flat (\pi_v,s)
\end{equation}
where, again, if we have to modify the left-hand-side of the conjecture to make sense of the global Tamagawa measure, then we also have to modify $L_X^{\flat\,\,(T)}$ by the appropriate factors.

\section{Examples} \label{sec:examples}

We finally outline some examples where the period conjectures of the prior section can be verified.
We do not make any claim to originality:  many of the
results are known to experts. The material of \S \ref{ssprincipal} is related to regularization of Eisenstein periods,
a topic which has been developed in the works \cite{JLR,LapidRogawski1,LapidRogawski2}.  The results about the Whittaker period are established
already in the paper \cite{Lapid-Mao-conjecture} of Lapid and Mao, and the results of Theorem \ref{holdsbyunfolding} concern periods whose Euler factorization is already known. Thus, our main concern has been to show that the local factors are equal to the ``Plancherel'' factors predicted by the Period Conjecture \ref{periodconjecture}, thus illustrating the compatibility of known methods with the framework of this paper. In particular -- see \S \ref{FoldingCompatibility} -- the formulation of ``unfolding''
as an isometry between local $L^2$-spaces arises naturally in the evaluation of the global period.

We will use $\int^*$  to denote a regularized integral. 
This notation will often be omitted for integral expressions which depend meromorphically on a parameter in some region of convergence and are meromorphically continued to other values of the parameter; those will generally be denoted by $\int$. 

An expression of the form $\lim_{s\to 0} I(s)$, where $I(s)$ is an expression which literally makes sense only for $\Re s \gg 0$, means the value at $s=0$ of the meromorphic continuation of $I(s)$, if it exists.

For normalizations of Tamagawa measure, etc.  we refer to \S \ref{Tamagawameasure}.

 \subsection{Principal Eisenstein periods}\label{ssprincipal}

Let $\XX=\HH\backslash\GG$ be a multiplicity-free spherical variety, and let us assume, for simplicity, that $\PP(\XX)=\BB$. Assume that $\BB$ and $\HH$ have been chosen so that $\BB\HH$ is open in $\GG$; the multiplicity-free assumption, together with the assumption that $\PP(\XX)=\BB$, implies that $\HH\cap \BB$ is a torus and that there is a unique open $\BB(k_v)$-orbit on $\XX(k_v)$, for every completion $v$; moreover, that $\XX(k_v)$ is a unique $\GG(k_v)$-orbit, cf.\ \cite{SaSpc}. We will assume, as we have done throughout, that the connected component of the center of $\GG$ acts faithfully on $\XX$.

Let $\pi=I(\chi)=I_{\BB(\adele)}^{\GG(\adele)}(\chi)$ be a unitary principal series representation with $\XX$-distinguished parameter, i.e.\ the idele class character $\chi$ corresponds by class field theory to a homomorphism with bounded image: $\mathcal W_K\to A_X^*$ (where $\mathcal W_K$ denotes the Weil group of $K$). For the purpose of regularization, however, we should at first drop the requirement of ``bounded image'', i.e.\ the assumption that $\chi$ is unitary, and consider all idele class characters $\chi$ of $\BB$ which are trivial on $(\BB\cap\HH)(\adele)$. We would like to compute the (regularized) period integral of $\mathcal E_B(u)$, for every $u\in \pi$, where $\mathcal E_B$ denotes, as before, the Eisenstein series morphism.

Consider the operator:
\begin{equation}\label{Deltachi}\Delta_\chi: I(\chi)\ni u\mapsto \Phi(\HH g)= \int_{(\HH\cap\BB)\backslash\HH(\adele)} u(hg) dh \in  C^\infty(\XX(\adele))
 \end{equation}
which was called ``unnormalized Eisenstein integral'' in section \ref{sec:explicit}; we take our measures to be given by volume forms with the understanding, as was explained in \S \ref{ssexplanation}, that if they have to be modified by convergence factors to make sense of them, then the same modification will be applied to the results. The operator $\Delta_\chi$ converges absolutely when $\chi^{-1}\gg 0$ in the notation introduced after Corollary \ref{corollaryconvergence}.

It is reasonable to postulate that for almost all unitary $\chi$ the correct normalization of the integral:
$$ \int_{[\HH]} \mathcal E_B(h) dh$$
is obtained as the analytic continuation (assuming it exists) of the evaluation at $\HH \cdot 1$ of the operator $\Delta_\chi$. This is easier to justify in the case that $\HH\cap \BB$ is trivial, where the period integral of a pseudo-Eisenstein series $\Phi=\int \mathcal E_B(u_\chi) d\chi$ (where $u_\chi\in I(\chi)$ and the integral is taken over a suitable translate of the set of unitary idele class characters) over $[\HH]$ has an expression whose main term is:
$$ \int \Delta_\chi (u_\chi)(1) d\chi,$$
integrated over the same set, cf.\ \cite[Section 10]{SaSph}. In the general case the analogous expression for $\int_{[\HH]} \Phi$ is over a smaller set of characters (corresponding to $X$-distinguished Arthur parameters), so the period integral should not give a function on the space of Eisenstein series, but a distribution (or, rather, a generalized function). Thus, the value of the period in this case should be thought to be ``infinity''; this is indeed the case with $\Delta_X$, if it is globally defined by invariant volume forms (or volume forms modified by the appropriate local factors for the measure on $[\HH]$ to make sense, as mentioned before): the volume of $(\HH\cap\BB)\backslash \HH (\adele)$ is infinite. Thus, in this case our calculations should be seen as formal manipulations -- both the period and the result will be ``infinite'' but ``with the same order of $\zeta(1)$ appearing''. One could dwell on the issue of how to make a rigorous statement out of this (how to describe the period as a 
generalized function on the space of $\chi$'s, for example), but we will not get into that now (again, cf. the literature on regularized Eisenstein periods, in particular \cite{JLR, LapidRogawski1, LapidRogawski2}). 

Let us therefore explain how this matches the Period Conjecture \ref{periodconjecture}. We fix a completion $k_v$, and start denoting by regular font the points of various varieties over $k_v$. Instead of the local factor $\Delta_{\chi_v}$ of the operator $\Delta_\chi$, we might consider the adjoint:
$$C_c^\infty(X)\ni \Phi_v \mapsto \Delta_{\chi_v}^*(\Phi)(g) = \int_{B\cap H\backslash B} \Phi_v(Hbg) \chi_v\delta^{-\frac{1}{2}}(b) db \in I_B^G(\chi_v^{-1})$$
in order to show that the corresponding Hermitian form:
$$\Vert\Phi_v\Vert^2_{\chi_v}:= \Vert\Delta^*_{\chi_v}(\Phi_v)\Vert^2,$$
where the norm on the right hand side is that on $I_{B^-}^G(\chi_v^{-1})$, is the form $\mathcal P_v^\Planch$ predicted by the Period Conjecture \ref{periodconjecture}; we will recall what this means. In fact, this will not quite be the case: what we will show is that there are local factors $\gamma_v$ (depending on $\chi_v$) with the properties:

\begin{enumerate}
 \item $\gamma_v \Vert\Delta^*_{\chi_v}(\Phi_v)\Vert^2 = \mathcal P_v^\Planch(\Phi_v)$; 
 \item for almost all $v$, $\gamma_v$ can be identified with a quotient of local $L$-values, and:
 $$\prod'_v \gamma_v = 1.$$
 The product here is taken over all places, and understood as in \S \ref{ssexplanation}, i.e.\ as a partial $L$-value times a finite number of factors.
\end{enumerate}
We drop the index $v$ from $\chi$ and $\Phi$ from now on.

Recall that both the norm on $I_{B^-}^G(\chi^{-1})$ and the form $\mathcal P_v^\Planch$ depend on $k$-rational volume forms used to define measures on $X, B\backslash G$ and $G_X$ (the split group with dual $\check G_X$); we will see that these volume forms can be chosen compatibly.

First of all, fix a $k$-rational $\GG$-eigen-volume form $\omega$ on $\XX$ that will be used to define measures on the points $\XX(k_v)$ over each completion. For simplicity, let us actually assume that the form is $\GG$-invariant. We recall from Proposition \ref{propmeasure} that this induces an invariant volume form on each boundary degeneration $\XX_\Theta$, and the latter was used to fix a measure on the points $\XX_\Theta(k_v)$. Clearly, the volume form on $\XX_\Theta$ provided by Proposition \ref{propmeasure} is $k$-rational if $\omega$ is so.

Recall that $\XX_\emptyset$ is the ``most degenerate'' boundary degeneration of $\XX$. We will recall the explicit Plancherel Theorem \ref{explicitPlancherel} for the most continuous part $L^2(X)_\emptyset$ of the spectrum (where $X=\XX(k_v)$ for some completion), in a formulation that is suitable for our present purposes. The variety $\XX_\emptyset$ here is isomorphic to: $\TT\UU^-\backslash\GG$, where $\UU^-$ is the maximal unipotent subgroup of $\GG$ (taken opposite to the chosen Borel $\BB$) and $\TT$ is the subtorus of $\AA$ such that $\AA/\TT = \AA_X$. We fix such an isomorphism over $k$. We have a Plancherel decomposition for $L^2(X_\emptyset)$:
$$L^2(X_\emptyset) = \int_{\widehat{A_X}} \mathcal H_\chi \nu(\chi),$$
where $\nu(\chi)$ is in the class of Haar measure. The precise measure $\nu$ and the square of the norm on $\mathcal H_\chi$ are not canonical, of course, but their product is. 

There is an action of the little Weyl group $W_X$ on the unitary dual $\widehat{A_X}$, and in this (multiplicity-free) case it identifies almost every point of the (set-theoretic) quotient $\widehat{A_X}/W_X$ with a subset of the unitary dual $\widehat G$ of $G$. At the same time, it is identified with a subset of the unitary dual of $G_X$, the ($k_v$-points of the) split group with dual $\check G_X$. Hence we have maps, defined off a set of measure zero:
$$\widehat{A_X}/W_X \dashrightarrow \widehat G,$$
$$\widehat{A_X}/W_X \dashrightarrow \widehat{G_X}.$$

Let us fix a measurable subset $S$ of $\widehat{A_X}$ where these maps are injective; hence, we may identify $S$ as a subset of $\widehat G$. Theorem \ref{explicitPlancherel} states that the corresponding part of the Plancherel formula for $L^2(X)$ is given by:
\begin{itemize}
 \item normalized adjoint Eisenstein integrals: $E_{\emptyset,\chi}^*: C_c^\infty(X)\to \mathcal H_\chi $;
 \item the restriction of the measure $\nu(\chi)$ to $S$.
\end{itemize}

It is particularly easy in this case to describe the normalized adjoint Eisenstein integrals $E_{\emptyset,\chi}^*$: by identification of $\mathring \XX/\UU$ with $\mathring \XX_\emptyset /\UU$ over $k$ \eqref{Borbitident}, choosing corresponding $k$-points on each of them we may identify a character of $A_X$ as a function on $X/U$ or $X_\emptyset/U$. The following diagram, then, where the integrals are obtained by analytic continuation and by restriction of the fixed measures on $X$ and $X_\emptyset$ obtained from the aforementioned volume forms, should commute, cf. \eqref{Radoncommutes2}:

$$\xymatrix{
 C_c^\infty(X) \ar[dr]_{E_{\emptyset,\chi}^*} \ar[drr]^{\int_{\mathring X} \chi^{-1}\delta^{-\frac{1}{2}}}   \\
& C_c^\infty(X_\emptyset)_\chi \ar[r] & I(\chi)\\
 C_c^\infty(X_\emptyset) \ar[ur] \ar[urr]_{\int_{\mathring X_\emptyset} \chi^{-1}\delta^{-\frac{1}{2}}}
}$$

Recall that $C_c^\infty(X_\emptyset)_\chi$ denotes simply the quotient through which the lower arrow factors, which in this (multiplicity-free) case coincides with the space of smooth vectors of $\mathcal H_\chi$, for almost all $\chi$. Notice that the top arrow is the morphism $\Delta_{\chi^{-1}}^*$ that we encountered above.

Finally, we may fix $k$-rational identifications: $\XX_\emptyset = \TT\UU^-\backslash\GG$, $\mathring\XX_\emptyset = \AA_X\times\UU$, hence $\AA_X\backslash\XX_\emptyset=\BB^-\backslash\GG$. We fix a corresponding factorization of the volume form on $\mathring\XX_\emptyset$ into a product of invariant volume forms on $\AA_X,\UU$, inducing Haar measures $da, du$ on the $k_v$-points of these spaces, as well as a $\delta^{-1}$-valued measure on the quotient $A_X\backslash X$, and let $d\chi$ denote the corresponding dual measure on $\widehat{A_X}$. Then we can also identify:
$$ \mathcal H_\chi^\infty\simeq I_{B^-}^G(\chi)$$
where the quotient $C_c^\infty(X)\to I_{B^-}^G(\chi)$ is given by the integral:
$$ \Phi\mapsto \int_{A_X} \Phi(a\bullet) \chi^{-1}\delta^{-\frac{1}{2}} da$$
with norm on $I_{B^-}^G(\chi)$ obtained from the aforementioned measure on $A_X\backslash X = B^-\backslash G$:
$$ \Vert u\Vert^2=\int_{B^-\backslash G} |u^2(g)| dg.$$
The Plancherel measure corresponding to this norm is the Haar measure $d\chi$ dual to $da$. This is \emph{not} yet the Plancherel measure that we need to use by the Period Conjecture \ref{periodconjecture}, but now the adjoint normalized Eisenstein morphism can be identified with the composition of the maps:
$$\xymatrix{ E_{\emptyset,\chi}^*:  C_c^\infty(X) \ar[r]^{\Delta_\chi^*}  & I_B^G(\chi) \ar[r]^{T_\emptyset}& I_{B^-}^G(\chi)},$$
where both arrows (the second represents the ``standard'' intertwining operator) \emph{are defined using global volume forms}, and so is the norm on $I_{B^-}^G(\chi)$. Thus, we have shown that for the application of the explicit Plancherel theorem \ref{explicitPlancherel} we can use fixed volume forms defined over $k$.

If the hermitian forms $\Phi\mapsto \Vert E^*_{\emptyset,\chi}(\Phi)\Vert^2_{I_{B^-}^G(\chi)}$ correspond to Haar Plancherel measure $d\chi$ on the set $S$, then the forms $\Phi\mapsto \Vert\Delta^*_\chi(\Phi)\Vert^2_{I_B^G(\chi)}$ correspond to Plancherel measure:
$$ c(\chi)^{-1} d\chi,$$
where $c(\chi)$ was defined in \eqref{ctau}. We emphasize once more that everthing here is defined by measures obtained from globally defined, $k$-rational volume forms.

We will now compare this measure with the Plancherel measure that corresponds to the Plancherel forms $\mathcal P_v^\Planch$ of Conjecture \ref{periodconjecture}. We will see that the quotient of the two is given by scalars $\gamma_v$ with the aforementioned properties. Already, we notice that the scalars $c(\chi)$ (let us write $c_v(\chi_v)$ now to distinguish global from local) have the properties stated for the factors $\gamma_v$: their ``global Euler product'' is trivial for unitary idele class characters $\chi$. We put quotation marks here because there is no convergent Euler product, not even in a certain region for $\chi$; instead, all but finitely many factors can be interpreted as quotients of (abelian, here) $L$-factors, and we replace the infinite product by the corresponding values of $L$-functions. The statement of triviality of the global $c(\chi)$ is the statement that \emph{for unitary idele class characters, the intertwining operator $T_\emptyset$ is an isometry}, as long as global volume 
forms are used to define them and the norms on principal series.
Thus, for the purpose of factorizing the global $\HH$-period on principal Eisenstein series, \emph{there is no difference whether we use the normalized or the unnormalized Eisenstein integrals} $E_{\emptyset,\chi}$ resp.\ $\Delta_\chi$ -- or whether we use local Plancherel measures $d\chi_v$ or $c_v(\chi_v)d\chi_v$.

The Plancherel measure that corresponds to the Plancherel forms $\mathcal P_v^\Planch$ of Conjecture \ref{periodconjecture} is the restriction to $S\subset\widehat{G_X}$ of standard Plancherel measure for $G_X$. This standard Plancherel measure is the one corresponding to a Haar measure obtained by a global invariant volume form on $\GG_X$. Let us see how this Plancherel measure compares to the measure $d\chi$ that we discussed before; again, we fix a completion $k_v$ and drop the index $v$ when not necessary. Also, recall that the definition of $d\chi$ arises from a fixed invariant volume form on $\AA_X$.

The volume form on $\GG_X$ induces, again, an invariant volume form on $\GG_{X,\emptyset}$ (the most degenerate boundary degeneration of $\GG_X$). To describe the most continuous part of the Plancherel formula of $G_X$, one could use again normalized Eisenstein integrals and the measure $d\chi$ on $\widehat{A_X}$ (everything with volume forms defined over $k$, just by replacing $\XX$ in the above discussion by $\GG_X$), or matrix coefficients and the measure $c_{G_X}(\chi)^{-1} d\chi$, cf.\ Theorem \ref{groupPlancherel}. (We introduced the index $G_X$ here, in order to distinguish from the factor $c(\chi)$ above: while the $c(\chi)$ are defined by the diagram \eqref{ctau} for the group $G$, $c_{G_X}(\chi)$ is defined by the same diagram for the group $G_X$.)

The important point here is the observation made on p.\ \pageref{keepthislabel}, that in order to define the intertwining operator used to define $c_{G_X}(\chi)$ one needs to use the measure on $N^-$ (there: $U^-$) which corresponds to the chosen measure on $P\backslash G$. In particular, \emph{if we use a global, $k$-rational volume form to define local norms on the principal series, the measure on $N^-$ is also defined by a global, $k$-rational volume form}. Hence, again, \emph{the constants $c_{G_X}(\chi)$ are globally trivial}.

To summarize:
\begin{itemize}
 \item the hermitian forms $\Phi_v\mapsto \Vert\Delta^*_{\chi_v}(\Phi_v)\Vert^2_{I_B^G(\chi_v)}$ correspond to Plancherel measure $c(\chi_v)^{-1} d\chi_v$;
 \item the hermitian forms $\mathcal P^\Planch_v$ correspond to Plancherel measure $c_{G_X}(\chi_v)^{-1}d\chi_v$.
\end{itemize}

Thus: $$\gamma_v \Vert\Delta^*_{\chi_v}(\Phi_v)\Vert^2_{I_B^G(\chi_v)} = \mathcal P_v^\Planch(\Phi_v)$$
with:
$$\gamma_v = \frac{c_{G_X}(\chi_v)}{c(\chi_v)},$$
which are globally trivial, in the above sense. 

This shows that the analytic continuation of the integral \eqref{Deltachi} (evaluated at $g=1$) satisfies the Period Conjecture \ref{periodconjecture}.

\subsection{Parabolic periods}\label{ssparabolic}

The computation of period integrals in the examples that follow is based on the ``trick'' of representing the constant function on $[\HH]$ as the residue of an Eisenstein series on $[\HH]$, thus effectively replacing the $\HH$-period integral by a period integral over a parabolic subgroup of $\HH$. Therefore, we develop here the basic result that we will use.
(We do not actually prove   any instances of the conjecture in the current subsection; but this basic result will be applied in 
\S \ref{Whitsubsection} and \S \ref{FoldingCompatibility} to prove instances of it).

 \begin{quote}
Let $\HH$ be a semisimple group and $\PP$ a parabolic subgroup of $\HH$;
let $\Delta_P$ be the set of simple roots of $\HH$ belonging to the unipotent radical of $\PP$. For any automorphic function $\phi$ or rapid decay\footnote{This can be relaxed, of course} on $[\HH]$, we have: 

\begin{equation} \label{parabolicintegration} \int_{[\HH]} \phi(h) dh =  \prod_{\alpha \in \Delta_P} \langle \delta_P, \alpha^{\vee} \rangle
\lim_{s \to 0} \frac{ \int_{[\PP]} \phi(p) \delta_P^s(p)  \prod_{v}  \zeta_v(1)^{\# \Delta_P}  d_v p}{ \zeta(1+s)^{\# \Delta_P} }. \end{equation}
 where the local measures $d_v p$ are the measures defined by a right-invariant differential form on $\PP$.
 Note that the abelianization of $\PP$ has rank equal to $\# \Delta_P$, and hence the Euler product of measures $\prod \zeta_v(1)^{\# \Delta_P} d_v p$ converges to a nonzero right-invariant measure on $\PP(\adele)$. 
 As usual, $\delta_P $ denotes the modular character of $P$. The pairing $\left<\, , \,\right>$ is the canonical linear pairing between the vector spaces spanned by roots and coroots; for example, $\left<\delta_P , \check\alpha\right>=n$ when $\PP$ is the parabolic of type $\GGL_{n-1}\times\GGm$ in $\GGL_n$ and $\alpha$ is the unique simple root in its unipotent radical. Note that the integral
$\int_{[\PP]} \phi(p) \delta^s(p)$ is absolutely convergent for $s > 0$
if $\phi$ is of rapid decay\footnote{  Indeed, 
choose any  linear  algebraic representation of $\HH$ 
on a $K$-vector space $V$, and let $v_0 \in V$; then the definition
of rapid decays shows that $| \varphi(h) | \leq \| h v_0\|_{\AA}^{-N}$.
for any ``adelic norm'' on $V$. But it is possible to 
such a vector $v_0$ that is a $\PP$-eigencharacter, and moreover
the eigencharacter may be any dominant character of $\PP$; that is to say, 
$|\varphi(p)| \ll_N \| \chi(p)\|_{\AA}^{-N}$ for any dominant character $\chi$. 
That is to say, $\varphi$ decays in all directions when $|\chi| > 1$
for some dominant character $\chi$; and if $|\chi| \leq 1$
for all dominant characters, then in particular $|\delta_P| \leq 1$. }.
\end{quote}
We will denote the right-hand side of \eqref{parabolicintegration} by $\int_{[\PP]}^*$, the ``regularized integral over $\PP$,''
so that the formula asserts simply that
\begin{equation} \label{starintdef} \int_{[\HH]} \phi(h) dh =  \int_{[\PP]}^* \phi(p) . \end{equation}  

We will actually prove the statement in a somewhat more intrinsic formulation;
in particular, in place of $\zeta(1+s)$, we will use a general meromorphic function 
of the same general type to regularize. 

   Let $c(s)$ be any meromorphic function of $s$ which has a pole at $s=0$ of order 
   $\# \Delta_P$,  and for $\Re s\gg 0$ admits an Euler product decomposition: 
$\prod_v c_v(s),$
with $c_v(0)\in \CC^\times$ and the property that the Euler product of measures: 
\begin{equation} \label{dpdef} d'p := \prod_v c_v(0) dp_v\end{equation} (where $dp_v$ is the measure obtained from a $K$-rational right invariant volume form on $\PP$) is convergent.

Then the measure: 
\begin{equation} \label{dgdef} 
d'g :=\prod_v c_v(0)^{-1} dg_v\end{equation}  on $\PP\backslash\HH(\adele)$, valued in the line bundle defined by $\delta_P^{-1}$ (so that the composition of the two is Tamagawa measure on $\HH(\adele)$), will also converge.

There is another invariant integral we can define on sections of this line bundle, i.e.
on the (unnormalized) induced space $ \Ind_{\PP(\adele)}^{\HH(\adele)}(\delta_P^{1})$:
Given a section $f_0$, we extend it to a continuous section 
 $f_s\in \Ind_{\PP(\adele)}^{\HH(\adele)}(\delta_P^{1+s})$ 
 and form $ \lim_{s\to 0} c(s)^{-1}\int_{\UU^-(\adele)} f_s(u) du$
 where $du$ is Tamagawa measure on  $\UU^-(\adele)$, 
 and $\UU^-$ is the unipotent radical of a $K$-rational parabolic opposite to $\PP$.

Let $C$ be the scalar quotient of the functionals that we just defined: 
$$ C =  \frac{ f_0\mapsto \int_{\PP\backslash\GG(\adele)} f_0(g) \prod_v c_v(0)^{-1} dg_v}{f_0\mapsto \lim_{s\to 0} c(s)^{-1}\int_{\UU^-(\adele)} f_s(u) du},$$
where $f_s\in \Ind_{\PP(\adele)}^{\HH(\adele)}(\delta_P^{1+s})$ (unnormalized induction) is any continuous section that specializes to $f_0$, $\UU^-$ is the unipotent radical of a $K$-rational parabolic opposite to $\PP$ and $du$ is Tamagawa measure on $\UU^-(\adele)$.

\begin{proposition}\label{parabolic}
$$ \int_{[\HH]} \phi(h) dh = \lim_{s\to 0} C\cdot c(s)^{-1} \cdot \int_{[\PP]} \phi(p) \delta^s(p) \prod_v c_v(0) dp_v,$$
where the integrals are with respect to Tamagawa measure, and the right-hand integral is regularized
as in the prior section. 
\end{proposition}
If we take, for example, $c_v(s) = \zeta_v(1+s)^{\# \Delta_P}$,
then one may evaluate $C = \prod_{\alpha \in \Phi} \langle \delta_P, \alpha^{\vee} \rangle$, thus obtaining \eqref{parabolicintegration}.

\begin{proof}
  Let $f_s\in \Ind_{\PP(\adele)}^{\HH(\adele)}(\delta_P^{1+s})$ be any continuous section; then $\lim_{s\to 0} c(s)^{-1} \mathcal E_P(f_s)$ is a constant function on $[\HH]$,
  where $\mathcal E_P$ denotes the Eisenstein series intertwiner whose value at the identity is given by $\sum_{\PP(K) \backslash \GG(K)}$. 
  Moreover, taking
  constant term along $\UU^-$ and taking residue as $s \rightarrow 0$, 
  we deduce  $$\lim_{s\to 0} \frac{\mathcal E_P(f_s)(g)}{\int_{\UU^-(\adele)} f_s(ug) du} =1$$ for any $g$,
 and therefore:
\begin{equation} \label{a26} \int_{[\HH]} \phi(h) dh = \frac{\lim_{s\to 0} c(s)^{-1} \int_{[\HH]}\phi \mathcal E_P(f_s)}{\lim_{s\to 0} c(s)^{-1}\int_{\UU^-(\adele)} f_s(ug) du}. \end{equation}
 
The numerator is equal to:
$$  \lim_{s\to 0} c(s)^{-1} \int_{\PP(K)\backslash\HH(\adele)} \phi f_s =\int_{\PP\backslash\HH(\adele)} f_0(g) \left(  \lim_{s\to 0} c(s)^{-1}  \int_{[\PP]} \phi(pg) \delta_P^s(p) d'p \right)d'g .$$
Here we have denoted by $d'p, d'g$ the modified measures defined as in the statement of the proposition.
Now, as we just saw, this defines a functional in $f_s$ 
 which is in fact invariant under $\H(\adele)$-translation; 
that shows that the inner expression $ \lim_{s\to 0} c(s)^{-1}  \int_{[\PP]} \phi(pg) \delta_P^s(p)$ is in fact {\em constant} as a function of $g \in \HH(\adele)$.

Therefore the numerator of \eqref{a26} equals $$ \int_{\PP\backslash\HH(\adele)} f_0 (g) d'g \cdot \lim_{s\to 0} c(s)^{-1}\int_{[\PP]} \phi(p) \delta_P^s(p) d'p$$
and the claim follows.
\end{proof}

\subsection{The Whittaker case for $\GL_n$} \label{Whitsubsection}

We denote by $\PP_n$ the mirabolic subgroup of $\GG=\GGL_n$, i.e.\ the stabilizer of a vector under the standard representation.
We now repeat the argument of Jacquet \cite{JEuler} to precisely compute the absolute value of the Whittaker period; this result has appeared already in \cite{Lapid-Mao-conjecture}, although the regularization of local periods has a slightly different definition there. Notice that the constant $q$ is nontrivial here, as opposed to \cite{Lapid-Mao-conjecture}, because we are computing the norms of cusp forms as integrals over $[\PPGL_n]$, as opposed to $\GGL_n(\adele)^1$.

\begin{theorem}\label{WhittakerGLn} 
 Conjecture \ref{periodconjecture} is true for the Whittaker period of cuspidal representations of $\GL_n$, with the constant $q$ equal to $n^{-1}$.
\end{theorem}

Notice that the local Plancherel formula for the Whittaker model was established in Theorem \ref{Whittakerplancherel}, thus the forms $P_v^\Planch$ of Conjecture \ref{periodconjecture} are given by the regularized integrals of matrix coefficients of Corollary \ref{corollaryonnormalization}. This, of course, assumes that the local components of the cuspidal representation are tempered, as we have assumed the Arthur conjectures for the formulation of the Period Conjecture, but even without assuming this the following proof can give some Euler factorization by ``meromorphically continuing'' these local factors to the nontempered spectrum.  

\begin{proof} 
For simplicity, we choose a factorization of Tamagawa measure for $\adele$ (which, we recall, is define so that the measure of $\adele/K$ is $1$) into self-dual measures on $K_v$ with respect to the characters $\psi_v$.

By Proposition \ref{parabolic}, we have:
$$\int_{[\PPGL_n]} |\phi(g)|^2 dg = \int_{[\PP_n]}^* |\phi(p)|^2 dp,$$
where by a slight abuse of notation, we will also denote by $\PP_n$ the image of the mirabolic inside $\PPGL_n$;
this image is a parabolic subgroup.

For any nonzero invariant measure $d'p$ on $[\PP_n]$ and any $s$ with $\Re s\gg 0$ the integral: 
$$\int_{[\PP_n]} |\phi(p)|^2 \delta_{P_n}^s(p) d'p$$
``unfolds'' to the Whittaker model, i.e.\ is equal to:
$$\int_{\UU\backslash \PP_n(\adele)} |W_\phi(p)|^2 \delta_{P_n}^s(p) d'p,$$
where $W_\phi(g) = \int_{[\UU]} \phi(u) \psi^{-1}(u) du$. This unfolding process is a sequence of inverse Fourier transforms, and is compatible with the Tamagawa (hence self-dual with respect to the given characters) measures that we are using on the adelic points of additive groups.

Hence:
\begin{equation}\label{altWhittaker}\int_{[\PPGL_n]} |\phi|^2 = \int_{\UU\backslash \PP_n(\adele)}^* |W_\phi(p)|^2 dp,\end{equation}
where the regularization should be understood 
exactly as in \eqref{starintdef}, but with integrals over $[\PP_n]$ replaced by integrals over $\UU \backslash \PP_n (\adele)$.

Now we will write the right-hand side as an Euler product. Fix local measures $dp_v'=\zeta_v(1) dp_v$, so that their Euler product is convergent, and factorize $W_\phi=\prod_v W_{\phi_v}$. The local factors:
$$\int_{U\backslash P_n(K_v)} \left|W_{\phi_v}(p_v)\right|^2 \delta_P^s(p_v) \zeta_v(1) dp_v$$
are, by Rankin-Selberg theory, almost everywhere equal to $\zeta_v(1+s)$ times a factor whose Euler product is analytic at $s=0$. 

Therefore, by \eqref{parabolicintegration}:
$$\int_{\UU\backslash \PP_n(\adele)}^* |W_\phi(p)|^2 dp = n\cdot \lim_{s\to 0} \prod_v \int_{U\backslash P_n(K_v)} |W_{\phi_v}|^2(p_v) \delta_P^{s}(p_v) dp_v .$$
Notice that there is no factor $\zeta_v(0)$ in front of the measure on the right hand side. 

We will now see that the local Euler factor:
\begin{equation}\label{localEuler}
\int_{U\backslash P_n(K_v)} \left|W_{\phi_v}(p_v)\right|^2 dp_v, 
\end{equation}
is as predicted by Conjecture \ref{periodconjecture}, that is: the ``adjoint'' of the regularized form:
\begin{equation} \label{rak2} \int^*_{\UU(K_v)}\left<\pi_v(u) \phi_v,\phi_v\right> \psi^{-1}(u) du \end{equation}
that we constructed in \S \ref{ssWhittaker}.\footnote{Literally speaking, the regularization
constructed there made use of the fact that we were working over a nonarchimedean field --
see after Corollary \ref{corollaryonnormalization}; but it is simple to extend it to the archimedean case.} 
More precisely, we would like to show that if $\phi_v\in\pi_v$ and $\phi_v\mapsto W_{\phi_v}(1)$ is a Whittaker functional with the property that $\Vert\phi_v\Vert^2$ is given by \eqref{localEuler}, then $|W_{\phi_v}(1)|^2$ is given by \eqref{rak2}.

We show this first for the case of $\GL_2$, which is simpler. By the definition of $\int^*$ as a Fourier transform, we have for every vector $\phi_v\in \pi_v$:
$$\left<\phi_v,\phi_v\right> =  \int_{\UU\backslash\PP_n(K_v)} \int^*_{\UU(K_v)}\left<\pi(u) \phi_v,\phi_v\right> \psi(pup^{-1}) du \delta_{P_n}(p) dp.$$
Here the inner integral is regularized, and the outer integral is in fact absolutely convergent.

To see this, 
recall   that our measures are always supposed to be given by invariant differential forms, defined globally; 
notice that in the case of $\GL_2$ we have $\UU\backslash\PP_n(K_v)\simeq K_v^\times$, and the measure $\delta_{P_n}(p) dp$ can be thought of as additive measure on $K_v$, restricted to $K_v^\times$. Moreover, it is easy to see that the restriction of the matrix coefficient to $\UU(K_v)$ is $L^2$ -- in particular, its Fourier transform is a function and does not include any distribution supported at $0\in K_v$. Hence, the above equation is just duality for the Fourier transform, using the self-duality of the chosen factorization of measures.

This can be re-written as:
$$ \int_{\UU\backslash\PP_n(K_v)} \int^*_{\UU(K_v)}\left<\pi(u) \pi(p)\phi_v,\pi(p)\phi_v\right> \psi(u) du  dp.$$
If the image of $\phi_v\otimes\bar\phi_v$ under the morphism: $\pi_v\otimes\overline{\pi_v}\to C^\infty(\UU\backslash\GG(K_v), \psi) \otimes C^\infty(\UU\backslash\GG(K_v), \psi^{-1})$ defined by the regularized integral of the matrix coefficient is denoted by $W(g)\otimes \overline{W(g)}$ then the last integral can be written:
$$ \int_{\UU\backslash\PP_n(K_v)} |W(p)|^2 dp,$$
and we have finished the proof for $G=\GL_2$.

The general case will be proven by an inductive application of the above argument. To state it, let $\UU_i$ denote the unipotent radical of the parabolic corresponding to the first $i$ nodes of the Dynkin diagram (hence, with Levi $\GGL_i \times \GGm^{n-i}$), and let $\tilde\PP_i$ be the preimage, in that parabolic, of the mirabolic subgroup of $\GGL_i$ under the natural quotient map. In particular, $\UU_1=\UU=\tilde\PP_1$ and $\tilde\PP_n=\PP_n$. Denote by $\NN$ the commutator subgroup of $\UU$, and by $\NN_i$ its intersection with $\UU_i$; in particular, $\NN_1=\NN_2=\NN$. 

Denote by $f(g)$ the tempered matrix coefficient $\left<\pi(g)\phi_v(g),\phi_v(g)\right>$. We start denoting $K_v$-points of algebraic groups by regular font, as we have done in previous sections.
As in Proposition \ref{Whittakerconvergence}, one can show that $f$ is integrable over $N$; indeed, that Proposition established integrability over the kernel $H_0$ of a generic character, but one has:
$$\int_{H_0} f(h) dh = \int_{H_0/N} \int_N f(hn) dn dh,$$
which means that the inner integral is finite for almost all $h$; but it is also locally constant in $h$, which shows that it is finite for every $h$. The same argument shows that $f$ is integrable over $N_i$ for every $i$; the function $f(g) = \int_{N_i} f(n_i g) dn_i$ will be denoted by $f_{N_i}$ (typically considered as a function on $g\in U_i$ or $U_{i-1}$).

Similarly, consider the integral of $|W|^2$ over $U\backslash P_n$, where $W$ is a Whittaker function for a tempered representation. The asymptotics of Whittaker functions make it easy to see that the integral is absolutely convergent. But this integral can be written as a consecutive application of integrals:
$$ \int_{U\backslash P_n} |W|^2= \int_{\tilde P_{n-1}\backslash P_n}\int_{\tilde P_{n-2}\backslash\tilde P_{n-1}} \cdots \int_{U\backslash \tilde P_2} |W|^2,$$
and by the same argument all of the integrals are absolutely convergent.

\begin{lemma}\label{1epsilon}
 The restriction of $f$ to $U_i$ is in $L^{1+\varepsilon}$, for every $i$ and every $\varepsilon>0$.
\end{lemma}

\begin{proof}
We will use the Cartan and Iwasawa decompositions, with $K=\GGL_n(\mathfrak o_v)$. 

Recall that $f(k_1 a k_2) \ll \delta^{-\frac{1}{2}}(a)$ when $k_1, k_2\in K$ and $a \in A^+\subset A$, i.e.\  is \emph{$B$-anti-dominant}.

For an element $g\in G$ write it in terms of the Cartan and Iwasawa decompositions with respect to the \emph{opposite} Borel $B^- = A U^-$:

$$ g = k_1 a_c(g) k_2,\,\,\,\ g = u^- a_i(g) k,$$
with $a_c(g)\in A^+/A_0$, $u^- \in U^-$, $a_i(g) \in A/A_0$, where $A_0$ denotes the maximal compact subgroup of $A$. 

It is known that $\log a_i(g) - \log a_c(g)$ (the same $\log$ maps as in \eqref{logmaps}) is in the cone spanned by positive coroots. In particular, $\delta(a_i)(g) \le \delta(a_c)(g)$. Hence:

$$ \int_{U_i} |f(u)|^{1+\varepsilon} du \ll \int_{U_i} \delta^{-\frac{1+\varepsilon}{2}}(a_c(u)) du \le $$
$$ \le \int_{U_i} \delta^{-\frac{1+\varepsilon}{2}}(a_i(u)) du.$$

The last integral represents the value on the spherical vector for the standard intertwining operator:
$$ I_{B^-} (\delta^{-\frac{\varepsilon}{2}}) \to I_{B'} (\delta^{-\frac{\varepsilon}{2}}),$$
where $B'$ is the Borel subgroup obtained from $B^-$ by inverting the opposites of the roots in the Lie algebra $\mathfrak u_i$. This intertwining operator is known to converge absolutely for $\varepsilon>0$.
\end{proof}

Recall that: 
 $$|W(1)|^2 = \int^*_U \left<\pi_v(u) \phi_v,\phi_v\right> \psi^{-1}(u) du $$
 was defined as the value at $\psi$ of Fourier transform of $f_N$, the latter considered as a function on $U/N$:
$$|W(1)|^2 = \widehat{f_N}(\psi).$$
To be precise, it was defined in terms of Fourier transform on $U/H_0$ (where $H_0$ denotes the kernel of $\psi$), but since the Fourier transform of $f_N$ is locally constant on the nondegenerate locus, its value on $\psi$ coincides with the one previously defined. More generally, we will make use of the following:

\begin{lemma}\label{Fourier}
 Let $V\subset W$ be two vector spaces and $f$ a function on $W$ which is integrable over preimages of compact subsets of $W/V$, and such that its product with Lebesgue measure is a tempered distribution (in the archimedean case). Then:
 \begin{equation}\label{clear} \hat f|_{V^\perp} = \widehat{f_V}, \end{equation}
 where: 
 \begin{itemize}
  \item $f_V(w) = \int_V f(w+v) dv$, considered as a function on $W/V$;
  \item measures have been chosen compatibly on $V, W$ and $W/V$ for defining Fourier transform between (tempered generalized) functions and for the definition of $f_V$;
  \item the meaning of ``restriction'' of Fourier transform to $V^\perp\subset W^*$ is the following: Let $W=V\oplus V'$ be a decomposition, and consider an approximation of the delta measure at the identity on $(V')^\perp$ by Schwartz measures $\mu_n$ on $(V')^\perp$. Then, by definition: $$\hat f|_{V^\perp} = \lim_n \left.\left(\mu_n\star \hat f|_{V^\perp}\right)\right|_{V^\perp}$$ as tempered generalized functions, provided that this limit is independent of choices (i.e.\ the independence is part of the above assertion).
 \end{itemize}
\end{lemma}

The proof of this fact is easy and left to the reader.

Going back to the notation that we introduced above, denote by $V_i$ the quotient of $U_i$ by $U_{i+1}$. 
As in the case of $\GGL_2$, the integral of $|W|^2$ over $U\backslash \tilde P_2$ can be identified with the integral of $\widehat{f_N}$ over the subset $\psi+V_1^* \subset \widehat{(U/N)}$ -- indeed, by Lemma \ref{1epsilon}, $\widehat{f_N}$ is a locally integrable function, so again the integral over an open dense subset of $\psi+V_1^*$, which is represented by the $U\backslash \tilde P_2$ integration, is the same as the integral over the whole set.

Applying Lemma \ref{Fourier} to the function $\widehat{f_N}$ with $V$ replaced by $V_1^*$, the function $\left.\left(\widehat{f_N}\right)\right|_{V_1^*}$, which is locally constant around $\psi$, is the Fourier transform of the restriction of $f_N$ to $U_2/N$. (Notice that this restriction is smooth under translation by elements of $U$, which proves that $\widehat{f_N}$ is indeed integrable over preimages of compact subsets of $(U_2/N)^*$, as required by the Lemma.)

Now we repeat the same step, this time over the vector space $U_2/N_3$ on which the group $\tilde P_3$ (and, in fact, its normalizer -- a parabolic subgroup) acts. The orbit under this parabolic of the restriction of a nondegenerate character of $U$ is open in $(U_2/N_3)$; this shows that the Fourier transform of $\widehat{f_{N_3}}$ is locally constant around nondegenerate elements of the subspace $(U_2/N)^*$, and by applying Lemma \ref{Fourier} we see: 
$$ \widehat{f_{N_3}}(\psi) = \widehat{f_N}(\psi)$$
for such characters $\psi$. Now, the integral of $\widehat{f_{N_3}}(\psi) $ over the action of $\tilde P_2\backslash \tilde P_3$ can be identified with its integral over $\psi+$characters of $U_2/U_3$ so we get, as before, the Fourier transform of the restriction of $f_{N_3}$ to $U_3$. 

Thus, inductively, in the end we see that $\int_{U\backslash P_n} |W(p)|^2 dp = f(1) = \left<\phi_v,\phi_v\right>$. This completes the proof of the theorem.

\end{proof}

\subsection{Compatibility of the conjecture with unfolding} \label{FoldingCompatibility}

Our local interpretation of the ``unfolding'' process as an isometry between different $L^2$-spaces (Theorem \ref{unfoldingisomorphism}) allows us to prove Conjecture \ref{periodconjecture} when one period integral ``unfolds'' to another, e.g.:

\begin{theorem}\label{holdsbyunfolding}
 Conjecture \ref{periodconjecture} is true for cuspidal representations, with the given value of $q$, for the following spherical varieties:
\begin{itemize}
\item $\SSL_n^\diag\backslash \GGL_n\times\GGL_{n+1}$ under the action of $\GGm\times \GGL_n\times\GGL_{n+1}$ ($q^{-1}=n\cdot (n+1)$); 
\item $\PP_n^\diag\backslash \GGL_n\times\GGL_n$ (the classical Rankin-Selberg integral, $q^{-1}=n^2$);
\item $\SSL_n\times\PP_n\backslash\GGL_{2n}$ under the action of $\GGm\times\GGL_{2n}$ where $\GGm=\GGL_n^\ab$ (the Rankin-Selberg integral of Bump-Friedberg \cite{BF}, $q^{-1}=2n$).
 \end{itemize}
\end{theorem}

As before, $\PP_n$ denotes the mirabolic subgroup. We should clarify the meaning of the conjecture when the stabilizer $\HH$ of a point on $\XX$ is not reductive, at least in the cases above. Notice that the above examples actually correspond to periods against \emph{characters} of spherical subgroups, but with the character expressed as a character of the group; for instance, in the first case the automorphic representation is of the form $\pi=\chi \otimes \pi_1\otimes \pi_2$, where $\chi$ denotes an idele class character of $\GGm=\GGL_n^\ab$, and the period of an element of $\pi$ over $[\HH]$ is the same as the $[\GGL_n]^\diag$-period of an element of $\pi_1\otimes\pi_2$ against the character $\chi$.

Since $\XX$ is always quasi-affine, we will say that a character $\chi$ of $\GG(\adele)$ which is trivial on $\HH(\adele)$ is ``sufficiently $X$-positive'' if it is of the form:
$$ \omega = \chi\cdot | \mathfrak c|^s,$$
where $\chi$ is unitary, $\mathfrak c$ is an algebraic character vanishing on $\overline{\XX}^\aff\smallsetminus\XX$ (where $\overline{\XX}^\aff$ denotes the affine closure) and $\Re(s)$ is sufficiently large. The same notion will be applied to characters of $G(K_v)$, of course. 

We will also apply this notion to \emph{central} characters of $\GG(\adele)$, provided they are of the form (a unitary character) $\times$ (a sufficiently $X$-positive character of $\GG(\adele)$).

Now, for a unitary cuspidal representation $\pi$ of $\GG$ the integral over $[\HH]$ is not in general convergent; the reason is, of course, that its elements are not rapidly decaying, but only rapidly decaying modulo the center, and $[\HH]$ does not have finite volume. By the way, if $\HH$ is not reductive then $[\HH]$ is not necessarily closed in $[\GG]$, which is another way to see the lack of convergence. However, it is easy to show that the $[\HH]$-period is convergent on elements of $\pi\otimes\omega$, where $\omega$ is any sufficiently $X$-positive idele class character of $\GG$. Thus, given a cusp form $\phi$ with unitary central character we can interpret:
$$\int_{[\HH]} \phi(h) dh = \lim_{s\to 0} \int_{[\HH]} \phi\cdot \omega^s(h) dh,$$
where $\omega$ is a sufficiently $X$-positive idele class character and the limit denotes, as before, the value at $s=0$ of the meromorphic continuation of the given expression. The ability to continue meromorphically, of course, comes in question, but in the cases above it is not an issue as the above period integrals can be interpreted as inner products against Eisenstein series; we leave the details to the reader.

Finally, we should remind that since the invariant measure on $[\HH]$ defined by a right-invariant volume form $dh$ is not well-defined (does not correspond to a convergent Euler product) when $\HH^\ab$ has nonzero split rank, one should heed the conventions of \S \ref{ssexplanation} in order to make sense of the conjecture: both the local measures and the local Euler factors on the right hand side of \eqref{periodEuler} should be multiplied by the same local factors so that the Euler products become convergent.

\begin{proof}
Fix an invariant volume form $dh$ on $\HH$ and let $d'h_v= \zeta_v(1)^{\rk \HH^\ab}dh_v$ so that the Euler product of measures converges. In the first two examples above $\rk\HH^\ab=1$, while in the third $\rk\HH^\ab=2$. We let $d'h=\prod_v d'h_v$, a measure on $\HH(\adele)$ and $[\HH]$.

Fix a unitary cuspidal representation $\pi$ of $\GG(\adele)$, identified with a space of functions on $[\GG]$, and let $\phi\in \pi$. In all of these cases the period integral of cusp forms ``unfolds'' to the Whittaker model, i.e.:
$$ \int_{[\HH]} \phi\cdot \omega^s(h) d'h = \int_{\UU\cap\HH\backslash \HH (\adele)} \int_{[\UU]} \phi\cdot \omega^s(nh) \psi^{-1}(n) du d'h $$
for sufficiently $\XX$-positive idele class characters $\omega^s$ of $\GG$.

As in the previous proof,
choose a morphism $\phi_v\mapsto W_v$ into the Whittaker space of $\pi_v$
in such a way that: 
\begin{equation}\label{Bonn}|W_v(1)|^2 =  \int^*_{\UU(K_v)}\left<\pi_v(u) \phi_v,\phi_v\right> \psi^{-1}(u) du.\end{equation}
 Since the conjecture holds for the Whittaker model, Theorem \ref{WhittakerGLn}, we can write:
$$ \left| \int_{[\HH]} \phi\cdot \omega^s(h)  d'h \right|^2 = q 
\left| \int_{ \UU\cap\HH\backslash \HH (\adele)}  \prod_{v} W_v\cdot \omega^s_v(h) d'h \right|^2,$$
where $q$ is as stated in the theorem.

Now we can take analytic continuation of both sides to $s=0$; this is compatible with the Period Conjecture, where Euler products where interpreted by means of analytic continuation. 

Thus, it remains to verify that the local factors:
\begin{equation}\label{localfactors} \phi_v\mapsto \lim_{s\to 0}\left|\int_{U\cap H\backslash H} W_v\cdot \omega^s(h) d'h_v\right|^2 
\end{equation}
(where we started again denoting by $U, H$ etc.\ the $K_v$-points of the corresponding groups)
are equal to $\zeta_v(1)^{2\rk \HH^\ab}$ (because we modified the measures) times the ``Plancherel'' hermitian forms $\mathcal P_v^\Planch$ of \eqref{periodEuler}. Equivalently, that the same local factors with measures $dh_v$ are equal to $\mathcal P_v^\Planch$.

In other words, we need to verify that if $\Phi\in C_c^\infty(X)$ and we consider the adjoints of the maps \eqref{localfactors} as morphisms:
$$I_{\pi_v}: C_c^\infty(X)\to \pi_v$$
then we have a Plancherel formula:
\begin{equation}\label{desired}\Vert \Phi\Vert^2_{L^2(X)} = \int_{\widehat{\GG(K_v)}} \Vert I_{\pi_v}(\Phi)\Vert^2 \mu_{G_v}(\pi_v),
\end{equation}
where $\mu_{G_v}$ is Plancherel measure on $\widehat{\GG(K_v)}$. Here we remind that all measures (including the measure on ${\GG(K_v)}$ and hence the Plancherel measure on $\widehat{\GG(K_v)}$) are chosen by global volume forms, which we factorize at will. For any choice of local measures on $G$ and $U$ (and, compatibly, on $U\backslash G$), the Plancherel formula of Theorem \ref{Whittakerplancherel} holds for the Whittaker model -- i.e.\ the analog of \eqref{desired} when $L^2(X)$ is replaced by $L^2$ of the Whittaker model and $\Vert I_{\pi_v}\Vert^2$ is replaced by the adjoint of the map $\phi_v\mapsto W_v$. 

Now recall our local interpretation of ``unfolding'' in \S \ref{ssunfolding} as an isometry:
\begin{equation}\label{unf} \Unf: L^2(X)\to L^2(U\backslash G,\psi).
\end{equation}
This isometry, restricted to $C_c^\infty(X)$, was given by a series of Fourier transforms over subgroups of $U$. We claim that \emph{we can factorize global Tamagawa measures to obtain local measures on $G, U, H$ etc.\ }(and compatibly on the corresponding homogeneous spaces) \emph{so that this series of Fourier transforms does indeed give rise to an isometry \eqref{unf} when these measures are used}.

Let us use the notation preceding Theorem \ref{unfoldingisomorphism}, according to which a step of the unfolding process consists in applying Fourier transform between sections of suitable complex line bundles along the fibers of the maps: $F\to Y$ and $V^*\to Y$. We remind that $\F$ is the total space of an affine bundle over a variety $\YY$ and $\VV^*$ is the total space of its dual vector bundle; all of $\F$, $\VV^*$ and $\YY$ carry compatible, homogeneous (or almost homogeneous, for $\VV^*$) actions of our group $\GG$. In our example at the beginning of the unfolding process we have $\F\simeq \XX$, while at the end of the unfolding process we have $\VV^*\simeq$ a partial compactification of $\UU\backslash\GG$. (And, more generally, by ``folding back'' after each step we may 
assume that at every step we have a homogeneous space $\F$ for $\GG$ and a ``Whittaker-type'' space $\VV^*$ where $\GG$ acts with an open orbit $\VV^{*+}$.)

Theorem \ref{unfoldingisomorphism} holds, locally, for measures on $F, V^*$ which can be written as the composition of dual Haar measures on the fibers of $F\to Y$ and $V^*\to Y$ with a $G$-eigenmeasure on $Y$ valued in a suitable line bundle. If we take Haar measures coming from invariant differential forms on the fibers, over a fixed $k$-point, of the map: $\F(\adele)\to \YY(\adele)$ and $\VV^*(\adele)\to \YY(\adele)$ then these measures \emph{are dual to each other} by global additive duality. Notice that the measure on $\YY$ obtained from suitable eigenforms over $K$ is infinite for the adelic points of $\YY$; as a result, the resulting measure on $\F(\adele)$ does not make sense, although the measure on $\VV^{*+}(\adele)$ does.

We return to the proof of \eqref{desired}. Since the unfolding map \eqref{unf} is an isometry, we have $\Vert \Phi\Vert_{L^2(X)}^2= \Vert\Phi\Vert_{L^2(U\backslash G,\psi)}$. It is easy to see that the inverse of the unfolding map preserves compact support:
$$\Unf^{-1}: C_c^\infty(U\backslash G,\psi)\to C_c^\infty(X)$$
(in the nonarchimedean case; in the arhimedean the image will lie in the Schwartz space). Let $\mathcal U$ denote the image; for $\Phi$ lying in the image the decomposition \eqref{desired} holds.  

To show the validity for all $\Phi\in C_c^\infty(X)$ there are several ways to proceed, none of which is very pleasant: One needs to describe the correct morphisms: $C_c^\infty(X)\to \pi_v$ which restrict to the morphisms $I_{\pi_v}$ on $\mathcal U$ and appear in the Plancherel formula. First, one can appeal to a multiplicity-one statement for $\mu_{G_v}$-almost all representations, whenever it is available, or try to prove it by analyzing the precise image of $C_c^\infty(X)$ under the unfolding map. Secondly, one can try to show that the image of $C_c^\infty(X)$ under $\Unf$ is in the Harish-Chandra Schwartz space $\mathscr C(U\backslash G,\psi)$ of the Whittaker model, and hence the morphisms $I_{\pi_v}$ are the correct ones as the continuous extension of the Plancherel norms from $C_c^\infty(U\backslash G,\psi)$ to $\mathscr C(U\backslash G,\psi)$ (with respect to the topology of the latter; in particular: the integral \eqref{localfactors} is convergent for $s=0$). We will appeal to a more direct but less 
informative argument, proving directly:

\begin{lemma}
 The integral \eqref{localfactors} is convergent for $s=0$, uniformly in terms of the asymptotics of the Whittaker function.
\end{lemma}
In particular, we may approximate an element of $C_c^\infty(X)$ in the $L^2$-norm by elements of $\mathcal U$, and then the integrands on the right hand side of \eqref{desired} converge uniformly, proving the validity of the formula.

The lemma itself is an easy application of the Iwasawa decomposition, the asymptotics of Whittaker functions, and straightforward volume computations, and is left to the reader. This completes the proof of the theorem.

\end{proof}

 \newpage
 
\appendix

\makeatletter
\def\thesection{\@Alph\c@section}
\def\section{\@startsection{section}{1}
 \z@{.7\linespacing\@plus\linespacing}{.5\linespacing}
 {\normalfont\bfseries\centering\appendixname~}}
\makeatother

\section{Prime rank one spherical varieties}\label{primerankone}

\subsection{Goals} In this appendix, we compute the dual group and the commuting (Arthur) $\SL_2$
for every affine spherical variety of rank one, proving Propositions \ref{prrkone}
and \ref{rankonesl2commute}. 
  
As in the text, we use the following terminology: for a morphism $f: \SL_2 \rightarrow \check{G}$,
we shall call the restriction $f |_{\Gm}$ (or to the Lie algebra $\mathfrak g_m$) the {\em weight} of $f$.

Let us recall that the notion of a normalized spherical root is defined in \ref{ssrootdatum}: It is 
a character of $\AA_X$ that is either a root or a sum of two superstrongly orthogonal roots; it is denoted below by $\gamma$. In particular, for each spherical variety of rank one we have well-defined maps: 
$$\Gm \xrightarrow{\gamma} A_X^* \to A^*.$$
Recall also that $A_{X, GN}^*$ is, by definition, the image of $A_X^*$ in $A^*$. Moreover, the sum $2\rho_{L(X)}$ of positive roots of the Levi $\LL(X)$ defines another map ${2\rho_{L(X)}}:\Gm\to A^*$.

We shall check two assertions for rank one spherical varieties, which we term ``existence'' and ``uniqueness.''

\begin{description}
 \item[\textbf{Existence:}]
For each spherical variety $\XX$ of rank one { and normalized spherical root $\gamma$},   there exists
a morphism:
 \begin{equation} f_X \times f_A: \SL_2 \times \SL_2 \longrightarrow  \check{G}\end{equation} 
 such that the weight of $f_X$ is $\gamma$, and $f_A$ is principal into $\check{L}(X)$ with weight equal to $2\rho_{L(X)}$.

To describe the uniqueness assertion, recall that Gaitsgory and Nadler have associated to 
every affine spherical variety a group $\check{G}_{X, GN}$, which we suppose   
to satisfy axioms (GN1) -- (GN5)  from \S \ref{GNaxioms}.  In the case of rank one it is necessarily the image of a morphism: $f_{GN}:\left(A_{X,GN}^*\right)^{W_X} \times\SL_2\to \check G$ which is the identity on  $\left(A_{X,GN}^*\right)^{W_X}$ and has weight positively proportional to $\gamma$ (by (GN2)). We will show:

\item[\textbf{Uniqueness:}] Assuming (GN1) -- (GN5), the restriction of $f_{GN}$ to $\SL_2$ has weight $\gamma$.

\end{description}

\subsection{Lie algebra versions}

\label{Lie algebra remark} 

Note that both the ``Existence'' and ``Uniqueness'' assertions can be checked at the level of Lie algebras, namely:

\begin{description}
 \item[\textbf{Existence -- Lie:}] For each spherical variety $\XX$ of rank one { and normalized spherical root $\gamma$},   there exists
a morphism:
 \begin{equation} f_X \times f_A: \ssl_2 \times \ssl_2 \longrightarrow  \check{\mathfrak g}\end{equation} 
 such that the weight of $f_X$ is $\gamma$, and $f_A$ is principal into $\check{\mathfrak l}(X)$ (with weight equal to $2\rho_{L(X)}$).

 \item[\textbf{Uniqueness -- Lie:}] Let $\XX$ be an affine spherical variety of rank one, and let $f_{GN}:\left(\mathfrak a_{X,GN}^*\right)^{w_\gamma} \times\ssl_2\to \check{\mathfrak g}$ be the Gaitsgory-Nadler morphism having image equal to $\check g_{X,GN}$ and weight positively proportional to $\gamma$ (which is possible by (GN2)). Assuming (GN1)--(GN5), it actually has weight equal to $\gamma$.
\end{description}

 \subsection{Reductions for the existence statement}

We first make several reductions:

{  
\subsubsection{$\XX$ is homogeneous}

The existence statement depends only on the data $\gamma$, $\check L(X)$, and hence only on the open $\GG$-orbit (which may not be affine).	
}

\subsubsection{Parabolic induction}  \label{red1} 
 Suppose that $\PP$ is a parabolic subgroup of $\GG$, with Levi quotient $\LL$,
and $\XX_1$ is a spherical variety for $\LL$. Suppose  moreover that $\XX$ is isomorphic to $\XX_1\times^{\PP^-}\GG$.

The Levi $\check L(X)$, as well as the normalized spherical root, coincide for $\XX$ and $\XX_1$
(more precisely: they are related by means of the canonical inclusion $\check{L} \rightarrow \check{G}$) therefore the existence statement is reduced to the case where $\XX$ is not parabolically induced.

\subsubsection{Group surjection} \label{red2}

Let $\mathfrak g_1 $ be a (Lie algebra) direct summand of $\mathfrak g$ such that the corresponding normal, connected subgroup $\GG_1$ acts trivially on $\XX$. The normalized spherical root is independent of whether we consider $\XX$  as  a $\GG$-variety or as a $\GG/\GG_1$-variety, and $\check l(X) = \check l'(X)\oplus \check g_1$, where $\check l'(X)$ is the analog of $\check l(X)$ when $\XX$ is considered as a $\GG/\GG_1$-variety. (Here $\check g_1$ does not quite make sense as a subalgebra of $\check g$, but its sum with the center does.) It is clear that the existence statement is equivalent whether we are talking about $\check l(X)$ or $\check l'(X)$, which reduces the problem to the case where the action is infinitesimally faithful.

\subsubsection{Quotient by the connected component of the center}

Let $\overline{\XX}$ be the   quotient of $\XX$ by $\mathcal Z(\XX)$. The spherical roots of $\XX$ and $\overline{\XX}$ coincide, as do the associated Levi subgroups.  Therefore, the existence statement is reduced to the case { that $\mathcal Z(\XX)$ is trivial}.
 In this case, $\XX$ admits a wonderful compactification.

\subsubsection{Passage to a simply connected cover}

{ By the previous two reductions, $\GG$ is semisimple. Notice that the existence statement for Lie algebras does not depend on whether we consider $\XX$ as a spherical variety for $\GG$ or $\GG^\scon$, the simply connected cover of $\GG$. Therefore, altogether, we are reduced to the case of a pair $(\GG,\XX)$ such that:
\begin{itemize}
 \item $\GG$ is semisimple simply connected;
 \item $\XX$ is homogeneous, not parabolically induced (in any nontrivial way) and with $\mathcal Z(\XX)$;
 \item there is no direct factor of $\GG$ acting trivially.
\end{itemize}

In particular, since $\XX$ has rank one, its isotropy groups are \emph{prime} in the sense of \cite[Definition 2.3]{Wa}, and its wonderful embedding is included in \cite[Table 1]{Wa}.}

\subsubsection{Quotient by a finite automorphism group} \label{FiniteQuotient}

The normalized spherical root and the associated Levi subgroup do not change if we divide $\XX$ by a finite group of $\GG$-automorphisms. Therefore, it is enough to prove the existence statement for a class of representatives for the varieties of Table 1 of \cite{Wa} modulo the operation of taking quotients by finite subgroups of the automorphism group. We present here a table of such representatives: 

  \begin{table}[ht]
\begin{tabular}{|c| |c | c|  c|  c| c | }
\hline
 & $X = H\backslash G$ &    $\LL(\X)$ &   $\gamma$& Root type. \\
\hline
1. & $\GL_n\backslash\SL_{n+1}$ &  $ [2, n-1]$ &$\sum_1^n \alpha_i$&     T \\

\hline

2.
 &
 $\SL_2\backslash \SL_2^2$ &
$\emptyset$ &
   $\alpha_1 + \alpha_1'$  &
     G
 \\
\hline

3. & $\Sp_4\backslash \SL_4$ &  \{1,3\} & $ \alpha_1 + 2 \alpha_2 + \alpha_3$. & G \\
\hline

4. & $\Spin_{2n}\backslash \Spin_{2n+1} $ & $[2,n]$  & $\sum_{i=1}^n \alpha_i$ &   T \\ 
\hline

5. & $\Spin_{2n-1}\backslash \Spin_{2n}$ ($n \geq 3$) & $[2,n]$ & $2 \sum_1^{n-2} \alpha_i+ \alpha_{n-1}+\alpha_n$   &    G  \\ \hline

6.  & $\Spin_7\backslash \Spin_8$ &  $\{1,3\}$ & $\alpha_1+ 2 \alpha_2 + \alpha_3$ & G  \\ \hline

7. & $\SL_n \times \Gm \ltimes \wedge^2 \Ga^n \backslash$& $[2,n]$ & 
$\sum_{i=1}^n \alpha_i $&  T
\\
 & $\Spin_{2n+1}$ &  &  &  \\
 \hline
8. & $ G_2\backslash \Spin_7$ & $[1,2]$ & $ \alpha_1 + 2 \alpha_2 + 3 \alpha_3$ 
& G \\ \hline

9. &$\SL_2 \times \Sp_{2n-2}\backslash $
 &
  $\{1\} \cup [3,n]$ 
 &
  $\alpha_1 +   \alpha_n +  $&     T
    \\ 
  & $  \Sp_{2n} $  ($n \geq 2$)&  & $2 \sum_2^{n-1} \alpha_i $ &   \\ \hline

10.  & $\Spin_9\backslash F_4$ &  $[1,3]$ & $\alpha_1 + 2 \alpha_2 + 3 \alpha_3 + 2 \alpha_4$  &  T \\ \hline
11.  & $\SL_3\backslash G_2$ & $\{2\}$ & $2 \alpha_1 + \alpha_2$ &  T \\ \hline
12 .   & $ \Gm \times \SL_2 \ltimes (\Ga \oplus \Ga^2)\backslash $  & $\emptyset $ &$ \alpha_1 + \alpha_2 $ &   T \\ 
  & $G_2$&   &    & \\ 
  \hline
\end{tabular}
\end{table}

{\em Notational conventions:} Simple roots are labelled according to Bourbaki. We parameterize 
$\LL(\X)$ by using the standard numbering of simple roots and giving the set of simple roots of $\check{G}$ (= simple cooots of $\GG$) contained in it: the notation $[a,b]$ is taken to mean the set of all integers $i$ with $a \leq i \leq b$. 
(Recall that,  in the quasi-affine case, 
these are exactly the same as those simple coroots which are orthogonal to $\gamma$).

A routine verification shows that the (Existence) assertion holds for
all the varieties in the above list. 
We refer to \S \ref{existence-more} for further discussion of  this and the following result:

\begin{lemma} \label{29lemma}
In cases $(2)-(11)$, as well as case $(1)$ with $n$ odd, 
there is a unique morphism $\SL_2 \rightarrow \check{G}$
which commutes
with $2 \rho_{L(X)}$ and has weight proportional to $\gamma$. 
\end{lemma}

\subsection{Reductions for the uniqueness statement}

In the case of the uniqueness statement, we should like to repeat the same reductions; but we face the problem that, a priori, the Gaitsgory-Nadler dual group depends on the choice of affine embedding.

Thus, for example, in the case when $\XX^+$ (the open $\GG$-orbit) is parabolically induced, we may not immediately replace $\XX$ by a spherical variety for the Levi -- (GN4) only gives us some information about the image of the dual group (namely, that it lies in the corresponding Levi subgroup).

We know, by the prior discussion, that to any rank one variety $\XX$
we may find: 
\begin{itemize}
\item a parabolic subgroup $\PP^-$, with Levi quotient $\LL$;
\item a spherical $\LL$-variety $\XX_1$, such that the action of $\LL$
on $\XX$ factors through some quotient $\LL'$ whose derived group is simple (since this is, by inspection, the case of all wonderful varieties of rank one); 
\item write: $\LL'' = \mbox{the simply connected cover of the adjoint group of } \LL',$
then the  variety $\XX_1/ \mathcal{Z}(\XX_1)$, as a variety under $\LL''$, appears in the prior table (up to the operation of taking a finite quotient, cf. \S \ref{FiniteQuotient}.);
\end{itemize}
in such a way that such that the open orbit on $\XX$ is isomorphic to $\XX_1 \times^{\PP^-} \GG$.
In this way, we  regard $\XX$ as belonging to one of twelve ``types'' indexed by 
the table above.

Using notation as in \S \ref{red1}, by (GN4) the dual group of $\XX$ is contained in $\check L$, the dual Levi corresponding to the class opposite to $\PP^-$.   In fact, again by (GN4),  $\check g_{X,GN}$ lives in the dual Lie algebra of $\mathfrak l'$ (which is canonically a subalgebra of $\check{\mathfrak l}$).

  We call a spherical variety $(\GG, \XX)$ of rank one {\em good} 
  if there exists data as above and a \emph{unique} morphism $f: \ssl_2 \to \check {\mathfrak l}'$
  with weight proportional to $\gamma$ that commutes both
  with the image of $2 \rho_{L(X)}$ and  $(\mathfrak{a}_X^*)^{W_X}$.

The uniqueness statement is evidently valid for any good $(\GG, \XX)$ (because of the
existence statement, which is already proven in the prior section.) That the uniqueness statement holds for any $\XX$ of rank one follows from the two Lemmas that follow:

\begin{lemma}
Any affine spherical variety of rank one, except possibly type (1) for $n$ even are good. 
\end{lemma}
\begin{proof}
Note the following sufficient criterion for $(\GG, \XX)$ to be good: for $(\GG', \XX')$  the corresponding entry in the table of types: 
there exists a unique morphism $\ssl_2 \rightarrow \check{G'}$ 
 commuting with $2 \rho_{L(X')}$.
 In particular, in all cases but (12), the assertion follows from  Lemma \ref{29lemma}.

For case (12):  Since $\XX$ is affine, 
  its homogeneous part $\XX_1\times^{\PP^-}\GG$ must be quasi-affine. 
  This forces $\XX_1$ to be quasi-affine, also.\footnote{This follows from the fact 
  that an orbit of a linear algebraic group on an affine variety is quasi-affine.}

Now we claim that $\mathfrak{a}_X^*$, considered inside $\check{\mathfrak{l}}$,
actually projects to a full Cartan subalgebra of $\check{\mathfrak{l}}''$; the result follows easily from there. This is equivalent to saying that if $\mathfrak a_1\subset \mathfrak a_{L'}$ is the Lie algebra of the stabilizer of a generic point of $\XX_1$ in the Borel of $\LL'$ modulo its unipotent radical (and $\mathfrak a_{L'}$ denotes the universal Cartan algebra of $\LL'$), then $a_1\cap [\mathfrak l', \mathfrak l']$ is trivial. 

Notice that in this case $\XX_1$ is, up to finite quotient, the quotient of $\LL'$ (whose derived group is $\GG_2$) by $\TT\cdot \SSL_2 \cdot  (\GG_a \oplus \GG_a^2)$, where $\TT$ is a torus in $\LL'$ commuting with $\SSL_2 \cdot  (\GG_a \oplus \GG_a^2)$. From this it is easy to see that $a_1\cap [\mathfrak l', \mathfrak l']$ is the Lie algebra of $\TT\cap [\LL',\LL']$, which is at most equal to $\GGm$. But it is easy to see that $\GGm\cdot \SSL_2 \cdot  (\GG_a \oplus \GG_a^2)\backslash\GG_2$ is not quasiaffine, which implies what we want.
  \end{proof}

Finally, if $\XX$ is of type (1) we show in \S \ref{uniqueness-more}:

  \begin{lemma}\label{torusbundle} 
Suppose $\XX$ of type (1). 
Then, with the above notation, there is a good $\LL$-spherical variety $\YY_1$, in fact a torus bundle over $\XX_1$, and an affine spherical embedding $\YY$ of $\YY^+:= \YY_1\times^{\PP^-}\GG$ such that $k[\XX] = k[\YY]^\TT$ (where $\TT$ is the torus of automorphisms of that bundle).
\end{lemma}
  
Then the uniqueness statement follows, in this final case,  from (GN5).

\subsection{Further discussion of the existence result} \label{existence-more}

 We now give a few details, or at least a table of some useful data,
 related to Lemma  \ref{29lemma} and the prior existence assertion. 
  
 We describe, in the case where $G$ is a classical group, 
the representation of $\mathfrak{sl}_2  \times 
\mathfrak{sl}_2$ arising as the composite
of $f_X \times f_A$ with the classical representation of 
$\check{\mathfrak{g}}$. 
{\em We denote the $n$-dimensional irreducible representation of $\SL_2$ by $\rho_n$.}

\begin{enumerate}
\item[1.]   $\rho_2 \otimes \rho_1 \oplus \rho_1 \otimes \rho_{n-1}$.

\item[2.]$\rho_A$ is trivial  and $\rho_{GN}$ is the diagonal morphism
$\mathfrak{sl}_2 \rightarrow \mathfrak{sl}_2^2$.

\item[3.]   $\rho_2 \otimes\rho_2$. 

\item[4.]   $\rho_{1} \otimes \rho_{2n-2} \oplus \rho_2 \otimes \rho_1$.

\item[5.]  $ \rho_1 \otimes\rho_{2n-3} \oplus \rho_{3} \otimes \rho_1$. 
\item[6.] 
Related to case 5 via triality of $\mathfrak{spin}_8$. 
\item[7.] Identical to (4). 
\item[8.]  $\rho_2 \otimes \rho_3$.

\item[9.]    $ \rho_1 \otimes \rho_{2n-3} \oplus \rho_2 \otimes \rho_2$.

\item[10.
]
 In a suitable system of coordinates,
  we have $\alpha_1 = e_1-e_2, \alpha_2= e_2-e_3, \alpha_3 = e_3, \alpha_4 = \frac{-e_1-e_2-e_3-e_4}{2}$.  (The roots consist of all vectors of norm $1$ or $2$
  in $\Z^4 \cup \left(\Z^4 + \frac{1}{2} (1,1,1,1)\right)$.) 
  
  The normalized spherical root is $\gamma = -e_4$, and $\rho_{L(X)} = e_3 + 3 e_2 + 5 e_1$. 

Now $\check{\mathfrak l}(X)$ is the Levi subalgebra of $\check{\mathfrak{f}_4} = \mathfrak{f}_4$
obtained by deleting the left-most vertex of the Dynkin diagram.  
The centralizer $\mathfrak{s}$ of $2 \rho_{L(X)}$, considered as a cocharacter into $\check{\mathfrak{l}}(X)$,
has semisimple rank $1$. 
This shows that  $(F_4, \mathrm{Spin}_9\backslash F_4)$ is {\em good} in the sense previously discussed. 
But also $[\mathfrak{s}, \mathfrak{s}]$ {\em commutes with} $\check{\mathfrak{l}}(X)$. 

Taking the principal $\mathfrak{sl}_2$ inside $\check{\mathfrak{l}}(X)$
gives a morphism
$$\mathfrak{sl}_2 \times \mathfrak{sl}_2 \rightarrow [\mathfrak{s}, \mathfrak{s}] \times \check{\mathfrak{l}}(X) \rightarrow \mathfrak{f}_4,$$
which verifies uniqueness.

\item[11.]
$\gamma$ corresponds to a {\em short coroot} of $\check{\mathfrak g}$, and
the $\SL_2 \times \SL_2$ is that associated to the orthogonal pair (long root, short root).

\item[12.]   Again $\gamma$ is a short coroot.

\end{enumerate}

\subsection{Proof of Lemma \ref{torusbundle}} \label{uniqueness-more} 

We recall the situation:

$\XX$ is an affine spherical $\GG$-variety whose open orbit is isomorphic to $\XX^+\simeq \XX_1\times^{\PP^-}\GG$, where $\XX_1$ is a torus bundle over $\GGL_n\backslash\SSL_{n+1}$ (or a finite quotient thereof). By (GN4), the image of $\mathfrak{sl}_2$ under the map $(\mathfrak a_X^*)^{W_X}\times\mathfrak{sl}_2\to \check{\mathfrak g}$ is contained in a well-defined direct summand of a Levi subalgebra $\check{\mathfrak l}$ of $\check{\mathfrak g}$ which is isomorphic to $\mathfrak{pgl}_{n+1}$. The weight of such a morphism satisfying axioms (GN1) and (GN2) is \emph{not} uniquely defined by the requirement that it commutes with $(\mathfrak a_X^*)^{W_X}$ \emph{if and only if} $n$ is even and one of the following equivalent conditions hold:
\begin{itemize}
 \item $(\mathfrak a_X^*)^{W_X}$ has trivial image under the projection to the summand: $\check{\mathfrak l}\to \mathfrak{pgl}_n$;
 \item the stabilizer of a point in $\XX_1$ under the action of $\SSL_{n+1}$ contains $\GGL_n$;
 \item the valuations induced (on $k(\XX)^{(\BB)}$) by the two colors (=$\BB$-stable divisors) contained in $\XX_1\cdot \BB$ are equal.
\end{itemize}

\begin{lemma*} 
In the above setting, there is a torus bundle $\YY_1 \rightarrow \XX_1$ which does not satisfy the equivalent conditions above, and an affine spherical embedding $\YY$ of $\YY^+:= \YY_1\times^{\PP^-}\GG$ such that $k[\XX] = k[\YY]^\TT$ (where $\TT$ is the torus of automorphisms of that bundle).
\end{lemma*}

\begin{proof}
 There are clearly many possible choices for $\YY_1$. Choosing any of them, we can describe the isomorphism class of a simple spherical embedding $\YY$ of $\YY^+:= \YY_1\times^{\PP^-}\GG$ by a pair $(\mathcal C(\YY),\mathcal F(\YY))$, where $\mathcal F(\YY)$ is the set of colors of $\YY^+$ which contain the closed $\GG$-orbit in their closure and $\mathcal C(\YY)$ is the cone in $\Hom(\varchi(\YY),\QQ)$ generated by the valuations induced by all $\BB$-invariant (including the $\GG$-invariant) divisors in $\YY$ containing the closed $\GG$-orbit \cite[Theorem 3.1]{KnLV}. The embedding is affine if there is a hyperplane containing $\mathcal C(\XX)$ and strictly separating $\mathcal V_{\XX}\cup \mathcal C(\XX)$ from the set of valuations induced by colors not in $\mathcal F(\XX)$ \cite[Theorem 6.7]{KnLV}.

Since there is a bijection between colors of $\XX$ and colors of $\YY$, we may choose for $\mathcal F(\YY)$ the preimage of $\mathcal F(\XX)$. Moreover, for any extremal ray of $\mathcal C(\XX)$ which does not contain the image of a color (and hence is generated by an element of $\mathcal V_X$), we can choose a non-zero element of $\mathcal V_Y$ in its preimage, and hence obtain a cone $\mathcal C(\YY)$ generated by those and the images of elements in $\mathcal F(\YY)$. The cone is strictly convex, since $\mathcal C(\YY)$ was, and the pair $(\mathcal C(\YY),\mathcal F(\YY))$ satisfies the criterion for affinity, since the corresponding pair for $\XX$ does. Thus, we have an affine spherical variety $\YY$, and the map $\YY^+\to\XX^+$ extends to: $\YY\to \XX$ \cite[Theorem 4.1]{KnLV}.

If we decompose $k[\YY]$ into highest weight spaces, the weights that will appear are precisely the elements of $\varchi(\YY)$ which are $\ge 0$ on $\mathcal C(\YY)$ and the valuations induced by colors. Those which restrict to the trivial character for $\TT$ are precisely the elements of $\varchi(\XX)$ which are $\ge 0$ on $\mathcal C(\XX)$ and the valuations induced by colors, hence: $k[\XX]=k[\YY]^\TT$. This proves the claim.
\end{proof}

\backmatter

\bibliographystyle{alphaurl}
\bibliography{biblio}

\end{document}